\documentclass[a4paper, 10pt]{article}
\usepackage{amsmath}
\numberwithin{equation}{section}
\usepackage{amssymb,esint,hyperref}
\usepackage{amscd}
\usepackage{xspace}
\usepackage{fancyhdr}
\usepackage{color}
\usepackage{verbatim}
\usepackage{graphicx}
\usepackage{cite}
\usepackage{stmaryrd}
\usepackage{bbm}
\usepackage{xcolor}
\usepackage{mathrsfs}
\usepackage{srcltx} 

\usepackage{marginnote}
\let\oldmarginpar\marginpar
\renewcommand\marginpar[1]{\-\oldmarginpar[\raggedleft\footnotesize #1]%
	{\raggedright\footnotesize\color{red} #1}} 
\newtheorem{theorem}{Theorem}[section]

\newtheorem{definition}[theorem]{Definition}

\newtheorem{lemma}[theorem]{Lemma}

\newtheorem{proposition}[theorem]{Proposition}
\newtheorem{remark}[theorem]{Remark}

\def\da{\mathsf{Data}}
\def\B{\mathsf{B}}
\def\A{\mathsf{A}}
\def\M{\mathsf{M}}

\def\E{\mathsf{E}}
\def\F{\mathsf{F}}
\def\G{\mathsf{G}}
\def\H{\mathsf{H}}
\def\I{\mathsf{I}}
\def\J{\mathsf{J}}
\def\K{\mathsf{K}}

\def\N{\mathsf{N}}

\def\R{\mathsf{R}}

\def\T{\mathsf{T}}

\def\b{\mathbf{b}}
\def\x{\mathbf{x}}
\def\y{\mathbf{y}}

\newcommand{\eps}{\varepsilon}
\newenvironment{proof}[1][Proof]{\textbf{#1.} }{\hfill\rule{0.5em}{0.5em}}
{\catcode`\@=11\global\let\AddToReset=\@addtoreset
	\AddToReset{equation}{section}
	
	\AddToReset{theorem}{section}

	\title{Well-posedness for local and
		nonlocal quasilinear evolution equations in fluids and geometry }
	\author{
		{{\bf Ke Chen,\thanks{E-mail address: k1chen@polyu.edu.hk, Department of Applied Mathematics, The Hong Kong Polytechnic University, Kowloon, Hong Kong, PR China.}
				~~Ruilin Hu\thanks{E-mail address: huruilin16@mails.ucas.ac.cn, Academy of Mathematics and Systems Science, Chinese Academy of Sciences, Beijing, 100190, China.},
				~~Quoc-Hung Nguyen\thanks{E-mail address: qhnguyen@amss.ac.cn, Academy of Mathematics and Systems Science, Chinese Academy of Sciences, Beijing, 100190, China.}}}}
	
\begin{document}
		\maketitle
		\begin{abstract}
			We establish a Schauder-type estimate for general local and non-local linear parabolic system $$\partial_tu+\mathcal{L}_su=\Lambda^\gamma f+g$$ in $(0,\infty)\times\mathbb{R}^d$ where $\Lambda=(-\Delta)^{\frac{1}{2}}$, $0<\gamma\leq s$, $\mathcal{L}_s$ is the Pesudo-differential operator defined by  $$\mathcal{L}_su(t,x)=(2\pi)^{-\frac{d}{2}}\int_{\mathbb{R}^d}\A(t,x,\xi)\hat u(t,\xi)e^{ix\cdot\xi}d\xi,\quad\quad \A(t,x,\xi)\sim |\xi|^s.$$
			To prove this, we develop a new freezing coefficient method for kernel, where we freeze the coefficient at $x_0$, then derive a representation formula of the solution, and finally we take $x_0=x$ when estimating the solution.
			
			By applying our Schauder-type estimate to suitably chosen differential operators $\mathcal{L}_s$, we obtain critical well-posedness results of various  local and non-local nonlinear evolution equations in geometry and fluids, including hypoviscous Navier--Stokes equations, the surface quasi-geostrophic equation, mean curvature equations, Willmore flow, surface diffusion flow, Peskin equations, thin-film equations and Muskat equations.
		\end{abstract}
		\section{Introduction}\label{secintro}
		In this paper, we develop Schauder-type estimates for evolution equations. We start from the elliptic equation 
		\begin{align*}
			\Delta u=\operatorname{div} F ~~\text{in}~~\mathbb{R}^d.
		\end{align*}
		The following Schauder estimate is well-known \cite{GTbook},
		\begin{align*}
			\|\nabla u\|_{\dot{C}^{a}} \lesssim\left(\sup _{k \in \mathbb{N}} \frac{1}{|a-k|}\right)\|F\|_{\dot{C}^{a}},\quad \quad \forall a \in (0, \infty)  \backslash \mathbb{N}.
		\end{align*}
		If we consider the elliptic equation with non-constant coefficient
		$$
		\operatorname{div}(\A(x) \nabla u)=\operatorname{div}(F),
		$$
		where the coefficient matrix $\A:\mathbb{R}^d\to \mathbb{R}^{d\times d}$  satisfies $C^{-1 }\mathrm{Id}\leq \A(x)\leq C\mathrm{Id}$. Applying the freezing coefficient method (e.g. \cite{Evans98}), we obtain the corresponding Schauder estimate 
		\begin{align*}
			\|\nabla u\|_{\dot{C}^{a}} \lesssim\left(\sup _{k \in \mathbb{N}} \frac{1}{|a-k|}\right)\left(\|F\|_{\dot{C
				}^{a}}+\left(1+\|\A\|_{\dot{C}^{a}}\right)\|u\|_{L^{\infty}}\right),\ \  \forall a \in(0, \infty) \backslash \mathbb{N}.
		\end{align*}
		It is interesting to extend this estimate to   parabolic systems. A simple example is the heat equation
		\begin{align*}
			&\partial_t u-\Delta u=\operatorname{div} F ,\ \ \text{in} \ (0,\infty)\times\mathbb{R}^d ,\\
			&u|_{t=0}=u_0.
		\end{align*}
		Here $u:(0,\infty)\times \mathbb{R}^d\to\mathbb{R}$ is unknown, the force term $F:(0,\infty)\times\mathbb{R}^d\to\mathbb{R}^{d}$ is known. 
		The solution has formula
		\begin{align*}
			u(t)=e^{t\Delta}u_0+\int_0^t e^{(t-\tau)\Delta}\operatorname{div}(F(\tau))d\tau,
		\end{align*}
		where $e^{t\Delta}$ denotes the semigroup associated to the heat kernel $K(t,x)=\frac{1}{(4\pi t )^\frac{d}{2}}\exp\left(-\frac{|x|^2}{4t}\right)$. By the elementary estimates of the heat kernel
		\begin{align*}
			&\|K(t,\cdot)\|_{L^1}=1,\ \ \ \ \|\nabla_x^n K(t,\cdot)|\cdot|^\beta\|_{L^1}\lesssim t^{-\frac{n-\beta}{2}},\ \ \ n\in\mathbb{N},\ \beta >-d,
		\end{align*}
		one can obtain that 
		for any $a\in(0,1)$, 
		\begin{align}
			&\sup _{t>0}(\|u(t)\|_{L^{\infty}}+ t^{\frac{1+a}{2}}\|\nabla u(t)\|_{\dot{C}^{a}} )\lesssim \left\|u_{0}\right\|_{L^{\infty}}+\frac{1}{a(1-a)}\sup _{t>0} t^{\frac{1+a}{2}}\|F(t)\|_{\dot{C}^{a}}.\label{heat}
		\end{align}
		The constant $\frac{1}{a(1-a)}$ is optimal in the estimate. We refer to Theorem \ref{lemmain}, where we use a new freezeing coefficient method, for the proof of the estimate in a more general setting.
		The estimate \eqref{heat} is useful to establish global well-posedness for quasi-linear parabolic equations. Here we give a simple example to illustrate this. \\
		Consider the mean curvature flow system 
		\begin{equation}\label{eqmcf}
			\begin{aligned}
				&\partial_t f=\A[\nabla f]:\nabla ^2f, \ \ \ \text{in}\ (0,\infty)\times \mathbb{R}^d,\\
				&f|_{t=0}=f_0.
			\end{aligned}
		\end{equation}
		Here $f:(0,\infty)\times \mathbb{R}^d\to\mathbb{R}^N$, and 
		\begin{align}\label{defamcf}
			\A[\nabla f]=\left[\mathrm{Id}+\sum_{i=1}^{N}\nabla f^i\otimes \nabla f^i\right]^{-1}\ \ \ \ \text{and}\ \ \ \A[\nabla f]:\nabla ^2 f=\sum_{i,j}(\A[\nabla f])_{ij}\partial_i\partial_j f.
		\end{align}
		It is easy to check that the Lipschitz space $\dot W^{1,\infty}$ is critical for \eqref{eqmcf}. We consider initial data satisfying $\|\nabla f_0\|_{L^\infty}\ll 1.$ Note that we can approximate $A[\nabla f]$ by $\mathrm{Id}$ when $\|\nabla f\|_{L^\infty}\ll1$, hence the main part of the system is semi-linear. We rewrite \eqref{eqmcf} as 
		\begin{align*}
			&\partial_t f-\Delta f=(\A[\nabla f]-\mathrm{Id}):\nabla ^2f, \ \ \ \text{in}\ (0,\infty)\times \mathbb{R}^d,\\
			&f|_{t=0}=f_0.
		\end{align*}
		We prove the global well-posedness by the contraction mapping theorem. Fix $a\in(0,1)$, define the norm 
		\begin{align*}
			\interleave h\interleave=\sup_{t>0}(\|h(t)\|_{L^\infty}+t^\frac{1+a}{2}\|\nabla h(t)\|_{\dot C^a}).
		\end{align*}
		Let $\sigma>0$, denote the set
		\begin{align*}
			\mathcal{X}^\sigma=\left\{g\in L_t^\infty\dot W^{1,\infty}:g|_{t=0}=f_0, \interleave\nabla g\interleave\leq \sigma\right\}.
		\end{align*}
		It is easy to check that $	\mathcal{X}^\sigma$ is a convex closed subset of a Banach space.
		We can construct a suitable contraction map $\mathcal{S}: 	\mathcal{X}^\sigma\to 	\mathcal{X}^\sigma$ whose fixed point is a solution to \eqref{eqmcf} (see section \ref{secmeanc} for detail). 
		This implies the global well-posedness in $\mathcal{X}^\sigma$.  In fact, the above result holds for any $a\in (0,+\infty)\backslash \mathbb{N}$, then we obtain the existence of a classical solution. 
		
		When considering non-small initial data, we can not approximate the coefficient of the main term by a constant. In this case, we impose ``continuity conditions" in the form 
		\begin{align}\label{ccgene}
			\liminf_{\eps\to 0}\|u_0-u_0\ast \rho_\eps\|_{\dot W^{s,\infty}}\ll 1.
		\end{align}
		Here $\rho_\eps$ is the standard mollifier in $\mathbb{R}^d$. These conditions allow us to approximate the variable coefficients (which depend on solutions themselves and may not be regular enough) by smooth functions. This essentially indicates that the main parts of the quasi-linear equation are semi-linear. Moreover, we remark that these conditions can be relaxed to a Besov space version (i.e $B^1_{\infty,\infty}$) if the variable coefficients are well controlled (see Theorem \ref{thmPesB} for example). To obtain local well-posedness,
		we study parabolic equations with smooth variable coefficients. Consider
		\begin{align*}
			&\partial_t u-\operatorname{div}(\A(x)\nabla u) =\operatorname{div}(F),\ \ \text{in} \ (0,\infty)\times \mathbb{R}^d,\\
			&u|_{t=0}=u_0,
		\end{align*}
		where the coefficient function $\A$ is assumed to be H\"{o}lder continuous. By the fundamental observation that equations with H\"{o}lder continuous coefficients can be treated locally as a 
		perturbation of constant coefficient equations, we can use the classical freezing coefficient method (to multiply a cutoff function $\chi_{x_0}$ to the equation) to deal with coefficients that depend also on time ($\A=\A(t,x)$). 
		In this paper, we use a new freezing kernel method to approximate the solution and derive a Schauder-type estimate. Precisely, we freeze the coefficient at $x_0$, then derive a representation formula of the solution, and finally we take $x_0=x$ when estimating the solution. We note that this method is somewhat equivalent to the classical freezing coefficient method. 
		
		More generally, consider the following non-local parabolic system for $u:(0,+\infty)\times \mathbb{R}^d\to \mathbb{R}^N$.
		\begin{equation}\label{eqpara}
			\begin{aligned}
				&	\partial_{t} u(t, x)+\mathcal{L}_{s} u(t, x)=\mathcal{P}_{\gamma} f(t, x)+g(t,x) \quad \text { in } (0, \infty)\times\mathbb{R}^{d},\\
				&	u|_{t=0}=u_0,
			\end{aligned} 
		\end{equation}
		with constants $0<\gamma\leq s$. \\
		In order to explain our setup more clearly and to better state the main results, we fix some constants and impose some conditions to the operators and data in \eqref{eqpara} as follows.\vspace{0.2cm}\\
		1) We fix $m\in\mathbb{N}$ and $\kappa\in(0,s)$, $\kappa\notin\mathbb{N}$ such that $\kappa_0=\kappa-s+\gamma>0$. And we take 
		\begin{align}\label{defbet}
			\beta=100(d+m+s+1),
		\end{align}
		which is large enough to dominate some polynomial exponents encountered in our proof.\\
		2) In \eqref{eqpara}, the functions $f,g:(0,+\infty)\times \mathbb{R}^d\to \mathbb{R}^N$ are known force terms satisfying $\da_\infty(f,g)<\infty$, where
		\begin{align*}
			\da_T(f,g)=  \sup _{t \in[0,T]} \left(t^{\frac{\kappa}{s}}\|f(t)\|_{\dot{C}^{\kappa_0}(\mathbb{R}^d)}+t^{\frac{m+\kappa }{s}}\|f(t)\|_{\dot{C}^{m+\kappa_0}(\mathbb{R}^d)}+t^{\frac{m}{s}+1}\|g(t)\|_{\dot C^m(\mathbb{R}^d)}\right)
			+\|g\|_{L_T^1L^\infty(\mathbb{R}^d)}.
		\end{align*}
		3) The Pseudo-differential operators $\mathcal{L}_{s}$  of order $s>0$ is defined by
		\begin{align}\label{defop}
			&\mathcal{L}_su(t,x)=(2\pi)^{-\frac{d}{2}}\int_{\mathbb{R}^d}\A(t,x,\xi)\hat u(t,\xi)e^{ix\cdot\xi}d\xi,\quad\quad\quad
		\end{align}
		where $\A(t,x,\xi)\in \mathbb{R}^{N\times N}$ is a  symmetric matrix satisfying
		\begin{equation}\label{condop}
			\begin{aligned}
				&\quad\A(t,x,\xi)\geq c_0|\xi|^s\mathrm{Id}, \\
				&\sum_ {j,l\leq d+ m+s+2 } {|\xi|^{l-s}}{\left|\nabla _x^j\nabla^{l}_\xi \A(t,x,\xi)\right| }\leq \M,
				\ \ \forall  \xi \neq 0, (t,x)\in(0,\infty)\times\mathbb{R}^d,
			\end{aligned} 
		\end{equation}
		for some constants $c_0\in(0,1)$, and $\M>1$. Here $\mathrm{Id}$ is the $N\times N$ identity matrix. Moreover, for any  $Q_1, Q_2\in\mathbb{R}^{N\times N}$, we write  $Q_1\leq Q_2$ if $Q_2-Q_1$ is semi-positive definite. \\
		The differential operator $\mathcal{P}_{\gamma}$ is defined by 
		\begin{align}\label{defPg}
			\mathcal{P}_{\gamma} f(t,x)=\int_{\mathbb{R}^d} \B(\xi)\hat f(t,\xi)e^{ix\cdot\xi}d\xi,
		\end{align}
		with $\B(\xi)\in \mathbb{R}^{N\times N}$ satisfying 
		\begin{align}\label{proB}
			\left|\nabla^{j}_\xi \B(\xi)\right| \lesssim |\xi|^{\gamma-j},\ \ \ \forall 0\leq j\leq d+m+s+2.
		\end{align}
		4) Define the constant 
		\begin{align}\label{deftilM}
			\tilde \M=c_0^{-1}\M,
		\end{align}
		where $c_0$ and $\M$ are the constants in \eqref{condop}.
		
		Now we are in position to state the main theorems of this paper. First, if the coefficient matrix $\A$ is independent in the space variable $x$, we have the following result, which holds for any positive time $T>0$.
		\begin{theorem}\label{thmcons} Assume $\A(t,x,\xi)$ is independent of $x$.
			If $u\in C([0,T],L^\infty_x(\mathbb{R}^d))\cap L^\infty_{loc}((0,T],C^{m+\kappa }_x(\mathbb{R}^d))$ \footnote{By $u\in L^\infty_{loc}((0,T],X)$ we mean $u\in L^\infty(K,X)$ for any compact set $K\subset(0,T]$.} for $T>0$ is a solution to the Cauchy problem \eqref{eqpara},  then we have \\
			i)
			\begin{equation}\label{main1}
				\begin{aligned}
					&\sup_{t\in[0,T]} (\|u(t)\|_{L^\infty}+t^\frac{m+\kappa }{s}\|u(t)\|_{\dot C^{m+\kappa }})\lesssim\tilde \M^\beta\left(\|u_0\|_{L^\infty}+
					\da_T(f,g)\right). 
				\end{aligned}
			\end{equation}
			ii)
			Fix $0<\eta\ll 1$, then 
			\begin{equation}\label{111}
				\begin{aligned}
					&\sup_{t\in[0,T]} (t^\frac{\eta}{s}\|u(t)\|_{\dot C^\eta}+t^\frac{m+\kappa }{s}\|u(t)\|_{\dot C^{m+\kappa }})\lesssim \tilde \M^{\beta}\left(\|u_0\|_{\dot B^0_{\infty,\infty}}+  \da_T(f,g)\right). 
				\end{aligned}
			\end{equation}
		\end{theorem}
		More generally, if the coefficient matrix $\A$ also varies with respect to the space variable $x$, we use a new freezing coefficient method to obtain the following result.
		\begin{theorem}\label{lemmain}For any $T>0$, if $u\in C([0,T],L^\infty_x(\mathbb{R}^d))\cap L^\infty_{loc}((0,T],C^{m+\kappa }_x(\mathbb{R}^d))$ is a solution to the Cauchy problem \eqref{eqpara}, then \\
			i)
			\begin{equation}\label{main2}
				\begin{aligned}
					&\sup_{t\in[0,T]} (\|u(t)\|_{L^\infty}+t^\frac{m+\kappa }{s}\|u(t)\|_{\dot C^{m+\kappa }})\lesssim  e^{CT\log (T+2)}\left(\|u_0\|_{L^\infty}+  \da_T(f,g)\right),
				\end{aligned}
			\end{equation}
			where $C>0$ is a constant depending only on $\tilde \M$, $m, s$ and $\beta$.\\
			ii) 
			Let $0<\eta\ll 1$, then 
			\begin{equation}\label{112}
				\begin{aligned}
					&\sup_{t\in[0,T]} (t^\frac{\eta}{s}\|u(t)\|_{\dot C^\eta}+t^\frac{m+\kappa }{s}\|u(t)\|_{\dot C^{m+\kappa }})\lesssim e^{CT\log(T+2)}\left(\|u_0\|_{\dot B^0_{\infty,\infty}}+  \da_T(f,g)\right).
				\end{aligned}
			\end{equation}
		\end{theorem}
		\begin{remark} By the interpolation inequality in Lemma \ref{maininterpo}, \eqref{main2} implies 
			\begin{equation}\label{main2'}
				\begin{aligned}
					\sup_{0\leq \kappa'\leq m+\kappa  }\sup_{t\in[0,T]} t^\frac{\kappa'}{s}\|u(t)\|_{\dot C^{\kappa'}}\lesssim  e^{CT\log (T+2)}\left(\|u_0\|_{L^\infty}+  \da_T(f,g)\right);
				\end{aligned}
			\end{equation}
			and \eqref{112} implies 
			\begin{equation}\label{112'}
				\begin{aligned}
					\sup_{\eta\leq \eta'\leq m+\kappa }	\sup_{t\in[0,T]} t^\frac{\eta'}{s}\|u(t)\|_{\dot C^{\eta'}}\lesssim e^{CT\log(T+2)}\left(\|u_0\|_{\dot B^0_{\infty,\infty}}+  \da_T(f,g)\right).
				\end{aligned}
			\end{equation}
		\end{remark}
		\begin{remark}\label{rmk1} ~\\
			i)
			Theorem \ref{lemmain} can be applied to obtain local well-posedness of some quasilinear equations with initial data satisfying some ``continuity condition" in the form \eqref{ccgene} (which indicates that the variable coefficient can be approximated by smooth function).\\ 
			Moreover, if $\A(t,x,\xi)$ is independent of $x$, Theorem \ref{thmcons} holds, where the right hand side does not have an exponential growth in time. This will be applied to obtain global well-posedness of some quasilinear equations with small initial data (which indicates that the variable coefficient can be approximated by constant).\\
			ii) The implicit constants in \eqref{main1} and \eqref{main2} are equivalent to $\frac{1}{\kappa_0(s-\kappa)}$. The implicit constants in \eqref{111} and \eqref{112} are equivalent to $\frac{1}{\eta\kappa_0(s-\kappa)}$. We can drop this constant when $\da_T(f,g)=0$.\\
			iii) The estimates in Theorem \ref{thmcons} and Theorem \ref{lemmain} are optimal, since 
			they are   consistent with the scaling of system \eqref{eqpara}.\\
			iv) We can extend the estimates in Theorem \ref{thmcons} and Theorem \ref{lemmain} to ultra-parabolic equations, where the operator $\mathcal{L}_s$ is hypoelliptic and the system \eqref{eqpara} also exhibits smoothing effect. (see \cite{Hormand,wz1,wz2}).
		\end{remark}
		\begin{remark}
			Theorem \ref{lemmain} can be extend to $u:(0,+\infty)\times\mathcal{M}\to \mathbb{R}^N$, where $\mathcal{M}$ is a smooth manifold without boundary, see Theorem \ref{thmmani}.
		\end{remark}
		With the a priori estimates in Theorem \ref{thmcons} and Theorem \ref{lemmain} in hand, we can obtain the following existence result.
		\begin{theorem}\label{thmexis} Let $u_0\in L^\infty(\mathbb{R}^d)$, suppose that
			\begin{align*}
				\|u_0\|_{L^\infty}+  \da_T(f,g)<\infty.
			\end{align*}
			i) 
			The system \eqref{eqpara} admits a unique solution $u\in C((0,T],L^\infty_x(\mathbb{R}^d))\cap L^\infty_{loc}((0,T],C_x^{m+\kappa }(\mathbb{R}^d))$ for any $T>0$ such that \eqref{main2} and \eqref{112} hold. If the coefficient $\A(t,x,\xi)$ is independent in $x$, then \eqref{main1} and \eqref{111} hold.\vspace{0.2cm}\\
			ii) If $u_0\in L^\infty(\mathbb{R}^d)\cap C(\mathbb{R}^d)$ and $\lim_{t\to 0}  \da_t(f,g)=0$, then the solution in i) satisfies $\lim_{t\rightarrow 0}\|u(t)-u_0\|_{L_x^\infty}=0$, hence $u\in C([0,T],L^\infty_x(\mathbb{R}^d))$.
		\end{theorem}
		\begin{remark}~~\\
			i)   The solution in Theorem \ref{thmexis} \textit{i)} satisfies $(1-\Delta)^{-\varkappa}u\in C([0,T],L^\infty)$ for any $\varkappa>0$.\\
			ii) If we further assume the  initial data $u_0\in W^{-2\varkappa,\infty}$, which means that $u_0=(-\Delta)^\varkappa\tilde u_0$ with $\tilde u_0\in L^\infty$, then the solution in Theorem \ref{thmexis} \textit{i)} satisfies $u\in C([0,T],W^{-2\varkappa,\infty})$.
		\end{remark}
		Theorem \ref{thmexis} is applicable to a variety of parabolic systems. Before stating the applications, we introduce some notations that will be frequently used throughout the paper.\vspace{0.3cm}\\
		\textbf{Notations}\\
		\begin{itemize}
			\item{Fourier transform}
			\begin{equation*}
				\begin{aligned}
					&\mathcal{F}(f)(\xi)=\frac{1}{(2\pi)^{\frac{d}{2}}}\int_{\mathbb{R}^d}f(x)e^{-ix\cdot \xi}dx,\\
					&\mathcal{F}^{-1}(g)(x)=\frac{1}{(2\pi)^{\frac{d}{2}}}\int_{\mathbb{R}^d}g(\xi)e^{ix\cdot \xi}d\xi.
				\end{aligned}
			\end{equation*}
			We also write $\hat f(\xi)=\mathcal{F}(f)(\xi)$ for short. 
			It is easy to check that $\mathcal{F}\mathcal{F}^{-1}(f)=\mathcal{F}^{-1}\mathcal{F}(f)=f$.
			Moreover, 
			\begin{align*}
				\mathcal{F}(f\ast g)(\xi)=(2\pi)^\frac{d}{2}\hat f(\xi)\hat g(\xi).
			\end{align*}
			
			\item{Finite difference operators} 
			\begin{align*}
				&\delta_\alpha f(x)=f(x)-f(x-\alpha),\\&
				\Delta_\alpha f(x)=\frac{\delta_\alpha f(x)}{\alpha}~~\text{in } ~\mathbb{R},\\ &\Delta_\alpha f(x)=\frac{\delta_\alpha f(x)}{|\alpha|}~~\text{in} ~\mathbb{R}^l, l\geq 2.
			\end{align*}
			The finite difference operators appeared in this paper are all on the spatial variable, \textit{i.e.} for $h(t,x)$ defined on $[0,T]\times \mathbb{R}^d$, we denote $\delta_\alpha h(t,x)=(\delta_\alpha h(t,\cdot))(x)$. Specifically, for domain $\Omega\subset \mathbb{R}^d$, we define finite difference in $\Omega$ as
			\begin{equation*}
				\delta_\alpha f(x)=f(x)-f(x-\alpha), \quad \text{for}\ x,x-\alpha\in\Omega.
			\end{equation*}
			\item{Fractional Laplacian operators}
			\begin{equation}\label{deffracla}
				\Lambda=(-\Delta)^{\frac{1}{2}},\quad\quad\quad\Lambda^a=(-\Delta)^{\frac{a}{2}}.
			\end{equation}
			\item{Nabla operator}\\
			For $f=(f^1,\cdots f^n)$ with $f^i\in C^1(\mathbb{R}^d)$, we denote $\nabla f$ to be the matrix with
			$$
			(\nabla f)_{ij}=\partial_{i}f^j.
			$$
			\item{Multi-index derivative} \\
			For any function $f$, and any $\beta=(\beta_1,...\beta_d)\in\mathbb{N}^d$, we denote $\partial_x^\beta f=\partial_1^{\beta_1}...\partial_d^{\beta_d}f$. With a slight abuse of notation, for $m\in\mathbb{N}$, we denote $\nabla^mf=(\partial^\beta_x f)_{|\beta|=m}$, where $|\beta|=\sum_{i=1}^d\beta_i$.
			\item Frobenius inner product: for any matrix $\M, \N\in \mathbb{R}^{d\times d}$, the Frobenius inner product of $
			\M:\N$ and the induced matrix norm are defined by
			\begin{align}\label{matrxnorm}
				\M:\N=\sum_{i,j=1}^d\M_{ij}\N_{ij},\quad\quad\quad |\M|=\sqrt{\M:\M}.
			\end{align}
			\item Friedrichs mollifiers: we denote $\rho_\eps$ the standard Friedrichs mollifier with parameter $\eps>0$.
			\item For $r\in\mathbb{R}^d$, we denote the bracket $\langle r\rangle=\sqrt{1+|r|^2}$.
			\item For $s\in\mathbb{R}^+$, denote $[s]=\max\{n:n\in\mathbb{N},n\leq s\}$ the integer part of $s$.
			\item For $f:[0,T]\times \mathbb{R}^d\to \mathbb{R}^m$, we denote $\|f\|_{L^\infty(I;X)}=\sup_{t\in I}\|f(t)\|_X$, where $I\subset[0,T]$ and $X$ is a Banach space equipped with norm $\|\cdot\|_X$. Specially, when $I=[0,T]$, we simply write $\|\cdot\|_{L^\infty_TX}$.
			\item Denote $\mathbf{1}_A$ the indicator function:
			\begin{equation}\label{indicaf}
				\begin{aligned}
					\mathbf{1}_A(x):=\begin{cases}
						1,\quad\quad x\in A,\\
						0,\quad\quad x\notin A.
					\end{cases}
				\end{aligned}
			\end{equation}
			
			\vspace{0.3cm}
		\end{itemize}
		This paper is concerned with four classes of parabolic systems: semi-linear parabolic system, local quasi-linear parabolic system, nonlocal quasi-linear parabolic system, and quasi-linear parabolic system with Dirichlet boundary condition. \vspace{0.1cm}\\
		\subsection{Semi-linear parabolic system}
		We start with semi-linear parabolic systems. Let initial data $u_0$ belong to  critical Besov space $\dot{B}^a_{\infty,\infty}$ and satisfy
		\begin{align}\label{conbesov}
			\liminf_{\eps\to 0}\|u_0-u_0\ast \rho_\eps\|_{\dot B^a_{\infty,\infty}}\ll 1,
		\end{align}
		where $a$ is the critical Besov exponent that depends on the structure of the system, we prove local well-posedness results of some semi-linear parabolic systems. For initial data with small critical Besov norms, we prove global well-posedness results.\\
		\begin{remark}
			If $u_0\in (C^\infty)^{\dot B_{\infty,\infty}^a}$, which means there is a sequence of $C^\infty$ functions $\{\phi_\eps\}$ such that $\lim_{\eps\rightarrow 0}\|u_0-\phi_\eps\|_{\dot B_{\infty,\infty}^a}=0$, then \eqref{conbesov} holds naturally.
		\end{remark}
		We consider the following three systems.\\
		i) For the Aggregation-Diffusion equations
		\begin{equation}\label{eqag}
			\begin{aligned}	   &\partial_t u=\Delta u+\operatorname{div} ( u\nabla \Lambda^{-s} u),\ \ \ 0<s<1,\\
				&u|_{t=0}=u_0,
			\end{aligned}
		\end{equation}
		we prove the local and global  well-posedness in the critical Besov space $\dot B^{-s}_{\infty,\infty}$ (see Proposition \ref{propade}).\vspace{0.2cm}\\
		ii) For the hypoviscous Navier--Stokes equations
		\begin{equation}\label{eqns}
			\begin{aligned}
				&\partial_t u+\Lambda^s u+\operatorname{div}(u\otimes u)+\nabla p=0,\ \ \ \ 1<s<2,\\
				&\operatorname{div}u=0,\ \ \ \quad\quad\quad u|_{t=0}=u_0,
			\end{aligned}
		\end{equation}
		we prove the local and global  well-posedness in the critical Besov space $\dot B^{1-s}_{\infty,\infty}$ (see Proposition \ref{propns}).\\
		We remark that our method also implies the well-posedness in $L^\infty$ for $s=1$. Moreover, our method is not available for the Navier-Stokes equation $s=2$ (see Remark \ref{remNS} for details). In this case, the problem is proved to be strongly ill-posed in $\dot B^{-1}_{\infty,\infty}$ by Bourgain-Pavlovic \cite{Bourgain}, and globally well-posed for small data in $BMO^{-1}$ by Koch-Tataru\cite{KochTar}. The main reason is that $\|e^{t\Delta}u_0\|_{L^\infty}^2\sim Ct^{-1}\|u_0\|_{\dot B_{\infty,\infty}^{-1}}^2$ for small $t$, and $\|e^{t\Delta}u_0\|_{L^\infty}^2\notin L^1([0,1])$. We also refer interested readers to \cite{NSbesov} for well-posedness of \eqref{eqns} in other critical spaces close to $\dot B^{1-s}_{\infty,\infty}$.\vspace{0.2cm}\\
		iii) For the SQG equations
		\begin{equation}\label{eqsqg}
			\begin{aligned}
				&\partial_t \theta+\Lambda^s \theta+ \nabla^\perp \Lambda^{-1}\theta\cdot\nabla \theta=0,\ \ \ \ 1<s<2,\\
				&\theta|_{t=0}=\theta_0,
			\end{aligned}
		\end{equation}
		we prove the local and global well-posedness in the critical Besov space $\dot B^{1-s}_{\infty,\infty}$ (see Proposition \ref{propsqg}).\vspace{0.1cm}\\
		
		\begin{remark}
			We obtain well-posedness in Besov spaces because the nonlinear terms can be controlled by non-endpoint norms.  For example, we apply $\Lambda^{-s}$ to \eqref{eqag}, and consider the equation of $\Lambda^s u$. Then the right-hand -side nonlinear term reads $\Lambda^{-s}\operatorname{div}(u\nabla \Lambda^{-s}u)$.  By Theorem \ref{thmcons}, one needs suitable H\"{o}lder norm of $u\nabla \Lambda^{-s}u$, which can be controlled by some H\"{o}lder norms of $\Lambda ^{-s} u$, without the end-point norm $\|\Lambda^{-s}u\|_{L^\infty}$. Hence we obtain well-posedness results in Besov spaces by \eqref{111}.
		\end{remark}
		\subsection{Local quasi-linear parabolic system}
		For quasi-linear parabolic systems, we prove local and global well-posedness in critical Sobolev spaces $\dot W^{a,\infty}$. More precisely, for any initial data in $\dot W^{a,\infty}$ satisfying the continuity condition 
		\begin{align}\label{contlwp}
			\liminf_{\eps\to 0}\|u_0-u_0\ast \rho_\eps\|_{\dot W^{a,\infty}}\ll 1,
		\end{align}
		we obtain local well-posedness results. For initial data with small $\dot W^{a,\infty}$ norm:
		\begin{align}
			\label{smallgwp}
			\|u_0\|_{\dot W^{a,\infty}}\ll 1,
		\end{align}
		we prove the global well-posedness. \vspace{0.2cm}\\
		Compared with results for semi-linear systems, we consider the initial data in Sobolev spaces rather than Besov spaces, simply because we need the variable coefficients to be bounded and non-degenerate. In fact, if the boundedness and non-degeneracy of coefficients are provided, we can relax the continuity condition to a Besov version (see section \ref{secbesov} and Theorem \ref{thmPesB} for example, it is possible to improve the results of other quasi-linear systems in this paper to such a Besov version, but we will not discuss this in more detail). 
		We consider the following parabolic systems.\vspace{0.2cm}\\
		i) Consider the mean curvature flow in graph case, where
		$f:(0,T)\times\mathbb{R}^d\rightarrow\mathbb{R}^N$ satisfies
		\eqref{eqmcf}.
		We prove the local and global well-posedness in the critical space $\dot W^{1,\infty}$ (see Theorem \ref{thmmc}). \vspace{0.2cm}\\
		ii) Consider the equation of Willmore flow, where  $u:(0,T)\times \mathbb{R}^2\to \mathbb{R}$ satisfies
		\begin{equation}\label{eqwma}
			\begin{aligned}
				&u_{t}+\langle\nabla u\rangle \operatorname{div}\left(\frac{1}{\langle\nabla u\rangle}\left(\left(\mathrm{Id}-\frac{\nabla u \otimes \nabla u}{\langle\nabla u\rangle^{2}}\right) \nabla\left(\langle\nabla u\rangle H\right)-\frac{1}{2} H^{2} \nabla u\right)\right)=0 \quad \text { in } \quad (0, T)\times\mathbb{R}^{2},\\
				&u|_{t=0}=u_0,
			\end{aligned}
		\end{equation}
		with $H=\nabla\cdot \left(\frac{\nabla u}{\langle\nabla u\rangle}\right)$.
		We prove the local and global well-posedness in the critical space $\dot W^{1,\infty}$ (see Theorem \ref{thmwill}).\vspace{0.2cm} \\
		iii) Consider the equation for surface diffusion flow, where $u:[0,T]\times \mathbb{R}^d\to\mathbb{R}$ satisfies \begin{equation}
			\begin{aligned}\label{eqsurfdiffa}
				&\partial_t u+\operatorname{div }\left(\left(\langle \nabla u \rangle \mathrm{Id}-\frac{\nabla u\otimes \nabla u}{\langle \nabla u \rangle}\right)\nabla H\right)=0,\ \ \text{in}\ (0,T)\times \mathbb{R}^d,\\
				&u|_{t=0}=u_0.
			\end{aligned}
		\end{equation}
		We prove the local and global well-posedness in critical space $\dot W^{1,\infty}$ (see Theorem \ref{thmsdf}).\vspace{0.2cm}\\
		iv)
		Consider diffusion of axisymmetric surfaces, which can be written as an equation for $h:[0,T]\times \mathbb{R}\to(0,\infty)$: 
		\begin{equation}
			\begin{aligned}
				&\partial_th=\frac{1}{h}\left(\frac{h}{\langle h_x\rangle}(\mathcal{H}(h))_x\right)_x,\ \ \text{in}\ (0,T)\times\mathbb{R},\\
				&h|_{t=0}=h_0,
			\end{aligned}
		\end{equation}
		with 
		\begin{align*}
			\mathcal{H}(h)=\frac{1}{h\langle h_x\rangle}-\frac{h_{xx}}{\langle h_x\rangle^3}.
		\end{align*}
		Note that the above equation is not homogeneous. We prove the local well-posedness in the space $\dot W^{1,\infty}$. Moreover, we prove the global well-posedness when considering the equation in the periodic setting (see Theorem \ref{thmdaslo}).\vspace{0.2cm}\\
		v) 
		Consider the following gravity-driven thin film equation
		\begin{equation}\label{eqGtfa}
			\begin{aligned}
				&\partial_t \zeta+\partial_x \left((1+\varepsilon\zeta)(\partial_x \zeta+\frac{\sqrt{\mu}}{bo}\partial_x ^3\zeta)\right)+\frac{\mu}{3}\partial_x^2\left((1+\varepsilon\zeta)^3\partial_x^2\zeta\right)=0,\ \ \text{in}\ (0,T)\times\mathbb{R},\\
				&\zeta|_{t=0}=\zeta_0.
			\end{aligned}
		\end{equation}
		We prove the local well-poseness for any continuous and bounded initial data. (see Theorem \ref{thmgd} for detail).\vspace{0.2cm}\\
		vi) Consider the thin-film equation,
		\begin{equation}
			\begin{aligned}
				&\partial_tu=\Delta(e^{-\Delta u}),\quad \operatorname{in}\ (0,T)\times\mathbb{R}^d,\\
				&u|_{t=0}=u_0,\quad \operatorname{in}\ \mathbb{R}^d.
			\end{aligned}
		\end{equation}
		We prove the local well-posedness for any initial data $u_0$ satisfying $\|\Delta u_0\|_{L^\infty}<\infty$ and $\lim_{\eps\rightarrow 0}\|\Delta u_0-\Delta u_0\ast \rho_{\eps}\|_{L^\infty}\ll 1$. Moreover, for any initial data satisfying smallness condition $\|\Delta u_0\|_{L^\infty}\ll 1$, we prove the global well-posedness (see Theorem \ref{thmtf} for detail).
		\vspace{0.3cm}\\
		\subsection{Nonlocal quasi-linear parabolic system}
		The main contribution of this paper is to establish a Schauder-type estimate for general non-local 
		and non-linear parabolic equations. Non-local and non-linear parabolic operators are common in both fluid and geometry equations, which bring great difficulties to the analysis of these equations. In this paper, we prove well-posedness theory to the following non-local quasi-linear equations.\vspace{0.2cm}\\
		\textbf{Non-local mean curvature flow}\\
		Consider the non-local counterpart of the classical mean curvature flow, which is  a geometric flow driven by the fractional mean curvature. In the graph case, the system can be simplified to the following evolution equation 
		
		\begin{align*}
			&\partial_t u(t,x)=-\langle\nabla u(t,x)\rangle \mathrm{P.V.}\int_{\mathbb{R}^{d-1}}\frac{\mathcal{G}(\Delta_\alpha u(t,x))}{|\alpha|^{d-1+a}}d\alpha,\ \ \ \ \ \text{in}\ \ (0,T)\times \mathbb{R}^{d-1},\\ 
			&u|_{t=0}=u_0,
		\end{align*}
		with $a\in(0,1)$, and
		\begin{align*}
			\Delta_\alpha u(x):=\frac{u(x)-u(x-\alpha)}{|\alpha|}, \ \ \ \ \ \mathcal{G}(\rho):=\int_{-\rho}^{\rho}\frac{d\tau}{(1+\tau^2)^\frac{d+a}{2}}. \end{align*}
		We remark that $\dot W^{1,\infty}$ is a critical space of \eqref{eqnmc}. In section \ref{secnmc},  we establish both global and local well-posedness for Lipschitz initial data. The main result is the following (see section \ref{secnmc} for detail).
		\begin{theorem}\label{thmnmc}
			Fix $\kappa$ such that $0<(1+a)-\kappa\ll1$, and $m\in\mathbb{N}$. There exists $\varepsilon_0>0$ such that 	~\\
			i)	For any  initial data $u_0$ with $\| u_0\|_{\dot W^{1,\infty}}\leq \varepsilon_0$, there exists a unique global solution $u$ to \eqref{eqnmc} satisfying
			\begin{align*}
				\sup_{t>0}(\|\nabla u(t)\|_{L^\infty}+t^\frac{m+\kappa }{1+a}\|\nabla u(t)\|_{\dot C^{m+\kappa }})\leq C \| \nabla u_0\|_{L^\infty}.
			\end{align*}
			ii) For any initial data $u_0\in \dot W^{1,\infty}$, if there exists $\eps_1>0$ such that 
			\begin{align}\label{nmccon}
				\|u_0-u_0\ast \rho_{\eps_1}\|_{\dot W^{1,\infty}}\leq \varepsilon_0,
			\end{align}
			then   there exists $T=T(\eps_0,\eps_1,\|u_0\|_{\dot W^{1,\infty}})>0$ such that \eqref{eqnmc} admits a unique solution  $u$ satisfying
			\begin{align*}
				\sup_{t\in[0,T]}(\|\nabla u(t)\|_{L^\infty}+t^\frac{m+\kappa }{1+a}\|\nabla u(t)\|_{\dot C^{m+\kappa }})\leq C \| \nabla u_0\|_{L^\infty}.
			\end{align*}
			for some $C>0$.
		\end{theorem}
		\textbf{The Peskin Problem}\\
		Consider a 1D elastic string/a 2D elastic membrane immersed in a 2D/3D Stokes fluid, which is known as the 2D/3D Peskin problem. We refer readers to section \ref{secpeskin} for physical background and related results for this model. 
		The model can be written as the following boundary integral form 
		\begin{itemize}
			\item  \textbf{2D Peskin}: The unknown function $X:[0,T]\times \mathbb{S}\to \mathbb{R}^2$ satisfies
			\begin{equation}
				\begin{aligned}\label{peskin}
					&\partial_{t} X(x) =\int_{\mathbb{S}} \mathbf{G}_2(X(x)-X(\sigma)) \partial_x \left(\mathcal{T}(|\partial_x X|)\frac{\partial_x X}{|\partial_x X|}\right)(\sigma) \mathrm{d} \sigma, \\
					&\mathbf{G}_2(x)=\frac{1}{4\pi}\left(-\log |x|         \mathrm{Id}_2+\frac{x\otimes x}{|x|^2}\right),\ \ x\in\mathbb{R}^2\backslash\{0\}\\
					&X|_{t=0}=X_0.
				\end{aligned}
			\end{equation}
			\item  \textbf{3D Peskin}: The unknown function $X:[0,T]\times \mathbb{S}^2\to \mathbb{R}^3$ satisfies
			\begin{equation}\label{eqpes3d}
				\begin{aligned}
					& \frac{\partial {X}}{\partial t}(\widehat{\boldsymbol{x}}) =\int_{\mathbb{S}^2} \mathbf{G}_3({X}(\widehat{\boldsymbol{x}})-{X}(\widehat{\boldsymbol{y}})) \nabla_{\mathbb{S}^2} \cdot\left(\mathcal{T}\left(\left|\nabla_{\mathbb{S}^2} {X}\right|\right) \frac{\nabla_{\mathbb{S}^2} {X}}{|\nabla_{\mathbb{S}^2} {X}|}\right)(\widehat{\boldsymbol{y}}) d \mu_{\mathbb{S}^2}(\widehat{\boldsymbol{y}}),\\
					&\mathbf{G}_3(x)=\frac{1}{8 \pi}\left(\frac{1}{|{x}|} \mathrm{Id}_3+\frac{{x} \otimes {x}}{|{x}|^3}\right),\ \ x\in\mathbb{R}^3\backslash\{0\},\\
					&X|_{t=0}=X_0,
				\end{aligned}
			\end{equation}
			where $\nabla_{\mathbb{S}^2}$ denotes the surface gradient operator on the unit sphere, $\mu_{\mathbb{S}^2}$ is the standard measure on the unit sphere, and $\mathcal{T}:\mathbb{R}^+\to\mathbb{R}^+$ is the smooth elastic tension.
		\end{itemize}
		
		In this paper, we will apply Theorem \ref{lemmain} to obtain local well-posedness of \eqref{peskin} in the critical Lipschitz space. To state the main result, we first introduce the notion of well-stretched condition. For 2D Peskin problem \eqref{peskin}, define 
		\begin{equation}\label{kap}
			\mathbf{\Theta}_Y=\sup_{x_1,x_2\in \mathbb{S},x_1\neq x_2}\frac{|x_1-x_2|}{|Y(x_1)-Y(x_2)|},
		\end{equation}
		where $|x_1-x_2|=\inf_{k\in \mathbb{Z}}|x_1-x_2-2k\pi|$ is the distance between $x_1$ and $x_2$ on the torus. For simplicity, we denote
		\begin{equation*}
			\mathbf{\Theta}(t)=	\sup_{\tau\in[0,t]}\mathbf{\Theta}_{X(\tau)},\quad \mathbf{\Theta}_0=	\mathbf{\Theta}_{X_0}.
		\end{equation*}
		We say the initial data $X_0$ satisfies the well-stretched condition if $	\mathbf{\Theta}_0<\infty$. 
		
		Analogously, we also impose the well-stretched condition for the 3D Peskin problem:
		\begin{align}\label{well3D}
			\mathbf{\Theta}_0:=\sup_{\widehat{\boldsymbol{x}},\widehat{\boldsymbol{y}}\in\mathbb{S}^2,\widehat{\boldsymbol{x}}\neq\widehat{\boldsymbol{y}}}\frac{|\widehat{\boldsymbol{x}}-\widehat{\boldsymbol{y}}|}{|{X}_0(\widehat{\boldsymbol{x}})-{X}_0(\widehat{\boldsymbol{y}})|}<+\infty,
		\end{align}
		where $|\widehat{\boldsymbol{x}}-\widehat{\boldsymbol{y}}|$ means the standard metric in $\mathbb{R}^3$.
		\begin{remark}
			The well-stretched conditions are necessary, since otherwise the integral in the right hand of \eqref{peskin} and \eqref{eqpes3d} are not well-defined.
		\end{remark}
		The main results are the following theorems.
		\begin{theorem}\label{thmPes} Consider $X_0\in \dot W^{1,\infty}$ with $\mathbf{\Theta}_0<\infty$. Fix $\kappa\in(0,1)$ satisfying $1-\kappa\ll 1$, and $m\in\mathbb{N}^+$.
			There exists $0<\varepsilon_0\ll 1$ such that  if 
			\begin{align}\label{pescon}
				\| X_0-X_0\ast \rho_{\eps_1}\|_{\dot W^{1,\infty}}\leq \varepsilon_0, 
			\end{align} for some $\eps_1>0$,  then there exists $T>0$ such that the 2D Peskin problem \eqref{peskin} admits a unique solution  $X$ satisfying
			\begin{align*}
				&\sup_{t\in[0,T]}	(\|\partial_x X(t)\|_{L^\infty}+t^{m+\kappa}\|\partial_x X(t)\|_{\dot C^{m+\kappa}})\leq C \|\partial_x X_0\|_{L^\infty}.\\
				&\mathbf{\Theta}(T)\leq 2\mathbf{\Theta}_0.
			\end{align*}
		\end{theorem}
		
		\begin{theorem}\label{thmPes3d}
			Fix $\kappa$ such that $0<1-\kappa\ll 1$. There exists $0<\varepsilon_0\ll 1$ such that, if the initial data ${X}_0\in \dot W^{1,\infty}(\mathbb{S}^2)$ satisfying $\mathbf{\Theta}_0<\infty$ and there exists smooth $\Phi\in C^\infty(\mathbb{S}^2)$, such that
			\begin{align}\label{con3dpes}
				\|X_0- \Phi\|_{\dot W^{1,\infty}(\mathbb{S}^2)}\leq \varepsilon_0,
			\end{align}
			then \eqref{eqpes3d} has a unique solution ${X}$ in $[0,T]$ satisfying
			\begin{align*}
				&\sup_{t\in[0,T]}\left(\| X(t)\|_{W^{1,\infty}(\mathbb{S}^2)  }+t^{m+\kappa}\|X(t)\|_{C^{m+1+\kappa}(\mathbb{S}^2)}\right)\leq C\|X_0\|_{W^{1,\infty}(\mathbb{S}^2)},\\
				&\mathbf{\Theta}(T)\leq 2\mathbf{\Theta}_0.
			\end{align*}
		\end{theorem}
		The H\"{o}lder norm on $\mathbb{S}^2$ is defined in \eqref{normsph}.
		
		We remark that the condition \eqref{pescon} in Theorem \ref{thmPes} can be relaxed to a Besov version if the boundedness and non-degeneracy of variable coefficients are provided, see section \ref{secbesov} for detail. More precisely, we can prove the following theorem.
		\begin{theorem}\label{thmPesB} 
			For any $X_0\in \dot W^{1,\infty}\cap (C^2)^{\dot B^1_{\infty,\infty}}$ with \begin{align}\label{con111}
				\liminf_{\vartheta\to 0} \mathbf{\Theta}_{X_0\ast \rho_{\vartheta}}\leq \mathbf{\Theta}_0<\infty,
			\end{align}
			there exists $T>0$ such that \eqref{peskin} admits a unique solution  $X$ in $[0,T]$ satisfying
			\begin{align*}
				&\sup_{t\in[0,T]}	(\|\partial_x X(t)\|_{L^\infty}+t^{m+\kappa}\|\partial_x X(t)\|_{\dot C^{m+\kappa}})\leq C \|\partial_x X_0\|_{L^\infty}.\\
				&\mathbf{\Theta}(T)\leq 2\mathbf{\Theta}_0.
			\end{align*}
		\end{theorem}
		Note that $(C^2)^{\dot B^1_{\infty,\infty}}\backslash \dot W^{1,\infty}\neq\emptyset$, see \cite[V.4.3.1]{Steinbook} and \cite{DLN14}. It is well-known that  $C^1\subsetneqq VMO^1\subsetneqq \dot B^1_{\infty,\infty}\cap (C^2)^{\dot B^1_{\infty,\infty}}$.
		As we remarked before, here we need $X_0\in \dot W^{1,\infty}$ to control the variable coefficients. Moreover, the condition $X_0\in (C^2)^{\dot B^1_{\infty,\infty}}$ helps to write the main part of the system as a semi-linear system. 
		\vspace{0.3cm}\\
		{\bf The 2D Muskat equation with surface tension}\\ Consider the 2D Muskat problem in the setting where the two fluids
		have equal viscosities and different densities. Moreover, we take the effect of surface tension into consideration. The free boundary can be described by the evolution equation (see section \ref{intomus1} for details)
		\begin{align}\label{eqmst0}
			\partial_t f(x)=\frac{\mathbf{k}}{2\pi\mu}\mathrm{P.V.} \int_{\mathbb{R}}\frac{1+\partial_x f(x)\Delta_\alpha f(x)}{\langle\Delta_\alpha f(x)\rangle^2}\partial_x\left(\sigma'\kappa(f)- \varrho_0 f\right)(x-\alpha)\frac{d\alpha}{\alpha},
		\end{align}
		for $(t,x)\in (0,+\infty)\times\mathbb{R}$.
		The main result is the following.
		\begin{theorem}\label{thmmst}
			Fix $\kappa$ such that $0<3-\kappa\ll1$ and $m\in\mathbb{N}$. There exists $\varepsilon_0>0$ such that 	~\\
			i) ($ \varrho_0\geq 0$)	For any initial data  $f_0$ with $\| f_0'\|_{L^\infty}\leq \varepsilon_0$, there exists a unique global solution $f$ to \eqref{eqmst0} such that
			\begin{align*}
				\sup_{t>0}(\|f'(t)\|_{L^\infty}+t^\frac{m+\kappa}{3}\|f'(t)\|_{\dot C^{m+\kappa}})\leq C \| f_0'\|_{L^\infty}.
			\end{align*}
			ii) ($ \varrho_0\in\mathbb{R}$) For any initial data $f_0\in W^{1,\infty}$, if there exists $\eps_1>0$ such that 
			\begin{align}\label{stcon}
				\| f_0-f_0\ast \rho_{\eps_1}\|_{W^{1,\infty}}\leq \varepsilon_0,
			\end{align}
			then   there exists $T_0>0$ such that \eqref{eqmst0} admits a unique solution  $f$ satisfying
			\begin{align*}
				\sup_{t\in[0,T_0]}(\|f'(t)\|_{L^\infty}+t^\frac{m+\kappa}{3}\|f'(t)\|_{\dot C^{m+\kappa}})\leq C \| f_0'\|_{L^\infty}.
			\end{align*}
		\end{theorem}
		\begin{remark}
			When $\varrho_0<0$ we can't get the global result, since in low frequency $\Lambda f$ cannot be controlled by $\Lambda^3f$.
		\end{remark}
		
		We remark that our method is not applicable to the Muskat equation without surface tension (when $\sigma'=0$ and $\varrho_0>0$). The main difficulty  comes from the fact that we cannot move the derivative from the term $\partial_x f(x)-\Delta_\alpha f(x)$ to other lower order terms, hence we can not close the a priori estimate with Lipschitz data. This can be  overcame using a different  method in \cite{KeC1}. For the same reason, our method cannot be applied to the Muskat equation with different viscosity. However, we can do for $C^{1,\log^{\varkappa}}$ data with $\varkappa>1.$ \vspace{0.3cm}\\
		{\bf The Muskat equation in general setting}\\
		Consider the Muskat equation in graph case, which describes the dynamics of two immiscible fluids in a porous medium with different densities $\rho^{\pm}$ and different
		viscosities $\mu_{\pm}$. 
		\begin{figure}[h]\label{figure1}
			\centering
			\includegraphics[width=6cm,height=6cm]{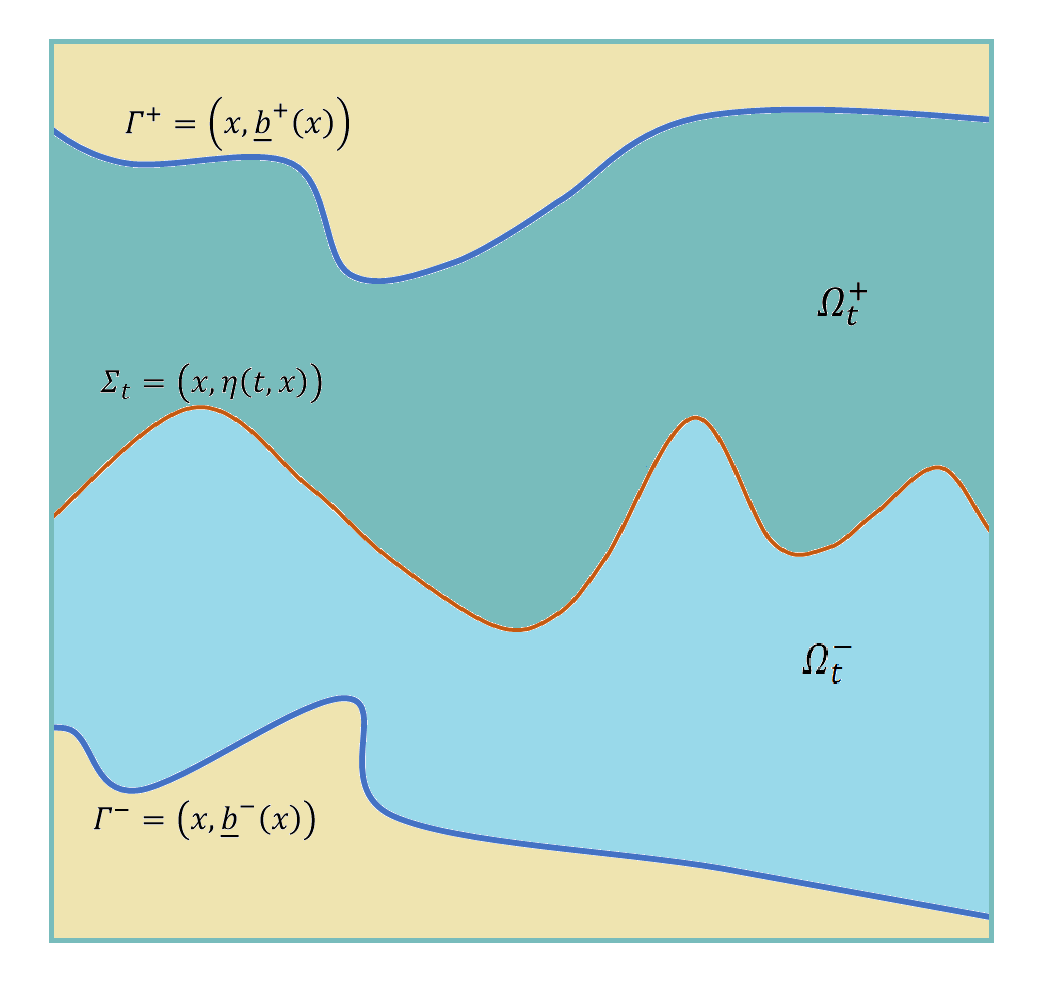}
			\caption{The Muskat problem}
		\end{figure}Assume the interface between the two fluids is a graph of a time-dependent function $\eta(t,x)$, which can be described by 
		\begin{align*}
			\Sigma_t=\{(x,\eta(t,x)):x\in\mathbb{R}^{d}\}.
		\end{align*}
		The domains are given by
		\begin{align*}
			&\Omega_{t}^{+}=\left\{(x, z) \in \mathbb{R}^{d} \times \mathbb{R}: \eta(t, x)<z<\underline{b}^{+}(x)\right\}, \\
			&\Omega_{t}^{-}=\left\{(x, z) \in \mathbb{R}^{d} \times \mathbb{R}: \underline{b}^{-}(x)<z<\eta(t, x)\right\},
		\end{align*}
		where $\underline{b}^\pm$ are the parametrizations of the rigid boundaries
		$$
		\Gamma^{\pm}=\left\{\left(x, \underline{b}^{\pm}(x)\right): x \in \mathbb{R}^{d}\right\}.
		$$
		The incompressible fluid velocity $u^\pm$ in each region is governed by Darcy’s law
		\begin{align*}
			\mu_{\pm} u^{\pm}+\nabla_{x, z} p^{\pm}=-\rho^{\pm} \mathfrak{g} \vec{e}_{d+1}, \quad \operatorname{div}_{x, z} u^{\pm}=0 \quad \text { in } \Omega^\pm(t).
		\end{align*}
		Here $p^\pm$ is pressure, $\mathfrak{g}$ is the Earth’s gravity. In this section, we simply denote $\nabla=\nabla_x$, and $\nabla_{x,z}=(\nabla,\partial_z)$.
		The normal velocity is continuous at the interface 
		\begin{align*}
			u^+\cdot n=	u^-\cdot n.
		\end{align*}
		where $n=\frac{1}{\sqrt{1+|\nabla \eta|^2}}(-\nabla \eta, 1)$. 
		The pressure jump at the interface is proportional to the mean curvature:
		\begin{align*}
			p^--p^+=-\operatorname{div}\left(\frac{\nabla\eta}{\langle\nabla \eta\rangle}\right),\ \ \ \text{on}\ \Sigma_t.
		\end{align*}
		The interface moves with the fluid
		\begin{equation*}
			\partial_{t} \eta=\left.\sqrt{1+|\nabla \eta|^{2}} u^{-} \cdot n\right|_{\Sigma_{t}}= \left. \mathbf{n}\cdot u^{-}\right|_{\Sigma_{t}},
		\end{equation*}
		where we denote $\mathbf{n}=\sqrt{1+|\nabla \eta|^{2}} n$ for simplicity.\vspace{0.1cm}\\
		At the two rigid boundaries, the no-penetration boundary conditions are imposed
		\begin{align*}
			u^\pm\cdot \nu^\pm=0,\ \ \text{on}\ \Gamma^\pm,
		\end{align*}
		where $\nu^\pm=\frac{1}{\langle\nabla \underline{b}^\pm\rangle}(-\nabla\underline{b}^\pm,1)$. \\
		Let $ q^\pm(x,z)=p(x,z)+\rho^\pm\mathfrak{g}z$ and denote  $\varrho_0=(\rho^--\rho^+)\mathfrak{g}$, We rewrite the system as 
		\begin{equation}\label{sysq}
			\begin{aligned}
				&\Delta_{x, z} q^{\pm}=0, \ \ \text{in}\ \Omega^\pm(t),\\
				&q^--q^+=-\operatorname{div}\left(\frac{\nabla \eta}{\langle\nabla \eta\rangle}\right)+\varrho_0\eta,\ \ \text{on}\ \Sigma_t,\\
				&\mathbf{n}\cdot\nabla_{x,z}\left(\frac{1}{\mu_+}q^+-\frac{1}{\mu_-}q^-\right) =0,\ \ \text{on}\ \Sigma_t,\\
				& \nu^\pm\cdot\nabla_{x,z}q^\pm=0,\ \ \text{on}\ \Gamma^\pm.
			\end{aligned} 
		\end{equation}
		The evolution equation of the interface reads 
		\begin{align}\label{evo}
			\partial_t\eta=-\frac{1 }{\mu_-}\left.\mathbf{n}\cdot\nabla_{x,z}q^-\right|_{\Sigma_t}.
		\end{align}
		\begin{theorem}\label{thmgm} Let $\kappa$ be such that $0<3-\kappa\ll 1$, and $m\in\mathbb{N}^+$.
			Take $\eta_0\in L^2(\mathbb{R}^d)\cap\dot C^{1,\log^\varkappa  }(\mathbb{R}^d)$ and $\varkappa>1$ such that \begin{align*}
				\operatorname{dist}(\eta_0,\Gamma^\pm)>2\mathbf{r}>0.
			\end{align*}
			There exist $0<\varepsilon_0,\varepsilon_1\ll 1$ such that if $\eta_0$ further satisfies 
			\begin{align}
				\label{gmcon}
				\|\eta_0-\eta_0\ast \rho_{\eps_1}\|_{L^2\cap \dot C^{1,\log^\varkappa  }}\leq \varepsilon_0,
			\end{align}
			for some $\eps_1>0$,
			where the definition of $\|\cdot\|_{\dot C^{1,\log ^\varkappa}}$ can be found in \eqref{ctlog}, then there exist $T>0$, depending on $\|\eta_0\|_{L^2\cap \dot C^{1,\log^\varkappa  }}, \mathbf{r}, \varepsilon_0,\eps_1$, $m$,  and a unique solution $\eta$ in $[0,T]$ satisfying 
			\begin{align}
				&	\sup_{t\in[0,T]}(\|\eta(t)\|_{L^2}+\|\nabla \eta(t)\|_{\dot C^{\log^\varkappa  }}+t^\frac{m+\kappa}{3}|\log t|^\varkappa   \|\nabla \eta(t)\|_{\dot C^{m+\kappa}})\leq C(\|\eta_0\|_{L^2}+\|\nabla\eta_0\|_{\dot C^{\log^\varkappa  }}),\label{fes1}\\
				&\inf_{t\in[0,T]}\operatorname{dist}(\eta(t),\Gamma^\pm)>\mathbf{r}.\label{dist}
			\end{align}
		\end{theorem}
		We remark that the result \eqref{dist} shows that the interface $\Sigma_t$ does not touch the rigid boundaries $\Gamma^\pm$ (see Figure \ref{figure1}). Compared with the critical Lipschitz space, we need extra $\log^\varkappa  $ regularity for initial data. 
		We remark that the condition $\varkappa>1$ is necessary and we can not obtain well-posedness results in critical Lipschitz space by our method. To see this, we can have a look at the equation in 2D without rigid boundary, which can be written as (see \cite{Matioctams}): 
		\begin{equation}\label{mus2d}
			\begin{aligned}
				& \partial_t f(x)=\int_{\mathbb{R}}  \frac{1+\partial_x f(x) \Delta_\alpha f(x)}{\left\langle\Delta_\alpha f(x)\right\rangle^2} \omega(x-\alpha) \frac{d \alpha}{\alpha},
			\end{aligned}
		\end{equation}
		where $f$ and $\omega$ are further coupled by the relation
		\begin{align}
			& \tilde \sigma \partial_x\left(\frac{\partial_x^2f}{\langle\partial_x f\rangle^3}\right)+\varrho \partial_x f(x)=\omega(x)+(\mu_--\mu_+) K(f)(\omega)(x),\label{ome}\\
			&
			K(f)(\omega)(x)=\frac{1}{\pi} \int_{\mathbb{R}} \frac{\partial_x f(x)-\Delta_\alpha f(x)}{\left\langle\Delta_\alpha f(x)\right\rangle^2} \omega(x-\alpha) \frac{d \alpha}{\alpha},\nonumber
		\end{align}
		where $\tilde \sigma>0$ is the surface tension coefficient at the free boundary, $\varrho\in\mathbb{R}$ denotes the difference between the upper and lower fluid densities.  
		When considering two fluids with the same viscosity, i.e. $\mu_--\mu_+=0$, we don't have the term $K(f)(w)$ in \eqref{ome} and the equation goes back to the one in section \ref{intomus1}. In the presence of different viscosities, we substitute \eqref{ome} into \eqref{mus2d}. 
		The main dissipative term in the evolution equation still comes from the effect of surface tension, i.e. the term associate  to $\partial_x\left(\frac{\partial_x^2f}{\langle\partial_x f\rangle^3}\right)$. The term associate to $K(f)(\omega)$ is nonlinear error term. By Theorem \ref{lemmain}, to obtain a priori estimate, one needs to control the $C^a$ norm of $\int_{\mathbb{R}} \frac{1+\partial_x f(x) \Delta_\alpha f(x)}{\left\langle\Delta_\alpha f(x)\right\rangle^2} K(f)(\omega)(x-\alpha) \frac{d \alpha}{\alpha}$. This in turn requires the boundedness of $\int_{\mathbb{R}}\frac{\partial_xf-\Delta_\alpha f}{\langle\Delta_\alpha f\rangle^2}\frac{d\alpha}{\alpha}$. It is easy to check that the $\dot C^{1,\log^\varkappa  }$ regularity with $\varkappa>1$ is the least to guarantee the desired  boundedness.\vspace{0.3cm}\\

		\subsection{Quasi-linear parabolic system with boundary condition}
		It is more subtle to study quasi-linear parabolic systems in a bounded domain because the solution may behave badly near the boundary. For example,  
		it is well known that the 3d Navier–Stokes equations with Dirichlet boundary condition has non $L^\infty_tC^\alpha$ solution for some $\alpha>0$ since the Leray operator $\mathbb{P}=\mathrm{Id}-\nabla\Delta_{D}^{-1}\operatorname{div}$, where $\Delta_{D}$ is the Neumann Laplacian, is irregular at the boundary, see \cite{ChangKang,sesv}. Moreover, the regularity of SQG equation in the whole space is well-studied in \cite{CWsqg}, but the problem in bounded domain was open until the recent result 
		by  Constantin, Ignatova, and the third author \cite{Constantin2023}, which proved that the solution of the  critical dissipative SQG equation on bounded domain is smooth using a very special structure of nonlinear terms.\vspace{0.1cm}\\
		As the last part of this paper, we study the well-posedness of the mean curvature flow system in domain with homogeneous Dirichlet boundary condition, which can be written as 
		\begin{equation}\label{meancur1}
			\begin{aligned}
				&\partial_t f=\A[\nabla f]:\nabla^2 f,\ \ \ \text{in}\ [0,T]\times \Omega,\\
				&f|_{t=0}=f_0,\ \ \ \text{in}\ \Omega,\\
				&f=0,\ \ \text{on}\ [0,T]\times\partial \Omega.
			\end{aligned}
		\end{equation}
		Here $\Omega\subset \mathbb{R}^d$ is $\mathbb{R}_+^d$ or a bounded domain with $\partial\Omega\in C^{2m+3}$, the coefficient $\A$ is defined in \eqref{defamcf}. 
		The main results for $\Omega=\mathbb{R}_+^d$ is the following.
		\begin{theorem}\label{eqmcls} Let $m\in \mathbb{N}$ with $m\geq 10d$, $\kappa>0$ with $1-\kappa \ll 1$.
			Assume $\Omega=\mathbb{R}_+^d$. There exists $\eps_0>0$ such that \\
			i) For any initial data $f_0\in \dot{W}^{1,\infty}$ with $f_0|_{\partial\mathbb{R}_+^d}=0$, if there exists a smooth function $\phi\in C^\infty(\bar{\mathbb{R}}_+^d)$ with $\phi|_{\partial\mathbb{R}_+^d}=0$, and 
			\begin{equation*}
				\|f_0-\phi\|_{\dot{W}^{1,\infty}}\leq\eps_0,
			\end{equation*}
			then there exists $T_0>0$ such that \eqref{meancur1} admits a unique solution $f$ in $[0,T_0]$ with initial data $f_0$ such that
			\begin{equation*}
				\sup_{t\in[0,T]}\left(\|\nabla f(t)\|_{L^\infty}+t^{m+\frac{1+\kappa}{2}}\|\nabla^{2m}(\nabla^2f,\partial_t f)(t)\|_{\dot C^{\kappa}}\right)\leq C\|f_0\|_{\dot{W}^{1,\infty}}.
			\end{equation*}
			and 
			\begin{equation}
				\Delta^{k}f(t)\vert_{\partial\mathbb{R}_+^d}=0,~~k=1,...,m,~~\forall ~t\in (0,T],
			\end{equation}
			for some constant $C>0.$ \vspace{0.3cm}\\
			ii) For any initial data $f_0\in\dot {W}^{1,\infty}(\mathbb{R}_+^d) $ with $f_0|_{\partial\mathbb{R}_+^d}=0$, if 
			\begin{equation*}
				\|f_0\|_{\dot{W}^{1,\infty}}\leq \eps_0,
			\end{equation*}
			then \eqref{meancur1} admits a unique solution f in $[0,\infty)$ with initial data $f_0$, such that
			\begin{equation*}
				\sup_{t>0}\left(\|\nabla f(t)\|_{L^\infty}+t^{m+\frac{1+\kappa}{2}}\|\nabla^{2m}(\nabla^2f,\partial_t f)(t)\|_{\dot C^{\kappa}}\right)\leq C\|f_0\|_{W^{1,\infty}},
			\end{equation*}
			and 
			\begin{equation*}
				\Delta^{k}f(t)\vert_{\partial\mathbb{R}_+^d}=0,~~k=1,...,m,~~\forall ~t>0,
			\end{equation*}
			for some constant $C>0$.
		\end{theorem}
		Moreover, for the case where $\Omega$ is a bounded domain with $\partial\Omega\in C^{2m+3}$, we have the following theorem.
		\begin{theorem}\label{mcbdglo}
			Let $\Omega$ be a bounded domain with $C^{2m+3}$ boundary. There exists $\eps_0>0$ such that,\\
			i) For any initial data $f_0\in W^{1,\infty}(\Omega)$ with $f_0|_{\partial\Omega}=0$, if there exists $\eps_1$ such that 
			\begin{equation*}
				\|f_0-e^{\eps_1\Delta}f_0\|_{W^{1,\infty}}\leq\eps_0,
			\end{equation*}
			where $e^{t\Delta}$ is the heat semigroup for $\Omega$ with Dirichlet boundary condition, then there exists $T_0>0$ such that \eqref{meancur1} admits a unique solution f in $[0,T_0]$ with initial data $f_0$ such that \begin{equation}\label{mcbdhd}
				\sup_{t\in[0,T]}\|\nabla f\|_{L^\infty}+\sup_{t\in[0,T]}t^{m+\frac{\kappa}{2}}\|\nabla^{2m+1}f(t)\|_{C^{\kappa}}\leq C\|f\|_{W^{1,\infty}},
			\end{equation}
			for some constant $C>0.$\vspace{0.3cm}\\
			ii) For any initial data $f_0\in W^{1,\infty}(\Omega)$ with $f_0|_{\partial\Omega}=0$, if 
			\begin{equation*}
				\|f_0\|_{W^{1,\infty}}\leq\eps_0,
			\end{equation*}
			then \eqref{meancur1} admits a unique solution f in $[0,\infty)$ such that 
			\begin{equation}\label{bddmcgloe}
				\sup_{t>0}\sum_{k=1}^m e^{ct}t^{-\frac{k-1}{2}}\|f(t)\|_{C^k}\leq C\|f_0\|_{W^{1,\infty}},
			\end{equation}
			for some constants $C>0$ and $c>0$.\\
			
		\end{theorem}
		Note that the property  $  \Delta^{k}f(t)\vert_{\partial\Omega}=0$ for $k\geq 1$ is not available when $\Omega$ is a bounded domain,  see Remark \ref{bdydiff}.\vspace{0.1cm}\\
		It is very interesting to extend this result to the Willmore flow \eqref{eqwma}, the surface diffusion flow \eqref{eqsurfdiffa} and other quasi-linear parabolic systems. We remark that this is not trivial, as both our proof of mean curvature flow in Section \ref{secbdd} and the proof of SQG equation in \cite{Constantin2023} needs to take advantage of the special structure of nonlinear terms. We leave further applications of our method to other equations for future research.
		\vspace{0.2cm}\\
		{\bf{Data availability statement.}} This research does not have any associated data.\\
		{\bf {Conflict of Interest.}} The authors declare that there are no conflicts of interest.\\
		{\bf{Acknowledgments.}}  Quoc-Hung Nguyen  is supported by the Academy of Mathematics and Systems Science, Chinese Academy of Sciences startup fund; CAS Project for Young Scientists in Basic Research, Grant No. YSBR-031;  and the National Natural Science Foundation of China (No. 12288201);  and the National Key R$\&$D Program of China under grant 2021YFA1000800. Ke Chen is supported by the Research Centre for Nonlinear Analysis at The Hong Kong Polytechnic University.
		\section{Proof of the main theorems}\label{secproof}
		When dealing with equations with constant coefficients, it is easy to write the equation in terms of Fourier transform and find the fundamental solution. We split the proof of Theorem \ref{thmcons} and Theorem \ref{lemmain} into four steps. First, we fix the coefficient, and transform the differential equation into an integral equation by Fourier method. Then we give some elementary estimates of the fundamental solution. Subsequently, we estimate the solution and show that the effect of smooth variable coefficients are lower-order error terms that can be ignored in short time estimate. Finally, we obtain a global-in-time estimate by a simple iterative procedure and complete the proof. \vspace{0.3cm}\\
		Now we start the proof of Theorem \ref{thmcons} and \ref{lemmain}.\\
		\textbf{Step 1. Transform the differential equation to an integral equation.}\\
		Fix $x_0\in\mathbb{R}^d$, we rewrite \eqref{eqpara} as 
		\begin{equation}\label{eqpafi}
			\begin{aligned}
				&\partial_tu(t,x)+(2\pi)^{-\frac{d}{2}}\int _{\mathbb{R}^d}\A(t,x_0,\xi)\hat{u}(t,\xi)e^{ix\cdot\xi}d\xi=\mathcal{P}_\gamma f(t,x)+g(t,x)+R_{x_0}[u](t,x),
			\end{aligned}
		\end{equation}
		with 
		\begin{align*}
			R_{x_0}[u](t,x)=(2\pi)^{-\frac{d}{2}}\int _{\mathbb{R}^d}(\A(t,x_0,\xi)-\A(t,x,\xi))\hat{u}(t,\xi)e^{ix\cdot\xi} d\xi.
		\end{align*}
		We remark that if the coefficient matrix $\A(t,x,\xi)$ is independent in $x$, then $R_{x_0}[u](t,x)\equiv 0$.\\
		Let $\K_{
			x_0}(t,\tau,x)$ be the fundamental solution of the adjoint system, satisfying 
		\begin{equation}\label{defk0}
			\begin{aligned}
				&-\partial_\tau \K_{x_0}(t,\tau,x)+(2\pi)^{-\frac{d}{2}}\int _{\mathbb{R}^d}\hat \K_{x_0}(t,\tau,\xi)\A(\tau,x_0,\xi)e^{ix\cdot \xi}d\xi=0, \ \ (\tau,x) \in (0,t)\times \mathbb{R}^d,\\
				&\lim_{\tau\to t^-}\K_{x_0}(t,\tau,x)=\delta(x)\mathrm{Id},\ \ \ \ \ x\in\mathbb{R}^d.
			\end{aligned}
		\end{equation}
		Here $\mathrm{Id}$ is the $N\times N$ identity matrix,  $\hat \K_{x_0}(t,\tau,\xi)$ is the Fourier transform of  $\K_{
			x_0}(t,\tau,x)$ in $x$-variable, and $\delta(x)$ is the Dirac delta function.
		Then we transform the differential equation \eqref{eqpafi} to the following form,
		\begin{align}
			u(t,x)=&\int_{\mathbb{R}^d} \K_{x_0}(t,0,x-z)u_0(z)dz+\int_0^t \int _{\mathbb{R}^d}\K_{x_0}(t,\tau,x-z)(\mathcal{P}_\gamma f+g)(\tau,z)dzd\tau\nonumber \\
			&+\int_0^t\int_{\mathbb{R}^d} \K_{x_0}(t,\tau,x-z)  R_{x_0}[u](\tau,z)dzd\tau\nonumber\\
			:=	&u_{L,x_0}(t,x)+u_{N,x_0}(t,x)+u_{R,{x_0}}(t,x),\label{uform}
		\end{align}
		which holds for any $x_0\in\mathbb{R}^d$. Note that when estimate the H\"{o}lder norm of $u$, we will first take derivatives in the formula \eqref{uform}, and then fix $x_0=x$ such that the remainder term $u_{R,x_0}$ can be absorbed by the main terms, see \eqref{1.1hol}. Hence our method do not need the regularity of $\K_{x_0}$ in $x_0$. Remark that 
		\begin{equation*}
			\big(\nabla^n_x(u_{L,x_0}(t,x)+u_{N,x_0}(t,x)+u_{R,{x_0}}(t,x))\big)|_{x_0=x}\neq \nabla^n_x\big((u_{L,x_0}(t,x)+u_{N,x_0}(t,x)+u_{R,{x_0}}(t,x))|_{x_0=x} \big).
		\end{equation*}
		\vspace{0.1cm}\\
		{\bf Step 2. Properties of the fundamental solution}\\
		To estimate the kernel $\K_{x_0}$ and its derivatives, we first prove the following lemma.
		\begin{lemma}\label{esthtker}
			Let $\lambda>0.$	Assume  $\F(\xi)$ satisfies
			\begin{equation}\label{condb}
				|\nabla^n_\xi\F(\xi)|\leq |\xi|^{\sigma-n}\min\{1,\lambda^{-1}|\xi|^{-1}\}^{[d+\sigma]+1},\quad n\leq [d+\sigma]+1,
			\end{equation}
			and define
			\begin{equation*}
				\G(x)=\int_{\mathbb{R}^d}\F(\xi)e^{ix\cdot\xi}d\xi.
			\end{equation*}
			If $\sigma>-d$, then there holds
			\begin{align}\label{leme}
				\left|\G(x)\right|\lesssim \frac{1}{(|x|+\lambda)^{d+\sigma}}.
			\end{align}
			Moreover, if $\sigma>0$, then for $n\in\mathbb{N}$,  $0\leq\sigma'<\sigma$, and for any $\alpha\in\mathbb{R}^d$,
			\begin{align}\label{leme1}
				\int_{\mathbb{R}^d}|\delta_\alpha\nabla^n \G(x)||x|^{\sigma'}dx\lesssim \frac{|\alpha|^{\sigma'}}{\lambda^{n+\sigma}}\min\left\{1,\frac{|\alpha|}{\lambda}\right\}^{1-\sigma'}.
			\end{align}
		\end{lemma}
		\begin{proof}
			We split 
			\begin{align*}
				\G(x)	&=\int_{\mathbb{R}^d} \F(\xi)\chi_{\delta}(\xi)e^{ix\cdot \xi}d\xi+\int_{\mathbb{R}^d} \F(\xi)(1-\chi_{\delta}(\xi))e^{ix\cdot \xi}d\xi\\
				&\quad\quad:=\text{I}+\text{II},
			\end{align*}
			where $\chi_\delta$ is a smooth cut-off function satisfying $\mathbf{1}_{|\xi|\leq \delta}\leq \chi_\delta(\xi)\leq \mathbf{1}_{|\xi|\leq 2\delta}$ and $\|\nabla^k\chi_\delta\|_{L^\infty}\lesssim_k\delta^{-k}$ for any $k\in \mathbb{N}.$
			Then 
			\begin{align*}
				|	\text{I}|\lesssim \int_{|\xi|\leq 2\delta} |\xi|^\sigma d\xi\lesssim \delta ^{d+\sigma}. 
			\end{align*}
			Then we deal with $\text{II}$. 	Define the operator $\hat{\nabla}=-\frac{ix\cdot\nabla_\xi}{|x|^2}$, then one has $\hat{\nabla}(e^{ix\cdot \xi})=e^{ix\cdot\xi}$. 
			Hence 
			\begin{align*}
				\text{II}=\int_{\mathbb{R}^d} \F(\xi)(1-\chi_{\delta}(\xi))\hat{\nabla}^l\left(e^{ix\cdot \xi}\right)d\xi=(-1)^l\int_{\mathbb{R}^d} \hat{\nabla}^l\left(\F(\xi)(1-\chi_{\delta}(\xi))\right) e^{ix\cdot \xi} d\xi.
			\end{align*}
			By \eqref{condb}, there holds for $l=[d+\sigma]+1$, 
			\begin{align*}
				|\text{II}|
				&\lesssim\sum_{k=0}^l\int_{\mathbb{R}^d}\left|\hat{\nabla}^k\F(\xi)\hat{\nabla}^{l-k}(1-\chi_{ \delta}(\xi)) \right|d\xi\\
				&\lesssim|x|^{-l}\left( \sum_{k=0}^{l-1}\delta^{k-l}\int_{\delta\leq|\xi|\leq 2\delta}|\xi|^{\sigma-k}+\int_{|\xi|\geq \delta}|\xi|^{\sigma-l}d\xi \right)\\
				&\lesssim \frac{ \delta^{d+\sigma-l}}{|x|^{l}}.
			\end{align*}
			Taking $\delta=|x|^{-1}$ to obtain that 
			\begin{align}\label{11}
				\left|	 \G(x)\right|\lesssim 		|	\text{I}|+	|\text{II}|\lesssim |x|^{-(d+\sigma)}.
			\end{align}
			On the other hand, we have 
			\begin{align*}
				\left|	 \G(x)\right|\lesssim \int_{\mathbb{R}^d} |\xi|^{\sigma}\min\{1,\lambda^{-1}|\xi|^{-1}\}^{[d+\sigma]+1}d\xi\lesssim \lambda^{-(d+\sigma)}.
			\end{align*}
			Combining this with \eqref{11}, we obtain \eqref{leme}. 
			A direct result of \eqref{leme} is 
			\begin{align*}
				|\nabla^n\G(x)|\lesssim \frac{1}{(|x|+\lambda)^{d+n+\sigma}},\ \ \ \forall\ n\in\mathbb{N}.
			\end{align*}
			Hence, for any $\alpha\in\mathbb{R}^d$, if $|\alpha|\geq \lambda$, it is easy to obtain that 
			\begin{equation}\label{d111}
				\begin{aligned}
					& \int_{\mathbb{R}^d}  |\delta_\alpha \nabla^n\G(x)||x|^{\sigma'}dx\\&\ \ \ \ \lesssim	  \int_{\mathbb{R}^d}  (| \nabla^n\G(x)|+| \nabla^n\G(x-\alpha)|)|x|^{\sigma'}dx\\
					&\ \ \ \	  \lesssim \int_{\mathbb{R}^d} \left(\frac{1}{(|x|+\lambda)^{d+n+\sigma}}+\frac{1}{(|x-\alpha|+\lambda)^{d+n+\sigma}}\right)|x|^{\sigma'}dx\\
					&\ \ \ \ \lesssim \int_{\mathbb{R}^d} \frac{|x|^{\sigma'}+|\alpha|^{\sigma'}}{(|x|+\lambda)^{d+n+\sigma}}dx\\
					&\ \ \ \ \lesssim\frac{|\alpha|^{\sigma'}}{\lambda^{n+\sigma}},
				\end{aligned}
			\end{equation}
			provided $0\leq \sigma'<\sigma$. Moreover, if $|\alpha|<\lambda$,	we can write
			\begin{align*}
				|\delta_\alpha \nabla^n\G(x)|\lesssim |\alpha|\int_0^1|\nabla^{n+1}\G(x-\tau\alpha)|d\tau.
			\end{align*}
			Then we obtain 
			\begin{align*}
				\int_{\mathbb{R}^d}  |\delta_\alpha \nabla^n\G(x)||x|^{\sigma'}dx\lesssim& |\alpha|\int_0^1  \int_{\mathbb{R}^d} |\nabla^{n+1}\G(x-\tau\alpha)||x|^{\sigma'}dxd\tau\\
				\lesssim&|\alpha|\int_0^1  \int_{\mathbb{R}^d} \frac{|x|^{\sigma'}}{(|x-\tau\alpha|+\lambda)^{d+n+1+\sigma}}dxd\tau\\
				\lesssim&  |\alpha|\int_0^1  \int_{\mathbb{R}^d} \frac{|x|^{\sigma'}+|\alpha|^{\sigma'}}{(|x|+\lambda)^{d+n+1+\sigma}}dxd\tau\\
				\lesssim &\frac{|\alpha|}{\lambda^{n+1+\sigma-\sigma'}}.
			\end{align*}
			Combining this with \eqref{d111}, we  obtain \eqref{leme1}.
			This completes the proof of Lemma \ref{esthtker}.
		\end{proof}\\
		\begin{remark} It is well-known that \begin{align}
				|x|^{d+\sigma}\left|\int_{\mathbb{R}^d}|\xi|^\sigma \chi(\xi) e^{ix\cdot\xi}d\xi\right|\sim 1
			\end{align}
			for any $|x|\geq L \gg1$ when $\sigma\not \in \mathbb{Z}$. Thus,  the estimate \eqref{leme} is optimal  for the decay in $x$.
		\end{remark}
		\begin{remark}
			Lemma \ref{esthtker} is available for a family of functions $\{F_{\lambda}(\xi)\}_{\lambda>0}$, where in some cases, the functions $F_{\lambda}(\xi)$ cannot be presented as $F(\lambda,\xi)$. For example, in Remark \ref{lemfourierK}, we apply  Lemma \ref{esthtker} to a family of functions  $F_\lambda(\xi)=\hat K_{x_0}(t_1,t_2,\xi)$ with $\lambda=(t_1-t_2)^\frac{1}{s}$.
		\end{remark}
		\begin{lemma}\label{lemadk}
			Let $\A(t,x,\xi)\in\mathbb{R}^{N\times N}$ satisfy \eqref{condop}, consider $\K_{x_0}(t,\tau,x)$ satisfying \eqref{defk0}.
			Then for the Fourier transform $\hat{\K}_{x_0}(t,\tau,\xi)$, we have the following estimates
			\begin{align}
				&\left|\hat{\K}_{x_0}(t,\tau,\xi) \right|\lesssim e^{-c_0(t-\tau)|\xi|^s},	 \label{expA}\\
				&\left|\nabla_{\xi}^l\hat{\K}_{x_0}(t,\tau,\xi) \right|\lesssim  \tilde\M^{l}(t-\tau)|\xi|^{{s}-l}e^{-\frac{c_0}{2}(t-\tau)|\xi|^s},\quad\forall l\in\mathbb{N}_+,  l\leq d+m+s+2,\label{derA}
			\end{align}
			where $\tilde \M$ is defined in \eqref{deftilM}.
		\end{lemma}
		\begin{proof}
			We take Fourier transform of \eqref{defk0} to get
			\begin{equation*}
				\begin{aligned}
					&-\partial_\tau\hat{\K}_{x_0}(t,\tau,\xi)+\hat{\K}_{x_0}(t,\tau,\xi)\A(\tau,x_0,\xi)=0,\ \ \tau\in(0,t),\\
					&\lim_{\tau\rightarrow t^-}\hat{\K}_{x_0}(t,\tau,\xi)=\mathrm{Id}.
				\end{aligned}
			\end{equation*}
			Remark that $\hat{\K}_{x_0}(t,\tau,\xi)$ is real-valued because $\A(t,x,\xi)\in\mathbb{R}^{N\times N}$.
			By a change of variable $\hat {\mathcal{K}}_{x_0}(t,\tau,\xi)=\hat{\K}_{x_0}(t,t-\tau,\xi)$, we transform the backward evolution system to a forward one:
			\begin{equation}\label{forwK}
				\begin{aligned}
					&\partial_\tau \hat {\mathcal{K}}_{x_0}(t,\tau,\xi)+\hat {\mathcal{K}}_{x_0}(t,\tau,\xi)\A(\tau,x_0,\xi)=0,\ \ \tau\in(0,t),\\
					&\lim_{\tau\to 0^+}\hat {\mathcal{K}}_{x_0}(t,\tau,\xi)=\mathrm{Id}.
				\end{aligned}
			\end{equation}
			Recalling the definition of Frobenius inner product and matrix norm in \eqref{matrxnorm}, we have 
			\begin{align*}
				\frac{1}{2}\partial_\tau |\hat {\mathcal{K}}_{x_0}(t,\tau,\xi)|^2+\left(\hat {\mathcal{K}}_{x_0}(t,\tau,\xi)\A(\tau,x_0,\xi)\right):\hat {\mathcal{K}}_{x_0}(t,\tau,\xi)=0.
			\end{align*}
			By the uniform elliptic condition \eqref{condop}, one has 
			\begin{align*}
				\left(\hat {\mathcal{K}}_{x_0}(t,\tau,\xi)\A(\tau,x_0,\xi)\right):\hat {\mathcal{K}}_{x_0}(t,\tau,\xi)\geq c_0|\xi|^s|\hat {\mathcal{K}}_{x_0}(t,\tau,\xi)|^2.
			\end{align*}
			Hence, we get 
			\begin{align*}
				&\partial_\tau |\hat {\mathcal{K}}_{x_0}(t,\tau,\xi)|^2+2 c_0|\xi|^s|\hat {\mathcal{K}}_{x_0}(t,\tau,\xi)|^2\leq 0,\\
				&\lim_{\tau\to 0^+} |\hat {\mathcal{K}}_{x_0}(t,\tau,\xi)|^2=N.
			\end{align*}
			Then it follows from Gronwall's inequality that 
			\begin{align}\label{huaK}
				|\hat {\mathcal{K}}_{x_0}(t,\tau,\xi)|^2\leq N e^{-2c_0\tau|\xi|^s},\ \ \ \forall \tau\in(0,t).
			\end{align}
			This implies that 
			\begin{align*}
				|\hat \K_{x_0}(t,\tau,\xi)|=|\hat {\mathcal{K}}_{x_0}(t,t-\tau,\xi)|\lesssim e^{-c_0(t-\tau)|\xi|^s},\ \ \ \forall \tau\in(0,t).
			\end{align*}
			This completes the proof of \eqref{expA}.
			
			For the estimates of derivatives, we prove by induction. We first prove that \eqref{derA} holds for $l=1$. Indeed, taking one derivative in \eqref{forwK} to obtain 
			\begin{align*}
				&\partial_\tau \nabla_\xi\hat {\mathcal{K}}_{x_0}(t,\tau,\xi)+\nabla_\xi\hat {\mathcal{K}}_{x_0}(t,\tau,\xi)\A(\tau,x_0,\xi)=-\hat {\mathcal{K}}_{x_0}(t,\tau,\xi)\nabla_\xi\A(\tau,x_0,\xi),\ \ \tau\in(0,t),\\
				&\lim_{\tau\to 0^+}\nabla_\xi\hat {\mathcal{K}}_{x_0}(t,\tau,\xi)=0.
			\end{align*}
			Then taking Frobenius inner product with $\nabla_\xi\hat {\mathcal{K}}_{x_0}(t,\tau,\xi)$, and using \eqref{huaK} we can deduce that 
			\begin{align*}
				&\partial_\tau |\nabla_\xi\hat {\mathcal{K}}_{x_0}(t,\tau,\xi)|+ c_0|\xi|^s|\nabla_\xi\hat {\mathcal{K}}_{x_0}(t,\tau,\xi)|\leq |\nabla_\xi\A(\tau,x_0,\xi)||\hat {\mathcal{K}}_{x_0}(t,\tau,\xi)|\leq \M |\xi|^{s-1}e^{-c_0\tau|\xi|^s}.
			\end{align*}
			This implies that 
			\begin{align*}
				|\nabla_\xi\hat {\mathcal{K}}_{x_0}(t,\tau,\xi)|\lesssim \M|\xi|^{s-1}\tau e^{-c_0\tau|\xi|^s},\ \ \ \tau\in(0,t).
			\end{align*}
			From this we can recover the estimate of $\nabla_\xi\hat \K_{x_0}$, which leads to \eqref{derA} in the case $l=1$.\\
			For $1<l\leq d+m+s+2$, we denote $\nabla_{\xi}^l\hat{\mathcal{K}}_{x_0}(t,\tau,\xi)=\hat{\mathcal{K}}_{x_0}^{l}(t,\tau,\xi)$, and we assume that the cases for $1,2,\cdots,l-1$ have been proved. We can write the equation of $\hat{\mathcal{K}}_{x_0}^{l}(t,\tau,\xi)$ as 
			\begin{equation*}
				\begin{aligned}
					&\partial_{\tau}\hat{\mathcal{K}}_{x_0}^{l}(t,\tau,\xi)+\hat{\mathcal{K}}_{x_0}^{l}(t,\tau,\xi)\A(\tau,x_0,\xi)=R^l(t,\tau,\xi),\\
					&\lim_{\tau\rightarrow 0_+}\hat{\mathcal{K}}_{x_0}^{l}(t,\tau,\xi)=0,
				\end{aligned}
			\end{equation*}
			with the remainder term satisfying
			\begin{equation*}
				\left|R^l(t,\tau,\xi)\right|\lesssim \sum_{i=1}^l\left|\nabla_\xi^i\A(\tau,\xi)\right|\left|\hat{\mathcal{K}}_{x_0}^{l-i}(t,\tau,\xi)\right|\lesssim \M^{l}\tau|\xi|^{2s-l}e^{-\frac{c_0}{2}\tau|\xi|^s}.
			\end{equation*}
			Then taking Frobenius inner product, and applying \eqref{condop} we have
			\begin{equation*}
				\begin{aligned}
					&\partial_\tau|\hat{\mathcal{K}}_{x_0}^{l}(t,\tau,\xi)|+c_0|\xi|^s|\hat{\mathcal{K}}_{x_0}^{l}(t,\tau,\xi)|\leq \left|R^l(t,\tau,\xi) \right|.
				\end{aligned}
			\end{equation*}
			This implies that 
			\begin{align*}
				\partial_\tau\left(e^{c_0\tau|\xi|^s}|\hat{\mathcal{K}}_{x_0}^{l}(t,\tau,\xi)|\right)\lesssim e^{c_0\tau|\xi|^s}\left|R^l(t,\tau,\xi) \right|\lesssim \M^{l}\tau|\xi|^{2s-l}e^{\frac{c_0}{2}\tau|\xi|^s}.
			\end{align*}
			Integrating the inequality in $\tau$, and using the fundamental inequality 
			\begin{align*}
				\int_0^a z e^{bz}dz\lesssim \frac{ae^{ba}}{b},\ \ \ \forall a,b>0,
			\end{align*}
			we can deduce that 
			\begin{equation*}
				|\hat{\mathcal{K}}_{x_0}^{l}(t,\tau,\xi)|\lesssim \tilde\M^{l} \tau|\xi|^{s-l}e^{-\frac{c_0}{2}\tau|\xi|^s}.
			\end{equation*}
			So we have proved
			\begin{equation*}
				\left|\nabla_{\xi}^l\hat{\mathcal{K}}_{x_0}(t,\tau,\xi) \right|\lesssim \tilde\M^{l} \tau|\xi|^{{s}-l}e^{-\frac{c_0}{2}\tau|\xi|^s},\quad\forall l\in\mathbb{N}_+, l\leq d+m+s+2.
			\end{equation*}
			This completes the proof.
		\end{proof}
		
		Combining Lemma \ref{esthtker} with Lemma \ref{lemadk}, we get the following estimates of the kernel $\K_{x_0}$.
		\begin{remark}\label{lemfourierK}It is easy to check that the kernel $\hat{\K}_{x_0}(t_1,t_2,\xi)$ satisfies
			\begin{equation*}
				|\nabla_{\xi}^n\hat{\K}_{x_0}(t_1,t_2,\xi)|\lesssim \tilde\M^n (t_1-t_2)|\xi|^{s-n}\min\{1,(t_1-t_2)^{-\frac{1}{s}}|\xi|^{-1}\}^{[d+s]+1},\quad \forall n\in\mathbb{N}_+, n\leq d+m+s+2.
			\end{equation*}
			Applying Lemma \ref{esthtker} with $\lambda=(t_1-t_2)^\frac{1}{s}$, we can get the point-wise estimate for $\K(t_1,t_2,x)$ that
			\begin{align}
				\sup_{x_0\in\mathbb{R}^d}|\K_{x_0}(t_1,t_2,x)|\lesssim \frac{\tilde\M^{d+s+1}(t_1-t_2)}{((t_1-t_2)^\frac{1}{s}+|x|)^{d+s}},\quad\sup_{x_0\in\mathbb{R}^d}|\nabla_x^n\K_{x_0}(t_1,t_2,x)|\lesssim \frac{\tilde\M^{d+s+n+1}}{((t_1-t_2)^\frac{1}{s}+|x|)^{d+n}}.\label{ptKd}
			\end{align}
			As a direct result of \eqref{ptKd}, we obtain 
			\begin{align}\label{kerL1}
				\sup_{x_0\in\mathbb{R}^d}\|\nabla _x^n\K_{x_0}(t_1,t_2)|\cdot|^{\sigma}\|_{L^1_x}\lesssim\tilde\M^{d+n+s+1}(t_1-t_2)^{-\frac{n-\sigma}{s}},\ \ \forall n\in\mathbb{N},\ \sigma<n.
			\end{align} 
			We remark that by interpolation, \eqref{kerL1} implies that for any $\alpha\neq 0$,
			\begin{align}\label{delKL1}
				\sup_{x_0\in\mathbb{R}^d}	\|\nabla _x^n\delta_{\alpha}\K_{x_0}(t_1,t_2)\|_{L^1_x}\lesssim \tilde\M^{d+s+n+1}\min\left\{1,\frac{|\alpha|}{(t_1-t_2)^{\frac{1}{s}}}\right\} (t_1-t_2)^{-\frac{n}{s}},\ \ \ \ \forall n\in\mathbb{N}, n\leq d+m+s+2.
			\end{align}  
		\end{remark}
		\textbf{Step 3: Estimate the solution: short time estimate}\\
		In this step, we assume $T\in [0,1].$\\
		In the following of this section, we recall $\kappa_0=\kappa-s+\gamma$, and we define the following norm,
		\begin{align}\label{defxt}
			\|u\|_{X_T}=	\sup_{t\in[0,T]} \left(\|u(t)\|_{L^\infty}+t^\frac{m+\kappa }{s}\|u(t)\|_{\dot C^{m+\kappa }}\right)
			.\end{align}
		By \eqref{uform} we have for any $\alpha \in \mathbb{R}^d\backslash\{0\}$, and any $n\in\mathbb{N}$,
		\begin{align}\label{duform}
			\delta_\alpha \nabla^n u(t,x)=\delta_\alpha \nabla ^nu_{L,x_0}(t,x)+\delta_\alpha \nabla ^nu_{N,x_0}(t,x)+\delta_\alpha \nabla ^nu_{R,x_0}(t,x),
		\end{align}
		holds for any fixed $x_0\in\mathbb{R}^d$. We fix $x_0=x$ in \eqref{uform}, and also in \eqref{duform}. Then we get 
		\begin{align}\label{1.1infty}
			\|u(t)\|_{L^\infty}\leq \|u_{L,x_0}(t,x)|_{x_0=x}\|_{L^\infty_x}+\|u_{N,x_0}(t,x)|_{x_0=x}\|_{L^\infty_x}+\|u_{R,x_0}(t,x)|_{x_0=x}\|_{L^\infty_x},
		\end{align}
		and 
		\begin{equation}\label{1.1hol}
			\begin{aligned}
				\|\delta_\alpha \nabla^n u(t)\|_{L^\infty}&\leq \|\big(\delta_\alpha \nabla^nu_{L,x_0}(t,x)\big)|_{x_0=x}\|_{L^\infty_x}+\|\big(\delta_\alpha \nabla^nu_{N,x_0}(t,x)\big)|_{x_0=x}\|_{L^\infty_x}\\
				&\quad\quad+\|\big(\delta_\alpha \nabla^nu_{R,x_0}(t,x)\big)|_{x_0=x}\|_{L^\infty_x}.
			\end{aligned}
		\end{equation}
		We first consider the linear term. By \eqref{kerL1} and \eqref{delKL1}, we can prove by Young's inequality that
		\begin{equation}\label{uL}
			\begin{aligned}
				&|u_{L,x_0}(t,x)|\lesssim \|\K_{x_0}(t,0)\|_{L^1}\|u_0\|_{L^\infty}\lesssim \tilde\M^\beta \|u_0\|_{L^\infty},\\
				&	\sup_\alpha\frac{|\delta_\alpha \nabla^{m+[\kappa]}u_{L,x_0}(t,x)|}{|\alpha|^{\kappa-[\kappa]}}\lesssim \sup_\alpha \frac{\|\delta_\alpha \nabla^{m+[\kappa]}\K_{x_0}(t,0)\|_{L^1}}{|\alpha|^{\kappa-[\kappa]}}\|u_0\|_{L^\infty}\lesssim \tilde\M^{\beta}t^{-\frac{m+\kappa}{s}}\|u_0\|_{L^\infty},
			\end{aligned}
		\end{equation}
		where $\beta$ is defined in \eqref{defbet}.
		On the other hand, by Proposition \ref{normequ}, one has
		\begin{align*}
			\sup_{t>0}\sup_\alpha t^{\frac{m+\kappa}{s}}\frac{|\delta_\alpha \nabla^{m+[\kappa]}u_{L,x_0}(t,x)|}{|\alpha|^{\kappa-[\kappa]}}\lesssim \tilde\M^\beta\|u_0\|_{\dot B^0_{\infty,\infty}}.
		\end{align*}
		The above estimates hold for any $x_0\in \mathbb{R}^d$, in particular for $x_0=x$. 
		\\
		So we have proved that
		\begin{equation}\label{uLf}
			\sup_{t\in[0,T]} \left(\|u_{L,x_0}(t)|_{x_0=x}\|_{L^\infty}+t^\frac{m+\kappa }{s}\frac{\left\|\left(\delta_\alpha\nabla^{m+[\kappa]}u_{L,x_0}(t)\right)|_{x_0=x}\right\|_{L^\infty}}{|\alpha|^{\kappa-[\kappa]}}\right)\lesssim \tilde\M^{\beta}\|u_0\|_{L^\infty},
		\end{equation}
		and 
		\begin{align}\label{bes1}
			\sup_{t\in[0,T]} t^\frac{m+\kappa }{s}\frac{\|\left(\delta_\alpha\nabla^{m+[\kappa]}u_{L,x_0}(t)\right)|_{x_0=x}\|_{L^\infty}}{|\alpha|^{\kappa-[\kappa]}}\lesssim \tilde\M^\beta\|u_0\|_{\dot B^0_{\infty,\infty}}.
		\end{align}
		For $u_{N,x_0}$, note that 
		\begin{equation}\label{formuun}
			\begin{aligned}
				u_{N,x_0}(t,x)&=\int_0^t \int_{\mathbb{R}^d} \K_{x_0}(t,\tau,x-z) \mathcal{P}_\gamma f(\tau,z)dz d\tau+\int_0^t \int_{\mathbb{R}^d} \K_{x_0}(t,\tau,x-z) g(\tau,z)dz d\tau\\
				:&=u_{N,1,x_0}(t,x)+u_{N,2,x_0}(t,x).
			\end{aligned}
		\end{equation}
		Recall the definition of $\mathcal{P}_\gamma$ in \eqref{defPg}, we have 
		\begin{align*}
			u_{N,1,x_0}(t,x)&=	\int_0^t \int_{\mathbb{R}^d}  L_{x_0}(t,\tau,x-z) \Lambda^{\vartheta}f(\tau,z)dz d\tau.
		\end{align*}
		Here we take $\kappa_0>\vartheta>\kappa_0-\kappa+[\kappa]$, and $L_{x_0}(t,x)$ is defined by  
		\begin{equation}\label{mthmdefl}
			L_{x_0}(t,\tau,x) =\frac{1}{(2\pi)^\frac{d}{2}}\int _{\mathbb{R}^d}\frac{\B(\xi)}{|\xi|^{\vartheta}} \hat{\K}_{x_0}(t,\tau,\xi) e^{ix\cdot \xi}d\xi.
		\end{equation}
		Applying \eqref{leme} and \eqref{leme1} in Lemma \ref{esthtker} with 
		$	   \F_{\sigma}(\xi)=\frac{\B(\xi)}{|\xi|^{\vartheta}}\hat{\K}_{x_0}(t,\tau,\xi)$, $\sigma=\gamma-\vartheta$ and $\sigma'=\kappa_0-\vartheta$ to obtain 
		\begin{align}
			&\sup_{x_0\in\mathbb{R}^d}\int_{\mathbb{R}^d}|L_{x_0}(t,\tau,z )|z|^{\kappa_0-\vartheta}dz\lesssim \tilde\M^\beta (t-\tau)^{-\frac{\gamma-\kappa_0}{s}}, \label{esL}\\
			&\sup_{x_0\in\mathbb{R}^d}\int_{\mathbb{R}^d}|\delta_\alpha \nabla^n L_{x_0}(t,\tau,z)||z|^{\kappa_0-\vartheta}dz\lesssim  \tilde\M^\beta \frac{|\alpha|^{\kappa_0-\vartheta}}{(t-\tau)^\frac{n+\gamma-\vartheta}{s}}\min\left\{1,\frac{|\alpha|}{(t-\tau)^\frac{1}{s}}\right\}^{1-\kappa_0+\vartheta}, \label{esdelL}
		\end{align}
		for any $n\in\mathbb{N}$, $n\leq m+\kappa$.
		Moreover,  observe that 
		\begin{align}\label{uNvancd}
			\frac{1}{(2\pi)^\frac{d}{2}}\int_{\mathbb{R}^d}L_{x_0}(t,\tau,z)dz=\left.\frac{\B(\xi)}{|\xi|^{\vartheta}} \hat{\K}_{x_0}(t,\tau,\xi)\right|_{\xi=0} =0.
		\end{align}
		Hence,
		\begin{align}
			&u_{N,1,x_0}(t,x)=	\int_0^t \int_{\mathbb{R}^{d}} L_{x_0}(t,\tau,x-z) (\Lambda^{\vartheta}f(\tau,z)-\Lambda^{\vartheta}f(\tau,x))dzd\tau,\nonumber\\
			&\delta_\alpha\nabla^n u_{N,1,x_0}(t,x)=\int_0^t \int_{\mathbb{R}^{d}}  \delta_\alpha\nabla^nL_{x_0}(t,\tau,x-z) (\Lambda^{\vartheta}f(\tau,z)-\Lambda^{\vartheta}f(\tau,x))dz d\tau.\label{un1x0}
		\end{align}
		Then we can apply \eqref{esL}  to obtain that 
		\begin{align*}
			|u_{N,1,x_0}(t,x)|
			&\lesssim \int_0^t\int_{\mathbb{R}^d}| L_{x_0}(t,\tau,x-z)||x-z|^{\kappa_0-\vartheta}dz\|f(\tau)\|_{\dot{C}^{\kappa_0}}d\tau\\
			&\lesssim \tilde\M^\beta \int_0^t(t-\tau)^{-\frac{s-\kappa}{s}}\tau^{-\frac{\kappa}{s}}d\tau \sup _{t \in[0, T]} t^{\frac{\kappa}{s}}\|f(t)\|_{\dot{C}^{\kappa_0}}\\
			&\lesssim \tilde\M^\beta\sup _{t \in[0, T]} t^{\frac{\kappa}{s}}\|f(t)\|_{\dot{C}^{\kappa_0}}.
		\end{align*}
		Moreover, by \eqref{un1x0} we obtain 
		\begin{equation}\label{deun1}
			\begin{aligned}
				|\delta_\alpha \nabla^{m+[\kappa]}u_{N,1,x_0}(t,x)|\lesssim 
				&\int _0^\frac{t}{2} \int_{\mathbb{R}^d}|\delta_\alpha \nabla^{m+[\kappa]}L_{x_0}(t,\tau,x-z)||\Lambda^\vartheta f(\tau,z)-\Lambda^\vartheta f(\tau,x)|dzd\tau\\
				&\quad\quad  +\int _\frac{t}{2}^t\int_{\mathbb{R}^d} |\delta_\alpha \nabla^{[\kappa]}L_{x_0}(t,\tau,x-z)||\Lambda^\vartheta\nabla^mf(\tau,z)-\Lambda^\vartheta\nabla^mf(\tau,x)|dzd\tau\\
				&\lesssim \int _0^\frac{t}{2} \int_{\mathbb{R}^d}|\delta_\alpha \nabla^{m+[\kappa]}L_{x_0}(t,\tau,x-z)||x-z|^{\kappa_0-\vartheta}dz\|f(\tau)\|_{\dot C^{\kappa_0}}d\tau\\
				&\quad\quad  +\int _\frac{t}{2}^t \int_{\mathbb{R}^d}|\delta_\alpha \nabla^{[\kappa]}L_{x_0}(t,\tau,x-z)||x-z|^{\kappa_0-\vartheta}dz\|f(\tau)\|_{\dot C^{m+\kappa_0}}d\tau\\
				& \lesssim \tilde\M^\beta |\alpha|^{\kappa-[\kappa]} \sup_{t\in[0,T]}(t^{\frac{\kappa}{s}}\|f(t)\|_{\dot C^{\kappa_0}}+t^{\frac{m+\kappa }{s}}\|f(t)\|_{\dot C^{\kappa_0+m}}).
			\end{aligned}
		\end{equation}
		Hence, for any $t\in[0,T]$, we get
		\begin{equation}\label{uN1}
			\begin{aligned}
				\|u_{N,1,x_0}(t,x)|_{x_0=x}\|_{L^\infty}&+\sup_\alpha t^{\frac{m+\kappa}{s}}\frac{\|\left(\delta_\alpha \nabla^{m+[\kappa]}u_{N,1,x_0}(t,x)\right)|_{x_0=x}\|_{L^\infty}}{|\alpha|^{\kappa-[\kappa]}}\\
				&\quad\quad\quad\lesssim \tilde\M^\beta\sup_{t\in[0,T]}(t^{\frac{\kappa}{s}}\|f(t)\|_{\dot C^{\kappa_0}}+t^{\frac{m+\kappa }{s}}\|f(t)\|_{\dot C^{\kappa_0+m}}).
			\end{aligned}
		\end{equation}
		On the other hand, for
		\begin{equation*}
			\begin{aligned}
				u_{N,2,x_0}(t,x)&=\int_0^t\int _{\mathbb{R}^d} \K_{x_0}(t,\tau,x-z)g(\tau,z)dzd\tau,
			\end{aligned}
		\end{equation*}
		by Lemma \ref{esthtker} it is easy to check that 
		\begin{align*}
			|u_{N,2,x_0}(t,x)|\lesssim \int_0^t \|\K_{x_0}(t,\tau)\|_{L^1}\|g(\tau)\|_{L^\infty}d\tau\overset{\eqref{kerL1}}{\lesssim} \tilde\M^\beta  \|g\|_{L^1_TL^\infty}.
		\end{align*}
		For H\"{o}lder norm, by \eqref{delKL1} we have
		\begin{equation}\label{deun2}
			\begin{aligned}
				&|\delta_\alpha\nabla^{m+[\kappa]} u_{N,2,x_0}(t,x)|\\
				&\lesssim\int_0^\frac{t}{2}\int _{\mathbb{R}^d} \left|\delta_\alpha\nabla^{m+[\kappa]}\K_{x_0}(t,\tau,x-z)\right| |g(\tau,z)|dzd\tau+\int_\frac{t}{2}^t\int _{\mathbb{R}^d} \left|\delta_\alpha\nabla^{[\kappa]} \K_{x_0}(t,\tau,x-z)\right| \left|\nabla^mg(\tau,z)\right|dzd\tau\\
				&\lesssim \int_0^\frac{t}{2}\frac{\tilde\M^\beta}{(t-\tau)^{\frac{m+[\kappa]}{s}}}\min\left\{1,\frac{|\alpha|}{(t-\tau)^{\frac{1}{s}}}\right\}\|g(\tau)\|_{L^\infty}d\tau+\int_{\frac{t}{2}}^{t}\frac{\tilde\M^\beta}{(t-\tau)^{\frac{[\kappa]}{s}}}\min\left\{1,\frac{|\alpha|}{(t-\tau)^{\frac{1}{s}}}\right\}\|g(\tau)\|_{\dot C^m}d\tau\\
				&\lesssim \tilde\M^\beta|\alpha|^{\kappa-[\kappa]} t^{-\frac{m+\kappa }{s}}\left(\|g\|_{L^1_TL^\infty}+\sup_{\tau'\in[0,T]}\tau'^{\frac{m}{s}+1}\|g(\tau')\|_{\dot C^m}\right).
			\end{aligned}
		\end{equation}
		Hence 
		\begin{equation}\label{uN2}
			\begin{aligned}
				\|u_{N,2,x_0}(t,x)|_{x_0=x}\|_{L^\infty}+&\sup_{\alpha} t^{\frac{m+\kappa }{s}}\frac{\|\delta_\alpha\nabla^{m+[\kappa]} u_{N,2,x_0}(t,x)|_{x_0=x}\|_{L^\infty}}{|\alpha|^{\kappa-[\kappa]}}\\
				&\quad\quad\lesssim \tilde\M^\beta\left(\|g\|_{L^1_TL^\infty}+\sup_{\tau'\in[0,t]}\tau'^{\frac{m}{s}+1}\|g(\tau')\|_{\dot C^m}\right).
			\end{aligned}
		\end{equation}
		So by \eqref{uN1} and \eqref{uN2}, we get 
		\begin{equation}\label{uNf}
			\begin{aligned}
				&\sup_{t\in[0,T]}\left( \|u_{N,x_0}(t)|_{x_0=x}\|_{L^\infty}+ t^\frac{m+\kappa }{s}\sup_\alpha\frac{\|\left(\delta_\alpha\nabla^{m+[\kappa]}u_{N,x_0}(t)\right)|_{x_0=x}\|_{L^\infty}}{|\alpha|^{\kappa-[\kappa]}} \right) \lesssim \tilde\M^\beta\da_T(f,g).
			\end{aligned}
		\end{equation}
		Considering the fact that $u_{R,x_0}\equiv 0$ for any $x_0\in\mathbb{R}^d$, if the coefficient $\A(t,x,\xi)$ is independent of $x$, we complete the proof of Theorem \ref{thmcons}, where \eqref{main1} follows from \eqref{1.1infty}, \eqref{1.1hol}, \eqref{uLf}, \eqref{uNf}, and \eqref{111} follows from \eqref{1.1infty}, \eqref{1.1hol}, \eqref{bes1}, \eqref{uNf}.
		\vspace{0.2cm}\\
		In the case where $\A(t,x,\xi)$ varies with $x$, it remains to consider
		$$u_{R,{x_0}}(t,x)= -\int_0^t \int_{\mathbb{R}^d}\K_{x_0}(t,\tau,x-z)R_{x_0}[u](\tau,z)dzd\tau,$$
		where 
		\begin{align}\label{defRx0}
			R_{x_0}[u](t,x)=(2\pi)^{-\frac{d}{2}}\int _{\mathbb{R}^d}(\A(t,x_0,\xi)-\A(t,x,\xi))\hat{u}(t,\xi)e^{ix\cdot\xi} d\xi.
		\end{align}
		First we estimate $\|(u_{R,x_0}(t))|_{x_0=x}\|_{L^\infty}$. Note that $\widehat{\Lambda^{\sigma}u}(\tau,\xi)=|\xi|^\sigma \hat u(\tau,\xi)$, where $\Lambda$ is the fractional Laplacian operator defined in \eqref{deffracla}. Let $\sigma\in(\max\{s-d,0\},s)$, we can write
		\begin{equation*}
			\begin{aligned}
				u_{R,x_0}(t,x)&=(2\pi)^{-\frac{d}{2}}\int_0^t\iint_{\mathbb{R}^{2d}}\K_{x_0}(t,\tau,x-z)\frac{\A(\tau,x_0,\xi)-\A(\tau,z,\xi)}{|\xi|^{\sigma}}\widehat{\Lambda^{\sigma}u}(\tau,\xi)e^{iz\cdot \xi}dzd\xi d\tau\\
				&=(2\pi)^{-\frac{d}{2}}\int_0^t\int_{\mathbb{R}^d} \tilde{\K}_{x_0}^\sigma(t,\tau,x,\xi)\widehat{\Lambda^{\sigma}u}(\tau,\xi)e^{ix\cdot\xi}d\xi d\tau,
			\end{aligned}
		\end{equation*}
		where
		\begin{equation*}
			\tilde{\K}_{x_0}^\sigma(t,\tau,x,\xi)=\int_{\mathbb{R}^{d}}\K_{x_0}(t,\tau,z)\frac{\A(\tau,x_0,\xi)-\A(\tau,x-z,\xi)}{|\xi|^{\sigma}}e^{-iz\cdot \xi}dz.
		\end{equation*}
		Denote 
		\begin{align*}
			\J_{x_0}(t,\tau,x,y)=\mathcal{F}^{-1}(\tilde{\K}_{x_0}^\sigma(t,\tau,x,\cdot))(y).
		\end{align*}
		Then 
		\begin{align*}
			u_{R,x_0}(t,x)=\int_0^t\int_{\mathbb{R}^d}\J_{x_0}(t,\tau,x,y)\Lambda^\sigma u(\tau,x-y)dy d\tau.
		\end{align*}
		We need to prove the condition \eqref{condb}. In general, we denote
		\begin{equation*}
			\J_{x_0}^{a,b}(t,\tau,x,y)=\mathcal{F}^{-1}(\tilde{\K}^{\sigma,a,b}_{x_0}(t,\tau,x,\cdot))(y)=\mathcal{F}^{-1}(\xi^{\otimes b}\nabla_{x}^a\tilde{\K}^\sigma_{x_0}(t,\tau,x,\cdot))(y),\ \ a,b\in\mathbb{N}, a,b\leq d+m+s+2.
		\end{equation*}
		In the following, take $n=d+[b+s-\sigma]+1$. \\
		For the case of $a\neq 0$, we denote the operator $\mathcal{D}_z(f)=\frac{\xi\cdot\nabla_z}{i|\xi|^2}f$, then $\mathcal{D}_z(e^{iz\cdot \xi})=e^{iz\cdot\xi}$. By  \eqref{condop} and \eqref{kerL1}, we have for $j\in\{0,n\}$,
		\begin{equation*}
			\begin{aligned}
				\left|\tilde{\K}_{x_0}^{\sigma,a,b}(t,\tau,x,\xi)\right|&\lesssim \left|\int_{\mathbb{R}^d}\left(\xi^{\otimes b}\K_{x_0}(t,\tau,z)\frac{\nabla_x^{a}\A(\tau,x-z,\xi)}{|\xi|^{\sigma}}\right)\mathcal{D}_z^j(e^{iz\cdot\xi})dz\right|\\
				&\lesssim \left|\int_{\mathbb{R}^d}\xi^{\otimes b}\mathcal{D}_z^j\left(\K_{x_0}(t,\tau,z)\frac{\nabla_x^a\A(\tau,x-z,\xi)}{|\xi|^{\sigma}}\right)e^{iz\cdot\xi}dz\right|\\
				&\lesssim |\xi|^{b-j}\sum_{0\leq l\leq j}\int_{\mathbb{R}^d}\left|(\nabla_z^l\K_{x_0}(t,\tau,z))\right| \left|\nabla_x^{j-l+a}\left(\frac{\A(\tau,x-z,\xi)}{|\xi|^{\sigma}}\right)\right|dz\\
				&\lesssim \tilde\M^\beta|\xi|^{b+s-\sigma-j}(t-\tau)^{-\frac{j}{s}}.
			\end{aligned}
		\end{equation*}
		This implies that 
		\begin{align}\label{a0K}
			&\left|\tilde{\K}_{x_0}^{\sigma,a,b}(t,\tau,x,\xi)\right|\lesssim \tilde\M^\beta|\xi|^{b+s-\sigma}\min\{1,(t-\tau)^{-\frac{1}{s}}|\xi|^{-1}\}^{n},\ \ a\in \mathbb{N}, 0<a\leq d+m+s+2.
		\end{align}
		Furthermore, for $a=0$, we have for $j\in\{0,n\}$,
		\begin{equation*}
			\begin{aligned}
				&\left|\tilde{\K}_{x_0}^{\sigma,0,b}(t,\tau,x,\xi)\right|\\
				&\quad\lesssim \left|\int_{\mathbb{R}^d}\xi^{\otimes b}\left(\K_{x_0}(t,\tau,z)\left(\frac{\A(\tau,x_0,\xi)-\A(\tau,x,\xi)}{|\xi|^{\sigma}}+\frac{\A(\tau,x,\xi)-\A(\tau,x-z,\xi)}{|\xi|^{\sigma}}\right)\right)\mathcal{D}_z^j(e^{iz\cdot\xi})dz\right|\\
				&\quad\lesssim\left|\int_{\mathbb{R}^d}\xi^{\otimes b}\mathcal{D}_z^j\left(\K_{x_0}(t,\tau,z)\left(\frac{\A(\tau,x_0,\xi)-\A(\tau,x,\xi)}{|\xi|^{\sigma}}+\frac{\A(\tau,x,\xi)-\A(\tau,x-z,\xi)}{|\xi|^{\sigma}}\right)\right)e^{iz\cdot\xi}dz\right|\\
				&\quad\lesssim |\xi|^{b-j}\int_{\mathbb{R}^d}\left|\nabla_z^j\K_{x_0}(t,\tau,z)\right|\left|\frac{\A(\tau,x_0,\xi)-\A(\tau,x,\xi)}{|\xi|^{\sigma}}\right|dz\\
				&\quad\quad+|\xi|^{b-j}\sum_{0\leq l\leq j}\int_{\mathbb{R}^d}\left|\nabla_z^l\K_{x_0}(t,\tau,z)\right|\left|\nabla_z^{j-l}\left(\frac{\A(\tau,x,\xi)-\A(\tau,x-z,\xi)}{|\xi|^{\sigma}}\right)\right|dz.
			\end{aligned}
		\end{equation*}
		Using \eqref{kerL1} and the fact that $\frac{|\A(\tau,x,\xi)-\A(\tau,y,\xi)|}{|\xi|^\sigma}\lesssim \M |x-y|^\eps|\xi|^{s-\sigma}$,  we have
		\begin{equation*}
			\int_{\mathbb{R}^d}\left|(\nabla_z^j\K_{x_0}(t,\tau,z))\right|\left|\frac{\A(\tau,x_0,\xi)-\A(\tau,x,\xi)}{|\xi|^{\sigma}}\right|dz\lesssim \tilde\M^\beta|x_0-x|^\eps(t-\tau)^{-\frac{j}{s}}|\xi|^{s-\sigma}.
		\end{equation*}
		Moreover, 
		\begin{equation*}
			\begin{aligned}
				&\sum_{0\leq l\leq j}\int_{\mathbb{R}^d}\left|\nabla_z^l\K_{x_0}(t,\tau,z)\right|\left|\nabla_z^{j-l}\left(\frac{\A(\tau,x,\xi)-\A(\tau,x-z,\xi)}{|\xi|^{\sigma}}\right)\right|dz\\
				&\quad\quad\lesssim\M\int_{\mathbb{R}^d}\left|\nabla_z^j\K_{x_0}(t,\tau,z)\right||z|^\eps|\xi|^{s-\sigma}dz+\sum_{1\leq l\leq j}\int_{\mathbb{R}^d}\left|(\nabla_z^l\K_{x_0}(t,\tau,z))\right|\left|\nabla_z^{j-l}\left(\frac{\A(\tau,x-z,\xi)}{|\xi|^{\sigma}}\right)\right|dz\\
				&\quad\quad\lesssim\tilde\M^\beta(t-\tau)^{-\frac{j-\eps}{s}}|\xi|^{s-\sigma},
			\end{aligned}
		\end{equation*}
		where the last inequality holds because $0<\tau<t\leq 1$ and therefore the second term in the second line is lower order term. 
		So we have
		\begin{equation}\label{an0k}
			\begin{aligned}
				&\left|\tilde{\K}_{x_0}^{\sigma,0,b}(t,\tau,x,\xi)\right|\lesssim \tilde\M^\beta\left((t-\tau)^{\frac{\eps}{s}}+|x_0-x|^\eps\right)|\xi|^{b+s-\sigma}\min\{1,(t-\tau)^{-\frac{1}{s}}|\xi|^{-1}\}^{n}.
			\end{aligned}
		\end{equation}
		Similarly, for $\nabla_\xi^k\tilde{\K}_{x_0}^{\sigma,a,b}(t,\tau,x,\xi)|_{x_0=x}$, for the case when $a\neq 0$, we have
		\begin{equation}\label{an0nak}
			\begin{aligned}
				&\left|(\nabla_\xi^k\tilde{\K}_{x_0}^{\sigma,a,b}(t,\tau,x,\xi))\right|\lesssim \left|\sum_{l=0}^k\int_{\mathbb{R}^d}z^{\otimes(k-l)}\left(\K_{x_0}(t,\tau,z)\nabla_\xi^l\left(\xi^{\otimes b}\frac{\nabla_x^a\A(\tau,x-z,\xi)}{|\xi|^{\sigma}}\right)\right)\mathcal{D}_z^{k-l+n}(e^{iz\cdot\xi})dz\right|\\
				&\quad\quad\quad\quad\lesssim \sum_{l=0}^k\int_{\mathbb{R}^d}|\xi|^{l-k-n}\left|\nabla_z^{k-l+n}\left(z^{\otimes(k-l)}\K_{x_0}(t,\tau,z)\nabla_\xi^l\left(\xi^{\otimes b}\frac{\nabla_x^a\A(\tau,x-z,\xi)}{|\xi|^{\sigma}}\right)\right)\right|dz\\
				&\quad\quad\quad\quad\lesssim\tilde \M^\beta |\xi|^{s-\sigma-k-n}(t-\tau)^{-\frac{n}{s}}.
			\end{aligned}
		\end{equation}
		For the case when $a=0$, since $|\nabla_\xi^l(\xi^{\otimes b}\frac{\A(\tau,x,\xi)-\A(\tau,x-z,\xi)}{|\xi|^{\sigma}})|\leq \M |\xi|^{b+s-\sigma-l}|z|^\eps$, we have
		\begin{equation}\label{a0naK}
			\begin{aligned}
				&\left|(\nabla_\xi^k\tilde{\K}_{x_0}^{\sigma,0,b}(t,\tau,x,\xi))\right|\\
				&\lesssim \left|\sum_{l=0}^k\int_{\mathbb{R}^d}z^{\otimes(k-l)}\left(\K_{x_0}(t,\tau,z)\nabla_\xi^l\left(\xi^{\otimes b}\frac{(\A(\tau,x_0,\xi)-\A(\tau,x-z,\xi))}{|\xi|^{\sigma}}\right)\right)\mathcal{D}_z^{k-l+n}(e^{iz\cdot\xi})dz\right|\\
				&\lesssim \left|\sum_{l=0}^k\int_{\mathbb{R}^d}|\xi|^{l-k-n}\nabla_z^{k-l+n}\left(z^{\otimes(k-l)}\K_{x_0}(t,\tau,z)\right)\nabla_\xi^l\left(\xi^{\otimes b}\frac{\A(\tau,x_0,\xi)-\A(\tau,x,\xi)}{|\xi|^{\sigma}}\right)e^{iz\cdot\xi}dz\right|\\
				&\quad+\left|\sum_{l=0}^k\int_{\mathbb{R}^d}|\xi|^{l-k-n}\nabla_z^{k-l+n}\left(z^{\otimes(k-l)}\K_{x_0}(t,\tau,z)\nabla_\xi^l\left(\xi^{\otimes b}\frac{\A(\tau,x,\xi)-\A(\tau,x-z,\xi)}{|\xi|^{\sigma}}\right)\right)e^{iz\cdot\xi}dz\right|\\
				&\lesssim\tilde \M^\beta |\xi|^{s-\sigma-k-n}(t-\tau)^{-\frac{n}{s}}\left((t-\tau)^{\frac{\eps}{s}}+|x_0-x|^{\eps}\right).
			\end{aligned}
		\end{equation}
		We conclude from \eqref{a0K}, \eqref{an0k}, \eqref{an0nak} and \eqref{a0naK} that 
		\begin{equation*}
			\begin{aligned}
				&|\nabla_\xi^k\tilde{\K}_{x_0}^{\sigma,a,b}(t,\tau,x,\xi)|\lesssim\tilde \M^\beta |\xi|^{b+s-\sigma-k}\min\{1,(t-\tau)^{-\frac{1}{s}}|\xi|^{-1}\}^n,\\
				&|\nabla_\xi^k\tilde{\K}_{x_0}^{\sigma,0,b}(t,\tau,x,\xi)|\lesssim\tilde\M^\beta \left((t-\tau)^{\frac{\eps}{s}}+|x_0-x|^{\eps}\right)|\xi|^{b+s-\sigma-k}\min\{1,(t-\tau)^{-\frac{1}{s}}|\xi|^{-1}\}^n.
			\end{aligned}
		\end{equation*}
		By Lemma \ref{esthtker}, we have
		\begin{equation*}
			\begin{aligned}
				&\left|\mathcal{F}^{-1}(\tilde{\K}_{x_0}^{\sigma,a,b}(t,\tau,x,\cdot))(y) \right|\lesssim\tilde \M^{\beta}\frac{1}{((t-\tau)^{\frac{1}{s}}+|y|)^{d+b+s-\sigma}},\\
				&\left|\mathcal{F}^{-1}(\tilde{\K}_{x_0}^{\sigma,0,b}(t,\tau,x,\cdot))(y) \right|\lesssim\tilde \M^{\beta}\frac{(t-\tau)^{\frac{\eps}{s}}+|x_0-x|^\eps}{((t-\tau)^{\frac{1}{s}}+|y|)^{d+b+s-\sigma}}.
			\end{aligned}
		\end{equation*}
		So we have
		\begin{equation}\label{mthkl1}
			\begin{aligned}
				&\|J_{x_0}^{a,b}(t,\tau,x,y)\|_{L_y^1}\lesssim  \tilde\M^\beta (t-\tau)^{-\frac{b+s-\sigma}{s}},\\
				&\|\J_{x_0}^{0,b}(t,\tau,x,y)\|_{L_y^1}\lesssim  \tilde\M^\beta (t-\tau)^{-\frac{b+s-\sigma}{s}}\left((t-\tau)^{\frac{\eps}{s}}+|x-x_0|^{\eps}\right).
			\end{aligned}
		\end{equation}
		Also, by Lemma \ref{esthtker}, we have
		\begin{equation}\label{mthkh}
			\begin{aligned}
				&\int_{\mathbb{R}^d}|\delta_\alpha^y\J_{x_0}^{a,b}(t,\tau,x,y)||y|^{\sigma'}dy\lesssim \tilde\M^\beta\frac{|\alpha|^{\sigma'}}{(t-\tau)^\frac{b+s-\sigma}{s}}\min\left\{1,\frac{|\alpha|}{(t-\tau)^{\frac{1}{s}}} \right\}^{1-\sigma'},\\
				&\int_{\mathbb{R}^d}|\delta_\alpha^y\J_{x_0}^{0,b}(t,\tau,x,y)||y|^{\sigma'}dy\lesssim \tilde\M^\beta\frac{|\alpha|^{\sigma'}}{(t-\tau)^\frac{b+s-\sigma}{s}}\left((t-\tau)^{\frac{\eps}{s}}+|x_0-x|^{\eps}\right)\min\left\{1,\frac{|\alpha|}{(t-\tau)^{\frac{1}{s}}} \right\}^{1-\sigma'},
			\end{aligned}
		\end{equation}
		for any $\sigma'<b+s-\sigma$. By \eqref{mthkl1} with $a=b=0$, we have
		\begin{equation}\label{uRlinf}
			\sup_{t\in[0,T]}\|(u_{R,x_0})|_{x_0=x}(t)\|_{L^\infty}\lesssim \tilde\M^{\beta}\sup_{t\in[0,T]}\int_0^t (t-\tau)^{-\frac{s-\sigma+\eps}{s}}s^{-\frac{\sigma}{s}}ds\sup_{t\in[0,T]}t^{\frac{\sigma}{s}}\|\Lambda^\sigma u(t)\|_{L^\infty}\lesssim\tilde\M^{\beta} T^{\frac{\eps}{s}}\|u\|_{T}.
		\end{equation}
		We need to estimate $\dot C^{m+\kappa}$ norm. To do so, we write
		\begin{equation*}
			\begin{aligned}
				\delta_\alpha\nabla^{m+[\kappa]} u_{R,{x_0}}(t,x)&=\delta_\alpha\nabla^{m+[\kappa]}\left(\int_0^\frac{t}{2} \int_{\mathbb{R}^d} \J_{x_0}(t,\tau,x,z)\Lambda^\sigma u(\tau,x-z)dzd\tau\right)\\&\quad\quad+\delta_\alpha\nabla^{m+[\kappa]}\left(\int_\frac{t}{2}^t \int_{\mathbb{R}^d} \J_{x_0}(t,\tau,x,z)\Lambda^\sigma u(\tau,x-z)dzd\tau\right)\\
				&=u_{R,{x_0}}^{1,1}(t,x)+u_{R,{x_0}}^{1,2}(t,x)+u_{R,{x_0}}^{2,1}(t,x)+u_{R,{x_0}}^{2,2}(t,x).
			\end{aligned}
		\end{equation*}
		where 
		\begin{align*}
			& u_{R,{x_0}}^{1,1}(t,x)=\nabla^{m+[\kappa]}\left(\int_0^\frac{t}{2} \int_{\mathbb{R}^d} \delta_\alpha^x\J_{x_0}(t,\tau,x,z)\Lambda^\sigma u(\tau,x-z)dzd\tau\right),\\&
			u_{R,{x_0}}^{1,2}(t,x)=\nabla^{m+[\kappa]}\left(\int_0^\frac{t}{2} \int_{\mathbb{R}^d} \J_{x_0}(t,\tau,x,z)\delta_\alpha^x\Lambda^\sigma u(\tau,x-z)dzd\tau\right),\\&
			u_{R,{x_0}}^{2,1}(t,x)=\nabla^{m+[\kappa]}\left(\int_\frac{t}{2}^t \int_{\mathbb{R}^d} \delta_\alpha^x\J_{x_0}(t,\tau,x,z)\Lambda^\sigma u(\tau,x-z)dzd\tau\right),\\&
			u_{R,{x_0}}^{2,2}(t,x)=\nabla^{m+[\kappa]}\left(\int_\frac{t}{2}^t \int_{\mathbb{R}^d} \J_{x_0}(t,\tau,x,z)\delta_\alpha^x\Lambda^\sigma u(\tau,x-z)dzd\tau\right).
		\end{align*}
		For $u_{R,x_0}^{1,1}$, we have 
		\begin{equation*}
			\begin{aligned}
				\left|u_{R,x_0}^{1,1}\right|(t,x)&\lesssim\sum_{j=0}^{m+[\kappa]}\left|\int_0^{\frac{t}{2}}\int_{\mathbb{R}^d}\nabla_x^{m+[\kappa]-j}\delta_\alpha^x\J_{x_0}(t,\tau,x,z)\nabla_x^j\Lambda^{\sigma}u(\tau,x-z)dzd\tau\right|\\
				&\lesssim\sum_{j=0}^{m+[\kappa]}\left|\int_0^{\frac{t}{2}}\int_{\mathbb{R}^d}\nabla_x^{m+[\kappa]-j}\nabla_z^{j}\delta_\alpha^x\J_{x_0}(t,\tau,x,z)\Lambda^{\sigma}u(\tau,x-z)dzd\tau\right|\\
				&\lesssim\sum_{j=0}^{m+[\kappa]}\left|\int_0^{\frac{t}{2}}\int_{\mathbb{R}^d}\delta_\alpha^x\J_{x_0}^{m+[\kappa]-j,j}(t,\tau,x,z)\Lambda^{\sigma}u(\tau,x-z)dzd\tau\right|\\
				&\lesssim \tilde\M^{\beta}|\alpha|^{\kappa-[\kappa]}t^{-\frac{m+[\kappa]}{s}}(1+t^{\frac{\eps}{s}}+|x_0-x|^\eps)\|u\|_{X_T},
			\end{aligned}
		\end{equation*}
		where we use interpolation and \eqref{mthkl1} in the last inequality. Similarly, for $u_{R,x_0}^{1,2}$, we have
		\begin{equation*}
			\begin{aligned}
				\left|u_{R,x_0}^{1,2}\right|(t,x)&\lesssim\sum_{j=0}^{m+[\kappa]}\left|\int_0^{\frac{t}{2}}\int_{\mathbb{R}^d}\nabla_x^{m+[\kappa]-j}\J_{x_0}(t,\tau,x,z)\nabla_x^j\delta_\alpha^z\Lambda^{\sigma}u(\tau,x-z)dzd\tau\right|\\
				&\lesssim\sum_{j=0}^{m+[\kappa]}\left|\int_0^{\frac{t}{2}}\int_{\mathbb{R}^d}\nabla_x^{m+[\kappa]-j}\nabla_z^{j}\delta_\alpha^z\J_{x_0}(t,\tau,x,z)\Lambda^{\sigma}u(\tau,x-z)dzd\tau\right|\\
				&\lesssim\sum_{j=0}^{m+[\kappa]}\left|\int_0^{\frac{t}{2}}\int_{\mathbb{R}^d}\delta_\alpha^z\J_{x_0}^{m+[\kappa]-j,j}(t,\tau,x,z)\left(\Lambda^{\sigma}u(\tau,x-z)-\Lambda^{\sigma}u(\tau,x)\right)dzd\tau\right|\\
				&\lesssim \tilde\M^{\beta}|\alpha|^{\kappa-[\kappa]}t^{-\frac{m+\kappa}{s}}(t^{\frac{\eps}{s}}+|x_0-x|^\eps)\|u\|_{X_T},
			\end{aligned}
		\end{equation*}
		where we apply \eqref{mthkh} in the last inequality. For $u_{R,x_0}^{2,1}$ we have
		\begin{equation*}
			\begin{aligned}
				\left|u_{R,x_0}^{2,1}\right|(t,x)&\lesssim\sum_{j=0}^{m+[\kappa]}\left|\int_{\frac{t}{2}}^t\int_{\mathbb{R}^d}\nabla_x^{m+[\kappa]-j}\delta_\alpha^x\J_{x_0}(t,\tau,x,z)\nabla_x^j\Lambda^{\sigma}u(\tau,x-z)dzd\tau\right|\\
				&\lesssim\sum_{j=0}^{m+[\kappa]}\left|\int_0^{\frac{t}{2}}\int_{\mathbb{R}^d}\delta_\alpha^x\J_{x_0}^{m+[\kappa]-j,0}(t,\tau,x,z)\nabla_x^j\Lambda^{\sigma}u(\tau,x-z)dzd\tau\right|\\
				&\lesssim \tilde\M^{\beta}|\alpha|^{\kappa-[\kappa]}t^{-\frac{m+[\kappa]}{s}}(1+t^{\frac{\eps}{s}}+|x_0-x|^\eps)\|u\|_{X_T},
			\end{aligned}
		\end{equation*}
		where we use interpolation and \eqref{mthkl1} in the last inequality. For the last term $u_{R,x_0}^{2,2}$, we have
		\begin{equation*}
			\begin{aligned}
				\left|u_{R,x_0}^{2,2}\right|(t,x)&\lesssim\sum_{j=0}^{m+[\kappa]}\left|\int_{\frac{t}{2}}^t\int_{\mathbb{R}^d}\nabla_x^{m+[\kappa]-j}\J_{x_0}(t,\tau,x,z)\nabla_x^j\delta_\alpha^z\Lambda^{\sigma}u(\tau,x-z)dzd\tau\right|\\
				&\lesssim\sum_{j=0}^{m+[\kappa]}\left|\int_0^{\frac{t}{2}}\int_{\mathbb{R}^d}\delta_\alpha^z\J_{x_0}^{m+[\kappa]-j,0}(t,\tau,x,z)\nabla_x^j\Lambda^{\sigma}u(\tau,x-z)dzd\tau\right|\\
				&\lesssim \tilde\M^{\beta}|\alpha|^{\kappa-[\kappa]}t^{-\frac{m+\kappa}{s}}(t^{\frac{\eps}{s}}+|x_0-x|^\eps)\|u\|_{X_T},
			\end{aligned}
		\end{equation*}
		where we apply \eqref{mthkh} in the last inequality.\\
		Combining all the estimates for $u_{R,x_0}^{i,j}$ and taking value at $x_0=x$, we have for any $\alpha\in\mathbb{R}^d$,
		\begin{equation}\label{uRHol}
			\|(\delta_\alpha\nabla^{m+[\kappa]}u_{R,x_0})|_{x_0=x}\|_{L^\infty}\lesssim  \tilde\M^{\beta}t^{-\frac{m+\kappa-\eps}{s}}|\alpha|^{\kappa-[\kappa]}\|u\|_{X_T},\quad t\leq T.
		\end{equation}
		We conclude from \eqref{1.1infty}, \eqref{uLf}, \eqref{uNf}, \eqref{uRlinf} and \eqref{uRHol} that
		\begin{align*}
			\sup_{t\in[0,T]}\|u(t)\|_{L^\infty}&\lesssim \tilde\M^\beta\Bigg(\|u_0\|_{L^\infty}+ \sup _{t \in[0, T]} (t^{\frac{\kappa}{s}}\|f(t)\|_{\dot{C}^{\kappa_0}}+t^{\frac{m+\kappa }{s}}\|f(t)\|_{\dot{C}^{\kappa_0+m}})\\
			&\quad\quad\quad\quad+\int_0^T\|g(\tau)\|_{L^\infty}d\tau+\sup_{t\in[0,T]}t^{\frac{m}{s}+1}\|g(t)\|_{\dot C^m}+T^\frac{\varepsilon}{s}\|u\|_{X_T}\Bigg).
		\end{align*}
		Moreover, combining \eqref{1.1hol}, \eqref{uLf}, \eqref{uNf} and \eqref{uRHol} to obtain that,
		\begin{align*}
			\sup_{t\in[0,T]} t^{\frac{m+\kappa}{s}}\sup_\alpha\frac{\|\delta_\alpha \nabla^{m+[\kappa]} u\|_{L^\infty}}{|\alpha|^{\kappa-[\kappa]}}&\lesssim \M^\beta\Bigg(\|u_0\|_{L^\infty}+ \da_T(f,g)+T^\frac{\varepsilon}{s}\|u\|_{X_T}\Bigg).
		\end{align*}
		Then we get
		\begin{equation*}
			\|u\|_{X_{T}}\leq C \tilde\M^\beta\left(\|u_0\|_{L^\infty}+ \da_{T}(f,g)+T^{\frac{\eps}{s}}\|u\|_{X_T}\right).
		\end{equation*}
		We can take $T_0=T_0(\tilde\M)$ small enough such that $\tilde\M^{\beta+10}T_0^\frac{\varepsilon}{s}\ll1$, then 
		\begin{equation}\label{eslocal}
			\begin{aligned}
				\|u\|_{X_{T_0}}\leq &C \tilde\M^\beta\left(\|u_0\|_{L^\infty}+ \da_{T_0}(f,g)\right).
			\end{aligned}
		\end{equation}~\vspace{0.2cm}\\
		{\bf Step 4: Estimate the solution: long time estimate}\\
		From \eqref{eslocal}, if we start from $\frac{T_0}{2}$, there holds 
		\begin{align*}
			\sup_{t\in[\frac{T_0}{2},\frac{3T_0}{2}]}(\|u(t)\|_{L^\infty}+(t-T_0/2)^\frac{m+\kappa}{s}\|u(t)\|_{\dot C^{m+\kappa}})\leq& C \tilde\M^\beta\left(\|u(T_0/2)\|_{L^\infty}+\da_{3T_0/2}(f,g)\right).
		\end{align*}
		Specially, for $t\in[T_0,3T_0/2]$, one has $t-T_0/2\sim t$, then we obtain 
		\begin{align*}
			\|u\|_{X_{3T_0/2}}\leq C \tilde\M^\beta\left(\|u(T_0/2)\|_{L^\infty}+\da_{3T_0/2}(f,g))\right)\leq C^2 \tilde\M^{2\beta}(\|u_0\|_{L^\infty}+\da_{3T_0/2}(f,g)).
		\end{align*}
		Repeat the procedure $n-1$ times, by Stirling formula, one has 
		\begin{align*}
			\|u\|_{X_{(n+1)T_0/2}}\leq C^n(n!)^\frac{m+\kappa}{s} \tilde\M^{n\beta}(\|u_0\|_{L^\infty}+\da_{(n+1)T_0/2}(f,g)).
		\end{align*}
		Take $n=[\frac{2T}{T_0}]+1$, then we have
		\begin{align*}
			\|u\|_{X_T}\lesssim e^{\tilde CT\log T}(\|u_0\|_{L^\infty}+\da_T(f,g)),\ \ \ \forall T>0.
		\end{align*}
		Here $\tilde C$ depends only on $\tilde\M$, $\beta$, $\kappa$ and $s$.
		This completes the proof of \eqref{main2}. Similarly, combining \eqref{1.1hol}, \eqref{bes1}, \eqref{uNf} and \eqref{uRHol}, one obtains \eqref{112}. This completes the proof of Theorem \ref{lemmain}.\vspace{0.3cm}\\
		\begin{remark}
			In contrast to general methods in other literature, we don't multiply cut-off function to the equation when freezing coefficient. Instead, we freeze constant $x_0$ in the coefficient then find the fundamental solution of the linear equation. By the fundamental solution we derive a representation formula of the solution, then we take derivatives and fix $x_0=x$ in the formula of solution, and prove that the effect of freezing coefficient only leads to some remainder terms of lower order.
		\end{remark}
		In the following, we give a proof of Theorem \ref{thmexis}.\\
		\begin{proof}[Proof of Theorem \ref{thmexis}]
			To prove the existence, we use the compactness method. Consider a regularized approximation system, for which the Cauchy problem is easily studied, and whose solutions are expected to converge to solutions of the original system \eqref{eqpara}. More precisely, let $\varepsilon_1,\varepsilon_2\in(0,1)$, we  consider for $u:(0,\infty)\times \mathbb R^d\rightarrow\mathbb R^N$, \begin{equation}\label{eqpara1}
				\begin{aligned}
					&  \partial_tu(t,x)+\mathcal{L}_su(t,x)=(1-\chi_{\varepsilon_1}(t))(\mathcal{P}_\gamma f(t,x)+g(t,x))\quad \text{in}\ (0,\infty)\times\mathbb{R}^d,\\
					& u|_{t=0}=u_0\ast \rho_{\varepsilon_2}.
				\end{aligned}
			\end{equation}
			Here $\rho_{\eps_2}$ is the standard mollifier in $\mathbb{R}^d$, $\chi_{\eps_1}$ is a smooth temporal cut-off function satisfying $\mathbf{1}_{[0,\eps_1]}\leq \chi_{\eps_1}\leq \mathbf{1}_{[0,2\eps_1]}$ and the force term is supported away from time $0$ and belongs to $C_c^\infty((0,\infty)\times\mathbb{R}^d)$, there exist $\tilde T=\tilde T(\eps_1,\eps_2)\in(0,1)$ and a unique classical solution $u_\eps=u_{\eps_1,\eps_2}\in C([0,\tilde T],L^\infty(\mathbb{R}^d))\cap L^\infty_{loc}((0,\tilde T],C^{m+\kappa }(\mathbb{R}^d))$ (see \cite{Fried} or Remark \ref{mfdhk} for a brief proof). Moreover, by Theorem \ref{lemmain}, there holds 
			\begin{align*}
				&\sup_{t\in[0,T]} (\|u_\eps(t)\|_{L^\infty}+t^\frac{m+\kappa }{s}\|u_\eps(t)\|_{\dot C^{m+\kappa }})\leq Ce^{TlogT}(\|u_0\|_{L^\infty}+ \da_T(f,g)),\quad\forall T<\tilde T.
			\end{align*}
			By continuation argument, the solution can be uniquely extended beyond time $T$. By taking $\varepsilon_1,\eps_2\to 0$, we obtain $u_{\varepsilon}$ converges to $u\in C((0,T],L^\infty(\mathbb{R}^d))\cap L^\infty_{loc}((0,T],C^{m+\kappa }(\mathbb{R}^d))$, which is a solution to the original system \eqref{eqpara}. This completes the proof of i).\vspace{0.1cm}\\
			Next we show that if $u_0\in C(\mathbb{R}^d)$ and $\lim_{t\to 0}\da_t(f,g)=0$, then 
			\begin{align*}
				\lim_{t\to 0}\|u(t)-u_0\|_{L^\infty}=0.
			\end{align*}
			By the superposition principle, we can write the solution to system \eqref{eqpara} as $u=u^1+u^2$, where $u^1,u^2$ satisfy
			\begin{align*}
				&\partial_t u^1(t,x)+\mathcal{L}_s u^1(t,x)=\mathcal{P}_\gamma f(t,x)+g(t,x),\ \ \ \ \text{in}\ (0,T)\times \mathbb{R}^d,\\
				&u^1|_{t=0}=0,
			\end{align*}
			and 
			\begin{align*}
				&\partial_t u^2(t,x)+\mathcal{L}_s u^2(t,x)=0,\ \ \ \ \text{in}\ (0,T)\times \mathbb{R}^d,\\
				&u^2|_{t=0}=u_0.
			\end{align*}
			Then apply the estimate \eqref{main1} to $u^1$, we have 
			\begin{align*}
				\|u^1\|_{L^\infty_T L^\infty}\lesssim \tilde \M^\beta \da_T(f,g).
			\end{align*}
			Hence 
			\begin{align}\label{liu1}
				\lim_{t\to 0}\|u^1(t)\|_{L^\infty}\lesssim\tilde \M^\beta\lim_{t\to 0}\da_t(f,g)=0.
			\end{align}
			For $u^2$, we further consider the approximation system with smooth initial data. For $0<\eps\ll 1$, let $u_\varepsilon^2\in C((0,T],L^\infty(\mathbb{R}^d))\cap L^\infty_{loc}((0,T],C^{m+\kappa }(\mathbb{R}^d))$ be the solution to the system
			\begin{align*}
				&\partial_t u^2_\eps(t,x)+\mathcal{L}_s u^2_\eps(t,x)=0,\ \ \ \ \text{in}\ (0,T)\times \mathbb{R}^d,\\
				&u^2_\eps|_{t=0}=u_0\ast\rho_\eps:=u_{0,\eps}.
			\end{align*}
			Then one has 
			\begin{align*}
				&\partial_t(u_\eps^2-u_{0,\eps})+\mathcal{L}_s (u^2_\eps-u_{0,\eps})=-\mathcal{L}_s u_{0,\eps},\\
				&(u^2_\eps-u_{0,\eps})|_{t=0}=0.
			\end{align*}
			By \eqref{main1}, for any $t\in[0,T]$,
			\begin{align*}
				\|u_\eps^2-u_{0,\eps}\|_{L^\infty_t L^\infty}\lesssim t\|\mathcal{L}_su_{0,\eps}\|_{L^\infty}\lesssim t \eps^{-s}\|u_0\|_{L^\infty}.
			\end{align*}
			Hence for any $\eps\in(0,1)$,
			\begin{align}\label{ue2u0}
				\lim_{t\to 0}\|u_\eps^2(t)-u_{0,\eps}\|_{L^\infty}=0.
			\end{align}
			On the other hand, one has 
			\begin{align*}
				&\partial_t (u^2-u^2_\eps)+\mathcal{L}_s(u^2-u^2_\eps)=0,\\
				&(u^2-u^2_\eps)|_{t=0}=u_0-u_{0,\eps}.
			\end{align*}
			Then 
			\begin{align}\label{usmt21}
				\|u^2-u^2_\eps\|_{L^\infty_TL^\infty}\lesssim \|u_0-u_0\ast\rho_{\eps}\|_{L^\infty}.
			\end{align}
			By \eqref{ue2u0} and \eqref{usmt21}, we have 
			\begin{align*}
				\lim_{t\to 0}\|u^2(t)-u_0\|_{L^\infty}&\leq  \lim_{t\to 0}\left( \|(u^2-u^2_\eps)(t)\|_{L^\infty}+\|u^2_\eps(t)-u_{0, \eps}\|_{L^\infty}+\|u_0-u_0\ast\rho_{\eps}\|_{L^\infty}\right)\\
				&\lesssim \|u_0-u_0\ast\rho_{\eps}\|_{L^\infty}.
			\end{align*}
			Remark that if $u_0\in C(\mathbb{R}^d)$, then 
			\begin{align*}
				\lim_{\varepsilon\to 0}\|u_0-u_0\ast\rho_{\eps}\|_{L^\infty}=0.
			\end{align*}
			Hence, we get 
			\begin{align}
				\lim_{t\to 0}\|u^2(t)-u_0\|_{L^\infty}=0.
			\end{align}
			Combining this with \eqref{liu1}, we get 
			\begin{align*}
				\lim_{t\to 0}\|u(t)-u_0\|_{L^\infty}\leq    \lim_{t\to 0}\|u^1(t)\|_{L^\infty}+ \lim_{t\to 0}\|u^2(t)-u_0\|_{L^\infty}=0.
			\end{align*}
			This completes the proof.
		\end{proof}
		\vspace{1cm}\\
		The result of Theorem \ref{lemmain} in Euclidean spaces can be extended to manifold without boundary. We consider $\mathcal{M}$ to be a smooth manifold without boundary, and assume we have local charts $\{\mathcal{U}_i, \phi_i \}_{i=1}^n$ such that for any $i$ there exists $\tilde{\mathcal{U}}_i\Subset\mathcal{U}_i$; and $\{ \tilde{\mathcal{U}}_i\}_{i=1}^n$ is an open cover of $\mathcal{M}$; and   $\phi_i:\mathcal{U}_i\rightarrow\mathbb{R}^d$ are local homeomorphisms. We take characteristic functions $\{\chi_i\}$ such that $\text{Supp}(\chi_i)\subset\mathcal{U}_i$, $\chi_i(\tilde{\mathcal{U}}_i)=1$. 
		
		Define the H\"{o}lder norm on $\mathcal{M}$ as
		\begin{equation}\label{mfdhd}
			\|u\|_{ C^{m}(\mathcal{M})}=\sum_{i=1}^n\|(u\chi_i)\circ\phi_i^{-1}\|_{ C^m(\mathbb{R}^d)}.
		\end{equation}
		We denote  $\mathcal{L}_{\mathcal{M}}^s$ a Pseudo-differential operator on the manifold $\mathcal{M}$, which is defined by 
		\begin{equation}\label{mfdcdl}
			(\chi_i\mathcal{L}_{\mathcal{M}}^su)\circ\phi_i^{-1}=\mathcal{L}_{\mathbb{R}^d}^{s,i}((\chi_iu)\circ\phi_i^{-1})+R_i(u),
		\end{equation}
		where $\mathcal{L}_{\mathbb{R}^d}^{s,i}$, $i=1,\cdots,n$ are Pseudo-differential operators in $\mathbb{R}^d$ as defined in \eqref{defop} and \eqref{condop}, and the remainder term $R_i$ satisfies that
		\begin{equation}\label{mtmfcd1}
			\begin{aligned}
				\|R_i(u)\|_{L^\infty}\lesssim C_{\mathcal{M}}\|u\|_{C^{s-\zeta_0}(\mathcal{M})},\quad\|R_i(u)\|_{\dot C^{k+\alpha}}\lesssim C_{\mathcal{M}}\|u\|_{C^{s+k+\alpha-\zeta_0}(\mathcal{M})},\quad \forall k\leq m,\alpha\in[0,1),
			\end{aligned}
		\end{equation}
		for some $\zeta_0\in (0,\gamma)$ small. Also, we denote the Pseudo-differential operator $\mathcal{P}_{\mathcal{M}}^\gamma$ as 
		\begin{equation}\label{mfdcdp}
			(\chi_i\mathcal{P}_{\mathcal{M}}^\gamma F)\circ\phi_i^{-1}=\mathcal{P}_{\mathbb{R}^d}^{\gamma,i}(F_i)+\tilde{R}_i,
		\end{equation}
		where $\mathcal{P}_{\mathbb{R}^d}^{\gamma,i}$ is defined in \eqref{defPg}, $\|F_i\|_{\dot C^{\varkappa}}\lesssim \|F\|_{C^\varkappa(\mathcal{M})}$, for any $\varkappa\geq 0$, and the remainder term $\tilde {R}_i$ satisfies that
		\begin{equation}\label{mtmfcd2}
			\|\tilde{R}_i\|_{L^\infty}\lesssim C_{\mathcal{M}}\|F\|_{C^{\gamma-\zeta_0}(\mathcal{M})},\quad \|\tilde{R}_i\|_{\dot C^{k+\alpha}}\lesssim C_{\mathcal{M}}\|F\|_{C^{k+\gamma+\alpha-\zeta_0}(\mathcal{M})},\quad \forall k\leq m,\alpha\in[0,1).
		\end{equation}
		Assume that 
		\begin{equation}\label{gardineq}
			(\mathcal{L}_{\mathcal{M}}^su,u)\geq c_0\|u\|_{H^{s_0}(\mathcal{M})}^2-c_0^{-1}\|u\|_{L^2(\mathcal{M})}^2,\quad\forall u\in C^\infty(\mathcal{M}),
		\end{equation}
		for some $s_0\in (0,\frac{s}{2}]$ and $c_0\in (0,1).$ Here $(\cdot,\cdot)$ is the inner product on $\mathcal{M}$. \\
		We consider the following equation on $\mathcal{M}$,
		\begin{equation}\label{maineqM}
			\begin{aligned}
				&\partial_tu(t,x)+\mathcal{L}_{\mathcal{M}}^su(t,x)=\mathcal{P}_{\mathcal{M}}^{\gamma}F(t,x)+G(t,x),\quad\text{in}\ (0,\infty)\times \mathcal{M},\\
				&u|_{t=0}=u_0,
			\end{aligned}
		\end{equation}
		with data $F$ and $G$. The existence and uniqueness of solution to the linear parabolic equation \eqref{maineqM} follow from the method in \cite{JJS} provided the G\aa{rding}'s inequality \eqref{gardineq}, see also \cite{Fried,Polden}.
		\begin{remark}\label{mfdhk}
			We mention an iterative method to solve \eqref{maineqM}. The solution is given by 
			\begin{align}\label{forsolma}
				u(t)=S(t,0)u_0+\int_0^t S(t,\tau)(\mathcal{P}^\gamma_{\mathcal{M}}F(\tau)+G(\tau))d\tau,
			\end{align}
			where 
			\begin{equation*}
				S(t,\tau)=\lim_{n\rightarrow \infty}\Pi_{k=1}^n\left(I+\frac{t-\tau}{n}\mathcal{L}^s_{\mathcal{M}}\left(\tau+\frac{k(t-\tau)}{n}\right)\right)^{-1}. 
			\end{equation*}
			To see this, we divide the time interval $[0,T]$ into $n$ parts, and denote $\Delta t=\frac{t-\tau}{n}$, $t_j=j\Delta t$, then consider the approximate system 
			\begin{equation*}
				\begin{aligned}
					\frac{u_{j+1}-u_j}{\Delta t}+\mathcal{L}^s_{\mathcal{M}}(t_j)u_{j+1}=\mathcal{P}^\gamma_{\mathcal{M}} F_j+G_j,\quad j=0,1,\cdots,n-1,
				\end{aligned}
			\end{equation*}
			which is equivalent to 
			$$           \frac{u_{j+1}}{\Delta t}+\mathcal{L}^s_{\mathcal{M}}(t_j)u_{j+1}=\frac{u_{j}}{\Delta t}+\mathcal{P}^\gamma_{\mathcal{M}} F_j+G_j,\quad j=0,1,\cdots,n-1.$$        Using $\frac{1}{\Delta t}\gg 1$ and   the G\aa{rding} inequality \eqref{gardineq}, we can solve inductively
			\begin{align*}
				u_{j+1}=S_ju_j+\Delta t S_j(\mathcal{P}^\gamma_{\mathcal{M}} F_j+G_j),
			\end{align*}
			where 
			\begin{align*}
				S_j=\left(\mathrm{Id}+\Delta t \mathcal{L}^s_{\mathcal{M}}(t_j)\right)^{-1}.
			\end{align*}
			This implies that 
			\begin{align*}
				u_n=\prod_{j=1}^nS_ju_0+\Delta t\sum_{j=0}^n\prod_{k=j}^nS_k(\mathcal{P}^\gamma_{\mathcal{M}} F_j+G_j).
			\end{align*}
			Take $n\to\infty$, we get \eqref{forsolma}.
		\end{remark}
		We will prove the following a priori estimate for evolution equation \eqref{maineqM}. 
		\begin{theorem}\label{thmmani}
			Let $0<T<\infty$. Assume that $u\in C([0,T],L^\infty(\mathcal{M}))\cap L_{loc}^\infty((0,T], C^{m+\kappa}(\mathcal{M}))$ is a solution of \eqref{maineqM} with $\mathcal{L}_\mathcal{M}^s$ and $\mathcal{P}_{\mathcal{M}}^\gamma$ satisfying \eqref{mfdcdl} and \eqref{mfdcdp}.\\
			i) If $u_0\in L^\infty(\mathcal{M})$, then
			\begin{equation}\label{pemfd}
				\sup_{t\in[0,T]}\left(\|u(t)\|_{L^\infty(\mathcal{M})}+t^{\frac{m+\kappa}{s}}\|u(t)\|_{\dot C^{m+\kappa}(\mathcal{M})}\right)\lesssim e^{CT\log(T+2)}(\|u_0\|_{L^\infty(\mathcal{M})}+\da'_T(F,G)),
			\end{equation}
			where
			\begin{align*}
				\da'_T(f,g)=  {\sup _{t \in[0,T]} \left(t^{\frac{\kappa}{s}}\|f(t)\|_{C^{\kappa_0}(\mathcal{M})}+t^{\frac{m+\kappa }{s}}\|f(t)\|_{C^{m+\kappa_0}(\mathcal{M})}+t^{\frac{m}{s}+1}\|g(t)\|_{  C^m(\mathcal{M})}\right)}
				+\|g\|_{L^1_TL^\infty(\mathcal{M})},
			\end{align*}
			with $\kappa_0=\kappa-s+\gamma$.\\
			ii) If $u_0\in C^n(\mathcal{M})$ with $n\leq s$, and $R_i$ satisfies \eqref{mtmfcd1} with $k\leq m+n$,
			then
			\begin{equation}\label{pehmfd}
				\sup_{t\in [0,T]}(\|u(t)\|_{C^n(\mathcal{M})}+t^{\frac{m+\kappa}{s}}\|u(t)\|_{C^{m+n+\kappa}(\mathcal{M})})\lesssim e^{CT\log(T+2)}( \|u_0\|_{C^n}+\da_{T,n}'(F,G)),
			\end{equation}
			where  
			\begin{align*}
				\da'_{T,n}(f,g)=  &\sup _{t \in[0,T]} \left(t^{\frac{\kappa}{s}}\|f(t)\|_{C^{\kappa_0+n}(\mathcal{M})}+t^{\frac{m+\kappa }{s}}\|f(t)\|_{C^{m+\kappa_0+n}(\mathcal{M})}\right)\\
				&+\|g\|_{L_T^1L^\infty(\mathcal{M})}+\sup _{t \in[0,T]}\left(t^{\frac{\kappa}{s}}\|g(t)\|_{  C^{\kappa-s+n}(\mathcal{M})}
				+t^{\frac{m+\kappa}{s}}\|g(t)\|_{  C^{m+\kappa-s+n}(\mathcal{M})}\right).
			\end{align*}
		\end{theorem}
		\begin{proof}
			Denote $v_j=(u\chi_j)\circ\phi_{j}^{-1}$, then by \eqref{mfdcdl} and \eqref{mfdcdp}, we can write the equation of $v_j$ as 
			\begin{equation*}
				\begin{aligned}
					&\partial_tv_j(t,x)+\mathcal{L}_{\mathbb{R}^d}^{s,j}v_j(t,x)=\mathcal{P}_{\mathbb{R}^d}^{\gamma,j}(F_j)+G_j(t,x)+R_j(t,x)+\tilde{R}_j(t,x),\quad\text{in}\ (0,\infty)\times\mathbb{R}^d,\\
					&v_j|_{t=0}=(u_0\chi_j)\circ\phi_j,
				\end{aligned}
			\end{equation*}
			where $\|F_j\|_{  C^{\varkappa}}\lesssim\|F\|_{C^\varkappa(\mathcal{M})}$ for any $\varkappa\geq 0$, $G_j=(G\chi_j)\circ\phi_j^{-1}$, and
			\begin{align*}
				R_j(t,x)=\mathcal{L}_{\mathbb{R}^2}^{s,j}(v_j)-(\chi_j\mathcal{L}_{\mathcal{M}}^su)\circ\phi_j^{-1},\quad\tilde R_j(t,x)=(\chi_j\mathcal{P}_{\mathcal{M}}^\gamma F)\circ\phi_j^{-1}-\mathcal{\mathcal{P}}_{\mathbb{R}^2}^{\gamma,j}(F_j).
			\end{align*}
			For simplicity, denote 
			\begin{align*}
				&\|f\|_{X_T}:=  \sup_{t\in[0,T]}\left(\|f(t)\|_{L^\infty}+t^{\frac{m+\kappa }{s}}\|f(t)\|_{\dot C^{m+\kappa }}\right),\\
				&\|f\|_{X_{T,n}}:=  \sup_{t\in[0,T]}\left(\|f(t)\|_{C^n}+t^{\frac{m+\kappa }{s}}\|f(t)\|_{\dot C^{m+n+\kappa }}\right),
			\end{align*}
			where we drop the domain of the function, which can be $\mathcal{M}$ or $\mathbb{R}^d$, whenever it is clear from context.\\
			Note that by definition \eqref{mfdhd}, 
			\begin{equation}\label{h1}
				\sum_{j=1}^n\|v_j\|_{X_{T}}\sim \|u\|_{X_{T}},\quad \sum_{j=1}^n\|v_j\|_{X_{T,n}}\sim \|u\|_{X_{T,n}},\quad\sum_{j=1}^n\|v_j(0)\|_{C^n}\sim \|u_0\|_{C^n(\mathcal{M})},
			\end{equation}
			and by the property of $F_j$,
			\begin{equation}\label{h2}
				\sum_{j=1}^n\da'_T(F_j,G_j)\sim \da'_T(F,G),\quad\sum_{j=1}^n\da'_{T,n}(F_j,G_j)\sim \da'_{T,n}(F,G).
			\end{equation}
			By Theorem \ref{lemmain}, one has for any $T_0\in (0,\min\{T,1\}]$,
			\begin{equation*}
				\|v_j\|_{X_{T_0}}\lesssim \|v_j(0)\|_{L^\infty}+\da'_{T_0}(F_j,G_j)+\da'_{T_0}(0,R_j)+\da'_{T_0}(0,\tilde{R}_j),
			\end{equation*}
			where the constant $C$ depends on the coefficient in the Pseudo-differential operator $\mathcal{L}^{s,j}_{\mathbb{R}^d}$. 
			By assumption \eqref{mtmfcd1} and \eqref{mtmfcd2}, one has
			\begin{equation}\label{mfdrt}
				\begin{aligned}
					&\|R_j\|_{L^\infty}\lesssim C_{\mathcal{M}}\|u\|_{C^{s-\zeta_0}(\mathcal{M})},\quad\|R_j\|_{\dot C^{m+\alpha}}\lesssim C_{\mathcal{M}}\|u\|_{C^{s+m+\alpha-\zeta_0}(\mathcal{M})},\\
					&\|\tilde{R}_j\|_{L^\infty}\lesssim C_{\mathcal{M}}\|F\|_{C^{\gamma-\zeta_0}(\mathcal{M})},\quad\|\tilde{R}_j\|_{\dot C^{m+\alpha}}\lesssim C_{\mathcal{M}}\|F\|_{ C^{m+\gamma+\alpha-\zeta_0}(\mathcal{M})},
				\end{aligned}
			\end{equation}
			for some $\zeta_0\in(0,\gamma)$ small and any $\alpha\in[0,1)$. So
			\begin{equation*}
				\begin{aligned}
					\da'_{T_0}(0,R_j)+\da'_{T_0}(0,\tilde{R}_j)\lesssim {T_0}^{\frac{\zeta_0}{s}}\left(\|u\|_{X_{T_0}}+\da'_{T_0}(F,G)\right).
				\end{aligned}
			\end{equation*}
			So one has
			\begin{equation*}
				\|v_j\|_{X_{T_0}}\lesssim \|v_j(0)\|_{L^\infty}+\da'_{T_0}(F_j,G_j)+{T_0}^{\frac{\zeta_0}{s}}\left(\|u\|_{X_{T_0}}+\da'_{T_0}(F,G)\right).
			\end{equation*}
			Combining this with \eqref{h1} and \eqref{h2} to get  the estimate of $u$,
			\begin{align}\label{manile}
				\|u\|_{X_{T_0}}\lesssim \|u_0\|_{L^\infty(\mathcal{M})}+\da'_{T_0}(F,G)+T_0^{\frac{\zeta_0}{s}}\|u\|_{X_{T_0}}.
			\end{align}
			We can take $T_0$ small enough such that the remainder term in the right hand side can be absorbed by the left hand side. Then by a similar argument as Step 4 in the proof of Theorem \ref{lemmain}, we get 
			\begin{align}\label{mfdetm}
				\|u\|_{X_{T}}\lesssim e^{C'T\log(T+2)}(\|u_0\|_{L^\infty(\mathcal{M})}+\da'_T(F,G)),\ \ \ \forall T>0.
			\end{align}
			This completes the proof for \eqref{pemfd}.\\
			Then we prove \eqref{pehmfd}. Denote $v_j^{n}:=\nabla^nv_j$, then the equation of $v_j^{n}$ can be written as
			\begin{equation*}
				\begin{aligned}
					&\partial_tv_j^{n}(t,x)+\mathcal{L}_{\mathbb{R}^d}^{s,j}v_j^{n}(t,x)=\nabla^n\mathcal{P}_{\mathbb{R}^d}^{\gamma,j}(F_j)+\nabla^n(G_j(t,x)+R_j(t,x)+\tilde{R}_j(t,x))+R_j^n,\quad\text{in}\ (0,\infty)\times\mathbb{R}^d\\
					&v_j^{n}|_{t=0}=\nabla^n((u_0\chi_j)\circ\phi_j^{-1}),
				\end{aligned}
			\end{equation*}
			where
			\begin{equation*}
				R_j^n=\mathcal{L}_{\mathbb{R}^d}^{s,j}v_j^{n}(t,x)-\nabla^n(\mathcal{L}_{\mathbb{R}^d}^{s,j}v_j).
			\end{equation*}
			Let $T_0\in (0,\min\{T,1\}]$.
			By \eqref{condop}, we have
			\begin{equation}\label{mfdrt2}
				\|R_j^n\|_{L_{T_0}^1L^\infty}+\sup_{t\in[0,T_0]}t^{\frac{m}{s}+1}\|R_j^n(t)\|_{\dot C^{m}}\lesssim T_0^{\frac{\zeta_0}{s}}\|v_j\|_{X_{T_0,n}}.
			\end{equation}
			So by Theorem \ref{lemmain}, we get
			\begin{equation*}
				\begin{aligned}
					\|v_j^{n}\|_{X_{T_0}}\lesssim&\|v_j^{n}(0)\|_{L^\infty}+\sup_{t\in[0,T_0]}(t^{\frac{\kappa}{s}}\|F_j(t)\|_{\dot C^{\kappa_0+n}}+t^{\frac{m+\kappa}{s}}\|F_j(t)\|_{\dot C^{\kappa_0+m+n}})\\
					&+\sup_{t\in[0,T_0]}(t^{\frac{\kappa}{s}}\|(G_j+R_j+\tilde{R}_j)(t)\|_{\dot C^{\kappa-s+n}}+t^{\frac{m+\kappa}{s}}\|(G_j+R_j+\tilde{R}_j)(t)\|_{\dot C^{\kappa-s+n+m}})\\
					&+\|R_j^n\|_{L_{T_0}^1L^\infty}+\sup_{t\in[0,T_0]}t^{\frac{m}{s}+1}\|R_j^n(t)\|_{\dot C^{m}}\\
					\lesssim& \|v_j^{n}(0)\|_{L^\infty}+\da'_{T_0,n}(F_j,G_j)+\da'_{T_0,n}(0,R_j)+\da'_{T_0,n}(0,\tilde{R}_j)+\da'_{T_0}(0,R_j^n).
				\end{aligned}
			\end{equation*}
			By \eqref{h2}, \eqref{mfdrt}, and the fact that
			\begin{equation*}
				\da'_{T_0,n}(0,R_j)+\da'_{T_0,n}(0,\tilde{R}_j)+\da'_{T_0}(0,R_j^n)\lesssim T_0^{\frac{\zeta_0}{s}}\left(\|v_j\|_{X_{T_0,n}}+\da'_{T_0,n}(F,G)\right),
			\end{equation*}
			we have
			\begin{equation*}
				\|v_j^{n}\|_{X_{T_0}}\lesssim\|v_j^{n}(0)\|_{L^\infty}+\da'_{T_0,n}(F,G)+T_0^{\frac{\zeta_0}{s}}\|v_j^{n}\|_{X_{T_0}}.
			\end{equation*}
			For $L^\infty$ norm, by \eqref{manile} and the fact that $\da'_{T_0}(F,G)\lesssim \da'_{T_0,n}(F,G)$, we have
			\begin{equation}\label{manihe}
				\|v_j\|_{X_{T_0}}\lesssim\|v_j(0)\|_{L^\infty}+\da'_{T_0,n}(F,G)+T_0^{\frac{\zeta_0}{s}}\|u\|_{X_{T_0}}.
			\end{equation}
			By \eqref{h1}, \eqref{mfdrt2}, \eqref{manihe}, and $\|u\|_{T_0}\lesssim\|u\|_{T_0,n}$, one has
			\begin{equation*}
				\|u\|_{X_{T_0,n}}\lesssim \|u_0\|_{C^n(\mathcal{M})}+\da'_{T,n}(F,G)+T_0^{\frac{\zeta_0}{s}}\|u\|_{X_{T_0,n}}.
			\end{equation*}
			Then by same procedure of \eqref{mfdetm}, one has
			\begin{equation*}
				\|u\|_{X_{T,n}}\lesssim e^{C'T\log(T+2)}\left(\|u_0\|_{C^n(\mathcal{M})}+\da'_{T,n}(F,G)\right),\ \ \ \forall T>0.
			\end{equation*}
			This completes the proof.
		\end{proof}\\
		
		\section{Applications}

		\subsection{Semi-linear parabolic equation}\label{secsemi}
		
		\begin{proposition}\label{propade}
			Fix $\gamma, \kappa$ such that $0<\gamma<s$, $1+s<\kappa<2$, and $m\in\mathbb{N}$. There exists $\varepsilon_0>0$ such that \\
			i)	For any $ u_0\in\dot B^{-s}_{\infty,\infty}$, if $$\| u_0\|_{\dot B^{-s}_{\infty,\infty}}\leq \varepsilon_0,
			$$ then \eqref{eqag} admits a unique global solution $u$ in $[0,\infty]$ satisfying
			\begin{align*}
				\sup_{t>0}\left(t^\frac{\gamma}{2}\|\Lambda^{-s} u(t)\|_{\dot C^{\gamma}}+t^\frac{m+\kappa}{2}\|\Lambda^{-s} u(t)\|_{\dot C^{m+\kappa}}\right)\leq C_m \| u_0\|_{\dot B^{-s}_{\infty,\infty}}.
			\end{align*}
			ii) For any initial data $ u_0\in \dot B^{-s}_{\infty,\infty}$, if there exists $\eps_1>0$ such that
			\begin{align}\label{aggcon}
				\|  u_0-u_0\ast \rho_{\varepsilon_1}\|_{\dot B^{-s}_{\infty,\infty}}\leq \varepsilon_0,
			\end{align}
			then   there exists $T=T(\eps_0,\eps_1,\|u_0\|_{\dot B^{-s}_{\infty,\infty}})>0$ such that \eqref{eqag} admits a unique solution  $u$ in $[0,T]$ satisfying
			\begin{align*}
				\sup_{t\in[0,T]}\left(t^\frac{\gamma}{2}\|\Lambda^{-s} u(t)\|_{\dot C^{\gamma}}+t^\frac{m+\kappa}{2}\|\Lambda^{-s} u(t)\|_{\dot C^{m+\kappa}}\right)\lesssim \| u_0\|_{\dot B^{-s}_{\infty,\infty}}.
			\end{align*}
		\end{proposition}
		\begin{proposition}\label{propns}
			Fix $\gamma,\kappa$ such that $0<\gamma<s-1$, $\max\{1,2s-2\}<\kappa<s$, and $m\in\mathbb{N}$. There exists $\varepsilon_0>0$ such that \\
			i)  For any initial data $u_0\in\dot B^{1-s}_{\infty,\infty}$, if
			$$
			\|u_0\|_{\dot B^{1-s}_{\infty,\infty}}\leq \varepsilon_0,
			$$
			then \eqref{eqns} admits a unique global solution $u$ in $[0,\infty]$ such that
			\begin{align*}
				\sup_{t>0}\left(t^{\frac{\gamma}{s}}\|\Lambda^{1-s}u(t)\|_{\dot{C}^{\gamma}}+t^\frac{m+\kappa}{s}\|\Lambda^{1-s}u(t)\|_{\dot C^{m+\kappa}}\right)\lesssim \|u_0\|_{\dot B^{1-s}_{\infty,\infty}}.
			\end{align*}
			ii) For any initial data $u_0\in \dot B^{1-s}_{\infty,\infty}$, if there exists $\eps_1>0$ such that
			\begin{align}\label{nscon}
				\|u_0-u_0\ast\rho_{\eps_1}\|_{\dot B^{1-s}_{\infty,\infty}}\leq\varepsilon_0,
			\end{align}
			then	there exists $T=T(\eps_0,\eps_1,\|u_0\|_{\dot B^{1-s}_{\infty,\infty}})>0$ such that \eqref{eqns} admits a unique solution $u$ in $[0,T]$ such that
			\begin{align*}
				\sup_{t\in[0,T]}\left(t^{\frac{\gamma}{s}}\|\Lambda^{1-s}u(t)\|_{\dot{C}^{\gamma}}+	t^\frac{m+\kappa}{s}\|\Lambda^{1-s}u(t)\|_{\dot C^{m+\kappa}}\right)\lesssim \|u_0\|_{\dot B^{1-s}_{\infty,\infty}}.
			\end{align*}	
		\end{proposition}
		
		\begin{proposition}\label{propsqg}
			Let $0<\gamma<s-1, \max\{1,2s-2\}<\kappa<s$, and $m\in\mathbb{N}$. There exists $\varepsilon_0>0$ such that \\
			i)  For any initial data $\theta_0\in \dot B^{1-s}_{\infty,\infty}$, satisfying $$
			\|\theta_0\|_{\dot B^{1-s}_{\infty,\infty}}\leq \varepsilon_0,
			$$
			then \eqref{eqsqg} admits a unique global solution $\theta$ in $[0,\infty]$ such that
			\begin{align*}
				\sup_{t>0}\left(t^{\frac{\gamma}{s}}\|\Lambda^{1-s}\theta(t)\|_{\dot{C}^{\gamma}}+t^\frac{m+\kappa}{s}\|\Lambda^{1-s}\theta(t)\|_{\dot C^{m+\kappa}}\right)\lesssim \|\theta_0\|_{\dot B^{1-s}_{\infty,\infty}}.
			\end{align*}
			ii) For any initial data $\theta_0 \in \dot B^{1-s}_{\infty,\infty}$, if there exists $\eps_1>0$ such that
			\begin{align}\label{sqgcon}
				\|\theta_0-\theta_0\ast\rho_{\eps_1}\|_{\dot B^{1-s}_{\infty,\infty}}\leq\varepsilon_0,
			\end{align}
			then there exists $T=T(\eps_0,\eps_1,\|\theta_0\|_{\dot B^{1-s}_{\infty,\infty}})>0$ such that \eqref{eqsqg} admits a unique solution $\theta$ in $[0,T]$ such that 
			\begin{align*}
				\sup_{t\in[0,T]}\left(t^{\frac{\gamma}{s}}\|\Lambda^{1-s}\theta(t)\|_{\dot{C}^{{\gamma}}}+t^\frac{m+\kappa}{s}\|\Lambda^{1-s}\theta(t)\|_{\dot C^{m+\kappa}}\right)\lesssim\|\theta_0\|_{\dot B^{1-s}_{\infty,\infty}}.
			\end{align*}
		\end{proposition}
		\begin{proof}[Proof of Proposition \ref{propade}]\\
			Denote 
			\begin{align*}
				\|h\|_T=\sup_{t\in[0,T]}(t^\frac{\gamma}{2}\|h(t)\|_{\dot C^\gamma}+t^\frac{m+\kappa}{2}\|h(t)\|_{\dot C^{m+\kappa}}).
			\end{align*}
			Let $\phi\in C^\infty$, $\sigma,T>0$ that will be fixed later. We consider the set 
			\begin{align*}
				\mathcal{X}^\sigma_{T,\phi}=\left\{v:\|\Lambda^{-s}v-\phi\|_T\leq \sigma,\ v|_{t=0}=u_0\right\}.
			\end{align*}
			For $v\in 	\mathcal{X}^\sigma_{T,\phi}$, we define $u=\mathcal{S}v$ as a solution to 
			\begin{align*}
				&			\partial_t u-\Delta u=\nabla \cdot ( v\nabla \Lambda^{-s} v),\\
				&u|_{t=0}=u_0.
			\end{align*}
			We apply $\Lambda^{-s}$ to both sides of the equation.
			\begin{align*}
				\partial_t \Lambda^{-s} u-\Delta\Lambda^{-s} u=\Lambda^{-s}\nabla \cdot F[v],\ \ \ \text{where}\ \ F[v]= v\nabla \Lambda^{-s} v.
			\end{align*}
			Then we have 
			\begin{align*}
				\partial_t (\Lambda^{-s} u-\phi)-\Delta(\Lambda^{-s} u-\phi)=\Lambda^{-s}\nabla \cdot( F[v]+\Lambda^s\nabla \phi). 
			\end{align*}
			Applying \eqref{111} in Theorem \ref{thmcons} with $F= F[v]+\Lambda^s\nabla \phi$ and $G=0$,  denote $\kappa_0=\kappa-1-s$, we have for any $T>0$,
			\begin{equation}\label{esagg}
				\begin{aligned}
					\|\Lambda^{-s} u-\phi\|_T\lesssim & \|\Lambda^{-s} u_0-\phi\|_{\dot B^0_{\infty,\infty}}+\sup _{t \in[0, T]} (t^{\frac{\kappa}{s}}\|F[v](t)\|_{\dot{C}^{\kappa_0}}+t^{\frac{m+\kappa}{s}}\|F[v](t)\|_{\dot{C}^{m+\kappa_0}})\\
					&+T^{\frac{\kappa}{s}}\|\Lambda^s\nabla \phi\|_{\dot{C}^{\kappa_0}}+T^{\frac{m+\kappa}{s}}\|\Lambda^s\nabla \phi\|_{\dot{C}^{m+\kappa_0}}.
				\end{aligned}
			\end{equation}
			By Lemma \ref{Hol}, for any $n\in\mathbb{N}$, $n\leq m$, 
			\begin{align*}
				\|\nabla^n F[v](t)\|_{\dot{C}^{\kappa_0}}&\lesssim \sum_{n_1+n_2=n}(\| v(t)\|_{\dot{C}^{n_1+\kappa_0}}\|\Lambda^{-s} v(t)\|_{\dot C^{n_2+1}}+\| v(t)\|_{\dot C^{n_1}}\|\Lambda^{-s}v(t)\|_{\dot C^{n_2+\kappa_0+1}})\\
				&\lesssim t^{-\frac{n+\kappa}{s}}\|\Lambda^{-s} v\|_T^2\lesssim t^{-\frac{n+\kappa}{s}}\left( \|\Lambda^{-s}  v-\phi\|_T^2+\|\phi\|_T(\|\Lambda^{-s}  v-\phi\|_T+\|\phi\|_T)\right).
			\end{align*}
			Hence
			\begin{align*}
				&		\sup _{t \in[0, T]} (t^{\frac{\kappa}{s}}\|F[v](t)\|_{\dot{C}^{\kappa_0}}+t^{\frac{m+\kappa}{s}}\|F[v](t)\|_{\dot{C}^{m+\kappa_0}})\\
				&\quad\quad\quad\lesssim \|\Lambda^{-s}  v-\phi\|_T^2+\|\phi\|_T(\|\Lambda^{-s}  v-\phi\|_T+\|\phi\|_T)\\
				&\quad\quad\quad\lesssim \sigma^2+(T^\frac{\gamma}{2}+T^\frac{m+\kappa}{2})\|\phi\|_{C^{m+2}}(\sigma+(T^\frac{\gamma}{2}+T^\frac{m+\kappa}{2})\|\phi\|_{C^{m+2}}).
			\end{align*}
			Combining this with \eqref{esagg} we obtain 
			\begin{align}\label{finagg}
				\|\Lambda^{-s} u-\phi\|_T\leq C_0( \|\Lambda^{-s} u_0-\phi\|_{\dot B^0_{\infty,\infty}}+\sigma^2+(T^\frac{\gamma}{2}+T^\frac{m+\kappa}{2})\|\phi\|_{C^{m+2}}(\sigma+(T^\frac{\gamma}{2}+T^\frac{m+\kappa}{2})\|\phi\|_{C^{m+2}})).
			\end{align}
			Moreover, for any $v_1,v_2\in \mathcal{X}^\sigma_{T,\phi}$, denote $u_1=\mathcal{S}v_1$, $u_2=\mathcal{S}v_2$, then we have 
			\begin{align*}
				&	\partial_t \Lambda^{-s}(u_1-u_2)-\Delta \Lambda^{-s}(u_1-u_2)=\Lambda^{-s}\nabla\cdot (F[v_1]-F[v_2]),\\
				&\Lambda^{-s}	(u_1-u_2)|_{t=0}=0.
			\end{align*}
			Then by Theorem \ref{thmcons} with $F=(F[v_1]-F[v_2])$ and $G=0$, one has for any $T>0$,
			\begin{align}\label{nsdif}
				\|\Lambda^{-s}(u_1-u_2)\|_T\lesssim \sup _{t \in[0, T]} (t^{\frac{\kappa}{s}}\|(F[v_1]-F[v_2])(t)\|_{\dot{C}^{\kappa_0}}+t^{\frac{m+\kappa}{s}}\|(F[v_1]-F[v_2])(t)\|_{\dot{C}^{m+\kappa_0}}).
			\end{align}
			Note that 
			\begin{align*}
				F[v_1]-F[v_2]=(v_1-v_2)\nabla\Lambda^{-s}v_1+v_2\nabla\Lambda^{-s}(v_1-v_2).
			\end{align*}
			By Lemma \ref{Hol},	for any $n\in\mathbb{N}$, $n\leq m$, 
			\begin{align*}
				&\|\nabla^n(F[v_1]-F[v_2])(t)\|_{\dot C^{\kappa_0}}\\
				&\quad\lesssim
				\sum_{n_1+n_2=n}(\| (v_1-v_2)(t)\|_{\dot{C}^{n_1+\kappa_0}}\|\Lambda^{-s} v_1(t)\|_{\dot C^{n_2+1}}+\| (v_1-v_2)(t)\|_{\dot C^{n_1}}\|\Lambda^{-s}v_1(t)\|_{\dot C^{n_2+\kappa_0+1}})\\
				&\quad\quad+\sum_{n_1+n_2=n}(\| v_2(t)\|_{\dot{C}^{n_1+\kappa_0}}\|\Lambda^{-s} (v_1-v_2)(t)\|_{\dot C^{n_2+1}}+\|v_2(t)\|_{\dot C^{n_1}}\|\Lambda^{-s}(v_1-v_2)(t)\|_{\dot C^{n_2+\kappa_0+1}})\\
				&\quad\lesssim t^{-\frac{n+\kappa}{s}}\|\Lambda^{-s}(v_1-v_2)\|_T\|\Lambda^{-s}(v_1,v_2)\|_T.
			\end{align*}
			Hence, we obtain 
			\begin{align*}
				\sup _{t \in[0, T]} (t^{\frac{\kappa}{s}}\|(F[v_1]-F[v_2])(t)\|_{\dot{C}^{\kappa_0}}+t^{\frac{m+\kappa}{s}}\|(F[v_1]-F[v_2])(t)\|_{\dot{C}^{m+\kappa_0}})\lesssim \|\Lambda^{-s}(v_1-v_2)\|_T\|\Lambda^{-s}(v_1,v_2)\|_T.
			\end{align*}
			Combining this with \eqref{nsdif}, one has 
			\begin{align}\label{sqgdif}
				\|\Lambda^{-s}(u_1-u_2)\|_T\leq  C_1\|\Lambda^{-s}(v_1-v_2)\|_T\|\Lambda^{-s}(v_1,v_2)\|_T.
			\end{align}
			With the estimates \eqref{finagg} and \eqref{sqgdif}, we can prove the global/local well-posedness of the Cauchy problem \eqref{eqsqg} in $\mathcal{X}^\sigma_{T,\phi}$ with suitable choice of $\sigma,T$ and $\phi$.\vspace{0.2cm}\\
			\textit{i) Global existence.}\\	Consider $ u_0$ with $\| u_0\|_{\dot B^{-s}_{\infty,\infty}}\leq \varepsilon_0$, where $0<\varepsilon_0<\frac{1}{100(C_0+C_1+1)^{10}}$. For $v\in\mathcal{X}^\sigma_{\infty,0}$ with $\phi\equiv 0$, $\sigma=2C_0\varepsilon_0$, and $T=\infty$, we obtain from \eqref{finagg} that 
			\begin{align*}
				\|\Lambda^{-s} u\|_\infty\leq C_0 \|\Lambda^{-s} u_0\|_{\dot B^0_{\infty,\infty}}+C_0\sigma^2\leq C_0\varepsilon_0+4C_0^3\varepsilon_0^2\leq \sigma.
			\end{align*}
			This implies that $\mathcal{S}$ maps $\mathcal{X}^\sigma_{\infty,0}$ to itself.\\
			Moreover, for $v_1,v_2\in \mathcal{X}^\sigma_{\infty,0}$, we obtain from \eqref{sqgdif} that 
			\begin{align*}
				\|\Lambda^{-s}(u_1-u_2)\|_\infty\leq  2C_1\sigma\|\Lambda^{-s}(v_1-v_2)\|_\infty\leq \frac{1}{2}\|\Lambda^{-s}(v_1-v_2)\|_T.
			\end{align*}
			Hence $\mathcal{S}$ is a contraction map from $\mathcal{X}^\sigma_{\infty,0}$ to itself.\vspace{0.2cm}\\
			\textit{ii) Local existence.}\\ If the  initial data satisfies \eqref{aggcon},  take $\phi= \Lambda^{-s}  u_0\ast \rho_{\eps_1}$,  then we have
			\begin{align*}
				\|\phi\|_T\lesssim (T^\frac{\gamma}{2} \eps_1^{-\gamma}+T^{\frac{m+\kappa}{2}}\eps_1^{-(m+\kappa)})\| u_0\|_{\dot B^{-s}_{\infty,\infty}},\ \ \ \ \|\Lambda^s \nabla \phi\|_{\dot{C}^{m+\kappa-1-s}}\lesssim \| u_0\|_{\dot B^{-s}_{\infty,\infty}}\eps_1^{-(m+\kappa)}.
			\end{align*}
			Now we fix $T_1=\left(\frac{\varepsilon_0\eps_1^{m+\kappa}}{100+C_0+\| u_0\|_{\dot B^{-s}_{\infty,\infty}}}\right)^\frac{100}{\gamma}$, then for any $v\in\mathcal{X}^\sigma_{T_1,\Lambda^{-s}  u_0\ast \rho_{\eps_1}}$,
			\eqref{finagg} leads to
			\begin{align*}
				\|\Lambda^{-s} u-\phi\|_T
				&\leq C_0\varepsilon_0^2+C_0\sigma^2+C_0T^\frac{\gamma}{2}\| u_0\|_{\dot B^{-s}_{\infty,\infty}}(2C_0\varepsilon_0+T^\frac{\gamma}{2}\| u_0\|_{\dot B^{-s}_{\infty,\infty}})+C_0T^{\frac{\kappa}{s}}\eps_1^{-2}\| u_0\|_{\dot B^{-s}_{\infty,\infty}}\\
				&\leq \sigma.
			\end{align*}
			Hence $\mathcal{S}$ maps $\mathcal{X}^\sigma_{T_1,\Lambda^{-s}  u_0\ast \rho_{\eps_1}}$ to itself.\\
			Moreover, for $v_1,v_2\in \mathcal{X}^\sigma_{T_1,\Lambda^{-s}  u_0\ast \rho_{\eps_1}}$, we have 
			\begin{align*}
				\|\Lambda^{-s}(u_1-u_2)\|_{T_1}\leq& C_2 \|\Lambda^{-s}(v_1-v_2)\|_{T_1}(\|\Lambda^{-s}v_1-\phi\|_{T_1}+\|\Lambda^{-s}v_2-\phi\|_{T_1}+\|\phi\|_{T_1})\\
				\leq &C_2(2\sigma+T_1^\frac{\gamma}{2}\eps_1^{-2}\| u_0\|_{\dot B^0_{\infty,\infty}})\|\Lambda^{-s}(v_1-v_2)\|_{T_1}\leq 3C_2\sigma \|\Lambda^{-s}(v_1-v_2)\|_{T_1}.
			\end{align*}
			We take $\varepsilon_0=\frac{1}{100(C_0+C_1+C_2+1)^{10}}$, which implies that $2C_1\sigma+3C_2\sigma\leq \frac{1}{2}$. 
			Hence we obtain that 
			\begin{align*}
				&	\|\Lambda^{-s}(u_1-u_2)\|_{T_1}\leq\frac{1}{2}\|\Lambda^{-s}(v_1-v_2)\|_{T_1}, \ \  \forall v_1,v_2\in\mathcal{X}^\sigma_{T_1,\Lambda^{-s}  u_0\ast \rho_{\eps_1}}.
			\end{align*}
			Hence we conclude that 
			\begin{itemize}
				\item If $\| u_0\|_{\dot B^{-s}_{\infty,\infty}}\leq \varepsilon_0$, then there exists a unique $u\in\mathcal{X}^\sigma_{\infty,\phi} $ such that $\mathcal{S}u=u$ with $\phi\equiv 0$.\\
				\item If \eqref{aggcon} holds, then there exists a unique $u\in\mathcal{X}^\sigma_{T_1,\phi}$ such that $\mathcal{S}u=u$ with $\phi=\Lambda^{-s}  u_0\ast \rho_{\eps_1}$.
			\end{itemize}
			This completes the proof.
			
		\end{proof}
		
		Then we prove Proposition \ref{propns}.	Applying the Leray operator $\mathbb{P}=\mathrm{Id}-\nabla \Delta^{-1}\operatorname{div}$ to both sides of \eqref{eqns}, we obtain 
		\begin{align}\label{nseq}
			\partial_t u+\Lambda^s u=-\mathbb{P}\operatorname{div}(u\otimes u),\ \ \ \operatorname{div } u=0.
		\end{align}
		\begin{proof}[Proof of Proposition \ref{propns}]
			Define
			\begin{align*}
				&\|h\|_{T}=\sup_{t\in[0,T]}\left(t^{\frac{\gamma}{s}}\|h(t)\|_{\dot{C}^{\gamma}}+t^{\frac{m+\kappa}{s}}\|h(t)\|_{\dot{C}^{m+\kappa}}\right),\quad \forall T\in(0,\infty],
			\end{align*}
			for $0<\gamma<s-1$, $\max\{1,2s-2\}<\kappa<s$.
			Let $\phi\in C^\infty$ that will be fixed later.	We consider the set 
			\begin{align*}
				\mathcal{X}^\sigma_{T,\phi}=\left\{v:\|\Lambda^{1-s}v-\phi\|_T\leq \sigma, \ v|_{t=0}=u_0\right\}.
			\end{align*}
			For any $v\in 	\mathcal{X}^\sigma_{T,\phi}$, we define $u=\mathcal{S}v$ as a solution to 
			\begin{align*}
				\partial_t u+\Lambda^s u=-\mathbb{P}\operatorname{div}(v\otimes v),\\ \operatorname{div } u=0,\ \:u|_{ t=0}=u_0.
			\end{align*}
			Applying $\Lambda^{1-s}$ to \eqref{nseq}, we obtain 
			\begin{align*}
				\partial_t(\Lambda^{1-s}u-\phi)+\Lambda^s(\Lambda^{1-s}u-\phi)=-\Lambda^{1-s}\mathbb{P}\operatorname{div}(v\otimes v)-\Lambda^s \phi.
			\end{align*}
			By Theorem \ref{thmcons} with $F=v\otimes v$ and $G=-\Lambda^s\phi$, if we denote $\kappa_0=\kappa+2-2s$, then we have for any $T>0$, 
			\begin{align*}
				\|\Lambda^{1-s}u-\phi\|_T&\lesssim  \|\Lambda^{1-s}u_0-\phi\|_{\dot B^0_{\infty,\infty}}+\sup_{t\in[0,T]} \left(t^\frac{\kappa}{s}\|(v\otimes v)(t)\|_{\dot C^{\kappa_0}}+t^\frac{m+\kappa}{s}\|(v\otimes v)(t)\|_{\dot C^{m+\kappa_0}}\right)\\
				&+T^{\frac{m}{s}+1}\|\Lambda^s\phi\|_{C^m}+T\|\Lambda^s\phi\|_{L^\infty}.
			\end{align*}
			Note that 
			\begin{align*}
				\sup_{t\in[0,T]} t^\frac{m+\kappa}{s}\|\partial_x^m(v\otimes v)(t)\|_{\dot C^{\kappa_0}}\lesssim\|\Lambda^{1-s}v\|_{T}^2\lesssim  (\|\Lambda^{1-s}v-\phi\|_{T}+\|\phi\|_T)^2.
			\end{align*}
			Hence we obtain
			\begin{equation}\label{nses}
				\begin{aligned}
					\|\Lambda^{1-s}u-\phi\|_T&\leq C_0 \left(\|\Lambda^{1-s}u_0-\phi\|_{\dot B^0_{\infty,\infty}}+(\|\Lambda^{1-s}v-\phi\|_{T}+\|\phi\|_T)^2\right)\\
					&+C_0\left(T^{\frac{m}{s}+1}\|\Lambda^s\phi\|_{C^{m}}+T\|\Lambda^s\phi\|_{L^\infty}\right).
				\end{aligned}
			\end{equation}
			If $\|u_0\|_{\dot B^{1-s}_{\infty,\infty}}\leq \varepsilon_0$, we take $\phi=0$, $\sigma =2C_0\varepsilon_0$. then   \eqref{nses} implies that for any $v\in \mathcal{X}^\sigma_{\infty,0}$,
			\begin{align*}
				\|\Lambda^{1-s}u\|_\infty&\leq C_0 \|u_0\|_{\dot B^{1-s}_{\infty,\infty}}+C_0\|\Lambda^{1-s}v\|_{T}^2\leq C_0\varepsilon_0+4C_0^3\varepsilon_0^2\leq \sigma,
			\end{align*}
			where we take $0<\varepsilon_0\leq \frac{1}{100(C_0+1)^{10}}$. Hence we obtain $u\in \mathcal{X}^\sigma_{\infty,0}$.\\
			We consider another case. Suppose \eqref{nscon} holds, then take $\phi= \Lambda^{1-s}  u_0\ast \rho_{\eps_1}$,  then we have
			$$\|\Lambda^{1-s}u_0-\phi\|_{\dot B^0_{\infty,\infty}}\leq \frac{5}{4}\varepsilon_0, \ \ \ \ \ \ \ \|\phi\|_T\leq C(T^\frac{\gamma}{s}\eps_1^{-\gamma}+T^\frac{m+\kappa}{s}\eps_1^{-(m+\kappa)})\|u_0\|_{\dot B^{1-s}_{\infty,\infty}}.$$ 
			Now we fix $T_1=\left(\frac{\varepsilon_0\eps_1^{m+\kappa}}{100+C_0+\| u_0\|_{\dot B^{1-s}_{\infty,\infty}}}\right)^\frac{100}{\gamma}$, then for any $v\in\mathcal{X}^\sigma_{T_1,\Lambda^{1-s}  u_0\ast \rho_{\eps_1}}$,
			we obtain from \eqref{nses} that 
			\begin{align*}
				\|\Lambda^{1-s}u-\phi\|_{T_1}&\leq C_0 (\|\Lambda^{1-s}u_0-\phi\|_{\dot B^0_{\infty,\infty}}+(\|\Lambda^{1-s}v-\phi\|_{T_1}+T_1^\frac{\gamma}{s}\eps_1^{-2}\| u_0\|_{\dot B^{1-s}_{\infty,\infty}})^2+T_1^\frac{\kappa}{s}\eps_1^{-2}\| u_0\|_{\dot B^{1-s}_{\infty,\infty}})\\
				&\leq \frac{5}{4}C_0\varepsilon_0+C_0(2C_0\varepsilon_0+T_1^\frac{\gamma}{s}\eps_1^{-2}\| u_0\|_{\dot B^{1-s}_{\infty,\infty}})^2+C_0T_1^\frac{\kappa}{s}\eps_1^{-2}\| u_0\|_{\dot B^{1-s}_{\infty,\infty}}\\
				&\leq \sigma.
			\end{align*}
			Hence we obtain that $v\in\mathcal{X}^\sigma_{T_1,\Lambda^{1-s}  u_0\ast \rho_{\eps_1}}$ implies that $\mathcal{S}v\in\mathcal{X}^\sigma_{T_1,\Lambda^{1-s}  u_0\ast \rho_{\eps_1}}$ under the assumption \eqref{nscon}.
			Then we do the contraction estimates. Let $v_1,v_2\in \mathcal{X}^\sigma_{\infty,0}$, $u_k=\mathcal{S}v_k, k=1,2$. One has 
			\begin{align*}
				\partial_t \Lambda^{1-s}(u_1-u_2)+\Lambda(u_1-u_2)=-\Lambda^{1-s}\mathbb{P}\operatorname{div}(v_1\otimes v_1-v_2\otimes v_2).
			\end{align*}
			Then apply Theorem \ref{thmcons} with $F=(v_1\otimes v_1-v_2\otimes v_2)$ and $G=0$ to get
			\begin{align}\label{contrans}
				\|\Lambda^{1-s}(u_1-u_2)\|_\infty\lesssim \sup_{t>0}\left(t^\frac{\kappa}{s}\|(v_1\otimes v_1-v_2\otimes v_2)(t)\|_{\dot C^{\kappa_0}}+t^\frac{m+\kappa}{s}\|(v_1\otimes v_1-v_2\otimes v_2)(t)\|_{\dot C^{m+\kappa_0}}\right).
			\end{align}
			By Lemma \ref{Hol} we have 
			\begin{align*}
				&\|\nabla_x^m(v_1\otimes v_1-v_2\otimes v_2)\|_{\dot C^{\kappa_0}}\\
				&\quad\quad\lesssim\sum_{m_1+m_2=m}\Big\{ \|\nabla_x^{m_1}(v_1-v_2)\|_{\dot C^{\kappa_0}}(\|\nabla_x^{m_2}v_1\|_{L^\infty}+\|\nabla_x^{m_2}v_2\|_{L^\infty})\\
				&\quad\quad\quad\quad +\|\nabla_x^{m_1}(v_1-v_2)\|_{L^\infty}(\|\nabla_x^{m_2}v_1\|_{\dot C^{\kappa_0}}+\|\nabla_x^{m_2}v_2\|_{\dot C^{\kappa_0}})\Big\}.
			\end{align*}
			Hence 
			$$
			\sup_{t>0}t^\frac{m+\kappa}{s}\|\nabla_x^m(v_1\otimes v_1-v_2\otimes v_2)(t)\|_{\dot C^{\kappa_0}}\lesssim \|\Lambda^{1-s}(v_1-v_2)\|_\infty(\|\Lambda^{1-s}v_1\|_\infty+\|\Lambda^{1-s}v_2\|_\infty).
			$$
			Substitute this into \eqref{contrans}, we obtain 
			\begin{align*}
				\|\Lambda^{1-s}(u_1-u_2)\|_\infty\leq C_1\|\Lambda^{1-s}(v_1-v_2)\|_\infty(\|\Lambda^{1-s}v_1\|_\infty+\|\Lambda^{1-s}v_2\|_\infty)\leq 2C_1\sigma \|\Lambda^{1-s}(v_1-v_2)\|_\infty.
			\end{align*}
			Similarly, for $v_1,v_2\in\mathcal{X}^\sigma_{T_1,\Lambda^{1-s}  u_0\ast \rho_{\eps_1}} $, we have 
			\begin{align*}
				\|\Lambda^{1-s}(u_1-u_2)\|_{T_1}&\leq C_2 \|\Lambda^{1-s}(v_1-v_2)\|_{T_1}(\|\Lambda^{1-s}(v_1-\phi)\|_{T_1}+\|\Lambda^{1-s}(v_2-\phi)\|_{T_1}+\|\phi\|_{T_1})\\
				&\leq C_2(2\sigma+T_1^\frac{\gamma}{s}\eps_1^{-2}\|u_0\|_{\dot B^{1-s}_{\infty,\infty}})\|\Lambda^{1-s}(v_1-v_2)\|_{T_1}\\
				&\leq 3C_2\sigma |\Lambda^{1-s}(v_1-v_2)\|_{T_1}.
			\end{align*}
			We take $\varepsilon_0=\frac{1}{100(C_0+C_1+C_2+1)^{10}}$, which implies that $2C_1\sigma+3C_2\sigma\leq \frac{1}{2}$. 
			Hence we obtain that 
			\begin{align*}
				&\|\Lambda^{1-s}(u_1-u_2)\|_\infty\leq \frac{1}{2}\|\Lambda^{1-s}(v_1-v_2)\|_\infty,\ \  \forall v_1,v_2\in \mathcal{X}^\sigma_{\infty,0},\\
				&	\|\Lambda^{1-s}(u_1-u_2)\|_{T_1}\leq\frac{1}{2}\|\Lambda^{1-s}(v_1-v_2)\|_{T_1}, \ \  \forall v_1,v_2\in\mathcal{X}^\sigma_{T_1,\Lambda^{1-s}  u_0\ast \rho_{\eps_1}}.
			\end{align*}
			Hence 
			\begin{itemize}
				\item If $\| u_0\|_{\dot B^{1-s}_{\infty,\infty}}\leq \varepsilon_0$, then there exists a unique $u\in\mathcal{X}^\sigma_{\infty,0} $ such that $\mathcal{S}u=u$.\\
				\item If \eqref{nscon} holds, then there exists a unique $u\in\mathcal{X}^\sigma_{T_1,\Lambda^{1-s}  u_0\ast \rho_{\eps_1}}$ such that $\mathcal{S}u=u$.
			\end{itemize}
			This completes the proof.
			
		\end{proof}
		\begin{remark}
			\label{remNS} When $s=2$, the nonlinear term is $-\Lambda^{-1}\mathbb{P}\operatorname{div}(v\otimes v)$. Note that the operator $\Lambda^{-1}\mathbb{P}\operatorname{div}$ is a $0$-order operator, hence the nonlinear term loses the ``derivative form". In this case, we can not get a closed estimate by Theorem \ref{thmcons}. Precisely, $\|(v\otimes v)(t)\|_{L^\infty}\sim t^{-1} $, which is not integrable in $t$.
		\end{remark}
		\begin{proof}[Proof of Proposition \ref{propsqg}]
			Observe that $\nabla \cdot \nabla ^\perp \theta=0$, we can rewrite the equation as 
			\begin{align*}
				\partial_t \theta+\Lambda^s \theta+\nabla \cdot ( \theta\nabla^\perp \Lambda^{-1}\theta )=0.
			\end{align*}
			The proof is parallel to Proposition \ref{propns}, we omit details here.
		\end{proof}
		~~~\vspace{0.3cm}\\
		In the following section \ref{secmeanc}- section \ref{secdas}, we consider some geometric flows. By geometric flow we mean the evolution of a geometric structure on a manifold under a differential equation, usually associated with curvature. The aim of these problems is to find a family $\Gamma=\{ \Gamma(t):t\geq 0\}$ of immersed orientable hypersurfaces in $\mathbb{R}^d$ which evolves according to a law
		\begin{align*}
			V=f(\kappa_1,\kappa_2),\ \ \ \Gamma|_{t=0}=\Gamma_0.
		\end{align*}
		Here $V$ denotes the normal velocity of $\Gamma$, $f$ is known relation, and $\kappa_1,\kappa_2$ are principal curvatures which describe the local geometry of the surface $\Gamma(t)$. The classical notions of curvature of a surface include 
		\begin{align*}
			\text{mean curvature:} \ H=\frac{\kappa_1+\kappa_2}{2},\ \ \ \quad\quad\quad \text{Gauss curvature:}\  K=\kappa_1\kappa_2.
		\end{align*}
		One obtains the mean curvature flow by setting $f=H$, the surface diffusion flow by setting $f=\Delta_{\Gamma}H$, and Willmore flow by setting $f=\Delta_{\Gamma}H-2H(H^2-K)$, where $\Delta_\Gamma$ is the Laplace--Beltrami operator of $\Gamma(t)$. We provide a succinct introduction to each of these equations in the following sections. For a more in-depth exploration, readers may refer to \cite{bookgeo,DDE2005,ES2001} and references therein.
		
		\subsection{Mean Curvature Flow}\label{secmeanc}
		Assume $\mathcal{M}$ is a $d$-dimensional orientable manifold, $F_0:\mathcal{M}\rightarrow \mathbb{R}^{d+N}$ is an immersion, then we consider the immersions $F:(0,T)\times \mathcal{M}\rightarrow \mathbb{R}^{d+N}$ to be the solutions to the following system
		\begin{equation}\label{mcfe}
			\begin{aligned}
				&\partial_tF=\mathbf{H},\ \operatorname{on}\ (0,T)\times \mathcal{M}, \\
				&F|_{t=0}=F_0,
			\end{aligned}
		\end{equation}
		where $\mathbf{H}$ is the mean curvature vector of $\mathcal{M}_t=F(t,\mathcal{M})$ at $F(t,x)$. We are mostly interested in the case when  $F_0(x)=(x,f_0(x))$ is the graph of a function $f_0:\mathbb{R}^d\rightarrow\mathbb{R}^N$. In this paper, we will solve the following system, where $f:(0,T)\times\mathbb{R}^d\rightarrow\mathbb{R}^m$ satisfies
		\begin{equation}\label{mcf}
			\begin{aligned}
				&\partial_tf=\A[\nabla f]:\nabla^2f,\ \operatorname{in} \ (0,T)\times\mathbb{R}^d,\\
				&f|_{t=0}=f_0,
			\end{aligned}
		\end{equation}
		with
		\begin{align*}
			\A[\nabla f]=\left[\mathrm{Id}+\sum_{i=1}^{N}\nabla f^i\otimes \nabla f^i\right]^{-1}\ \ \ \ \text{and}\ \ \A[\nabla f]:\nabla ^2=\sum_{i,j=1}^d(\A[\nabla f])_{ij}\partial_i\partial_j.
		\end{align*}
		The relation between \eqref{mcfe} and \eqref{mcf} can be seen in \cite{W2004a}. For any solution of \eqref{mcfe}, there exists a family of diffeomorphisms $r:(0,\infty)\times\mathbb{R}^d\rightarrow\mathbb{R}^d$ such that $\tilde{F}(t,x)=F(t,r(t,x))$ can be written as $\tilde{F}(t,x)=(x,f(t,x))$, where $f(t,x)$ is a solution of \eqref{mcf}. Conversely, for any solution of \eqref{mcf}, $\tilde{F}(t,x)=(x,f(t,x))$ is a solution to \eqref{mcfe}. 
		
		We make a short review about the works before. In \cite{EH1989}, Ecker and Huisken proved that when $N=1$, system \eqref{mcfe} admits a global solution for any initial data which is locally Lipschitz. They also proved that the solution is smooth for any positive time. The case when $N>1$ is more subtle. Lawson and Osserman \cite{LO1977} constructed a solution $f:\mathbb{R}^4\rightarrow\mathbb{R}^3$,  which is Lipschitz but not $C^1$, to the minimal surface system (which is a stationary phase of the mean curvature flow). Hence it is necessary to impose some smallness conditions on the initial data to obtain global results in high codimensions. In \cite{W2004b}, Wang showed the existence of a local smooth solution in arbitrary codimension, under certain smallness assumptions. Other results for map between constant curvature manifolds can be found in \cite{W2002}. Also in \cite{CCW2009}, the authors constructed a global solution when $d=N$ and the initial data $f_0$ has the form that $f_0=\nabla u_0$ for some $u_0$. We also refer readers to \cite{mcsing} and references therein for results about singularities in mean curvature flow.\\
		The result in this paper is somehow inspired by \cite{Koch2011}, where for initial data $f_0:\mathbb{R}^d\rightarrow\mathbb{R}^N$ satisfying 
		\begin{equation*}
			\|\nabla f_0\|_{L^\infty(\mathbb{R}^d,\mathbb{R}^N)}<\delta
		\end{equation*} 
		for some $\delta$ small enough, Koch and Lamm \cite{Koch2011} proved that there exists a unique global solution $f$ satisfying
		\begin{equation*}
			\sup_{t} \left(\|\nabla f(t)\|_{L^\infty}+\sup_{x} \sup_{R} R^{\frac{2}{N+4}}\|\nabla^2 f(t)\|_{L^{N+4}(B_R(x)\times(\frac{R^2}{2},R^2))}\right)<C\|\nabla f_0\|_{L^\infty(\mathbb{R}^d,\mathbb{R}^N)}.
		\end{equation*}
		We will prove the existence of local solution with large initial data, and global solution with small initial data. Precisely, we will prove
		\begin{theorem}\label{thmmc}
			Fix $\kappa$ such that $0<2-\kappa\ll1$, and $m\in\mathbb{N}$. There exists $\varepsilon_0>0$ such that 	~\\
			i)	For any  initial data $f_0$ with $\| \nabla f_0\|_{L^\infty}\leq \varepsilon_0$, there exists a unique solution $f$ to \eqref{mcf} satisfying
			\begin{align*}
				\sup_{t>0}(\|\nabla f(t)\|_{L^\infty}+t^\frac{m+\kappa}{2}\|\nabla f(t)\|_{\dot C^{m+\kappa}})\leq C \| \nabla f_0\|_{L^\infty}.
			\end{align*}
			ii) For any initial data $f_0\in \dot W^{1,\infty}$, if there exists $\eps_1>0$ such that
			\begin{align}\label{mccon}
				\| \nabla f_0-\nabla f_0\ast \rho_{\eps_1}\|_{L^\infty}\leq \varepsilon_0,
			\end{align}
			then   there exist $T=T(\eps_0,\eps_1,\|f_0\|_{\dot W^{1,\infty}})>0$ and a unique solution  $f$ to \eqref{mcf} in $[0,T]$ satisfying
			\begin{align*}
				\sup_{t\in[0,T]}(\|\nabla f(t)\|_{L^\infty}+t^\frac{m+\kappa}{2}\|\nabla f(t)\|_{\dot C^{m+\kappa}})\leq C \| \nabla f_0\|_{L^\infty}.
			\end{align*}
		\end{theorem}
		\begin{remark}
			Compared with the result in \cite{Koch2011}, the result in Theorem \ref{thmmc} includes estimates of higher order derivatives, and  also indicates that the global solution decays in time.  
		\end{remark}
		\begin{proof} Define
			\begin{align*}
				\|h\|_{T}=\sup_{t\in[0,T]}\left(\|\nabla h(t)\|_{L^\infty}+t^{\frac{m+\kappa}{2}}\|\nabla h(t)\|_{\dot C^{m+\kappa}}\right).
			\end{align*}
			Let $\phi\in C^\infty(\mathbb{\mathbb{R}}^d)$, $\sigma,T>0$ that will be fixed later, we consider the set
			$$
			\mathcal{X}_{T,\phi}^\sigma=\left\{g\in L_T^\infty\dot W^{1,\infty}: g|_{t=0}=f_0, \ \|g-\phi\|_{T}\leq \sigma\right\},
			$$ 
			where $g|_{t=0}=g_0$ in the weak sense. For any $g\in	\mathcal{X}_{T,\phi}^\sigma$, we define a map $\mathcal{S}g=f$, where $f$ is the solution to the Cauchy problem
			\begin{align}\label{mceq}
				\partial_t f-\A[\nabla \phi]:\nabla^2f=(\A[\nabla g]-\A[\nabla \phi]):\nabla ^2 g, \ \ \ \ f|_{t=0}=f_0.
			\end{align}
			In the following, we prove that $\mathcal{S}$ is a contraction map from $	\mathcal{X}_{T,\phi}^\sigma$ to itself with suitable choice of $\phi,T,\sigma$. Hence there exists a fixed point in $	\mathcal{X}_{T,\phi}^\sigma$, which is a solution to  equation \eqref{mcf}.
			We first prove that
			\begin{itemize}
				\item i. If $\|\nabla f_0\|_{L^\infty}\leq \varepsilon_0$, then there exists $\sigma_1>0$ such that $\mathcal{S}$ maps $\mathcal{X}_{\infty,0}^{\sigma_1}$ to itself.
				\item ii. If \eqref{mccon} holds, then there exist $l,\sigma_2>0$ and $T>0$ such that $\mathcal{S}$ maps $\mathcal{X}_{T,f_0\ast \rho_{\eps_1}}^{\sigma_2}$ to itself.
			\end{itemize} 
			Denote  $\tilde f(t,x)=f(t,x)-\phi(x)$, $\tilde g(t,x)=g(t,x)-\phi(x)$. We have 
			\begin{align*}
				\partial_t \partial_i\tilde f-\A[\nabla\phi]: \nabla ^2 \partial_i\tilde f=\partial_i [(\A[\nabla g]-\A[\nabla \phi] ):\nabla ^2g]+\big(\partial_i (\A[\nabla\phi]): \nabla^2  f+\A[\nabla \phi]:\nabla ^2\partial_i\phi\big):=\partial_i F_1+F_2.
			\end{align*}	
			Applying Theorem \ref{lemmain} with $F=F_1$ and $G=F_2$, we obtain that there exists $T=T(\|\phi\|_{C^{m+3}})>0$ such that 
			\begin{equation}\label{mcmain}
				\begin{aligned}
					\|\tilde f\|_T\lesssim &\|\nabla ( f_0-\phi)\|_{L^\infty}+ \sup_{t\in[0,T]}\left(t^\frac{\kappa}{2}\|F_1(t)\|_{\dot C^{\kappa-1}}+t^\frac{m+\kappa}{2}\|F_1(t)\|_{\dot C^{m+\kappa-1}}\right)\\
					&+\sup_{t\in[0,T]}t^\frac{m}{2}\|F_2(t)\|_{\dot C^m}+\|F_2\|_{L^1_TL^\infty}.
				\end{aligned}
			\end{equation}
			We first compute terms associated with $F_1$. By Lemma \ref{Hol}, we have 
			\begin{equation}\label{fhol}
				\begin{aligned}
					&\|F_1\|_{\dot C^{\kappa-1}}\lesssim \|\A[\nabla g]-\A[\nabla \phi]\|_{\dot C^{\kappa-1}}\|\nabla ^{2} g\|_{L^\infty}+\|\A[\nabla g]-\A[\nabla \phi]\|_{L^\infty}\|\nabla ^{2} g\|_{\dot C^{\kappa-1}},\\
					&\|F_1\|_{\dot C^{m+\kappa-1}}\lesssim \sum_{j=0}^m\left(\|\nabla ^j(\A[\nabla g]-\A[\nabla \phi])\|_{\dot C^{\kappa-1}}\|\nabla ^{2+m-j} g\|_{L^\infty}+\|\nabla ^j(\A[\nabla g]-\A[\nabla \phi])\|_{L^\infty}\|\nabla ^{2+m-j} g\|_{\dot C^{\kappa-1}}\right).
				\end{aligned}
			\end{equation}
			We obtain from Lemma \ref{lemcom} that 
			\begin{equation}\label{aahol}
				\begin{aligned}
					&\|(\A[\nabla g]-\A[\nabla \phi])(t)\|_{L^\infty}\lesssim \|\nabla\tilde g(t)\|_{L^\infty}\lesssim \|\tilde g\|_T,\\
					&\|(\A[\nabla g]-\A[\nabla \phi])(t)\|_{\dot C^{\kappa-1}}\lesssim \|\nabla \tilde g(t)\|_{\dot C^{\kappa-1}}+\|\nabla\tilde g(t)\|_{L^\infty}(\|\nabla g(t)\|_{\dot C^{\kappa-1}}+\|\nabla \phi\|_{\dot C^{\kappa-1}})\\
					&\quad\quad\quad\quad\quad\quad\quad\quad\quad\quad\quad\ \lesssim t^{-\frac{\kappa-1}{2}}\|\tilde g\|_T(1+\|\tilde g\|_T+T^\frac{1}{10}\|\nabla \phi\|_{\dot C^{\kappa-1}}).
				\end{aligned}
			\end{equation}
			Moreover, by Lemma \ref{lemcom}, one has
			\begin{equation}\label{aahol1}
				\begin{aligned}
					&\|\nabla ^m(\A[\nabla g]-\A[\nabla \phi])(t)\|_{L^\infty}\\
					&\quad\quad\lesssim \|\nabla\tilde g(t)\|_{\dot C^{m}}+\sum_{m_1=1}^m \|\nabla\tilde g(t)\|_{\dot C^{m-m_1}}\left(\|\nabla g(t)\|_{\dot C^{m_1}}+\|\nabla \phi\|_{\dot C^{m_1}}+(\|\nabla g(t)\|_{\dot C^1}+\|\nabla \phi\|_{\dot C^1})^{m_1}\right)\\
					&\quad\quad\lesssim t^{-\frac{m}{2}}\|\nabla\tilde g\|_T(1+\|\nabla\tilde g\|_T+T^\frac{1}{10}\|\nabla \phi\|_{C^m})^{m+10},\\
					&\|\nabla ^m(\A[\nabla g]-\A[\nabla \phi])(t)\|_{\dot C^{\kappa-1}}\\
					&\quad\quad \lesssim \|\nabla\tilde g(t)\|_{\dot C^{m+\kappa-1}}+\sum_{m_1=1}^m \Big\{\|\nabla\tilde g(t)\|_{\dot C^{m-m_1+\kappa-1}}\left(\|\nabla g(t)\|_{\dot C^{m_1}}+\|\nabla \phi\|_{\dot C^{m_1}}+(\|\nabla g(t)\|_{\dot C^1}+\|\nabla \phi\|_{\dot C^1})^{m_1}\right)\\
					&\quad\quad\quad +\|\nabla\tilde g(t)\|_{\dot C^{m-m_1}}\left(\|\nabla g(t)\|_{\dot C^{m_1+\kappa-1}}+\|\nabla \phi\|_{\dot C^{m_1+\kappa-1}}+(\|\nabla g(t)\|_{\dot C^{\kappa-1}}+\|\nabla \phi\|_{\dot C^{\kappa-1}})^\frac{m_1+\kappa-1}{\kappa-1}\right)\Big\}\\
					&\quad\quad \lesssim t^{-\frac{m+\kappa-1}{2}}\|\tilde g\|_T(1+\|\tilde g\|_T+T^\frac{1}{10}\|\nabla \phi\|_{C^m})^{m+10}.
				\end{aligned}
			\end{equation}
			Then we collect \eqref{fhol}, \eqref{aahol} and \eqref{aahol1} to obtain
			\begin{equation}\label{F1}
				\begin{aligned}
					&\sup_{t\in[0,T]}(t^\frac{\kappa}{2}\|F_1(t)\|_{\dot C^{\kappa-1}}+t^\frac{m+\kappa }{2}\|F_1(t)\|_{\dot C^{m+\kappa -1}})\\
					&\quad\quad\quad\quad\quad\quad\quad\quad\quad\lesssim \|\tilde g\|_T(\|\tilde g\|_{T}+{T}^\frac{1}{10}\|\nabla \phi\|_{C^{m+2}})(1+\|\tilde g\|_T+T^\frac{1}{10}\|\nabla \phi\|_{C^{m+2}})^{m+10}.
				\end{aligned}
			\end{equation}
			Now we estimate $F_2$, 
			\begin{align*}
				&\|\partial_i (\A[\nabla\phi]) \nabla^2  f\|_{L^\infty}\lesssim \|\nabla \phi\|_{\dot C^1}\|\nabla f\|_{\dot C^1},\ \ \ 
				\|\A[\nabla \phi]\nabla ^2\partial_i\phi\|_{L^\infty}\lesssim \|\nabla \phi\|_{\dot C^2},\\
				&	\|\partial_i (\A[\nabla\phi]) \nabla^2  f\|_{\dot C^m}\lesssim \sum_{l=0}^m (\|\nabla \phi\|_{\dot C^{l+1}}+\|\nabla \phi\|_{\dot C^1}^{m+1})\|\nabla ^2f \|_{\dot C^{m-l}},\\
				&	\|\A[\nabla \phi]\nabla ^2\partial_i\phi\|_{ C^m}\lesssim \|\nabla \phi\|_{C^{m+3}}(1+\|\nabla \phi\|_{C^{m+3}})^{m+3}.
			\end{align*}
			Hence 
			\begin{align}\label{F2}
				{\|F_2\|_{L^1_TL^\infty}} +\sup_{t\in[0,T]}t^\frac{m}{2}\|F_2(t)\|_{\dot C^m}\lesssim T^\frac{1}{10}(1+\|\tilde f\|_T+\|\nabla \phi\|_{C^{m+3}})^{m+10}\|\nabla \phi\|_{C^{m+3}}.
			\end{align}
			Combining \eqref{mcmain}-\eqref{F2}, we obtain 
			\begin{equation} \label{mcm}\begin{aligned}
					\|\tilde f\|_T\lesssim &\|\nabla (f_0-\phi)\|_{L^\infty}+\|\tilde g\|_T(\|\tilde g\|_{T}+T^\frac{1}{10}\|\nabla \phi\|_{C^{m+3}})(\|\tilde g\|_T+T^\frac{1}{10}\|\nabla \phi\|_{ C^{m+3}}+1)^{m+10}\\
					&\quad\quad\quad+T^\frac{1}{10}(\|\tilde f\|_T+1+\|\nabla \phi\|_{C^{m+3}})^{m+10}\|\nabla \phi\|_{C^{m+3}}.
				\end{aligned}
			\end{equation}
			Suppose $\|\nabla f_0\|_{L^\infty}\leq \varepsilon_0$ for some $\varepsilon_0\ll 1 $ that will be fixed later. Take $\phi\equiv 0$ in \eqref{mcm}, then $A[\nabla \phi]$ is constant. By Theorem \ref{thmcons}, we know that \eqref{mcm} holds for any $T>0$, hence
			\begin{align*}
				\| f\|_\infty\leq C_0\|\nabla  f_0\|_{L^\infty}+C_0\| f\|_{\infty}\|  g\|_\infty(\|g\|_\infty+1)^{m+10}.
			\end{align*}
			Take $\varepsilon_0=\frac{1}{100^m(2C_0+1)^{100}}$ and $\sigma_1=2C_0 \varepsilon_0$. Then for any $g\in \mathcal{X}_{\infty,0}^{\sigma_1}$, we have $\|g\|_{\infty}\leq 2C_0 \varepsilon_0$. Then 
			\begin{align*}
				\|  f\|_\infty\leq C_0\|\nabla  f_0\|_{L^\infty}+C_0\| f\|_{\infty}\|  g\|_\infty(\| g\|_\infty+1)^{m+10}\leq C_0\varepsilon_0+2C_0^2\varepsilon_0(2C_0\varepsilon_0+1)	^{m+10}\|  f\|_\infty.
			\end{align*}
			Hence we get 
			\begin{align*}
				\| f\|_\infty\leq \frac{C_0\varepsilon_0}{1-2C_0^2\varepsilon_0(2C_0\varepsilon_0+1)^{m+10}}\leq \sigma_1.
			\end{align*}
			We obtain that $f\in \mathcal{X}_{\infty,0}^{\sigma_1}$.\\
			If condition \eqref{stcon} holds, we take $\phi=f_0\ast \rho_{\eps_1}$ in \eqref{mcm}, then there exists $T_0=T_0(\| f_0\|_{W^{1,\infty}},\eps_0,\eps_1)>0$ such that 
			\begin{align*}
				&\|\tilde f\|_{T_0}\leq C_0\left(\|\nabla(f_0-\phi) \|_{L^\infty}+{T_0}^\frac{1}{100}(\|\tilde f\|_{T_0}+1)\eps_1^{-(m+10)}\|\nabla f_0\|_{L^\infty}\right.\\
				&\left.\quad\quad\quad\quad\quad\quad+\|\tilde g\|_{T_0}(\|\tilde g\|_{T_0}+{T_0}^\frac{1}{100}\eps_1^{-(m+10)}\|\nabla f_0\|_{L^\infty})(\|\tilde g\|_{T_0}+{T_0}^\frac{1}{100}\eps_1^{-(m+10)}\|\nabla f_0\|_{ L^\infty}+1)^{m+10}\right).
			\end{align*}
			Let \begin{align}\label{cot1si}
				T_1=\min\left\{T_0,\left(\frac{\varepsilon_0\eps_1}{1+C_0+\|\nabla f_0\|_{L^\infty}}\right)^{10^4m}\right\}, \quad\quad\sigma_2=2C_0\varepsilon_0.
			\end{align} 
			If $g\in  \mathcal{X}_{T_1,f_0\ast \rho_{\eps_1}}^{\sigma_2}$, then 
			\begin{align*}
				\|\tilde f\|_{T_1}\leq &\frac{5}{4}C_0\varepsilon_0+2C_0\varepsilon_0(\|\tilde f\|_{T_1}+{T_1}^\frac{1}{100}\eps_1^{-(m+10)}\|\nabla f_0\|_{L^\infty})(2C_0\varepsilon_0+{T_1}^\frac{1}{100}\eps_1^{-(m+10)}\|\nabla f_0\|_{ L^\infty}+1)^{m+10}\\
				&+{T_1}^\frac{1}{100}(\|\tilde f\|_{T_1}+1)\eps_1^{-(m+10)}\|\nabla f_0\|_{L^\infty}\\
				\leq &\frac{5}{4}C_0\varepsilon_0+2C_0\varepsilon_0^2(3C_0\varepsilon_0+1)+2C_0\varepsilon_0(3C_0\varepsilon_0+1)	\|\tilde f\|_{T_1}.
			\end{align*}
			Then it is easy to check that 
			\begin{align*}
				\|\tilde f\|_{T_1}\leq \sigma_2.
			\end{align*}
			This implies that 
			$f\in  \mathcal{X}_{T_1,f_0\ast \rho_{\eps_1}}^{\sigma_2}$.\vspace{0.3cm}\\
			In the following we do the contraction estimates. Consider $g_1,g_2\in \mathcal{X}^\sigma_{T,\phi}$ for some $\sigma,T>0$ and $\phi\in C^{\infty}(\mathbb{R}^d)$. Denote  $\mathbf{g}=g_1-g_2$, $\mathbf{f}=f_1-f_2=\mathcal{S}g_1-\mathcal{S}g_2$, and $(\tilde g_1,\tilde g_2)=(g_1-\phi,g_2-\phi)$. We have
			\begin{align*}
				&\partial_t  \partial_i\mathbf{f}-\A[\nabla \phi]:\nabla ^2 \partial_i\mathbf{f}\\
				&=\partial_i \left( (\A[\nabla g_1]-\A[\nabla g_2]):\nabla^2 g_1+(\A[\nabla g_2]-\A[\nabla \phi]):\nabla^2 (g_1-g_2)\right)+\partial_i\A[\nabla \phi]:\nabla ^2(f_1-f_2)\\
				:&=\partial_i\tilde F_1+\tilde F_2.
			\end{align*} 
			By Theorem \ref{lemmain} with $F=\tilde F_1$ and $G=\tilde F_2$, there exists $T>0$ such that 
			\begin{align*}
				\|\mathbf{f}\|_T\lesssim&\sup_{t\in[0,T]}\left(t^\frac{\kappa}{2}\|\tilde F_1(t)\|_{\dot C^{\kappa-1}}+t^\frac{m+\kappa}{2}\|\tilde F_1(t)\|_{\dot C^{m+\kappa-1}}+t^\frac{m}{2}\|\tilde F_2(t)\|_{\dot C^m}\right)+{\|\tilde F_2\|_{L^1_TL^\infty}}.
			\end{align*}
			Following the estimates for $F_1$ in \eqref{F1} and $F_2$ in \eqref{F2}, we obtain that 
			\begin{equation}\label{hahaha}
				\begin{aligned}
					\|\mathbf{f}\|_T\lesssim&\|\mathbf{g}\|_T(\|(\tilde g_1,\tilde g_2)\|_T+{T}^\frac{1}{10}\|\nabla \phi\|_{C^{m+10}})(1+\|(\tilde g_1,\tilde g_2)\|_T+{T}^\frac{1}{10}\|\nabla \phi\|_{C^{m+10}})^{m+10}\\
					&+{T}^\frac{1}{10}\|\nabla \phi\|_{C^{m+10}}(1+\|\nabla \phi\|_{C^{m+10}})^{m+10}	\|\mathbf{f}\|_T.
				\end{aligned}
			\end{equation}
			If $g_1,g_2\in \mathcal{X}^{\sigma_1}_{\infty,0}$, then $\|(\tilde g_1,\tilde g_2)\|_\infty\leq 2\sigma_1$. Take $\phi\equiv0$ in \eqref{hahaha}, one has  
			\begin{align*}
				\|\mathbf{f}\|_\infty\leq & C_0\|\mathbf{g}\|_\infty \|(\tilde g_1,\tilde g_2)\|_\infty(1+\|(\tilde g_1,\tilde g_2)\|_\infty)^{m+10}
				\leq  2C_0\sigma_1(1+2\sigma_1)^{m+10}\|\mathbf{g}\|_\infty\leq \frac{1}{2}\|\mathbf{g}\|_\infty.
			\end{align*}
			Hence $\mathcal{S}:\mathcal{X}^{\sigma_1}_{\infty,0}\to \mathcal{X}^{\sigma_1}_{\infty,0}$ is a contraction map.\\
			\vspace{0.3cm}\\
			If $g_1,g_2\in  \mathcal{X}_{T_1,f_0\ast \rho_{\eps_1}}^{\sigma_2}$, one has  $\|(\tilde g_1,\tilde g_2)\|_T\leq 2 \sigma_2$. Take $\phi=f_0\ast \rho_{\eps_1}$ in \eqref{hahaha},  one has 
			\begin{align*}
				\|\mathbf{f}\|_{T_1}\leq & C_0(\|\mathbf{g}\|_{T_1}+\|\mathbf{f}\|_{T_1}) (4\sigma_2+2{T_1}^\frac{1}{100}\eps_1^{-(m+10)}\|\nabla f_0\|_{L^\infty})(1+4\sigma_2+2{T_1}^\frac{1}{100}\eps_1^{-(m+10)}\|\nabla f_0\|_{L^\infty})^{m+10}\\
				:= &\Theta(\sigma_2, T,\|\nabla f_0\|_{L^\infty})(\|\mathbf{g}\|_{T_1}+\|\mathbf{f}\|_{T_1}).
			\end{align*}
			We choose $\sigma_2$ and $T_1$ as in \eqref{cot1si} to ensure that $\Theta(\sigma_2, T,\|\nabla f_0\|_{L^\infty})\leq\frac{1}{100}$, which implies that 
			\begin{align*}
				\|\mathbf{f}\|_{T_1}\leq \frac{1}{2}\|\mathbf{g}\|_{T_1}.
			\end{align*}
			Hence $\mathcal{S}:\mathcal{X}_{T_1,f_0\ast \rho_{\eps_1}}^{\sigma_2}\to \mathcal{X}_{T_1,f_0\ast \rho_{\eps_1}}^{\sigma_2}$ is a contraction map. This completes the proof.
		\end{proof}

		\subsection{Non-local mean curvature flow}\label{secnmc}
		We consider the non-local counterpart of the classical mean curvature flow, which is  a geometric flow driven by the fractional mean curvature.\\
		Let $E_0\subset \mathbb{R}^d$, $\mathcal{M}_0=\partial E_0$. For any fixed $a\in(0,1)$, we consider the family of immersions $F:[0,T)\times \mathcal{M}_0\to \mathbb{R}^d$ which satisfies
		\begin{align*}
			\begin{cases}
				\partial_t F(t,z)=-H_a(t,z)\nu(t,z),\ \ \ z\in \mathcal{M}_0,\ t>0,\\
				F(0,z)=z,\ \ \ z\in \mathcal{M}_0.
			\end{cases}
		\end{align*}
		Here $\nu(t,z)$ denotes the  outward  unit  normal vector to  the  hypersurface $\mathcal{M}_t=F(t,\mathcal{M}_0)$, and $H_a(t,z)$ is the fractional mean curvature of order $a$, which is given by 
		\begin{align*}
			H_a(t,z):=c_a \mathrm{P.V}.\int_{\mathbb{R}^d}\frac{\mathbf{1}_{\mathbb{R}^d\backslash E_t}(y)-\mathbf{1}_{E_t}(y)}{|z-y|^{d+a}}dy.
		\end{align*}
		The study of the non-local mean curvature flow has attracted much interest in recent years (see the survey paper \cite{Cinti} and references therein). We refer readers to \cite{Imbert,Ca19,CN21,SV19} for existence, uniqueness, and regularity of viscosity solutions. In \cite{JM20}, the authors prove short time existence and uniqueness of classical
		solution to the fractional mean curvature flow with a bounded $C
		^{1,1}$
		initial set. More recently, the authors in \cite{AFW} established short time existence, uniqueness and optimal H\"{o}lder regularity of classical solutions when the sets $E_t$ are given by subgraphs of functions $u(t,\cdot)\in C_{loc}^{1+\beta}(\mathbb{R}^{d-1})$, where $\beta>a$.

		In the present paper, we consider the case when the sets  $\{E_t\}_{t\in[0,T)}$  are given by subgraphs of functions $u(t,\cdot)\in \dot W^{1,\infty}(\mathbb{R}^{d-1})$.  In this setting, we have 
		\begin{align*}
			H[u](t,x):=H_a(t,x)=\mathrm{P.V.}\int_{\mathbb{R}^{d-1}}\frac{\mathcal{G}(\Delta_\alpha u(t,x))}{|\alpha|^{d-1+a}}d\alpha,
		\end{align*}
		with 
		\begin{align*}
			\Delta_\alpha u(t,x):=\frac{u(t,x)-u(t,x-\alpha)}{|\alpha|}, \ \ \ \ \ \mathcal{G}(\rho):=\int_{-\rho}^{\rho}\frac{1}{\langle \tau\rangle^{d+a}}d\tau. \end{align*}
		The evolution of $u$ is given by the flow associated to the
		quasi-linear evolution equation (see \cite{AFW} for details)
		\begin{equation}\label{eqnmc}
			\begin{aligned}
				&\partial_t u(t,x)=-\langle\nabla u(t,x)\rangle H[u](t,x),\ \ \ \ \ (t,x)\in (0,T)\times \mathbb{R}^{d-1},\\ 
				&u|_{t=0}=u_0.
			\end{aligned}
		\end{equation}
		\begin{proof}[Proof of Theorem \ref{thmnmc}]
			By changes of variables and the fact that $\mathcal{G}(\rho)$ is odd, we have 
			\begin{align*}
				H[u](x)&=\int_{\mathbb{R}^{d-1}}\frac{\mathcal{G}(\Delta_\alpha u(x))}{| \alpha |^{d-1+a}}d \alpha =\frac{1}{2}\int_{\mathbb{R}^{d-1}}\frac{\mathcal{G}(\Delta_ \alpha u(x))+\mathcal{G}(\Delta_{- \alpha }u(x))}{| \alpha |^{d-1+a}}d\alpha\\
				&=\frac{1}{2}\int_{\mathbb{R}^{d-1}}\frac{\mathcal{G}(\Delta_ \alpha u(x))-\mathcal{G}(-\Delta_{- \alpha }u(x))}{| \alpha |^{d-1+a}}d \alpha .
			\end{align*}
			By fundamental calculus we have $\mathcal{G}'(\rho)={2}{\langle\rho\rangle^{-{(d+a)}}}$, and 
			\begin{align*}
				\mathcal{G}(\Delta_\alpha u(x))-\mathcal{G}(-\Delta_{-\alpha}u(x))=&\mathcal{O}_\alpha u(x)\int_0^1\mathcal{G}'\left(\tau \Delta_\alpha u(x)-(1-\tau)\Delta_{-\alpha}u(x)\right)d\tau\\
				=&2\mathcal{O}_\alpha u(x)\int_0^1\langle\tau \Delta_\alpha u(x)-(1-\tau)\Delta_{-\alpha}u(x)\rangle^{-{(d+a)}}d\tau.
			\end{align*}
			Here we denote $\mathcal{O}_ \alpha  u(x)=\Delta_ \alpha u(x)+\Delta_{- \alpha }u(x)$ for simplicity.  
			Hence we obtain 
			\begin{align*}
				H[u](x)&=\int_{\mathbb{R}^{d-1}}\frac{\mathcal{O}_\alpha u(x)}{\langle\hat \alpha\cdot \nabla u(x)\rangle^{d+a}}\frac{d\alpha}{|\alpha|^{d-1+a}}+\int_{\mathbb{R}^{d-1}}\mathcal{O}_\alpha u(x)P[u](x,\alpha)\frac{d\alpha}{|\alpha|^{d-1+a}},
			\end{align*}
			where 
			\begin{align}\label{defP}
				P[u](x,\alpha)=\int_0^1\langle\tau \Delta_\alpha u(x)-(1-\tau)\Delta_{-\alpha}u(x)\rangle^{-{(d+a)}}d\tau-\langle\hat \alpha\cdot \nabla u(x)\rangle^{-(d+a)},\quad\hat{\alpha}=\frac{\alpha}{|\alpha|}.	
			\end{align}
			We rewrite the equation as 
			\begin{align*}
				\partial_t u(x) +&\langle\nabla u(x)\rangle\int_{\mathbb{R}^{d-1}}\frac{\mathcal{O}_\alpha u(x)}{\langle\hat \alpha\cdot \nabla u(x)\rangle^{d+a}}\frac{d\alpha}{|\alpha|^{d-1+a}}\\
				&\quad\quad=-\langle\nabla u(x)\rangle\int_{\mathbb{R}^{d-1}}\left(\mathcal{O}_\alpha u(x)\right)P[u](x,\alpha)\frac{d\alpha}{|\alpha|^{d-1+a}}.
			\end{align*}
			Let  $\phi:{\mathbb{R}^{d-1}}\to\mathbb{R}$ be a smooth function that will be fixed later. Denote 
			\begin{align*} 
				\mathcal{L}_\phi u(x)= \langle\nabla \phi(x)\rangle\int_{\mathbb{R}^{d-1}}\frac{\mathcal{O}_\alpha u(x)}{\langle\hat \alpha\cdot \nabla \phi(x)\rangle^{d+a}}\frac{d\alpha}{|\alpha|^{d-1+a}}.
			\end{align*}
			We can write 
			\begin{align*}
				\mathcal{L}_\phi u(x)=(2\pi)^{-\frac{d-1}{2}}\int_{\mathbb{R}^{d-1}} A(x,\xi)\hat{u}(\xi)e^{i\xi\cdot x}d\xi.
			\end{align*}
			Here
			\begin{equation}\label{nmcab}
				\begin{aligned}
					&A(x,\xi)=\int _{\mathbb{R}^{d-1}}B[\phi](x,\alpha)(2-e^{i\alpha\cdot \xi}-e^{-i\alpha\cdot \xi})\frac{d\alpha}{|\alpha|^{d+a}},\\
					&B[\phi](x,\alpha):=\frac{\langle\nabla \phi(x)\rangle}{\langle\hat \alpha\cdot \nabla \phi(x)\rangle^{d+a}}.
				\end{aligned}
			\end{equation}
			We remark that if $\phi\equiv 0$, then $\mathcal{L}_\phi=\Lambda^{1+a}$, which is the fractional Laplacian operator defined in \eqref{deffracla}.
			It is easy to check that there exists $c_0=c_0(\|\nabla \phi\|_{L^\infty})>0$ such that  
			\begin{align*}
				c_0\leq \frac{A(x,\xi)}{|\xi|^{1+a}} \leq  c_0^{-1},\ \ \ \ |\nabla_x^{j}\nabla _\xi^{l} A(x,\xi)|\leq (c_0^{-1} +\|\nabla \phi\|_{C^j})^j|\xi|^{1+a-l}, \ \ \forall j,l\in \mathbb{N}.
			\end{align*}
			Hence the operator $\mathcal{L}_\phi u$ satisfies the condition \eqref{defop}-\eqref{condop} with $s=1+a$ and $\M=(c_0^{-1} +\|\nabla \phi\|_{C^{m+3}})^{m+4}$.\\
			Now we rewrite the equation as 
			\begin{align}\label{dedeF}
				\partial_t u(x)+\mathcal{L}_{\phi}u(x)=-\mathcal{L}_{u}u(x)+\mathcal{L}_{\phi}u(x)-\langle\nabla u(x)\rangle\int_{\mathbb{R}^{d-1}}\mathcal{O}_\alpha u(x)P[u](x,\alpha)\frac{d\alpha}{|\alpha|^{d-1+a}}:=F[u,\phi](x).
			\end{align}
			Fix $\kappa, \varepsilon$ such that $0<1+a-\kappa\ll1$, $0<\varepsilon\leq \frac{1+a-\kappa}{10}$, and $m\in\mathbb{N}$. Denote
			\begin{equation}\label{nornmc}
				\begin{aligned}
					&\|f\|_{T}=\sup_{t\in[0,T]}\left(\|\nabla f(t)\|_{L^\infty}+t^\frac{m+\kappa}{1+a}\|\nabla f(t)\|_{\dot C^{m+\kappa}}\right),\\
					&\|f\|_{T,*}=\sup_{t\in[0,T]}\left(t^\frac{1}{5(1+a)}\|\nabla f(t)\|_{\dot C^{\frac{1}{5}}}+t^\frac{m+\kappa}{1+a}\|\nabla f(t)\|_{\dot C^{m+\kappa}}\right).
				\end{aligned}
			\end{equation}
			Let $\sigma, T>0$, we define 
			\begin{align*}
				\mathcal{X}^\sigma_{\phi,T}:=\left\{v\in L^\infty([0,T];\dot W^{1,\infty}):v|_{t=0}=u_0, \| v-\phi\|_T\leq \sigma\right\}.
			\end{align*}
			For any $v\in 	\mathcal{X}^\sigma_{\phi,T}$, we define a map $\mathcal{S}:v\to u$, where $u$ solves the Cauchy problem 
			\begin{align*}
				\partial_t u(x)+\mathcal{L}_{\phi}u(x)=F[v,\phi](x),\ \ \ \ u|_{t=0}=u_0.
			\end{align*}
			The remaining part of the proof is devoted to prove that $\mathcal{S}$ has a unique fixed point in $	\mathcal{X}^\sigma_{\phi,T}$ for some $\sigma,\phi,T$ that will be specified later.  For simplicity, denote $\tilde u=u-\phi$, $\tilde v=v-\phi$. We have 
			\begin{align*}
				\partial_t \nabla \tilde u(x)+\mathcal{L}_{\phi}\nabla \tilde u(x)=\nabla F[v,\phi](x)-\nabla (\mathcal{L}_\phi\phi)(x)+ R[u,\phi](x),
			\end{align*}
			where 
			\begin{align}\label{dedeR}
				R[u,\phi](x)=\int_{\mathbb{R}^{d-1}} \nabla_x B[\phi](x,\alpha)\mathcal{O}_\alpha \tilde u(x)\frac{d\alpha}{|\alpha|^{d-1+a}}.
			\end{align}
			Applying Theorem \ref{lemmain} with $F=F[v,\phi]$ and $G=R[u,\phi]-\nabla\mathcal{L}_\phi(\phi)$ to obtain that there exists $T=T(\|\nabla\phi\|_{C^{m+3}})>0$ such that 
			\begin{equation}\label{nmcmm}
				\begin{aligned}
					\|\tilde u\|_{T}\lesssim &\|\nabla (u_0-\phi)\|_{L^\infty}+\sup_{t\in[0,T]}\left(t^\frac{\kappa}{1+a}\|F[v,\phi](t)\|_{\dot C^{\kappa-a}}+t^{\frac{m+\kappa}{1+a}}\|F[v,\phi](t)\|_{\dot C^{m+\kappa-a}} \right)\\
					&+ \|R[u,\phi]\|_{L^1_TL^\infty}+\sup_{t\in[0,T]}t^{\frac{m}{1+a}+1}\|R[u,\phi](t)\|_{\dot C^m}+T\|\nabla (\mathcal{L}_\phi\phi)\|_{L^\infty}+T^{\frac{m}{1+a}+1}\|\nabla (\mathcal{L}_\phi\phi)\|_{\dot C^m}.
				\end{aligned}
			\end{equation}
			Applying Lemma \ref{lemnmc} below with $v_1=v,v_2=\phi$ in \eqref{df1}, and $v_1=u,v_2=\phi$ in \eqref{dr1}, one has for any $ n\in\mathbb{N}, \ n\leq m,$
			\begin{equation}\label{nmcfc}
				\begin{aligned}
					&\sup_{t\in[0,T]}\left(t^\frac{\kappa}{1+a}\|(F[v,\phi]-F[\phi,\phi])(t)\|_{\dot C^{\kappa-a}}+t^\frac{m+\kappa}{1+a}\|(F[v,\phi]-F[\phi,\phi])(t)\|_{\dot C^{m+\kappa-a}} \right)\\
					&\quad\quad\quad\quad\quad\quad\quad\quad\quad\quad\quad\quad\quad\quad\quad\lesssim \|\tilde v\|_T\|v\|_{T,*}(1+\|(v,\phi)\|_T)^{n+10},\\
					&\|R[u,\phi]\|_{L^1_TL^\infty}+\sup_{t\in[0,T]}t^{\frac{m}{1+a}+1}\|R[u,\phi](t)\|_{\dot C^m}\lesssim (T^{\frac{1}{1+a}}+T^{\frac{1+m}{1+a}})\|\nabla\phi\|_{C^{m+2}}(1+\|\nabla\phi\|_{C^{m+2}})^{m+3}\|\tilde u\|_T.
				\end{aligned}
			\end{equation}
			Moreover, by the smoothness of $\phi$, one can check that 
			\begin{equation}\label{nmcf}
				\begin{aligned}
					&	T^\frac{\kappa}{1+a}\|F[\phi,\phi]\|_{\dot C^{\kappa-a}}+T^\frac{m+\kappa}{1+a}\|F[\phi,\phi]\|_{\dot C^{m+\kappa-a}}+T\|\nabla (\mathcal{L}_\phi\phi)\|_{L^\infty}+T^\frac{m+\kappa}{1+a}\|\nabla (\mathcal{L}_\phi\phi)\|_{\dot C^m}\\
					&\quad\quad\quad\quad\quad\quad\quad\quad\quad\lesssim (	T^\frac{\kappa}{1+a}+T^\frac{m+\kappa}{1+a})\|\nabla\phi\|_{C^{m+2}}(1+\|\nabla\phi\|_{C^{m+2}})^{m+2}.
				\end{aligned}
			\end{equation}
			Take \eqref{nmcfc} and \eqref{nmcf} into \eqref{nmcmm} to get 
			\begin{equation}\label{esmap}
				\begin{aligned}
					\|\tilde u\|_{T}\lesssim &\|\nabla (u_0-\phi)\|_{L^\infty}+\|\tilde v\|_T\|(v,\phi)\|_{T,*}(1+\|(v,\phi)\|_T)^{n+10}+(T^{\frac{1}{1+a}}+T^{\frac{1+m}{1+a}})(1+\|\nabla\phi\|_{C^{m+2}})^{m+3}\|\tilde u\|_T\\
					&\quad\quad+ (	T^\frac{\kappa}{1+a}+T^\frac{m+\kappa}{1+a})\|\nabla\phi\|_{C^{m+2}}(1+\|\nabla\phi\|_{C^{m+2}})^{m+2}.
				\end{aligned}
			\end{equation}
			For  $v_1,v_2\in \mathcal{X}^\sigma_{\phi,T}$, denote 
			\begin{align}\label{notavu}
				\mathbf{v}=v_1-v_2, \quad \mathbf{u}= u_1-u_2=\mathcal{S}v_1-\mathcal{S}v_2,\quad\tilde v_i=v_i-\phi,\quad\tilde u_i=u_i-\phi,\ i=1,2.
			\end{align} Then
			\begin{align*}
				&\partial_t \nabla \mathbf{u}+\mathcal{L}_\phi \nabla  \mathbf{u}=\nabla  \mathsf{f}+\mathsf{g},\\
				&\nabla \mathbf{u}|_{t=0}=0,
			\end{align*}
			where 
			\begin{align*}
				&\mathsf{f}=F[v_1,\phi]-F[v_2,\phi],\\&
				\mathsf{g}=R[u_1,\phi]-R[u_2,\phi].
			\end{align*}
			Applying Theorem \ref{lemmain}, we obtain that there exists $T=T(\|\nabla \phi \|_{C^{m+3}})>0$ such that 
			\begin{equation}\label{escontra}
				\begin{aligned}
					\|\mathbf{u}\|_T\lesssim &\sup_{t\in[0,T]}t^\frac{\kappa}{1+a}\|\mathsf{f}(t)\|_{\dot C^{\kappa-a}}+\sup_{t\in[0,T]}t^\frac{m+\kappa}{1+a}\|\mathsf{f}(t)\|_{\dot C^{m+\kappa-a}}+ \|\mathsf{g}\|_{L^1_TL^\infty}+\sup_{t\in[0,T]}t^{\frac{m}{1+a}+1}\|\mathsf{g}(t)\|_{\dot C^{m}}.
				\end{aligned}
			\end{equation}
			By Lemma \ref{lemnmc}, we get 
			\begin{align*}
				&\sup_{t\in[0,T]}t^\frac{\kappa}{1+a}\|\mathsf{f}(t)\|_{\dot C^{\kappa-a}}+\sup_{t\in[0,T]}t^\frac{m+\kappa}{1+a}\|\mathsf{f}(t)\|_{\dot C^{m+\kappa-a}}+ \|\mathsf{g}\|_{L^1_TL^\infty}+\sup_{t\in[0,T]}t^{\frac{m}{1+a}+1}\|\mathsf{g}(t)\|_{\dot C^{m}}\\
				&\lesssim \|\mathbf{v}\|_T(\|v_1\|_{T,*}+\|\tilde v_2\|_{T})(1+\|(v_1,v_2)\|_T)^{m+10}+(T^{\frac{1}{1+a}}+T^{\frac{1+m}{1+a}})\|\nabla \phi\|_{C^{n+2}}(1+\|\nabla \phi\|_{C^{n+2}})^{n+3}\|\mathbf{u}\|_T.
			\end{align*}
			Hence we obtain 
			\begin{align}\label{recont}
				\|\mathbf{u}\|_T\lesssim \|\mathbf{v}\|_T(\|v_1\|_{T,*}+\|\tilde v_2\|_{T})(1+\|(v_1,v_2)\|_T)^{m+10}+(T^{\frac{1}{1+a}}+T^{\frac{1+m}{1+a}})\|\nabla \phi\|_{C^{m+2}}(1+\|\nabla \phi\|_{C^{m+2}})^{m+3}\|\mathbf{u}\|_T.
			\end{align}
			With \eqref{esmap} and \eqref{recont}, we can prove that $\mathcal{S}$ has a fixed point in $\mathcal{X}^\sigma_{\phi,T}$ with suitable choice  of  $\sigma, T$ and $\phi$.
			
			We consider two cases:
			\begin{itemize}
				\item i) If $\|\nabla u_0\|_{L^\infty}\leq \varepsilon_0$ for some $\varepsilon_0\ll1$. We take $ \phi\equiv 0$, then by Theorem \ref{thmcons}, \eqref{esmap} and \eqref{recont} hold for $T=\infty$, which can be written as 
				\begin{align*}
					&\|u\|_\infty\leq C_0(\|\nabla u_0\|_{L^\infty}+\|v\|_\infty^2(1+\|v\|_\infty)^{m+10})\leq  C_0(\varepsilon_0+\sigma^2(1+\sigma)^{m+10}),\ \ \ \forall v\in \mathcal{X}^\sigma_{0,\infty},\\
					&\|\mathbf{u}\|_\infty\leq C_0\|\mathbf{v}\|_\infty\|(v_1,v_2)\|_\infty(1+\|(v_1,v_2)\|_\infty)^{m+10},\ \ \ \forall  v_1,v_2\in \mathcal{X}^\sigma_{0,\infty}.
				\end{align*}
				Then
				by taking $\varepsilon_0\leq \frac{1}{100(1+2C_0)^{m+100}}$ and $\sigma =2C_0\varepsilon_0$, we have 
				$$\| u\|_{\infty}\leq \sigma,\ \ \ \ \ \text{and}\ \ \ \|\mathbf{u}\|_\infty\leq \frac{1}{2}\|\mathbf{v}\|_\infty.$$
				Hence $\mathcal{S}$ is a contraction map from $\mathcal{X}^\sigma_{0,\infty}$ to itself.
				\item ii) If \eqref{nmccon} holds, we take $\phi=u_0\ast \rho_{\eps_1}$, then there exists $T_0=T_0(\|\nabla u_0\|_{L^\infty},\eps_0,\eps_1)\in(0,1)$ such that \eqref{esmap} and \eqref{recont} hold for any $0<T\leq T_0$. With such choice of $\phi$, we have 
				\begin{align*}
					&\|\nabla\phi\|_{C^{m+2}}\lesssim \eps_1^{-(m+2)}\|u_0\|_{\dot W^{1,\infty}},\\
					&\|\phi\|_{T,*}\lesssim T^\frac{\varepsilon}{1+a}\|\nabla\phi\|_{C^{m+\kappa}}\lesssim T^\frac{1}{5(1+a)}\eps_1^{-(m+2)}\|u_0\|_{\dot W^{1,\infty}},\\
					&\|\phi\|_{T}\lesssim \|\nabla\phi\|_{L^\infty}+T^\frac{m+\kappa}{1+a}\|\nabla\phi\|_{\dot C^{m+\kappa}}\lesssim (1+T^\frac{m+\kappa}{1+a}\eps_1^{-(m+2)})\|u_0\|_{\dot W^{1,\infty}}.
				\end{align*}
				We can take $T$ small enough such that $T^\frac{1}{100}\eps_1^{-(m+10)^2}(1+\|u_0\|_{\dot W^{1,\infty}})\leq 1$. Then one has $\|\phi\|_{T,*}\lesssim T^\frac{1}{20}$ and $\|\phi\|_{T}\lesssim \|u_0\|_{\dot W^{1,\infty}}$.
				Put these estimates into \eqref{esmap} and \eqref{recont}, we get
				\begin{align}
					&\|\tilde u\|_{T}\leq \frac{5C_1\varepsilon_0}{4}+C_1\|\tilde v\|_T(\|\tilde v\|_T+T^\frac{1}{20})(1+\|\tilde v\|_T+\|u_0\|_{\dot W^{1,\infty}})^{m+10}+C_1T^\frac{1}{10}(1+\|\tilde u\|_T),\label{tut}\\
					&\|\mathbf{u}\|_T\leq C_1(1+\|(\tilde v_1, \tilde v_2)\|_T+T^\frac{1}{20})^{m+10}(\|(\tilde v_1, \tilde v_2)\|_T+T^\frac{1}{20})\|\mathbf{v}\|_T+C_1T^\frac{1}{10}\|\mathbf{u}\|_T.\label{but}
				\end{align}
				Let $\varepsilon_0=\frac{1}{(10+C_1+\|u_0\|_{\dot W^{1,\infty}})^{m+10}}$, $T_1=\min\left\{T_0,\left(\frac{\varepsilon_0\eps_1}{10+C_1+\|u_0\|_{\dot W^{1,\infty}}}\right)^{(m+100)^2}\right\}$, $\sigma'=2C_1\varepsilon_0$.\\
				If $v\in  \mathcal{X}_{T_1,u_0\ast \rho_{\eps_1}}^{\sigma'}$, then one has $\|\tilde v\|_{T_1}\leq \sigma'$, hence we obtain from \eqref{tut} that 
				\begin{align*}
					\|\tilde u\|_{T_1}\leq &\frac{3}{2}C_1\varepsilon_0+C_1T_1^\frac{1}{10}\|\tilde u\|_{T_1},
				\end{align*}
				which implies that 
				\begin{align*}
					\|\tilde u\|_{T_1}\leq\frac{\frac{5}{4}C_1\varepsilon_0}{1-C_1T_1^\frac{1}{10}}\leq \sigma'.	
				\end{align*}
				Hence
				$\mathcal{S}v=u\in  \mathcal{X}_{T_1, u_0\ast \rho_{\eps_1}}^{\sigma'}$.\\
				If $v_1,v_2\in  \mathcal{X}_{T_1, u_0\ast \rho_{\eps_1}}^{\sigma_2}$. Take $\phi= u_0\ast \rho_{\eps_1}$ in \eqref{mapwm},  one has $\|(v_1-\phi, v_2-\phi)\|_T\leq 2\sigma'$. Hence we obtain from \eqref{but} that 
				\begin{align*}
					\|\mathbf{u}\|_{T_1}\leq & 2^{m+7}C_1^2\varepsilon_0\|\mathbf{v}\|_{T_1}+C_1T_1^\frac{1}{10}\|\mathbf{u}\|_{T_1}.
				\end{align*}
				which leads to 
				\begin{align*}
					\|\mathbf{u}\|_{T_1}\leq &\frac{2^{m+7}C_1^2\varepsilon_0}{1-C_1T_1^\frac{1}{10}}\|\mathbf{v}\|_{T_1}\leq \frac{1}{2}\|\mathbf{v}\|_{T_1}.
				\end{align*}
				Hence $\mathcal{S}:\mathcal{X}_{T_1,u_0\ast \rho_{\eps_1}}^{\sigma'}\to \mathcal{X}_{T_1, u_0\ast \rho_{\eps_1}}^{\sigma'}$ is a contraction map. 
			\end{itemize}
			This completes the proof.
		\end{proof}
		\begin{lemma}\label{lemnmc}
			Let $\|\cdot \|_T$ and $\|\cdot\|_{T,*}$ be defined in \eqref{nornmc}.  Consider $F[u,\phi]$, $R[u,\phi]$ defined in \eqref{dedeF} and \eqref{dedeR}, respectively. Then for any integer $ n\in [0,m]$, and any $v_1,v_2$, there holds 
			\begin{align}
				\sup_{t\in[0,T]}t^\frac{n+\kappa}{1+a}\|(F[v_1,\phi]-F[v_2,\phi])(t)\|_{\dot C^{n+\kappa-a}}&\lesssim \|v_1-v_2\|_T(\|v_1\|_{T,*}+\|v_2-\phi\|_{T})(1+\|(v_1,v_2)\|_T)^{n+10},\label{df1}\\
				\sup_{t\in[0,T]}t^{\frac{n}{1+a}+1}\|(R[v_1,\phi]-R[v_2,\phi])(t)\|_{\dot C^n}+& \|(R[v_1,\phi]-R[v_2,\phi])\|_{L^1_TL^\infty}\nonumber\\
				&\lesssim (T^\frac{1}{1+a}+T^\frac{1+n}{1+a})\|\nabla\phi\|_{C^{n+2}}(1+\|\nabla\phi\|_{C^{n+2}})^{n+3}\|v_1-v_2\|_{T}.\label{dr1}
			\end{align}
		\end{lemma}	
		\begin{proof} In this proof, we still use the notations in \eqref{notavu}.
			Note that 
			\begin{align*}
				F[v_1,\phi]-F[v_2,\phi]=&(\mathcal{L}_{v_1}v_1-\mathcal{L}_{\phi}v_1)-(\mathcal{L}_{v_2}v_2-\mathcal{L}_{\phi}v_2)\\
				&-\int_{\mathbb{R}^{d-1}}\Big(\mathcal{O}_\alpha v_1(x)\langle\nabla v_1(x)\rangle P[v_1](x,\alpha)-\mathcal{O}_\alpha v_2(x)\langle\nabla v_2(x)\rangle P[v_2](x,\alpha)\Big)\frac{d\alpha}{|\alpha|^{d-1+a}}\\
				:=&J_1+J_2.
			\end{align*}
			We can further write 
			\begin{align*}
				J_1(x)=&\int_{\mathbb{R}^{d-1}}\left(B[v_1](x,\alpha)-B[v_2](x,\alpha)\right)\mathcal{O}_\alpha v_1(x)\frac{d\alpha}{|\alpha|^{d-1+a}}\\
				&\quad\quad\quad+\int_{\mathbb{R}^{d-1}}\left(B[v_2](x,\alpha)-B[\phi](x,\alpha)\right)\mathcal{O}_\alpha  \mathbf{v}(x)\frac{d\alpha}{|\alpha|^{d-1+a}}\\
				:=&J_{1,1}(x)+J_{1,2}(x),
			\end{align*}
			where $B[\cdot]$ is defined in \eqref{nmcab}. Elementary calculus yields that for any integer $n\in[0,m]$, 
			\begin{align*}
				\|\nabla_x^n(B[v_1]-B[v_2])\|_{L^\infty_{x,y}}&\lesssim \sum_{j=0}^n\|\nabla  \mathbf{v}\|_{\dot C^j}(\|\nabla v\|_{\dot C^{n-j}}+\|(\nabla v_1,\nabla v_2)\|_{\dot C^1}^{n-j})\lesssim t^{-\frac{n}{1+a}}\| \mathbf{v}\|_T(1+\|(v_1,v_2)\|_T)^n.
			\end{align*}
			Similarly, 
			\begin{align*}
				\|\nabla_x^n(B[v_1]-B[v_2])\|_{L^\infty_\alpha \dot C^{\kappa-a}_x}&\lesssim \sum_{j=0}^n\left(\|\nabla  \mathbf{v}\|_{\dot C^{j+\kappa-a}}(\|\nabla (v_1,v_2)\|_{\dot C^{n-j}}+\|\nabla (v_1,v_2)\|_{\dot C^1}^{n-j})\right.\\
				&\quad\quad+\left.\|\nabla  \mathbf{v}\|_{\dot C^{j}}(\|\nabla (v_1,v_2)\|_{\dot C^{n-j+\kappa-a}}+\|\nabla (v_1,v_2)\|_{\dot C^1}^{n-j+\kappa-a})\right)\\
				&\lesssim t^{-\frac{n+\kappa-a}{1+a}}\| \mathbf{v}\|_T(1+\|(v_1,v_2)\|_T)^n.
			\end{align*}
			Moreover, fix $0<\varepsilon<\frac{1-a}{2}$, then for any $\alpha \in \mathbb{R}^{d-1}$,
			\begin{align}\label{oav1}
				&\|\mathcal{O}_\alpha\nabla^nv_1\|_{L^\infty}\lesssim \min\left\{\|\nabla v_1\|_{\dot C^{n+a+\varepsilon}}|\alpha|^{a+\varepsilon},\|\nabla v_1\|_{\dot C^{n+a-\varepsilon}}|\alpha|^{a-\varepsilon}\right\}.
			\end{align}
			This implies that 
			\begin{align*}
				\int_{\mathbb{R}^{d-1}}\frac{\|\mathcal{O}_\alpha\nabla^nv_1\|_{L^\infty}}{|\alpha|^{d-1+a}}d\alpha\lesssim \|\nabla v_1\|_{\dot C^{n+a+\varepsilon}}^\frac{1}{2}\|\nabla v_1\|_{\dot C^{n+a-\varepsilon}}^\frac{1}{2}\lesssim t^{-\frac{n+a}{1+a}}\|v_1\|_{T,*}.
			\end{align*}
			Finally,  for any $\beta\in \mathbb{R}^{d-1}\backslash\{0\}$, by Lemma \ref{douint}, one has 
			\begin{align*}
				\|\delta_\beta \mathcal{O}_\alpha\nabla^n v_1\|_{L^\infty}\lesssim |\alpha|^{a}|\beta|^{\kappa-a} \min\left\{|\alpha| |\beta|^{-1},|\alpha|^{-1} |\beta|\right\}^\eps\|\nabla v_1\|_{\dot C^{n+\kappa}}.
			\end{align*}
			Hence,
			\begin{align*}
				\int_{{\mathbb{R}^{d-1}}}\frac{\|\delta_\beta\mathcal{O}_\alpha \nabla^nv_1(t)\|_{L^\infty}}{|\alpha|^{d-1+a}}d\alpha&\lesssim |\beta|^{\kappa-a}\|\nabla v_1(t)\|_{\dot C^{n+\kappa}}\int_{\mathbb{R}^{d-1}}\frac{\min\left\{|\alpha| |\beta|^{-1},|\alpha|^{-1} |\beta|\right\}^\eps}{|\alpha|^{d-1}} d\alpha\\
				&\lesssim|\beta|^{\kappa-a}\|\nabla v_1(t)\|_{\dot C^{n+\kappa}}\lesssim t^{-\frac{n+\kappa}{1+a}}|\beta|^{\kappa-a}\|v_1\|_{T,*}.
			\end{align*}
			Then we get 
			\begin{align*}
				\sup_{\beta}|\beta|^{-(\kappa-a)} \int_{\mathbb{R}^{d-1}}\frac{\|\delta_\beta\mathcal{O}_\alpha \nabla^lv_1(t)\|_{L^\infty}}{|\alpha|^{d-1+a}}d\alpha&\lesssim t^{-\frac{l+\kappa}{1+a}}\|v_1\|_{T,*}.
			\end{align*}
			We collect all the estimates above to get 
			\begin{align*}
				\|\nabla ^n J_{1,1}(t)\|_{\dot C^{\kappa-a}}\lesssim &\sum_{l=0}^n\left(\|\nabla _x^l\left(B[v_1]-B[v_2]\right)(t)\|_{L^\infty_\alpha\dot C^{\kappa-a}_x}\int_{\mathbb{R}^{d-1}}\frac{\|\mathcal{O}_\alpha\nabla ^{n-l} v_1(t)\|_{L^\infty}}{|\alpha|^{d-1+a}}d\alpha\right.\\
				&\quad\quad+\left.\|\nabla _x^l\left(B[v_1]-B[v_2]\right)(t)\|_{L^\infty_{x,\alpha}}\sup_{\beta}|\beta|^{-(\kappa-a)}\int_{\mathbb{R}^{d-1}}\frac{\|\delta_\beta\mathcal{O}_\alpha\nabla ^{n-l} v_1(t)\|_{L^\infty}}{|\alpha|^{d-1+a}}d\alpha\right)\\
				\lesssim &t^{-\frac{n+\kappa}{1+a}}\| \mathbf{v}\|_T\|v_1\|_{T,*}(1+\|(v_1,v_2)\|_T)^n.
			\end{align*}
			Similarly, we have 
			\begin{align*}
				\|\nabla ^n J_{1,2}(t)\|_{\dot C^{\kappa-a}}
				\lesssim &t^{-\frac{n+\kappa}{1+a}}\|v_2-\phi\|_T\| \mathbf{v}\|_{T,*}(1+\|(v_1,v_2)\|_T)^n.
			\end{align*}
			Then we obtain that 
			\begin{align}\label{j1}
				\|\nabla ^n J_{1}(t)\|_{\dot C^{\kappa-a}}
				\lesssim &t^{-\frac{n+\kappa}{1+a}}\| \mathbf{v}\|_{T}(\|v_1\|_{T,*}+\|v_2-\phi\|_T)(1+\|(v_1,v_2,\phi)\|_T)^n.
			\end{align}
			Then we deal with $J_2$, which can be further written as 
			\begin{align*}
				J_2(x)=&-\int_{\mathbb{R}^{d-1}}\Big(\mathcal{O}_\alpha v_1(x)\langle\nabla v_1(x)\rangle-\mathcal{O}_\alpha v_2(x)\langle\nabla v_2(x)\rangle\Big) P[v_1](x,\alpha)\frac{d\alpha}{|\alpha|^{d-1+a}}\\
				&-\int_{\mathbb{R}^{d-1}}\mathcal{O}_\alpha v_2(x)\langle\nabla v_2(x)\rangle \Big(P[v_1](x,\alpha)-P[v_2](x,\alpha)\Big)\frac{d\alpha}{|\alpha|^{d-1+a}}\\
				:=&J_{2,1}(x)+J_{2,2}(x).
			\end{align*} 
			For $J_2=J_{2,1}+J_{2,2}$, we have
			\begin{equation}\label{j21f}
				\begin{aligned}
					\|\nabla^nJ_{2,1}\|_{\dot C^{\kappa-a}}\lesssim &\sum_{l=0}^n\int_{\mathbb{R}^{d-1}}\|\nabla^{n-l}(\mathcal{O}_\alpha v_1\langle\nabla v_1\rangle-\mathcal{O}_\alpha v_2\langle\nabla v_2\rangle)\|_{L^\infty}\|\nabla^lP[v_1](\cdot,\alpha)\|_{\dot C^{\kappa-a}}\frac{d\alpha}{|\alpha|^{d-1+a}}\\
					&+\sum_{l=0}^n\int_{\mathbb{R}^{d-1}}\|\nabla^{n-l}(\mathcal{O}_\alpha v_1\langle\nabla v_1\rangle-\mathcal{O}_\alpha v_2\langle\nabla v_2\rangle)\|_{\dot C^{\kappa-a}}\|\nabla^lP[v_1](\cdot,\alpha)\|_{L^\infty}\frac{d\alpha}{|\alpha|^{d-1+a}}.
				\end{aligned}
			\end{equation}
			\begin{equation}\label{j22f}
				\begin{aligned}
					\|\nabla^nJ_{2,2}\|_{\dot C^{\kappa-a}}\lesssim &\sum_{l=0}^n\int_{\mathbb{R}^{d-1}}\|\nabla^{n-l}(\mathcal{O}_\alpha v_2\langle\nabla v_2\rangle)\|_{L^\infty}\|\nabla^l(P[v_1]-P[v_2])(\cdot,\alpha)\|_{\dot C^{\kappa-a}}\frac{d\alpha}{|\alpha|^{d-1+a}}\\
					&+\sum_{l=0}^n\int_{\mathbb{R}^{d-1}}\|\nabla^{n-l}(\mathcal{O}_\alpha v_2\langle\nabla v_2\rangle)\|_{\dot C^{\kappa-a}}\|\nabla^l(P[v_1]-P[v_2])(\cdot,\alpha)\|_{L^\infty}\frac{d\alpha}{|\alpha|^{d-1+a}}. 
				\end{aligned}
			\end{equation}
			By Lemma \ref{Hol},  we deduce that 
			\begin{equation}\label{esooo1}
				\begin{aligned}
					&\|\nabla^{n}(\mathcal{O}_\alpha v_2\langle\nabla v_2\rangle)(t)\|_{L^\infty}\lesssim t^{-\frac{n}{1+a}}\|v_2\|_{T,*}(1+\|v_2\|_{T})^{n+1}\min\{1,|\alpha| t^{-\frac{1}{1+a}}\}^\kappa,\\
					&\|\nabla^{n}(\mathcal{O}_\alpha v_2\langle\nabla v_2\rangle)(t)\|_{\dot C^{\kappa-a}}\lesssim t^{-\frac{n+\kappa-a}{1+a}}\|v_2\|_{T,*}(1+\|v_2\|_{T})^{n+2}\min\{1,|\alpha| t^{-\frac{1}{1+a}}\}^a,
				\end{aligned}
			\end{equation}
			and 
			\begin{equation}\label{esooo}
				\begin{aligned}
					&\|\nabla^n(\mathcal{O}_\alpha v_1\langle\nabla v_1\rangle-\mathcal{O}_\alpha v_2\langle\nabla v_2\rangle)(t)\|_{L^\infty}\lesssim t^{-\frac{n}{1+a}}\| \mathbf{v}\|_{T}(1+\|(v_1,v_2)\|_{T})^{n+1}\min\{1,|\alpha| t^{-\frac{1}{1+a}}\}^\kappa,\\
					&\|\nabla^{n}(\mathcal{O}_\alpha v_1\langle\nabla v_1\rangle-\mathcal{O}_\alpha v_2\langle\nabla v_2\rangle)(t)\|_{\dot C^{\kappa-a}} \lesssim t^{-\frac{n+\kappa-a}{1+a}}\| \mathbf{v}\|_{T}(1+\|(v_1,v_2)\|_{T})^{n+2}\min\{1,|\alpha| t^{-\frac{1}{1+a}}\}^a.
				\end{aligned}
			\end{equation}
			Recalling the definition of $P[v](x,\alpha)$ in \eqref{defP}, for any fixed $\alpha\in\mathbb{R}^d$ and  $\tau\in(0,1)$, we  apply Lemma \ref{lemcom} with $ f(x)=\langle x\rangle^{-(d+a)}$, and  
			\begin{align*}
				&  	g_1(x)=\tau \Delta_\alpha v_1(x)-(1-\tau)\Delta_\alpha v_1(x),\ \ \ g_2(x)=\tau \Delta_\alpha  v_2(x)-(1-\tau)\Delta_\alpha v_2(x),\\
				&	g_3(x)=\hat \alpha\cdot \nabla v_1(x),\quad\quad\quad\quad\ \ \ g_4(x)=\hat \alpha\cdot \nabla v_2(x).
			\end{align*}
			Then we get for any $\alpha\in\mathbb{R}^d\backslash\{0\}$,
			\begin{equation}\label{dfp}
				\begin{aligned}
					&\|\nabla_x^n(P[v_1]-P[v_2])(t,\cdot ,\alpha)\|_{L^\infty}\\
					&\lesssim \sum_{l=0}^n\|\Delta_{\pm \alpha} \mathbf{v}(t)\mp\hat \alpha \cdot \nabla  \mathbf{v}(t)\|_{\dot C^l}(\|(\nabla v_1,\nabla v_2)(t)\|_{\dot C^1_x}^{n-l}+\|(\nabla v_1,\nabla v_2)(t)\|_{\dot C^{n-l}_x})\\
					&\quad\quad+\sum_{n_1+n_2+n_3=n}\|\nabla  \mathbf{v}(t)\|_{\dot C^{n_1}_x}(\|\Delta_{\pm \alpha}v_1(t)\mp\hat \alpha \cdot \nabla v_1(t)\|_{\dot C^{n_2}_x}+\|\Delta_{\pm \alpha}v_2(t)\mp\hat \alpha \cdot \nabla v_2(t)\|_{\dot C^{n_2}_x})\\
					&\quad\quad\quad \times (\|(\nabla v_1,\nabla v_2)(t)\|_{\dot C^1_x}^{n_3}+\|(\nabla v_1,\nabla v_2)(t)\|_{\dot C^{n_3}_x})\\
					&\lesssim t^{-\frac{n}{1+a}}(\| \mathbf{v}\|_{T,*}+\| \mathbf{v}\|_T\|(v_1,v_2)\|_{T,*})(1+\| (v_1,v_2)\|_T)^n\min\{1,|\alpha| t^{-\frac{1}{1+a}}\}^\kappa,
				\end{aligned}
			\end{equation}
			and 
			\begin{equation}\label{dfph}
				\begin{aligned}
					&\|\nabla^n_x(P[v_1]-P[v_2])(t,\cdot,\alpha)\|_{\dot C^{\kappa-a}}\\
					&\lesssim \sum_{l=0}^n\left[\|\Delta_{\pm \alpha} \mathbf{v}(t)\mp\hat \alpha\cdot \nabla  \mathbf{v}(t)\|_{\dot C^{l+\kappa-a}}(\|(\nabla v_1,\nabla v_2)(t)\|_{\dot C^1}^{n-l}+\|(\nabla v_1,\nabla v_2)(t)\|_{\dot C^{n-l}})\right.\\
					&\quad\quad+\left.\|\Delta_{\pm \alpha} \mathbf{v}(t)\mp\hat \alpha\cdot \nabla  \mathbf{v}(t)\|_{\dot C^{l}}(\|(\nabla v_1,\nabla v_2)(t)\|_{\dot C^{\kappa-a}}^{\frac{n-l}{\kappa-a}}+\|(\nabla v_1,\nabla v_2)(t)\|_{\dot C^{n-l+\kappa-a}})\right]\\
					&\quad\quad+\sum_{n_1+n_2+n_3=n}\left(\|\nabla  \mathbf{v}(t)\|_{\dot C^{n_1+\kappa-a}}\|\Delta_{\pm \alpha}v_2(t)\mp\hat \alpha\cdot \nabla v_2(t)\|_{\dot C^{n_2}}(\|(\nabla v_1,\nabla v_2)(t)\|_{\dot C^{n_3}}+\|(\nabla v_1,\nabla v_2)(t)\|_{\dot C^{1}}^{n_3})\right.\\
					&\quad\quad+\|\nabla  \mathbf{v}(t)\|_{\dot C^{n_1}}\|\Delta_{\pm \alpha}v_2(t)\mp\hat y\cdot \nabla v_2(t)\|_{\dot C^{n_2+\kappa-a}}(\|(\nabla v_1,\nabla v_2)(t)\|_{\dot C^{n_3}}+\|(\nabla v_1,\nabla v_2)(t)\|_{\dot C^{1}}^{n_3})\\
					&\quad\quad+\left.\|\nabla  \mathbf{v}(t)\|_{\dot C^{n_1}}\|\Delta_{\pm y}v_2(t)\mp\hat \alpha\cdot \nabla v_2(t)\|_{\dot C^{n_2}}(\|(\nabla v_1,\nabla v_2)(t)\|_{\dot C^{n_3+\kappa-a}}+\|(\nabla v_1,\nabla v_2)(t)\|_{\dot C^{\kappa-a}}^{\frac{n_3}{\kappa-a}})\right)\\
					&\lesssim t^{-\frac{n+\kappa-a}{1+a}}(\| \mathbf{v}\|_{T,*}+\| \mathbf{v}\|_{T}\|(v_1,v_2)\|_{T,*})(1+\| (v_1,v_2)\|_T)^{n+5}\min\{1,|\alpha| t^{-\frac{1}{1+a}}\}^a.
				\end{aligned}
			\end{equation}
			Specially, \eqref{dfp} and \eqref{dfph} imply that 
			\begin{equation}\label{pv1}
				\begin{aligned}
					&\|\nabla_x^nP[v_1](t,\cdot ,\alpha)\|_{L^\infty}\lesssim  t^{-\frac{n}{1+a}}\|v_1\|_{T,*}(1+\| v_1\|_T)^{n+1}\min\{1,|\alpha| t^{-\frac{1}{1+a}}\}^\kappa,\\
					&\|\nabla^n_xP[v_1](t,\cdot,\alpha)\|_{\dot C^{\kappa-a}}\lesssim t^{-\frac{n+\kappa-a}{1+a}}\|v_1\|_{T,*}(1+\| v_1\|_T)^{n+8}\min\{1,|\alpha| t^{-\frac{1}{1+a}}\}^a.
				\end{aligned}
			\end{equation}
			Combining \eqref{j21f}, \eqref{esooo} with \eqref{pv1} to obtain that 
			\begin{align*}
				\|\nabla^nJ_{2,1}(t)\|_{\dot C^{\kappa-a}}&\lesssim t^{-\frac{n+\kappa-a}{1+a}} \int_{\mathbb{R}^{d-1}}\min\{1,|\alpha|t^{-\frac{1}{1+a}}\}^{\kappa+a}\frac{d\alpha}{|\alpha|^{d-1+a}}\| \mathbf{v}\|_T\|v_1\|_{T,*}(1+\|(v_1,v_2)\|_T)^{n+10}\\
				&\lesssim t^{-\frac{n+\kappa-a}{1+a}}\| \mathbf{v}\|_T\|v_1\|_{T,*}(1+\|(v_1,v_2)\|_T)^{n+10}.
			\end{align*}
			Moreover, by \eqref{j22f}, \eqref{esooo1}, \eqref{dfp} and \eqref{dfph} we obtain 
			\begin{align*}
				\|\nabla^mJ_{2,2}(t)\|_{\dot C^{\kappa-a}}&\lesssim t^{-\frac{m+\kappa-a}{1+a}} \int_{\mathbb{R}^{d-1}}\min\{1,|\alpha|t^{-\frac{1}{1+a}}\}^{\kappa+a}\frac{d\alpha}{|\alpha|^{d-1+a}}\| \mathbf{v}\|_T\|v_2\|_{T,*}(1+\|(v_1,v_2)\|_T)^{n+10}\\
				&\lesssim t^{-\frac{m+\kappa-a}{1+a}}\| \mathbf{v}\|_T\|v_2\|_{T,*}(1+\|(v_1,v_2)\|_T)^{n+10}.
			\end{align*}
			Hence, 
			\begin{align*}
				\|\nabla^nJ_{2}(t)\|_{\dot C^{\kappa-a}}
				&\lesssim t^{-\frac{n+\kappa-a}{1+a}}\| \mathbf{v}\|_T\|(v_1,v_2)\|_{T,*}(1+\| (v_1,v_2)\|_T)^{n+10},\ \ \forall n\in[0,m], n\in\mathbb{N}.
			\end{align*}
			Combining this with \eqref{j1},  we get \eqref{df1}.\\
			It remains to estimate the lower order term \begin{align*}
				R[v_1,\phi]-R[v_2,\phi]=\int_{\mathbb{R}^{d-1}}\nabla_x B[\phi](x,\alpha)\mathcal{O}_\alpha(v_1-v_2)(x)\frac{d\alpha}{|\alpha|^{d-1+a}}.
			\end{align*}
			Note that 
			\begin{align}
				\|\nabla_x^l B[\phi](\cdot,\alpha)\|_{L^\infty}\lesssim\|\nabla\phi\|_{C^{n+2}}^{l+2},\ \ \ \forall l\in\mathbb{N}^+.
			\end{align}
			Moreover, by \eqref{oav1}, for any $0\leq l\leq m$, we get 
			\begin{align*}
				\|\mathcal{O}_\alpha\nabla^{l} \mathbf{v}(t)\|_{L^\infty}\lesssim |\alpha|^at^{-\frac{l+a}{1+a}}\| \mathbf{v}\|_T\min\left\{|\alpha|t^{-\frac{1}{1+a}}, (|\alpha|t^{-\frac{1}{1+a}})^{-1}\right\}^\eps.
			\end{align*}
			Then for any $n\in \mathbb{N}$, $n\leq m$, we have 
			\begin{align*}
				& \|\nabla^n(R[v_1,\phi]-R[v_2,\phi])(t)\|_{L^\infty}\\&\quad\quad\lesssim \sum_{l=0}^n\int_{\mathbb{R}^{d-1}}\|\nabla_x^{l+1}B[\phi](\cdot,\alpha)\|_{L^\infty}\|\mathcal{O}_\alpha\nabla^{n-l} \mathbf{v}(x)\|_{L^\infty}\frac{d\alpha}{|\alpha|^{d-1+a}}\\
				&\quad\quad\lesssim (t^{-\frac{n+a}{1+a}}+t^{-\frac{a}{1+a}})\|\nabla\phi\|_{C^{n+2}}(1+\|\nabla\phi\|_{C^{n+2}})^{n+3}\| \mathbf{v}\|_{T}\int_{\mathbb{R}^{d-1}}\min\left\{|\alpha|t^{-\frac{1}{1+a}}, (|\alpha|t^{-\frac{1}{1+a}})^{-1}\right\}^\eps\frac{d\alpha}{|\alpha|^{d-1}}\\
				&\quad\quad\lesssim (t^{-\frac{n+a}{1+a}}+t^{-\frac{a}{1+a}})\|\nabla\phi\|_{C^{n+2}}(1+\|\nabla\phi\|_{C^{n+2}})^{n+3}\| \mathbf{v}\|_{T}.
			\end{align*}
			Hence
			\begin{align*}
				\sup_{t\in[0,T]} t^{\frac{n}{1+a}+1} \|\nabla^n(R[v_1,\phi]-R[v_2,\phi])(t)\|_{L^\infty}\lesssim (T^\frac{1}{1+a}+T^\frac{1+n}{1+a})\|\nabla\phi\|_{C^{n+2}}(1+\|\nabla\phi\|_{C^{n+2}})^{n+3}\| \mathbf{v}\|_{T}.
			\end{align*}
			And 
			\begin{align*}
				\|R[v_1,\phi]-R[v_2,\phi]\|_{L^1_TL^\infty}\lesssim T^\frac{1}{1+a}\|\nabla\phi\|_{C^{3}}(1+\|\nabla\phi\|_{C^1})^3\| \mathbf{v}\|_T.
			\end{align*}
			Then we obtain \eqref{dr1}. This completes the proof.
		\end{proof}
		
		\subsection{The Willmore flow}
		In this section we will consider the Willmore flow equation. Assume $f:\Sigma\rightarrow \mathbb{R}^2$ to be an immersion, and denote the Willmore fractional as follows
		\begin{equation*}
			\mathcal{W}(f)=\int_{\Sigma}|H|^2dA,
		\end{equation*}
		where $H$ is the mean curvature vector of $f$, and $dA$ is the area element. The Willmore surface consists of the critical points of the Willmore fractional $\mathcal{W}$ and is described by the following Euler-Lagrange equation
		\begin{equation*}
			V=\Delta_\Sigma H+\frac{1}{2}H^3-2HK,
		\end{equation*}
		where $K$ is the Gauss curvature of $\Sigma$, $\Delta_\Sigma$ is the Laplace-Beltrami operator of $\Sigma$. We consider the case where $\Sigma$ is a graph of a function on $\mathbb{R}^2$, which means that $\Sigma=\{(x,u(x))|x\in\mathbb{R}^2\}$, and recall $\langle\nabla u\rangle$ is Japanese bracket, then we can see that 
		\begin{align*}
			V=-\frac{u_t}{\langle\nabla u\rangle},\ H=\operatorname{div}\left(\frac{\nabla u}{\langle\nabla u\rangle}\right),\ K=\frac{\operatorname{det} \nabla^2u}{\langle\nabla u\rangle^4},\\
			\Delta_\Sigma H=\frac{1}{\langle\nabla u\rangle}\cdot \left((\langle\nabla u\rangle \mathrm{Id}-\frac{\nabla u\otimes\nabla u}{\langle\nabla u\rangle})\nabla H\right).
		\end{align*}
		The Willmore flow equation can be rewritten as the following equation for $u:(0,T)\times \mathbb{R}^2\to \mathbb{R}$ (see \cite{DD2006} for details).
		\begin{equation}\label{eqwm}
			\begin{aligned}
				&u_{t}+F[u]=0 \quad \text { in } \quad (0, T)\times\mathbb{R}^{2},\\
				&u|_{t=0}=u_0,
			\end{aligned}
		\end{equation}
		where the nonlinear term reads
		\begin{align*}
			F[u]=\langle\nabla u\rangle \operatorname{div}\left(\frac{1}{\langle\nabla u\rangle}\left(\left(\mathrm{Id}-\frac{\nabla u \otimes \nabla u}{\langle\nabla u\rangle^{2}}\right) \nabla\left(\langle\nabla u\rangle H\right)-\frac{1}{2} H^{2} \nabla u\right)\right).
		\end{align*}
		It should be noticed that we want to show the existence and regularity of solutions with critical initial data. Since our method doesn't rely on dimension, we will consider more general case for any dimension $d$ in the following. \\
		We give a short review about the work before. In \cite{ES2001} and \cite{ES2004}, Kuwert and Sch\"{a}tzle proved that if the Willmore fractional of initial data is small, then the Willmore flow equation has a smooth global solution, which converges exponentially to sphere. In their another paper \cite{ES2002}, they proved the existence of local solution, and gived the lower bound of the time interval which only relies on the concentration of the curvature of initial data. In \cite{CFS2009}, Chill, Fa\v{s}angov\'{a} and Sch\"{a}tzle showed that for initial data which is close to a $C^k(k\geq 2)$ local minimizer in $W^{2,2}\cap C^1$ norm, then there exists a global solution to Willmore flow equation, and the solution converges to a local minimizer(up to a diffeomorphism). For singularity formation of the Willmore flow, Blatt proved that there are surfaces that become singular under the influence of Willmore flow in \cite{B2009}, and Mayer and Simonett \cite{MS2002} present similar result by numerical methods.\\
		We want to emphasize the result in  \cite{Koch2011}, where Koch and Lamm proved that for initial data $f_0$ satisfying
		\begin{equation*}
			\|f_0\|_{C^{0,1}}<\delta
		\end{equation*}
		for some small $\delta$, there exists an analytic solution $f$ to the Willmore flow equation which satisfies that
		\begin{equation}\label{Xinfty}
			\|f\|_{X_\infty}:=\sup_{t}\|\nabla f\|_{L^\infty}+\sup_{x,R}R^{\frac{2}{d+6}}\|\nabla^2f\|_{L^{d+6}(B_R(x)\times (\frac{R^4}{2},R^4))}<C\|f_0\|_{C^{0,1}}.
		\end{equation}
		In our paper, we can prove both local existence of solution for large initial data and global existence for small initial data. The main result in this section is as follows. 
		\begin{theorem} \label{thmwill}~Fix $\kappa$ such that $0<4-\kappa\ll1$, and $m\in\mathbb{N}$. There exists $\varepsilon_0>0$ such that \\
			i)For any $u_0$ with $\|\nabla u_0\|_{L^\infty}\leq \varepsilon_0$, there exists a unique solution $u$ to \eqref{eqwm} such that
			\begin{align*}
				\sup_{t>0}	(\|\nabla u(t)\|_{L^\infty}+t^\frac{m+\kappa}{4}\|\nabla u(t)\|_{\dot C^{m+\kappa}})\leq C\|\nabla u_0\|_{L^\infty}.
			\end{align*}
			ii) For any initial data $u_0\in \dot W^{1,\infty}$, if there exists $\eps_1>0$ such that 
			\begin{align}\label{wmcon}
				\|\nabla u_0-\nabla u_0\ast \rho_{\eps_1}\|_{L^\infty}\leq \varepsilon_0,
			\end{align} then there exists $T>0$ such that \eqref{eqwm} admits a unique solution  $u$ satisfying
			\begin{align*}
				\sup_{t\in[0,T]}	(\|\nabla u(t)\|_{L^\infty}+t^\frac{m+\kappa}{4}\|\nabla u(t)\|_{\dot C^{m+\kappa}})\leq C \|\nabla u_0\|_{L^\infty}.
			\end{align*}
		\end{theorem}	
		\begin{proof}[Proof of Theorem \ref{thmwill}]
			Note that 
			\begin{align*}
				\langle\nabla u\rangle H=\langle\nabla u\rangle\operatorname{div}\left(\frac{\nabla u}{\langle\nabla u\rangle}\right)=M(\nabla u):\nabla^2 u,
			\end{align*}
			where we denote $M(b)=\mathrm{Id}-\frac{b\otimes b}{\langle \b\rangle^2}$ for $b\in \mathbb{R}^d$.
			Then we can rewrite 
			\begin{align*}
				F[u]&=\langle\nabla u\rangle \operatorname{div}\left(\frac{1}{\langle\nabla u\rangle}\left(\left(\mathrm{Id}-\frac{\nabla u \otimes \nabla u}{\langle\nabla u\rangle^{2}}\right) \nabla\left(\langle\nabla u\rangle H\right)-\frac{1}{2} H^{2} \nabla u\right)\right)\\
				&=A[\nabla u]\nabla^4u+F_1[u]+F_2[u]+F_3[u],
			\end{align*}
			with 
			\begin{align*}
				&A[b]\nabla^4=A^{ijkl}[b]\partial_{ijkl},\ \ \ \ A^{ijkl}(b)=M_{ij}(b)M_{kl}(b),\\
				&    F_1[u]=-\langle\nabla u\rangle \operatorname{div}\left(\frac{ H^{2} \nabla u}{2\langle\nabla u\rangle}\right),\\
				&  F_2[u]=\langle\nabla u\rangle \partial_i\left(\frac{M_{ij}(\nabla u)}{\langle\nabla u\rangle}\right)\partial_j(M_{kl}(\nabla u)\partial_{kl} u),
				\\& F_3[u]=M_{ij}(\nabla u)\partial_{ij}(M_{kl}(\nabla u))\partial_{kl}u.
			\end{align*}
			Here we used the Einstein summation notation.\\
			Note that by Cauchy-Schwarz inequality,  $$A^{ijkl}[b]\xi_i\xi_j\xi_k\xi_l\geq \frac{1}{\langle \b\rangle^4}|\xi|^4, $$ for any $\xi\in \mathbb{R}^d$. \\
			For simplicity, fix $m\in\mathbb{N}$ and $\kappa$ such that $0<4-\kappa\ll 1$ as stated in Theorem \ref{thmwill}. Denote 
			\begin{equation}\label{defnowi}
				\begin{aligned}
					&		\|h\|_{T}=\sup_{t\in[0,T]}(\|\nabla h(t)\|_{L^\infty}+t^{\frac{m+\kappa }{4}}\|\nabla h(t)\|_{\dot C^\kappa}),\\
					&\|h\|_{T,*}=\sup_{t\in[0,T]}(t^\frac{1}{12}\|\nabla h(t)\|_{\dot C^\frac{1}{3}}+t^{\frac{m+\kappa }{4}}\|\nabla h(t)\|_{\dot C^\kappa}).
				\end{aligned}
			\end{equation}
			Let $\phi\in C^\infty(\mathbb{\mathbb{R}}^d,\mathbb{\mathbb{R}})$, $\sigma,T>0$, we consider the set
			$$
			\mathcal{X}_{T,\phi}^\sigma=\left\{v\in L_T^\infty\dot W^{1,\infty}: v|_{t=0}=v_0, \ \| v-\phi\|_{T}\leq \sigma\right\}.
			$$ 
			For any $v\in	\mathcal{X}_{T,\phi}^\sigma$, where $T,\sigma,\phi$ will be fixed later, we define a map $\mathcal{S}v=u$, where $u$ is the weak solution to the Cauchy problem
			\begin{equation}\label{mapwm}
				\begin{aligned}
					&		\partial_tu+A[\nabla\phi]\nabla ^4 u=(A[\nabla\phi]-A[\nabla v])\nabla ^4 v-\sum_{k=1}^3F_k[v],\\
					&u|_{t=0}=u_0.
				\end{aligned} 
			\end{equation}
			Denote $\tilde u=u-\phi$, $\tilde v=v-\phi$.	
			We have the equation 
			\begin{equation*}
				\begin{aligned}
					&\partial_t \partial_j\tilde u+A[\nabla\phi]\nabla ^4 \partial_j\tilde u=\partial_j (N[v,\phi]+\tilde F[v])+R_j[u,\phi],\ \ \ j=1,2,\\
					&\partial_j\tilde u|_{t=0}=\partial_j(u_0-\phi),
				\end{aligned} 
			\end{equation*}
			where 
			\begin{equation}\label{defrhs}
				\begin{aligned}
					&N[v,\phi]=(A[\nabla\phi]-A[\nabla v])\nabla ^4 v,\ \ \ 
					\tilde F[v]=\sum_{k=1}^3F_k[v],\\
					&R_j[u,\phi]=-\partial_j(A[\nabla\phi])\nabla ^4u-A[\nabla\phi]\nabla^4\partial_j\phi.
				\end{aligned}
			\end{equation}
			By Theorem \ref{lemmain} with $F=N[v,\phi]+\tilde F[v]$ and $G=\tilde F[v])+R_j[u,\phi]$, there exists $T=T(\|\nabla\phi\|_{C^{m+4}})>0$ such that 
			\begin{equation}\label{bla1}
				\begin{aligned}
					\|\tilde u\|_T\lesssim &\|\nabla(u_0-\phi)\|_{L^\infty}+ \sup_{t\in[0,T]}(t^\frac{\kappa}{4}\|N[v,\phi](t)\|_{\dot C^{\kappa-3}}+t^\frac{m+\kappa}{4}\|N[v,\phi](t)\|_{\dot C^{m+\kappa-3}})\\
					&\quad+\sup_{t\in[0,T]}(t^\frac{\kappa}{4}\|\tilde F[v](t)\|_{\dot C^{\kappa-3}}+t^\frac{m+\kappa}{4}\|\tilde F[v](t)\|_{\dot C^{m+\kappa-3}})\\
					&\quad+ \|R[u,\phi]\|_{L^1_TL^\infty}+\sup_{t\in[0,T]}t^{\frac{m}{4}+1}\|\nabla^mR[u,\phi](t)\|_{L^\infty}.
				\end{aligned} 
			\end{equation}
			Observe that $N[\phi,\phi]\equiv 0$. Then applying Lemma \ref{lemwi} with $(v_1,v_2)=(v,\phi)$,  we obtain 
			\begin{align}
				&\sup_{t\in[0,T]}(t^\frac{\kappa}{4}\|N[v,\phi](t)\|_{\dot C^{\kappa-3}}+t^\frac{m+\kappa}{4}\|N[v,\phi](t)\|_{\dot C^{m+\kappa-3}})\lesssim \|\tilde v\|_T\|v\|_{T,*}(1+\|(v,\phi)\|_{T,*})^{m+5},\label{nnn}\\
				&\sup_{t\in[0,T]}(t^\frac{\kappa}{4}\|(\tilde F[v]-\tilde F[\phi])(t)\|_{\dot C^{\kappa-3}}+t^\frac{m+\kappa}{4}\|(\tilde F[v]-\tilde F[\phi])(t)\|_{\dot C^{m+\kappa-3}})\lesssim \|\tilde v\|_T\|(v,\phi)\|_{T,*}(1+\|(v,\phi)\|_{T,*})^{m+5}.\label{fff1}
			\end{align}
			Moreover, we have 
			\begin{align*}
				\|\tilde F[\phi]\|_{\dot C^{\kappa-3}}+\|\tilde F[\phi]\|_{\dot C^{m+\kappa-3}}\lesssim (1+\|\nabla\phi\|_{C^{m+4}})^{m+5}. 
			\end{align*}
			This yields that 
			\begin{align*}
				\sup_{t\in[0,T]}(t^\frac{\kappa}{4}\|\tilde F[\phi](t)\|_{\dot C^{\kappa-3}}+t^\frac{m+\kappa}{4}\|\tilde F[\phi](t)\|_{\dot C^{m+\kappa-3}})\lesssim (T^\frac{\kappa}{4}+T^\frac{m+\kappa}{4})(1+\|\nabla\phi\|_{C^{m+4}})^{m+5}.
			\end{align*}
			Combining this with \eqref{fff1} to get 
			\begin{equation}\label{fff}
				\begin{aligned}
					&\sup_{t\in[0,T]}(t^\frac{\kappa}{4}\|\tilde F[v](t)\|_{\dot C^{\kappa-3}}+t^\frac{m+\kappa}{4}\|\tilde F[v](t)\|_{\dot C^{m+\kappa-3}})\\
					&\quad\quad\quad\quad\lesssim \|\tilde v\|_T\|(v,\phi)\|_{T,*}(1+\|(v,\phi)\|_{T,*})^{m+5}+(T^\frac{\kappa}{4}+T^\frac{m+\kappa}{4})(1+\|\nabla\phi\|_{C^{m+4}})^{m+5}.
				\end{aligned}
			\end{equation}
			Finally, we have $R[u,\phi]\equiv 0$ if $\phi\equiv 0$. Moreover,
			\begin{align*}
				&\|R[u,\phi]\|_{L^\infty}\lesssim \|\nabla \phi\|_{C^1}\|\nabla^4u\|_{L^\infty}+\|\nabla \phi\|_{C^4},\\
				&\|\nabla^mR[u,\phi]\|_{L^\infty}\lesssim \|\nabla \phi\|_{C^{m+4}}(\|\nabla^4u\|_{C^m}+(1+\|\nabla \phi\|_{C^{m+4}})^{m+5}).
			\end{align*}
			Hence
			\begin{equation}\label{rrr1}
				\begin{aligned}
					&\|R[u,\phi]\|_{L^1_TL^\infty}+\sup_{t\in[0,T]}t^{\frac{m}{4}+1}\|\nabla^mR[u,\phi](t)\|_{L^\infty}\\
					&\quad\quad\quad\quad\lesssim (T^\frac{1}{4}+T^{\frac{m}{4}+1})\|\nabla \phi\|_{C^{m+4}}(1+\|\nabla \phi\|_{C^{m+4}})^{m+5}(1+\|u\|_T).
				\end{aligned}
			\end{equation}
			We collect \eqref{nnn}, \eqref{fff} and \eqref{rrr1} into \eqref{bla1} to obtain 
			\begin{equation}\label{wmre}
				\begin{aligned}
					\|\tilde u\|_T\lesssim &\|\nabla(u_0-\phi)\|_{L^\infty}+\|\tilde v\|_T\|(v,\phi)\|_{T,*}(1+\|(v,\phi)\|_{T,*})^{m+5}\\
					&\ \ +(T^\frac{1}{4}+T^{\frac{m}{4}+1})\|\nabla \phi\|_{C^{m+4}}(1+\|\nabla\phi\|_{C^{m+4}})^{m+5}(1+\|u\|_T).
				\end{aligned}
			\end{equation}
			Then we do the contraction estimates. Consider $v_1,v_2\in \mathcal{X}^\sigma_{T,\phi}$ for some $\sigma,T>0$ and $\phi\in C^\infty(\mathbb{R}^d,\mathbb{R})$. Denote $\mathbf{v}=v_1-v_2$, $\mathbf{u}=u_1-u_2=\mathcal{S}v_1-\mathcal{S}v_2$.
			For $j=1,2$, We have
			\begin{align*}
				&\partial_t \partial_j\mathbf{u}+A [\nabla\phi]\nabla ^4 \partial_j\mathbf{u}
				=\partial_j (N[v_1,\phi]-N[v_2,\phi]+\tilde F[v_1]-\tilde F[v_2])+R_{j}[u_1,\phi]-R_{j}[u_2,\phi],\ \ \text{in}\ (0,T)\times\mathbb{R}^d, \\
				& \mathbf{u}|_{t=0}=0,\quad\text{in}\ \mathbb{R}^d.
			\end{align*} 
			By Theorem \ref{lemmain} with $F=N[v_1,\phi]-N[v_2,\phi]+\tilde F[v_1]-\tilde F[v_2]$, $G=R_{j}[u_1,\phi]-R_{j}[u_2,\phi]$, there exists $T=T(\|\nabla\phi\|_{C^{m+4}})>0$ such that 
			\begin{equation}\label{difu1}
				\begin{aligned}
					\|\mathbf{u}\|_T\lesssim&\sup_{t\in[0,T]}(t^\frac{\kappa}{4} \|(N[v_1,\phi]-N[v_2,\phi])(t)\|_{\dot C^{\kappa-3}}+t^\frac{m+\kappa}{4} \|(N[v_1,\phi]-N[v_2,\phi])(t)\|_{\dot C^{m+\kappa-3}})\\
					&+\sup_{t\in[0,T]}(t^\frac{\kappa}{4} \|(\tilde F[v_1]-\tilde F[v_2])(t)\|_{\dot C^{\kappa-3}}+t^\frac{m+\kappa}{4} \|(\tilde F[v_1]-\tilde F[v_2])(t)\|_{\dot C^{m+\kappa-3}})\\
					&+\|R[u_1,\phi]-R[u_2,\phi]\|_{L^1_TL^\infty}+\sup_{t\in[0,T]}t^{\frac{m}{4}+1}\|\nabla^m (R[u_1,\phi]-R[u_2,\phi])(t)\|_{L^\infty}.
				\end{aligned}
			\end{equation}
			Applying Lemma \ref{lemwi} to obtain 
			\begin{equation}\label{es111}
				\begin{aligned}
					& \sup_{t\in[0,T]}(t^\frac{\kappa}{4} \|(N[v_1,\phi]-N[v_2,\phi])(t)\|_{\dot C^{\kappa-3}}+t^\frac{m+\kappa}{4} \|(N[v_1,\phi]-N[v_2,\phi])(t)\|_{\dot C^{m+\kappa-3}})\\
					&\quad\quad+\sup_{t\in[0,T]}(t^\frac{\kappa}{4} \|(\tilde F[v_1]-\tilde F[v_2])(t)\|_{\dot C^{\kappa-3}}+t^\frac{m+\kappa}{4} \|(\tilde F[v_1]-\tilde F[v_2])(t)\|_{\dot C^{m+\kappa-3}})\\
					&\quad\quad\quad\quad\quad\quad\quad\lesssim (1+\|(v_1,v_2,\phi)\|_{T,*})^{m+5}(\|(v_1,v_2)\|_{T,*}+\|v_2-\phi\|_T)\|\mathbf{v}\|_T.
				\end{aligned}
			\end{equation}
			Moreover, we have 
			\begin{align*}
				R[u_1,\phi]-R[u_2,\phi]=-\nabla(A[\nabla\phi])\nabla ^4\mathbf{u}.
			\end{align*}
			Hence 
			\begin{align*}
				&\|(R[u_1,\phi]-R[u_2,\phi])(t)\|_{L^\infty}\lesssim t^{-\frac{3}{4}}\|\nabla \phi\|_{C^1}\|\mathbf{u}\|_T,\\
				& \|\nabla^m(R[u_1,\phi]-R[u_2,\phi])(t)\|_{L^\infty}\lesssim (t^{-\frac{3}{4}}+t^{-\frac{m+3}{4}})(1+\|\nabla \phi\|_{C^{m+4}})^{m+5}\|\nabla \phi\|_{C^{m+4}}\|\mathbf{u}\|_T.
			\end{align*}
			Then we obtain
			\begin{equation}\label{rerrr}
				\begin{aligned}
					& \|R[u_1,\phi]-R[u_2,\phi]\|_{L^1_TL^\infty}+\sup_{t\in[0,T]}t^{\frac{m}{4}+1}\|\nabla^m (R[u_1,\phi]-R[u_2,\phi])(t)\|_{L^\infty}\\
					&\quad\quad\quad\quad\quad\quad\quad\quad\quad\quad\quad\quad\lesssim (T^\frac{1}{4}+T^{\frac{m}{4}+1})(1+\|\nabla \phi\|_{C^{m+4}})^{m+5}\|\nabla \phi\|_{C^{m+4}}\|\mathbf{u}\|_T.
				\end{aligned}
			\end{equation}
			Combining \eqref{es111}, \eqref{rerrr} with \eqref{difu1} to get 
			\begin{equation}\label{hahawm}
				\begin{aligned}
					\|\mathbf{u}\|_T\lesssim &(1+\|(v_1,v_2,\phi)\|_{T,*})^{m+5}(\|(v_1,v_2)\|_{T,*}+\|v_2-\phi\|_T)\|\mathbf{v}\|_T\\
					&+(T^\frac{1}{4}+T^{\frac{m}{4}+1})\|\nabla \phi\|_{C^{m+4}}(1+\|\nabla \phi\|_{C^{m+4}})^{m+5}\|\mathbf{u}\|_T.
				\end{aligned}
			\end{equation}	
			With \eqref{wmre} and \eqref{hahawm}, we can prove that the map $\mathcal{S}$ has a unique fixed point in $\mathcal{X}^\sigma_{T,\phi}$ by suitable choice of $\sigma,\phi$ and $T$.
			We consider two cases:
			\begin{itemize}
				\item 	Suppose $\|\nabla u_0\|_{L^\infty}\leq \varepsilon_0$ for some $\varepsilon_0\ll 1 $ that will be fixed later. Take $\phi\equiv 0$ in \eqref{mapwm}, then $A [\nabla\phi]$ is constant.  By Theorem \ref{thmcons}  we obtain that \eqref{wmre} and \eqref{hahawm} hold for  $T=\infty$. After taking $\phi\equiv 0$ and $T=\infty$, we rewrite \eqref{wmre} and \eqref{hahawm} as follows.
				\begin{align*}
					&	\|u\|_\infty\leq C_0\|\nabla  u_0\|_{L^\infty}+C_0\|(v_1,v_2)\|_\infty^2(1+\|(v_1,v_2)\|_\infty)^{m+5},\\
					&\|u_1-u_2\|_\infty\leq C_0 (1+\|(v_1,v_2)\|_{\infty})^{m+5}\|(v_1,v_2)\|_{\infty}\|v_1-v_2\|_\infty.
				\end{align*}
				Take $\varepsilon_0=\frac{1}{100(2C_0+1)^{m+10}}$ and $\sigma_1=2C_0 \varepsilon_0$. Then for any $v\in \mathcal{X}_{\infty,0}^{\sigma_1}$, we have $\|v\|_{\infty}\leq 2C_0 \varepsilon_0$. Then 
				\begin{align*}
					\|  u\|_\infty\leq C_0\varepsilon_0+4C_0^3\varepsilon_0^2(2C_0\varepsilon_0+1)^{m+5}\leq \sigma_1.
				\end{align*}
				Hence $\mathcal{S}v=u\in \mathcal{X}_{\infty,0}^{\sigma_1}$.
				
				Moreover, for any $v_1,v_2\in \mathcal{X}^{\sigma_1}_{\infty,0}$,   
				\begin{align*}
					\|\mathcal{S}v_1-\mathcal{S}v_2\|_\infty=\|u_1-u_2\|_\infty\textbf{}\leq 2C_0\sigma(1+2\sigma)^{m+5}\|v_1-v_2\|_\infty\leq \frac{1}{2}\|v_1-v_2\|_\infty.
				\end{align*}
				Hence $\mathcal{S}:\mathcal{X}^{\sigma_1}_{\infty,0}\to \mathcal{X}^{\sigma_1}_{\infty,0}$ is a contraction map. Hence there exists a unique $u\in \mathcal{X}^{\sigma_1}_{\infty,0}$ such that $\mathcal{S}u=u$.

				\item 	If the condition \eqref{wmcon} holds, take $\phi= u_0\ast \rho_{\eps_1}$ in \eqref{mapwm}, then there exists $T_0=T_0(\|\nabla u_0\|_{L^\infty},\eps_1)\in(0,1)$ such that \eqref{wmre} and \eqref{hahawm} hold for any $0<T\leq T_0$. With such choice of $\phi$, we have 
				\begin{align*}
					&\|\nabla\phi\|_{C^{m+4}}\lesssim \eps_1^{-(m+4)}\|u_0\|_{\dot W^{1,\infty}},\\
					&\|\phi\|_{T,*}\lesssim T^\frac{1}{12}\|\nabla\phi\|_{C^{m+\kappa}}\lesssim T^\frac{1}{12}\eps_1^{-(m+4)}\|u_0\|_{\dot W^{1,\infty}},\\
					&\|\phi\|_{T}\lesssim \|\nabla\phi\|_{L^\infty}+T^\frac{m+\kappa}{4}\|\nabla\phi\|_{\dot C^{m+\kappa}}\lesssim (1+T^\frac{m+\kappa}{4}\eps_1^{-(m+\kappa)})\|u_0\|_{\dot W^{1,\infty}}.
				\end{align*}
				We can take $T$ small enough such that $T^\frac{1}{100}\eps_1^{-(m+10)^2}(1+\|u_0\|_{\dot W^{1,\infty}})\leq 1$. Then one has $\|\phi\|_{T,*}\lesssim T^\frac{1}{20}$ and $\|\phi\|_{T}\lesssim \|u_0\|_{\dot W^{1,\infty}}$.
				Put these estimates into \eqref{wmre} and \eqref{hahawm}, we obtain 
				\begin{align}
					\|\tilde u\|_T&\leq \frac{5}{4}C_1 \varepsilon_0+C_1\|\tilde v\|_T(\|\tilde v\|_T+T^\frac{1}{20})(1+\|\tilde v\|_T+T^\frac{1}{20})^{m+5}+C_1T^\frac{1}{10}
					(1+\|\tilde u\|_T+\|u_0\|_{\dot W^{1,\infty}}),\label{eeee1}\\
					\|\mathbf{u}\|_T&\leq C_1(1+\|(\tilde v_1,\tilde v_2)\|_T+T^\frac{1}{20})^{m+5}(\|\tilde v_1,\tilde v_2\|_{T}+T^\frac{1}{20})\|\mathbf{v}\|_T+C_1T^\frac{1}{10}\|\mathbf{u}\|_T.\label{eee2}
				\end{align}
				Let $\varepsilon_0=\frac{1}{(10+C_1)^{m+10}}$, $T_1=\min\left\{T_0,\left(\frac{\varepsilon_0\eps_1}{10+C_1+\|u_0\|_{\dot W^{1,\infty}}}\right)^{(m+100)^2}\right\}$, $\sigma_2=2C_1\varepsilon_0$. If $v\in  \mathcal{X}_{T_1,u_0\ast \rho_{\eps_1}}^{\sigma_2}$, then one has $\|\tilde v\|_{T_1}\leq \sigma_2$, hence we obtain from \eqref{eeee1} that 
				\begin{align*}
					\|\tilde u\|_{T_1}\leq &\frac{3}{2}C_1\varepsilon_0+C_1T_1^\frac{1}{10}\|\tilde u\|_{T_1},
				\end{align*}
				which implies that 
				\begin{align*}
					\|\tilde u\|_{T_1}\leq\frac{\frac{5}{4}C_1\varepsilon_0}{1-C_1T_1^\frac{1}{10}}\leq \sigma_2.	
				\end{align*}
				Hence
				$\mathcal{S}v=u\in  \mathcal{X}_{T_1, u_0\ast \rho_{\eps_1}}^{\sigma_2}$.\\
				If $v_1,v_2\in  \mathcal{X}_{T_1, u_0\ast \rho_{\eps_1}}^{\sigma_2}$,  one has $\|(\tilde v_1,\tilde v_2)\|_T\leq 2\sigma_2$. Hence we obtain from \eqref{eee2} that 
				\begin{align*}
					\|\mathbf{u}\|_{T_1}\leq & 2^{m+7}C_1^2\varepsilon_0\|\mathbf{v}\|_{T_1}+C_1T_1^\frac{1}{10}\|\mathbf{u}\|_{T_1}.
				\end{align*}
				which leads to 
				\begin{align*}
					\|\mathbf{u}\|_{T_1}\leq &\frac{2^{m+7}C_1^2\varepsilon_0}{1-C_1T_1^\frac{1}{10}}\|\mathbf{v}\|_{T_1}\leq \frac{1}{2}\|\mathbf{v}\|_{T_1}.
				\end{align*}
				Hence $\mathcal{S}:\mathcal{X}_{T_1,u_0\ast \rho_{\eps_1}}^{\sigma_2}\to \mathcal{X}_{T_1, u_0\ast \rho_{\eps_1}}^{\sigma_2}$ is a contraction map. 
			\end{itemize}
			This completes the proof.
			
		\end{proof}
		\begin{lemma}\label{lemwi}
			Let $v_1,v_2:[0,T]\times\mathbb{R}^d\to \mathbb{R}$, and $\phi\in C^{m+4}$. Then for $N[v,\phi], \tilde F[v]$ defined in \eqref{defrhs}, and $\|\cdot\|_T, \|\cdot\|_{T,*}$ defined in \eqref{defnowi},there hold
			\begin{align}
				&\sup_{t\in[0,T]}t^{\frac{n+\kappa}{4}}\|(N[v_1,\phi]-N[v_2,\phi])(t)\|_{\dot C^{n+\kappa-3}}\lesssim(1+\|(v_1,v_2,\phi)\|_{T,*})^{m+5} \|v_1-v_2\|_T(\|v_1\|_{T,*}+\|v_2-\phi\|_T),\label{difN}\\
				&\sup_{t\in[0,T]}t^{\frac{n+\kappa}{4}}\|(\tilde F[v_1]-\tilde F[v_2])(t)\|_{\dot C^{n+\kappa-3}}\lesssim(1+\|(v_1,v_2)\|_{T,*})^{m+5}\|v_1-v_2\|_{T}\|(v_1,v_2)\|_{T,*}, \label{diftf}
			\end{align}
			for any $n\in\mathbb{N}$, $n\leq m$. 
		\end{lemma}
		\begin{proof}
			For simplicity, we 
			first prove \eqref{difN}. Note that 
			\begin{align*}
				N[v_1,\phi]-N[v_2,\phi]=(A[\nabla v_2]-A[\nabla v_1])\nabla^4v_1+(A[\nabla \phi]-A[\nabla v_2])\nabla^4(v_1-v_2).
			\end{align*}
			Applying Lemma \ref{lemcom}, one has for any $n\in\mathbb{N}$, $n\leq m$,
			\begin{align*}
				\|\nabla^n(A[\nabla v_2]-A[\nabla v_1])\|_{L^\infty}\lesssim &\sum_{l=0}^{n-1}\|\nabla(v_1-v_2)\|_{\dot C^l}(\|\nabla(v_1,v_2)\|_{\dot C^1}^{n-l}+\|\nabla(v_1,v_2)\|_{\dot C^{n-l}})\\
				&+\|\nabla (v_1-v_2)\|_{\dot C^n}\\
				\lesssim & t^{-\frac{n}{4}}\|v_1-v_2\|_T(1+\|(v_1,v_2)\|_{T,*})^n.
			\end{align*}
			Moreover, for $a\in(0,1)$,
			\begin{align*}
				&\|\nabla^n(A[\nabla v_2]-A[\nabla v_1])\|_{\dot C^a}\\
				&\quad\quad\lesssim \|\nabla(v_1-v_2)\|_{\dot C^{n+a}}+\sum_{l=0}^{n-1}\|\nabla(v_1-v_2)\|_{\dot C^{l+a}}(\|\nabla(v_1,v_2)\|_{\dot C^1}^{n-l}+\|\nabla(v_1,v_2)\|_{\dot C^{n-l}})\\
				&\quad\quad\quad\quad\quad+\sum_{l=0}^{n}\|\nabla(v_1-v_2)\|_{\dot C^{l}}\left(\|\nabla(v_1,v_2)\|_{\dot C^a}^{\frac{n-l+a}{a}}+\|\nabla(v_1,v_2)\|_{\dot C^{n-l+a}}\right)\\
				&\quad\quad\lesssim  t^{-\frac{n+a}{4}}\|v_1-v_2\|_T(1+\|(v_1,v_2)\|_{T,*})^\frac{n+a}{a}.
			\end{align*}
			Hence we obtain 
			\begin{align*}
				\|\nabla^n((A[\nabla v_2]-A[\nabla v_1])\nabla^4v_1)\|_{\dot C^{\kappa-3}}&\lesssim \sum_{k=0}^n\left(\|\nabla^k(A[\nabla v_2]-A[\nabla v_1])\|_{\dot C^{\kappa-3}}\|\nabla^{4+n-k}v_1\|_{L^\infty}\right.\\
				&\quad\quad +\left.\|\nabla^k(A[\nabla v_2]-A[\nabla v_1])\|_{L^\infty}\|\nabla^{4+n-k}v_1\|_{\dot C^{\kappa-3}}\right)\\
				&\lesssim t^{-\frac{n+\kappa}{4}}\|v_1-v_2\|_T\|v_1\|_{T,*}(1+\|(v_1,v_2)\|_{T,*})^{n+5}.
			\end{align*}
			Similarly, one has 
			\begin{align*}
				\|\nabla^n((A[\nabla \phi]-A[\nabla v_2])\nabla^4(v_1-v_2))\|_{\dot C^{\kappa-3}}&\lesssim \sum_{k=0}^n\left(\|\nabla^k(A[\nabla v_2]-A[\nabla \phi])\|_{\dot C^{\kappa-3}}\|\nabla^{4+n-k}(v_1-v_2)\|_{L^\infty}\right.\\
				&\quad\quad +\left.\|\nabla^k(A[\nabla v_2]-A[\nabla \phi])\|_{L^\infty}\|\nabla^{4+n-k}(v_1-v_2)\|_{\dot C^{\kappa-3}}\right)\\
				&\lesssim t^{-\frac{n+\kappa}{4}}\|v_2-\phi\|_T\|v_1-v_2\|_{T,*}(1+\|(v_2,\phi)\|_{T,*})^{n+5}.
			\end{align*}
			Then we get 
			\begin{align*}
				&\sup_{t\in[0,T]}t^{\frac{n+\kappa}{4}}\|(N[v_1,\phi]-N[v_2,\phi])(t)\|_{\dot C^{n+\kappa-3}}\lesssim \|v_1-v_2\|_T(\|v_1\|_{T,*}+\|v_2-\phi\|_T)(1+\|(v_1,v_2,\phi)\|_{T,*})^{n+5}.
			\end{align*}
			This completes the proof of \eqref{difN}. Then we consider \eqref{diftf}. Note that $\tilde F[v]$ only contains up to the third order derivatives of $v$. By Lemma \ref{maininterpo}, it suffices to prove that 
			\begin{equation}\label{hhdF}
				\begin{aligned}
					\sup_{t\in[0,T]}t^{\frac{k}{4}}&\|\nabla^k(\tilde F[v_1]-\tilde F[v_2])(t)\|_{L^\infty}\\
					&\quad\quad\quad\quad\lesssim\left(1+\|(v_1,v_2)\|_{T,*}\right)^{k+5}\|v_1-v_2\|_{T}\|(v_1,v_2)\|_{T,*}, \quad\quad\forall k\in\mathbb{N}, k\leq m+1. 
				\end{aligned}
			\end{equation}
			Indeed, by Lemma \ref{lemcom}, we have 
			\begin{align*}
				&\|\nabla^k(\tilde F[v_1]-\tilde F[v_2])(t)\|_{L^\infty}
				\lesssim\sum_{l=0}^{k+2}\|\nabla^{l+1}(v_1-v_2)\|_{L^\infty}\left(\|\nabla^2(v_1,v_2)\|_{L^\infty}^{k-l+4}+\|\nabla^{k-l+3} (v_1,v_2)\|_{L^\infty}^\frac{k-l+4}{k-l+2}\right).
			\end{align*}
			We remark that the right hand side contains up to $k+2$ order derivatives of $\nabla(v_1,v_2)$, and $k+2<m+\kappa$ for any $k\leq m+1$. Hence the right hand side terms can be controlled by $\|(v_1,v_2)\|_{T}$. Specially, terms with higher order derivatives can be controlled by the non-endpoint norm $\|(v_1,v_2)\|_{T,*}$. This implies \eqref{hhdF} and completes the proof of \eqref{diftf}.
			
		\end{proof}
		\subsection{Surface diffusion flow}\label{secdas}
		In this section we consider the surface diffusion equation for the graph of a function $\Gamma=\{(x,u(x))\}$, where $u(t,x):(0,+\infty)\times \mathbb{R}^d\to \mathbb{R}$. The equation can be written as
		\begin{align*}
			V=\Delta_\Gamma H,
		\end{align*}
		where the notation and meaning of $V=-\frac{u_t}{\langle \nabla u\rangle}$ and $H=\nabla\cdot(\frac{\nabla u}{\langle\nabla u\rangle})$ are just the same as what we denote in Willmore flow, and $\Delta_\Gamma H$ is Laplace-Beltrami operator of $\Gamma$ applied on $H$. We can see that the equation is equivalent to the following
		\begin{equation}
			\begin{aligned}\label{eqsurfdiff}
				&\partial_t u+\operatorname{div }\left((\langle \nabla u \rangle \mathrm{Id}-\frac{\nabla u\otimes \nabla u}{\langle \nabla u \rangle})\nabla H\right)=0,\ \ \text{on}\ (0,T)\times \mathbb{R}^d,\\
				&u|_{t=0}=u_0.
			\end{aligned}
		\end{equation}
		The readers can review \cite{DDE2005} to find further details of the derivation. We need to mention some works before. In \cite{Koch2011}, Koch and Lamm proved that for initial data $f_0\in C^{0,1}$ with $\|f_0\|_{C^{0,1}}<\eps$ for some small $\eps$, then there exists a unique solution $f$ in $X_\infty$ with $\|f\|_{X_\infty}\leq C\|f_0\|_{C^{0,1}}$. The definition of $X_\infty$ can be seen in \eqref{Xinfty}. For anisotropic cases, in \cite{DDE2005} and \cite{DDE2005b}, they used time discretization method, and proved error bounds by numerical scheme.\\
		Our main result is the following.
		\begin{theorem}\label{thmsdf} Fix $\kappa$ such that $0<4-\kappa\ll1$, and $m\in\mathbb{N}$. There exists $\varepsilon_0>0$ such that \\
			1)For any $u_0$ with $\|\nabla u_0\|_{L^\infty}\leq \varepsilon_0$, there exists a unique solution $u$ to \eqref{eqsurfdiff} such that
			\begin{align*}
				\sup_{t>0}	(\|\nabla u(t)\|_{L^\infty}+t^\frac{m+\kappa}{4}\|\nabla u(t)\|_{\dot C^{m+\kappa}})\leq C\|\nabla u_0\|_{L^\infty}.
			\end{align*}
			2) For any initial data $u_0\in \dot W^{1,\infty}$, if there exists $\eps_1>0$ such that 
			\begin{align}\label{sdfcon}
				\|\nabla u_0-\nabla u_0\ast \rho_{\eps_1}\|_{L^\infty}\leq \varepsilon_0,
			\end{align} then there exists $T>0$ such that \eqref{eqsurfdiff} admits a unique solution  $u$ satisfying
			\begin{align*}
				\sup_{t\in[0,T]}	(\|\nabla u(t)\|_{L^\infty}+t^\frac{m+\kappa}{4}\|\nabla u(t)\|_{\dot C^{m+\kappa}})\leq C \|\nabla u_0\|_{L^\infty}.
			\end{align*}
		\end{theorem}
		\begin{proof}
			By a direct calculation, we can write \eqref{eqsurfdiff} as 	\begin{align*}
				u_t+\tilde A^{ijkl}(u)\partial_{ijkl} u=f_0[u]+\operatorname {div} (f_1[u]),
			\end{align*} 
			where 
			$$
			\tilde A^{ijkl}(u)=\delta_{ij}\delta_{kl}-2\delta_{ij}\frac{\partial_k u\partial_l u}{\langle\nabla u\rangle^2}+\frac{\partial_i u\partial_{j} u\partial_k u\partial_l u}{\langle\nabla u\rangle^4}.
			$$
			It is easy to check that $A^{ijkl}(u)\xi_i\xi_j\xi_k\xi_l\geq \frac{1}{\langle\nabla u\rangle^4}|\xi|^4$ for any $\xi\in \mathbb{R}^n$. Then we can follow the proof for Theorem \ref{thmwill} line by line to obtain the results.
		\end{proof}
		
		\subsection{Diffusion of axisymmetric surfaces}
		We consider the axisymmetric surfaces diffusion equations (ASD) in this section. In this section, we recall the derivation of the equation as the way in \cite{LS2013}. To begin with, we consider the following axisymmetric surface
		\begin{equation*}
			\Gamma=\{(x,h(x)\nu)|x\in\mathbb{R},\nu\in\mathbb{S}^n\},
		\end{equation*}
		where $h:\mathbb{R}\rightarrow (0,\infty)$ is the profile function. We denote $\Delta_\Gamma$ to be the Laplace-Beltrami operator, which is the surface Laplace operator on the surface $\Gamma$. If we denote $\mathcal{H}$ to be the mean curvature vector of $\Gamma$, and $V$ be the velocity in the normal direction of $\Gamma$, then we have
		\begin{equation*}
			V(\Gamma(t))=\Delta_{\Gamma(t)}\mathcal{H}(\Gamma(t)).
		\end{equation*} 
		For axisymmetric case, the Laplace-Beltrami operator has the form 
		\begin{align*}
			&\Delta_{\Gamma(t)}=\frac{1}{h\sqrt{1+h_x^2}}\partial_x\left(\frac{1}{\sqrt{1+h_x^2}}\partial_x\right),\\
			&V(\Gamma(t))=\frac{h_t}{\sqrt{1+h_x^2}},
		\end{align*}
		and the mean curvature is given by 
		\begin{align*}
			\mathcal{H}(h)=\kappa_1+\kappa_2=\frac{1}{h(1+h_x^2)^{\frac{1}{2}}}-\frac{h_{xx}}{(1+h_x^2)^{\frac{3}{2}}},
		\end{align*}
		where $\kappa_1$ and $\kappa_2$ are azimuthal and axial principle curvatures respectively on $\Gamma$. Substituting the  formulas into above to obtain the equation
		\begin{equation}\label{eqdas}
			\begin{aligned}
				&\partial_th=\frac{1}{h}\left(\frac{h}{(1+h_x^2)^{\frac{1}{2}}}(\mathcal{H}(h))_x\right)_x,\quad\text{in}\ [0,T]\times \mathbb{R},\\
				&h|_{t=0}=h_0,\quad \text{in}\ \mathbb{R}.
			\end{aligned}
		\end{equation}
		We will make a short review about the work before. The equation was derived from \cite{Mul1957}, and more details about the derivation and the model can be seen in \cite{CT1996}. In \cite{EMS1997}, the authors proved the existence and uniqueness of local smooth solution, and the authors also proved that for any initial data close to sphere, the solution exists globally and converges exponentially in time to the sphere. Results about stability by numerical investigation can be seen on \cite{CFM1996}. In \cite{WVMD1998}, the authors use self-similar coordinates to study the pinchoff of a cylinder and get supportive result with \cite{CFM1996}. \cite{JBW1998} has very important contribution to the study of surface discussion. In \cite{JBW1998}, the authors studied linear stability under evolution via surface diffusion of axisymmetric equilibrium surfaces(the cylinder, the Delaunay unduloid, and the sphere) by analytic and numerical methods. Also, for self-similar solutions, \cite{JBW1998} considered the possible behavior of pinchoff and the relationship between pinchoff and stability under perturbations in axisymmetric motion by surface diffusion. \\
		
		We apply Theorem \ref{lemmain} to obtain the local well-posedness of \eqref{eqdas} in the critical space $W^{1,\infty}$. Moreover, due to the non-homogeneity nature of the equation, the presence of some low-order terms prevents us from obtaining global well-posedness. However, the global well-posedness is still available if we consider the equation in periodic domain, where the lower order terms are dominated by higher order terms.
		Our main theorem is as follows.
		\begin{theorem}\label{thmdaslo}
			Fix $\kappa$ such that $0<4-\kappa\ll 1$, and $m\in\mathbb{N}$. 	 \\
			\textit{i)} For any initial data $h_0\in \dot W^{1,\infty}(\mathbb{R})$  satisfying $\inf_{x\in\mathbb{R}}|h_0(x)|\geq\lambda>0$, there exists $\varepsilon_0>0$ such that if
			\begin{equation}\label{condas}
				\|h_0-h_0\ast \rho_{\eps_1}\|_{ W^{1,\infty}}\leq \varepsilon_0,
			\end{equation}
			for some $\eps_1>0$,
			then there exists $T>0$ such that \eqref{eqdas} admits a unique solution $h$ in $[0,T]$ satisfying
			\begin{equation*}
				\sup_{t\in [0,T]}(\|\partial_xh(t)\|_{L^\infty}+t^{\frac{m+\kappa}{4}}\|\partial_x h(t)\|_{\dot C^{m+\kappa}})\lesssim \|h_0\|_{\dot W^{1,\infty}}.
			\end{equation*}
			\textit{ii)} Consider \eqref{eqdas} in the periodic setting. For any initial data $h_0\in \dot W^{1,\infty}(\mathbb{S})$, if 
			\begin{align}\label{conglo}
				\bar h_0:=\frac{1}{2\pi}\int_{-\pi}^\pi h_0(x)dx>1+\delta_0.
			\end{align}
			Then there exists $0<\varepsilon_0=\varepsilon_0(\delta_0)\ll 1$ such that if $	\|\partial_x h_0\|_{L^\infty(\mathbb{S})}\leq \varepsilon_0$, then \eqref{eqdas} admits a constant $c_0>0$ and a unique global solution $h$ satisfying 
			\begin{align*}
				\sup_{t\geq 0}(e^{c_0t}\|\partial_x h(t)\|_{L^\infty}+t^\frac{m+\kappa }{4}e^{c_0t}\|\partial_x h(t)\|_{\dot C^{m+\kappa }}+\|h(t)-\bar h_0\|_{L^\infty})\leq \frac{\bar h_0}{2}.
			\end{align*}
		\end{theorem}
		\begin{remark}
			The proof of ii) cannot be extended to the whole line $\mathbb{R}$ even if $\frac{1}{\bar {h}_0}\ll 1$.
		\end{remark}
		We first prove the local well-posedness result in the whole space.\\
		\begin{proof}[Proof of Theorem \ref{thmdaslo} \textit{i)}]
			We rewrite \eqref{eqdas} as 
			\begin{align*}
				&\partial_t h+\frac{\partial_x
					^4h}{\langle h_x\rangle^4}=F[h],\quad\text{in}\ [0,T]\times\mathbb{R},\\
				&h|_{t=0}=h_0,\quad\text{in}\ \mathbb{R}.
			\end{align*}
			Here the nonlinear remainder term can be written as 
			\begin{equation}\label{nff}
				\begin{aligned}
					F[h]&=\left(\frac{(\mathcal{H}(h))_x}{\langle h_x\rangle} \right)_x+\frac{\partial_x^4h}{\langle h_x\rangle^4}+\frac{1}{h}\frac{h_x(\mathcal{H}(h))_x}{\langle h_x\rangle}\\
					&=G_1(h_x)\partial_{x}^3h\partial_{x}^2h+G_2(h_x)(\partial_{x}^2h)^3+\sum\frac{1}{h^r}G_{r,n_1,n_2}(h_x)(\partial_{x}^3h)^{n_1}(\partial_{x}^2h)^{n_2},
				\end{aligned}
			\end{equation}
			where the sum of the second line is taken along $r=1,2,3$, $(n_1,n_2)=(1,0),(0,i),i=0,1,2$, and $G_i,G_{r,n_1,n_2}'$ bounded smooth functions satisfying
			\begin{equation}\label{dascdcoe}
				\|\nabla^jG_1(\cdot)\|_{L^\infty}+\|\nabla^jG_2(\cdot)\|_{L^\infty}+\|\nabla^jG_{r,n_1,n_2}(\cdot)\|_{L^\infty}\leq C_j,\quad j=1,\cdots,m.
			\end{equation}
			For simplicity, fix $m\in\mathbb{N}$ and take $0<4-\kappa\ll 1$. We denote
			\begin{equation}\label{normdas}
				\begin{aligned}
					&\|f\|_{T}:=\sup_{t\in[0,T]}(\|f(t)\|_{L^\infty}+t^{\frac{m+\kappa }{4}}\|f(t)\|_{\dot C^{m+\kappa }}),\\
					&\|f\|_{T,*}:=\sup_{t\in[0,T]}(t^{\frac{\kappa-3}{4}}\|f(t)\|_{\dot C^{\kappa-3}}+t^{\frac{m+\kappa }{4}}\|f(t)\|_{\dot C^{m+\kappa }}),\\
					&\|f\|_{X_T}:=\sup_{t\in[0,T]}\|f(t)\|_{L^\infty}+\|\partial_xf\|_{T}.
				\end{aligned}
			\end{equation}
			By \eqref{condas}, take $\phi=h_0\ast \rho_{\eps_1}$, one has
			\begin{align}\label{inh-phi}
				\|h_0-\phi\|_{W^{1,\infty}}\leq \varepsilon_0. 
			\end{align}
			Moreover, by taking $\varepsilon_0\leq \frac{\lambda}{100}$, one has 
			\begin{align}\label{infphi}
				\inf_x|\phi(x)|\geq \inf_x|h_0(x)|-\|h_0-\phi\|_{L^\infty}\geq\frac{3\lambda}{4}.
			\end{align}
			Let $\sigma,T>0$ that will be fixed later, we define a space $\mathcal{X}_{T,\phi}^{\sigma}$ as follows 
			\begin{equation*}
				\mathcal{X}_{T,\phi}^{\sigma}=\left\{f\in L_T^\infty L^\infty:\inf_{t\in[0,T]}\inf_x|f(t,x)|\geq \frac{\lambda}{2},\|f-\phi\|_{X_T}\leq\sigma\right\}.
			\end{equation*}
			For $f\in  \mathcal{X}_{T,\phi}^{\sigma}$, we define a map $\mathcal{S}f=h$, such that $h$ solves the Cauchy problem
			\begin{equation}\label{daseq}
				\begin{aligned}
					& \partial_th+\frac{\partial_x^4h}{\langle\phi_x\rangle^4}=F[f]-\partial_x^4f\left(\frac{1}{\langle f_x\rangle^4}-\frac{1}{\langle\phi_x\rangle^4}\right),\\
					& h|_{t=0}=h_0.
				\end{aligned}
			\end{equation}
			Denote $\tilde f=f-\phi$, $\tilde  h=h-\phi$, then $\tilde { h}_x=\partial_x (h-\phi)$ satisfies
			\begin{equation}\label{dasdeq}
				\begin{aligned}
					&\partial_t \tilde {h}_x+\frac{\partial_x
						^4\tilde h_x}{\langle \phi_x\rangle^4}=\partial_x(F[f]+N[f,\phi])+R[h,\phi],\\
					&\tilde h_x|_{t=0}=\partial_x(h_0-\phi),
				\end{aligned}	
			\end{equation}
			where
			\begin{align}
				&N[f,\phi]=-\partial_x^4f\left(\frac{1}{\langle f_x\rangle^4}-\frac{1}{\langle\phi_x\rangle^4}\right),\label{nfp}\\
				&R[h,\phi]=-\partial_x\left(\frac{\partial_x^4\phi}{\langle\phi_x\rangle^4}\right)+\partial_x^4\tilde h_x\partial_x\left(\frac{1}{\langle\phi_x\rangle^4}\right).
			\end{align}
			Apply Theorem \ref{lemmain} for $F=F[f]+N[f,\phi]$ and $G=R[h,\phi]$, there exists time $T=T(\|\phi_x\|_{C^{m+4}})>0$ such that 
			\begin{equation}\label{dasmain}
				\begin{aligned}
					&\|\tilde h_x\|_{T}\lesssim \|h_0-\phi\|_{W^{1,\infty}}+\sup_{t\in[0,T]}(t^{\frac{\kappa}{4}}\|F[f](t)\|_{\dot C^{\kappa-3}}+t^{\frac{m+\kappa}{4}}\|F[f](t)\|_{\dot C^{m+\kappa-3}})\\
					&\quad\quad+\sup_{t\in[0,T]}(t^{\frac{\kappa}{4}}\|N[f,\phi](t)\|_{\dot C^{\kappa-3}}+t^{\frac{m+\kappa}{4}}\|N[f,\phi](t)\|_{\dot C^{m+\kappa-3}}+t^{\frac{m}{4}+1}\|R[h,\phi](t)\|_{\dot C^m})+\|R[h,\phi]\|_{L^1_TL^\infty}.
				\end{aligned}
			\end{equation}
			Observe that $N[\phi,\phi]=0$,  then take $f_1=f$, $f_2=\phi$ in Lemma \ref{dasnl}, we have
			\begin{equation}\label{dasN1}
				\begin{aligned}
					&\sup_{t\in[0,T]}\left(t^{\frac{\kappa}{4}}\|N[f,\phi](t)\|_{\dot C^{\kappa-3}}+t^{\frac{m+\kappa}{4}}\|N[f,\phi](t)\|_{\dot C^{m+\kappa-3}}\right)\\
					&\quad\quad\quad\quad  \lesssim \|\tilde f_x\|_{T}(\|\tilde f_x\|_{T}+\|\phi_x\|_{T,*})(1+\|\tilde f_x\|_{T}+\|\phi_x\|_T)^{m+10},
				\end{aligned}
			\end{equation}
			and
			\begin{equation*}
				\begin{aligned}
					&\sup_{t\in[0,T]}\left(t^{\frac{\kappa}{4}}\|(F[f]-F[\phi])(t)\|_{\dot C^{\kappa-3}}+t^{\frac{m+\kappa}{4}}\|(F[f]-F[\phi])(t)\|_{\dot C^{m+\kappa-3}}\right)\\
					&\quad\quad\quad\quad \lesssim \|\tilde f\|_{X_T}\left(\|(f,\phi)_{x}\|_{T,*}+T^{\frac{1}{5}}\right)(1+\|\tilde f_x\|_T+\|\phi_x\|_T)^{m+10}.
				\end{aligned}
			\end{equation*}
			Combining this with the fact that 
			\begin{equation*}
				\sup_{t\in[0,T]}\left(t^{\frac{\kappa}{4}}\|F[\phi]\|_{\dot C^{\kappa-3}}+t^{\frac{m+\kappa}{4}}\|F[\phi]\|_{\dot C^{m+\kappa-3}}\right)\lesssim T^{\frac{\kappa-3}{4}}(1+\|\phi\|_{C^{m+5}})^{m+2},
			\end{equation*}
			we can deduce that
			\begin{equation}\label{dasF1}
				\begin{aligned}
					&\sup_{t\in[0,T]}\left(t^{\frac{\kappa}{4}}\|F[f](t)\|_{\dot C^{\kappa-3}}+t^{\frac{m+\kappa}{4}}\|F[f](t)\|_{\dot C^{m+\kappa-3}}\right)\\
					&\quad\quad\lesssim \|\tilde f\|_{X_T}\left(\|(f,\phi)_{x}\|_{T,*}+T^{\frac{1}{5}}\right)(1+\|\tilde f_x\|_T+\|\phi_x\|_T)^{m+10}+T^{\frac{1}{5}}(1+\|\phi\|_{C^{m+5}})^{m+2}.
				\end{aligned}
			\end{equation}
			By the definition $\phi=h_0\ast \rho_{\eps_1}$, we have 
			\begin{equation}\label{esph}
				\begin{aligned}
					&\|\phi_x\|_{T,*}\lesssim T^\frac{1}{5}\|\phi\|_{C^{m+5}}\lesssim T^\frac{1}{5}\eps_1^{-(m+4)}\|h_0\|_{\dot W^{1,\infty}}\lesssim T^\frac{1}{10},\\
					&\|\phi_x\|_T\lesssim \|h_0\|_{\dot W^{1,\infty}}, \quad\quad\quad\quad T^{\frac{1}{5}}(1+\|\phi\|_{C^{m+5}})^{m+2}\lesssim T^\frac{1}{10},
				\end{aligned}
			\end{equation}
			provided $T<\eps_1^{10(m+10)^2}(1+\|h_0\|_{\dot W^{1,\infty}})^{-(m+10)}$.
			Combining this with \eqref{dasF1} and \eqref{dasN1}, we obtain 
			\begin{equation}\label{dasnf}
				\begin{aligned}
					& \sup_{t\in[0,T]}\left(t^{\frac{\kappa}{4}}\|F[f](t)\|_{\dot C^{\kappa-3}}+t^{\frac{m+\kappa}{4}}\|F[f](t)\|_{\dot C^{m+\kappa-3}}\right)\\
					&\quad\quad+\sup_{t\in[0,T]}\left(t^{\frac{\kappa}{4}}\|N[f,\phi](t)\|_{\dot C^{\kappa-3}}+t^{\frac{m+\kappa}{4}}\|N[f,\phi](t)\|_{\dot C^{m+\kappa-3}}\right)\\
					&\lesssim \|\tilde f\|_{X_T}(\|\tilde f_x\|_{T}+T^\frac{1}{10})(1+\|\tilde f_x\|_{T}+\|h_0\|_{\dot W^{1,\infty}})^{m+20}+T^\frac{1}{10}.
				\end{aligned}
			\end{equation}
			Finally, for $R[h,\phi]$,
			\begin{equation*}
				\begin{aligned}
					&\|R[h,\phi](t)\|_{L^\infty}\lesssim t^{-\frac{3}{4}}\|\tilde h_x\|_{T,*}\|\phi\|_{C^5}+\|\phi\|_{C^5}(1+\|\phi\|_{C^5}),\\
					&\|\partial_x^{m}R[h,\phi](t)\|_{L^\infty}\lesssim t^{-\frac{m+3}{4}}\|\tilde h_x\|_{T,*}\|\phi\|_{C^{m+5}}(1+\|\phi\|_{C^{m+5}})^{m+4}+\|\phi\|_{C^{m+5}}(1+\|\phi\|_{C^{m+5}})^{m+5},
				\end{aligned}
			\end{equation*}
			which infer
			\begin{equation}\label{dasR}
				\begin{aligned}
					\|R[h,\phi]\|_{L^1_TL^\infty}+\sup_{t\in[0,T]}t^{\frac{m}{4}+1}\|\partial_x^{m}R[h,\phi](t)\|_{L^\infty}&\lesssim T^{\frac{1}{4}}(\|\tilde h_x\|_{T}+1)\|\phi\|_{C^{m+5}}(1+\|\phi\|_{C^{m+5}})^{m+5}\\
					&\lesssim T^\frac{1}{10}(\|\tilde h_x\|_{T}+1).
				\end{aligned}
			\end{equation}
			We collect \eqref{inh-phi}, \eqref{dasmain}, \eqref{dasnf},  and \eqref{dasR} to get 
			\begin{align*}
				\|\tilde h_x\|_T\leq C_0\left( \frac{3\varepsilon_0}{2}+\sigma(\sigma+T^\frac{1}{10})(1+\sigma+\|h_0\|_{\dot W^{1,\infty}})^{m+20}+T^\frac{1}{10}(1+\|\tilde h_x\|_T)\right).
			\end{align*}
			Here the constant $C_0$ depends only on $m,\kappa$ and $\lambda$.
			Take 
			\begin{equation}\label{dascoe}
				\eps_0=\frac{1}{(10+C_0+\|h_0\|_{\dot W^{1,\infty}})^{(m+10)^2}},\quad \sigma=10C_0\eps_0,\quad T=\frac{\eps_1^{10(m+10)^2}}{(10+C_0+\|h_0\|_{\dot W^{1,\infty}})^{m+10}},
			\end{equation}
			then we obtain
			\begin{equation*}
				\|\tilde h_x\|_T\leq\frac{\sigma}{2}.
			\end{equation*}
			Finally, by \eqref{daseq}, we have
			\begin{equation*}
				\|\tilde h\|_{L_T^\infty L^\infty}\leq \|h_0-\phi\|_{L^\infty}+\int_0^T\|\partial_t\tilde h(\tau)\|_{L^\infty}d\tau\leq \frac{\sigma}{4}\leq \frac{\lambda}{40},
			\end{equation*}
			which infers 
			\begin{equation}\label{dasinfh}
				\inf_{t,x}|h|\geq\frac{\lambda}{2},\quad\quad\quad   \|\tilde h\|_{X_T}\leq\frac{3\sigma}{4}.
			\end{equation}
			Hence $h=\mathcal{S}f\in\mathcal{X}_{T,\phi}^{\sigma}$.\\
			Now we only need to prove the contraction property. For $f_1,f_2\in\mathcal{X}_{T,\phi}^{\sigma}$, denote $\mathbf{f}=f_1-f_2$, $\mathbf{h}=h_1-h_2$. We have
			\begin{equation}\label{dasdeq1}
				\begin{aligned}
					&\partial_t \mathbf{h}-\frac{\partial_x^4\mathbf{h}}{\langle\phi_x\rangle^4}=F[f_1]-F[f_2]-(N[f_1,\phi]-N[f_1,\phi]),\quad\text{in}\ [0,T]\times\mathbb{R},\\
					&\mathbf{h}|_{t=0}=0,\quad\text{in}\ \mathbb{R}.
				\end{aligned}
			\end{equation}
			Taking one spatial derivative to obtain 
			\begin{equation}\label{dasdeq}
				\begin{aligned}
					&\partial_t\mathbf{h}_x-\frac{\partial_x^4\mathbf{h}_x}{\langle\phi_x\rangle^4}=\partial_x(N[f_1,\phi]-N[f_2,\phi]+F[f_1]-F[f_2])+R[h_1,\phi]-R[h_2,\phi],\quad\text{in}\ [0,T]\times\mathbb{R},\\
					&\mathbf{h}_x|_{t=0}=0,\quad\text{in}\ \mathbb{R}.
				\end{aligned}
			\end{equation}
			Then by Theorem \ref{lemmain}, there exists $T=T(\|\phi\|_{C^{m+4}})$ such that 
			\begin{equation}\label{dasd}
				\begin{aligned}
					\|\mathbf{h}_x\|_{T}\lesssim  &\sup_{t\in[0,T]}\left(t^{\frac{\kappa}{4}}\|(N[f_1,\phi]-N[f_2,\phi])(t)\|_{\dot C^{\kappa-3}}+t^{\frac{m+\kappa}{4}}\|(N[f_1,\phi]-N[f_2,\phi])(t)\|_{\dot C^{m+\kappa-3}}\right)\\
					&+\sup_{t\in[0,T]}\left(t^{\frac{\kappa}{4}}\|(F[f_1]-F[f_2])(t)\|_{\dot C^{\kappa}}+t^{\frac{m+\kappa}{4}}\|(F[f_1]-F[f_2])(t)\|_{\dot C^{m+\kappa}}\right)\\
					&+\|R[h_1,\phi]-R[h_2,\phi]\|_{L^1_TL^\infty}+\sup_{t\in[0,T]}t^{\frac{m}{4}+1}\|(R[h_1,\phi]-R[h_2,\phi])(t)\|_{\dot C^m}.
				\end{aligned}
			\end{equation}
			By Lemma \ref{dasnl} and \eqref{esph},  we have
			\begin{equation}\label{dasdN}
				\begin{aligned}
					&\sup_{t\in[0,T]}\left(t^{\frac{\kappa}{4}}\|(N[f_1,\phi]-N[f_2,\phi])(t)\|_{\dot C^{\kappa-3}}+t^{\frac{m+\kappa}{4}}\|(N[f_1,\phi]-N[f_2,\phi])(t)\|_{\dot C^{m+\kappa-3}}\right)\\
					&\quad+\sup_{t\in[0,T]}\left(t^{\frac{\kappa}{4}}\|(F[f_1]-F[f_2])(t)\|_{\dot C^{\kappa-3}}+t^{\frac{m+\kappa}{4}}\|(F[f_1]-F[f_2])(t)\|_{\dot C^{m+\kappa-3}}\right)\\
					&\quad\quad\quad\lesssim \|\mathbf{f}\|_{X_T}(\|(\tilde f_1,\tilde f_2)_x\|_{T}+T^{\frac{1}{10}})(1+\| (\tilde f_1,\tilde f_2)_x\|_{T}+\|h_0\|_{\dot W^{1,\infty}})^{m+10},
				\end{aligned}
			\end{equation}
			where $\tilde f_i=f_i-\phi$.
			For $R[h,\phi]$, one has
			\begin{equation*}
				\begin{aligned}
					&\|(R[h_1,\phi]-R[h_2,\phi])(t)\|_{L^\infty}\lesssim t^{-\frac{3}{4}}\|\mathbf{h}_x\|_{T},\\
					&\|(R[h_1,\phi]-R[h_2,\phi])(t)\|_{\dot C^m}\lesssim t^{-\frac{m+3}{4}}\|\mathbf{h}_x\|_{T}(1+\|\phi\|_{C^{m+5}})^{m+5},
				\end{aligned}
			\end{equation*}
			which yield that 
			\begin{equation}\label{dasdR}
				\begin{aligned}
					\|R[h_1,\phi]-R[h_2,\phi]\|_{L^1_TL^\infty}&+\sup_{\in[0,T]}t^{\frac{m}{4}+1}\|(R[h_1,\phi]-R[h_2,\phi])(t)\|_{\dot C^m}\\
					&\lesssim T^\frac{1}{4} \|\mathbf{h}_x\|_{T}(1+\|\phi\|_{C^{m+5}})^{m+5}\lesssim T^\frac{1}{10} \|\mathbf{h}_x\|_{T}.
				\end{aligned}
			\end{equation}
			Combining \eqref{dasd}, \eqref{dasdN}, and \eqref{dasdR}, since $f_1,f_2\in\mathcal{S}_{T,\phi}^\sigma$, one has
			\begin{equation*}
				\|\mathbf{h}_x\|_T\leq C(1+T)^{\frac{m}{4}+1}(1+\sigma+\|\phi\|_{C^{m+5}})(\sigma+T^{\frac{1}{10}})(\|\mathbf{h}_x\|_T+\|\mathbf{f}_x\|_T).
			\end{equation*}
			By taking T satisfying \eqref{dascoe}, we will have
			\begin{equation}\label{dasdT}
				\|\mathbf{h}_x\|_T\leq\frac{1}{4}\|\mathbf{f}\|_{X_T}.
			\end{equation}
			By \eqref{dasdeq1}, if we take T the same as above, then we have 
			\begin{equation*}
				\|\mathbf{h}\|_{L^\infty_TL^\infty}\leq\int_0^T\|\partial_t\mathbf{h}(t)\|_{L^\infty}dt\leq \frac{1}{4}\|\mathbf{f}\|_{X_T},
			\end{equation*}
			which combines with \eqref{dasdT} gives
			\begin{equation}
				\|\mathbf{h}\|_{X_T}\leq\frac{1}{2}\|\mathbf{f}\|_{X_T}.
			\end{equation}
			Hence $\mathcal{S}$ is a contraction map from $\mathcal{X}^\sigma_{T,\phi}$ to itself. This completes the proof.
		\end{proof}
		\begin{lemma}\label{dasnl}
			Let $T\in (0,1)$,  $f_1,f_2\in \mathcal{X}_{T,\phi}^{\sigma}$ and $\phi\in C^{m+5}$, denote $f_{12}=f_1-f_2$. Then for $N[f_,\phi], F[f]$ defined in \eqref{nfp}, \eqref{nff}, and norms  $\|\cdot\|_{T}, \|\cdot\|_{T,*}, \|\cdot\|_{X_T}$ defined in \eqref{normdas}, there hold
			\begin{align}
				& \sup_{t\in[0,T]}t^{\frac{n+\kappa}{4}}  \|(N[f_1,\phi]-N[f_2,\phi])(t)\|_{\dot C^{n+\kappa-3}}\nonumber\\
				&\quad\quad\quad\quad\quad\quad\lesssim \|\partial_xf_{12}\|_{T}(\|(f_1-\phi)_x\|_{T}+\|f_{2x}\|_{T,*})(1+\|(f_1,f_2,\phi)_x\|_T)^{n+10},\label{dasn}\\
				& \sup_{t\in[0,T]}t^{\frac{n+\kappa}{4}} \|(F[f_1]-F[f_2])(t)\|_{\dot C^{n+\kappa-3}}\nonumber\\
				&\quad\quad\quad\quad\quad\quad\lesssim \|f_{12}\|_{X_T}(\|(f_1, f_2)_x\|_{T,*}+T^\frac{1}{5})(1+\|(f_1,f_2)_x\|_{T})^{n+10},\label{dasf}
			\end{align}
			for any $n\in\mathbb N$, $n\leq m$.
		\end{lemma}
		\begin{proof}
			Note that 
			\begin{align*}
				N[f_1,\phi]-N[f_2,\phi]=\partial_x^{4}f_{12}\left(\frac{1}{\langle f_{1x}\rangle^4}-\frac{1}{\langle \phi_{x}\rangle^4}\right)+\partial_x^4f_2\left(\frac{1}{\langle f_{1x}\rangle^4}-\frac{1}{\langle f_{2x}\rangle^4}\right).
			\end{align*}
			For simplicity, we will drop time variable $t$ in the following proof. Applying Lemma \ref{lemcom}, one has
			\begin{align*}
				\left\|\partial_x^k\left(\frac{1}{\langle f_{1x}\rangle^4}-\frac{1}{\langle f_{2x}\rangle^4}\right)\right\|_{L^\infty}&\lesssim  \sum_{j=0}^{k}\|\partial_xf_{12}\|_{\dot C^{j}}(\|(f_1,f_2)_x\|_{\dot C^{1}}^{n-j}+\|(f_1,f_2)_x\|_{\dot C^{k-j}})\\
				&\lesssim t^{-\frac{k}{4}}\|\partial_xf_{12}\|_{T}(1+\|(f_1,f_2)_x\|_T)^k,\ \ \ \ \ \ \forall k\in\mathbb{N}, k\leq m+1.
			\end{align*}
			By interpolation, for $\gamma=\kappa-3\in(0,1)$, one has
			\begin{align*}
				\left\|\partial_x^n\left(\frac{1}{\langle f_{1x}\rangle^4}-\frac{1}{\langle f_{2x}\rangle^4}\right)\right\|_{\dot C^\gamma}
				\lesssim &t^{-\frac{n+\gamma}{4}}\|\partial_xf_{12}\|_{T}(1+\|(f_1,f_2)_x\|_T)^{n+1},\ \ \ \forall n\in\mathbb{N},n\leq m.
			\end{align*}
			Hence we obtain 
			\begin{align*}
				&\|\partial_x^n(N[f_1,\phi]-N[f_2,\phi])\|_{\dot C^\gamma}\\
				&\lesssim\sum_{j=0}^n\left( \|\partial_xf_{12}\|_{\dot C^{j+3+\gamma}}\left\|\partial_x^{n-j}\left(\frac{1}{\langle f_{1x}\rangle^4}-\frac{1}{\langle \phi_{x}\rangle^4}\right)\right\|_{L^\infty}+\|\partial_xf_{12}\|_{\dot C^{j+3}}\left\|\partial_x^{n-j}\left(\frac{1}{\langle f_{1x}\rangle^4}-\frac{1}{\langle \phi_{x}\rangle^4}\right)\right\|_{\dot C^\gamma}\right.\\
				&\quad\left.+\|f_{2x}\|_{\dot C^{j+3+\gamma}}\left\|\partial_x^{n-j}\left(\frac{1}{\langle f_{1x}\rangle^4}-\frac{1}{\langle f_{2x}\rangle^4}\right)\right\|_{L^\infty}+\|f_{2x}\|_{\dot C^{j+3}}\left\|\partial_x^{n-j}\left(\frac{1}{\langle f_{1x}\rangle^4}-\frac{1}{\langle f_{2x}\rangle^4}\right)\right\|_{\dot C^\gamma}\right)\\
				&\lesssim t^{-\frac{n+\gamma}{4}}\|\partial_xf_{12}\|_{T}(\|(f_1-\phi)_x\|_{T}+\|f_{2x}\|_{T,*})(1+\|(f_1,f_2,\phi)_x\|_T)^{n+\gamma}.
			\end{align*}
			This completes the proof of \eqref{dasn}.
			For $F[f]$, it suffices to estimate
			\begin{equation*}
				\begin{aligned}
					\|\partial_x^k(F[f_1]-F[f_2])\|_{L^\infty},\quad \forall k\in\mathbb{N},\ k\leq m+1.
				\end{aligned}
			\end{equation*}
			In fact, by Lemma \ref{maininterpo}, one has 
			\begin{equation}\label{dasint}
				\begin{aligned}
					\|\partial_x^k(F[f_1]-F[f_2])\|_{\dot C^\gamma}\lesssim \|\nabla^k(F[f_1]-F[f_2])\|_{L^\infty}^{1-\gamma}\|\nabla^{k+1}(F[f_1]-F[f_2])\|_{L^\infty}^\gamma.
				\end{aligned}
			\end{equation}
			By direct computation, it is easy to find that $F[f]$ can be divide into two parts, the first part $F_1[f]$ includes terms which have the same scaling with $\partial_x^4 f$, while another part $F_2[f]$ includes lower order terms. More precisely, we can write 
			\begin{align*}
				&F_1[f]=G_1(h_x)\partial_{x}^3h\partial_{x}^2h+G_2(h_x)(\partial_{x}^2h)^3,\\
				&F_2[f]=\sum\frac{1}{h^r}G_{r,n_1,n_2}(h_x)(\partial_{x}^3h)^{n_1}(\partial_{x}^2h)^{n_2},
			\end{align*}
			where the sum of the second line is taken along $r=1,2,3$, and $(n_1,n_2)=(1,0),(0,i),i=0,1,2$, and $G_i,G_{r,n_1,n_2}$ satisfying \eqref{dascdcoe}. By Lemma \ref{lemcom}, we obtain
			\begin{equation}\label{dasf1d}
				\|\partial_x^k(F_1[f_1]-F_1[f_2])(t)\|_{L^\infty}\lesssim t^{-\frac{k+3}{4}}\|f_{12}\|_{X_T}\|(f_1, f_2)_x\|_{T,*}(1+\|(f_1,f_2)_x\|_{T})^{k+10},\ \ \ \forall k\in\mathbb{N}, k\leq m+1.
			\end{equation}
			For lower order term $F_2[f]$, by \eqref{dasinfh}, we can prove that
			\begin{equation}\label{dasf2d}
				\|\partial_x^k(F_2[f_1]-F_2[f_2])(t)\|_{L^\infty}\lesssim t^{-\frac{k+2}{4}}\|f_{12}\|_{X_T}(1+\|(f_1,f_2)_x\|_{T})^{k+10},\ \ \ \forall k\in\mathbb{N}, k\leq m+1.
			\end{equation}
			Combining this with \eqref{dasf1d} to obtain that for any $k\in\mathbb{N}, k\leq m+1$,
			\begin{align*}
				\sup_{t\in[0,T]}   t^{\frac{k+3}{4}}\|\partial_x^k(F[f_1]-F[f_2])(t)\|_{L^\infty}\lesssim \|f_{12}\|_{X_T}(\|(f_1, f_2)_x\|_{T,*}+T^\frac{1}{4})(1+\|(f_1,f_2)_x\|_{T})^{k+10}.
			\end{align*}
			By \eqref{dasint}, we obtain \eqref{dasf}.
			This completes the proof.
		\end{proof}

		\vspace{0.3cm}
		Before starting the proof of the global result, we first estimate the parabolic kernel associate with \eqref{eqdas} in the periodic setting. We study the following model:
		\begin{equation}\label{model}
			\begin{aligned}
				&\partial_t h+{\partial_x^4 h}+\frac{\partial_x ^2 h}{\bar h_0^2}
				=F,\quad\quad\text{in}\ (0,T)\times \mathbb{S},\\
				&h|_{t=0}= h_0.
			\end{aligned}
		\end{equation}
		Here $\mathbb{S}=[-\pi,\pi)$.
		We write \eqref{model} in terms of Fourier series. Let $h(t,x)=\sum_{k=-\infty}^\infty h_k(t)e^{ikx}$, $F(t,x)=\sum_{k=-\infty}^\infty F_k(t)e^{ikx}$. Then for each $k\in \mathbb{Z}$, we get the ODE 
		\begin{align*}
			&\frac{d}{dt} h_k(t)+\left(|k|^4-\frac{|k|^2}{\bar h_0^2}\right)h_k(t)=F_k(t),\\
			&h_k|_{t=0}=h_{0,k}.
		\end{align*}
		Denote $A(k)=|k|^4-\frac{|k|^2}{\bar h_0^2}$. Then 
		\begin{align*}
			h_k=h_{0,k}e^{-A(k)t}+\int_0^tF_k(\tau)e^{-A(k)(t-\tau)}d\tau.
		\end{align*}
		Define  $	h_{\neq 0}(t,x)=\sum_{n\neq 0}h_ne^{inx}$. We can write 
		\begin{align}\label{defh}
			h_{\neq 0}(t,x)=\int_{-\pi}^\pi K_{\neq 0}(t,x-y)h_0(y)dy+\int_0^t \int_{-\pi}^\pi K_{\neq 0}(t-\tau,x-y)F[f](\tau,y)dyd\tau,
		\end{align}
		where \begin{align}\label{defKn0}
			K_{\neq 0}(t,x)=\frac{1}{2\pi}\sum_{n\in\mathbb{N}\backslash\{0\}}e^{-A(n)t}e^{inx}.
		\end{align}
		
		\begin{lemma}\label{lemkerper}
			For the kernel $K_{\neq 0}$ defined in \eqref{defKn0}, there holds
			\begin{align}\label{lemrK}
				\|\partial_x^lK_{\neq 0}(t)\|_{L^1(\mathbb{S})}\lesssim e^{-c_0t}t^{-\frac{l}{4}},
			\end{align}
			where $c_0=\frac{1}{4}\left(1-\frac{1}{\bar h_0^2}\right)$.
		\end{lemma}
		\begin{proof}
			We can write 
			\begin{align*}
				K_{\neq 0}(t,x)=\frac{e^{-c_0t}}{2\pi}\sum_{n\in\mathbb{N}_+}e^{-\tilde A(n)t}e^{inx},\ \ \ \text{where}\ \tilde A(n)=A(n)-c_0,\ \ \ c_0=\frac{1}{4}\left(1-\frac{1}{\bar h_0^2}\right).
			\end{align*}
			It is easy to check that 
			\begin{align*}
				\tilde A(n)\geq \frac{1}{4}\left(1-\frac{1}{\bar h_0^2}\right),\ \ \forall n\in\mathbb{N}_+.
			\end{align*}
			We claim that for $\tilde K(t,x)=\frac{1}{2\pi}\sum_{n\in\mathbb{N}_+}e^{-\tilde A(n)t}e^{inx}$, there holds
			\begin{align}\label{ptesK}
				|	\partial_x^l\tilde K(t,x)|\lesssim_l \frac{t}{(t^\frac{1}{4}+|x|)^{5+l}},
			\end{align}
			for any non-negative integer $l$. As a result, we have 
			\begin{align}\label{esL1K}
				\|\partial_x^l\tilde K(t,\cdot)\|_{L^1}\lesssim t^{-\frac{l}{4}}.
			\end{align}
			Then 
			\begin{align*}
				\left\|\partial_x^l\left(\sum_{n\in\mathbb{N}_+}e^{-\tilde A(n)t}e^{inx}\right)\right\|_{L^1}=\left\|\partial_x^n(\tilde K-(2\pi)^{-1}e^{-{c_0t}})\right\|_{L^1}\lesssim t^{-\frac{l}{4}}.
			\end{align*}
			This yields \eqref{lemrK}. 
			So it suffices to prove the point-wise estimate \eqref{ptesK}.
			For simple, we only prove \eqref{ptesK} for the case $l=0$, other cases can be done analogously. Indeed, 
			we can write 
			\begin{align*}
				\tilde K(t,x)&=\frac{1}{2\pi}\frac{1}{1-e^{-ix}}\sum_{n\in\mathbb{N}_+}e^{-\tilde A(n)t}(e^{inx}-e^{i(n-1)x})\\
				&=\frac{1}{2\pi}\frac{1}{1-e^{-ix}}\sum_{n\in\mathbb{N}_+}(e^{-\tilde A(n)t}-e^{-\tilde A(n+1)t})e^{inx}:=\sum_{n\in\mathbb{N}} \hat \nabla  \left(e^{-\tilde A(n)t}\right)e^{inx}.
			\end{align*}
			Here we 
			denote $\hat \nabla  f(k)=\frac{1}{2\pi}\frac{f(k)-f(k+1)}{1-e^{-ix}}$. \\
			We repeat the above procedure $5$ times to obtain 
			\begin{align*}
				\tilde	K(t,x)&=\sum_{n\in\mathbb{N}_+} \hat \nabla^5  \left(e^{-\tilde A(n)t}\right)e^{inx}.
			\end{align*}
			By the fact \begin{align*}
				&	\left|\hat \nabla ^5 \left(e^{-\tilde A(n)t}\right)\right|\lesssim \frac{t^5|n|^{15}+t^2|n|^3}{|x|^5}e^{-|n|^4t}\\
				&\ \ \lesssim \frac{t^\frac{5}{4}}{|x|^5}\left((|n|^4t)^\frac{15}{4}+(|n|^4t)^\frac{3}{4}\right)e^{-|n|^4t}\lesssim \frac{t^\frac{5}{4}}{|x|^5}e^{-|n|^4t/2},
			\end{align*}
			we have 
			\begin{align}\label{esK1}
				|\tilde K(t,x)|\lesssim \frac{t^\frac{5}{4}}{|x|^5}\sum_{n\in\mathbb{N}_+}e^{-|n|^4t/2}\lesssim \frac{t}{|x|^5},
			\end{align}
			where we use the estimate $\sum_{n\in\mathbb{N}_+}e^{-|n|^4t/2}\lesssim t^{-\frac{1}{4}}.$
			On the other hand, we have 
			\begin{align*}
				|\tilde K(t,x)|\lesssim \sum_{n\in\mathbb{N}_+}e^{-|n|^4t}\lesssim t^{-\frac{1}{4}}.
			\end{align*}
			Combining this with \eqref{esK1}, we get 
			$$
			|\tilde K(t,x)|\lesssim \frac{t}{(t^\frac{1}{4}+|x|)^5}.
			$$
			This verifies \eqref{ptesK} and completes the proof.
		\end{proof}
		~\vspace{0.3cm}\\
		With the estimate \eqref{esL1K}, the proof of Theorem \ref{lemmain} implies that the solution to \eqref{model} satisfies
		\begin{equation}\label{mainper}
			\begin{aligned}
				&\sup_{t\in[0,T]}(e^{c_0t}\|\partial_x h(t)\|_{L^\infty}+t^\frac{m+\kappa}{4}e^{c_0t}\|\partial_x h(t)\|_{\dot C^{m+\kappa}})\\
				&\quad\quad\quad\lesssim \|\partial_x h_0\|_{L^\infty}+\sup_{t\in[0,T]}\left(t^{\frac{\kappa}{4}}e^{c_0t}\|F(t)\|_{\dot C^{\kappa-3}}+t^{\frac{m+\kappa }{4}}e^{c_0t}\|F(t)\|_{\dot C^{m+\kappa -3}}\right),\ \ \ \ \ \forall\ T>0.
			\end{aligned}
		\end{equation}
		\begin{proof}[Proof of Theorem \ref{thmdaslo} \textit{ii)}]\\
			We rewrite the equation \eqref{eqdas} as 
			\begin{align*}
				\partial_t h+{\partial_x^4 h}+\frac{\partial_x ^2 h}{\bar h_0^2}
				=\tilde F[h]+N[h].\end{align*}
			Here 
			\begin{equation}\label{nfn2}
				\begin{aligned}
					\tilde F[h]=&G_1(\partial_xh)\partial_{x}^3h\partial_x^2h+G_2(\partial_xh)(\partial_x^2h)^3+\sum_{r=1}^3\frac{1}{h^r}\sum_{\substack{(n_1,n_2,n_3)\\
							2n_1+n_2< 3\\
							n_1+n_2+n_3\geq 2}}G_{r,n_1,n_2,n_3}(\partial_xh)(\partial_x^3h)^{n_1}(\partial_x^2h)^{n_2}(\partial_xh)^{n_3},\\
					N[h]=&\left(1-\frac{1}{\langle \partial_xh\rangle^4}\right)\partial_x^4 h+\left(\frac{1}{\bar h_0^2}-\frac{1}{h^2\langle \partial_xh\rangle^2}\right)\partial_x^2h=N_1[h]+N_2[h],
				\end{aligned}
			\end{equation}
			with $G_i, G_{r,n_1,n_2,n_3}$ satisfy \eqref{dascdcoe}. For $T>0$, define 
			\begin{equation}\label{defyt}
				\begin{aligned}
					&\|f\|_{Y_T,*}:=\sup_{t\in[0,T]}(e^{c_0t}\|\partial_x f(t)\|_{L^\infty}+t^\frac{m+\kappa }{4}e^{c_0t}\|\partial_x f\|_{\dot C^{m+\kappa }}),\\
					&\|f\|_{Y_T}:=\|f\|_{Y_T,*}+\|f\|_{L^\infty_TL^\infty}.
				\end{aligned}
			\end{equation}
			Let $\sigma, T>0$, define the space 
			\begin{align*}
				\mathcal{Y}_{\sigma,T}:=\left\{f\in L^\infty_T\dot W^{1,\infty}(\mathbb{S}): f|_{t=0}=h_0,\|f-\bar h_0\|_{Y_T}\leq \sigma\right\}.
			\end{align*}
			For any $f\in \mathcal{Y}_{\sigma,T}$, we define a map $\mathcal{S}f=h$, where $h$ solves the Cauchy problem
			\begin{equation}\label{map}
				\begin{aligned}
					&\partial_t h+{\partial_x^4 h}+\frac{\partial_x ^2 h}{\bar h_0^2}
					=\tilde F[f]+N[f],\\
					&h|_{t=0}= h_0.
				\end{aligned}
			\end{equation}
			We assume $\sigma\leq \frac{\bar h_0}{2}$, then $\inf_{x}|f(t,x)|\geq \frac{\bar h_0}{2}$ for any $f\in \mathcal{Y}_{\sigma,T}$. By \eqref{mainper}, we get 
			\begin{equation}\label{gaa1}
				\begin{aligned}
					\|h\|_{Y_{T,*}}\lesssim &\|\partial_x h_0\|_{L^\infty}+\sum_{j\in\{0,m\}}\sup_{t\in[0,T]}\left(e^{c_0t}t^{\frac{j+\kappa}{4}}\|N[f](t)\|_{\dot C^{j+\kappa-3}}+e^{c_0t}t^{\frac{j+\kappa}{4}}\|\tilde F[f](t)\|_{\dot C^{j+\kappa-3}}\right).
				\end{aligned}
			\end{equation}
			Observing that $N[\bar h_0]=\tilde F[\bar h_0]=0$, then apply Lemma \ref{dasgsml} with $f_1=f$, $f_2=\bar{h}_0$ to obtain
			\begin{equation*}
				\begin{aligned}
					&\sum_{j\in\{0,m\}}\sup_{t\in[0,T]}\left(e^{c_0t}t^{\frac{j+\kappa}{4}}\|N[f](t)\|_{\dot C^{j+\kappa-3}}+e^{c_0t}t^{\frac{j+\kappa}{4}}\|\tilde F[f](t)\|_{\dot C^{j+\kappa-3}}\right)\\
					&\quad \quad\quad\lesssim\|f-\bar h_0\|_{Y_T}^2(1+\|f-\bar h_0\|_{Y_T}+\bar h_0)^{m+4}.
				\end{aligned}
			\end{equation*}
			Hence 
			\begin{align}\label{yts}
				\|h\|_{Y_{T,*}}\leq &C_1\|\partial_x h_0\|_{L^\infty}+C_1\|f-\bar h_0\|_{Y_T}^2(1+\|f-\bar h_0\|_{Y_T}+\bar h_0)^{m+4}.
			\end{align}
			Then we estimate  $\|h(t)-\bar h_0\|_{L^\infty}$. We first estimate $\bar h(t)=\frac{1}{2\pi}\int_{-\pi}^{\pi}h(t,x)dx$.  Integrate \eqref{map} in $[0,T]\times \mathbb{S}$, we obtain 
			\begin{align*}
				\bar h(t)=\bar h_0+\frac{1}{2\pi}\int_0^t \int_{-\pi}^{\pi}(\tilde F[f]+N[f])(\tau,x)dxd\tau.
			\end{align*}
			We have 
			\begin{align*}
				\|\tilde F[f](t)\|_{L^\infty}+\|N[f](t)\|_{L^\infty}\lesssim t^{-\frac{3}{4}}e^{-2c_0t}\|f-\bar h_0\|_{Y_T}^2(1+\|f-\bar h_0\|_{Y_T})^5.
			\end{align*}
			Hence 
			\begin{align*}
				|	\bar h(t)-\bar h_0|\lesssim \int_0^t  \tau^{-\frac{3}{4}}e^{-2c_0\tau}d\tau \|f-\bar h_0\|_{Y_T}^2(1+\|f-\bar h_0\|_{Y_T}^5)\lesssim \|f-\bar h_0\|_{Y_T}^2(1+\|f-\bar h_0\|_{Y_T})^5.
			\end{align*}
			Combining this with the fact that $\sup_{x\in(-\pi,\pi)}|h(t,x)-\bar h(t)|\lesssim \|\partial_x h(t)\|_{L^\infty}\lesssim \|h\|_{Y_{T,*}}$ to obtain 
			\begin{equation}\label{est2}
				\begin{aligned}
					\|h(t)-\bar h_0\|_{L^\infty}&\leq \|h(t)-\bar h(t)\|_{L^\infty}+|\bar h(t)-\bar h_0|\\
					&\leq C_2\| h\|_{Y_{T,*}}+C_2\|f-\bar h_0\|_{Y_T}^2(1+\|f-\bar h_0\|_{Y_T})^5.
				\end{aligned}
			\end{equation}
			Finally,  we take $\varepsilon_0=\frac{\sigma}{(10+C_1)^2}$,  $\sigma<\min\left\{\frac{\bar h_0}{10},\frac{1}{((100+C_1+C_2+\bar h_0))^{m+10}} \right\}$. Then for any $f\in \mathcal{Y}_{\sigma,T}$, we obtain from \eqref{yts} that 
			\begin{equation}\label{dasgsder}
				\|h\|_{Y_{T,*}}\leq C_1\varepsilon_0+C_1\sigma^2(1+\sigma+\bar h_0)^{m+4}\leq \frac{\sigma}{4(1+C_2)}.
			\end{equation}  
			From this and \eqref{est2}, we get
			\begin{align*}
				\|h-\bar h_0\|_{L^\infty_TL^\infty}&\leq \frac{C_2\sigma}{4(1+C_2)}+C_2\sigma^2(1+\sigma)^5\leq \frac{\sigma}{4}.
			\end{align*}
			Combining this with \eqref{dasgsder} and \eqref{defyt} to obtain \begin{align*}
				\|h-\bar h_0\|_{Y_T}\leq 	\|h\|_{Y_{T,*}}+	\|h-\bar h_0\|_{L^\infty_TL^\infty}\leq \frac{\sigma}{2}.
			\end{align*}
			Hence $\mathcal{S}$ maps $\mathcal{Y}_{\sigma,T}$ to itself.\\ In the following, we prove that $\mathcal{S}$ is a contraction map. Let $f_1,f_2\in \mathcal{Y}_{\sigma,T}$, and $h_i=\mathcal{S}f_i$, $i=1,2$. We have 
			\begin{equation}\label{mapdif}
				\begin{aligned}
					&\partial_t (h_1-h_2)+{\partial_x^4 (h_1-h_2)}+\frac{\partial_x ^2 (h_1-h_2)}{\bar h_0^2}
					=\tilde F[f_1]-\tilde F[f_2]+N[f_1]-N[f_2],\\
					&(h_1-h_2)|_{t=0}=0.
				\end{aligned}
			\end{equation}
			Applying \eqref{mainper} again to obtain 
			\begin{align*}
				\|h_1-h_2\|_{Y_{T,*}}	\lesssim \sum_{j\in\{0,m\}}\sup_{t\in[0,T]}t^\frac{j+\kappa}{4}e^{c_0t}(\|(\tilde F[f_1]-\tilde F[f_2])(t)\|_{\dot C^{j+\kappa-3}}+\|(N[f_1]-N[f_2])(t)\|_{\dot C^{j+\kappa-3}}).
			\end{align*}
			Apply Lemma \ref{dasgsml} below to obtain 
			\begin{equation*}
				\begin{aligned}
					&\sum_{j\in\{0,m\}}\sup_{t\in[0,T]}t^\frac{j+\kappa}{4}e^{c_0t}(\|(\tilde F[f_1]-\tilde F[f_2])(t)\|_{\dot C^{j+\kappa-3}}+\|(N[f_1]-N[f_2])(t)\|_{\dot C^{j+\kappa-3}})\\ &\quad\quad\lesssim\|f_1-f_2\|_{Y_T}\|(f_1,f_2)\|_{Y_T,*}(1+\|(f_1-\bar h_0,f_2-\bar h_0)\|_{Y_T}+\bar h_0)^{m+4}.
				\end{aligned}
			\end{equation*}
			So we obtain 
			\begin{equation}\label{ythhh}
				\begin{aligned}
					\|h_1-h_2\|_{Y_{T,*}}	\leq C_3 \|f_1-f_2\|_{Y_T}\|(f_1,f_2)\|_{Y_T,*}(1+\|(f_1-\bar h_0,f_2-\bar h_0)\|_{Y_T}+\bar h_0)^{m+4}.
				\end{aligned}
			\end{equation}
			Then we estimate $\|h_1-h_2\|_{L^\infty}$. Integrate \eqref{mapdif} in time, we obtain 
			\begin{equation}\label{dfinf}
				\begin{aligned}
					\|(h_1-h_2)(t)\|_{L^\infty}\lesssim &\int_0^t \left(\|\partial_x^4(h_1-h_2)(\tau)\|_{L^\infty}+\frac{\|\partial_x^2(h_1-h_2)(\tau)\|_{L^\infty}}{\bar h_0}\right)d\tau\\
					&+\int_0^t\left(\|(\tilde F[f_1]-\tilde F[f_2])(\tau)\|_{L^\infty}+ \|(N[f_1]-N[f_2])(\tau)\|_{L^\infty}\right)d\tau.
				\end{aligned}
			\end{equation}
			Note that for any $\tau\in[0,T]$,
			\begin{align*}
				&	\|\partial_x^4(h_1-h_2)(\tau)\|_{L^\infty}+\|\partial_x^2(h_1-h_2)(\tau)\|_{L^\infty}\lesssim (\tau^{-\frac{1}{4}}+\tau^{-\frac{3}{4}})e^{-c_0\tau}\|h_1-h_2\|_{Y_{T,*}},\\
				&\|(\tilde F[f_1]-\tilde F[f_2])(\tau)\|_{L^\infty}+ \|(N[f_1]-N[f_2])(\tau)\|_{L^\infty}\\
				&\quad\quad\quad\lesssim  \tau^{-\frac{3}{4}}e^{-c_0\tau} \|f_1-f_2\|_{Y_T}\|(f_1-\bar h_0,f_2-\bar h_0)\|_{Y_T}(1+\|(f_1-\bar h_0,f_2-\bar h_0)\|_{Y_T})^3.
			\end{align*} 
			Combining this with \eqref{dfinf}, and the fact that 
			$$\int_0^\infty(\tau^{-\frac{1}{4}}+\tau^{-\frac{3}{4}})e^{-c_0\tau} d\tau \lesssim 1,
			$$
			we get 
			\begin{align}\label{h12inf}
				\|h_1-h_2\|_{L^\infty_TL^\infty}\leq C_4 \|h_1-h_2\|_{Y_{T,*}}+C_4\|f_1-f_2\|_{Y_T}\|(f_1-\bar h_0,f_2-\bar h_0)\|_{Y_T}(1+\|(f_1-\bar h_0,f_2-\bar h_0)\|_{Y_T})^3.
			\end{align}
			Now we fix $\sigma=\min\left\{\frac{\bar h_0}{10},\frac{1}{((100+C_1+C_2+C_3+C_4+\bar h_0))^{m+10}} \right\}$. Then it follows from \eqref{ythhh} that 
			\begin{align}\label{rh12yts}
				\|h_1-h_2\|_{Y_{T,*}}\leq 2 C_3\sigma(1+2\sigma+\bar h_0)^{m+4}\|f_1-f_2\|_{Y_T}\leq \frac{1}{4(1+C_4)}\|f_1-f_2\|_{Y_T}.
			\end{align}
			Combining this with \eqref{h12inf} to get
			\begin{align}\label{h1212}
				\|h_1-h_2\|_{L^\infty_TL^\infty}\leq\frac{C_4}{4(1+C_4)}\|f_1-f_2\|_{Y_T}+C_4\sigma(1+2\sigma)^3\|f_1-f_2\|_{Y_T}\leq \frac{1}{4}\|f_1-f_2\|_{Y_T}.
			\end{align}
			Hence we conclude from \eqref{rh12yts}, \eqref{h1212}, and the definition \eqref{defyt} that 
			\begin{equation*}
				\|h_1-h_2\|_{Y_T}\leq\frac{1}{2}\|f_1-f_2\|_{Y_T}.
			\end{equation*}
			Thus $\mathcal{S}:\mathcal{Y}_{\sigma,T}\rightarrow\mathcal{Y}_{\sigma,T}$ is a contraction map. This completes the proof.
		\end{proof}
		\begin{lemma}\label{dasgsml}
			Let $f_1,f_2:[0,\infty]\times\mathbb R^d\rightarrow\mathbb R$, then for $N[f]$, $\tilde F[f]$ defined in \eqref{nfn2}, and norms  $\|\cdot\|_{Y_T}$, $\|\cdot\|_{Y_{T,*}}$ defined in \eqref{defyt}, there holds
			\begin{equation*}
				\begin{aligned}
					&\|(N[f_1]-N[f_2])(t)\|_{\dot C^{n+\kappa-3}}+\|(\tilde F[f_1]-\tilde F[f_2])(t)\|_{\dot C^{n+\kappa-3}}\\
					&\quad\quad\quad\lesssim e^{-c_0t}t^{-\frac{n+\kappa}{4}}\|f_1-f_2\|_{Y_T,*}\|(f_1,f_2)\|_{Y_T,*}(1+\|(f_1,f_2)\|_{Y_T,*})^{n+2},\ \ \forall t\in[0,T],
				\end{aligned}
			\end{equation*}
			for $n\in\mathbb{N}$, $n\leq m$.
		\end{lemma}
		\begin{proof}
			For $N_1[f]$, denote $\gamma:=\kappa-3$. Then applying Lemma \ref{lemcom}, we can obtain for any $n\leq m$,
			\begin{equation*}
				\begin{aligned}
					\|\nabla^nN_1[f_1]-\nabla^nN_1[f_2]\|_{\dot C^{\gamma}}\lesssim& \sum_{l=0}^{n}\left(\|\partial_x(f_1-f_2) \|_{\dot C^{l+\gamma}}(\|\partial_x (f_1,f_2) \|_{\dot C^1}^{n-l}+\|\partial_x(f_1,f_2)\|_{\dot C^{n-l}})\right.\\
					&\left.+\|\partial_x(f_1-f_2) \|_{\dot C^l}(\|\partial_x(f_1,f_2)\|_{\dot C^\gamma}^{\frac{n-l+\gamma}{\gamma}}+\|\partial_x(f_1,f_2)\|_{\dot C^{n-l+\gamma}})\right)\\
					\lesssim& e^{-2c_0t}t^{-\frac{n+\kappa}{4}}\|f_1-f_2\|_{Y_T,*}\|(f_1,f_2)\|_{Y_T,*}(1+\|(f_1,f_2)\|_{Y_T,*})^{n+2}.
				\end{aligned}
			\end{equation*}
			For $N_2[f]$, note that since $f_i-\bar{h}_0$ are periodic functions, the lower order derivatives can be controlled by higher order derivatives, which means that $\|\partial_x^2h\|_{\dot C^\gamma}\lesssim\|\partial_x^4h\|_{\dot C^\gamma}$. Similarly by Lemma \ref{lemcom}, for any $n\leq m$, we will get
			\begin{equation*}
				\begin{aligned}
					\|\nabla^nN_2[f_1]-\nabla^nN_2[f_2]\|_{\dot C^{\gamma}}\lesssim e^{-c_0t}t^{-\frac{n+\kappa}{4}}\|f_1-f_2\|_{Y_T,*}\|(f_1-\bar{h}_0,f_2-\bar{h}_0)\|_{Y_T}(1+\|(f_1,f_2)\|_{Y_T,*})^{n+2}.
				\end{aligned}
			\end{equation*}
			For $F[f]$, we can decompose it into two parts, which are the critical one $F_3[f]$ and the lower order one $F_4[f]$. We can write as 
			\begin{align*}
				&F_3[f]=G_1(\partial_xf)\partial_{x}^3f\partial_x^2f+G_2(\partial_xf)(\partial_x^2f)^3,\\
				&F_4[f]=\sum_{r=1}^3\frac{1}{f^r}\sum_{\substack{(n_1,n_2,n_3)\\
						2n_1+n_2< 3\\
						n_1+n_2+n_3\geq 2}}G_{r,n_1,n_2,n_3}(\partial_xf)(\partial_x^3f)^{n_1}(\partial_x^2f)^{n_2}(\partial_xf)^{n_3}.
			\end{align*}
			By interpolation, we only need to estimate
			\begin{equation*}
				\|\nabla^{n}(F_i[f_1]-F_i[f_2])\|_{L^\infty},\quad i=3,4,\ n\leq m+1.
			\end{equation*}
			In fact, denote $\gamma:=\kappa-3$, then by Lemma \ref{maininterpo}, one has
			\begin{equation}\label{dasgloint}
				\|\nabla^{k}(F_i[f_1]-F_i[f_2])\|_{\dot C^{\gamma}}\lesssim\|\nabla^k(F_i[f_1]-F_i[f_2])\|_{L^\infty}^{1-\gamma}\|\nabla^{k+1}(F_i[f_1]-F_i[f_2])\|_{\dot C^{1}}^{\gamma}.
			\end{equation}
			For $F_3[f]$, we obtain 
			\begin{equation*}
				\|\nabla^n(F_3[f_1]-F_3[f_2])\|_{L^\infty}\lesssim e^{-2c_0t}t^{-\frac{n+2}{4}}\|f_1-f_2\|_{Y_T}\| (f_1,f_2)\|_{Y_T,*}(1+\|(f_1,f_2)\|_{Y_T})^{n+3}.
			\end{equation*}
			For $F_4[f]$, since $f_1,f_2\in \mathcal{X}_{T,\phi}^\sigma$, by \eqref{dasinfh} we obtain 
			\begin{equation*}
				\|\nabla^n(F_4[f_1]-F_4[f_2])\|_{L^\infty}\lesssim h_0^{-4-n}e^{-2c_0t}t^{-\frac{n+2}{4}}(1+t^{\frac{n+2}{4}})\|f_1-f_2\|_{Y_T}\|(f_1,f_2)\|_{Y_T,*}(1+\|(f_1,f_2)\|_{Y_T})^{n+5}.
			\end{equation*}
			Combine the estimates of $\|\nabla^nF_i[f_1]-\nabla^nF_i[f_2]\|_{L^\infty}$ with \eqref{dasgloint}, we obtain the result.
		\end{proof}

		\subsection{The Peskin problem in 2D}\label{secpeskin}
		We consider the problem of an elastic filament immersed in a two-dimensional Stokes fluid. It is inspired by the numerical immersed boundary method introduced by Peskin \cite{PeskinFlow1972,Peskin1972Thesis} to study the flow patterns around heart valves.  The numerical study for such fluid-structure interaction problems has attracted a lot of interest, which gives birth to wide applications in physics, biology, and medical sciences\cite{PeskinImmersed2002,MittalImmersed2005,HouNumerical2012}. The Peskin problem is named after Peskin in honor of his seminal contributions.

		Let $\Gamma$ be a simple closed curve which partitions $\mathbb{R}^2$ into two regions, the interior of the curve, $\Omega^+$ and the exterior $\Omega^-=\mathbb{R}^2\backslash\Omega^+$.
		Let $\Gamma$ be parameterized by vector valued function $X(t,x)=(X_1(t,x),X_2(t,x))\in\mathbb{R}^2$. Here $x\in\mathbb{S}:=\mathbb{R}/(2\pi\mathbb{Z})$ is the material coordinate and $t\geq 0$ denotes time. For fixed $x$, $X(t,x)$ moves with the local fluid velocity. Suppose further that the elastic structure has force density $F(X(t,x))$ with the form 
		\[
		F(X)=\partial_x(\mathcal{T}(|\partial_x X|)\tau(X)),
		\]
		where $\tau(X)=\frac{\partial_x X}{|\partial_x X|}$ is the unit tangent of the boundary $\Gamma$, and $\mathcal{T}:\mathbb{R}^+\to \mathbb{R}^+$ is the elastic tension which is assumed to be smooth and satisfies the structure condition $\mathcal{T}'>0$. Denote $u$ the fluid velocity and $p$ the pressure. The Peskin problem reads
		\begin{align*}
			\begin{cases}
				-\Delta u =-\nabla p~~~~~&\text{in} ~\mathbb{R}^2\backslash\Gamma(t),  \\
				\operatorname{div} u =0 ~~~~\quad\quad&\text{in} ~\mathbb{R}^2\backslash\Gamma(t),\\
				\llbracket u \rrbracket =0 ~~~~\quad\quad&\text{on} ~\Gamma(t),\\
				\llbracket\left(\nabla u+\nabla u^{T}-p \mathrm{Id}\right) n\rrbracket=\frac{F(X)}{|\partial_x X|}~~~&\text{on}~\Gamma(t),\\
				\partial_tX=u~~~~\quad\quad&\text{on}~ \Gamma(t).
			\end{cases}	
		\end{align*}
		Here $n$ is the outward unit normal to the free boundary $\Gamma(t)$ and $\llbracket\cdot\rrbracket$ denotes the jump across $\Gamma$:
		\[\llbracket U\rrbracket=U^+-U^-,\]
		where $U^{\pm}$ denotes the limiting value of $U$ evaluated on $\Gamma$ from the $\Omega^\pm$ side.
		The Peskin problem can be equivalently written as the following contour equations\cite{LinTongSolvability2019,MoriWell2019}
		\begin{align*}
			\partial_{t} X(t,x) &=\int_{\mathbb{S}} \mathbf{G}(X(t,x)-X(t,\sigma)) F(X(t,\sigma)) \mathrm{d} \sigma, \\
			\mathbf{G}(Z)&=\frac{1}{4\pi}\left(-\log |Z|         \mathrm{Id}+\frac{Z\otimes Z}{|Z|^2}\right),\quad \quad Z=(z_1,z_2)\in \mathbb{R}^2\backslash\{0\},\nonumber
		\end{align*}
		where $\mathbf{G}$ is the fundamental solution of the $2D$ Stokes problem. It is easy to check that if $X(t,x)$ is a solution, then for any $\lambda>0$, $X_\lambda(t,x)=\lambda^{-1}X(\lambda t,\lambda x)$ is also a solution. Under this scaling, $\dot B^1_{\infty,\infty}\supset BMO^1\supset $ $\dot W^{1,\infty}, $ and $\dot H^\frac{3}{2}$ are  critial spaces.
		
		The analytical study of the Peskin problem was initiated in \cite{LinTongSolvability2019,MoriWell2019}. 
		Lin and Tong \cite{LinTongSolvability2019} proved the local well-posedness for arbitrary $H^\frac{5}{2}$ data. Their proof relies on energy arguments and an application of the Schauder fixed point theorem. They also proved the global existence result and exponential decay towards equilibrium when the initial configuration is sufficiently close to the equilibrium. Tong \cite{TongRegularized2021} also established global well-posedness of a regularized Peskin problem and proved the convergence as the regularization parameter diminishes. Mori, Rodenberg and Spirn \cite{MoriWell2019} established a local well-posedness result for initial data in $C^{1,\gamma}$ with $\gamma>0$ (see also \cite{Rodenberg2018}).
		These spaces are subcritical under the scaling of the Peskin equation. For the well-posedness in critical spaces, Garcia-Juarez, Mori and Strain \cite{GarciaViscosityContrast2020} proved the global well-posedness result with initial data in the Wiener space $\mathcal{F}^{1,1}$ and sufficiently close to the stationary states. Their result holds for two fluids with different viscosity. Gancedo, Belinch\'{o}n and Scrobogna \cite{GancedoGlobal2020} considered a toy model of the Peskin problem and proved a global existence result in the critical Lipschitz space. In the  work \cite{KN} of the first and the third author, the local and global well-posedness of the Peskin problem (for Hookean
		material) are established in critical Besov spaces $\dot B^1_{\infty,\infty}$. The new ingredient in \cite{KN} is the construction of a new norm based on the structure of the nonlinear terms of the equation.   More recently, Cameron and Strain \cite{Camepeskin} considered the problem with fully non-linear tension, and proved local well-posedness in the critical Besov space $\dot B^\frac{3}{2}_{2,1}$.
		
		For simplicity, denote $\mathbf{T}(\lambda)=\frac{\mathcal{T}(\lambda)}{\lambda}$.		The Peskin equation \eqref{peskin} can be written as 
		\begin{equation}\label{eqpesk}
			\begin{aligned}
				\partial_t X(x)&=\int G(\delta_\alpha X(x))\partial_x\left({\mathbf{T}(|\partial_x X|)}\partial_x X\right)(x-\alpha)d\alpha\\
				&=-\int \partial_\alpha G(\delta_\alpha X(x))\delta_\alpha\left({\mathbf{T}(|\partial_x X|)}\partial_xX\right)d\alpha:=\mathrm{N}(X(x)).
			\end{aligned}
		\end{equation}
		Following \cite{KN}, we can reformulate the equation as 
		\begin{align}\label{eqpes}
			\partial_t X(x)+\frac{1}{4}\mathcal{H}({\mathbf{T}(|\partial_x X|)}\partial_x X)(x)=\mathcal{N}(X(x)),
		\end{align}
		where $\mathcal{H}$ is the Hilbert transform on torus, which is defined by
		$$
		\mathcal{H}f(x)=\frac{1}{2\pi}\int_{\mathbb{S}}f(x-\alpha)\cot\left(\frac{\alpha}{2}\right)d\alpha,
		$$
		and the nonlinear term reads
		\begin{equation}\label{defnonpes}
			\begin{aligned}
				\mathcal{N}(X)
				&=\frac{1}{4\pi}\int_{\mathbb{R}}  \frac{\tilde\Delta_\alpha X\cdot \tilde{\E}^\alpha X}{|\tilde\Delta_\alpha X|^2}\delta_\alpha ({\mathbf{T}(|\partial_x X|)}\partial_x X)\frac{d\alpha}{\alpha}\\
				&\quad\quad\quad\quad-\frac{1}{4\pi}\int_{\mathbb{R}}\frac{ \tilde{\E}^\alpha X\otimes \tilde\Delta_\alpha X+\tilde\Delta_\alpha X\otimes \tilde{\E}^\alpha X}{|\tilde\Delta_\alpha X|^2} \delta_\alpha ({\mathbf{T}(|\partial_x X|)}\partial_x X)\frac{d\alpha}{\alpha}\\
				&\quad\quad\quad\quad+\frac{1}{2\pi}\int_{\mathbb{R}} \frac{\tilde\Delta_\alpha X\otimes \tilde\Delta_\alpha X}{|\tilde\Delta_\alpha X|^4}\left(\tilde\Delta_\alpha X\cdot \tilde{\E}^\alpha X\right) \delta_\alpha ({\mathbf{T}(|\partial_x X|)}\partial_x X)\frac{d\alpha}{\alpha}.
			\end{aligned}
		\end{equation}
		Here we denote 
		\begin{align}\label{notapes}
			\tilde\Delta_\alpha X(x)=\frac{\delta_\alpha X(x)}{\tilde \alpha},\quad \tilde{\E}^\alpha X(x)=X'(x-\alpha)-\tilde\Delta_\alpha X(x), \quad \tilde \alpha =\left(\frac{1}{2}\cot \left(\frac{\alpha}{2}\right)\right)^{-1}. 
		\end{align}
		Take one spatial derivative in \eqref{eqpes}, we obtain 
		\begin{align}\label{eqonede}
			\partial_t \partial_x X+\frac{1}{4}\Lambda({\mathbf{T}(|\partial_x X|)}\partial_x X)=\partial_x\mathcal{N}(X).	\end{align}
		Fix $a\in(0,1)$ such that $1-a\ll 1$.	We write 
		\begin{align*}
			\Lambda({\mathbf{T}(|\partial_x X|)}\partial_x X)=\Lambda^{1-a}\Lambda^a({\mathbf{T}(|\partial_x X|)}\partial_x X).
		\end{align*}
		It is easy to check that 
		\begin{align*}
			\Lambda^a({\mathbf{T}(|\partial_x X|)}\partial_x X)&=\Lambda^a(\mathbf{T}(|\partial_x X|))\partial_x X+\mathbf{T}(|\partial_x X|)\Lambda^a\partial_x X-\frac{1}{\pi}\int_{\mathbb{R}} \delta_\alpha (\mathbf{T}(|\partial_x X|))\delta_\alpha \partial_x X\frac{d\alpha}{|\alpha|^{1+a}}\\
			&=\A(\partial_x X)\Lambda^a \partial_x X+ \mathsf{M}(\partial_xX),
		\end{align*}
		where we denote 
		\begin{align}
			&\A(\partial_x X)=\frac{1}{4}\mathbf{T}(|\partial_x X|)\mathrm{Id}+\frac{\mathbf{T}'(|\partial_x X|)}{4}\frac{\partial_x X\otimes \partial_x X}{|\partial_x X|},\nonumber\\
			&\mathsf{M}(\partial_xX)=-\frac{1}{\pi}\int_{\mathbb{R}} \delta_\alpha (\mathbf{T}(|\partial_x X|))\delta_\alpha \partial_x X\frac{d\alpha}{|\alpha|^{1+a}}+\frac{1}{\pi} \partial_x X\int_{\mathbb{R}} \delta_\alpha (\mathbf{T}(|\partial_x X|))-\delta_\alpha \partial_xX\cdot \nabla \big(\mathbf{T}(|\partial_x X|)\big)\frac{d\alpha}{|\alpha|^{1+a}}\label{defpeR}.
		\end{align}
		Note that 
		\begin{align*}
			\A(b)=\frac{1}{4}\frac{\mathcal{T}(|b|)}{|b|}\left(\mathrm{Id}-\frac{b\otimes b}{|b|^2}\right)+\frac{1}{4}\mathcal{T}'(|b|)\frac{b\otimes b}{|b|^2}.
		\end{align*}
		Hence
		\begin{align*}
			c(\mathcal{T},b)\mathrm{Id} \leq \A(b)\leq C(\mathcal{T},b)\mathrm{Id},\ \ \ \text{with}\   c(\mathcal{T},b)=\frac{1}{4}\min\{\mathcal{T}'(|b|),\frac{\mathcal{T}(|b|)}{|b|}\},\ \ C(\mathcal{T},b)=\frac{1}{4}\left(\mathcal{T}'(|b|)+\frac{\mathcal{T}(|b|)}{|b|}\right).
		\end{align*}
		Hence we can rewrite \eqref{eqonede} as
		\begin{align}\label{pesre}
			\partial_t \partial_x X+\Lambda^{1-a}(\A(\partial_x X)\Lambda^a \partial_x X)=\partial_x\mathcal{N}(X)-\Lambda^{1-a}\mathsf{M}(\partial_xX).
		\end{align}
		\begin{proof}[Proof of Theorem \ref{thmPes}]
			Take $\Phi=\partial_x X_0\ast \rho_{\eps_1}$.	By \eqref{pescon}, 
			\begin{align}\label{inism}
				\| \partial_x X_0-\Phi\|_{L^\infty}\leq \varepsilon_0,
			\end{align}
			We can rewrite \eqref{pesre} as  
			\begin{align*}
				\partial_t \partial_x X+\A(\Phi) \Lambda \partial_x X=\partial_x\mathcal{N}(X(x))+ \Lambda ^{1-a}F(\Phi,X)+R(\Phi,X),
			\end{align*}
			where 
			\begin{align*}
				&F(\Phi,X)=\mathsf{M}(\partial_xX)+(\A(\Phi)-\A(\partial_x X)) \Lambda ^a \partial_x X,\\
				&R(\Phi,X)=\A(\Phi) \Lambda \partial_x X- \Lambda ^{1-a}(\A(\Phi) \Lambda ^a \partial_x X).
			\end{align*}
			We introduce the norm 
			
			\begin{equation*}
				\begin{aligned}
					&\|h\|_{T}=\sup_{t\in[0,T]}(\|h(t)\|_{L^\infty}+t^{m+\kappa}\|h(t)\|_{\dot C^{m+\kappa}}),\\
					&\|h\|_{T,*}=\sup_{t\in[0,T]}(t^\frac{1}{10}\|h(t)\|_{\dot C^\frac{1}{10}}+t^{m+\kappa}\|h(t)\|_{\dot C^{m+\kappa}}).
				\end{aligned}
			\end{equation*}
			Let  $\sigma,T>0$. We consider the set 
			\begin{align}\label{defset}
				\mathcal{X}_{T,\Phi}^\sigma=\left\{Y\in L_T^\infty\dot W^{1,\infty}: Y|_{t=0}=X_0, \ \|\partial_x Y-\Phi\|_{T}\leq \sigma\right\}.
			\end{align}
			By definition and \eqref{inism}, it is easy to check that $|\Phi(x)|$ has the following upper and lower bound
			\begin{align*}
				\frac{1}{2\mathbf{\Theta}_0}\leq\frac{1}{\mathbf{\Theta}_0}-\frac{5\varepsilon_0}{4}\leq |\partial_x X_0(x)|-\|\partial_x X_0-\Phi\|_{L^\infty}\leq \frac{3}{2}\|\partial_x X_0\|_{L^\infty}, \ \ \ \ \forall x\in[-\pi,\pi),
			\end{align*}
			provided $\varepsilon_0\leq \min\{\frac{1}{100\mathbf{\Theta}_0},\frac{1}{2}\|\partial_xX_0\|_{L^\infty}\}$. Moreover, by taking $\sigma\leq \min\{\frac{1}{100\mathbf{\Theta}_0},\frac{1}{2}\|\partial_xX_0\|_{L^\infty}\}$, we have for any $Y\in\mathcal{X}^\sigma_{T,\phi}$,
			\begin{align*}
				\frac{1}{2\mathbf{\Theta}_0}\leq |\partial_x Y(t,x)|\leq 2\|\partial_x X_0\|_{L^\infty},\ \ \forall t\in[0,T], x\in[-\pi,\pi).
			\end{align*}
			Hence, we have for any $n\in\mathbb{N}$,
			\begin{align}
				\sup_{\tau\in[0,1]}|(\nabla^n\A)(\tau\partial_x Y+(1-\tau)\Phi)|\lesssim C_{\mathcal{T},n}(1+10\mathbf{\Theta}_0)^{n+1},
			\end{align}
			with 
			\begin{align}\label{consct}
				C_{\mathcal{T},n}=\sum_{k=0}^{n+1}\sup\{|\mathcal{T}^{(k)}(\lambda)|:\lambda\in[(10\mathbf{\Theta}_0)^{-1}, 2\|\partial_x X_0\|_{L^\infty}]\}.\end{align}
			For any $Y\in \mathcal{X}_{T,\Phi}^\sigma$, where $T, \sigma$ will be fixed later, we define a map $\mathcal{S}Y=X$, where $X$ is a solution to the Cauchy problem
			\begin{align*}
				&\partial_t \partial_x X+\A(\Phi) \Lambda \partial_x X=\partial_x\mathcal{N}(Y)+ \Lambda ^{1-a}F(\Phi,Y)+R(\Phi,Y),\\
				&\partial_xX|_{t=0}=\partial_xX_0.
			\end{align*}
			We first prove that there exist $\sigma,T>0$ such that $\mathcal{S}$ maps $\mathcal{X}_{T,\Phi}^\sigma$ to itself.
			Denote $\tilde X'=\partial_xX-\Phi$,  one has 
			\begin{align}\label{peseqr}
				\partial_t \tilde  X'+\A(\Phi) \Lambda \tilde  X'=\partial_x\mathcal{N}(Y)+ \Lambda ^{1-a}F(\Phi,Y)+R(\Phi,Y)-\A(\Phi) \Lambda \Phi.
			\end{align}
			Applying Theorem \ref{lemmain} with $G=R(\Phi,Y)-\A(\Phi) \Lambda \Phi$, $F=\mathcal{N}(Y)$ and $F(\Phi,Y)$ separately. we obtain that  there exists $T_0=T_0(\|\partial_x X_0\|_{L^\infty},\eps_0,\eps_1)>0$ such that  for any $0<T\leq T_0$,
			\begin{equation}
				\begin{aligned}\label{pesma}
					\|\tilde X'\|_T&\lesssim \|\tilde X'_0\|_{L^\infty}+\sum_{j=0,m} \sup_{t\in[0,T]}t^{j+\kappa}\left(\|\mathcal{N}(Y)(t)\|_{\dot C^{j+\kappa}}+\|F(\Phi,Y)(t)\|_{\dot C^{j+\kappa-a}}\right)\\
					&+\|R(\Phi,Y)\|_{L_T^1L^\infty}+T\|\A(\Phi) \Lambda \Phi\|_{L^\infty}+\sup_{t\in[0,T]}t^m\|R(\Phi,Y)(t)\|_{\dot C^m}+T^m\|\A(\Phi) \Lambda \Phi\|_{\dot C^m}.
				\end{aligned}
			\end{equation}
			By Lemma \ref{lemnonpes} we have
			\begin{equation}	\label{pesN}
				\begin{aligned}
					&\sup_{t\in[0,T]}(t^{\kappa}\|\mathcal{N}(Y)(t)\|_{\dot C^\kappa}+t^{m+\kappa}\|\mathcal{N}(Y)(t)\|_{\dot C^{m+\kappa}})\\
					&\quad\quad\quad\lesssim (1+\mathbf{\Theta}_Y(T))^{m+5}(\|\tilde Y'\|_{T}+T^\frac{1}{10}\|\Phi\|_{C^{m+1}})^2(1+\|\partial_x Y\|_T)^{m+5}.
				\end{aligned}
			\end{equation}
			Here we also used the fact that $\|\partial_x Y\|_{T,*}\leq \|\tilde Y'\|_{T,*}+\|\Phi\|_{T,*}\lesssim \|\tilde Y'\|_{T}+T^\frac{1}{10}\|\Phi\|_{C^{m+1}}$, where $\|\cdot \|_{T,*}$ is defined in \eqref{normstar}.\vspace{0.1cm}\\
			Then we estimate $F(\Phi,Y)$. We first deal with $\mathsf{M}(\partial_xY)$.
			Denote $\tilde Y'=\partial_x Y-\Phi$.	By Lemma \ref{lempesR} we have for any $0<T<\eps_1^2$,
			\begin{align*}
				\sup_{t\in[0,T]}(t^{\kappa}\|\mathsf{M}(\partial_xY)(t)-\mathsf{M}(\Phi)\|_{\dot C^{\kappa-a}}&+t^{m+\kappa}\|\mathsf{M}(\partial_x Y)(t)-\mathsf{M}(\Phi)\|_{\dot C^{m+\kappa-a}})\\&
				\lesssim \|\tilde Y'\|_T(\|\tilde Y'\|_{T}+T^\frac{1}{10}\|\Phi\|_{C^{m+1}})(1+\|\tilde Y'\|_T+\|\Phi\|_T)^{m+5}.
			\end{align*}
			Moreover, it is easy to check that 
			$$
			\sup_{t\in[0,T]}(t^{\kappa}\|\mathsf{M}(\Phi)\|_{\dot C^{\kappa-a}}+t^{m+\kappa}\|\mathsf{M}(\Phi)\|_{\dot C^{m+\kappa-a}})\lesssim T^\frac{1}{10}\|\Phi\|_{C^{m+3}}(1+\|\Phi\|_{C^{m+3}})^{m+3}.
			$$
			Then we obtain
			\begin{equation}\label{RRR}
				\begin{aligned}
					&\sup_{t\in[0,T]}(t^{\kappa}\|\mathsf{M}(\partial_xY)(t)\|_{\dot C^{\kappa-a}}+t^{m+\kappa}\|\mathsf{M}(\partial_x Y)(t)\|_{\dot C^{m+\kappa-a}})\\
					&\quad\quad\quad\quad\lesssim \|\tilde Y'\|_T(\|\tilde Y'\|_{T}+T^\frac{1}{10}\|\Phi\|_{C^{m+1}})(1+\|\tilde Y'\|_T+\|\Phi\|_T)^{m+5}+T^\frac{1}{10}\|\Phi\|_{C^{m+3}}(1+\|\Phi\|_{C^{m+3}})^{m+3}.
				\end{aligned}
			\end{equation}
			Moreover, by Lemma \ref{Hol}, we find that for any $0\leq n\leq m$,
			\begin{align*}
				&	\|(\A(\Phi)-\A(\partial_x Y)) \Lambda ^a \partial_x Y(t)\|_{\dot C^{n+\kappa-a}}\\
				&\quad	\lesssim \sum_{l=0}^n\left(\|\partial_x^l(\A(\Phi)-\A(\partial_x Y))(t)\|_{\dot C^{\kappa-a}}\| \Lambda ^a\partial_x^{1+n-l} Y(t)\|_{L^\infty}\right.\\
				&\quad\quad\quad\quad\left.+\|\partial_x^l(\A(\Phi)-\A(\partial_x Y))(t)\|_{L^\infty}\| \Lambda ^a\partial_x ^{1+n-l}Y(t)\|_{\dot C^{\kappa-a}}\right)\\
				&\quad	\lesssim C_{\mathcal{T},n}t^{-(n+\kappa)} (\|\tilde Y'\|_{T}+T^\frac{1}{10}\|\Phi\|_{C^{m+1}})^2(1+\|\tilde Y'\|_{T}+\|\Phi\|_{T})^{m+5},
			\end{align*}
			where $C_{\mathcal{T},n}$ is defined in \eqref{consct}.
			Hence  
			\begin{align*}
				&\sup_{t\in[0,T]}(t^\kappa\|(\A(\Phi)-\A(\partial_x Y)) \Lambda ^a \partial_x Y(t)\|_{\dot C^{\kappa-a}}+t^{m+\kappa}\|(\A(\Phi)-\A(\partial_x Y)) \Lambda ^a \partial_x Y(t)\|_{\dot C^{m+\kappa-a}})\\
				&\quad\quad\quad\lesssim C_{\mathcal{T},n}(\|\tilde Y'\|_{T}+T^\frac{1}{10}\|\Phi\|_{C^{m+1}})^2(1+\|\tilde Y'\|_{T}+\|\Phi\|_{T})^{m+5}.
			\end{align*}
			Combining this  with \eqref{RRR} to obtain that 
			\begin{equation}
				\begin{aligned}
					\sup_{t\in[0,T]}(t^\kappa\|F(\Phi,Y)(t)\|_{\dot C^{\kappa-a}}&+t^{m+\kappa}\|F(\Phi,Y)(t)\|_{\dot C^{m+\kappa-a}})\\
					&\lesssim (\|\tilde Y'\|_{T}+T^\frac{1}{10}\|\Phi\|_{C^{m+1}})^2(1+\|\tilde Y'\|_{T}+\|\Phi\|_{T})^{m+5}.
				\end{aligned}
			\end{equation}	
			It remains to estimate lower order terms $R(\Phi,Y)$ and $\A(\Phi) \Lambda \Phi$.
			Recall that 
			\begin{align*}
				R(\Phi,Y)=\A(\Phi) \Lambda \partial_x Y- \Lambda ^{1-a}(\A(\Phi) \Lambda ^a \partial_x Y)=c\int_{\mathbb{R}} \A(\Phi)\delta_\alpha ( \Lambda ^a \partial_x Y)-\delta_\alpha ( \A(\Phi) \Lambda ^a \partial_x Y)\frac{d\alpha}{|\alpha|^{2-a}}.
			\end{align*}
			By the identity $\delta_\alpha (fg)(x)=f(x)\delta_\alpha g(x)+\delta_\alpha f(x)g(x-\alpha)$, we have 
			\begin{equation}
				\begin{aligned}\label{R222}
					&\|	R(\Phi,Y)(t)\|_{L^\infty}\lesssim \int_{\mathbb{R}} \|\delta_\alpha \A(\Phi)\|_{L^\infty}\frac{d\alpha}{|\alpha|^{2-a}}\| \Lambda ^a\partial_x Y(t)\|_{L^\infty}\lesssim \|\Phi\|_{C^1}\|( t^{-a}\|\tilde Y'\|_T+\|\Phi\|_{C^3}),\\
					& \|\partial_x ^m R(\Phi,Y)(t)\|_{L^\infty}\lesssim (1+\|\Phi\|_{C^{m+1}})^{m+1}(t^{-(m+a)}\|\tilde Y'\|_T+\|\Phi\|_{C^{m+3}}).
				\end{aligned} 
			\end{equation}
			Finally, it is easy to check that $\|\A(\Phi) \Lambda \Phi\|_{L^\infty}\lesssim (1+\|\Phi\|_{L^\infty})\|\Phi\|_{C^{2}}$. Hence 
			\begin{equation}\label{peslow}
				\begin{aligned}
					\|R(\Phi,Y)\|_{L_T^1L^\infty}+T\|\A(\Phi) \Lambda \Phi\|_{L^\infty}&+\sup_{t\in[0,T]}t^{m+1}\|R(\Phi,Y)(t)\|_{\dot C^m}+T^{m+1}\|\A(\Phi) \Lambda \Phi\|_{\dot C^m}\\
					&\lesssim T^{1-a}(1+\|\tilde Y'\|_T+\|\Phi\|_{C^{m+3}})^{m+3}.
				\end{aligned}
			\end{equation}
			Combining \eqref{pesN}-\eqref{peslow} with \eqref{pesma},  we obtain for any $0<T\leq T_0$,
			\begin{equation}\label{maespes}
				\begin{aligned}
					\|\tilde X'\|_T\leq& C_0\|\tilde X'_0\|_{L^\infty}+C_0(1+\mathbf{\Theta}_Y(T)+\mathbf{\Theta}_0)^{m+5}(\|\tilde Y'\|_{T}+T^\frac{1}{10}\|\Phi\|_{C^{m+1}})^2(1+\|\tilde Y'\|_{T}+\|\Phi\|_{C^{m+1}})^{m+1}\\
					&+C_0(1+\mathbf{\Theta}_0)^{m+5}T^\kappa(1+\|\Phi\|_{C^{m+3}})^{m+5}+C_0T^{1-a}\|\Phi\|_{C^{2}}(1+\|\tilde Y'\|_T+\|\Phi\|_{C^{m+3}})^{m+3}.
				\end{aligned}
			\end{equation}
			Take $0<\varepsilon_0<\frac{1}{100(2C_0+2\mathbf{\Theta}_0+1)^{m+5}}$ and $\sigma=2C_0\varepsilon_0$.  Take $T_1=\min\left\{T_0,\left(\frac{\varepsilon_0\eps_1}{10+C_0+\|\partial_x X_0\|_{L^\infty}}\right)^{100m}\right\}$. Then it is easy to check that 
			for any  $Y\in\mathcal{X}^\sigma_{T_1,\Phi}$, we have $\|\tilde Y'\|_{T_1}\leq \sigma$ and
			\begin{align}\label{wes}
				\inf_{\alpha,x}|\Delta_\alpha Y(x)|\geq 	\inf_{\alpha,x}|\Delta_\alpha X_0(x)|-\|\partial_xX_0-\Phi\|_{L^\infty}-\|\Phi-\partial_x Y\|_{L^\infty}\geq \frac{1}{\mathbf{\Theta}_0}-\frac{5}{4}\varepsilon_0-\sigma\geq  \frac{1}{2\mathbf{\Theta}_0},
			\end{align}
			which implies that  $\mathbf{\Theta}_Y(T_1)\leq 2\mathbf{\Theta}_0$. Then \eqref{maespes} leads to 
			\begin{align*}
				\|\tilde X'\|_{T_1}\leq&\frac{5}{4}C_0\varepsilon_0+2C_0^2\varepsilon_0 (1+3\mathbf{\Theta}_0)^{m+5}(2C_0\varepsilon_0+T_1^\frac{1}{10}\eps_1^{-m-1}\|X_0'\|_{L^\infty})(1+2C_0\varepsilon_0+\eps_1^{-m-1}\|X_0'\|_{L^\infty})^{m+1}\\
				&+C_0(1+\mathbf{\Theta}_0)^{m+5}T_1^\kappa(1+\eps_1^{-2}\|X_0'\|_{L^\infty})^{m+5}+T_1^{1-a}\eps_1^{-m-3}\|X_0'\|_{L^\infty}(1+2C_0\varepsilon_0+\eps_1^{-m-3}\|X_0'\|_{L^\infty})^{m+3}\\
				\leq &\sigma.
			\end{align*}
			Hence $\mathcal{S}$ maps $\mathcal{X}_{T_1,\Phi}^\sigma$ to itself. \\
			Now we prove that $\mathcal{S}$ is a contraction map. Consider $Y_1,Y_2\in \mathcal{X}_{T_1,\Phi}^\sigma$. Denote $\mathbf{Y}'=\partial_x (Y_1-Y_2)$, and $\mathbf{X}'=\partial_x (X_1-X_2)=\partial_x (\mathcal{S}Y_1-\mathcal{S}Y_2)$. Then we have 
			\begin{align*}
				\partial_t \mathbf{X}'+\A(\Phi) \Lambda \mathbf{X}'=\partial_x(\mathcal{N}(Y_1)-\mathcal{N}(Y_2))+ \Lambda ^{1-a}(F(\Phi,Y_1)-F(\Phi,Y_2))+R(\Phi,Y_1)-R(\Phi,Y_2).
			\end{align*}
			Applying Theorem \ref{lemmain} again, with $G=R(\Phi,Y_1)-R(\Phi,Y_2)$, $F=\mathcal{N}(Y_1)-\mathcal{N}(Y_2)$ and $F(\Phi,Y_1)-F(\Phi,Y_2)$ separately, we have for any $0<T\leq T_0$,
			\begin{align*}
				\|\mathbf{X}'\|_{T}\lesssim& \sup_{t\in[0,T]}\left(t^{\kappa}\|(\mathcal{N}(Y_1)-\mathcal{N}(Y_2))(t)\|_{\dot C^\kappa}+t^{m+\kappa}\|(\mathcal{N}(Y_1)-\mathcal{N}(Y_2))(t)\|_{\dot C^{m+\kappa}}\right)\\
				&+\sup_{t\in[0,T]}\left(t^{\kappa}\|(F(\Phi,Y_1)-F(\Phi,Y_2))(t)\|_{\dot C^{\kappa}}+t^{m+\kappa}\|(F(\Phi,Y_1)-F(\Phi,Y_2))(t)\|_{\dot C^{m+\kappa-a}}\right)\\
				&+\|R(\Phi,Y_1)-R(\Phi,Y_2)\|_{L^1_TL^\infty}d\tau+\sup_{t\in[0,T]}t^{m+1}\|(R(\Phi,Y_1)-R(\Phi,Y_2))(t)\|_{\dot C^m}.
			\end{align*}
			We obtain from Lemma \ref{lemnonpes} that 
			\begin{align*}
				&\sup_{t\in[0,T]}\left(t^{\kappa}\|(\mathcal{N}(Y_1)-\mathcal{N}(Y_2))(t)\|_{\dot C^\kappa}+t^{m+\kappa}\|(\mathcal{N}(Y_1)-\mathcal{N}(Y_2))(t)\|_{\dot C^{m+\kappa}}\right)\\
				&\quad\lesssim C_{\mathcal{T},m} (1+\mathbf{\Theta}_{Y_1}(T)+\mathbf{\Theta}_{Y_2}(T))^{m+5} \|\mathbf{Y}'\|_{T}(\|\partial_x (\tilde{Y}_1, \tilde{Y}_2)\|_{T}+T^{\frac{1}{10}}\|\Phi\|_{C^{m+3}})(1+\|\partial_x (Y_1,Y_2)\|_T)^{5(m+1)},
			\end{align*}
			where $\tilde{Y}_i=Y_i-\Phi$. Then we estimate $F$. By Lemma \ref{lempesR}, we get 
			\begin{align*}
				&\sup_{t\in[0,T]} (t^{\kappa}\|(\mathsf{M}(\partial_xY_1)-\mathsf{M}(\partial_xY_2))(t)\|_{\dot C^{\kappa-a}}+t^{m+\kappa}\|(\mathsf{M}(\partial_xY_1)-\mathsf{M}(\partial_xY_2))(t)\|_{\dot C^{m+\kappa-a}})\\
				&\quad\quad\quad\quad\lesssim C_{\mathcal{T},m}\|\partial_x(Y_1-Y_2)\|_T\|\partial_x(Y_1, Y_2)\|_{T,*}(1+\|\partial_x(Y_1,Y_2)\|_T)^{m+5}.
			\end{align*}
			Moreover, observe that $R(\Phi,Y_1)-R(\Phi,Y_2)=R(\Phi,Y_1-Y_2)$. Following \eqref{R222} to obtain 
			\begin{align*}
				\|(R(\Phi,Y_1)-R(\Phi,Y_2))(t)\|_{L^\infty}\lesssim \int_{\mathbb{R}} |\delta_\alpha \A(\Phi)|\frac{d\alpha}{|\alpha|^{2-a}}\| \Lambda ^a\mathbf{Y}'(t)\|_{L^\infty}\lesssim t^{-a}\|\Phi\|_{C^1}\|\mathbf{Y}'\|_T.
			\end{align*}
			Hence 
			\begin{align*}
				\|\mathbf{X}'\|_{T_1}\leq& C_1 (1+\mathbf{\Theta}_{Y_1}+\mathbf{\Theta}_{Y_2})^{m+5} \|\mathbf{Y}'\|_{T_1}(\|(\tilde Y_1',\tilde Y_2')\|_{T_1}+T_1^\frac{1}{10}\|\Phi\|_{C^{m+3}})(1+\| (\tilde Y_1',\tilde Y_2')\|_{T_1}+\|\Phi\|_{C^{m+3}})^{m+3}\\
				&+C_1T_1^{1-a}\|\Phi\|_{C^{m+3}}\|\mathbf{Y}'\|_{T_1}\\
				\leq&\left(4C_1(1+4\mathbf{\Theta}_0)^{m+5}\sigma(1+4\sigma+T_1^\frac{1}{10}\eps_1^{-m-3}\|\partial_x X_0\|_{L^\infty})^{m+3}+C_1T_1^{1-a}\eps_1^{-m-3}\|\partial_x X_0\|_{L^\infty}\right)\|\mathbf{Y}'\|_{T_1},
			\end{align*}
			where we use the fact that $\|\tilde  Y_1'\|_{T_1},\|\tilde  Y_2'\|_{T_1}\leq \sigma$ and  $\mathbf{\Theta}_{Y_1},\mathbf{\Theta}_{Y_2}\leq 2\mathbf{\Theta}_0$. We have 
			\begin{align*}
				4C_1(1+4\mathbf{\Theta}_0)^{m+5}\sigma(1+4\sigma+T_1^\frac{1}{10}\eps_1^{-m-3}\|\partial_x X_0\|_{L^\infty})^{m+3}+C_1T_1^{1-a}\eps_1^{-m-3}\|\partial_x X_0\|_{L^\infty}\leq \frac{1}{2}
			\end{align*}
			by taking $\varepsilon_0=\frac{1}{100(2C_0+2C_1+2\mathbf{\Theta}_0+1)^{m+10}}$ and the corresponding $T$ small enough. This implies that 
			$\mathcal{S}:\mathcal{X}_{T_1,\Phi}^\sigma \to \mathcal{X}_{T_1,\Phi}^\sigma$ is a contraction map. Hence there exists a unique $X\in \mathcal{X}_{T_1,\Phi}^\sigma$ such that $X=\mathcal{S}X$, which indicates that $X$ is a solution to \eqref{eqpes}. Moreover, one has $\mathbf{\Theta}_X(T_1)\leq 2\mathbf{\Theta}_0$. This completes the proof.
		\end{proof}

		\vspace{0.5cm}
		Theorem \ref{thmPes} is for local well-posedness result of the 2d Peskin \eqref{eqpes}. We can extend this theorem to global well-posedness result. Let us discuss about this.\vspace{0.1cm}\\
		It is well-known that the only stationary solutions to the 2D peskin problem \eqref{eqpesk} are the uniformly parametrized circles which form a four-dimensional vector space, see \cite{LinTongSolvability2019}. If we treat $X(x)$ as a complex-valued function in the sense that $X=X_1+iX_2$, then the set of stationary solutions is $\{X(x)=a_1+a_2e^{ix}:a_1,a_2\in\mathbb{C},a_2\neq 0\}:=\mathcal{V}\backslash\{0\}$.
		Let $\mathcal{P}$ be the $L^2$ projection on to the space $\mathcal{V}$:
		\begin{align*}
			\mathcal{P}{Z}(t,x)=z_0(t)+z_1(t)e^{ix},\ \ \ \text{where}\ z_0(t)=\langle Z(t),1\rangle_{L^2(\mathbb{S})}, z_1(t)=\left\langle Z(t),e^{ix}\right\rangle_{L^2(\mathbb{S})}.
		\end{align*}
		Recently, Garcia-Juarez and V. Haziot \cite{GHpeskin} showed that the following the linearized system
		\begin{align*}
			\partial_tX=\left.\frac{d}{d\eps}\frac{\N(Z+\eps X)}{d\eps}\right|_{\eps=0}:=\mathcal{L}X
		\end{align*}
		satisfies
		\begin{equation}
			\langle\mathcal{L}X_1,X_1\rangle_{L^2(\mathbb{S})} \leq - 2c_0||X_1||_{L^2(\mathbb{S})}^2~~\forall X_1\in \mathcal{V}^\perp.
		\end{equation}
		Therefore, it is not hard to check that its kernel satisfies the point-wise estimate
		\begin{align*}
			|\partial_x^l(\mathrm{I}-\mathcal{P})K(t,\cdot )(x)|\lesssim _l\frac{e^{-c_0t}}{t^{1+l}}\langle\frac{|x|}{t}\rangle^{-(l+3)},\ \ \quad\quad \forall l\in\mathbb{N}.
		\end{align*}
		Define
		\begin{align*}
			\|Z\|_T:=\sup_{t\in [0,T]}e^\frac{c_0t}{4}(\|\partial_tz_0\|_{L^\infty}+\|\partial_t z_1(t)\|_{L^\infty}+\|(\mathrm{I}-\mathcal{P})Z\|_{\dot W^{1,\infty}}+t^{m+\kappa}\|(\mathrm{I}-\mathcal{P})Z\|_{\dot C^{m+1+\kappa}}).
		\end{align*}
		Consider the Cauchy problem \eqref{eqpesk} with initial data $X_0$ satisfying $\|(\mathrm{I}-\mathcal{P})X_0\|_{\dot W^{1,\infty}}\ll 1$, if $X$ is a solution in $[0,T]$, then using \eqref{mainper}, Lemma \ref{lemnonpes} and Lemma \eqref{lempesR}, we can estimate 
		\begin{align*}
			\|X\|_T \leq C \|(\mathrm{I}-\mathcal{P})X_0\|_{\dot W^{1,\infty}}+C\|X\|_T^2(1+\|X\|_T)^{m+3}.
		\end{align*}
		By bootstrap argument one obtains that there exists a global solution satisfying
		\begin{align*}
			\|X\|_\infty\leq 2C\|(\mathrm{I}-\mathcal{P})X_0\|_{\dot W^{1,\infty}}.
		\end{align*}
		This indicates that the global solution converges
		exponentially in time to a stationary circle solution, i.e   $$\tilde X(t,x)=X(t,x)-\langle X(t=\infty),1\rangle_{L^2(\mathbb{S})}-\left\langle X(t=\infty),e^{ix}\right\rangle_{L^2(\mathbb{S})}e^{ix}$$
		has the following decay estimate:  
		\begin{align*}
			\|\tilde X(t)\|_{L^\infty}+\|\tilde X(t)\|_{\dot W^{1,\infty}}+t^{m+\kappa}\|\tilde X(t)\|_{\dot C^{m+1+\kappa}}\leq C'\|(\mathrm{I}-\mathcal{P})X_0\|_{\dot W^{1,\infty}} e^{-\frac{c_0t}{8}}
		\end{align*}
		for any $t>0.$
		\subsection{The Peskin problem in 3D}
		We introduce the Peskin problem in 3D, which describes  a two-dimensional
		elastic membrane immersed in a three-dimensional steady Stokes flow. 
		Compared to the 2D case, the Peskin problem in 3D is more complex due to the additional dimension and the resulting non-local effects. 
		We consider general nonlinear elastic laws, i.e., the
		fully nonlinear Peskin problem, and prove that the problem is locally well-posed in critical $C^1$ space.

		The authors in \cite{3Dpeskin} derived the following Boundary Integral formulation,
		\begin{equation*}
			\begin{aligned}
				&	\frac{\partial F}{\partial t}(\widehat{\boldsymbol{x}}) =\int_{\mathbb{S}^2} G(F(\widehat{\boldsymbol{x}})-F(\widehat{\boldsymbol{y}})) \nabla_{\mathbb{S}^2} \cdot\left(\mathcal{T}\left(\left|\nabla_{\mathbb{S}^2} F(\widehat{\boldsymbol{y}})\right|\right) \frac{\nabla_{\mathbb{S}^2} F(\widehat{\boldsymbol{y}})}{\left|\nabla_{\mathbb{S}^2}F(\widehat{\boldsymbol{y}})\right|}\right) d \mu_{\mathbb{S}^2}(\widehat{\boldsymbol{y}}), \\
				&	\left.F(\widehat{\boldsymbol{x}})\right|_{t=0} =F_0(\widehat{\boldsymbol{x}}),
			\end{aligned}
		\end{equation*}
		where $\nabla_{\mathbb{S}^2}$ denotes the surface gradient operator on the unit sphere, $\mu_{\mathbb{S}^2}$ is the standard measure on the unit sphere, $G(\widehat{\boldsymbol{x}})$ is the 3D Stokeslet tensor
		$$
		G(\widehat{\boldsymbol{x}})=\frac{1}{8 \pi}\left(\frac{1}{|\widehat{\boldsymbol{x}}|} \mathrm{Id}_3+\frac{\widehat{\boldsymbol{x}} \otimes \widehat{\boldsymbol{x}}}{|\widehat{\boldsymbol{x}}|^3}\right).
		$$
		In the following of this section, we shortly denote $\mathbf{T}(\left|\nabla_{\mathbb{S}^2} F(\widehat{\boldsymbol{y}})\right|)= \frac{\mathcal{T}\left(\left|\nabla_{\mathbb{S}^2} F(\widehat{\boldsymbol{y}})\right|\right)}{\left|\nabla_{\mathbb{S}^2}F(\widehat{\boldsymbol{y}})\right|}$.\\
		Due to the non-trivial geometry in the 3D setting, we need to work in local charts. This makes the problem more complicated in view of the non-local character of the equation. However, we mention that the essential structure of 3D problem is the same as that of 2D problem.\\
		Recall that the initial data $F_0$ satisfies
		\begin{equation}\label{3dpesnc}
			\sup_{\widehat{\boldsymbol{x}},\widehat{\boldsymbol{y}}\in\mathbb{S}^2}\frac{|\widehat{\boldsymbol{x}}-\widehat{\boldsymbol{y}}|}{|F_0(\widehat{\boldsymbol{x}})-F_0(\widehat{\boldsymbol{y}})|}\leq \mathbf \Theta_0<\infty.
		\end{equation}
		By the standard stereographic projection $\mathcal{X}:\mathbb{R}^2\rightarrow\mathbb{S}^2$, see \cite[Definition 3.9]{3Dpeskin}, we can  
		transform the equation from $\mathbb{S}^2$ to $\mathbb{R}^2$. Precisely, we define
		\begin{equation*}
			\mathcal{X}(\theta)=\left(\frac{2\theta_1}{1+|\theta|^2},\frac{2\theta_2}{1+|\theta|^2},\frac{-1+|\theta|^2}{1+|\theta|^2} \right).
		\end{equation*}
		We can see that for any compact set $K\Subset\mathbb{S}^2$ with $(0,0,1)\notin K$, $\mathcal{X}$ is a homeomorphism from $\mathcal{X}^{-1}(K)$ to $K$. By local chart, for $F\in C(\mathbb{S}^2)$ with $ (0,0,1)\notin \text{Supp}(F)$, we can define $\nabla_{\mathbb{S}^2}F$ as 
		\begin{equation*}
			(\nabla_{\mathbb{S}^2}F)\circ\mathcal{X}(\theta)=\left(\frac{1+|\theta|^2}{2}\right)^2\sum_i\frac{\partial(F\circ\mathcal{X})}{\partial\theta_i}\frac{\partial\mathcal{X}}{\partial\theta_i},
		\end{equation*}
		which naturally gives
		\begin{equation*}
			\left|\nabla_{\mathbb{S}^2}F\right|\circ\mathcal{X}=\frac{1+|\theta|^2}{2}\left|\nabla(F\circ\mathcal{X})\right|.
		\end{equation*}
		We will shortly denote $\widetilde\nabla=\nabla_{\mathbb{S}^2}$ in this section.\\
		Denote $\widehat{\boldsymbol{x}}_*=(0,0,-1)$. For any $\widehat{\boldsymbol{x}}_0\in\mathbb{S}^2$, we define the rotation $\mathcal R_{\widehat{\boldsymbol{x}}_0}:\mathbb{S}^2\rightarrow\mathbb{S}^2$ such that $\mathcal R_{\widehat{\boldsymbol{x}}_0}\widehat{\boldsymbol{x}}_*=\widehat{\boldsymbol{x}}_0$, and define the smooth cut-off functions $0\leq\chi_{\widehat{\boldsymbol{x}}_0},\tilde{\chi}_{\widehat{\boldsymbol{x}}_0}\leq 1$ such that $\mathrm{Supp}(\chi_{\widehat{\boldsymbol{x}}_0})\subset \tilde\chi_{\widehat{\boldsymbol{x}}_0}^{-1}(1)$, and
		\begin{equation*}
			\begin{aligned}
				&\text{Supp}(\chi_{\widehat{\boldsymbol{x}}_0}\circ\mathcal{R}_{\widehat{\boldsymbol{x}}_0})\subset B(\widehat{\boldsymbol{x}}_*,a_2)\cap\mathbb{S}^2,\quad \chi_{\widehat{\boldsymbol{x}}_0}\circ\mathcal{R}_{\widehat{\boldsymbol{x}}_0}(B(\widehat{\boldsymbol{x}}_*,a_1)\cap\mathbb{S}^2)=1,\\
				&\text{Supp}(\tilde\chi_{\widehat{\boldsymbol{x}}_0}\circ\mathcal{R}_{\widehat{\boldsymbol{x}}_0})\subset B(\widehat{\boldsymbol{x}}_*,a_4)\cap\mathbb{S}^2,\quad \tilde \chi_{\widehat{\boldsymbol{x}}_0}\circ\mathcal{R}_{x_0}(B(\widehat{\boldsymbol{x}}_*,a_3)\cap\mathbb{S}^2)=1,\\
				&B(\widehat{\boldsymbol{x}}_*,a_1)\cap\mathbb{S}^2\subset \mathcal{X}(B(0,R_0))\subset \mathcal{X}(B(0,2R_0))\subset B(\widehat{\boldsymbol{x}}_*,a_2)\cap\mathbb{S}^2,\\
				&B(\widehat{\boldsymbol{x}}_*,a_3)\cap\mathbb{S}^2\subset \mathcal{X}(B(0,3R_0))\subset \mathcal{X}(B(0,4R_0))\subset B(\widehat{\boldsymbol{x}}_*,a_4)\cap\mathbb{S}^2,
			\end{aligned}
		\end{equation*}
		where $0<a_1<a_2<a_3<a_4$, $B(a,r)$ denotes the sphere with center $a$ and radius $r$ on $\mathbb{R}^2$ or $\mathbb{R}^3$, and $R_0$ is a fix constant in this section. We can take finite $\{\widehat{\boldsymbol{x}}_i\}_{i=1}^n\subset \mathbb{S}^2$ and $\{\chi_{i}\}_{i=1}^n=\{\chi_{\widehat{\boldsymbol{x}}_i}\}_{i=1}^n$ such that $\cup_{i=1}^n\chi_{i}^{-1}(1)$ is a cover of $\mathbb{S}^2$, and we denote $\mathcal{R}_i=\mathcal{R}_{\widehat{\boldsymbol{x}}_i}$. \\
		We define the H\"{o}lder norm on $\mathbb{S}^2$ as the norm on manifold
		\begin{equation}\label{normsph}
			\begin{aligned}
				\|f\|_{ C^\alpha(\mathbb{S}^2)}:=\sum_{i=1}^n\|(f\chi_{i})\circ\mathcal{R}_{i}\circ\mathcal{X}\|_{ C^{\alpha}(\mathbb{R}^2)}.
			\end{aligned}
		\end{equation}
		It is easy to check that this definition is equivalent to the normal definition of H\"{o}lder norm on sphere.\\
		Let $\Phi\in C^\infty(\mathbb{S}^2)$ close to $F_0$ in the sense of $C^1$ that will be fixed later, we can rewrite the equation of $H=F-\Phi$ as
		\begin{equation}\label{3dpesb}
			\partial_t(F-\Phi)(\widehat{\boldsymbol{x}})+\int_{\mathbb{S}^2}G(\Phi(\widehat{\boldsymbol{x}})-\Phi(\widehat{\boldsymbol{y}}))\widetilde\nabla\cdot(J(\widetilde\nabla\Phi) \widetilde\nabla(F-\Phi))(\widehat{\boldsymbol{y}})d\mu_{\mathbb{S}^2}(\widehat{\boldsymbol{y}})=N_1+N_2+N_3:=N(F,\Phi),
		\end{equation}
		with
		\begin{equation}\label{3dpesnlt}
			\begin{aligned}
				&N_1(F,\Phi)=-\int_{\mathbb{S}^2}\big(G(F(\widehat{\boldsymbol{x}})-F(\widehat{\boldsymbol{y}}))-G(\Phi(\widehat{\boldsymbol{x}})-\Phi(\widehat{\boldsymbol{y}}))\big)\widetilde\nabla\cdot(\mathbf{T}(|\widetilde\nabla F|) \widetilde\nabla F)(\widehat{\boldsymbol{y}})d \mu_{\mathbb{S}^2}(\widehat{\boldsymbol{y}}),\\
				&N_2(F,\Phi)=-\int_{\mathbb{S}^2}G(\Phi(\widehat{\boldsymbol{x}})-\Phi(\widehat{\boldsymbol{y}}))\widetilde\nabla\cdot\big(\mathbf{T}(|\widetilde\nabla F|) \widetilde\nabla F-\mathbf{T}(|\widetilde\nabla\Phi|) \widetilde\nabla\Phi-J(\widetilde\nabla\Phi) \widetilde\nabla(F-\Phi)\big)(\widehat{\boldsymbol{y}})d \mu_{\mathbb{S}^2}(\widehat{\boldsymbol{y}}),\\
				&N_3(\Phi)=\int_{\mathbb{S}^2}G(\Phi(\widehat{\boldsymbol{x}})-\Phi(\widehat{\boldsymbol{y}}))\widetilde\nabla\cdot(\mathbf{T}(|\widetilde\nabla\Phi|) \widetilde\nabla\Phi)(\widehat{\boldsymbol{y}})d \mu_{\mathbb{S}^2}(\widehat{\boldsymbol{y}}),\\
				&(J(\widetilde\nabla\Phi)\widetilde\nabla(F-\Phi))^j=\mathbf{T}(|\widetilde\nabla\Phi|)\widetilde\nabla(F-\Phi)^j+\mathbf{T}'(|\widetilde\nabla\Phi|)\frac{\widetilde\nabla\Phi:\widetilde\nabla(F-\Phi)}{|\widetilde\nabla\Phi|}\widetilde\nabla\Phi^j\\
				&\quad\quad:=(J_1(\widetilde\nabla\Phi)\widetilde\nabla(F-\Phi))^j+(J_2(\widetilde\nabla\Phi)\widetilde\nabla(F-\Phi))^j.
			\end{aligned}
		\end{equation}
		We define the following norms and set,
		\begin{equation}
			\begin{aligned}
				&\|F\|_{T}=\sup_{t\in[0,T]}\left(\|F(t)\|_{L^\infty(\mathbb{S}^2)}+t^{m+\kappa}\|F(t)\|_{ C^{m+\kappa}(\mathbb{S}^2)}\right),\\
				&\|F\|_{Z_T}=\sup_{t\in[0,T]}\|F(t)\|_{L^\infty(\mathbb{S}^2)}+\|\widetilde\nabla F\|_{T},\\
				&\mathcal{Z}_{T,\Phi}^\sigma=\left\{F:\sup_{\substack{t<T\\ \widehat{\boldsymbol{x}},\widehat{\boldsymbol{y}}\in\mathbb{S}^2}}\frac{|\widehat{\boldsymbol{x}}-\widehat{\boldsymbol{y}}|}{|F(t,\widehat{\boldsymbol{x}})-F(t,\widehat{\boldsymbol{y}})|}<2\mathbf\Theta_0,\ \|F-\Phi\|_{Z_T}\leq\sigma \right\}.
			\end{aligned}
		\end{equation}
		\begin{theorem}\label{3dpesmthm}
			There exists $\eps_0>0$ such that for any $F_0\in W^{1,\infty}$, if there exists $\Phi\in C^{\infty}(\mathbb{S}^2)$ such that $\|F_0-\Phi\|_{C^1(\mathbb{S}^2)}\leq\eps_0$ and $F_0$ satisfies \eqref{3dpesnc}, then there exists $T=T(\eps_0,\|\Phi\|_{C^{m+4}},\mathbf\Theta_0)<1$ such that the equation \eqref{3dpesb} admits a unique solution with initial data $F_0$ in $[0,T]$, and 
			\begin{equation*}
				\|F-\Phi\|_{Z_T}\lesssim\|F_0-\Phi\|_{C^1(\mathbb{S}^2)}.
			\end{equation*}
		\end{theorem}
		\begin{proof}
			We define the map $\mathcal{S}(Q)=F$ to be the solution of 
			\begin{equation}\label{3dpesm}
				\begin{aligned}
					&\partial_t(F-\Phi)(\widehat{\boldsymbol{x}})+\int_{\mathbb{S}^2}G(\Phi(\widehat{\boldsymbol{x}})-\Phi(\widehat{\boldsymbol{y}}))\widetilde\nabla\cdot(J(\widetilde\nabla\Phi) \widetilde\nabla(F-\Phi))(\widehat{\boldsymbol{y}})d \mu_{\mathbb{S}^2}(\widehat{\boldsymbol{y}})=N(Q,\Phi)(\widehat{\boldsymbol{x}}),\\
					&F(0,\widehat{\boldsymbol{x}})=F_0(\widehat{\boldsymbol{x}}).
				\end{aligned}
			\end{equation}
			What we will prove in this subsection is that $\mathcal{S}$ is a contraction map on $\mathcal{Z}_{T,\Phi}^\sigma$ for some small $\sigma$ and $T$.
			We define the operator 
			\begin{equation}\label{3dpesdefl}
				\mathcal{L}_{\Phi}H(\hat{x})=\int_{\mathbb{S}^2}G(\Phi(\widehat{\boldsymbol{x}})-\Phi(\widehat{\boldsymbol{y}}))\widetilde\nabla\cdot(J(\widetilde\nabla\Phi) \widetilde\nabla H)(\widehat{\boldsymbol{y}})d \mu_{\mathbb{S}^2}(\widehat{\boldsymbol{y}}).
			\end{equation}
			So the equation \eqref{3dpesm} can be written as 
			\begin{equation*}
				\begin{aligned}
					&\partial_t H+\mathcal{L}_{\Phi}H=N(Q,\Phi),\\
					&H(0,\widehat{\boldsymbol{x}})=H_0(\widehat{\boldsymbol{x}}),
				\end{aligned}
			\end{equation*}
			where we will denote $H=F-\Phi$ in the following of the proof. We want to apply \eqref{pehmfd} in Theorem \ref{thmmani} with $n=1$, $F=0$, $G=N(Q,\Phi)$. Now what we need to prove is the condition \eqref{mfdcdl}.
			\begin{lemma}\label{3dpesmlem}
				There exists $\mathcal{L}_{\mathbb{R}^2}^i$ satisfying \eqref{defop} and \eqref{condop}, such that for $l=0,m$, we have
				\begin{equation}\label{3dpesml}
					\|(\chi_i\mathcal{L}_{\Phi}H)\circ\mathcal{R}_i\circ\mathcal{X}-\mathcal{L}^i_{\mathbb{R}^2}((\chi_iH)\circ\mathcal{R}_i\circ\mathcal{X})\|_{\dot C^{l+\kappa}}\lesssim \|H\|_{C^{l+\kappa+1-\zeta_0}},
				\end{equation}
				for some $\zeta_0\in(0,1)$. 
			\end{lemma}
			\begin{proof}
				Since $\|F_0-\Phi\|_{C^1(\mathbb{S}^2)}\leq \eps_0$, by \eqref{3dpesnc} and taking $\eps_0<\frac{1}{100(C+10)\mathbf{\Theta}_0}$ with some $C$ satisfying $\sup_{\widehat{\boldsymbol{x}},\widehat{\boldsymbol{y}}\in\mathbb{S}^2}\frac{|F(\widehat{\boldsymbol{x}})-F(\widehat{\boldsymbol{y}})|}{|\widehat{\boldsymbol{x}}-\widehat{\boldsymbol{y}}|}\leq C\|F\|_{C^1(\mathbb{S}^2)}$, we can see that
				\begin{equation}\label{3dpesPc}
					\inf_{\widehat{\boldsymbol{x}},\widehat{\boldsymbol{y}}\in\mathbb{S}^2}\frac{|\Phi(\widehat{\boldsymbol{x}})-\Phi(\widehat{\boldsymbol{y}})|}{|\widehat{\boldsymbol{x}}-\widehat{\boldsymbol{y}}|}\geq \mathbf{\Theta}_0^{-1}-2C\|F_0-\Phi\|_{C^1(\mathbb{S}^2)}\geq\frac{4\mathbf{\Theta}_0^{-1}}{5}. 
				\end{equation}
				Recall that we can split the integration region of \eqref{3dpesdefl} into the union of many balls $\cup_{i=1}^nB(\widehat{\boldsymbol{x}}_i,a_1)$, so we consider a small neighborhood of $\widehat{\boldsymbol{x}}_i$, so we define
				\begin{equation*}
					\mathcal{L}_{\Phi,i}H=\int_{\mathbb{S}^2\cap B(\hat{\boldsymbol{x}}_i,a_2)}\widetilde\nabla G(\Phi(\widehat{\boldsymbol{x}})-\Phi(\widehat{\boldsymbol{y}}))\cdot(J(\widetilde\nabla\Phi) \widetilde\nabla H)(\widehat{\boldsymbol{y}})d \mu_{\mathbb{S}^2}(\widehat{\boldsymbol{y}}).
				\end{equation*}
				We denote $H_i=H\chi_i$, $h_i=H_i\circ\mathcal{R}_i\circ\mathcal{X}$, $\phi_i=(\tilde{\chi}_i\Phi)\circ\mathcal{R}_i\circ\mathcal{X}$. Now we approximate $\chi_i(\mathcal{L}_\Phi H)$ by $\tilde\chi_i\mathcal{L}_{\Phi,i}H_i$, and denote the error term as
				\begin{equation*}
					N_4=\chi_i(\mathcal{L}_\Phi H)-\tilde\chi_i\mathcal{L}_{\Phi,i}H_i,\quad \tilde{N}_4=N_4\circ\mathcal{R}_i\circ\mathcal{X}.
				\end{equation*}
				So we only need to estimate $\tilde\chi_i\mathcal{L}_{\Phi,i}H_i$. We consider $\left(\widetilde\nabla G\cdot(J(\widetilde\nabla\Phi)\widetilde\nabla(H\chi_i))\right)\circ\mathcal{R}_i\circ\mathcal{X}$. By changing variable formula, generally for any $F:\mathbb{S}^2\rightarrow \mathbb{R}^3$, and $f=F\circ\mathcal{X}:\mathbb{R}^2\rightarrow\mathbb{R}^3$, we have
				\begin{equation}\label{3dpescvar}
					\int_{\mathbb{S}^2}F(\widehat{\boldsymbol{y}})d \mu_{\mathbb{S}^2}(\widehat{\boldsymbol{y}})=\int_{\mathbb{R}^2}\left(\frac{2}{1+|\eta|^2}\right)^2f(\eta)d\eta
				\end{equation}
				By classical calculus, 
				\begin{equation}\label{3dpescvar2}
					\begin{aligned}
						&\left(\widetilde\nabla G\cdot(J(\widetilde\nabla\Phi)\widetilde\nabla(H\chi_i))\right)\circ\mathcal{R}_i\circ\mathcal{X}=(\frac{1+|\theta|^2}{2})^2\frac{\partial\mathcal{X}^j}{\partial\theta_l}(\frac{\partial G}{\partial\theta_l}\circ\mathcal{R}_i\circ\mathcal{X})(J(\widetilde\nabla\Phi)\widetilde\nabla(H\chi_i))^j\circ\mathcal{R}_i\circ\mathcal{X}\\
						&\quad=\left(\frac{1+|\theta|^2}{2}\right)^2(\frac{\partial G}{\partial\theta_l}\circ\mathcal{R}_i\circ\mathcal{X})\left(\mathbf{T}(|\widetilde\nabla\Phi|\circ\mathcal{R}_i\circ\mathcal{X})\frac{\partial h_i^j}{\partial\theta_l}+\frac{1+|\theta|^2}{2}\mathbf{T}'(|\widetilde\nabla\Phi|\circ\mathcal{R}_i\circ\mathcal{X})\frac{\nabla\phi_i:\nabla h_i}{|\nabla\phi_i|}\frac{\partial\phi_i^j}{\partial\theta_l} \right)\\
						&\quad=\left(\frac{1+|\theta|^2}{2}\right)^2(\frac{\partial G}{\partial\theta_l}\circ\mathcal{R}_i\circ\mathcal{X})\left(\mathbf{T}(\frac{1+|\theta|^2}{2}|\nabla\phi_i|)\frac{\partial h_i^j}{\partial\theta_l}+\frac{1+|\theta|^2}{2}\mathbf{T}'(\frac{1+|\theta|^2}{2}|\nabla\phi_i|)\frac{\nabla\phi_i:\nabla h_i}{|\nabla\phi_i|}\frac{\partial\phi_i^j}{\partial\theta_l} \right).
					\end{aligned}
				\end{equation}
				We define 
				\begin{equation*}
					\begin{aligned}
						&\tilde{\mathbf{T}}(\theta,\beta)=\mathbf{T}(\rho(\theta)\beta)=\frac{\tilde{\mathcal{T}}(\theta,\beta)}{\rho(\theta)\beta}=\frac{\mathcal{T}(\rho(\theta)\beta)}{\rho(\theta)\beta},\\
						&\left(\tilde{J}(\eta,\nabla\phi_i)\nabla h_i\right)^j(\eta):=\tilde{\mathbf{T}}(\eta,\nabla\phi_i)\nabla h_i^j+\frac{\partial\tilde{\mathbf{T}}(\eta,\beta)}{\partial\beta}\Big|_{\beta=|\nabla\phi_i|}\frac{\nabla\phi_i:\nabla h_i}{|\nabla\phi_i|}\nabla\phi_i^j.
					\end{aligned}
				\end{equation*}
				for some smooth $\rho$ bounded up and below, such that $\rho(\theta)\equiv \frac{1+|\theta|^2}{2}$ for $|\theta|\leq 4R_0$. Then by \eqref{3dpescvar} and \eqref{3dpescvar2} we can write 
				\begin{equation*}
					\begin{aligned}
						&\left(\int_{\mathbb{S}^2}\widetilde\nabla G(\Phi_i(\widehat{\boldsymbol{x}})-\Phi_i(\widehat{\boldsymbol{y}}))\cdot(J((\widetilde\nabla\Phi))\widetilde\nabla H_i)(\widehat{\boldsymbol{y}})d\mu_{\mathbb{S}^2}(\widehat{\boldsymbol{y}})\right)\circ\mathcal{R}_{\widehat{\boldsymbol{x}}_i}\circ\mathcal{X}\\
						&\quad=\int_{\mathbb{R}^2} G(\phi_i(\theta)-\phi_i(\eta))\nabla\cdot \left(\tilde{J}(\eta,\nabla\phi_i)\nabla h_i(\eta)\right)d\eta,
					\end{aligned}
				\end{equation*}
				So the main term of the symbol of $\nabla\cdot \left(\tilde{J}(\eta,\nabla\phi_i)\nabla h_i(\eta)\right)$ is
				\begin{equation}\label{3dpesopJ}
					\frac{\tilde{\mathcal{T}}(\eta,|\nabla\phi_i|)}{\rho(\eta)|\nabla\phi_i|}\left(|\xi|^2\mathrm{Id}_3-\frac{(\nabla\phi_i\cdot\xi)\otimes(\nabla\phi_i\cdot\xi)}{|\nabla\phi_i|^2}\right)+\frac{\partial\tilde{\mathcal{T}}(\eta,\beta)}{\partial\beta}\Big|_{\beta=|\nabla\phi_i|}\frac{(\nabla\phi_i\cdot\xi)\otimes(\nabla\phi_i\cdot\xi))}{|\nabla\phi_i|^2},
				\end{equation}
				which is uniformly elliptic. So we approximate $G(\phi_i(\theta)-\phi_i(\eta))$ by $G(\nabla\phi_i(\theta)(\theta-\eta))$
				\begin{equation}\label{3dpesdefn5}
					\begin{aligned}
						&\left(\tilde\chi_i\mathcal{L}_{\Phi,i}H_i\right)\circ\mathcal{R}_i\circ\mathcal{X}=\left(\tilde\chi_i\circ\mathcal{R}_i\circ\mathcal{X}\right)(\theta)\int_{\mathbb{R}^2} G(\phi_i(\theta)-\phi_i(\eta))\nabla\cdot \left(\tilde{J}(\eta,\nabla\phi_i)\nabla h_i(\eta)\right)d\eta\\
						&\quad=\left(\tilde\chi_i\circ\mathcal{R}_i\circ\mathcal{X}\right)(\theta)\int_{\mathbb{R}^2} G(\nabla\phi_i(\theta)(\theta-\eta))\nabla\cdot \left(\tilde{J}(\eta,\nabla\phi_i)\nabla h_i(\eta)\right)d\eta+\tilde{N}_5,
					\end{aligned}
				\end{equation}
				with error term defined as
				\begin{equation*}
					\tilde{N}_5=(\tilde\chi_i\circ\mathcal{R}_i\circ\mathcal{X})(\theta)\int_{\mathbb{R}^2}\left(G(\nabla\phi_i(\theta)(\theta-\eta))-G(\phi_i(\theta)-\phi_i(\eta))\right)\nabla\cdot(\tilde{J}(\eta,\nabla\phi_i)\nabla h_i)(\eta)d\eta.
				\end{equation*}
				Now we need to estimate the first term of \eqref{3dpesdefn5}. For $G(\nabla\phi_i(\theta)(\theta-\eta))$, we denote
				\begin{align*}
					A(\theta)=\nabla\phi_i(\theta)\in \mathbb{M}_{3\times 2}(\mathbb{R}),\quad B(\theta)=\sqrt{A^T(\theta)A(\theta)}\in \mathbb{M}_{2\times 2}(\mathbb{R}),\quad Q(\theta)=A(\theta)B^{-1}(\theta)\in \mathbb{M}_{3\times 2}(\mathbb{R}),
				\end{align*}
				where $\mathbb{M}$ means matrix.\\
				By fundamental calculus, we know that 
				\begin{equation*}
					\begin{aligned}
						&\left|A\eta \right|^2=(A\eta)^TA\eta=\eta^T(A^TA)\eta=\left|B\eta \right|^2,\\
						&(A\eta)\otimes (A\eta)=(A\eta)(A\eta)^T=A\eta\eta^TA^T=QB\eta\eta^TB^TQ^T=Q(B\eta)\otimes(B\eta)Q^T,
					\end{aligned}
				\end{equation*}
				and
				\begin{equation*}
					\mathcal{F}\left(\frac{1}{|\theta|}\right)(\xi)=\frac{2\pi}{|\xi|},\quad\mathcal{F}\left(\frac{\theta_i\theta_j}{|\theta|^3}\right)(\xi)=2\pi(\frac{\delta_{ij}}{|\xi|}-\frac{\xi_i\xi_j}{|\xi|^3}).
				\end{equation*}
				So by changing variable, we can write $\mathcal{F}(G(\nabla\phi_i(\theta)\cdot))(\xi)$ as
				\begin{equation}\label{3dpesopG}
					\begin{aligned}
						\frac{2\pi\mathrm{Id}_3}{|\det B(\theta)||B^{-1}\xi|}+\frac{2\pi}{|\det B(\theta)|} Q\left(\frac{\mathrm{Id}_3}{|B^{-1}\xi|}-\frac{(B^{-1}\xi)\otimes( B^{-1}\xi)}{|B^{-1}\xi|^3}\right) Q^T=2\pi\frac{\mathrm{Id}_3+v(\theta,\xi)\otimes v(\theta,\xi)}{|\det B(\theta)||B^{-1}\xi|}
					\end{aligned}
				\end{equation}
				with $v(\theta,\xi)=Q\begin{pmatrix}
					0&1\\
					-1&0
				\end{pmatrix}\frac{B^{-1}\xi}{|B^{-1}\xi|}$.\\
				We define $\tilde \phi_i\in C^\infty(\mathbb{R}^2;\mathbb{R}^{3\times 2})$ as the expansion of $\nabla\phi_i$ by
				\begin{align*}
					&\tilde \phi_i(\theta)\equiv \nabla\phi_i(\theta),\quad \text{for}\ |\theta|\leq 4R_0,\\
					&0<C_1<|\tilde \phi_i|<C_2<\infty,\\
					&\| {\tilde{\phi}_i}\|_{C^{m+2}}\leq C_3\|\Phi\|_{C^{m+2}}<\infty.
				\end{align*}
				Combining \eqref{3dpesopJ} and \eqref{3dpesopG}, we define 
				\begin{equation}\label{3dpesdefop}
					\mathcal{L}_{\mathbb{R}^2}^if(\theta)=\int_{\mathbb{R}^2}\frac{(\mathrm{Id}_3+\tilde v(\theta,\xi)\otimes \tilde v(\theta,\xi))\big(\mathcal{ J}(\tilde\phi_i)(\theta,\xi)\big)}{|\det B(\theta)||B^{-1}\xi|}\hat{f}(\xi)e^{i\theta\cdot\xi}d\xi,
				\end{equation}
				with
				\begin{equation*}
					\begin{aligned}
						&\mathcal{J}(\tilde\phi_i)(\theta,\xi)=\frac{\mathcal{T}(\rho(\theta)|\tilde\phi_i|)}{\rho(\theta)|\tilde\phi_i|}(|\xi|^2\mathrm{Id}_3-\frac{(\tilde\phi_i\cdot\xi)\otimes (\tilde\phi_i\cdot\xi)}{|\tilde\phi_i|^2|}) +\mathcal{T}'(\rho(\theta)|\tilde\phi_i|)\frac{(\tilde\phi_i\cdot\xi)\otimes (\tilde\phi_i\cdot\xi)}{|\tilde\phi_i|^2},\\     
						&\tilde A(\theta)=\tilde{\phi}_i(\theta),\quad \tilde B(\theta)=\sqrt{\tilde A(\theta)\tilde A^T(\theta)},\quad \tilde Q(\theta)=\tilde A(\theta)\tilde B^{-1}(\theta), \quad \tilde v(\theta,\xi)=\tilde Q\begin{pmatrix}
							0&1\\
							-1&0
						\end{pmatrix}\frac{\tilde B^{-1}\xi}{|\tilde B^{-1}\xi|}.
					\end{aligned}
				\end{equation*}
				So we have
				\begin{equation*}
					\left(\tilde\chi_i\circ\mathcal{R}_i\circ\mathcal{X}\right)(\theta)\int_{\mathbb{R}^2} G(\nabla\phi_i(\theta)(\theta-\eta))\nabla\cdot \left(\tilde{J}(\eta,\nabla\phi_i)\nabla h_i(\eta)\right)d\eta=\mathcal{L}_{\mathbb{R}^2}^if(\theta)+\tilde{N}_6,
				\end{equation*}
				with error term $\tilde N_6$
				\begin{equation*}
					\tilde{N}_6=(\tilde\chi_i\circ\mathcal{R}_i\circ\mathcal{X})(\theta)\int_{\mathbb{R}^2}G(\nabla\phi_i(\theta)(\theta-\eta))\nabla\cdot(\tilde{J}(\eta,\nabla\phi_i)\nabla h_i)(\eta)d\eta-\mathcal{L}_{\mathbb{R}^2}^ih_i,
				\end{equation*}
				We define
				\begin{equation}
					\tilde{N}_0=\sum_{i=4}^6\tilde{N}_i=(\chi_i\mathcal{L}_{\Phi}H)\circ\mathcal{R}_i\circ\mathcal{X}-M(\tilde\phi_i,\Lambda,\theta)((\chi_iH)\circ\mathcal{R}_i\circ\mathcal{X}),
				\end{equation}
				then \eqref{3dpesml} is equivalent to
				\begin{equation*}
					\|\tilde{N}_0\|_{ C^{l+\kappa}}\lesssim\|H\|_{C^{l+\kappa+1-\zeta_0}},\quad l=0,m,
				\end{equation*}
				for some small $\zeta_0\in(0,1)$. By Lemma \ref{maininterpo}, we only need to estimate $\|\tilde N_0\|_{C^{n}}$ for $n\leq m+1$.\\
				For $\tilde{N}_4$, by definition of H\"{o}lder norm on $\mathbb{S}^2$ and interpolation, we only need to estimate the $C^1$ norm of
				\begin{equation*}
					\begin{aligned}
						N_4&=\tilde{\chi}_i(\widehat{\boldsymbol{x}})\int_{\mathbb{S}^2}(\chi_i(\widehat{\boldsymbol{x}})-\chi_i(\widehat{\boldsymbol{y}}))\nabla_{\mathbb{S}^2}(G(\Phi(\widehat{\boldsymbol{x}})-\Phi(\widehat{\boldsymbol{y}})))\cdot(J(\nabla_{\mathbb{S}^2}\Phi) \nabla_{\mathbb{S}^2}H)(\widehat{\boldsymbol{y}})d \mu_{\mathbb{S}^2}(\widehat{\boldsymbol{y}})\\
						&\quad\quad\quad\quad\quad-\tilde{\chi}_i(\widehat{\boldsymbol{x}})\int_{\mathbb{S}^2}\nabla_{\mathbb{S}^2}G(\Phi(\widehat{\boldsymbol{x}})-\Phi(\widehat{\boldsymbol{y}}))\cdot(HJ(\nabla_{\mathbb{S}^2}\Phi) \nabla_{\mathbb{S}^2}\chi_i(\widehat{\boldsymbol{y}}))(\widehat{\boldsymbol{y}})d \mu_{\mathbb{S}^2}(\widehat{\boldsymbol{y}})\\
						&=N_{41}+N_{42}.
					\end{aligned}
				\end{equation*}
				By \eqref{3dpesPc}, we have 
				\begin{equation*}
					\begin{aligned}
						&|\widetilde\nabla_{\widehat{\boldsymbol{x}},\widehat{\boldsymbol{y}}}\left(G(\Phi(\widehat{\boldsymbol{x}})-\Phi(\widehat{\boldsymbol{y}}))\right)|\lesssim\frac{\|\Phi\|_{C^3}(1+\|\Phi\|_{C^3})\|\chi_i\|_{C^2}}{|\widehat{\boldsymbol{x}}-\widehat{\boldsymbol{y}}|^2},\\
						&|\widetilde\nabla_{\widehat{\boldsymbol{x}}}\widetilde\nabla_{\widehat{\boldsymbol{y}}}\left(G(\Phi(\widehat{\boldsymbol{x}})-\Phi(\widehat{\boldsymbol{y}}))\right)|\lesssim\frac{\|\Phi\|_{C^3}(1+\|\Phi\|_{C^3})^2\|\chi_i\|_{C^2}}{|\widehat{\boldsymbol{x}}-\widehat{\boldsymbol{y}}|^3},\\
						&|(\widetilde\nabla_{\widehat{\boldsymbol{x}}}+\widetilde\nabla_{\widehat{\boldsymbol{y}}})\left(G(\Phi(\widehat{\boldsymbol{x}})-\Phi(\widehat{\boldsymbol{y}})\right)|\lesssim\frac{\|\Phi\|_{C^3}(1+\|\Phi\|_{C^3})}{|\widehat{\boldsymbol{x}}-\widehat{\boldsymbol{y}}|}.
					\end{aligned}
				\end{equation*}
				Integral by parts, and we can easily check that
				\begin{equation*}
					\|N_4(t)\|_{L^\infty(\mathbb{S}^2)}\lesssim \|\Phi\|_{C^3}(1+\|\Phi\|_{C^3})\|H(t)\|_{C^1},
				\end{equation*}
				Then we estimate $C^1$ norm. For $N_{41}$,
				\begin{equation*}
					\begin{aligned}
						&\|\widetilde\nabla_{\hat{\boldsymbol{x}}}N_{41}\|_{L^\infty}\lesssim\|\tilde{\chi}_i\|_{C^2}\|\widetilde\nabla H\|_{L^\infty}\\
						&\quad+\Big\|\int_{\mathbb{S}^2}\widetilde\nabla\tilde\chi_i(\widehat{\boldsymbol{x}})\widetilde\nabla_{\hat{\boldsymbol{y}}}G(\Phi(\widehat{\boldsymbol{x}})-\Phi(\widehat{\boldsymbol{y}}))\cdot(J(\widetilde\nabla\Phi) \widetilde\nabla H(\widehat{\boldsymbol{y}})-J(\widetilde\nabla\Phi) \widetilde\nabla H(\widehat{\boldsymbol{x}}))d \mu_{\mathbb{S}^2}(\widehat{\boldsymbol{y}})\Big\|_{L^\infty}\\
						&\quad+\Big\|\int_{\mathbb{S}^2}(\tilde\chi_i(\widehat{\boldsymbol{x}})-\tilde\chi_i(\widehat{\boldsymbol{y}}))\widetilde\nabla_{\hat{\boldsymbol{x}}}\widetilde\nabla_{\widehat{\boldsymbol{y}}}G(\Phi(\widehat{\boldsymbol{x}})-\Phi(\widehat{\boldsymbol{y}}))\cdot(J(\widetilde\nabla\Phi) \widetilde\nabla H(\widehat{\boldsymbol{y}})-J(\widetilde\nabla\Phi) \widetilde\nabla H(\widehat{\boldsymbol{x}}))d \mu_{\mathbb{S}^2}(\widehat{\boldsymbol{y}})\Big\|_{L^\infty}\\
						&\quad+\Big\|\int_{\mathbb{S}^2}(\widetilde\nabla_{ \widehat{\boldsymbol{x}}}+\widetilde\nabla_{ \widehat{\boldsymbol{y}}})G(\Phi(\widehat{\boldsymbol{x}})-\Phi(\widehat{\boldsymbol{y}}))\cdot\widetilde\nabla^2\tilde\chi_i(\widehat{\boldsymbol{y}})\left(J(\widetilde\nabla\Phi) \widetilde\nabla H(\widehat{\boldsymbol{x}})\right)d \mu_{\mathbb{S}^2}(\widehat{\boldsymbol{y}})\Big\|_{L^\infty}\\
						&\lesssim\|\Phi\|_{C^3}(1+\|\Phi\|_{C^3})\|\chi_i\|_{C^2}\|H\|_{C^{1+\eps}(\mathbb{S}^2)}.
					\end{aligned}
				\end{equation*}
				For $N_{42}$, we have
				\begin{equation*}
					\begin{aligned}
						&\|\nabla_{\mathbb{S}^2,\hat{\boldsymbol{x}}}N_{42}\|_{L^\infty}\lesssim\|\tilde{\chi}_i\|_{C^2}\|\nabla_{\mathbb{S}^2}H\|_{L^\infty}\\
						&\quad+\Big\|\int_{\mathbb{S}^2}(\widetilde\nabla_{\hat{\boldsymbol{y}}}+\widetilde\nabla_{\hat{\boldsymbol{x}}})G(\Phi(\widehat{\boldsymbol{x}})-\Phi(\widehat{\boldsymbol{y}}))\widetilde\nabla_{\widehat{\boldsymbol{y}}}\cdot(HJ(\widetilde\nabla\Phi) \widetilde\nabla\tilde\chi_i(\hat{\boldsymbol{y}}))d \mu_{\mathbb{S}^2}(\widehat{\boldsymbol{y}})\Big\|_{L^\infty}\\
						&\quad+\Big\|\int_{\mathbb{S}^2}\widetilde\nabla_{\widehat{\boldsymbol{y}}}G(\Phi(\widehat{\boldsymbol{x}})-\Phi(\widehat{\boldsymbol{y}}))\left(\widetilde\nabla\cdot(HJ(\widetilde\nabla\Phi) \widetilde\nabla\tilde\chi_i)(\widehat{\boldsymbol{y}})-\widetilde\nabla\cdot(HJ(\widetilde\nabla\Phi) \widetilde\nabla\tilde\chi_i)(\widehat{\boldsymbol{x}})\right)d \mu_{\mathbb{S}^2}(\widehat{\boldsymbol{y}})\Big\|_{L^\infty}\\
						&\lesssim\|\Phi\|_{C^3}(1+\|\Phi\|_{C^3})\|\chi_i\|_{C^2}\|H\|_{C^{1+\eps}(\mathbb{S}^2)}.
					\end{aligned}
				\end{equation*}
				So we have
				\begin{equation}\label{3dpesn4l}
					\|N_4(t)\|_{C^1(\mathbb{S}^2)}\lesssim \|\Phi\|_{C^3}(1+\|\Phi\|_{C^3})\|\tilde\chi\|_{C^2} \|H(t)\|_{C^{1+\eps}(\mathbb{S}^2)}.
				\end{equation}
				For higher order derivatives, we use the decomposition on both $\mathbb{S}^2$ and $\mathbb{R}^2$
				\begin{equation}\label{3dpeskerd}
					\begin{aligned}
						&\nabla_{\mathbb{S}^2,\widehat{\boldsymbol{x}}}^k=\sum_{l=0}^kC(k,l)(\nabla_{\mathbb{S}^2,\widehat{\boldsymbol{x}}}+\nabla_{\mathbb{S}^2,\widehat{\boldsymbol{y}}})^l\nabla_{\mathbb{S}^2,\widehat{\boldsymbol{y}}}^{k-l},\\
						&\nabla_{\theta}^k=\sum_{l=0}^kC(k,l)(\nabla_{\theta}+\nabla_{\eta})^l\nabla_{\eta}^{k-l}.
					\end{aligned}
				\end{equation}
				Then use \eqref{3dpesPc}, one can see that for any $l$, 
				\begin{equation*}
					\begin{aligned}
						&|(\widetilde\nabla_{\widehat{\boldsymbol{x}}}+\widetilde\nabla_{\widehat{\boldsymbol{y}}})^lG(\Phi(\widehat{\boldsymbol{x}})-\Phi(\widehat{\boldsymbol{y}}))|\lesssim \frac{\|\Phi\|_{C^{l+2}}(1+\|\Phi\|_{C^{l+2}})^{l+1}}{|\widehat{\boldsymbol{x}}-\widehat{\boldsymbol{y}}|},\\
						&|\widetilde\nabla_{\widehat{\boldsymbol{y}}}(\widetilde\nabla_{\widehat{\boldsymbol{x}}}+\widetilde\nabla_{\widehat{\boldsymbol{y}}})^lG(\Phi(\widehat{\boldsymbol{x}})-\Phi(\widehat{\boldsymbol{y}}))|\lesssim \frac{\|\Phi\|_{C^{l+3}}(1+\|\Phi\|_{C^{l+3}})^{l+2}}{|\widehat{\boldsymbol{x}}-\widehat{\boldsymbol{y}}|^2},\\
						&|\widetilde\nabla_{\widehat{\boldsymbol{x}}}\widetilde\nabla_{\widehat{\boldsymbol{y}}}(\widetilde\nabla_{\widehat{\boldsymbol{x}}}+\widetilde\nabla_{\widehat{\boldsymbol{y}}})^lG(\Phi(\widehat{\boldsymbol{x}})-\Phi(\widehat{\boldsymbol{y}}))|\lesssim \frac{\|\Phi\|_{C^{l+3}}(1+\|\Phi\|_{C^{l+3}})^{l+2}}{|\widehat{\boldsymbol{x}}-\widehat{\boldsymbol{y}}|^3}.
					\end{aligned}
				\end{equation*}
				For $\widetilde\nabla^kN_{41}$ with $k\geq 2$, we have
				\begin{equation*}
					\begin{aligned}
						&\|\widetilde\nabla^{k}N_{41}\|_{L^\infty}\\
						&\quad\lesssim \sum_{m_1\leq k}\Big\|\widetilde\nabla^{m_1}\tilde\chi_i(\widehat{\boldsymbol{x}})\int_{\mathbb{S}^2}(\tilde\chi_i(\widehat{\boldsymbol{y}})-\tilde\chi_i(\widehat{\boldsymbol{x}}))\widetilde\nabla_{\widehat{\boldsymbol{x}}}^{k-m_1}\widetilde\nabla_{\widehat{\boldsymbol{y}}}(G(\Phi(\widehat{\boldsymbol{x}})-\Phi(\widehat{\boldsymbol{y}})))\cdot(J(\widetilde\nabla\Phi) \widetilde\nabla H)(\widehat{\boldsymbol{y}})d \mu_{\mathbb{S}^2}(\widehat{\boldsymbol{y}})\Big\|_{L^\infty_{\widehat{\boldsymbol{x}}}}\\
						&\quad\quad+\sum_{\substack{m_2\neq 0\\
								m_1+m_2\leq k}}\Big\|\widetilde\nabla^{m_1}\tilde\chi_i(\widehat{\boldsymbol{x}})\widetilde\nabla^{m_2}\tilde\chi_i(\widehat{\boldsymbol{x}})\int_{\mathbb{S}^2}\widetilde\nabla_{\widehat{\boldsymbol{x}}}^{k-m_1-m_2}\widetilde\nabla_{\widehat{\boldsymbol{y}}}(G(\Phi(\widehat{\boldsymbol{x}})-\Phi(\widehat{\boldsymbol{y}})))\cdot(J(\widetilde\nabla\Phi) \widetilde\nabla H)(\widehat{\boldsymbol{y}})d \mu_{\mathbb{S}^2}(\widehat{\boldsymbol{y}})\Big\|_{L_{\hat{\boldsymbol{x}}}^\infty}\\
						&\quad=\mathcal{N}_{411}+\mathcal{N}_{412}
					\end{aligned}
				\end{equation*}
				For $\mathcal{N}_{412}$, we can see that
				\begin{equation*}
					\begin{aligned}
						&\mathcal{N}_{412}\lesssim\sum_{\substack{m_2\neq 0\\
								m_1+m_2\leq k\\
								l\leq k-m_1-m_2}}\Big\|\nabla_{\mathbb{S}^2}^{m_1}\tilde\chi_i(\widehat{\boldsymbol{x}})\widetilde\nabla^{m_2}\tilde\chi_i(\hat{\boldsymbol{x}})\int_{\mathbb{S}^2}(\widetilde\nabla_{ \widehat{\boldsymbol{x}}}+\widetilde\nabla_{ \widehat{\boldsymbol{y}}})^{k-m_1-m_2-l}\widetilde\nabla_{ \widehat{\boldsymbol{y}}}(G(\Phi(\widehat{\boldsymbol{x}})-\Phi(\widehat{\boldsymbol{y}}))) \\
						&\quad\quad\quad\quad\quad \cdot\left(\widetilde\nabla^l(J(\widetilde\nabla\Phi) \widetilde\nabla H)(\widehat{\boldsymbol{y}})-\widetilde\nabla^l(J(\widetilde\nabla\Phi) \widetilde\nabla H)(\widehat{\boldsymbol{x}})\right)d \mu_{\mathbb{S}^2}(\widehat{\boldsymbol{y}}) \Big\|_{L_{\widehat{\boldsymbol{x}}}^\infty}\\
						&\lesssim\|\tilde\chi_i\|_{C^{k}}(1+\|\tilde\chi_i\|_{C^{k}})\|\Phi\|_{C^{k}}(1+\|\Phi\|_{C^{k}})^{k}\|H\|_{C^{k+\eps}}.
					\end{aligned}
				\end{equation*}
				For $\mathcal{N}_{411}$, by \eqref{3dpeskerd}, we have
				\begin{equation*}
					\begin{aligned}
						&\mathcal{N}_{411}\lesssim\|\tilde\chi_i\|_{C^k}(1+\|\tilde\chi_i\|_{C^k})\|H\|_{C^{1}}\\
						&+\sum_{\substack{m_1\leq k-1\\
								l\leq k-m_1-1}}\|\tilde\chi_i\|_{C^{m_1}}\Big\|\int_{\mathbb{S}^2}(\tilde\chi_i(\widehat{\boldsymbol{x}})-\tilde\chi_i(\widehat{\boldsymbol{y}}))\widetilde\nabla_{ \widehat{\boldsymbol{y}}}\widetilde\nabla_{ \widehat{\boldsymbol{x}}}(\widetilde\nabla_{ \widehat{\boldsymbol{x}}}+\widetilde\nabla_{ \widehat{\boldsymbol{y}}})^{k-m_1-1-l}G(\Phi(\widehat{\boldsymbol{x}})-\Phi(\widehat{\boldsymbol{y}}))\\
						&\quad\quad\quad\quad\quad\cdot\left(\widetilde\nabla^{l}(J(\widetilde\nabla\Phi)\widetilde\nabla H)(\widehat{\boldsymbol{y}})-\widetilde\nabla^{l}(J(\widetilde\nabla\Phi)\widetilde\nabla H)(\widehat{\boldsymbol{x}})\right)d\mu_{\mathbb{S}^2}(\widehat{\boldsymbol{y}})\Big\|_{L^\infty_{\widehat{\boldsymbol{x}}}}\\
						&+\sum_{\substack{m_1\leq k-1\\
								1\leq m_2\leq k-m_1\\
								l\leq k-m_1-m_2}}\|\tilde\chi_i\|_{C^{m_1}}\Big\|\int_{\mathbb{S}^2}\widetilde\nabla_{ \widehat{\boldsymbol{x}}}^{m_2}\tilde\chi_i(\widehat{\boldsymbol{x}})\widetilde\nabla_{ \widehat{\boldsymbol{y}}}(\widetilde\nabla_{ \widehat{\boldsymbol{x}}}+\widetilde\nabla_{ \widehat{\boldsymbol{y}}})^{k-m_1-m_2-l}G(\Phi(\widehat{\boldsymbol{x}})-\Phi(\widehat{\boldsymbol{y}}))\\
						&\quad\quad\quad\quad\quad\cdot\left(\widetilde\nabla^{l}(J(\widetilde\nabla\Phi)\widetilde\nabla H)(\widehat{\boldsymbol{y}})-\widetilde\nabla^{l}(J(\widetilde\nabla\Phi)\widetilde\nabla H)(\widehat{\boldsymbol{x}})\right)d\mu_{\mathbb{S}^2}(\widehat{\boldsymbol{y}})\Big\|_{L^\infty_{\widehat{\boldsymbol{x}}}}\\
						&+\sum_{\substack{m_1\leq k-1\\
								l\leq k-m_1-1}}\|\tilde\chi_i\|_{C^{m_1}}\Big\|\int_{\mathbb{S}^2}(\widetilde\nabla_{ \widehat{\boldsymbol{x}}}+\widetilde\nabla_{ \widehat{\boldsymbol{y}}})^{k-m_1-l}G(\Phi(\widehat{\boldsymbol{x}})-\Phi(\widehat{\boldsymbol{y}}))\cdot\left(\widetilde\nabla^2\tilde\chi_i(\widehat{\boldsymbol{y}})\widetilde\nabla^{l}(J(\widetilde\nabla\Phi)\widetilde\nabla H)(\widehat{\boldsymbol{x}})\right)d\mu_{\mathbb{S}^2}(\widehat{\boldsymbol{y}})\Big\|_{L^\infty_{\widehat{\boldsymbol{x}}}}\\
						&+\sum_{\substack{m_1\leq k-2\\
								l\leq k-m_1-1\\
								1\leq j\leq l}}\|\tilde\chi_i\|_{C^{m_1}}\Big\|\int_{\mathbb{S}^2}\widetilde\nabla_{ \widehat{\boldsymbol{y}}}(\widetilde\nabla_{ \widehat{\boldsymbol{x}}}+\widetilde\nabla_{ \widehat{\boldsymbol{y}}})^{k-m_1-l}G(\Phi(\widehat{\boldsymbol{x}})-\Phi(\widehat{\boldsymbol{y}}))\\
						&\quad\quad\quad\quad\quad\cdot\left(\widetilde\nabla^j\tilde\chi_i\widetilde\nabla^{l-j}(J(\widetilde\nabla\Phi)\widetilde\nabla H)(\widehat{\boldsymbol{y}})-\widetilde\nabla^j\tilde\chi_i\widetilde\nabla^{l-j}(J(\widetilde\nabla\Phi)\widetilde\nabla H)(\widehat{\boldsymbol{x}})\right)d\mu_{\mathbb{S}^2}(\widehat{\boldsymbol{y}})\Big\|_{L^\infty_{\widehat{\boldsymbol{x}}}}\\
						&\lesssim\|\tilde\chi_i\|_{C^{k}}(1+\|\tilde\chi_i\|_{C^{k}})\|\Phi\|_{C^{k}}(1+\|\Phi\|_{C^{k}})^{k}\|H\|_{C^{k+\eps}}
					\end{aligned}
				\end{equation*}
				So we have proved
				\begin{equation*}
					\|N_{41}\|_{C^{m+1}}\lesssim\|\tilde\chi_i\|_{C^{m+2}}(1+\|\tilde\chi_i\|_{C^{m+2}})\|\Phi\|_{C^{m+2}}(1+\|\Phi\|_{C^{k}})^{k}\|H\|_{C^{m+1+\eps}}.
				\end{equation*}
				For $\widetilde\nabla^kN_{42}$, we have
				\begin{equation*}
					\begin{aligned}
						&\|\widetilde\nabla^kN_{42}\|_{L^\infty}\lesssim\|\tilde\chi_i\|_{C^k}(1+\|\tilde\chi_i\|_{C^k})\|\Phi\|_{C^k}(1+\|\Phi\|_{C^k})\|H\|_{C^{1}}\\
						&+\sum_{\substack{m_1\leq k-1\\
								l\leq k-m_1-1}}\left\|\widetilde\nabla^{m_1}\tilde{\chi}_i(\widehat{\boldsymbol{x}})\int_{\mathbb{S}^2}(\widetilde\nabla_{ \widehat{\boldsymbol{x}}}+\widetilde\nabla_{ \widehat{\boldsymbol{y}}})^{k-m_1-l}G(\Phi(\widehat{\boldsymbol{x}})-\Phi(\widehat{\boldsymbol{y}}))\cdot\left( \widetilde\nabla^{l+1}((HJ(\widetilde\nabla\Phi) \widetilde\nabla\chi_i))(\widehat{\boldsymbol{y}})\right)d \mu_{\mathbb{S}^2}(\widehat{\boldsymbol{y}})\right\|_{L^\infty}\\
						&+\sum_{m_1\leq k-1}\Big\|\widetilde\nabla^{m_1}\tilde{\chi}_i(\widehat{\boldsymbol{x}})\int_{\mathbb{S}^2}\widetilde\nabla_{ \widehat{\boldsymbol{y}}}G(\Phi(\widehat{\boldsymbol{x}})-\Phi(\widehat{\boldsymbol{y}})) \\
						&\quad\quad\quad\quad\quad \cdot\left(\widetilde\nabla^{k-m_1}((HJ(\widetilde\nabla\Phi) \widetilde\nabla\chi_i))(\widehat{\boldsymbol{y}})-\widetilde\nabla^{k-m_1}((HJ(\widetilde\nabla\Phi) \widetilde\nabla\chi_i))(\widehat{\boldsymbol{x}})\right)d \mu_{\mathbb{S}^2}(\widehat{\boldsymbol{y}})\Big\|_{L^\infty}\\
						&\lesssim \|\tilde\chi_i\|_{C^{k+1}}(1+\|\tilde\chi_i\|_{C^{k+1}})\|\Phi\|_{C^{k+1}}(1+\|\Phi\|_{C^{k+1}})^{k+1}\|H\|_{C^{k+\eps}}.
					\end{aligned}
				\end{equation*}
				So we can get that
				\begin{equation}\label{3dpesn4h}
					\|N_4\|_{C^{m+1}(\mathbb{S}^2)}\lesssim \|\Phi\|_{C^{m+2}}(1+\|\Phi\|_{C^{m+2}})^{m+1}\|\tilde\chi_i\|_{C^{k+1}}(1+\|\tilde\chi_i\|_{C^{k+1}})\|H(t)\|_{C^{m+1+\eps}(\mathbb{S}^2)}.
				\end{equation}
				For $\tilde{N}_5$, we also use interpolation, and estimate $ C^1$ norm. For $L^\infty$ norm, it is easy to check that
				\begin{equation*}
					\|\tilde{N}_5\|_{L^\infty}\lesssim \|\Phi\|_{C^2}\|\nabla h_i\|_{L^\infty}.
				\end{equation*}
				For $C^1$ norm, we shortly denote $\chi_i^R(\theta)=(\chi_i\circ\mathcal{R}_i\circ\mathcal{X})(\theta)$, $\tilde\chi_i^R(\theta)=(\tilde\chi_i\circ\mathcal{R}_i\circ\mathcal{X})(\theta)$ and $\tilde{G}(\theta,\eta)=G(\phi_i(\theta)-\phi_i(\eta))-G(\nabla\phi_i(\theta)(\theta-\eta))$. By \eqref{3dpesPc}, we can see that 
				\begin{align*}
					|\tilde\chi_i^R(\theta)\tilde\chi_i^R(\eta)\nabla_{\theta}\nabla_{\eta}\tilde{G}(\theta,\eta)|\lesssim \frac{\|\nabla\Phi\|_{C^2}}{|\theta-\eta|^2}.
				\end{align*}
				So 
				\begin{equation}\label{3dpesn5l}
					\begin{aligned}
						\|\tilde{N}_5\|_{\dot C^1}\lesssim& \|\tilde\chi_i^R\|_{C^1}\|\Phi\|_{C^3}(1+\|\Phi\|_{C^3})^2\|h_i\|_{C^1}\\ &+\Big\|\tilde\chi_i^R(\theta)\int_{\mathbb{R}^2} \tilde\chi_i^R(\eta)\nabla_\theta\nabla_\eta \tilde{G}(\theta,\eta)\cdot\left(\tilde{J}(\eta,\nabla\phi_i)\nabla h_i(\eta) \right)d\eta\Big\|_{L^\infty}\\
						\lesssim& \|\tilde\chi_i^R\|_{C^1}(1+\|\tilde\chi_i^R\|_{C^1})^2\|\Phi\|_{C^3}(1+\|\Phi\|_{C^3})^2\|h_i\|_{C^1}\\
						&+\Big\|\tilde\chi_i^R(\theta)\int_{\mathbb{R}^2} \tilde\chi_i^R(\eta)\nabla_\theta\nabla_\eta \tilde{G}(\theta,\eta)\cdot\left(\tilde{J}(\eta,\nabla\phi_i)\nabla h_i(\eta)-\tilde{J}(\theta,\nabla\phi_i)\nabla h_i(\theta) \right)d\eta\Big\|_{L^\infty}\\
						\lesssim&\|\tilde\chi_i\|_{C^1}(1+\|\tilde\chi_i\|_{C^1})^2\|\Phi\|_{C^3}(1+\|\Phi\|_{C^3})^2\|h_i\|_{C^{1+\eps}}.
					\end{aligned}
				\end{equation}
				For higher order derivatives, notice that 
				\begin{align*}
					|\tilde\chi_i^R(\theta)\tilde\chi_i^R(\eta)\nabla_{\theta}\nabla_{\eta}(\nabla_{\theta}+\nabla_{\eta})^l\tilde{G}(\theta,\eta)|\lesssim \frac{\|\nabla\Phi\|_{C^{l+2}}(1+\|\nabla\Phi\|_{C^{l+2}})^{l+1}}{|\theta-\eta|^2}.
				\end{align*}
				We apply \eqref{3dpeskerd} and integral by parts, then similar to \eqref{3dpesn4h}, for $k\geq 2$, we can get
				\begin{equation*}
					\begin{aligned}
						&\|\nabla^k\tilde{N}_5\|_{L^\infty}\lesssim \|\tilde\chi_i\|_{C^{k+1}}(1+\|\tilde\chi_i\|_{C^{k+1}})^{k+1}\|\Phi\|_{C^{3}}\|h_i\|_{C^1}\\
						&\quad+\sum_{\substack{2\leq l\leq k\\
								m\leq l-2}}\|\tilde\chi_i\|_{C^{k-l}}\Big\|\int_{\mathbb{R}^2}\nabla_{\theta}\nabla_\eta(\nabla_\theta+\nabla_\eta)^{l-m-1}\tilde{G}(\theta,\eta)\cdot\nabla^{m}(\tilde{J}(\eta,\nabla\phi_i)\nabla h_i)(\eta)d\eta\Big\|_{L^\infty_\theta}\\
						&\quad+\Big\|\int_{\mathbb{R}^2}\nabla_{\theta}\nabla_\eta\tilde{G}(\theta,\eta)\cdot\left(\nabla^{k-1}(\tilde{J}(\eta,\nabla\phi_i)\nabla h_i)(\theta)-\nabla^{k-1}(\tilde{J}(\eta,\nabla\phi_i)\nabla h_i)(\eta)\right)\tilde\chi_i(\eta)d\eta\Big\|_{L^\infty_\theta}\\
						&\quad+\Big\|\int_{\mathbb{R}^2}\nabla_{\theta}\tilde{G}(\theta,\eta)\cdot\nabla_\eta\tilde\chi_i(\eta)\nabla^{k-1}(\tilde{J}(\eta,\nabla\phi_i)\nabla h_i)(\theta)d\eta\Big\|_{L^\infty_\theta}\\
						&\lesssim \|\tilde\chi_i\|_{C^{k+1}}(1+\|\tilde\chi_i\|_{C^{k+1}})^{k+1}\|\Phi\|_{C^{k+2}}(1+\|\Phi\|_{C^{k+2}})^{k+2}\|h_i\|_{C^{k+\eps}},
					\end{aligned}
				\end{equation*}
				which gives 
				\begin{equation}\label{3dpesn5h}
					\|\tilde{N}_5(t)\|_{ C^{m+1}}\lesssim \|\tilde\chi_i\|_{C^{k+1}}(1+\|\tilde\chi_i\|_{C^{k+1}})^{k+1}\|\Phi\|_{C^{m+2}}(1+\|\Phi\|_{C^{m+2}})^{m+2}\|\nabla h_i(t)\|_{C^{m+\eps}}.
				\end{equation}
				For $\tilde{N}_6$, we can write $\tilde N_6$ as
				\begin{equation*}
					\begin{aligned}
						\tilde{N}_6=& \int_{\mathbb{R}^2}(G(\tilde \phi_i(\theta)(\theta-\eta)))\nabla\cdot\left(\big(\tilde{J}(\eta,\tilde\phi_i)-\tilde{J}(\theta,\tilde\phi_i)\big)\nabla h_i(\eta)\right)d\eta\\
						&+\int_{\mathbb{R}^2}(1-\tilde{\chi}_i\circ\mathcal{R}_i\circ\mathcal{X})(\theta)G(\tilde \phi_i(\theta)(\theta-\eta))\nabla\cdot\left(\tilde{J}(\eta,\tilde\phi_i)\nabla h_i(\eta)\right)d\eta\\
						=&\tilde{N}_{61}+\tilde{N}_{62}.
					\end{aligned}
				\end{equation*}
				For $\tilde{N}_{61}$ we use similar methods as $\tilde{N}_4$ and $\tilde{N}_5$
				to get
				\begin{equation*}
					\begin{aligned}
						&\|\tilde{N}_{61}\|_{C^1}\lesssim\|\tilde\chi_i\|_{C^{2}}(1+\|\tilde\chi_i\|_{C^{2}})^{2}\|\Phi\|_{C^{3}}(1+\|\Phi\|_{C^{3}})^{3}\|h_i\|_{C^{1+\eps}}\\
						&\|\tilde{N}_{61}\|_{C^k}\lesssim\|\tilde\chi_i\|_{C^{k+1}}(1+\|\tilde\chi_i\|_{C^{k+1}})^{k+1}\|\Phi\|_{C^{k+2}}(1+\|\Phi\|_{C^{k+2}})^{k+2}\|h_i\|_{C^{k+\eps}}.
					\end{aligned}
				\end{equation*}
				For $\tilde{N}_{62}$, we notice that $\text{Supp}(1-\tilde{\chi}_i\circ\mathcal{R}_i\circ\mathcal{X})\subset \{|\theta|\geq 3R_0\}$, and $\text{Supp}(\nabla h_i)\subset \{|\eta|\leq 2R_0\}$, so $\tilde{N}_{62}$ has no singularity. Precisely, we have
				\begin{equation*}
					\left|\nabla_{\theta}^l\left(G(\tilde \phi_i(\theta)(\theta-\eta))\right)\right|\lesssim(1+\|\Phi\|_{C^{l+2}})^{k+2} R_0^{-l},\quad \forall |\theta|\geq 3R_0,\ |\eta|\leq 2R_0.
				\end{equation*}
				So it is easy to see that
				\begin{equation*}
					\|\tilde{N}_{62}\|_{C^k}\lesssim \|\tilde\chi_i\|_{C^{k+1}}(1+\|\tilde\chi_i\|_{C^{k+1}})^{k+1}\|\Phi\|_{C^{k+2}}(1+\|\Phi\|_{C^{k+2}})^{k+2}R_0^{-l}\|h_i\|_{C^1}.
				\end{equation*}
				So we can conclude that
				\begin{equation}\label{3dpesn6}
					\begin{aligned}
						& \|\tilde N_6(t)\|_{C^1}\lesssim \|\tilde\phi_i\|_{C^3}(1+\|\tilde\phi_i\|_{C^3})^2\|\nabla h_i\|_{C^\eps},\\
						&\|\tilde N_6(t)\|_{C^{m+1}}\lesssim \|\tilde\phi_i\|_{C^{m+3}}(1+\|\tilde\phi_i\|_{C^{m+3}})^{m+2}\|\nabla h_i(t)\|_{C^{m+\eps}}.
					\end{aligned}
				\end{equation}
				Now combine \eqref{3dpesn4l}, \eqref{3dpesn4h}, \eqref{3dpesn5l}, \eqref{3dpesn5h} and \eqref{3dpesn6} to get \eqref{3dpesml}.
			\end{proof}\\
			We continue to prove Theorem \ref{3dpesmthm}. We will apply \eqref{pehmfd} of Theorem \ref{thmmani} by taking $\mathcal{L}_{\mathcal{M}}^{s}=\mathcal{L}_\Phi$, $\mathcal{L}_{\mathbb{R}^d}^{s,j}$ as \eqref{3dpesdefop}, $F=0$, $G=N(Q,\Phi)$, and by Lemma \ref{3dpesmlem}, for some small $T$ we can get
			\begin{equation*}
				\sup_{t\in[0,T]}(\|(F-\Phi)(t)\|_{C^1(\Omega)}+t^{m+\kappa}\|(F-\Phi)(t)\|_{C^{m+1+\kappa}(\Omega)})\leq C_0\|F_0-\Phi\|_{W^{1,\infty}}+C_0\da'_{T,1}(0,N(Q,\Phi)).
			\end{equation*}
			Apply Lemma \ref{3dpeslemn1}, Lemma \ref{3dpeslemn2} and Lemma \ref{3dpeslemn3}, we can see that
			\begin{equation}\label{3dpeshn}
				\begin{aligned}
					&\|N(Q,\Phi)\|_{L_T^1L^\infty}+\sup_{t\in[0,T]}t^{\kappa}\|N(Q,\Phi)(t)\|_{ C^{\kappa}}+\sup_{t\in[0,T]}t^{m+\kappa}\|N(Q,\Phi)(t)\|_{ C^{m+\kappa}}\\
					&\leq C_1(1+\|\Phi\|_{C^{1}})(\|Q-\Phi\|_{Z_T}+T^{\kappa}\|\Phi\|_{C^{m+3}(\mathbb{S}^2)})(1+\|Q-\Phi\|_{Z_T}+T\|\Phi\|_{C^{m+3}(\mathbb{S}^2)})^{m+1}.
				\end{aligned}
			\end{equation}
			Hence, since $Q\in\mathcal{Z}_{T,\Phi}^\sigma$, we can take $\eps_0<\frac{1}{100C_0C_1(1+\|F_0\|_{W^{1,\infty}})}$, $\sigma<2C_0\eps_0$ and $T<\frac{\sigma}{100C_1(1+\|\Phi\|_{C^{m+3}})^{m+1}}$ small enough to get
			\begin{align*}
				\|F-\Phi\|_{Z_T}\leq C_0\eps_0+ \frac{\sigma}{4}\leq\frac{\sigma}{2}.
			\end{align*}
			Furthermore, by \eqref{3dpesPc} one has
			\begin{equation*}
				\begin{aligned}
					\inf_{\substack{t<T\\ \widehat{\boldsymbol{x}},\widehat{\boldsymbol{y}}\in\mathbb{S}^2}}\frac{|F(t,\widehat{\boldsymbol{x}})-F(t,\widehat{\boldsymbol{y}})|}{|\widehat{\boldsymbol{x}}-\widehat{\boldsymbol{y}}|}&\geq \inf_{ \widehat{\boldsymbol{x}},\widehat{\boldsymbol{y}}\in\mathbb{S}^2}\frac{|\Phi(\widehat{\boldsymbol{x}})-\Phi(\widehat{\boldsymbol{y}})|}{|\widehat{\boldsymbol{x}}-\widehat{\boldsymbol{y}}|}-2C\sup_{t\in[0,T]}\|\nabla_{\mathbb{S}^2}(F-\Phi)\|_{L^\infty}(t)\\
					&\geq \frac{4}{3\mathbf{\Theta}_0}-2C^2\sigma.
				\end{aligned}
			\end{equation*}
			Hence 
			\begin{align*}
				\sup_{\substack{t<T\\ \widehat{\boldsymbol{x}},\widehat{\boldsymbol{y}}\in\mathbb{S}^2}}\frac{|\widehat{\boldsymbol{x}}-\widehat{\boldsymbol{y}}|}{|F(t,\widehat{\boldsymbol{x}})-F(t,\widehat{\boldsymbol{y}})|} &\leq \left(\frac{4}{3\mathbf{\Theta}_0}-2C^2\sigma\right)^{-1}\leq 2\mathbf{\Theta}_0.
			\end{align*}
			So we have proved that $\mathcal{S}:\mathcal{Z}_{T,\Phi}^\sigma\rightarrow\mathcal{Z}_{T,\Phi}^\sigma$. Now we only need to prove that $\mathcal{S}$ has the compression property. For $Q_1,Q_2\in\mathcal Z_{T,\Phi}^\sigma$, we denote $\mathbf{H}=H_1-H_2$, $\mathbf{Q}=Q_1-Q_2$, $\mathbf{N}=N(Q_1,\Phi)-N(Q_2,\Phi)$ and $\mathbf{N}_i=N_i(Q_1,\Phi)-N_i(Q_2,\Phi)$ for $i=1,2$. We can write the equation of $\mathbf{H}$ as 
			\begin{equation*}
				\begin{aligned}
					&\partial_t\mathbf{H}(\widehat{\boldsymbol{x}})+\int_{\mathbb{S}^2}G(\Phi(\widehat{\boldsymbol{x}})-\Phi(\widehat{\boldsymbol{y}}))\widetilde\nabla\cdot(J(\widetilde\nabla\Phi)\widetilde\nabla\mathbf{H})(\widehat{\boldsymbol{y}})d \mu_{\mathbb{S}^2}(\widehat{\boldsymbol{y}})=\mathbf{N}(\widehat{\boldsymbol{x}}),\\
					&\mathbf{H}(0,\widehat{\boldsymbol{x}})=0.
				\end{aligned}
			\end{equation*}
			Similarly we have
			\begin{equation*}
				\begin{aligned}
					\sup_{t\in[0,T]}(\|\mathbf{H}(t)\|_{C^1(\Omega)}+t^{m+\kappa}\|\mathbf{H}(t)\|_{C^{m+1+\kappa}(\Omega)})\lesssim \da'_{T,1}(0,\mathbf{N}).
				\end{aligned}
			\end{equation*}
			By Lemma \ref{3dpeslemn1} and Lemma \ref{3dpeslemn2}, we have
			\begin{equation}\label{3dpesnh}
				\begin{aligned}
					&\|\mathbf N\|_{L_T^1L^\infty}+\sup_{t\in[0,T]}t^{\kappa}\|\mathbf{N}(t)\|_{\dot C^{\kappa}}+\sup_{t\in[0,T]}t^{m+\kappa}\|\mathbf N(t)\|_{\dot C^{m+\kappa}}\\
					&\lesssim  C_2\|\Phi\|_{C^{m+2}}(1+\|\Phi\|_{C^{m+2}})^{m+1}\|\mathbf{Q}\|_{Z_T}(\|Q-\Phi\|_{Z_T}+T^{\kappa}\|\Phi\|_{C^3})^{m},
				\end{aligned}
			\end{equation}
			and we can take $\sigma<\frac{1}{2^{10m}C_2(1+\|\Phi\|_{C^{m+1}})^{10m}}$ and $T<\frac{\sigma}{100(1+\|\Phi\|_{C^{m+3}})^{10m}}$ be small enough to ensure that $\mathcal S:\mathcal{Z}_{T,\Phi}^\sigma\rightarrow\mathcal{Z}_{T,\Phi}^\sigma$, and by \eqref{3dpesnh}, one has
			\begin{equation}
				\|\mathcal{S}(Q_1)-\mathcal{S}(Q_2)\|_{Z_T}\leq\frac{1}{2}\|\mathbf{Q}\|_{Z_T},
			\end{equation}
			and by compression mapping theorem, we get the result.
		\end{proof}\\

		\subsection{Thin-Film equation}
		In this section, we consider the following equation
		\begin{equation}\label{eqtf}
			\begin{aligned}
				\partial_tu(t,x)=\Delta(e^{-\Delta u(t,x)})&,\quad \operatorname{in}\ \mathbb{R}^+\times\mathbb{R}^d,\\
				u|_{t=0}=u_0&,\quad \operatorname{in}\ \mathbb{R}^d,
			\end{aligned}
		\end{equation}
		which can be used to describe the evolution of crystal film. This problem has been focused in mathematics and physics since Burton, Cabrera and Frank studied crystal growth on solid surfaces in \cite{BCF51}. Some researchers concentrated on the precise motion prediction for crystal surface and developed related PDEs, for which can be seen in \cite{KDM95} and \cite{MW13}. \\
		Following the derivation of \cite{LR18}, the equation can be explained by the following. Take $u$ to be the height of the film, and we consider the dynamics of atom deposition, detachment and diffusion on a crystal surface in the epitaxy growth process. By Fick's law we have
		\begin{equation*}
			u_t+\nabla \cdot J=0,\quad J=-D_c\nabla\rho_c,
		\end{equation*}
		where $D_c$ is the surface diffusion constant and $\rho_c$ is the equilibrium density of adatoms on a substrate of the thin film. Up to a normalization constant, $\rho_c$ can be described as $\rho_c=e^{-\frac{E_c-\mu_c}{k_BT}}$, where $E_c$ is the energy of per adatom, $\mu_c$ is the chemical potential per adatom, $k_B$ is the Boltzmann constant and $T$ is the temperature. We absorb $e^{-\frac{E_c}{k_BT}}$ and the normalization constant into a reference density $\rho_c$, and it becomes $\rho_c=\varrho_0e^{\frac{\mu_c}{k_BT}}$. \\
		Define the energy functional $E(u)=\frac{1}{p}\int |\nabla u|^pdx$, and $\mu_c=\frac{\delta E}{\delta u}=\Delta_pu$, where $\Delta_p$ is the $p$-Laplacian operator. Then we can recast the equation into the following Cahn-Hilliard equation with curvature-dependent mobility as in \cite{LR18}
		\begin{equation*}
			u_t=\operatorname{div}\mathcal{M}\nabla\frac{\delta E}{\delta u},\quad\mathcal{M}=e^{-\Delta_pu},
		\end{equation*}
		and then up to a normalization constant, the equation becomes 
		\begin{equation*}
			u_t=\Delta(e^{-\Delta_pu}).
		\end{equation*}
		It should be mentioned that we will only consider the case when $p=2$ in this paper. In \cite{GLL19}, Gao, Liu and Lu proved the existence and uniqueness for global solutions by gradient flow methods, details can be seen in \cite{CZF13}. In \cite{LR18}, Liu and Strain proved the global existence, uniqueness and optimal large time decay rates of \eqref{eqtf} for initial data that is of medium size in the critical Wiener algebra $u\in \dot {\mathcal{F}}^{2,1}$. 
		In this paper, we can prove local well-posedness for large data in $\{u|\Delta u\in L^\infty\}$. Moreover, if $\|\Delta u\|_{L^\infty}$ is small, we obtain global existence, uniqueness and optimal large time decay rates. We remark that $\dot {\mathcal{F}}^{2,1}\subsetneqq\{u|\Delta u\in L^\infty\}$.  Precisely, we have the following result.
		\begin{theorem}\label{thmtf}~Fix $\kappa$ such that $0<4-\kappa\ll1$, and $m\in\mathbb{N}$. There exists $\varepsilon_0>0$ such that  \\
			i) For any initial data $u_0$ with $\|\Delta u_0\|_{L^\infty}\leq \eps_0$, there exists a unique solution $u$ to the Cauchy problem \eqref{eqtf} such that
			$$
			\sup_{t>0}(\|\Delta  u(t)\|_{L^\infty}+t^{\frac{m+\kappa}{4}}\| \Delta u(t)\|_{\dot{C}^{m+\kappa}})\lesssim\|\Delta u_0\|_{L^\infty}.
			$$
			ii) For any initial data $u_0$ with $\Delta u_0\in L^\infty$, if there exists $\eps_1>0$ such that
			\begin{equation}\label{contf}
				\|\Delta u_0-\Delta u_0\ast \rho_{\eps_1}\|_{L^\infty}\leq\eps_0,
			\end{equation}
			then there exists $T=T(\eps_0,\eps_1,\|\Delta u_0\|_{L^\infty})>0$ such that the Cauchy problem \eqref{eqtf} admits a unique solution $u$ in $[0,T]$ satisfying
			$$
			\sup_{t\in[0,T]}(\| \Delta u(t)\|_{L^\infty}+t^{\frac{m+\kappa}{4}}\| \Delta u(t)\|_{\dot{C}^{m+\kappa}})\lesssim \|\Delta u_0\|_{L^\infty}.
			$$
		\end{theorem}
		\begin{proof}
			We define the following norm
			$$
			\|u\|_{T}=\sup_{t\in[0,T]}(\| u(t)\|_{L^\infty}+t^{\frac{m+\kappa}{4}}\|  u(t)\|_{\dot{C}^{m+\kappa}}).
			$$
			By Lemma \ref{maininterpo}, it is easy to check that for any $0<\gamma<\frac{m+\kappa}{4}$, 
			$$
			\sup_{t\in[0,T]}t^\frac{\gamma}{4}\| u(t)\|_{\dot C^{\gamma}}\lesssim \|u\|_{T}.
			$$
			Let $\sigma,T>0$, $\phi\in C^\infty(\mathbb{R}^d)$
			that will be fixed later. Define 
			\begin{align*}
				\mathcal{X}^\sigma_{T,\phi}=\left\{u:u|_{t=0}=u_0,\Delta u\in L^\infty_TL^\infty, \|\Delta(u-\phi)\|_T\leq \sigma \right\}.
			\end{align*}
			We  write the equation as follows
			\begin{align*}
				&\partial_tu+e^{-\Delta\phi}\Delta^2u=e^{-\Delta u}|\nabla\Delta u|^2+(e^{-\Delta\phi}-e^{-\Delta u})\Delta^2 u,\\
				&u|_{t=0}=u_0.
			\end{align*}
			For any $v\in   \mathcal{X}^\sigma_{T,\phi}$, we define a map $\mathcal{S}:v\rightarrow u$, where $u$ is the unique solution to the Cauchy problem
			\begin{align*}
				& \partial_tu+e^{-\Delta\phi}\Delta^2u=e^{-\Delta v}|\nabla\Delta v|^2+(e^{-\Delta\phi}-e^{-\Delta v})\Delta^2 v,\ \ \text{in}\ (0,T)\times \mathbb{R}^d,\\
				&u|_{t=0}=u_0.
			\end{align*}
			Take $\Delta$ on both sides of the equation above, and denote $\tilde\phi=\Delta\phi$,  $f=\Delta u-\tilde\phi$, $g=\Delta v-\tilde\phi$, then one has
			\begin{align*}
				&\partial_t f+e^{-\tilde\phi}\Delta^2f=\Delta R_1[g,\phi]+R_2[f,\phi],\\
				&f|_{t=0}=\Delta u_0-\tilde\phi,
			\end{align*}
			where
			\begin{align*}
				R_1[g,\phi]=&e^{- g+\tilde\phi}|\nabla( g+\tilde\phi)|^2+e^{-\tilde\phi}(1-e^{- g})\Delta(g+\tilde\phi)-e^{-\tilde\phi}\Delta\tilde\phi\\
				R_2[f,\phi]=&-\nabla(e^{-\tilde\phi})\cdot \nabla \Delta f-\Delta(e^{-\tilde\phi})\Delta f.
			\end{align*}
			We prove that, with suitable choices of $\phi,T,\sigma$, the map $\mathcal{S}$ is a contraction map from $\mathcal{X}_{T,\phi}^{\sigma}$ to itself.\\
			By Theorem \ref{thmexis} for $F=R_1[g,\phi]$ and $G=R_2[f,\phi]$, there exists $T=T(\|\Delta\phi\|_{C^{m+4}})>0$ such that 
			\begin{equation}\label{apl}
				\begin{aligned}
					\|f\|_T
					&\lesssim\|\Delta u_0-\tilde\phi\|_{L^\infty}+\sup_{t\in[0,T]}(t^{\frac{\kappa}{4}}\|R_1[g,\phi](t)\|_{\dot{C}^{\kappa-2}}+t^{\frac{m+\kappa}{4}}\|R_1[g,\phi](t)\|_{\dot{C}^{m+\kappa-2}})\\
					&\quad\quad+\|R_2[f,\phi]\|_{L^1_TL^\infty}+\sup_{t\in[0,T]}t^{\frac{m}{4}+1}\|R_2[f,\phi](t)\|_{\dot{C}^{m}}.
				\end{aligned}
			\end{equation}
			Now we estimate the right hand side nonlinear terms $R_1[g,\phi]$ and $R_2[f,\phi]$. We remark that the main term is $e^{-\tilde\phi}(1-e^{-g})\Delta g$ appeared in $R_1[g,\phi]$, while all the other terms are lower order terms because the regular function $\phi$ shared derivatives. It is easy to check that 
			\begin{align*}
				&\sup_{t\in[0,T]}\left(t^\frac{\kappa}{4}\|R_1[g,\phi](t)\|_{\dot{C}^{\kappa-2}}+t^{\frac{m+\kappa}{4}}\|R_1[g,\phi](t)\|_{\dot{C}^{m+\kappa-2}}\right)\\
				&\quad\quad\quad\quad\quad\quad\lesssim \exp(\|\Delta\phi\|_{L^\infty}+\|g\|_{L^\infty})\left(\|g\|_T^2+\mathbf{1}_{\phi\not\equiv 0}T^\frac{1}{4}(1+\|\Delta\phi\|_{C^{m+4}})^{m+2}\right)(1+\|g\|_T)^m,
			\end{align*}
			and 
			\begin{align*}
				&\|R_2[f,\phi]\|_{L^1_TL^\infty}d\tau+\sup_{t\in[0,T]}t^{\frac{m}{4}+1}\|R_2[f,\phi](t)\|_{\dot{C}^{m+4}}\\
				&\quad\quad\quad\quad\lesssim \mathbf{1}_{\phi\not\equiv 0}T^\frac{1}{4}\exp(\|\Delta\phi\|_{L^\infty})(1+\|f\|_T+\|\Delta\phi\|_{C^{m+4}})^{m+2}.
			\end{align*}
			Here we use the notation $\mathbf{1}_{\phi\not\equiv 0}$, which indicates that term vanishes if $\phi\equiv 0$. 
			
			We consider the following two cases, respectively.\\
			$\bullet$ If $\|\Delta u_0\|_{L^\infty} \leq \varepsilon_0$,  take $\phi\equiv 0$. By Theorem \ref{lemmain}, the estimate \eqref{apl} holds for $T=\infty$. One has 
			\begin{align*}
				\|f\|_\infty\leq C(\|\Delta(u_0-\phi)\|_{L^\infty}+\|g\|_T^2(1+\|g\|_T)^m)\leq C(\varepsilon_0+\sigma^2(1+\sigma)^m).
			\end{align*}
			This implies that $\|f\|_\infty\leq \sigma$ by taking $\varepsilon_0$ small enough and $\sigma=\sigma_1(\varepsilon_0)$ small enough. Hence $\mathcal{S}v\in \mathcal{X}_{\infty,0}^{\sigma_1}$ provided $v\in \mathcal{X}_{\infty,0}^{\sigma_1}$.\\ 
			$\bullet$ If	\eqref{contf} holds for $0<\varepsilon_0\ll 1$ that will be fixed later, then take $\phi=u_0\ast \rho_{\eps_1}$. We have $\|\Delta(u_0-\phi)\|_{L^\infty}\leq 2\varepsilon_0$.  By Theorem \ref{lemmain}, the estimate \eqref{apl} holds for some $T=T(\|\Delta u_0\|_{L^\infty},\eps_1)>0$. One has 
			\begin{align*}
				\|f\|_T
				&\leq C\|\Delta(u_0-\phi)\|_{L^\infty}+C\exp(2\|\Delta u_0\|_{L^\infty})\left(\|g\|_T^2+T^\frac{1}{4}(1+\|\Delta\phi\|_{C^{m+4}})^{m+2}\right)(1+\|g\|_T+\|f\|_T)^{m+2}\\
				&\leq C\varepsilon_0+C\exp(2\|\Delta u_0\|_{L^\infty})(\sigma^2+T^\frac{1}{4}(1+\eps_1^{-(m+4)}\|\Delta u_0\|_{L^\infty})^{m+2})(1+\sigma)^m.
			\end{align*}
			This implies $\|g\|_{T^*}\leq \sigma$ by taking $\varepsilon_0$ small enough, $\sigma=\sigma_2(\varepsilon_0,\|\Delta u_0\|_{L^\infty})$ small enough, and finally $T^*=T^*(\varepsilon_0, \eps_1, \|\Delta u_0\|_{L^\infty})$ small enough. Hence $\mathcal{S} v\in \mathcal{X}^{\sigma_2}_{T^*,u_0\ast \rho_{\eps_1}}$ provided $ v\in \mathcal{X}^{\sigma_2}_{T^*,u_0\ast \rho_{\eps_1}}$.\vspace{0.3cm}\\
			It remains to prove that $\mathcal{S}$ is a contraction map in both settings. For $v_1,v_2\in\mathcal{X}_{T,\phi}^\sigma$, 
			denote $u_i=\mathcal{S}v_i$, $f_i=\Delta (v_i-\phi)$, $g_i=\Delta(u_i-\phi)$ for $i=1,2$, and $\mathbf{g}=g_1-g_2$, $\mathbf{f}=f_1-f_2$, then we have
			\begin{align*}
				&\partial_t\mathbf{g}+e^{-\Delta\phi}\Delta^2\mathbf{g}=\Delta \mathbf{R}_1+\mathbf{R}_2,\\
				&\mathbf{g}|_{t=0}=0,
			\end{align*}
			where 
			\begin{equation*}
				\begin{aligned}   
					\mathbf{R}_1=&R_1[g_1,\phi]-R_1[g_2,\phi]\\
					=&e^{-(g_1+\tilde\phi)}\nabla(g_1+g_2+2\tilde\phi):\nabla\mathbf{g}+e^{-(g_2+\tilde\phi)}(e^{-\mathbf{g}}-1)|\nabla(g_2+\tilde\phi)|^2\\
					&\quad\quad+e^{-(g_2+\tilde\phi)}(1-e^{-\mathbf{g}})\Delta(g_1+\tilde\phi)+e^{-\tilde\phi}(1-e^{-g_2})\Delta\mathbf{g},\\
					\mathbf{R}_2=&R_2[f_1,\phi]-R_2[f_2,\phi]=\nabla(e^{-\tilde\phi})\nabla\Delta\mathbf{f}-\Delta(e^{-\tilde\phi})\Delta\mathbf{f}.
				\end{aligned}
			\end{equation*}
			By Theorem \ref{thmexis} with $F=\mathbf{R}_1$ and $G=\mathbf{R}_2$, there exists \footnote{The existence time may be different with the one in \eqref{apl},  we use the symbol $T$ to denote the smaller existence time.} $T=T(\|\Delta\phi\|_{C^{m+4}})>0$ such that
			\begin{equation*}
				\|\mathbf{g}\|_{T}\lesssim \sup_{t\in[0,T]}(t^{\frac{\kappa}{4}}\|\mathbf{R}_1(t)\|_{\dot C^{\kappa-2}}+t^{\frac{m+\kappa}{4}}\|\mathbf{R}_1(t)\|_{\dot C^{m+\kappa-2}})+\|\mathbf{R}_2\|_{L^1_TL^\infty}+\sup_{t\in[0,T]}t^{\frac{m}{4}+1}\|\mathbf{R}_2(t)\|_{\dot C^{m}}.
			\end{equation*}
			Note that 
			\begin{align*}
				&\sup_{t\in[0,T]}(t^\frac{\kappa}{4}\|\mathbf{R}_1(t)\|_{\dot C^{\kappa-2}}+t^{\frac{m+\kappa}{4}}\|\mathbf{R}_1(t)\|_{\dot C^{m+\kappa-2}})\\
				&\quad\quad\quad\lesssim \|\mathbf{g}\|_T\left(\|(g_1,g_2)\|_T+\mathbf{1}_{\phi\not\equiv 0}T^\frac{1}{4}(1+\|\Delta\phi\|_{C^{m+4}})^{m+2}\right)(1+\|(g_1,g_2)\|_T)^m,
			\end{align*}
			and 
			\begin{align*}
				&\|\mathbf{R}_2\|_{L^1_TL^\infty}+\sup_{t\in[0,T]}t^{\frac{m}{4}}\|\mathbf{R}_2(t)\|_{\dot C^{m}}\lesssim \mathbf{1}_{\phi\not\equiv 0}T^\frac{1}{4}\|\mathbf{f}\|_T(1+\|\Delta\phi\|_{C^{m+4}})^{m+2}.
			\end{align*}
			Hence, 
			\begin{align*}
				\|\mathbf{g}\|_T&\leq C (\|\mathbf{f}\|_T+\|\mathbf{g}\|_T)\left(\|(g_1,g_2)\|_T+\mathbf{1}_{\phi\not\equiv 0}T^\frac{1}{4}(1+\|\Delta\phi\|_{C^{m+4}})^{m+2}\right)(1+\|(g_1,g_2)\|_T)^m\\
				&\leq C (\|\mathbf{f}\|_T+\|\mathbf{g}\|_T)(\sigma+\mathbf{1}_{\phi\not\equiv 0}T^\frac{1}{4}(1+\|\Delta\phi\|_{C^{m+4}})^{m+2}(1+\sigma)^m):=\Theta(\sigma,T,\Delta\phi)(\|\mathbf{f}\|_T+\|\mathbf{g}\|_T).
			\end{align*}
			Similar to previous arguments, for small $T$ such that $\Theta(\sigma,T,\Delta\phi)\leq\frac{1}{100}$, we obtain that 
			\begin{align*}
				\|\mathbf{g}\|_T\leq \frac{1}{2}\|\mathbf{f}\|_T,
			\end{align*}
			holds in both setting
			by taking $\sigma$ small enough. 
			
			We conclude that there exists $\varepsilon_0>0$ such that, $\mathcal{S}$ is a contraction map from $\mathcal{X}^{\sigma_1}_{\infty,0}$ to itself provided $\|\Delta u_0\|_{L^\infty}\leq \varepsilon_0$, and is a contraction map from $\mathcal{X}^{\sigma_2}_{T^*,u_0\ast \rho_{\eps_1}}$ to itself provided \eqref{contf}. This completes the proof.
			
		\end{proof}
		
		\subsection{Gravity-driven thin film equation}\label{secgtf}
		
		We consider the one-phase Muskat problem driven by gravity
		and surface tension. In particular, we study the unstable scenario in which an incompressible,
		homogeneous and viscous fluid lies on top of a dry region. In the mathematical literature, this scenario has attracted a
		lot of attention due to its wide applicability. 
		Without surface tension effects, we refer to \cite{SCH04,CG07}
		for instability results when a more dense fluid is on top of a less dense one,  a less viscous fluid penetrates a more viscous one,
		or combining both features. The presence of surface tension is crucial to get rid of instabilities for short time as it introduces high
		order parabolic character to the system. We also mention the scenario where different instabilities proved to exist for a very short time giving exponentially growing modes\cite{Escher11,GHS07}.
		Moreover, in \cite{GBS20},
		the authors studied the stabilization of the
		Rayleigh-Taylor instability by capillary effects. They showed that the gravity unstable scenarios have global existence and instant smoothing for near planar
		solutions. This physical situation has been widely
		studied by means of thin film approximations.
		In \cite{Gancedo}, the authors studied gravity-capillarity driven configuration, and derived the following asymptotic approximation for this scenario by using a depth-averaged formulation. Denote $\zeta:[0,T]\times \mathbb{R}$ the free surface, the asymptotic model describes the evolution of $\zeta$ as 
		\begin{equation}\label{eqGtf}
			\begin{aligned}
				&\partial_t \zeta+\partial_x \left((1+\varepsilon\zeta)(\partial_x \zeta+\frac{\sqrt{\mu}}{bo}\partial_x ^3\zeta)\right)+\frac{\mu}{3}\partial_x^2\left((1+\varepsilon\zeta)^3\partial_x^2\zeta\right)=0,\ \ (t,x)\in(0,T)\times\mathbb{R},\\
				&\zeta|_{t=0}=\zeta_0.
			\end{aligned}
		\end{equation}
		Here $\mu>0$ is the shallowness parameter, $0<\varepsilon\leq 1$ is the nonlinearity parameter and $bo>0$ is the Bond number describing the ratio between gravitational and surface forces. The authors rigorously justified the asymptotic model and proved optimal order of convergence in the shallowness parameter to the original Muskat problem. In this paper, we prove the local well-posedness of the model \eqref{eqGtf} in the non-periodic setting. The main result is the following.
		\begin{theorem}\label{thmgd} ~\\
			Fix $m\in\mathbb{N}$, $0<\kappa<4$ and $4-\kappa\ll 1$.	There exists $0<\varepsilon_0\ll 1$ such that for any initial data $\zeta_0\in L^\infty(\mathbb{R})$, if 
			\begin{align}\label{tfcon}
				\|\zeta_0-\zeta_0\ast \rho_{\eps_1}\|_{L^\infty}\leq \varepsilon_0,
			\end{align} holds for some $\eps_1>0$,then there exists $T(\eps_0,\eps_1,\|\zeta_0\|_{L^\infty})>0$ such that \eqref{eqGtf} admits a unique solution  $u$ satisfying
			\begin{align*}
				\sup_{t\in[0,T]}	(\|\zeta(t)\|_{L^\infty}+t^\frac{m+\kappa}{4}\|\zeta(t)\|_{\dot C^{m+\kappa}})\leq C \|\zeta_0\|_{L^\infty}.
			\end{align*}
		\end{theorem}
		\begin{proof}[Proof of Theorem \ref{thmgd}]
			Remark that our analysis does not depend on the specific values of positive parameters $\varepsilon,\mu$ and $bo$. Hence, without loss of generality, we take $\varepsilon=\mu=bo=1$.\\	Let $\phi\in C^\infty(\mathbb{R})$ that will be fixed later. We rewrite \eqref{eqGtf} as 
			\begin{align*}
				&\partial_t  \zeta+\left(1+\phi+\frac{(1+\phi)^3}{3}\right)\partial_x ^4 \zeta=\partial_xF(\zeta,\phi)+R(\zeta,\phi),\quad \text{in}\ [0,T]\times \mathbb{R},\\
				&\zeta|_{t=0}=\zeta_0,\quad\text{in}\ \mathbb{R},
			\end{align*}
			where 
			\begin{align*}
				F(\zeta,\phi)&=-\left((\zeta-\phi)\partial_x ^3\zeta+(1+\zeta)\partial_x \zeta\right)+\frac{1}{3}\partial_x\left(((1+\phi)^3-(1+\zeta)^3)\partial_x^2\zeta\right),\\
				R(\zeta,\phi)&=-\partial_x \phi\partial_x^3\zeta-\frac{1}{3}\Big({2}\partial_x\left((1+\phi)^3\right)\partial_x^3\zeta+\partial_x^2\left((1+\phi)^3\right)\partial_x^2\zeta\Big).
			\end{align*}	
			For simplicity, in this proof we denote 
			$$
			\|h\|_T:=	\sup_{t\in[0,T]}	\left(\|h(t)\|_{L^\infty}+t^\frac{m+\kappa}{4}\|h(t)\|_{\dot C^{m+\kappa}}\right).
			$$
			We consider the set 
			$$
			\mathcal{X}^\sigma_{T,\phi}=\left\{\psi\in L^\infty_TL^\infty_{x}:\psi|_{t=0}=\zeta_0,\|\psi-\phi\|_{T}\leq \sigma\right\},
			$$ 
			where $\phi, T,\sigma$ will be fixed later. For any $\psi\in \mathcal{X}^\sigma_{T,\phi}$, we define a map $\mathcal{S}\psi=\zeta$, where $\zeta$ solves the equation 
			\begin{align*}
				&\partial_t  \zeta+\left(1+\phi+\frac{(1+\phi)^3}{3}\right)\partial_x ^4 \zeta=\partial_xF(\psi,\phi)+R(\psi,\phi),\\
				&\zeta|_{t=0}=\zeta_0.
			\end{align*}
			Denote $\tilde \zeta=\zeta-\phi$, $\tilde \psi=\psi-\phi$. Then 
			\begin{align*}
				&\partial_t  \tilde\zeta+\left(1+\phi+\frac{(1+\phi)^3}{3}\right)\partial_x ^4 \tilde \zeta=\partial_x F(\psi,\phi)+\tilde R(\psi,\phi),\\
				&\tilde\zeta|_{t=0}=\zeta_0-\phi
			\end{align*}
			where 
			$$
			\tilde R(\psi,\phi)= R(\psi,\phi)-\left((1+\phi)+\frac{(1+\phi)^3}{3}\right)\partial_x ^4\phi.
			$$
			By Theorem \ref{lemmain} with $F=F(\psi,\phi)$ and $G=\tilde R(\psi,\phi)$, there exists $T=T(\|\phi\|_{C^{m+4}})>0$ such that 
			\begin{align*}
				\|\tilde \zeta\|_T\lesssim& \|\tilde \zeta_0\|_{L^\infty}+\sup_{t\in[0,T]}\left(t^\frac{\kappa}{4}\|F(\psi,\phi)(t)\|_{\dot C^{\kappa-3}}+t^\frac{m+\kappa}{4}\|F(\psi,\phi)(t)\|_{\dot C^{m+\kappa-3}}+t^{\frac{m}{4}+1}\|\tilde R(\psi,\phi)(t)\|_{\dot C^m}\right)\\
				&+\|\tilde R(\psi,\phi)\|_{L^1_TL^\infty}.
			\end{align*}
			Applying Lemma \ref{Hol} we obtain for any $j\leq m$, 
			\begin{align*}
				\|F(\psi,\phi)\|_{\dot C^{j+\kappa-3}}&\lesssim \|\tilde \psi\partial_x^3\psi\|_{\dot C^{j+\kappa-3}}+\|(1+\psi)\partial_x \psi\|_{\dot C^{j+\kappa-3}}+\|((1+\phi)^3-(1+\psi)^3)\partial_x^2\psi\|_{\dot C^{j+\kappa-2}}\\
				&\lesssim t^{-\frac{j+\kappa}{4}}(\|\tilde \psi\|_T+T^\frac{1}{10}(1+\|\phi\|_{C^{m+4}}))^2(1+\|\tilde \psi\|_T+T^\frac{1}{10}(1+\|\phi\|_{C^{m+4}}))^{4}.
			\end{align*}
			Moreover, it is easy to check that 
			\begin{align*}
				&	\|\partial_x^j\tilde R(\psi,\phi)(t)\|_{L^\infty}\lesssim (1+\|\phi\|_{C^{j+4}})^5(1+\|\psi\|_{C^{j+3}})\lesssim t^{-\frac{j+3}{4}}(1+\|\phi\|_{C^{m+4}})^5(1+\|\psi\|_T), \quad\quad \quad \forall j\leq m.
			\end{align*}
			This yields
			\begin{align*}
				\sup_{t\in[0,T]}t^{\frac{m}{4}+1}\|\tilde R(\psi,\phi)(t)\|_{\dot C^m}+\|\tilde R(\psi,\phi)\|_{L^1_TL^\infty}\lesssim T^\frac{1}{4}(1+\|\phi\|_{C^{m+4}})^5(1+\|\psi\|_T).
			\end{align*}
			Hence 
			\begin{align}\label{Gtfre}
				\|\tilde \zeta\|_T\lesssim& \|\tilde \zeta_0\|_{L^\infty}+\|\tilde \psi\|_T(\|\tilde \psi\|_T+T^\frac{1}{10}\|\phi\|_{C^{m+4}})(1+\|\tilde \psi\|_T+T^\frac{1}{10}\|\phi\|_{C^{m+4}})^4+T^\frac{1}{10}(\|\tilde\psi\|_T+T^\frac{1}{10}\|\phi\|_{C^{m+4}})^2\\
				&+T^\frac{1}{10}(\|\tilde \psi\|_T+T^\frac{1}{10}\|\phi\|_{C^{m+4}}+1)\|\phi\|_{C^{m+4}}(1+\|\phi\|_{C^{m+4}})^4.\nonumber
			\end{align}
			By condition \eqref{tfcon}, we can take $\phi=\zeta_0\ast \rho_{\eps_1}$ in \eqref{Gtfre}, 
			then we obtain that there exists $T_0=T_0(\|\zeta_0\|_{L^\infty},\eps_1)$ such that for any $0<T\leq T_0$, 
			\begin{align*}
				\|\tilde \zeta\|_T\leq& C_0\|\tilde \zeta_0\|_{L^\infty}+C_0\|\tilde \psi\|_T(\|\tilde \psi\|_T+T^\frac{1}{10}\eps_1^{-4}\|\zeta_0\|_{L^\infty})(1+\|\tilde \psi\|_T+T^\frac{1}{10}\eps_1^{-4}\|\zeta_0\|_{L^\infty})^4\\
				&+C_0T^\frac{1}{10}(\|\tilde \psi\|_T+T^\frac{1}{10}\eps_1^{-4}\|\zeta_0\|_{L^\infty}+1)\|\phi\|_{C^4}(1+\eps_1^{-4}\|\zeta_0\|_{L^\infty})^4.
			\end{align*}
			Let $T_1=\min\{T_0,\left(\frac{\varepsilon_0\eps_1}{10+C_0+\|\zeta_0\|_{L^\infty}}\right)^{10m}\}$, $\sigma=2C_0\varepsilon_0$. If $\psi\in \mathcal{X}^\sigma_{T_1,\zeta_0\ast \rho_{\eps_1}}$, then 
			\begin{align*}
				\|\tilde \zeta\|_{T_1}\leq& C_0\varepsilon_0+2C_0^2\varepsilon_0(2C_0\varepsilon_0+T_1^\frac{1}{100}\eps_1^{-m-4}\|\zeta_0\|_{L^\infty})(1+2C_0\varepsilon_0+T_1^\frac{1}{100}\eps_1^{-m-4}\|\zeta_0\|_{L^\infty})^4\\
				&+C_0T_1^\frac{1}{100}(2C_0\varepsilon_0+T_1^\frac{1}{100}\eps_1^{-m-4}\|\zeta_0\|_{L^\infty}+1)\eps_1^{-m-4}\|\zeta_0\|_{L^\infty}(1+\eps_1^{-m-4}\|\zeta_0\|_{L^\infty})^4\\
				\leq&2C_0\varepsilon_0=\sigma.
			\end{align*}
			Hence $\zeta=\mathcal{S}\psi\in \mathcal{X}^\sigma_{T_1,\zeta_0\ast \rho_{\eps_1}}$. Then we prove that $\mathcal{S}:\mathcal{X}^\sigma_{T_1,\zeta_0\ast \rho_{\eps_1}}\to \mathcal{X}^\sigma_{T_1,\zeta_0\ast \rho_{\eps_1}}$ is a contraction map. Consider $\psi_1,\psi_2\in \mathcal{X}^\sigma_{T_1,\zeta_0\ast \rho_{\eps_1}}$, denote $\Psi=\psi_1-\psi_2$, $\Xi=\zeta_1-\zeta_2=\mathcal{S}\psi_1-\mathcal{S}\psi_2$. We have the equation 
			\begin{align*}
				\partial_t\Xi+&\left(1+\phi+\frac{(1+\phi)^3}{3}\right)\partial_x ^4 \Xi\\&=-\partial_x \left(\left(\Xi\partial_x ^3\psi_1+(\psi_2-\phi)\partial_x^3\Xi\right)+\big(\Xi\partial_x \psi_1+\psi_2\partial_x\Xi\big)\right)\\
				&\quad+\frac{1}{3}\partial_x^2\left(((1+\psi_2)^3-(1+\psi_1)^3)\partial_x^2\psi_1+((1+\phi)^3-(1+\psi_2)^3)\partial_x^2\psi_2\right)+R(\psi_1,\phi)-R(\psi_2,\phi)\\
				&=\partial_x F_1+\partial_{x}^2F_2+R(\psi_1,\phi)-R(\psi_2,\phi).
			\end{align*} 
			Note that 
			\begin{align*}
				&\|F_1(t)\|_{\dot C^{\kappa-3}}\lesssim t^{-\frac{\kappa}{4}} \|\Xi\|_T(\|\tilde \psi_1\|_T+\|\tilde \psi_2\|_T+T^\frac{1}{10}\|\phi\|_{C^4}).\\
				&\|F_2(t)\|_{\dot C^{\kappa-2}}\lesssim t^{-\frac{\kappa}{4}}\|\Xi\|_T(\|\tilde \psi_1\|_T+\|\tilde \psi_2\|_T+T^\frac{1}{10}\|\phi\|_{C^4})(1+\|\tilde \psi_1\|_T+\|\tilde \psi_2\|_T+T^\frac{1}{10}\|\phi\|_{C^4})^4.\\
				&\|R(\psi_1,\phi)(t)-R(\psi_2,\phi)(t)\|_{L^\infty}\lesssim t^{-\frac{3}{4}}\|\Xi\|_T\|\phi\|_{C^4}(1+\|\phi\|_{C^4})^4.
			\end{align*}
			For higher order estimates, we have
			\begin{align*}
				&\|F_1(t)\|_{\dot C^{m+\kappa-3}}\lesssim t^{-\frac{m+\kappa}{4}}\|\Xi\|_T(\|\tilde\psi_1\|_T+\|\tilde\psi_2\|_T+T^{\frac{1}{10}}\|\phi\|_{C^{m+4}})^4,\\
				&\|F_2(t)\|_{\dot C^{m+\kappa-2}}\lesssim t^{-\frac{m+\kappa}{4}}\|\Xi\|_T(\|\tilde\psi_1\|_T+\|\tilde\psi_2\|_T+T^{\frac{1}{10}}\|\phi\|_{C^{m+4}})^4(1+\|\tilde \psi_1\|_T+\|\tilde \psi_2\|_T+T^\frac{1}{10}\|\phi\|_{C^{m+4}})^4,\\ &\|R(\psi_1,\phi)(t)-R(\psi_2,\phi)(t)\|_{\dot C^m}\lesssim t^{-\frac{m+3}{4}}\|\Xi\|_T\|\phi\|_{C^{m+4}}(1+\|\phi\|_{C^{m+4}})^4.
			\end{align*}
			Take $
			\phi=\zeta_0\ast \rho_{\eps_1}$, by Theorem \ref{lemmain} we obtain 
			\begin{align*}
				\|\Xi\|_{T_1}&\leq	C_0\|\Psi\|_{T_1} (\|(\tilde \psi_1,\tilde \psi_2)\|_{T_1}+{T_1}^\frac{1}{100}\eps_1^{-m-4}\|\zeta_0\|_{L^\infty})(1+\|(\tilde \psi_1,\tilde \psi_2)\|_{T_1}+{T_1}^\frac{1}{100}\eps_1^{-m-4}\|\zeta_0\|_{L^\infty})^4\\
				&\leq\frac{1}{2}\|\Psi\|_{T_1}.
			\end{align*}
			Hence $\mathcal{S}:\mathcal{X}^\sigma_{T_1,\zeta_0\ast \rho_{\eps_1}}\to \mathcal{X}^\sigma_{T_1,\zeta_0\ast \rho_{\eps_1}}$ is a contraction map. This completes the proof.
		\end{proof}

		\subsection{The Muskat equation with surface tension }
		\label{intomus1}
		The Muskat equation is a classical model in the analysis of free surface flows,
		which describes the dynamics of two incompressible and immiscible fluids in a porous medium. We consider here the particular case when the fluids
		have equal viscosities and we assume that the flows are two-dimensional. The free boundary separating the fluids
		is described by the evolution equation
		\begin{equation}\label{eqmst}
			\begin{aligned}
				&\partial_t f=\frac{\mathbf{k}}{2\pi\mu}\mathrm{P.V.} \int_{\mathbb{R}}\frac{1+\partial_x f(x)\Delta_\alpha f(x)}{\langle\Delta_\alpha f(x)\rangle^2}\partial_x\left(\sigma'\kappa(f)- \varrho_0 f\right)(x-\alpha)\frac{d\alpha}{\alpha},\\
				&f|_{t=0}=f_0.
			\end{aligned}
		\end{equation}
		We refer readers to \cite{1Matioc2018} for a proof of equivalence of \eqref{eqmst} to the classical formulation of the Muskat problem.
		For brevity, we write $f'$
		for the spatial derivative $\partial_x f$. Here $\mathbf{k}, \mu, \sigma'$ are positive constants that denote the permeability of the
		homogeneous porous medium,  the viscosity of the fluids and the surface tension
		coefficient, respectively. And 
		\begin{align}\label{defvr0}
			\varrho_0=\mathbf{g}(\rho_--\rho_+),
		\end{align}
		where $\mathbf{g}$ is the Earth’s gravity and $\rho_+, \rho_-$ denote the density of the upper fluid and lower fluid.  
		Moreover, to shorten the notation we denote $\langle A\rangle=(1+A^2)^\frac{1}{2}$, and 
		$$
		\delta_\alpha f(x)=f(x)-f(x-\alpha),\quad\quad \ \ \Delta_\alpha f(x)=\frac{\delta_\alpha f(x)}{\alpha}.
		$$
		Finally, $\kappa(f)=\frac{f''}{\langle f' \rangle^3}$ is the curvature of the graph $\{y=f(t,x)\}$. 
		The Muskat problem has attracted much interest in the last decades.  The well-posedness of the Muskat problem with surface tension 
		effects has been investigated in bounded geometries in \cite{Escher11,Escher12,Escher18}. For unbounded geometry, we refer to \cite{Tofts17,1Matioc2018,Matioctams,1HuyNguyen2020} 
		for results in subcritical $L^2$ based Sobolev space. More recently, the reference \cite{MATIOC2022} establish a local well-posedness
		result in all subcritical Sobolev spaces $W^{s,p}(\mathbb{R})$ with $p\in(1,\infty)$ and $s\in(1+\frac{1}{p},2)$. We remark that  $\varrho_0>0$ corresponds to the stable regime (heavier fluid below). And $\varrho_0<0$ corresponds to the unstable regime  (heavier fluid up). 
		
		For simplicity, we fixed the constants $\mathbf{k}=2, \mu=1$ and $\sigma'=1$. 
		Note that 
		\begin{align*}
			&	\frac{1}{\pi}\mathrm{P.V.}\int_{\mathbb{R}}\frac{1+\partial_x f(x)\Delta_\alpha f(x)}{\langle\Delta_\alpha f(x)\rangle^2}\partial_x\left(\kappa(f)- \varrho_0 f\right)(x-\alpha)\frac{d\alpha}{\alpha} = - \frac{\Lambda^3f(x)}{\langle\partial_x f(x)\rangle^3}+\N[f](x),
		\end{align*}
		where $\N[f]=\N_1[f]+\N_2[f]+\varrho_0\N_3[f]$ with 
		\begin{equation}\label{defNst}
			\begin{aligned}
				&\N_{1}[f](x)=\frac{1}{\pi} \int_{\mathbb{R}}\left(\frac{ \Delta_{\alpha} f(x)(\partial_{x} f(x)-\Delta_{\alpha} f(x))}{\left\langle\Delta_{\alpha} f(x)\right\rangle^{2}}\right) \partial_{x}\left(\frac{\partial_{x}^{2} f(x-\alpha)}{\left\langle\partial_{x} f(x-\alpha)\right\rangle^{3}}\right) \frac{d \alpha}{\alpha},\\
				&\N_2[f](x)=-\frac{1}{\pi} \int_{\mathbb{R}}\partial_x^2f(x-\alpha)\left(\frac{1}{\left\langle\partial_{x} f(x-\alpha)\right\rangle^{3}}-\frac{1}{\left\langle\partial_{x} f(x)\right\rangle^{3}}\right)\frac{d\alpha}{\alpha^2},\\
				&\N_3[f](x)=-\frac{1}{\pi} \int_{\mathbb{R}}\left(\frac{ \Delta_{\alpha} f(x)(\partial_{x} f(x)-\Delta_{\alpha} f(x))}{\left\langle\Delta_{\alpha} f(x)\right\rangle^{2}}\right) \partial_{x}f(x-\alpha)\frac{d \alpha}{\alpha}-\Lambda f.
			\end{aligned}
		\end{equation}
		Hence \eqref{eqmst} can be written as 
		\begin{align*}
			&\partial_t f(x)+\frac{\Lambda^3f(x)}{\langle \partial_x f(x)\rangle^3}=\N[f](x),\\
			&f|_{t=0}=f_0.
		\end{align*}
		\begin{proof}[Proof of Theorem \ref{thmmst}]
			For simplicity, denote 
			\begin{equation}\label{normmst}
				\begin{aligned}
					&	\|h\|_{T}:=\sup_{t\in[0,T]}(\| \partial_xh(t)\|_{L^\infty}+t^\frac{m+\kappa}{3}\| \partial_xh(t)\|_{\dot C^{m+\kappa}}).
				\end{aligned}
			\end{equation}
			We first prove \textit{1)}. Suppose $\varrho_0=0$ and $\|f_0'\|_{L^\infty}\leq \varepsilon_0$. Consider the set 
			\begin{align*}
				\mathcal{X}^\sigma=\mathcal{X}_{\infty,0}^\sigma=\left\{g\in L^\infty_T\dot W^{1,\infty}:g|_{t=0}=f_0,\|g\|_{\infty}\leq \sigma\right\},
			\end{align*}
			where $\sigma>0$ is a constant that will be fixed later. For $g\in \mathcal{X}^\sigma$, we define a map $\mathcal{S}g=f$, where $f$ is the solution to the equation 
			\begin{equation}\label{mapst}
				\begin{aligned}
					&\partial_t f+{\Lambda^3f(x)}=\N[g](x)+\Lambda^3g(x)\left(1-\frac{1}{\langle\partial_x g(x)\rangle^3}\right),\\
					&f|_{t=0}=f_0.
				\end{aligned}
			\end{equation}
			For short of notation, we denote 
			\begin{align}\label{defffg}
				\G[g_1,g_2]=	\left(\frac{1}{\langle\partial_xg_2\rangle^3}-\frac{1}{\langle\partial_xg_1\rangle^3}\right)\Lambda^3 g_1.
			\end{align}
			In the following, we prove that $\mathcal{S}$ is a contraction map from $	\mathcal{X}^\sigma$ to itself. Hence there exists a fix point in $	\mathcal{X}^\sigma$, which is a global solution to  equation \eqref{eqmst}.
			Take one space derivative in \eqref{mapst}, we get 
			\begin{align*}
				&\partial_t \partial_xf+{\Lambda^3\partial_xf(x)}=\partial_x\N[g](x)+\partial_x\G[g,0](x),\\
				&\partial_xf|_{t=0}=\partial_xf_0.
			\end{align*}
			Then we apply Theorem \ref{lemmain} for $F=\N[g](x)+\G[g,0](x)$ and $G=0$ to obtain that 
			\begin{align*}
				\|f\|_\infty\lesssim &\|\partial_x f_0\|_{L^\infty}+\sup_{t>0}\left(t^\frac{\kappa}{3}\|\N[g](t)\|_{\dot C^{\kappa-2}}+t^\frac{m+\kappa }{3}\|\N[g](t)\|_{\dot C^{m+\kappa -2}}\right)\\
				&\ \ +\sup_{t>0}\left(t^\frac{\kappa}{3}\|\G[g,0](t)\|_{\dot C^{\kappa-2}}+t^\frac{m+\kappa }{3}\|\G[g,0](t)\|_{\dot C^{m+\kappa -2}}\right).
			\end{align*}
			Applying Lemma \ref{nonst} to get 
			\begin{align*}
				\sup_{t>0}\left(t^\frac{\kappa}{3}\|\N[g](t)\|_{\dot C^{\kappa-2}}+t^\frac{m+\kappa }{3}\|\N[g](t)\|_{\dot C^{m+\kappa -2}}\right)\lesssim \|g\|_\infty^2(1+\|g\|_\infty)^{m+5}.
			\end{align*}
			Applying Lemma \ref{GGG} with $f_1=g$, $f_2\equiv 0$, and $\phi\equiv 0$ to obtain that
			\begin{align*}
				\sup_{t>0}\left(t^\frac{\kappa}{3}\|\G[g,0](t)\|_{\dot C^{\kappa-2}}+t^\frac{m+\kappa }{3}\|\G[g,0](t)\|_{\dot C^{m+\kappa -2}}\right)\lesssim \|g\|_\infty^2(1+\|g\|_\infty)^{m+5}.
			\end{align*}
			Hence it follows that 
			\begin{align*}
				\|f\|_\infty\leq C_0\|\partial_x f_0\|_{L^\infty}+C_0\|g\|_\infty^2(1+\|g\|_\infty)^{m+5},
			\end{align*}
			here the constant $C_0$ depends only on $m$. 	Take $\varepsilon_0=\frac{1}{100(2C_0+1)^{m+10}}$ and $\sigma=2C_0 \varepsilon_0$. Then for any $g\in \mathcal{X}^{\sigma}$, we have $\|g\|_{\infty}\leq 2C_0 \varepsilon_0$. Hence 
			\begin{align*}
				\|f\|_{\infty}\leq C_0\varepsilon_0+4C_0^3\varepsilon_0^2(2C_0 \varepsilon_0+1)^{m+5}\leq \sigma.
			\end{align*}
			Hence we obtain that $f\in \mathcal{X}^{\sigma}$.
			Moreover, for $g_1,g_2\in \mathcal{X}^\sigma$, 
			let $f_1=\mathcal{S}g_1, f_2=\mathcal{S}g_2$. Then 
			\begin{align*}
				&\partial_t \partial_x(f_1-f_2)+{\Lambda^3\partial_x(f_1-f_2)(x)}=\partial_x(\N[g_1]-\N[g_2])(x)+\partial_x(\G[g_1,0]-\G[g_2,0])(x),\\
				&\partial_x(f_1-f_2)|_{t=0}=0.
			\end{align*}
			Apply Theorem \ref{lemmain} again, we obtain
			\begin{align*}
				\|f_1-f_2\|_\infty\lesssim &\sup_{t>0}\left(t^\frac{\kappa}{3}\|(\N[g_1]-\N[g_2])(t)\|_{\dot C^{\kappa-2}}+t^\frac{m+\kappa }{3}\|(\N[g_1]-\N[g_2])(t)\|_{\dot C^{m+\kappa -2}}\right)\\
				&\ \ +\sup_{t>0}\left(t^\frac{\kappa}{3}\|(G[g_1,0]-G[g_2,0])(t)\|_{\dot C^{\kappa-2}}+t^\frac{m+\kappa }{3}\|(\G[g_1,0]-\G[g_2,0])(t)\|_{\dot C^{m+\kappa -2}}\right).
			\end{align*}
			By Lemma \ref{nonst} and Lemma \ref{GGG},
			\begin{align*}
				&\sup_{t>0}\left(t^\frac{\kappa}{3}\|(\N[g_1]-\N[g_2])(t)\|_{\dot C^{\kappa-2}}+t^\frac{m+\kappa }{3}\|(\N[g_1]-\N[g_2])(t)\|_{\dot C^{m+\kappa -2}}\right)\\
				&\quad\quad+\sup_{t>0}\left(t^\frac{\kappa}{3}\|(\G[g_1,0]-\G[g_2,0])(t)\|_{\dot C^{\kappa-2}}+t^\frac{m+\kappa }{3}\|(\G[g_1,0]-\G[g_2,0])(t)\|_{\dot C^{m+\kappa -2}}\right)\\
				&\lesssim \|g_1-g_2\|_\infty\|(g_1,g_2)\|_{\infty}(1+\|(g_1,g_2)\|_{\infty})^{m+5}.
			\end{align*}
			Then we get 
			\begin{align*}
				\|f_1-f_2\|_\infty&\leq C_1\|g_1-g_2\|_\infty\|(g_1,g_2)\|_{\infty}(1+\|(g_1,g_2)\|_{\infty})^{m+5}\leq C_1\sigma(1+2\sigma)^{m+5}\|g_1-g_2\|_\infty\\
				&\leq \frac{1}{2}\|g_1-g_2\|_\infty,
			\end{align*}
			provided $C_1\sigma(1+2\sigma)^{m+5}\leq \frac{1}{2}$.
			Hence there exists a unique $f\in \mathcal{X}^\sigma$ such that $\mathcal{S}f=f$. This completes the proof of \textit{i)}.
			\vspace{0.3cm}\\
			Now we prove \textit{ii)}. When $\varrho_0\neq 0$, the equation \eqref{eqmst} is not homogeneous. 
			We introduce the following norm 
			\begin{align}\label{musnol}
				&\|h\|_{X_T}:=\|h\|_{L^\infty_TL^\infty}+\|h\|_T.
			\end{align}
			If condition \eqref{stcon} holds, then
			take $\phi=f_0\ast \rho_{\eps_1}$, there holds 
			\begin{align}\label{esinit}
				\|f_0-\phi\|_{W^{1,\infty}}\leq\varepsilon_0,\ \ \ \|\phi\|_{C^{m+5}}\lesssim \eps_1^{-(m+5)}\|f_0\|_{W^{1,\infty}}.
			\end{align}
			Consider the set
			$$
			\mathcal{X}_{T,\phi}^\sigma=\{g\in L^\infty_T\dot W^{1,\infty}:g|_{t=0}=f_0,\|g-\phi\|_{X_T}\leq \sigma\}.
			$$ 
			For any $g\in	\mathcal{X}_{T,\phi}^\sigma$, where $T,\sigma$ will be fixed later, we define a map $\tilde {\mathcal{S}}g=f$, where $f$ is the solution to the equation 
			\begin{equation}\label{mapst}
				\begin{aligned}
					&\partial_t f(x)+\frac{\Lambda^3f(x)}{\langle\partial_x\phi(x)\rangle^3}=\N[g](x)+\G[g,\phi](x),\\
					&f|_{t=0}=f_0,
				\end{aligned}
			\end{equation}
			where $\G[g,\phi]$ is defined in \eqref{defffg}.
			In the following, we prove that $\tilde{\mathcal{S}}$ is a contraction map from $	\mathcal{X}_{T,\phi}^\sigma$ to itself. Hence there exists a fixed point in $	\mathcal{X}_{T,\phi}^\sigma$, which is a solution to  equation \eqref{eqmst}.\\
			Take one space derivative in \eqref{mapst}, and denote $\tilde f(t,x)=f(t,x)-\phi(x)$, $\tilde g(t,x)=g(t,x)-\phi(x)$, one has 
			\begin{align*}
				&	\partial_t \partial_x\tilde f+\frac{\Lambda^3\partial_x\tilde f}{\langle\partial_x\phi\rangle^3}= \partial_x(\N[g]+
				\G[g,\phi])+R[f,\phi],\\
				&\partial_x\tilde f|_{t=0}=\partial_x (f_0-\phi ),
			\end{align*}
			where
			\begin{equation}\label{defnonl}
				R[f,\phi]=\frac{3\Lambda^3 f\partial_x\phi\partial_{x}^2\phi}{\langle\partial_x\phi\rangle^5}-\frac{\Lambda^3\partial_x\phi}{\langle\partial_x\phi\rangle^3}.
			\end{equation}
			Applying Theorem \ref{lemmain} with $f=\N[g]+\G[g,\phi]$ and $g=R[f,\phi]$ to obtain that there exists $T=T(\|\phi\|_{C^{m+4}})>0$ such that 
			\begin{equation}\label{stm}
				\begin{aligned}
					\|\tilde f\|_{T}
					\lesssim&\|\partial_x(f_0-\phi)\|_{L^\infty}+	\sup_{t\in[0,T]} \left(t^\frac{\kappa}{3}\|(\N[g],\G[g,\phi])(t)\|_{\dot C^{\kappa-2}}+t^\frac{m+\kappa }{3}(\|(\N[g],\G[g,\phi])(t)\|_{\dot C^{m+\kappa -2}})\right)\\
					&+\|R[f,\phi]\|_{L^1_TL^\infty}+\sup_{t\in[0,T]}t^{\frac{m}{3}+1}\|R[f,\phi](t)\|_{\dot C^m}.
				\end{aligned}
			\end{equation}
			Applying Lemma \ref{nonst} to obtain  
			\begin{equation}\label{st1}
				\begin{aligned}
					&	\sup_{t\in[0,T]}(t^\frac{\kappa}{3}\| \N[g](t)\|_{\dot C^{\kappa-2}}+t^\frac{m+\kappa}{3}\| \N[g](t)\|_{\dot C^{m+\kappa-2}})\lesssim (\| g\|_{T,*}^2+T^\frac{2}{3}\|g\|_{X_T})(1+\| g\|_{X_T})^{m+5}.
				\end{aligned}
			\end{equation}
			See \eqref{defnn} for the definition of the non-endpoint norm $\|\cdot\|_{T,*}$.
			Observe that $\G[\phi,\phi]\equiv 0$, hence applying Lemma \ref{GGG} with $f_1=g$, $f_2=\phi$ to obtain that 
			\begin{align}\label{st2}
				\sup_{t\in[0,T]}(t^\frac{\kappa}{3}	\|\G[g,\phi](t)\|_{\dot C^{\kappa-2}}+t^\frac{m+\kappa}{3}	\|\G[g,\phi](t)\|_{\dot C^{m+\kappa-2}})\lesssim \|\tilde g\|_{T}(\|g\|_{T,*}+\|\tilde g\|_T)(1+\|(g,\phi)\|_{T})^{m+5}.
			\end{align}
			Moreover, 
			\begin{align*}
				&		\|R[f,\phi]\|_{L^\infty}\lesssim (\|\Lambda^3 f\|_{L^\infty}+1)\|\partial_x\phi\|_{C^4},\\
				&	\|R[f,\phi]\|_{\dot C^m}\lesssim (\|\Lambda^3 f\|_{C^m}+1)\|\partial_x\phi\|_{C^{m+4}}^m.
			\end{align*}
			Hence 
			\begin{equation}
				\begin{aligned}
					&		\|R[f,\phi](\tau)\|_{L_T^1L^\infty}+\sup_{t\in[0,T]}t^{\frac{m}{3}+1}\|R[f,\phi](t)\|_{\dot C^m}\lesssim T^\frac{1}{3}(\|\tilde f\|_T+1)\|\partial_x\phi\|_{C^{m+4}}^m.\label{Res}
				\end{aligned}
			\end{equation}
			We conclude from \eqref{stm}-\eqref{Res} that 
			\begin{equation}\label{stfin}
				\begin{aligned}
					\|\tilde f\|_{T}\lesssim  
					&\|\partial_x (f_0-\phi)\|_{L^\infty}+(\| g\|_{T,*}^2+T^\frac{1}{10}\|g\|_{X_T}+\|\tilde g\|_{X_T}^2)(1+\| (g,\phi)\|_{X_T})^{m+10}\\
					&\ \ \ +T^\frac{1}{3}(\|\tilde f\|_T+1)\|\partial_x\phi\|_{C^{m+4}}^m.
				\end{aligned}
			\end{equation}
			On the other hand, integrating \eqref{mapst} in time to obtain 
			\begin{align*}
				|f(t,x)-f_0(x)|&\leq \int_0^t \|\Lambda^3f(\tau)\|_{L^\infty}+\|\N[g](\tau)\|_{L^\infty}+\|\G[g,\phi](\tau)\|_{L^\infty}d\tau\\
				&\leq \tilde C_0t^\frac{1}{3} (1+\|f\|_{X_T}+\|(g,\phi)\|_{X_T})^5,\ \ \ \forall t\in[0,T].
			\end{align*}
			This implies that 
			\begin{align*}
				\|\tilde f\|_{L^\infty_TL^\infty}\leq \|f_0-\phi\|_{L^\infty}+T^\frac{1}{3}(1+\|f\|_{X_T}+\|(g,\phi)\|_{X_T})^5.
			\end{align*}
			Combining this with \eqref{stfin}, 
			and the facts that 
			\begin{align*}
				&   \|g\|_{T,*}\leq \|\tilde g\|_{X_T}+T^\frac{1}{15}\|\phi\|_{C^{m+4}}\leq  \|\tilde g\|_{X_T}+T^\frac{1}{15}\eps_1^{-(m+4)}\|f_0\|_{W^{1,\infty}},\\
				&    \|g\|_{X_T}\leq \|\tilde g\|_{X_T}+\|\phi\|_{X_T}\leq\|\tilde g\|_{X_T}+C\|f_0\|_{W^{1,\infty}}(1+T^\frac{m+4}{3}\eps_1^{-(m+4)}),
			\end{align*}
			we obtain 
			\begin{align*}
				\|\tilde f\|_{X_T}&\leq \tilde C_0\|f_0-\phi\|_{W^{1,\infty}}+\tilde C_0(\|\tilde g\|_{X_T}+T^\frac{1}{15}\eps_1^{-(m+4)}\|f_0\|_{W^{1,\infty}})^2(1+T^\frac{1}{15}\eps_1^{-(m+4)}\|f_0\|_{W^{1,\infty}})^{m+10}\\
				&\quad\quad+\tilde C_0T^\frac{1}{3}(1+\|\tilde f\|_{X_T}+\|\tilde g\|_{X_T}+\eps_1^{-(m+5)}\|f_0\|_{W^{1,\infty}})(1+T^\frac{1}{15}\eps_1^{-(m+4)}\|f_0\|_{W^{1,\infty}})^{m+10}.
			\end{align*}
			Here the constant $\tilde C_0$ depends only on $m$ and $\varrho_0$. 
			Let  
			\begin{align}\label{sma}
				\sigma'=2\tilde C_0\varepsilon_0, \quad\quad\text{and}\quad 0<T<\left(\frac{\varepsilon_0\eps_1}{1+\tilde C_0+\|f_0\|_{W^{1,\infty}}}\right)^{(m+10)^2}.
			\end{align}  For any $g\in  \mathcal{X}_{T,f_0\ast \rho_{\eps_1}}^{\sigma'}$, we have
			\begin{align*}
				\|\tilde f\|_{X_T}&\leq \frac{5}{4}\tilde C_0\varepsilon_0+\tilde C_0(2\tilde C_0\varepsilon_0+T^\frac{1}{15}\eps_1^{-(m+4)}\|f_0\|_{W^{1,\infty}})^2(1+T^\frac{1}{15}\eps_1^{-(m+4)}\|f_0\|_{W^{1,\infty}})^{m+10}\\
				&\leq 2\tilde C_0\varepsilon_0,
			\end{align*}
			provided $\varepsilon_0\leq \frac{1}{(10+\tilde C_0)^{10}}$.
			Hence $f\in  \mathcal{X}_{T,f_0\ast \rho_{\eps_1}}^{\sigma'}$.
			
			In the following we do the contraction estimates. Consider $g_1,g_2\in \mathcal{X}^\sigma_{T,\phi}$. Denote $\mathbf{g}=g_1-g_2$, $\mathbf{f}=f_1-f_2=\tilde {\mathcal{S}}g_1-\tilde {\mathcal{S}}g_2$. We have 
			\begin{equation}\label{mapst1}
				\begin{aligned}
					&\partial_t \mathbf{f}+\frac{\Lambda^3\mathbf{f}}{\langle\partial_x\phi\rangle^3}=\N[g_1]-\N[g_2]+\G[g_1,\phi]-\G[g_2,\phi],\\
					&\mathbf{f}|_{t=0}=0.
				\end{aligned}
			\end{equation}
			Take one spatial derivative, we get
			\begin{align*}
				&\partial_t \partial_x\mathbf{f}+\frac{\Lambda^3\partial_x\mathbf{f}}{\langle\partial_x\phi\rangle^3}
				=\partial_x\left(\N[g_1]-\N[g_2]+\G[g_1,\phi]-\G[g_2,\phi]\right)-\frac{3\partial_x\phi\partial_x^2\phi\Lambda^3\mathbf{f}}{\langle\partial_x\phi\rangle^5}:=\partial_x F_1+F_2,\\
				&\partial_x \mathbf{f}|_{t=0}=0.
			\end{align*}
			Apply Theorem \ref{lemmain} with $f=F_1$ and $g=F_2$, there exists $T>0$ such that 
			\begin{equation}\label{ft1}
				\begin{aligned}
					\|\mathbf{f}\|_T\lesssim\sup_{t\in[0,T]}(t^\frac{\kappa}{3} \left\|F_1(t)\right\|_{\dot C^{\kappa-2}}+t^\frac{m+\kappa}{3} \left\|F_1(t)\right\|_{\dot C^{m+\kappa-2}})+\|F_2\|_{L^1_TL^\infty}+\sup_{t\in[0,T]}t^{\frac{m}{3}+1}\|F_2[t]\|_{\dot C^m}.
				\end{aligned}
			\end{equation}
			Applying Lemma \ref{nonst} to obtain 
			\begin{align*}
				&\sup_{t\in[0,T]}(t^\frac{\kappa}{3}\|(\N[g_1]-\N[g_2])(t)\|_{\dot C^{\kappa-2}}+t^\frac{m+\kappa}{3}\|(\N[g_1]-\N[g_2])(t)\|_{\dot C^{m+\kappa-2}}\\
				&\quad\quad\lesssim \|\mathbf{g}\|_{X_T}(\|(g_1,g_2)\|_{T,*}+T^\frac{1}{10})(1+\|(g_1,g_2)\|_{X_T})^{m+4}.
			\end{align*}
			Moreover, by Lemma \ref{GGG},
			\begin{align*}
				&\sup_{t\in[0,T]}(t^\frac{\kappa}{3}\|(\G[g_1,\phi]-\G[g_2,\phi])(t)\|_{\dot C^{\kappa-2}}+t^\frac{m+\kappa}{3}\|(\G[g_1,\phi]-\G[g_2,\phi])(t)\|_{\dot C^{m+\kappa-2}})\\
				&\quad\quad\quad\lesssim \|\mathbf{g}\|_{T}(\|g_1\|_{T,*}+\|g_2-\phi\|_T)(1+\|(g_1,g_2,\phi)\|_{T})^{m+5}.
			\end{align*}
			Hence we conclude that 
			\begin{equation}\label{esF1}
				\begin{aligned}
					&\sup_{t\in[0,T]}(t^\frac{\kappa}{3} \left\|F_1(t)\right\|_{\dot C^{\kappa-2}}+t^\frac{m+\kappa}{3} \left\|F_1(t)\right\|_{\dot C^{m+\kappa-2}})\\
					&\quad\quad\lesssim\|\mathbf{g}\|_{X_T} \left(\|(g_1-\phi,g_2-\phi)\|_T+T^\frac{1}{15}\|\phi\|_{C^{m+5}}\right)(1+\|(g_1,g_2,\phi)\|_{X_T})^{m+5}.
				\end{aligned}
			\end{equation}
			Finally, for $F_2$, we have 
			\begin{align*}
				&\|F_2(t)\|_{L^\infty}\lesssim \|\phi\|_{C^2}\|\Lambda^3\mathbf{f}(t)\|_{L^\infty}\lesssim t^{-\frac{2}{3}}\|\phi\|_{C^2}\|\mathbf{f}\|_T,\\
				&\|F_2(t)\|_{\dot C^m}\lesssim \sum_{l=0}^m\|\partial_x ^l\Lambda^3\mathbf{f}\|_{L^\infty}(\|\partial_x^2\phi\|_{L^\infty}^{1+m-l}+\|\partial_x^{2+m-l}\phi\|_{L^\infty})\lesssim t^{-\frac{m+2}{3}}(1+\|\phi\|_{C^{m+2}})^{m+2}\|\mathbf{f}\|_T.
			\end{align*}
			Hence 
			\begin{align}\label{esF2}
				\|F_2\|_{L^1_TL^\infty}+\sup_{t\in[0,T]}t^{\frac{m}{3}+1}\|F_2(t)\|_{\dot C^m}\lesssim T^\frac{1}{3}(1+\|\phi\|_{C^{m+2}})^{m+2}\|\mathbf{f}\|_{T}.
			\end{align}
			We obtain from \eqref{ft1}, \eqref{esF1} and \eqref{esF2} that 
			\begin{equation}   \label{haha1}
				\begin{aligned}
					\|\mathbf{f}\|_T\lesssim & \|\mathbf{g}\|_{X_T} (\|(g_1-\phi,g_2-\phi)\|_T+T^\frac{1}{15}\|\phi\|_{C^{m+5}})(1+\|(g_1,g_2,\phi)\|_{X_T})^{m+5}\\
					&+T^\frac{1}{3}(1+\|\phi\|_{C^{m+2}})^{m+2}\|\mathbf{f}\|_{T}.
				\end{aligned}
			\end{equation}
			Finally, to control $\|\mathbf{f}\|_{X_T}$, it remains to estimate $\|\mathbf{f}\|_{L^\infty_TL^\infty}$. Integrate \eqref{mapst1} in time, we obtain that 
			\begin{align*}
				\|\mathbf{f}\|_{L^\infty_TL^\infty} &\leq \int_0^T \left(\|\Lambda^3 \mathbf{f}(\tau)\|_{L^\infty}+\|(\N[g_1]-\N[g_2])(\tau)\|_{L^\infty}+\|(\G[g_1,\phi]-
				\G[g_2,\phi])(\tau)\|_{L^\infty}\right)d\tau\\
				&\lesssim T^\frac{1}{3}(\|\mathbf{f}\|_T+\|\mathbf{g}\|_T)(1+\|(g_1,g_2)\|_{X_T})^{m+5}.
			\end{align*}
			Combining this with \eqref{haha1} to get 
			\begin{align*}
				\|\mathbf{f}\|_{X_T}\lesssim & \|\mathbf{g}\|_{X_T} (\|(g_1-\phi,g_2-\phi)\|_T+T^\frac{1}{15}(1+\|\phi\|_{C^{m+5}}))(1+\|(g_1,g_2,\phi)\|_{X_T})^{m+5}\\
				&+T^\frac{1}{3}\|\mathbf{f}\|_{T}(1+\|\phi\|_{C^{m+2}}+\|(g_1,g_2)\|_{X_T})^{m+5}.
			\end{align*}
			If $g_1,g_2\in  \mathcal{X}_{T,\phi}^{\sigma'}$, one has $\|(g_1-\phi,g_2-\phi)\|_{X_T}\leq 2 \sigma'$ and  $\|(f_1-\phi,f_2-\phi)\|_{X_T}\leq 2 \sigma'$. Combining this with \eqref{esinit} yields that 
			\begin{align*}
				\|\mathbf{f}\|_{X_T}\leq & \tilde C_1\|\mathbf{g}\|_{X_T} (2\sigma'+{T}^\frac{1}{100}\eps_1^{-(m+5)}\|f_0\|_{W^{1,\infty}})(1+2\sigma'+{T}^\frac{1}{100}\eps_1^{-(m+5)}\|f_0\|_{W^{1,\infty}})^{m+5}\\
				&+\tilde C_1T^\frac{1}{3}\|\mathbf{f}\|_{X_T}(2\sigma'+1)(2\sigma'+2{T}^\frac{1}{100}\eps_1^{-(m+5)}\|f_0\|_{W^{1,\infty}})\\
				:=&\Theta(T,\sigma',\eps_1,\|f_0\|_{W^{1,\infty}})(\|\mathbf{f}\|_{X_T}+\|\mathbf{g}\|_{X_T}).
			\end{align*}
			Then we can take $\varepsilon_0$ small enough and $\sigma', T$ defined in \eqref{sma} to let $\Theta(T,\sigma',\eps_1,\|f_0\|_{W^{1,\infty}})\leq \frac{1}{100}$, then
			\begin{align*}
				\|\mathbf{f}\|_{T}\leq \frac{1}{2}\|\mathbf{g}\|_{X_T}.
			\end{align*}
			Hence $\tilde{\mathcal{S}}:\mathcal{X}_{T,\phi}^{\sigma'}\to \mathcal{X}_{T,\phi}^{\sigma'}$ is a contraction map. This completes the proof.
		\end{proof}
		
		\subsection{The Muskat equation in  general setting}\label{secmusgene}
		As we mentioned in section \ref{secintro}, for the Muskat equation in general setting, we need to work in H\"{o}lder space with logarithm derivatives.
		Fix $\varkappa>1$. For $\gamma\in[0,1)$, $k\in\mathbb{N}$, define
		\begin{equation}\label{ctlog}
			\begin{aligned}
				&\|h\|_{\dot C^{\gamma,\log^\varkappa}}:=\sup_{x\neq y}|h(x)-h(y)|\log^\varkappa\left(2+\frac{1}{|x-y|}\right)|x-y|^{-\gamma},\\
				&\|h\|_{\dot C^{k+\gamma,\log ^\varkappa}}=\|\nabla^k h\|_{\dot C^{\gamma,\log^\varkappa}}.
			\end{aligned}
		\end{equation}
		And we denote $\|h\|_{\dot C^{\log^\varkappa}}=\|h\|_{\dot C^{0,\log^\varkappa}}$ for short. Note that
		\begin{equation*}
			\int_{|\alpha|<1}\|\delta_\alpha f\|_{L^\infty}\frac{d\alpha}{|\alpha|^d}\lesssim\frac{1}{\varkappa-1}\|f\|_{\dot C^{\log^\varkappa}}.
		\end{equation*}
		We introduce the following lemma, which is an analogy to Theorem \ref{lemmain} in  H\"{o}lder space with logarithm derivatives.
		\begin{lemma}
			\label{lemlog}
			Let $u$ be a solution to system \eqref{eqpara}.
			There exists $T=T(\tilde \M,\beta)\ll 1$ such that
			\begin{equation}\label{relog}
				\begin{aligned}
					&\sup_{t\in[0,T]}(\|u(t)\|_{L^\infty}+\|u(t)\|_{\dot C^{\log^\varkappa}}+|\log t|^{a}t^{\frac{m+\kappa}{s}}\|u(t)\|_{\dot C^{m+\kappa}})\\
					&\quad\quad\quad	\lesssim\tilde  \M^\beta \left(\|u_0\|_{L^\infty}+\|u_0\|_{\dot C^{\log^\varkappa}}+\sup_{t\in[0,T]}|\log t|^\varkappa (t^{\frac{\kappa}{s}}\|f(t)\|_{\dot C^{\kappa_0}}+t^{\frac{m+\kappa}{s}}\|f(t)\|_{\dot C^{m+\kappa _0}})\right)\\
					&\quad\quad\quad\quad\quad+\tilde \M^\beta \left(\sup_{t\in[0,T]}\int_0^t|\log (t-\tau)|^\varkappa\|g(\tau)\|_{L^\infty}d\tau+\sup_{t\in[0,T]}|\log t|^\varkappa t^{\frac{m}{s}+1}\|g(t)\|_{\dot C^{m}}\right),
				\end{aligned}
			\end{equation}
			where $\kappa_0=\kappa-s+\gamma$, and $\tilde \M,\beta$ are defined in \eqref{deftilM}, \eqref{defbet}, respectively.
		\end{lemma}
		\begin{proof} We already estimate $\sup_{t\in[0,T]}\|u(t)\|_{L^\infty}$ in Theorem \ref{lemmain}. It remains to estimate
			\begin{equation*}
				\|u\|_{X_T^\varkappa }=\sup_{t\in[0,T]}(\|u(t)\|_{\dot C^{\log^\varkappa}}+|\log t|^\varkappa  t^{\frac{\kappa}{s}}\|u(t)\|_{\dot C^{\kappa}}).
			\end{equation*}
			We still use the formula \eqref{uform}, from which we obtain \eqref{1.1hol}, we rewrite here for readers' convenience:
			\begin{equation}\label{deu}
				\begin{aligned}
					\|\delta_\alpha \nabla^n u\|_{L^\infty}\leq \|(\delta_\alpha \nabla^nu_{L,x_0})|_{x_0=x}\|_{L^\infty}+\|(\delta_\alpha \nabla^nu_{N,x_0})|_{x_0=x}\|_{L^\infty}+\|(\delta_\alpha \nabla^nu_{R,x_0})|_{x_0=x}\|_{L^\infty},
				\end{aligned} 
			\end{equation}
			for any $\alpha\in\mathbb{R}^d$ and any $n\in\mathbb{N}$.\\
			We first consider $u_{L,x_0}=\int \K_{x_0}(t,0,z)u_0(x-z)dz$.
			By \eqref{kerL1}, one has
			\begin{align*}
				\|(\delta_\alpha u_{L,x_0})|_{x_0=x}\|_{L^\infty}
				&\lesssim \sup_{x_0\in\mathbb{R}^d}\|\K_{x_0}(t,0,\cdot)\|_{L^1}\|\delta_\alpha u_0\|_{L^\infty}
				\lesssim \tilde\M^\beta \log^{-\varkappa}\left(2+|\alpha|^{-1}\right)\|u_0\|_{\dot C^{\log^\varkappa}}. 
			\end{align*}
			Hence
			\begin{align}\label{logL}
				\sup_{\alpha}	\|(\delta_\alpha u_{L,x_0})|_{x_0=x}\|_{L^\infty}\log^{\varkappa}\left(2+|\alpha|^{-1}\right)\lesssim \tilde\M^\beta \|u_0\|_{\dot C^{\log^\varkappa}}.
			\end{align}
			Moreover, we have 
			\begin{align*}
				|\delta_\alpha \nabla_{x}^{m+[\kappa]}u_{L,x_0}(t,x)|=&\left|\int_{\mathbb{R}^d}\delta_\alpha\nabla_{x}^{m+[\kappa]}\K_{x_0}(t,0,x-z)u_0(z)dz\right|\\
				\lesssim&\int_{\mathbb{R}^d}|\delta_\alpha \nabla_{x}^{m+[\kappa]}\K_{x_0}(t,0,x-z)||u_0(z)-u_0(x)|dz\\
				\lesssim& \int_{\mathbb{R}^d}|\delta_\alpha\nabla_{z}^{m+[\kappa]}\K_{x_0}(t,0,z)|\log^{-\varkappa}\left(2+|x-z|^{-1}\right)dz\|u_0\|_{\dot C^{\log^\varkappa }}\\
				\lesssim&\int_{\mathbb{R}^d}|\delta_\alpha\nabla_{z}^{m+[\kappa]}\K_{x_0}(t,0,z)|\log^{\varkappa}\left(2+(t^{-\frac{1}{s}}|x-z|)\right)dz\|u_0\|_{\dot C^{\log^\varkappa }}|\log(t^{-\frac{1}{s}})|^{-\varkappa},
			\end{align*}
			where in the last inequality we have applied the elementary inequality 
			\begin{equation}\label{ine2}
				|\log(b)|\lesssim\log (2+a^{-1}) \log(2+ab),\ \ \ \forall a>0,b>2.
			\end{equation}
			Observe that  $\log^\varkappa(2+(t^{-\frac{1}{s}}|x|))\lesssim 1+(t^{-\frac{1}{s}}|x|)^\eta$ for $0<\eta\ll 1$. Combining this with \eqref{delKL1} to deduce that
			\begin{align*}
				&\left.\int_{\mathbb{R}^d}|\delta_\alpha\nabla_{z}^{m+[\kappa]}\K_{x_0}(t,0, z)\log(2+(t^{-\frac{1}{s}}| z|))^{\varkappa}|dz\right|_{x_0=x}\\
				&\quad\quad\quad\lesssim \tilde {\M}^\beta \int_{\mathbb{R}^d}|\delta_\alpha\nabla_{z}^{m+[\kappa]}K_{x_0}(t,0, z)|(1+(t^{-\frac{1}{s}}| z|)^\eta)dz\\
				&\quad\quad\quad\lesssim \tilde {\M}^\beta t^{-\frac{m+[\kappa]+\eta}{s}}|\alpha|^{\eta}\min\{1,t^{-\frac{1}{s}}|\alpha| \}^{1-\eta}
				\lesssim \tilde {\M}^\beta |\alpha|^{\kappa-[\kappa]} t^{-\frac{m+\kappa}{s}},
			\end{align*}
			where we use \eqref{leme1}. So we get  
			\begin{align}\label{loghL}
				\sup_\alpha \frac{	\|(\delta_\alpha \nabla_{x}^{m+[\kappa]}u_{L,x_0}(t))|_{x_0=x}\|_{L^\infty}}{|\alpha|^{\kappa-[\kappa]}}\lesssim \tilde {\M}^\beta|\log t|^{-\varkappa}t^{-\frac{m+\kappa}{s}}\|u_0\|_{\dot C^{\log^\varkappa}}.
			\end{align}
			Then we consider $u_{N,x_0}=u_{N,1,x_0}+u_{N,2,x_0}$ as defined in \eqref{formuun}.  First, for $u_{N,1,x_0}$, by \eqref{esdelL} and \eqref{un1x0} we have 
			\begin{align*}
				|\delta_\alpha u_{N,1,x_0}(t,x)|
				\lesssim &\left|\int_0^t \int_{\mathbb{R}^d} \delta_\alpha L_{x_0}(t,\tau,x-z)(\Lambda^{\vartheta}f(\tau,z)-\Lambda^{\vartheta}f(\tau,x))d\tau dz\right|\\
				\lesssim &	\int_0^t\int_{\mathbb{R}^d}| \delta_\alpha L_{x_0}(t,\tau,z)||z|^{\kappa_0-\vartheta}dz\|f(\tau)\|_{\dot{C}^{\kappa_0}}d\tau\\
				\lesssim &\int_0^t \frac{\tilde {\M}^\beta|\alpha|^{\kappa_0-\vartheta}}{(t-\tau)^\frac{[\kappa]+\gamma-\vartheta}{s}|\log \tau|^\varkappa\tau^{\frac{\kappa}{s}}}\min\left\{1,\frac{|\alpha|}{(t-\tau)^\frac{1}{s}}\right\}^{1-(\kappa_0-\vartheta)}d\tau\sup _{t \in[0, T]} |\log t|^\varkappa t^{\frac{\kappa}{s}}\|f(t)\|_{\dot{C}^{\kappa_0}},
			\end{align*}
			where we use \eqref{uNvancd} in the first step and the definition of $L_{x_0}$ as \eqref{mthmdefl}. Note that $|\log t|\lesssim |\log \tau|$ for any $0<\tau<t<\frac{1}{2}$, then it is easy to check that 
			\begin{align*}
				&\int_0^t \frac{|\alpha|^{\kappa_0-\vartheta}}{(t-\tau)^\frac{[\kappa]+\gamma-\vartheta}{s}|\log \tau|^\varkappa\tau^{\frac{\kappa}{s}}}\min\left\{1,\frac{|\alpha|}{(t-\tau)^\frac{1}{s}}\right\}^{1-(\kappa_0-\vartheta)}d\tau\\
				&\quad\quad\quad\lesssim |\log t|^{-\varkappa}\min\left\{1,|\alpha| t^{-1/s}\right\}^\varepsilon\lesssim \log^{-\varkappa} (2+|\alpha|^{-1}).
			\end{align*}
			Hence we obtain that 
			\begin{align}\label{clog1}
				\log^{\varkappa} (2+|\alpha|^{-1}) \|(\delta_\alpha u_{N,1,x_0}(t,x))|_{x_0=x}\|_{L^\infty_x}
				\lesssim &\tilde {\M}^{\beta}\sup _{t \in[0, T]} |\log t|^\varkappa t^{\frac{\kappa}{s}}\|f(t)\|_{\dot{C}^{\kappa_0}}.
			\end{align}
			For $u_{N,2,x_0}$, we have 
			\begin{align*}
				\|\delta_\alpha u_{N,2,x_0}(t,x)|_{x_0=x}\|_{L^\infty_x}&\lesssim \int_0^t\sup_{z} \|\delta_\alpha \K_{z}(t,\tau)\|_{L^1} \|g(\tau)\|_{L^\infty}d\tau\\
				&\lesssim \tilde \M^\beta\log ^{-\varkappa}(2+|\alpha|^{-1}) \int_0^t |\log(t-\tau)|^{\varkappa}\|g(\tau)\|_{L^\infty}d\tau.
			\end{align*}
			Combining this with \eqref{clog1} to obtain 
			\begin{align*}
				&\sup_{t\in[0,T]}\sup_\alpha\log^{\varkappa} (2+|\alpha|^{-1}) \|(\delta_\alpha u_{N,x_0}(t,x))|_{x_0=x}\|_{L^\infty_x}\\
				&\quad\quad\quad
				\lesssim \tilde  \M^\beta\sup_{t\in[0,T]}\left(|\log t|^\varkappa t^{\frac{\kappa}{s}}\|f(t)\|_{\dot C^{\kappa_0}}+\int_0^t |\log(t-\tau)|^{\varkappa}\|g(\tau)\|_{L^\infty}d\tau\right).
			\end{align*}
			For higher order derivative, by \eqref{deun1} and \eqref{deun2}, 
			\begin{align*}
				&|\delta_\alpha\nabla^{m+[\kappa]} u_{N,1,x_0}(t,x)|+|\delta_\alpha\nabla^{m+[\kappa]} u_{N,2,x_0}(t,x)|\\
				&\lesssim \tilde  \M^\beta|\alpha|^{\kappa-[\kappa]} t^{-\frac{m+\kappa }{s}}\left(\int_0^t\|g(\tau)\|_{L^\infty}d\tau+\sup_{\tau\in[0,t]}(\tau^{\frac{m}{s}+1}\|g(\tau)\|_{\dot C^m}+\tau^{\frac{\kappa}{s}}\|f(\tau)\|_{\dot C^{\kappa_0}}+\tau^{\frac{m+\kappa}{s}}\|f(\tau)\|_{\dot C^{\kappa_0+m}})\right).
			\end{align*}
			The above estimates imply that 
			\begin{equation}\label{nlin}
				\begin{aligned}
					&\sup_{t\in[0,T]}\left(\sup_{\alpha}\log^\varkappa(2+|\alpha|^{-1}){\|(\delta_\alpha u_{N,x_0})|_{x_0=x}\|_{L^\infty}}+\sup_{\alpha}|\log t|^\varkappa t^{\frac{m+\kappa}{s}}\frac{\|(\delta_\alpha \nabla^{m+[\kappa]}u_{N,x_0})|_{x_0=x}\|_{L^\infty}}{|\alpha|^{\kappa-[\kappa]}}\right)\\
					&\lesssim \sup_{t\in[0,T]}\left(\int_0^t|\log(t-\tau)|^\varkappa\|g(\tau)\|_{L^\infty}d\tau+|\log t|^\varkappa (t^{\frac{m}{s}+1}\|g(t)\|_{\dot C^m}+t^{\frac{\kappa}{s}}\|f(t)\|_{\dot C^{\kappa_0}}+t^{\frac{m+\kappa }{s}}\|f(t)\|_{\dot C^{\kappa_0+m}})\right)\\
					&\quad\quad\quad\times  \tilde  \M^\beta.
				\end{aligned}
			\end{equation}
			Finally, for the remainder term $u_{R,x_0}$, we obtain from (with a minor change) \eqref{uRHol} that 
			\begin{equation}\label{ur2}
				\begin{aligned}
					\sup_{t\in[0,T]}\sup_\alpha&\left(\log^\varkappa(2+|\alpha|^{-1})\|(\delta_\alpha u_{R,x_0}(t,x))|_{x_0=x}\|_{L^\infty_x}+|\log t|^\varkappa   t^\frac{m+\kappa}{s}\frac{\|\delta_\alpha \nabla_x ^{m+[\kappa]}u_{R,x_0}(t,x)|_{x_0=x}\|_{L^\infty_x}}{|\alpha|^{\kappa-[\kappa]}}\right)\\
					&\lesssim \tilde \M^\beta T^\frac{\varepsilon}{2s} \|u\|_{X_T},
				\end{aligned}
			\end{equation}
			where the norm $\|\cdot\|_{X_T}$ is defined in \eqref{defxt}. The above estimate holds for any $x\in\mathbb{R}^d$. Note that 
			$$
			\|u\|_{X_T}\leq \|u(t)\|_{L_T^\infty L^\infty}+\|u\|_{X_T^\varkappa }.
			$$
			Combining \eqref{logL}, \eqref{loghL},  \eqref{ur2} with \eqref{deu} to obtain that 
			\begin{align*}
				&\|u(t)\|_{L_T^\infty L^\infty}+\|u\|_{X_T^\varkappa }\\
				&\lesssim 
				\tilde  \M^\beta \left(\|u_0\|_{L^\infty}+\|u_0\|_{\dot C^{\log^\varkappa}}+\sup_{t\in[0,T]}|\log t|^\varkappa (t^{\frac{\kappa}{s}}\|f(t)\|_{\dot C^{\kappa-s+\gamma}}+t^{\frac{m+\kappa}{s}}\|f(t)\|_{\dot C^{m+\kappa -s+\gamma}})\right)\\
				&\quad\quad+ \tilde  \M^\beta \left(\sup_{t\in[0,T]}\int_0^t|\log (t-\tau)|^\varkappa\|g(\tau)\|_{L^\infty}d\tau+\sup_{t\in[0,T]}|\log t|^\varkappa t^{\frac{m}{s}+1}\|g(t)\|_{\dot C^{m}}\right)\\
				&\quad\quad+\tilde  \M^\beta T^\frac{\varepsilon}{2s}(\sup_{t\in[0,T]}\|u(t)\|_{L^\infty}+\|u\|_{X_T^\varkappa }).
			\end{align*}
			This implies \eqref{relog} by taking $T$ small enough.
			\vspace{0.3cm}
		\end{proof}\\
		\begin{proof}[Proof of Theorem \ref{thmgm}]  \\
			\textbf{Step 1. Derive the contour dynamic equation} \\
			It is important to find the main terms of the evolution equation \eqref{evo}. We note that the effect of rigid boundaries $ \Gamma^\pm$ and density contrast $\varrho_0 \eta$ are lower order terms. The dominated term comes from  the effect of surface tension  $-\operatorname{div}\left(\frac{\nabla\eta}{\langle\nabla \eta\rangle}\right)$. To see the dominated term in the evolution equation, 
			we first reformulate the system \eqref{sysq}-\eqref{evo} by straightening the interface. More precisely, 
			define \begin{align}\label{relavq}
				v^\pm(x,z)=v^\pm_\eta(x,z)=q^\pm(x,z+\eta(x)).
			\end{align}
			Then 
			\begin{align*}
				&\nabla_{x,z}  q^\pm(x,z)=\M_0(\nabla\eta(x))(\nabla_{x,z}v^\pm)(x,z-\eta(x)),\\ &\Delta_{x,z} q^\pm(x,z)=\left(\operatorname{div}_{x,z}(\M(\nabla\eta)\nabla_{x,z}v^\pm)\right)(x,z-\eta(x)),
			\end{align*}
			with
			\begin{align}\label{ndefM}
				\M_0(\nabla\eta)=\left(\begin{array}{cc}
					\mathrm{	Id}&-\nabla \eta\\
					0&1
				\end{array}\right),\quad\quad\quad \M(\nabla\eta)=(\M_0(\nabla\eta))^\top\M_0(\nabla\eta)=\left(\begin{array}{cc}
					\mathrm{	Id}&-\nabla \eta\\
					-(\nabla \eta)^\top&1+|\nabla\eta|^2
				\end{array}\right).
			\end{align} We obtain
			\begin{equation}\label{elleqeta}
				\begin{aligned}
					&\operatorname{div}_{x,z}(\M(\nabla \eta)\nabla_{x,z} v^\pm )=0,\ \ \text{in}\ \tilde \Omega_\eta^\pm,\\
					&v^--v^+=-\operatorname{div}\left(\frac{\nabla \eta}{\langle\nabla \eta\rangle}\right)+\varrho_0\eta, \quad \text { on } \{z=0\},\\
					&e_{d+1}\cdot (\M(\nabla \eta)(\mu_+^{-1}\nabla_{x,z}v^+-\mu_-^{-1}\nabla_{x,z}v^-))=0, \quad \text { on } \{z=0\},\\
					&\tilde \nu^\pm\cdot (\M(\nabla \eta)\nabla_{x,z} v^\pm) =0,\ \ \ \ \ \ \ \ \ \text{on} \ \tilde \Gamma^\pm_\eta,
				\end{aligned}
			\end{equation}
			where $e_{d+1}$ denotes the (d+1)-th unit vector, $\tilde \nu^\pm$ denote the outward normal vector of $\tilde \Gamma_\eta ^\pm$,  and 
			\begin{equation}\label{defdomain}
				\begin{aligned}
					&\tilde \Omega_{\eta}^{+}=\left\{(x, z) \in \mathbb{R}^{d} \times \mathbb{R}: 0<z<\underline{b}^{+}(x)-\eta(t,x)\right\}, \\
					&\tilde \Omega_{\eta}^{-}=\left\{(x, z) \in \mathbb{R}^{d} \times \mathbb{R}: \underline{b}^{-}(x)-\eta(t, x)<z<0\right\},\\
					&\tilde \Gamma ^\pm_\eta=\{(x,z):z=\underline{b}^\pm(x)-\eta(t,x)\}.
				\end{aligned} 
			\end{equation}
			Note that we transform the effect of free interface $\Sigma_t=\{(x,z):z=\eta(x)\}$ in the system \eqref{sysq} to the variable coefficient $\M(\nabla\eta)$ in the system \eqref{elleqeta}. 
			We remark that $v^\pm$ is not differentiable at $z=0$. When evaluating a function with a ``$+$" or ``$-$" superscript at $z=0$, it implies approaching the limit from either the positive or negative direction. \textit{i.e.} $ v^\pm(x,0)=\lim_{z=0^\pm}v^\pm(x,z)$.\\
			The evolution equation \eqref{evo} becomes
			\begin{equation}\label{inievo}
				\begin{aligned}
					\partial_t \eta (x)&=-\frac{1}{\mu_-}e_{d+1}\cdot (\M(\nabla\eta)\nabla_{x,z}v^-(x,z)|_{z=0})
					\\&=\frac{1}{\mu_-}\nabla \eta(x)\cdot \nabla _xv^-(x,z)|_{z=0} -\frac{1}{\mu_-}\langle\nabla \eta(x)\rangle^2\partial_zv^-(x,z)|_{z=0},\\
					\eta|_{t=0}&=\eta_0.
				\end{aligned}
			\end{equation}
			We remark that \eqref{inievo} is well-defined because of the existence and uniqueness of weak solution to the system \eqref{elleqeta}.
			\begin{remark}
				The elliptic system \eqref{elleqeta} can be written as the following elliptic equation with jump coefficient and measure data:
				\begin{align*}
					&\operatorname{div}_{x,z}(\bar \M(\nabla\eta)\nabla_{x,z}v)=\delta_{z=0}f,\ \ \ \text{in}\ \tilde \Omega_\eta= \{(x,z):\underline{b}^{-}(x)-\eta(t,x)<z<\underline{b}^{+}(x)-\eta(t,x)\},\\
					&\tilde \nu^\pm\cdot (\bar \M(\nabla\eta)\nabla_{x,z}v)=0,\ \ \ \ \ \ \text{on}\ \tilde\Gamma _\eta^\pm.
				\end{align*}
				Here the coefficient matrix $\bar \M(\nabla\eta)\in L^\infty$ is defined by 
				\begin{align*}
					\bar \M(\nabla\eta)=\begin{cases}
						\mu_+^{-1}\M(\nabla\eta),\quad 0<z<\underline{b}^{+}(x)-\eta(t,x),\\
						\mu_-^{-1}\M(\nabla\eta),\quad \underline{b}^{-}(x)-\eta(t,x)<z<0.
					\end{cases}
				\end{align*}
				And $f(x)=-\mu_+^{-1}\M(\nabla\eta)\nabla v_0(x,0)$, where $v_0$ solves the equation 
				\begin{align*}
					&\operatorname{div}_{x,z}(\M(\nabla\eta)\nabla_{x,z}v_0^\pm)=0,\ \ \ \text{in}\ \tilde \Omega_\eta,\\
					&v_0=-\operatorname{div}\left(\frac{\nabla\eta}{\langle\nabla\eta\rangle}\right)+\varrho_0\eta,\ \ \ \text{on}\ \{z=0\},\\
					&\tilde \nu^+\cdot (\M(\nabla\eta)\nabla_{x,z}v_0)=0,\ \ \text{on}\ \tilde \Gamma_\eta^\pm.
				\end{align*}
				The existence of weak solution to elliptic systems has been well-studied in literature(see \cite{Ag}).  \vspace{0.2cm} 
			\end{remark}
			To see the main terms in the evolution equation, we consider the following system.
			For any fixed $\b\in\mathbb{R}^{d}$, let $\tilde v^\pm_\b=\tilde v^\pm_\b[\mathcal{G}]$ be a solution to the following system
			\begin{equation}\label{eqmain}
				\begin{aligned}
					&\operatorname{div}_{x,z}(\M(\b)\nabla_{x,z}\tilde v^\pm_\b )=0,\ \ \text{in}\ \mathbb{R}^{d+1}_\pm,\\
					&\tilde v^-_\b-\tilde v^+_\b=\mathcal{G} \quad \text { on } \{z=0\},\\
					&e_{d+1}\cdot (\M(\b)(\mu_+^{-1}\nabla_{x,z}\tilde v^+_\b-\mu_-^{-1}\nabla_{x,z}\tilde v^-_\b))=0 \quad \text { on } \{z=0\},
				\end{aligned}
			\end{equation}
			where we denote $\mathcal{G}=\mathcal{G}[\eta]=-\operatorname{div}\left(\frac{\nabla \eta}{\langle\nabla \eta\rangle}\right)$.
			Take Fourier transform with respect to $x$ variable, and denote 
			\begin{align*}
				V^\pm_\b(\xi,z)=(2\pi)^{-\frac{d}{2}}\int_{\mathbb{R}^d}\tilde v_\b ^\pm(x,z)e^{-i\xi\cdot x}dx,\ \ \ \xi\in\mathbb{R}^d,\ \ \ \quad z\in\mathbb{R}^\pm.
			\end{align*}
			One has 
			\begin{align*}
				&	(1+|\b|^2)\partial_z^2V^\pm_\b(\xi,z)-|\xi|^2V^\pm_\b(\xi,z)-2i\b\cdot \xi \partial_z V_\b^\pm(\xi,z)=0,\\
				&	V^-_\b(\xi,0)-V^+_\b(\xi,0)= \hat{\mathcal{G}}(\xi),\\
				&\left(-\frac{ V^+_\b}{\mu_+}+\frac{V^-_\b}{\mu_-}\right)i\xi\cdot \b+\left(\frac{\partial_z V^+_\b}{\mu_+}-\frac{\partial_zV^-_\b}{\mu_-}\right)(|\b|^2+1)=0,\ \ \ \text{on}\ \{z=0\}.
			\end{align*}
			By a direct calculation, we obtain the eigenvalues of the above ODE:
			\begin{align}\label{deflanb}
				\lambda^\pm=	\lambda^\pm(\xi,\b)=\frac{i\b\cdot \xi\pm \sqrt{\langle \b\rangle^2|\xi|^2-(\b\cdot \xi)^2}}{\langle \b\rangle^2}.
			\end{align}
			We consider weak solutions that decay at  $z=\infty$,
			hence	we have the formula
			\begin{align}\label{deffour}
				V^+_\b(\xi,z)=c^+\hat{\mathcal{G}}(\xi)e^{\lambda^-z}, \ \ \ \ \ V^-_\b(\xi,z)=c^-\hat{\mathcal{G}}(\xi)e^{\lambda^+z},
			\end{align}
			where 
			\begin{align*}
				c^\pm=\frac{\pm\mu_\pm  }{\mu_++\mu_-}.
			\end{align*}
			Taking inverse Fourier transform to obtain 
			\begin{align*}
				\tilde v^\pm_\b(x,z)=(2\pi)^{-\frac{d}{2}}c^\pm \int_{\mathbb{R}^d} \hat {\mathcal{G}}(\xi)e^{\lambda^\mp(\xi,\b) z}e^{ix\cdot \xi}d\xi,\quad\quad\quad z\in\mathbb{R}^\pm.
			\end{align*}
			More precisely, we have
			\begin{align*}
				\tilde v^\pm_{\b}(x,z)=c^\pm \mathcal{V}_\b[\mathcal{G}](x,z),\quad\quad z\in\mathbb{R}^\pm,
			\end{align*}
			where 
			\begin{align}\label{defvb}
				\mathcal{V}_\b[\mathcal{G}](x,z)=\int_{\mathbb{R}^d} K_\b(x-y,z) \mathcal{G}(y)dy,
			\end{align}
			with
			\begin{equation}\label{defKor}
				\begin{aligned}
					K_\b(x,z)=\begin{cases}
						\frac{1}{(2\pi)^\frac{d}{2}}\int_{\mathbb{R}^d}\exp\left(ix\cdot\xi+\lambda^-(\xi,\b)z\right)d\xi,& z\geq0,\\
						\frac{1}{(2\pi)^\frac{d}{2}}\int_{\mathbb{R}^d}\exp\left(ix\cdot\xi+\lambda^+(\xi,\b)z\right)d\xi,& z<0.
					\end{cases}	
				\end{aligned}
			\end{equation}
			See Lemma \ref{Z10} in Appendix for the explicit formula of $	K_\b(x,z)$.
			Moreover, it is easy to check that 
			\begin{align}\label{dx0}
				\frac{\left(	\nabla_x\tilde v_\b^-[\mathcal{G}]\right)(x,0)}{\mu_-}=\frac{\nabla \mathcal{G}(x)}{\mu_++\mu_-}.
			\end{align}
			To obtain $ \partial_z\tilde v^-_\b|_{(x,0)}$, 
			we apply Lemma \ref{lemfdz0} to obtain that 
			\begin{equation}\label{hz0}
				\begin{aligned}
					\lim_{z\to 0^-}	\int_{\mathbb{R}^d} \partial_zK_\b(x-y,z)\mathcal{G}(y)dy&=\frac{1}{(2\pi)^\frac{d}{2}}\int_{\mathbb{R}^d}\lambda^+(\xi,\b)\hat {\mathcal{G}}(\xi)\exp\left(ix\cdot\xi\right)d\xi \\
					&=\tilde c_d\int_{\mathbb{R}^d} \frac{\delta_\alpha \mathcal{G}(x)}{\langle \hat\alpha\cdot \b\rangle^{d+1}}\frac{d\alpha}{|\alpha|^{d+1}}+\frac{\b\cdot\nabla \mathcal{G}(x)}{\langle \b\rangle^2}.
				\end{aligned}
			\end{equation}
			Here $\tilde c_d$ is the constant in \eqref{contc}.
			This implies that
			\begin{align}\label{z0}
				\frac{	\left(\partial_z \tilde v_\b^-[\mathcal{G}]\right)(x,0)}{\mu_-}=\frac{1}{\mu_++\mu_-}\left(\tilde c_d \int _{\mathbb{R}^d}\frac{\delta_\alpha\mathcal{G}(x)}{\langle\hat \alpha \cdot \b\rangle^{d+1}}\frac{d\alpha}{|\alpha|^{d+1}}+\frac{\b\cdot \nabla \mathcal{G}(x)}{\langle \b\rangle^2}\right).
			\end{align}
			For simplicity, denote 
			\begin{align}\label{ndefB}
				B(\b)=\frac{\mathrm{Id}}{\langle \b\rangle}-\frac{\b\otimes \b}{\langle\b\rangle^3}, \ \ \ \forall \b\in\mathbb{R}^d.
			\end{align}
			It is easy to check that \begin{align}\label{Bellip}
				\frac{1}{\langle \b\rangle^3} \mathrm{Id}\leq   B(\b)\leq 	\frac{1}{\langle \b\rangle} \mathrm{Id}.
			\end{align} And   \begin{align}\label{ndefG}
				\mathcal{G}[f]=-\operatorname{div}\left(\frac{\nabla f }{\langle\nabla f \rangle}\right) =-B(\nabla f):\nabla ^2 f=-\frac{\Delta f}{\langle \nabla f\rangle}+\frac{(\nabla f)^\top  \nabla ^2 f\nabla f}{\langle \nabla f\rangle^3}.
			\end{align}
			
			We add subscript in order to emphasize the dependence of the solution to elliptic systems on data. More precisely,  we denote $v_h^\pm(x,z)$ the solution to system \eqref{elleqeta} with data $\mathcal{G}=\mathcal{G}[h]$, 
			for any $h: \mathbb{R}^d\to \mathbb{R}$.
			With a slight abuse of notation, we simply denote 
			\begin{align}\label{deftvb}
				\tilde v^\pm_{\b}[\eta](x,z)=c^\pm  {\mathcal{V}}_{\b}[\mathcal{G}[\eta]](x,z).
			\end{align}
			Denote
			\begin{align}\label{defphigm}
				\phi=\eta_0\ast \rho_{\eps_1}.
			\end{align}
			By the condition \eqref{gmcon}, we have $\phi\in C^\infty(\mathbb{R}^d)\cap L^2(\mathbb{R}^d)$, and 
			\begin{align}\label{gmcon1}
				\|\eta_0-\phi\|_{L^2\cap \dot C^{1,\log^\varkappa  }}\leq \varepsilon_0.
			\end{align}
			We will approximate $(\nabla_{x,z}v^-_\eta)(x,0)$ by $\left(\nabla_{x,z}\tilde v^-_{\b}[\eta])(x,0)\right|_{\b=\nabla \phi(x)}$ in the evolution equation \eqref{inievo}.
			
			Now we fix  $\b\in\mathbb{R}^d$.
			The evolution equation \eqref{inievo} can be written as 
			\begin{equation}\label{eqeta1}
				\begin{aligned}
					\partial_t \eta (x)=&-\frac{1}{\mu_-}e_{d+1}\cdot (\M(\nabla\eta(x))(\nabla_{x,z}\tilde v^-_{\b}[\eta])(x,z))|_{z=0}\\
					&\quad\quad\quad-\frac{1}{\mu_-}e_{d+1}\cdot\left. \left(\M(\nabla\eta(x))\big(\nabla_{x,z}v^-_\eta(x,z)-\nabla_{x,z}\tilde v^-_{\b}[\eta](x,z)\big)\right)\right|_{z=0}.
				\end{aligned}
			\end{equation}
			By \eqref{dx0}, \eqref{z0}, the first term can be written as
			\begin{equation}\label{spli4}
				\begin{aligned}
					&-\frac{1}{\mu_-}e_{d+1}\cdot (\M(\nabla\eta(x))(\nabla_{x,z}\tilde v^-_{\b}[\eta])(x,z))|_{z=0}\\
					&=\frac{1}{\mu_-}\left.\left(\nabla \eta(x)\cdot \nabla_x \tilde v^-_{\b}[\eta](x,z)-\langle \nabla \eta(x)\rangle^2\partial_z\tilde v^-_{\b}[\eta](x,z)\right)\right|_{z=0}\\
					&=\frac{\langle\nabla \eta(x)\rangle^2 }{\mu_++\mu_-}\left(-\tilde c_d\int_{\mathbb{R}^d}\frac{\delta_\alpha \mathcal{G}[\eta](x)}{\langle\hat \alpha \cdot \mathbf{b}\rangle^{d+1}}\frac{d\alpha}{|\alpha|^{d+1}}+\left(\frac{\nabla \eta(x)}{\langle\nabla \eta(x)\rangle^2 }-\frac{\mathbf{b}}{\langle\mathbf{b}\rangle^2}\right)\cdot \nabla\mathcal{G}[\eta](x)\right)\\
					&=\int_{\mathbb{R}^d}\A(\b,\alpha):\delta_\alpha \nabla^2\eta(x)\frac{d\alpha}{|\alpha|^{d+1}}+\int_{\mathbb{R}^d}\frac{(\mathbf{B}(\nabla \eta(x))-\mathbf{B}(\mathbf{b})):\delta_\alpha \nabla^2\eta(x)}{\langle\hat \alpha \cdot \mathbf{b}\rangle^{d+1}}\frac{d\alpha}{|\alpha|^{d+1}}\\
					&\ \quad\quad\quad\quad+\frac{\tilde c_d\langle\nabla \eta(x)\rangle^2 }{\mu_++\mu_-}\int_{\mathbb{R}^d}\frac{\delta_\alpha B(\nabla \eta): \nabla^2\eta(x-\alpha)}{\langle\hat \alpha \cdot \mathbf{b}\rangle^{d+1}}\frac{d\alpha}{|\alpha|^{d+1}}+\frac{\langle\nabla \eta(x)\rangle^2 }{\mu_++\mu_-}\left(\frac{\nabla \eta(x)}{\langle\nabla \eta(x)\rangle^2 }-\frac{\mathbf{b}}{\langle \mathbf{b}\rangle^2}\right)\cdot \nabla\mathcal{G}[\eta](x),
				\end{aligned}
			\end{equation}
			where we denote \begin{align}
				\label{ndefbfb}
				\A(\b,\alpha)=\frac{\mathbf{B}(\b)}{\langle\hat \alpha \cdot \mathbf{b}\rangle^{d+1}},\quad\quad\quad	\mathbf{B}(\mathbf{b})=\frac{\tilde c_d }{\mu_++\mu_-}{\langle \mathbf{b}\rangle^2}B(\mathbf{b}).
			\end{align}	
			Now we extract the dominated term of the evolution equation, which is the first term in the right hand side of \eqref{spli4} with $\b=\nabla\phi(x)$, where $\phi$ defined in \eqref{defphigm} is the mollification of the initial interface. More precisely,
			we rewrite the system \eqref{eqeta1} as
			\begin{align}\label{evoeta}
				\partial_t \eta (x)-\int _{\mathbb{R}^d}\A(\nabla\phi(x),\alpha):\delta_\alpha \nabla^2 \eta(x)\frac{d\alpha}{|\alpha|^{d+1}}=\H[\eta](x).
			\end{align} 
			Here
			\begin{equation}\label{defH}
				\begin{aligned}
					&	\left.\begin{aligned}	\H[\eta](x)=&\int _{\mathbb{R}^d}(\A(\b,\alpha)-\A(\nabla\phi(x),\alpha)):\delta_\alpha \nabla^2 \eta(x)\frac{d\alpha}{|\alpha|^{d+1}}\\
						&	+		\int_{\mathbb{R}^d}\frac{(\mathbf{B}(\nabla \eta(x))-\mathbf{B}(\mathbf{b})):\delta_\alpha \nabla^2\eta(x)}{\langle\hat \alpha \cdot \mathbf{b}\rangle^{d+1}}\frac{d\alpha}{|\alpha|^{d+1}}\\
						&	+\frac{\tilde c_d\langle\nabla \eta(x)\rangle^2 }{\mu_++\mu_-}\int_{\mathbb{R}^d}\frac{\delta_\alpha B(\nabla \eta): \nabla^2\eta(x-\alpha)}{\langle\hat \alpha \cdot \mathbf{b}\rangle^{d+1}}\frac{d\alpha}{|\alpha|^{d+1}}\\
						&	+\frac{\langle\nabla \eta(x)\rangle^2 }{\mu_++\mu_-}\left(\frac{\nabla \eta(x)}{\langle\nabla \eta(x)\rangle^2 }-\frac{\mathbf{b}}{\langle \mathbf{b}\rangle^2}\right)\cdot \nabla\mathcal{G}[\eta](x)
					\end{aligned}\ \quad\quad\ \ \ \quad\quad\quad \right\}:=	\H^{ex}_\b[\eta](x)\\
					&\quad\quad\quad\quad\ 	\left.-\frac{1}{\mu_-}e_{d+1}\cdot\left. \left(\M(\nabla\eta(x))\left((\nabla_{x,z}v^-_\eta)(x,z)-(\nabla_{x,z}\tilde v^-_{\b}[\eta])(x,z)\right)\right)\right|_{z=0}\right\}:=	\H^{im}_\b[\eta](x).
				\end{aligned}
			\end{equation}
			Note that the first four terms in  \eqref{defH} have explicit formula that only depend on $\eta$. And the last term depends on $\eta$ implicitly, which involves the solution $v^\pm_\eta$ to the elliptic system \eqref{elleqeta}. Moreover, $\H[\eta]$ itself is not related to the choice of $\b$.\\
			Denote 
			\begin{align}\label{defwb}
				\omega_\mathbf{b}^\pm[\eta]=v^\pm_\eta-\tilde v^\pm_{\b}[\eta].
			\end{align} Then 
			\begin{align}\label{defHim}
				\H^{im}_\b[\eta](x)=-\frac{1}{\mu_-}e_{d+1}\cdot (\M(\nabla\eta(x))(\nabla_{x,z}\omega^-_\b[\eta])(x,z))|_{z=0}.
			\end{align} 
			And $\omega_\mathbf{b}^\pm=\omega_\mathbf{b}^\pm[\eta]$ solves the elliptic system
			\begin{equation}	\label{elliw}
				\begin{aligned}
					&\operatorname{div}_{x,z}(\M(\nabla \eta)\nabla_{x,z} \omega^\pm_{\b})=\operatorname{div} F_{1,\b}^\pm,\ \ \text{in}\ \tilde\Omega_\eta^\pm,\\
					&\omega^-_{\b}-\omega^+_{\b}=F_2, \quad \text { on } \{z=0\},\\
					&e_{d+1}\cdot (\M(\nabla \eta)\nabla_{x,z} (\mu_+^{-1}\omega_\mathbf{b}^+-\mu_-^{-1}\omega^-_{\b}))=F_{3,\b}, \quad \text { on } \{z=0\},\\
					&\tilde \nu^\pm\cdot (\M(\nabla \eta)\nabla _{x,z} \omega^\pm_{\b}) =F_{4,\b},\ \ \quad\quad\quad\quad\ \text{on}\  \tilde \Gamma_\eta^\pm,
				\end{aligned}
			\end{equation}	
			where 
			\begin{equation*}
				\begin{aligned}
					&F_{1,\b}^\pm=\left(\M(\mathbf{b})-\M(\nabla \eta)\right)\nabla_{x,z} \tilde v^\pm_\b[\eta],\\
					&F_2=\varrho_0\eta,\\
					&F_{3,\b}=\left.e_{d+1}\cdot \left((\mu_+^{-1}F_{1,\b}^+-\mu_-^{-1}F_{1,\b}^-)\right)\right|_{z=0},\\
					&F_{4,\b}^\pm=\left.\tilde \nu^\pm\cdot (\M(\nabla \eta)(\nabla_{x,z}  v^\pm_\eta-\nabla_{x,z}\tilde v_{\mathbf{b}}^\pm[\eta])) \right|_{ \tilde \Gamma_\eta^\pm}.
				\end{aligned}
			\end{equation*}
			\begin{remark}
				Note  that
				\begin{align*}
					-\int _{\mathbb{R}^d}\A(\nabla\phi(x),\alpha):\delta_\alpha \nabla^2 \eta(x)\frac{d\alpha}{|\alpha|^{d+1}}=(2\pi)^{-\frac{d}{2}}\int_{\mathbb{R}^d} \tilde \A(x,\xi)\hat \eta(\xi)e^{ix\cdot \xi}d\xi,
				\end{align*}
				where 
				\begin{align*}
					\tilde 	\A(x,\xi)=&\int_{\mathbb{R}^d}\A(\nabla\phi(x),\alpha):(\xi\otimes\xi)(1-e^{-i\alpha\cdot\xi})\frac{d\alpha}{|\alpha|^{d+1}}\\
					=&\frac{1}{2}\int_{\mathbb{R}^d}\A(\nabla\phi(x),\alpha):(\xi\otimes\xi)(2-e^{i\alpha\cdot\xi}-e^{-i\alpha\cdot\xi})\frac{d\alpha}{|\alpha|^{d+1}}.
				\end{align*}
				Here we used the fact that $\A(\mathbf{b},\alpha)=\A(\mathbf{b},-\alpha)$ in the last equality. By the ellipticity property \eqref{Bellip},
				one can check that
				\begin{align*}
					&\frac{\tilde c_d}{(\mu_++\mu_-)\langle \nabla\phi(x)\rangle^{d+2}}\mathrm{Id}\leq \frac{\tilde \A(x,\xi)}{|\xi|^3}\leq 	\frac{\tilde c_d}{(\mu_++\mu_-)\langle \nabla\phi(x)\rangle}\mathrm{Id},\\ &|\nabla_x^{n_1}\nabla_{\xi}^{n_2}\tilde \A(x,\xi)|\lesssim (1+\|\nabla\phi\|_{C^{n_1}})^{n_1+10} |\xi|^{3-n_2},\ \forall n_1,n_2\in\mathbb{N}.
				\end{align*}
				Hence the operator on the left hand side of \eqref{evoeta} satisfies condition \eqref{condop} with $s=3$, $c_0=\frac{\tilde c_d}{(\mu_++\mu_-)(1+\|\nabla\phi\|_{L^\infty})^{d+2}}$ and $\M=(1+\|\nabla\phi\|_{C^{d+m+4}})^{d+m+10}$.
				
			\end{remark}\vspace{0.3cm}
			\textbf{Step 2. Construct the solution map}\\
			Fix the constants 
			\begin{align}\label{consgm}
				\varkappa>1,\ \ \ m\in\mathbb{N},\ \ \ \kappa\in(2,3)\ \text{such that}\ 0<3-\kappa\ll 1.
			\end{align}
			For $h:[0,T]\times \mathbb{R}^d\to \mathbb{R}$, define
			\begin{equation}	\label{defnorgm}
				\begin{aligned}
					&\|h\|_T:=\sup_{t\in[0,T]}\left(\|\nabla h(t)\|_{L^\infty}+\|\nabla h(t)\|_{ \dot C^{\log ^\varkappa}}+t^\frac{m+\kappa}{3}|\log t|^\varkappa  \|\nabla h(t)\|_{\dot C^{m+\kappa}}\right),\\
					&\|h\|_{X_T}:=\|h\|_{L^\infty_TL^2}+\| h\|_T.
				\end{aligned}
			\end{equation}
			Moreover, for a function $\phi$ on $\mathbb{R}^d$, we can define $\|\phi\|_T, \|\phi\|_{X_T}$ by considering it as a function defined on $[0,T]\times \mathbb{R}^d$ and is invariant with respect to the time variable.\\
			Let	$\sigma,T>0$ that will be fixed later.	Consider the set 
			\begin{align}\label{defsp}
				\mathcal{X}^\sigma_{T,\phi}=\left\{f\in L^\infty_T\dot W^{1,\infty}:\|f-\phi\|_{X_T}\leq \sigma, f|_{t=0}=\eta_0, \operatorname{dist}(f,\Gamma^\pm)>\mathbf{r}\right\}.
			\end{align}
			For a function $f\in	\mathcal{X}^\sigma_{T,\phi}$,
			we construct a map $g=\mathcal{S}f$, where $g$ is the solution to the Cauchy problem
			\begin{equation}\label{infaeq1}
				\begin{aligned}
					&\partial_t g(x) -\int_{\mathbb{R}^d} \A(\nabla\phi(x),\alpha):\delta_\alpha \nabla^2 g(x)\frac{d\alpha}{|\alpha|^{d+1}}=\H[f](x),\\
					&g|_{t=0}=\eta_0.
				\end{aligned}
			\end{equation}
			It is clear that a fixed point of $\mathcal{S}$ is a solution to \eqref{evoeta}. We will prove that the map $\mathcal{S}$ has a unique fixed point in $\mathcal{X}^\sigma_{T,\phi}$ with suitable choice of $\sigma, T$.\vspace{0.3cm}\\
			\textbf{Step 3. Existence and regularity}\\
			By adding some terms on both sides of \eqref{infaeq1}, and taking one spatial derivative $\partial_i=\partial_{x_i}$, $i=1,2,\cdots,d$, we obtain 
			\begin{equation}\label{infaeq11}
				\begin{aligned}
					&\partial_t \partial_i (g-\phi)(x)-\int_{\mathbb{R}^d}\A(\nabla \phi(x),\alpha):{\delta_\alpha \nabla ^2\partial_i}(g-\phi)(x)\frac{d\alpha}{|\alpha|^{d+1}}=\partial_i \H[f](x)+R_i[g](x),\\
					&\partial_i (g-\phi)|_{t=0}=\partial_i(\eta_0-\phi),
				\end{aligned}  
			\end{equation}
			where
			\begin{align}\label{defR}
				&R_i[g](x)=\int _{\mathbb{R}^d}\partial_{x_i}\left(\A(\nabla \phi(x),\alpha)\right):\delta_\alpha \nabla^2 (g-\phi)(x)\frac{d\alpha}{|\alpha|^{d+1}}+\partial_i\left(\int_{\mathbb{R}^d}\A(\nabla \phi(x),\alpha):{\delta_\alpha \nabla ^2}\phi(x)\frac{d\alpha}{|\alpha|^{d+1}}\right).
			\end{align}
			We first estimate $\| g-\phi\|_T$. Applying Theorem \ref{lemlog} with $f=\H[f]$ and $g=R_i[g]$ to obtain that there exists $T=T(\|\nabla\phi\|_{C^{m+4}})>0$ such that
			\begin{equation}\label{mapes}
				\begin{aligned}
					\|g-\phi\|_T\lesssim &\|\eta_0-\phi\|_{ C^{1,\log^\varkappa  }}+\sup_{t\in[0,T]}|\log t|^\varkappa  \left(t^\frac{\kappa}{3}\| \H[f](t)\|_{\dot C ^{\kappa-2}}+t^\frac{m+\kappa}{3}\| \H[f](t)\|_{\dot C ^{m+\kappa-2}}\right)\\
					&+\sup_{t\in[0,T]}\left(\int_0^t{|\log(t-\tau)|^\varkappa\|R[g](\tau)\|_{L^\infty}}d\tau+|\log t|^\varkappa  t^{\frac{m}{3}+1}\|R[g](t)\|_{\dot C^{m}}\right).
				\end{aligned}
			\end{equation}
			In the following we use the notation \eqref{defnormbb}.
			Then for $j=0,m$,
			\begin{align*}
				\|\H[f]\|_{\dot C^{j+\kappa-2}}\lesssim  |||\nabla^j\H^{ex}_\b[f]|||_{\kappa-2}+|||\nabla^j\H^{im}_\b[f]|||_{\kappa-2}. 
			\end{align*}
			We have explicit formula of $\H^{ex}$. The main difficulty is to deal with the implicit term $\H^{im}_\b[f]$. \\
			By Lemma \ref{lemH1} and Lemma \ref{lemR} we have 
			\begin{equation}\label{H11}
				\begin{aligned}
					&\sup_{t\in[0,T]}|\log t|^\varkappa  \left(t^\frac{\kappa}{3}||| 
					\H_\b^{ex}[f](t)|||_{\kappa-2}+t^\frac{m+\kappa}{3}||| 
					\nabla^m\H_\b^{ex}[f](t)|||_{\kappa-2}\right)\\
					&\quad\quad+\sup_{t\in[0,T]}\left(\int_0^t|\log(t-\tau)|^\varkappa\|R[g](\tau)\|_{L^\infty}d\tau+|\log t|^\varkappa t^{\frac{m}{3}+1}\|R[g](t)\|_{\dot C^m}\right)\\
					&\lesssim (\|f-\phi\|_{X_T}+\|g-\phi\|_{X_T})(\| f-\phi\|_{X_T}+T^\frac{1}{20}\|\nabla \phi\|_{C^{m+4}})(1+\|f-\phi\|_{X_T}+T^\frac{1}{20}\|\nabla \phi\|_{C^{m+4}})^{m+3},
				\end{aligned}
			\end{equation}
			and from Lemma \ref{lemh2} one has
			\begin{equation}\label{H22}
				\begin{aligned}
					&\sup_{t\in[0,T]}|\log t|^\varkappa  (t^{\frac{\kappa}{3}}|||\H_\b^{im}[f]|||_{\kappa-2}+t^{\frac{m+\kappa}{3}}|||\nabla^m\H_\b^{im}[f]|||_{\kappa-2})\\
					&\lesssim \sup_{t\in[0,T]}|\log t|^\varkappa  \Big\{t^{\frac{\kappa}{3}}(\|\nabla f\|_{\dot C^{\kappa-2}}\|\nabla _{x,z}\omega^-_\b(\cdot,0)\|_{L^\infty}+\|\nabla f\|_{L^\infty}|||\nabla _{x,z}\omega^-_\b(\cdot,0)|||_{\kappa-2})\\
					&\quad\quad+t^{\frac{m+\kappa}{3}}\sum_{m_1+m_2=m}(\|\nabla^{1+m_1} f\|_{\dot C^{\kappa-2}}\|\nabla_x^{m_2}\nabla _{x,z}\omega^-_\b(\cdot,0)\|_{L^\infty}+\|\nabla^{1+m_1} f\|_{L^\infty}|||\nabla_x^{m_2}\nabla _{x,z}\omega^-_\b(\cdot,0)|||_{\kappa-2})\Big\}\\
					&\lesssim ( |\log T|^{-\varkappa } +\|f-\phi\|_{X_T})\|f-\phi\|_{X_T}+|\log T|^{-\varkappa },
				\end{aligned}
			\end{equation}
			where the implicit constants depend on $\|\phi\|_{L^2}$, $\|\nabla\phi\|_{C^{d+m+4}}$, and $\mathbf{r}$.
			Substituting \eqref{H11} and \eqref{H22} into \eqref{mapes}, and apply \eqref{gmcon1} we obtain 
			\begin{equation}\label{gpht}
				\begin{aligned}
					\|g-\phi\|_T&\leq C_1\|\eta_0-\phi\|_{\dot C^{1,\log^\varkappa  }}+C_1|\log T|^{-\varkappa }(1+\|f-\phi\|_{X_T})\\
					&\quad\quad+C_1 (\|f-\phi\|_{X_T}+\|g-\phi\|_{X_T})(\| f-\phi\|_{X_T}+T^\frac{1}{20}\|\nabla \phi\|_{C^{m+4}})\\
					&\quad\quad\quad\quad\times(1+\|f-\phi\|_{X_T}+T^\frac{1}{20}\|\nabla \phi\|_{C^{m+4}})^{m+3}\\
					&\leq C_1(\eps_0+(\sigma+\|g-\phi\|_{X_T})(\sigma+T^\frac{1}{50})+|\log T|^{-\varkappa }).
				\end{aligned}
			\end{equation}
			Here we obtain the last inequality from $\|f-\phi\|_{X_T}\leq \sigma<1$ and by taking $T<\left(\frac{\eps_1}{100+\|\eta_0\|_{\dot W^{1,\infty}}}\right)^{100(d+m+10
				)}$ such that $T^\frac{1}{50}\|\nabla \phi\|_{C^{d+m+4}}\leq C T^\frac{1}{50}\eps_1^{-(d+m+4)}\|\eta_0\|_{\dot W^{1,\infty}}\leq 1$.\\
			This completes the control of $\|\cdot\|_T$ norm, it remains to estimate the lower order $L^2$ norm. Thanks to Lemma \ref{mapl2}, we have
			\begin{align}\label{gphl2}
				\|g-\phi\|_{L^\infty_TL^2}\leq  \|\eta_0-\phi\|_{L^2}+C_1T^\frac{1}{50}\leq \eps_0+C_1T^\frac{1}{50}.
			\end{align}
			We note that $C_1=C_1(\|\phi\|_{L^2},\|\nabla \phi\|_{C^{d+m+4}},\mathbf{r},)$.
			With the definition of $\phi$ in \eqref{defphigm}, we have $C_1=C_1(\|\eta_0\|_{L^2\cap \dot C^{1,\log^\varkappa  }},\eps_1,\mathbf{r})$.
			Then we obtain from \eqref{gpht} and \eqref{gphl2} that 
			\begin{equation}\label{ineeeq1}
				\begin{aligned}
					\|g-\phi\|_{X_T}&=	\|g-\phi\|_{L^\infty_TL^2}+\|g-\phi\|_T\\
					&\leq (1+C_1)\eps_0+C_1(\sigma+\|g-\phi\|_{X_T})(\sigma+T^\frac{1}{50})+C_1(|\log T|^{-\varkappa}+T^\frac{1}{50}).
				\end{aligned}
			\end{equation}
			We take $\varepsilon_0\leq \frac{1}{100 (1+C_1)^2}$, $\sigma=4(1+C_1)\eps_0$, and $T$ small enough such that $T^\frac{1}{50}+|\log T|^{-\varkappa } \leq \frac{\eps_0}{10}$. Then 
			\eqref{ineeeq1} yields
			\begin{align*}
				\|g-\phi\|_{X_T}\leq \frac{\sigma}{2}.
			\end{align*}
			Finally, we prove the propagation of distance  between the surface and the bottom. We have 
			\begin{align*}
				\|g(t)-\eta_0\|_{L^\infty_TL^\infty}\leq \|\eta_0-\phi\|_{L^\infty}+\|g(t)-\phi\|_{L^\infty_TL^\infty}\lesssim \|\eta_0-\phi\|_{L^2\cap \dot C^{1,\log^\varkappa  }}+\|g-\phi\|_{X_T}\leq 2\sigma.
			\end{align*}
			Then we obtain 
			$$
			\operatorname{dist}(g,\Gamma^\pm)\geq 	\operatorname{dist}(g,\Gamma^\pm)-\|g(t)-\eta_0\|_{L^\infty_TL^\infty}\geq 2\mathbf{r}-2\sigma\geq \mathbf{r},
			$$
			by taking $\sigma \leq \frac{\mathbf{r}}{100}$.
			This implies that $\mathcal{S}$ maps $\mathcal{X}^\sigma_{T,\phi}$ to itself.\\
			To prove the contraction estimate, we consider $f_1,f_2\in \mathcal{X}^\sigma_{T,\phi}$, and denote $g_1=\mathcal{S}f_1$, $g_2=\mathcal{S}f_2$. We have the equation
			\begin{align*}
				&\partial_t \partial_i (g_1-g_2)(x)-\int_{\mathbb{R}^d}\A(\nabla \phi(x),\alpha):{\delta_\alpha \nabla ^2\partial_i}(g_1-g_2)(x)\frac{d\alpha}{|\alpha|^{d+1}}\\
				&\ \ \ \ \quad\quad\quad\quad\quad\quad\quad\quad =\partial_i (\H[f_1,\phi]-\H[f_2,\phi])(x)+R_i[g_1,\phi](x)-R_i[g_2,\phi](x),\\
				&\partial_i (g_1-g_2)|_{t=0}=0.
			\end{align*}
			Applying Lemma \ref{lemlog} to obtain 
			\begin{align}
				\|g_1-g_2\|_T\lesssim &\sup_{t\in[0,T]}|\log t|^\varkappa  \left(t^\frac{\kappa}{3}\|( \H[f_1]-\H[f_2])(t)\|_{\dot C ^{\kappa-2}}+t^\frac{m+\kappa}{3}\|( \H[f_1]-\H[f_2])(t)\|_{\dot C ^{m
						+\kappa-2}}\right)\nonumber\\
				&+\sup_{t\in[0,T]}\left(\int_0^t{|\log (t-\tau)|^\varkappa\|(R[f_1]-R[f_2])(t)\|_{L^\infty}}d\tau+|\log t|^\varkappa  t^{\frac{m}{3}+1}\|(R[f_1]-R[f_2])(t)\|_{\dot C^m}\right).
			\end{align}
			By Lemma \ref{lemH1} and Lemma \ref{lemR},
			\begin{equation}\label{H1con}
				\begin{aligned}
					&\sup_{t\in[0,T]}|\log t|^\varkappa  \left(t^\frac{\kappa}{3}|
					|| (\H_\b^{ex}[f_1]-\H_\b^{ex}[f_2])(t)|||_{ {\kappa-2}}+t^\frac{m+\kappa}{3}|||\nabla^m (\H_\b^{ex}[f_1]-\H_\b^{ex}[f_2])(t)|||_{ {\kappa-2}}\right)\\
					&\quad+\sup_{t\in[0,T]} \left(\int_0^t|\log (t-\tau)|^\varkappa {\|(R[g_1]-R[g_2])(\tau)\|_{L^\infty}}d\tau+|\log t|^\varkappa  t^{\frac{m}{3}+1}\|(R[g_1]-R[g_2])(t)\|_{\dot C^m}\right)\\
					&\lesssim (\|f_1-f_2\|_{X_T}+\|g_1-g_2\|_{X_T})(\|f-\phi\|_{X_T}+\|g-\phi\|_{X_T}+T^\frac{1}{100}\|\nabla \phi\|_{C^{m+4}})\\
					&\ \ \ \ \ \ \times(1+\|f-\phi\|_{X_T}+\|g-\phi\|_{X_T}+T^\frac{1}{20}\|\nabla \phi\|_{C^{m+4}})^{m+4}.
				\end{aligned}
			\end{equation}
			And from Lemma \ref{lemh2} we obtain
			\begin{equation}\label{H2con}
				\begin{aligned}
					&\sup_{t\in[0,T]}|\log t|^\varkappa  \left(t^\frac{\kappa}{3}|
					|| (\H_\b^{im}[f_1]-\H_\b^{im}[f_2])(t)|||_{ {\kappa-2}}+t^\frac{m+\kappa}{3}||| \nabla^m(\H_\b^{im}[f_1]-\H_\b^{im}[f_2])(t)|||_{ {\kappa-2}}\right)\\&\lesssim ( |\log T|^{-\varkappa } +\|f_1-\phi\|_{X_T}+\|f_2-\phi\|_{X_T})\|f_1-f_2\|_{X_T},
				\end{aligned}
			\end{equation}
			Substituting \eqref{H1con} and \eqref{H2con} into the above inequality, we obtain 
			\begin{align*}
				\|g_1-g_2\|_T\leq C_2 (\|f_1-f_2\|_{X_T}+\|g_1-g_2\|_{X_T})(\sigma+T^\frac{1}{100}+|\log T|^{-\varkappa } ).
			\end{align*}
			Moreover, Lemma \ref{mapl2} implies that 
			\begin{align*}
				\|g_1-g_2\|_{L^\infty_TL^2}\leq C_2T^\frac{1}{50}\|f_1-f_2\|_{X_T}.
			\end{align*}
			We can further take $\sigma <\frac{1}{100(C_2+1)}$. Then 
			\begin{align*}
				\|g_1-g_2\|_{X_T}\leq \frac{1}{2}\|f_1-f_2\|_{X_T}.
			\end{align*}
			Hence $\mathcal{S}:\mathcal{X}^\sigma_{T,\phi}\to \mathcal{X}^\sigma_{T,\phi}$ is a contraction map. There exists a unique fixed point $\eta\in \mathcal{X}^\sigma_{T,\phi}$ satisfying $\mathcal{S}\eta=\eta$, which implies that $\eta$ is the solution to the Muskat equation \eqref{evoeta} with initial data $\eta_0$. This completes the proof.
		\end{proof}

		\section{Well-posedness in Besov space}\label{secbesov}
		In this section we introduce a general method to deal with quasi-linear parabolic systems in critical Besov spaces. This is a generalization of the results in section \ref{secsemi}, where we consider some semi-linear systems, and obtain the global/local well-posedness in critical Besov spaces in view of the fact that nonlinear remainder terms can be controlled by non-endpoint norms. In the following, we prove that this still holds true for quasi-linear systems as long as the variable coefficients are well-controlled. More precisely, the continuity condition \eqref{contlwp} for local well-posedness and the smallness condition \eqref{smallgwp} for global  well-posedness can be relaxed to a Besov version if the nonlinear remainder terms can be controlled by non-endpoint norms, and the variable coefficients are both bounded and non-degenerate.
		
		Consider a general quasi-linear evolution system: $U(t):\mathbb{R}^d\to \mathbb{R}^N$
		\begin{equation}\label{besovgen}
			\partial_t U(t,x)+A(\nabla U(t,x),\nabla)U(t,x)=\mathcal{N}(U)(t,x),
		\end{equation}
		where $A(Z, \xi)\sim_{Z}|\xi|^\sigma$, $\sigma>0$ for any $Z\in \mathbb{R}^{d\times N}$ and suppose that the Lipschitz space is a critical space of this problem and 
		\begin{equation}\label{besovno}
			\begin{aligned}
				&   \sup_{t\in[0,T]} (t^\frac{\kappa}{\sigma}\|\mathcal{N}(U)\|_{\dot C^{\kappa-\sigma+1}}+t^\frac{m+\kappa}{\sigma}\|\mathcal{N}(U)\|_{\dot C^{m+\kappa-\sigma+1}})\\
				& \quad\quad\quad \lesssim (1+||\nabla U||_{L^\infty_TL^\infty_x})^{c_1} ||\nabla U||_{T,*}^{1+\epsilon_0}(1+||\nabla U||_{L^\infty_TL^\infty_x}+||\nabla U||_{T,*})^{c_2m},
			\end{aligned}
		\end{equation}
		where $  \|f\|_{T,*}$ is the non-endpoint critical norm defined by
		\begin{align*}
			\|f\|_{T,*}=\sup_{t\in[0,T]}\left(t^{\frac{\mu}{\sigma}}\|f(t)\|_{\dot C^\mu}+t^{\frac{m+\kappa}{\sigma}}\|f(t)\|_{\dot C^{m+\kappa}}\right),
		\end{align*}
		with $0<\mu<\kappa<\sigma$ and $m\in \mathbb{N}, m\geq 10d$.
		For any fixed $x_0\in \mathbb{R}^d$, we can estimate the fundamental solution $K_{x_0}(t,\tau,x)$, which satisfies
		\begin{align*}
			&-\partial_\tau K_{x_0}(t,\tau,x)+A(\nabla U(\tau,x_0),\nabla)K_{x_0}(t,\tau,x)=0,\ \ \  (t,x) \in (\tau,+\infty)\times \mathbb{R}^d,\\
			&\lim_{\tau \to t^-}K_{x_0}(t,\tau,x)= \delta_{x=0}.
		\end{align*}
		Then, with some suitable assumptions on $A(\nabla U(t,x_0),\xi)$, and using argument in  Remark \ref{lemfourierK} to get something like the following estimates:
		\begin{equation}\label{esKbes}
			\begin{aligned}
				&\int_{\mathbb{R}^d}|\nabla_x^{n} K_{x_0}(t,\tau,x)|dx\lesssim (t-\tau)^{-\frac{n}{\sigma}}(1+\|\nabla U\|_{L^\infty([\tau,t]\times\mathbb{R}^d)})^{c_3}, 
			\end{aligned}
		\end{equation}
		for any $n\in\mathbb{N}$. 
		To control $\|\nabla U\|_{L^\infty_TL_x^\infty}$, we need to assume that $c_3<1. $ We would like to mention that all the models in this paper satisfy the above estimates with $c_3=0$, but essentially we only need $c_3<1$.\\
		By the fundamental solution, we can transform \eqref{besovgen} to the following integral equation of $U$:
		\begin{align*}
			U(t,x)
			=&\int_{\mathbb{R}^d} K_{x_0}(t,0,x-y)U_0(y)dy+\int_0^t\int_{\mathbb{R}^d} K_{x_0}(t,\tau,x-y)\mathcal{N}(U)(\tau,y)dyd\tau\\
			&+\int_0^t\int_{\mathbb{R}^d} K_{x_0}(t,\tau,x-y)(A(\nabla U(\tau,x_0),\nabla)-A(\nabla U(\tau,y),\nabla))U(\tau,y)dyd\tau\\
			:=&U_{L,x_0}(t,x)+U_{N,x_0}(t,x)+U_{R,x_0}(t,x).
		\end{align*}
		The above holds for any fixed $x_0\in\mathbb{R}^d$.
		In particular, 
		\begin{align}
			& \nabla U(t,x)=(\nabla_xU_{L,x_0}(t,x))\vert_{x_0=x}+(\nabla_xU_{N,x_0}(t,x))\vert_{x_0=x}+(\nabla_xU_{R,x_0}(t,x))\vert_{x_0=x}\label{hahaa1}.
		\end{align}
		By \eqref{besovno} and \eqref{esKbes}, we get 
		\begin{equation}\label{besovlinf}
			\begin{aligned}
				\|\nabla U\|_{L^\infty}&\leq \|(\nabla_xU_{L,x_0}(t,x))\vert_{x_0=x}\|_{L^\infty}+
				\|(\nabla_xU_{N,x_0}(t,x))\vert_{x_0=x}\|_{L^\infty}+\|(\nabla_xU_{R,x_0}(t,x))\vert_{x_0=x}\|_{L^\infty}\\
				&\lesssim 
				(1+\|\nabla U\|_{L^\infty_TL^\infty_x})^{c_3} \|\nabla U_0\|_{L^\infty}\\
				&\ \ \ +(1+||\nabla U||_{L^\infty_TL^\infty_x})^{c_1+c_4} ||\nabla U||_{T,*}^{1+\epsilon_0}(1+||\nabla U||_{L^\infty_TL^\infty_x}+||\nabla U||_{T,*})^{c_2m}.
			\end{aligned}
		\end{equation}
		Then we estimate the non-endpoint norm $\|\nabla U\|_{T,*}$. Denote $U_0^\eps=U_0\ast\rho_\eps$. Since 
		\begin{align*}
			(\nabla_xU_{L,x_0}(t,x))-\nabla U_0^\eps(x)=\int_{\mathbb{R}^d} (K_{x_0}(t,0,x-y))(\nabla U_0(y)-\nabla U_0^\eps(x))dy:=W_{L,x_0}(t,x), 
		\end{align*}
		we can write 
		\begin{align*}
			\nabla U(t,x)-\nabla U_0^\eps(x)
			&=W_{L,x_0}(t,x)+\nabla U_{N,x_0}(t,x)+\nabla U_{R,x_0}(t,x).
		\end{align*}
		For any $\alpha\in \mathbb{R}^d\backslash\{0\}$,
		\begin{align*}
			\delta_\alpha^x\nabla (U(t,x)-U_0^\eps(x))= (\delta_\alpha^x W_{L,x_0})(t,x)\vert_{x_0=x}+(\delta_\alpha^x\nabla U_{N,x_0})(t,x)\vert_{x_0=x}+(\delta_\alpha^x\nabla U_{R,x_0})(t,x)\vert_{x_0=x}.
		\end{align*}
		Then,  by Proposition \ref{normequ}, using \eqref{besovno} and arguments in the proof of Theorem \ref{lemmain}, we can obtain 
		\begin{equation}\label{bes11}
			\begin{aligned}
				&||\nabla (U- U_0^\eps)||_{T,*}\\
				&\lesssim 
				(1+||\nabla U||_{L^\infty_TL^\infty_x})^{c} \|U_0-U_0^\eps\|_{\dot B^1_{\infty,\infty}}+||\nabla U||_{T,*}^{1+\epsilon_0}(1+||\nabla U||_{L^\infty_TL^\infty_x}+||\nabla U||_{T,*})^{c_2m+c},
			\end{aligned}
		\end{equation}
		with $c=c_1+c_2+c_3+c_4$.
		Now we denote
		\begin{align*}
			B(T):=\|\nabla (U-U_0^\eps)\|_{T,*},\ \ \ \ M(T):=\|\nabla U\|_{L^\infty_TL^\infty}.
		\end{align*}
		Note that 
		\begin{align*}
			\|\nabla U\|_{T,*}\leq \|\nabla (U-U_0^\eps)\|_{T,*}+\|\nabla U_0^\eps\|_{T,*}\lesssim B(T)+T^\frac{\mu}{\sigma}\eps^{-(m+\kappa)}\|\nabla U_0\|_{L^\infty}.
		\end{align*}
		Hence, we write \eqref{besovlinf} and \eqref{bes11} as 
		\begin{align*}
			&M(T)\lesssim M(T)^{c_3} \|\nabla U_0\|_{L^\infty}+B(T)^{1+\eps_0}(1+M(T)+B(T)+T^\frac{\mu}{\sigma}\eps^{-(m+\kappa)}\|\nabla U_0\|_{L^\infty})^{c+c_2m+2},\\
			&B(T)\lesssim M(T)^{c_3}\|U_0-U_0^\eps\|_{\dot B^1_{\infty,\infty}}+B(T)^{1+\eps_0}(1+M(T)+B(T)+T^\frac{\mu}{\sigma}\eps^{-(m+\kappa)}\|\nabla U_0\|_{L^\infty})^{c+c_2m+2}.
		\end{align*}
		By Young's inequality, one has
		\begin{align*}
			M(T)^{c_3} \|\nabla U_0\|_{L^\infty}\lesssim \varepsilon M(T)+ C(\varepsilon)\|\nabla U_0\|_{L^\infty}^\frac{1}{1-c_3},
		\end{align*}
		provided $c_3\in(0,1)$. It is easy to check that under the assumption $(1+||\nabla U_0||_{L^\infty})^\frac{c_3}{1-c_3}\|\nabla(U_0-U_0^\eps)\|_{\dot B_{\infty,\infty}^0}\ll 1$, one can take
		$T$ small enough such that
		\begin{align*}
			&M(T)\leq C \|\nabla U_0\|_{L^\infty}^\frac{1}{1-c_3},\\
			&B(T)\leq C\|U_0-U_0^\eps\|_{\dot B^1_{\infty,\infty}}(1+\|\nabla U_0\|_{L^\infty})^\frac{c_3}{1-c_3}.
		\end{align*}
		This implies the local well-posedness under the condition 
		\begin{align*}
			\liminf_{\eps\to 0} (1+\|\nabla U_0\|_{L^\infty})^\frac{c_3}{1-c_3}\|\nabla(U_0-U_0^\eps)\|_{\dot B_{\infty,\infty}^0}\ll 1.
		\end{align*}
		\begin{proof}[Proof of Theorem \ref{thmPesB}]
			Note that $X_0\in \dot W^{1,\infty}\cap (C^2)^{\dot B^1_{\infty,\infty}}$, hence for any $0<\varepsilon_0\ll 1$, there exists $\eps_1>0$ such that 
			\begin{align*}
				\|X_0-X_{0}\ast \rho_{\eps_1}\|_{\dot B^1_{\infty,\infty}}\leq \varepsilon_0.
			\end{align*}
			Define
			\begin{align*}
				&\|f\|_{T,*}:=\sup_{t\in[0,T]}(t^\frac{1}{10}\|f(t)\|_{\dot C^\frac{1}{10}}+t^{m+\kappa }\|f(t)\|_{\dot C^{m+\kappa }}),\\
				&\|f\|_{T}:=\sup_{t\in[0,T]}(\|f(t)\|_{L^\infty}+t^{m+\kappa }\|f(t)\|_{\dot C^{m+\kappa }}).
			\end{align*}
			Here $m\in\mathbb{N}$ is a fixed integer. Denote $\tilde X=X-X\ast \rho_{\eps_1}$. 
			Applying Theorem \ref{lemmain} to \eqref{peseqr} again with $F=\mathcal{N}(X)+F(\Phi,X)$ and $G=R(\Phi,X)+A(\Phi)\Lambda\Phi$ to obtain that there exists $T_0=T_0(\|\partial_x X_0\|_{L^\infty},\eps_1)>0$ such that  for any $0<T\leq T_0$,
			\begin{equation}
				\begin{aligned}\label{pesma1}
					\|\tilde X'\|_{T,*}\lesssim &\|\tilde X'_0\|_{\dot B^0_{\infty,\infty}}+ \sup_{t\in[0,T]}\sum_{j\in\{0,m\}}\left(t^{j+\kappa}\left(\|\mathcal{N}(X)(t)\|_{\dot C^{j+\kappa}}+\|F(\Phi,X)(t)\|_{\dot C^{j+\kappa-a}}\right)\right)\\
					&+\|R(\Phi,X)\|_{L_T^1L^\infty}+\|A(\Phi)\Lambda\Phi\|_{L_T^1L^\infty}+\sup_{t\in[0,T]}t^{\frac{m}{s}+1}(\|R(\Phi,X)(t)\|_{\dot C^m}+\|A(\Phi)\Lambda\Phi\|_{L^\infty})\\
					&:=\|\tilde X'_0\|_{\dot B^0_{\infty,\infty}}+\da(\mathcal{N}(X)+F(\Phi,X),R(\Phi,X)+A(\Phi)\Lambda\Phi).
				\end{aligned}
			\end{equation}
			We get global solution if $\da(\mathcal{N}(X)+F(\Phi,X),R(\Phi,X)+A(\Phi)\Lambda\Phi)\leq (1+\tilde X'\|_{T,*})^{\eps}$ with $\eps\in(0,1)$. By Lemma \ref{lemnonpes},
			\begin{align*}
				&\sup_{t\in[0,T]}(t^\kappa\|\mathcal{N}(X)(t)\|_{\dot C^{\kappa}}+ t^{m+\kappa }	\|\mathcal{N}(X)(t)\|_{\dot C^{m+\kappa }})\lesssim C_{\mathcal{T},m}(1+\mathbf{\Theta}_X(T))^{m+5}\|\partial_x X\|_{T,*}^2 (1+\|\partial_x X\|_T)^{m+5}.
			\end{align*}
			Note that 
			\begin{align*}
				\|\partial_x X\|_{T,*}\lesssim \| \tilde X'\|_{T,*}+T^\frac{1}{10}\|\Phi\|_{C^{m+1}}\lesssim \| \tilde X'\|_{T,*}+T^\frac{1}{10}\eps_1^{-m}\|\partial_x X_0\|_{L^\infty}.
			\end{align*}
			Hence 
			\begin{align}\label{esTs}
				\|\tilde X'\|_{T,*}\lesssim&\|\tilde X'_0\|_{\dot B^0_{\infty,\infty}}+C_{\mathcal{T},m}(1+\mathbf{\Theta}_X(T))^{m+5}( \| \tilde X'\|_{T,*}+T^\frac{1}{10}\eps_1^{-m}\|\partial_x X_0\|_{L^\infty})^2 (1+\|\partial_x X\|_T)^{m+5}.
			\end{align}
			To show the existence of solution, we use the standard compactness argument. Denote $X_{0,\vartheta}=X_0\ast \rho_\vartheta$, and denote $X_{\vartheta}$ to be the solution of \eqref{peskin} with initial data $X_{0,\vartheta}$. By condition \eqref{con111}, we can take $\vartheta$ small enough such that for any $\alpha,s\in \mathbb{S}$, $\alpha\neq 0$,
			\begin{align*}
				\frac{|\delta_\alpha X_{0,\vartheta}(x)|}{|\alpha|}\geq \frac{1}{2\mathbf{\Theta}_0}.
			\end{align*}
			This implies that 
			$$
			\mathbf{\Theta}(X_{0,\vartheta})\leq 2 \mathbf{\Theta}_0.
			$$ 
			Moreover, 
			Then, by Theorem \ref{thmPes}, there exist $T_0=T_0(\|\partial_x X_0\|_{L^\infty},\vartheta)>0$, and a unique solution $X_\vartheta$ to the Peskin system in $[0,T_0]$ with initial data $X_{0,\vartheta}$, and 
			\begin{align*}
				&\sup_{t\in[0,T_0]}	(\|\partial_x X_\vartheta(t)\|_{L^\infty}+t^{m+\kappa}\|\partial_x X_\vartheta(t)\|_{\dot C^{m+\kappa}})\leq C \|\partial_x X_0\|_{L^\infty}.\\
				&\mathbf{\Theta}_{X_\vartheta}(T_0)\leq 4\mathbf{\Theta}_0.
			\end{align*}
			Denote 
			\begin{align*}
				T_1=\sup\{T\in\mathbb{R}^+:\|\tilde X'_\vartheta\|_{T,*}\leq 10\varepsilon_0,\|\partial_x X_\vartheta\|_T\leq 10 \|\partial_x X_0\|_{L^\infty}, \mathbf{\Theta}_{X_\vartheta}(T)\leq 10\Theta_0\}.
			\end{align*}
			We claim that $T_1\geq \frac{\varepsilon_0}{(1+\eps_1^{-m}+10\|\partial_x X_0\|_{L^\infty})^{100(m+5)}}$. If this is not true, then $T_1< \frac{\varepsilon_0}{(1+\eps_1^{-m}+10\|\partial_x X_0\|_{L^\infty})^{100(m+5)}}$. By the standard bootstrap argument, it suffices to show that 
			\begin{align}\label{boo}
				\|\tilde X'_\vartheta\|_{T_1,*}\leq 5\varepsilon_0,\ \ \ \|\partial_x X_\vartheta\|_{T_1}\leq 5 \|\partial_x X_0\|_{L^\infty},\ \ \  \mathbf{\Theta}_{X_\vartheta}(T_1)\leq 5\Theta_0.
			\end{align}
			This contradicts the definition of $T_1$. Hence we obtain a lower bound of $T_1$ that is independent in $\vartheta$. \\
			Now we only need to get the estimate \eqref{boo}. First we prove $\|\tilde X'_\vartheta\|_{T_1,*}\leq 5\varepsilon_0$. By \eqref{esTs}, we have 
			\begin{equation}\label{mmm}
				\begin{aligned}
					&	\|\tilde X'_\vartheta\|_{T_1,*}\\
					&	\leq C_{\mathcal{T},m}\|\tilde X'_0\|_{\dot B^0_{\infty,\infty}}+\tilde C_1(1+10\mathbf{\Theta}_0)^{m+5}( \| \tilde X'_\vartheta\|_{T,*}+T_1^\frac{1}{10}\eps_1^{-m}\|\partial_x X_0\|_{L^\infty})^2 (1+10\|\partial_x X_0\|_{L^\infty})^{2(m+5)}\\
					&	\leq	2\varepsilon_0+100\tilde C (1+10\mathbf{\Theta}_0)^{m+5}\varepsilon_0^2\leq 5\varepsilon_0.
				\end{aligned}
			\end{equation}
			Here the last inequality follows by fixing $\varepsilon_0\leq \frac{1}{(100+10\mathbf{\Theta}_0+\tilde C)^{2(m+10)}}$. 
			Moreover, applying Theorem \ref{lemmain} again we have 
			\begin{equation}\label{tnor}
				\begin{aligned}
					\|\partial_x X\|_{T_1}&\leq \|\partial_x X_0\|_{L^\infty}+\tilde C_2(1+10\mathbf{\Theta}_0)^{m+5}( \| \tilde X'\|_{T_1,*}+T_1^\frac{1}{10}\eps_1^{-m}\|\partial_x X_0\|_{L^\infty})^2 (1+\|\partial_x X\|_{T_1})^{m+5}\\
					&:=(1+\Theta_1(T,\eps_0,\eps_1,\|\partial_xX_{\vartheta}\|_{L^\infty}))\|\partial_x X_0\|_{L^\infty}+\Theta_2(T,\eps_0,\eps_1,\|\partial_xX_{\vartheta}\|_{L^\infty}).
				\end{aligned}
			\end{equation} 
			By taking $T$ and $\eps_0$ small enough, we can let $\Theta_1(T,\eps_0,\eps_1,\|\partial_xX_{\vartheta}\|_{L^\infty})\leq 1$ and $\Theta_2(T,\eps_0,\eps_1,\|\partial_xX_{\vartheta}\|_{L^\infty})\leq 2\|\partial_x X_0\|_{L^\infty}$, so we can prove $\|\partial_x X\|_{T_1}\leq 5\|\partial_x X_0\|_{L^\infty}$.\\ 
			Finally, we need to control $\mathbf{\Theta}_{X_\vartheta}(T_1)$. When considering $\dot B^1_{\infty,\infty}$ data, it is different from the estimate \eqref{wes} when considering $C^1$ data. We follow the idea in \cite{KN}, and introduce a quantity
			\begin{align*}
				Q_h(T)=\sup_{t\in[0,T]}\sup_{\alpha,s}\frac{|\alpha|^\varepsilon}{t^\varepsilon}\left|\frac{1}{|\Delta_\alpha h(t,x)|}-\frac{1}{|\Delta_\alpha h(0,x)|}\right|.
			\end{align*}
			By \eqref{eqpes}, we can write
			\begin{align*}
				\frac{1}{|\Delta_\alpha X(t,\cdot)(x)|}-\frac{1}{|\Delta_\alpha X_{0}(x)|}&=\int_0^t \partial_t\left(\frac{1}{|\Delta_\alpha X(\tau,\cdot)(x)|}\right)d\tau\\
				&\leq \mathbf{\Theta}_X(t)^2 \int_0^t |\Delta_\alpha \mathcal{H} (\mathbf{T}(|\partial_x X|)\partial_xX)(\tau,\cdot)(x)|+|\Delta_\alpha N(X(\tau,\cdot))(x)|d\tau.
			\end{align*}
			The integral in the right hand side can be controlled by $\|\partial_x X\|_{T,*}$.
			It can be checked that 
			\begin{align*}
				Q_{X_\vartheta}(T_1)&\leq \mathbf{\Theta}_{X_\vartheta}(T_1)^2 \sup_{t\in[0,T_1]} t^{1-\varepsilon}(\|\mathbf{T}(|\partial_x X_\vartheta|)\partial_xX_\vartheta(t)\|_{\dot C^{1-\varepsilon}}+\|N(X_\vartheta(t))\|_{\dot C^{1-\varepsilon}})\\
				&\lesssim (1+\mathbf{\Theta}_{X_\vartheta}(T_1))^{10} (1+\|\partial_x X_\vartheta\|_{T_1})^{10}\|\partial_x X_\vartheta\|_{T_1,*}\leq \tilde C_3(1+\mathbf{\Theta}_0)^{10} (1+\|\partial_x X_0\|_{L^\infty})^{10}\varepsilon_0.
			\end{align*}
			Then we can take $\varepsilon_0$ small enough, and apply \cite[Lemma 2.8]{KN} to obtain 
			\begin{align*}
				\mathbf{\Theta}_{X_\vartheta}(T_1)\leq 2\mathbf{\Theta}_0.
			\end{align*}
			Combining this with \eqref{mmm} and \eqref{tnor}, we obtain \eqref{boo}.
			Hence $T_1$ has a lower bound that is independent in $\vartheta$. Then we can pass to the limit $\vartheta\to 0$, the sequence $\{X_\vartheta\}_\vartheta$ will converge to a solution $X$, which is a solution to \eqref{eqpes} with initial data $X_0$. And 
			$$
			\|\tilde X'\|_{T,*}\leq 10\varepsilon_0,\ \  \|\partial_x X\|_T\leq 10 \|\partial_x X_0\|_{L^\infty},\ \  \mathbf{\Theta}_{X}(T)\leq 10\Theta_0.
			$$
			This completes the proof.

		\end{proof}
		\section{Well-posedness with boundary condition}\label{secbdd}
		In this section, we will prove  well-posedness results of the mean curvature system with Dirichlet boundary condition, which can be written as 
		\begin{equation}\label{eqmcbdd}
			\begin{aligned}
				&\partial_t f=\A[\nabla f]:\nabla^2 f,\ \ \ \text{in}\ [0,T]\times \Omega,\\
				&f|_{t=0}=f_0,\ \ \ \text{in}\ \Omega,\\
				&f=0,\ \ \text{on}\ [0,T]\times\partial \Omega,
			\end{aligned}
		\end{equation}
		where the coefficient $\A$ is defined in \eqref{defamcf}. We will consider two cases: (i) $\Omega=\mathbb{R}_+^d$, (ii) $\Omega\subset \mathbb{R}^d$ is a bounded domain with $\partial\Omega\in C^{2m+3}$.
		We fix two constants: $m\in\mathbb{N}^+$, $m\geq 10d$, $\kappa\in(\frac{1}{2},1)$, $1-\kappa\ll1$.	
		
		We remark that for general quasi-linear evolution equation on bounded domain, higher order derivative may blow up at boundary.  However, the special structure of the nonlinear term in the mean curvature system \eqref{eqmcbdd} allows us to improve the regularity. The main ingredients are the analysis of the behaviour of solution on the boundary. More precisely, when considering $\Omega=\mathbb{R}^d_+$, we proved Lemma \ref{lembdc}, which indicates that any solution $f\in L^\infty(0,T;C^{2m+2}(\bar {\mathbb{R}}_+^d))$ satisfies $\Delta^l f|_{(0,T)\times \partial\mathbb{R}_+^d}=0$ for any $0\leq l\leq m+1$. When considering general bounded domain, the situation is more subtle, we proved in Lemma \ref{mcbdlembdy} that normal derivatives of the solution can be reduced to lower order ones on the boundary. These good properties essentially help us to obtain the estimates of higher order derivatives.
		\subsection{Well-posedness of the Mean curvature flow on half space} We first consider the case $\Omega=\mathbb{R}_+^d$. For $f:[0,T]\times\mathbb{R}^d_+\to \mathbb{R}^N$, define 
		\begin{align*}
			&\|f\|_{X_T}:=\|\nabla f\|_{L_T^\infty L^\infty(\mathbb{R}^d_+)}+\sup_{t\in[0,T]}t^{\frac{1+\kappa}{2}}\|\left(\nabla^2f,\partial_tf\right)(t)\|_{\dot C^{\kappa}(\mathbb{R}^d_+)},\\
			&\|f\|_{Y_T}:=\|\nabla f\|_{L_T^\infty L^\infty(\mathbb{R}^d_+)}+\sum_{n\leq m}\sup_{t\in[0,T]}t^{n+\frac{\kappa+1}{2}}\|\nabla^{2n}\left(\nabla^2f,\partial_tf\right)(t)\|_{\dot C^{\kappa}(\mathbb{R}^d_+)}.
		\end{align*}
		Let $\phi\in C^{\infty}(\bar {\mathbb{R}}^d_+)$ satisfy  $\phi|_{\partial\mathbb{R}^d_+}=0$, and $\sigma\in(0,1)$. We define the following sets 
		\begin{equation*}
			\begin{aligned}
				\mathcal{X}_{T,\phi}^\sigma=\left\{f\in L^\infty_T C^2(\bar{\mathbb{R}}^d_+),f|_{[0,T]\times\partial\mathbb{R}^d_+}=0:\|f-\phi\|_{X_T}\leq\sigma\right\},
			\end{aligned}
		\end{equation*}
		\begin{equation*}
			\begin{aligned}
				&\mathcal{Y}_{T,\phi}^\sigma=\left\{f\in L^\infty_T C^{2m+2}(\bar{\mathbb{R}}_+^d), f|_{[0,T]\times \partial\mathbb{R}_+^d}=0, \Delta^k f|_{(0,T]\times\partial\mathbb{R}_+^d}=0,k=1,\cdots,m:\|f-\phi\|_{Y_T}\leq\sigma\right\}.
			\end{aligned}
		\end{equation*}
		We will prove Theorem \ref{eqmcls} first. To prove this, we will need the following lemma.
		\begin{lemma}\label{lemheat}
			Let $\B\in C^\infty(\mathbb{R}^d_+;\mathbb{R}^{d\times d}_{sym})$ be uniformly elliptic, which means there exists $c_0$ such that 
			\begin{equation*}
				c_0^{-1}|\xi|^2\leq \xi \B \xi^\top  \leq c_0|\xi|^2,\ \ \ \forall\ \xi\in\mathbb{R}^d.
			\end{equation*}
			Let $T>0$, and $u$ be a solution to the system 
			\begin{equation}\label{mcbd}
				\begin{aligned}
					&\partial_tu-\B:\nabla^2 u=F^1+F^2,\quad \text{in}\ [0,T]\times\mathbb{R}^d_+,\\
					&u|_{t=0}=u_0,\quad \text{in}\ \mathbb{R}^d_+,\\
					&u=0,\quad \text{on}\ [0,T]\times\partial\mathbb{R}^d_+,
				\end{aligned}
			\end{equation}
			with $F^1,F^2\in C^\infty([0,T]\times \bar {\mathbb{R}}^d_+)$, $F^2|_{\partial \mathbb{R}^d_+}=0$, and $u_0\in \dot W^{1,\infty}(\mathbb{R}^d_+)$, $u_0|_{\partial\mathbb{R}^d_+}=0$,
			then  we have
			\begin{align*}
				&  \sup_{t\in[0,T]}\left(\|\nabla u(t)\|_{L^\infty(\mathbb{R}^d_+)}+t^\frac{1+\kappa}{2}\|\nabla^2 u(t)\|_{\dot C^{\kappa}(\mathbb{R}^d_+)}\right)\\
				&\quad \lesssim \|u_0\|_{\dot W^{1,\infty}(\mathbb{R}^d_+)}+\sup_{t\in[0,T]}\left(t^{\frac{1}{2}}\|F^1(t)\|_{L^\infty(\mathbb{R}^d_+)}+t^\frac{1+\kappa}{2}\|F^1(t)\|_{\dot C^\kappa(\mathbb{R}^d_+)}+t^\frac{1+\kappa}{2}\|F^2(t)\|_{\dot C^\kappa(\mathbb{R}^d_+)}\right)\\
				&\quad+\sup_{0<\tau<t<T}\tau^{\frac{1+\kappa}{2}}\frac{\|F^1(t)-F^1(\tau)\|_{L^\infty(\mathbb{R}^d_+)}}{(t-\tau)^\frac{\kappa}{2}}+T^{\frac{1}{2}}\|\B\|_{\dot W^{1,\infty}}\sup_{t\in[0,T]}t^{\frac{1}{2}}\|\nabla^2u(t)\|_{L^\infty(\mathbb{R}^d_+)}.
			\end{align*}
		\end{lemma}
		\begin{proof}
			We freeze coefficient at $x_0\in \mathbb{R}^d_+$, and denote $\B_{0}=\B(x_0)$, then 
			\begin{align*}
				\partial_tu(t,x)-\B_{0}:\nabla^2 u(t,x)=F^1(t,x)+F^2(t,x)+\mathsf{R}_{x_0}(t,x),
			\end{align*}
			where
			\begin{equation*}
				\mathsf{R}_{x_0}(x)=(\B(x)-\B_{0}):\nabla^2 u(t,x).
			\end{equation*}
			There exists a linear  transform $\psi_{x_0}:\mathbb{R}^d_+\to \mathbb{R}^d_+$ such that 
			$(B_{0}:\nabla^2u)\circ\psi_{x_0}=-\Delta(u\circ\psi_{x_0})$. Since $\B$ is uniformly elliptic, one has $|\nabla\psi_{x_0}|+|\nabla\psi_{x_0}^{-1}|\leq C(c_0)$. Let $v(t,x)=u(t,\cdot)\circ\psi_{x_0}$, we obtain the heat equation on half space:
			\begin{equation}\label{mcbdheat}
				\begin{aligned}
					&\partial_t v-\Delta v=\tilde F^1_{x_0}+\tilde F^2_{x_0}+\tilde {\mathsf{R}}_{x_0}, \quad \text{in}\ [0,T]\times\mathbb{R}^d_+,\\
					&v|_{t=0}=v_{x_0}=u_0\circ\psi_{x_0},\quad \text{in}\ \mathbb{R}^d_+,\\
					&v(t,x)=0,\quad \text{on}\ [0,T]\times\partial\mathbb{R}^d_+,
				\end{aligned}
			\end{equation}
			with
			$$
			\tilde F^1_{x_0}(t,x)=(F^1(t)\circ\psi_{x_0})(x), \ \ \ \tilde F^2_{x_0}(t,x)=(F^2(t)\circ\psi_{x_0})(x),\ \ \ \tilde {\mathsf{R}}_{x_0}(t,x)=(\mathsf{R}_{x_0}(t)\circ\psi_{x_0})(x).
			$$
			We have  the following fundamental solution to the heat equation on half space with Dirichlet boundary condition:
			\begin{align}\label{defH}
				\H(t,x'-y',x_d,y_d)= \K(t,x'-y',x_d-y_d)-\K(t,x'-y',x_d+y_d),
			\end{align}
			where $\K(t,x)$ is the standard heat kernel on $\mathbb{R}^d$, and we denote $x=(x',x_d)\in\mathbb{R}^d$ with $x'\in\mathbb{R}^{d-1}$. Similarly we denote
			\begin{align}\label{deftilH}
				\tilde{\H}(t,x'-y',x_d,y_d)= \K(t,x'-y',x_d-y_d)+\K(t,x'-y',x_d+y_d).
			\end{align}
			Remark that 
			\begin{align}\label{heates}
				\|\nabla ^n \K(t,\cdot) \|_{L^1}\lesssim t^{-\frac{n}{2}},\ \ \ \text{which implies that}\ \   \|\delta_\alpha\nabla ^n \K(t,\cdot) \|_{L^1}\lesssim t^{-\frac{n}{2}}\min\{1,|\alpha|t^{-\frac{1}{2}}\},\ \ \forall \alpha\in\mathbb{R}^d.
			\end{align}
			From this we obtain, for any $x\in\mathbb{R}_+^d$, $n\in\mathbb{N}$,
			\begin{equation}\label{eskerha}
				\begin{aligned}
					&\int_{\mathbb{R}^d_+}|\nabla_x ^n \H(t,x'-y',x_d,y_d)|dy\lesssim t^{-\frac{n}{2}},\\  &\int_{\mathbb{R}^d_+}|\delta_\alpha\nabla_x ^n \H(t,x'-y',x_d,y_d)||x-y|^{\beta}dy\lesssim t^{-\frac{n}{2}}|\alpha|^{\beta}\min\{1,|\alpha|t^{-\frac{1}{2}}\}^{1-\beta},\ \ \forall \alpha\in\mathbb{R}^d, x-\alpha\in\mathbb{R}^d_+,\beta\in[0,1).
				\end{aligned}
			\end{equation}
			In fact, to get \eqref{eskerha}, we use a change of variable to obtain 
			\begin{align*}
				&\int_{\mathbb{R}^d_+}|\delta_\alpha\nabla_x ^n \H(t,x'-y',x_d,y_d)||x-y|^{\beta}dy\\
				&\quad\lesssim\int_{\mathbb{R}^d_+}|\delta_\alpha\nabla_x ^n \K(t,x-y)||x-y|^{\beta}dy+\int_{\mathbb{R}^d_-}|\delta_\alpha\nabla_x ^n \K(t,x-y)||(x'-y',x_d+y_d)|^{\beta}dy.
			\end{align*}
			By \eqref{leme1}, the first part can be estimated. For the second part, since $x\in\mathbb{R}_+^d$, and $y\in\mathbb{R}_-^d$, obviously $|(x'-y',x_d+y_d)|^{\beta}\leq |x-y|^{\beta}$, so also apply \eqref{leme1} to get the estimate. Similar results also holds for $\tilde{\H}$.
			\\
			Using the fundamental solution $\H$, we can	transform \eqref{mcbdheat} to the following form.
			\begin{equation}\label{mcbdfor}
				\begin{aligned}
					v(t,x)&=\int_{\mathbb{R}_+^d} \H(t,x'-y',x_d,y_d)v_{x_0}(y)dy+\sum_{k=1}^2\int_0^t\int_{\mathbb{R}^d_+} \H(t-\tau,x'-y',x_d,y_d)\tilde F^k_{x_0}(\tau,y)dyd\tau\\
					&\quad\quad+\int_0^t\int_{\mathbb{R}^d_+} \H(t-\tau,x'-y',x_d,y_d)\tilde {\mathsf{R}}_{x_0}(\tau,y)dyd\tau\\
					&:=v_{L,x_0}(t,x)+\sum_{k=1,2}v_{N,x_0}^k(t,x)+v_{R,x_0}(t,x).
				\end{aligned}
			\end{equation}
			Moreover,   \begin{align*}
				&\nabla v(t,x)=\nabla v_{L,x_0}(t,x)+\sum_{k=1,2}\nabla v_{N,x_0}^k(t,x)+\nabla v_{R,x_0}(t,x),\\ &\delta_\alpha\nabla^2v(t,x)=\delta_\alpha\nabla^2v_{L,x_0}(t,x)+\sum_{k=1,2}\delta_\alpha\nabla^2v_{N,x_0}^k(t,x)+\delta_\alpha\nabla^2v_{R,x_0}(t,x),\ \  \ \forall \alpha\in\mathbb{R}^d,\  x-\alpha\in\mathbb{R}^d_+,
			\end{align*}
			hold for any fixed $x_0\in\mathbb{R}^d_+$. Specifically, we can take $x_0=x$, then 
			\begin{equation}\label{spliv}
				\begin{aligned}
					& \|\nabla v(t)\|_{L^\infty}\leq \|(\nabla v_{L,x_0}(t))|_{\substack{x_0=x}}\|_{L^\infty}+\sum_{k=1,2}\|(\nabla v_{N,x_0}^k(t))|_{\substack{x_0=x}}\|_{L^\infty}+\|(\nabla v_{R,x_0}(t))|_{\substack{x_0=x}}\|_{L^\infty},\\
					&\|\delta_\alpha\nabla^2v(t)\|_{L^\infty}\leq\|(\delta_\alpha\nabla^2v_{L,x_0}(t))|_{\substack{x_0=x}}\|_{L^\infty}+\sum_{k=1,2}\|(\delta_\alpha\nabla^2v_{N,x_0}^k(t))|_{\substack{x_0=x}}\|_{L^\infty}\\
					&\quad\quad\quad\quad\quad\quad+\|(\delta_\alpha\nabla^2v_{R,x_0}(t))|_{x_0=x}\|_{L^\infty}.
				\end{aligned}
			\end{equation}
			We first estimate the linear part $v_{L,x_0}$. Using integration by parts, since $v_{x_0}|_{\partial\mathbb{R}_+^d}=0$, one has 
			\begin{align*}
				\nabla v_{L,x_0}(t,x)&=\int_{\mathbb{R}^d_+}\nabla_x \H(t,x'-y',x_d,y_d)v_{x_0}(y)dy=\int_{\mathbb{R}^d_+} \H(t,x'-y',x_d,y_d)\nabla v_{x_0}(y)dy.
			\end{align*}
			Hence by \eqref{eskerha} we obtain 
			\begin{equation}\label{llha}
				\begin{aligned}
					&\|\nabla v_{L,x_0}(t)\|_{L^\infty(\mathbb{R}^d_+)}\lesssim \|(\nabla v_{x_0})\|_{L^\infty(\mathbb{R}^d_+)}\lesssim \|\nabla u_0\|_{L^\infty(\mathbb{R}^d_+)},
					\\
					&\| \delta_\alpha\nabla^2 v_{L,x_0}(t)\|_{L^\infty(\mathbb{R}^d_+)}\lesssim \int_{\mathbb{R}^d_+} |\delta_\alpha\nabla_x\H(t,x'-y',x_d,y_d)|dy\|\nabla v_{x_0}\|_{L^\infty(\mathbb{R}^d_+)}\\
					&\quad\quad\quad\quad\quad\quad\quad\quad\quad\quad\quad\quad\lesssim \min\left\{1,\frac{|\alpha|}{t^{\frac{1}{2}}}\right\}t^{-\frac{1}{2}}\|\nabla u_0\|_{L^\infty(\mathbb{R}^d_+)},
				\end{aligned}
			\end{equation}
			for any $x_0$. Then we get
			\begin{align}\label{lllha}
				\sup_{t\in[0,T]}\left(\|(\nabla v_{L,x_0}(t))|_{x_0=x}\|_{L^\infty(\mathbb{R}^d_+)}+t^\frac{1+\kappa}{2}\sup_{\alpha}\frac{\| (\delta_\alpha\nabla^2 v_{L,x_0}(t))|_{x_0=x}\|_{L^\infty(\mathbb{R}^d_+)}}{|\alpha|^\kappa}\right)\lesssim \|u_0\|_{\dot W^{1,\infty}}
			\end{align}
			Then we consider $v_{N,x_0}^1(t,x)$. We have 
			\begin{align*}
				\nabla v_{N,x_0}^1(t,x)=\int_0^t\int_{\mathbb{R}^d_+} \nabla_x \H(t-\tau,x'-y',x_d,y_d)\tilde F^1_{x_0}(\tau,y)dyd\tau.
			\end{align*}
			Applying \eqref{eskerha} again to obtain 
			\begin{align}\label{unlinft}
				\|\nabla v_{N,x_0}^1(t)\|_{L^\infty}\lesssim  \int_0^t (t-\tau)^{-\frac{1}{2}}\|\tilde F^1_{ x_0}(\tau)\|_{L^\infty} d\tau\lesssim \sup_{t\in[0,T]}t^{\frac{1}{2}}\|F^1(t)\|_{ L^\infty}.
			\end{align}
			For higher order derivatives, it suffices to consider $\nabla_{x'}\nabla v_{N,x_0}^1$ and $\partial_{x_d}^2  v_{N,x_0}^1$ respectively. Note that
			\begin{align*}
				\delta_\alpha\nabla_{x'}\nabla v_{N,x_0}^1(t,x)&=\int_0^t\int_{\mathbb{R}^d_+}\delta_\alpha\nabla_{x'}\nabla_x \H(t-\tau,x'-y',x_d,y_d)\tilde F^1_{x_0}(\tau,y)dyd\tau\\
				&=\int_0^t\int_{\mathbb{R}^d_+}\delta_\alpha\nabla_{y'}\nabla_x \H(t-\tau,x'-y',x_d,y_d)(\tilde F^1_{x_0}(\tau,y)-\tilde F^1_{ x_0}(\tau,x))dyd\tau.
			\end{align*}
			Hence, by \eqref{eskerha},
			\begin{equation}\label{dxpdx}
				\begin{aligned}
					&\|\delta_\alpha\nabla_{x'}\nabla  v_{N,x_0}^1(t) \|_{L^\infty}\lesssim \int _0^t \int_{\mathbb{R}^d_+} |\delta_\alpha \nabla_{y'}\nabla_x\H(t-\tau,x'-y',x_d,y_d)||x-y|^\kappa dy  \|\tilde F^1_{x_0}(\tau)\|_{\dot C^\kappa(\mathbb{R}^d_+)} d\tau\\
					&\quad\quad\quad\lesssim  |\alpha|^\kappa\int_0^t (t-\tau)^{-\frac{1}{2}}\tau^{-\frac{1+\kappa}{2}}d\tau \sup_{t\in[0,T]}t^\frac{1+\kappa}{2}\|F^1(t)\|_{\dot C^\kappa(\mathbb{R}^d_+)}\\
					&\quad\quad\quad\lesssim  |\alpha|^\kappa t^{-\frac{1+\kappa}{2}}\sup_{t\in[0,T]}t^\frac{1+\kappa}{2}\|F^1(t)\|_{\dot C^\kappa(\mathbb{R}^d_+)}.
				\end{aligned}
			\end{equation}
			Finally, for $\partial_{x_d}^2  v_{N,x_0}^1$, observe that 
			$$\partial_{x_d}^2\H(t,x'-y',x_d,y_d)=\partial_t\H(t,x'-y',x_d,y_d)-\Delta_{x'}\H(t,x'-y',x_d,y_d),$$
			then we can write 
			\begin{align*}
				\partial_{x_d}^2 v_{N,x_0}^1(t,x)=\int_0^t\int_{\mathbb{R}^d_+} \partial_t\H(t-\tau,x'-y',x_d,y_d)\tilde F^1_{x_0}(\tau,y)dyd\tau- \Delta_{x'}v_{N,x_0}^1(t,x).
			\end{align*}
			By \eqref{dxpdx}, we know that 
			\begin{align*}
				\sup_{t\in[0,T]}t^\frac{1+\kappa}{2} \sup_{\alpha} \frac{\|\delta_\alpha\Delta_{x'}v_{N,x_0}^1(t)\|_{L^\infty}}{|\alpha|^\kappa}\lesssim \sup_{t\in[0,T]}t^\frac{1+\kappa}{2}\|{F}^1(t)\|_{\dot   C^\kappa(\mathbb{R}^d_+)}.
			\end{align*}
			Hence it suffices to consider 
			\begin{align*}
				& \int_0^t\int_{\mathbb{R}^d_+} \partial_t\H(t-\tau,x'-y',x_d,y_d)\tilde F^1_{ x_0}(\tau,y)dyd\tau\\
				&\quad\quad\quad\quad\quad=  \int_0^t\int_{\mathbb{R}^d_+} \partial_t\H(t-\tau,x'-y',x_d,y_d)(\tilde F^1_{ x_0}(\tau,y)-\tilde F^1_{x_0}(t,y))dyd\tau\\
				&\quad\quad\quad\quad\quad\quad\quad+ \int_0^t\int_{\mathbb{R}^d_+} \partial_t\H(t-\tau,x'-y',x_d,y_d)\tilde F^1_{ x_0}(t,y)dyd\tau\\
				&\quad\quad\quad\quad\quad=  \int_0^t\int_{\mathbb{R}^d_+} \partial_t\H(t-\tau,x'-y',x_d,y_d)(\tilde F^1_{ x_0}(\tau,y)-\tilde F^1_{x_0}(t,y))dyd\tau\\
				&\quad\quad\quad\quad\quad\quad\quad+ \int_{\mathbb{R}^d_+} \H(t,x'-y',x_d,y_d)\tilde F^1_{ x_0}(t,y)dy-\tilde F^1_{ x_0}(t,x)\\
				&\quad\quad\quad\quad\quad= I^1_{ x_0}(t,x)+I^2_{ x_0}(t,x)-\tilde F^1_{x_0}(t,x).
			\end{align*}
			By \eqref{eskerha} and the boundness of $\psi_{x_0}$, it is easy to check that
			\begin{align*}
				&\left|\delta_\alpha I^1_{ x_0}(t,x) \right|\\
				&\quad\lesssim \int_0^t\int_{\mathbb{R}^d_+}| \delta_\alpha\partial_t\H(t-\tau,x'-y',x_d,y_d)|dy|t-\tau|^\frac{\kappa}{2}\tau^{-\frac{1+\kappa}{2}}d\tau \sup_{0<\tau<t<T}\tau^{\frac{1+\kappa}{2}}\frac{\|\tilde{F}_{x_0}^1(t)-\tilde{F}^1_{x_0}(\tau)\|_{L^\infty}}{(t-\tau)^\frac{\kappa}{2}}\\
				&\quad\lesssim |\alpha|^\kappa t^{-\frac{1+\kappa}{2}}\sup_{0<\tau<t<T}\tau^{\frac{1+\kappa}{2}}\frac{\|F^1(t)- F^1(\tau)\|_{L^\infty}}{(t-\tau)^\frac{\kappa}{2}},
			\end{align*}
			and
			\begin{align*}
				\left|\delta_\alpha I^2_{ x_0}(t,x)\right|\lesssim |\alpha|^\kappa t^{-\frac{\kappa}{2}}\|F^1(t)\|_{L^\infty}.
			\end{align*}
			Hence we conclude that 
			\begin{align*}
				&\sup_{t\in[0,T]}t^\frac{1+\kappa}{2}\sup_\alpha\frac{ \|\delta_\alpha \partial_{x_d}^2 v_{N,x_0}^1(t)\|_{L^\infty}}{|\alpha|^\kappa}\\
				&\quad\quad\lesssim \sup_{0<\tau<t<T}\tau^{\frac{1+\kappa}{2}}\frac{\|F^1(t)-F^1(\tau)\|_{L^\infty(\mathbb{R}^d_+)}}{(t-\tau)^\frac{\kappa}{2}}+\sup_{t\in[0,T]}\left(t^{\frac{1}{2}}\|F^1(t)\|_{L^\infty(\mathbb{R}^d_+)}+t^\frac{1+\kappa}{2}\|F^1(t)\|_{\dot C^\kappa(\mathbb{R}^d_+)}\right).
			\end{align*}
			Combining this with \eqref{unlinft} and \eqref{dxpdx}, and notice that those estimates holds for any $x_0$, especially for $x_0=x$, to obtain that 
			\begin{equation}\label{unnha}
				\begin{aligned}
					& \sup_{t\in[0,T]}\left(\|(\nabla v_{N,x_0}^1(t))|_{\substack{x_0=x}}\|_{L^\infty(\mathbb{R}^d_+)}+t^\frac{1+\kappa}{2}\sup_{\alpha}\frac{\| (\delta_\alpha\nabla^2 v_{N,x_0}^1(t))|_{\substack{x_0=x}}\|_{L^\infty(\mathbb{R}^d_+)}}{|\alpha|^\kappa}\right)\\
					&\quad\lesssim \sup_{0<\tau<t<T}\tau^{\frac{1+\kappa}{2}}\frac{\|F^1(t)-F^1(\tau)\|_{L^\infty(\mathbb{R}^d_+)}}{(t-\tau)^\frac{\kappa}{2}}+\sup_{t\in[0,T]}\left(t^{\frac{1}{2}}\|F^1(t)\|_{L^\infty(\mathbb{R}^d_+)}+t^\frac{1+\kappa}{2}\|F^1(t)\|_{\dot C^\kappa(\mathbb{R}^d_+)}\right).
				\end{aligned}
			\end{equation}
			Next we consider $v_{N,x_0}^2$. Note that by assumption, $\tilde F^2_{ x_0}$ vanishes on boundary, hence we can make an odd extension of $\tilde F^2_{x_0}$, then the estimate of $v_{N,x_0}^2$ is similar to the nonlinear term in the whole space in section \ref{secproof}. We write the proof shortly in the following. Note that 
			\begin{align*}
				v_{N,x_0}^2(t,x)&=\int_0^t \int_{\mathbb{R}^d_+}\H(t-\tau,x'-y',x_d,y_d)\tilde F^2_{ x_0}(\tau,y)dyd\tau\\
				&=\int_0^t \int_{\mathbb{R}^d}\K(t-\tau,x-y)\bar F^2_{x_0}(\tau,y)dyd\tau.
			\end{align*}
			Here $\K$ is the standard heat kernel on $\mathbb{R}^d$, and 
			\begin{align*}
				\bar F^2_{ x_0}(t,x)=\begin{cases}
					\tilde F^2_{ x_0}(t,x',x_d),\ \ \ x\in\mathbb{R}^d_+,\\
					-\tilde F^2_{ x_0}(t,x',-x_d),\ \ \ x\in \mathbb{R}^d_-,
				\end{cases}
			\end{align*}
			is the odd extension of $\tilde F^2_{ x_0}$. Note that 
			\begin{align*}
				\|\bar F^2_{x_0}\|_{\dot C^\kappa(\mathbb{R}^d)}\lesssim \| \tilde F^2_{x_0}\|_{\dot C^\kappa(\mathbb{R}^d_+)}\lesssim \|F^2\|_{\dot C^\kappa(\mathbb{R}^d_+)},\ \ \  \forall x_0\in\mathbb{R}^d_+.
			\end{align*}
			We have 
			\begin{equation}\label{unlha}
				\begin{aligned}
					\|\nabla v_{N, x_0}^2(t)\|_{L^\infty}&\lesssim \left\|\int_0^t \int_{\mathbb{R}^d} \nabla \K(t-\tau,x-y)(\bar F^2_{ x_0}(\tau,y)-\bar F^2_{ x_0}(\tau,x))dyd\tau \right\|_{L^\infty}\\
					&\lesssim \int_0^t \int_{\mathbb{R}^d} |\nabla \K(t-\tau,y)||y|^\kappa dy \|\bar F^2_{x_0}(\tau)\|_{\dot C^\kappa} d\tau \\
					&\lesssim \int_0^t(t-\tau)^{-\frac{1-\kappa}{2}}\tau^{-\frac{1+\kappa}{2}}d\tau \sup_{t\in[0,T]}t^\frac{1+\kappa}{2}\| F^{2}(t)\|_{\dot C^\kappa}\\
					&\lesssim \sup_{\tau\in[0,T]}\tau^\frac{1+\kappa}{2}\| F^{2}(\tau)\|_{\dot C^\kappa},
				\end{aligned}
			\end{equation}
			since $0<\frac{1+\kappa}{2},\frac{1-\kappa}{2}<1$, the third line is integrable. On the other hand, by \eqref{deun1}, we have 
			\begin{equation}\label{unhha}
				\sup_{t\in[0,T]}\sup_{\alpha}t^{\frac{1+\kappa}{2}}\frac{\|(\delta_\alpha \nabla^2v_{N,x_0}^2(t))|_{x_0=x}\|_{L^\infty}}{|\alpha|^{\kappa-[\kappa]}}\lesssim \sup_{\tau\in[0,T]}\tau^{\frac{1+\kappa}{2}}\|\bar{F}^2(t)\|_{\dot C^\kappa}\lesssim\sup_{\tau\in[0,T]}\tau^{\frac{1+\kappa}{2}}\|F^2(t)\|_{\dot C^\kappa}.
			\end{equation}
			We conclude from  \eqref{unlha} and \eqref{unhha}, and notice that they hold for any $x_0$, so especially for $x_0=x$, so we have
			\begin{equation}\label{ununun2}
				\begin{aligned}
					& \sup_{t\in[0,T]}\left(\|(\nabla v_{N,x_0}^2(t))|_{\substack{x_0=x}}\|_{L^\infty(\mathbb{R}^d_+)}+t^\frac{1+\kappa}{2}\sup_{\alpha}\frac{\| (\delta_\alpha\nabla^2 v_{N, x_0}^2(t))|_{\substack{x_0=x}}\|_{L^\infty(\mathbb{R}^d_+)}}{|\alpha|^\kappa}\right)\lesssim \sup_{t\in[0,T]}t^\frac{1+\kappa}{2}\| F^2(t)\|_{\dot C^\kappa}.
				\end{aligned}
			\end{equation}
			For the remainder term $v_{R,x_0}(t,x)$, we follow the idea in Section \ref{secproof} and expect that its contribution to be some lower order terms by fixing $x_0=x$ in the formula, where we essentially gain regularity from the variable coefficient.  Since $\nabla\psi_{x_0}$ is uniformly bounded, we have
			\begin{equation*}
				|\tilde {\mathsf{R}}_{ x_0}(t,x)|\lesssim\left|\B_0-\B(x)\right|\|\nabla^2u(t)\|_{L^\infty}\lesssim |x_0-x|\|\B\|_{\dot {W}^{1,\infty}}\|\nabla^2u(t)\|_{L^\infty}.
			\end{equation*}
			So by \eqref{eskerha}, we take $x_0=x$, and obtain
			\begin{equation*}
				\begin{aligned}
					&\left|(\nabla v_{R,x_0}(t,x))|_{x_0=x}\right|
					\lesssim \int_0^t \int_{\mathbb{R}_+^d}|\nabla_x\H(t-\tau,x'-y',x_d,y_d)||x-y|dy\|\nabla^2u(\tau)\|_{L^\infty}d\tau\|\B\|_{\dot W^{1,\infty}}\\
					&\quad\lesssim T^{\frac{1}{2}}\|\B\|_{\dot W^{1,\infty}}\sup_{\tau\in[0,T]}\tau^{\frac{1}{2}}\|\nabla^2u(\tau)\|_{L^\infty},
					\\
					&\left|(\delta_\alpha\nabla^2 v_{R,x_0}(t,x))|_{x_0=x}\right|\lesssim\int_0^t \int_{\mathbb{R}_+^d}|\delta_\alpha\nabla_x^2\H(t-\tau,x'-y',x_d,y_d)||x-y| dy\tau^{-\frac{\kappa}{2}}\|\nabla^2u(\tau)\|_{L^\infty}d\tau\|\B\|_{\dot W^{1,\infty}}\\
					&\quad\lesssim \|\B\|_{\dot W^{1,\infty}}|\alpha|^{\kappa}t^{-\frac{\kappa}{2}}\sup_{\tau\in[0,T]}\tau^{\frac{1}{2}}\|\nabla^2u(\tau)\|_{L^\infty}.
				\end{aligned}
			\end{equation*}
			The interchange of taking limit $x_0\to x$ and integration is not trivial, and the proof above only holds formally. To get rigorous proof, we can change $\partial_{x_d}^2\H$ into $\lim_{\beta\rightarrow 0}\Delta_\beta\partial_{x_d}\H=\lim_{\beta\rightarrow 0}\frac{\delta_{\beta e_d} \partial_{x_d}\H}{|\beta|}$, and for every $\beta>0$, the interchange is legal.
			These imply that 
			\begin{equation}\label{vrhha}
				\begin{aligned}
					& \sup_{t\in[0,T]}\left(\|(\nabla v_{R,x_0}(t))|_{\substack{x_0=x}}\|_{L^\infty(\mathbb{R}^d_+)}+t^\frac{1+\kappa}{2}\sup_{\alpha}\frac{\|( \delta_\alpha\nabla^2 v_{R,x_0}(t))|_{\substack{x_0=x}}\|_{L^\infty(\mathbb{R}^d_+)}}{|\alpha|^\kappa}\right)\\
					&\quad\quad\quad\quad\quad\lesssim T^{\frac{1}{2}}\|\B\|_{\dot W^{1,\infty}}\sup_{t\in[0,T]}t^{\frac{1}{2}}\|\nabla^2u(t)\|_{L^\infty}.
				\end{aligned}
			\end{equation}
			Combining  \eqref{spliv}, \eqref{lllha}, \eqref{unnha}, \eqref{ununun2}, and \eqref{vrhha} to get
			\begin{equation*}
				\begin{aligned}
					&\sup_{t\in[0,T]}(\|\nabla u(t)\|_{L^\infty}+t^{\frac{1+\kappa}{2}}\|\nabla^2u(t)\|_{\dot C^\kappa})\lesssim \|u\|_{\dot W^{1,\infty}}+\sup_{t\in[0,T]}(t^{\frac{1}{2}}\sup_{t\in[0,T]}\|F^1(t)\|_{L^\infty}+t^{\frac{1+\kappa}{2}}\|F^1(t)\|_{\dot C^\kappa})\\
					&\quad+\sup_{0<\tau<t<T}\tau^{\frac{1+\kappa}{2}}\frac{\|F^1(t)- F^1(\tau)\|_{L^\infty}}{(t-\tau)^\frac{\kappa}{2}}+\sup_{t\in[0,T]}t^{\frac{1+\kappa}{2}}\|F^2(t)\|_{\dot C^{\kappa}}+T^{\frac{1}{2}}\|\B\|_{\dot W^{1,\infty}}\sup_{t\in[0,T]}t^{\frac{1}{2}}\|\nabla^2u(t)\|_{L^\infty}.
				\end{aligned}
			\end{equation*}
			This completes the proof.
		\end{proof}\\
		We prove the following interpolation lemma to control temporal Holder norm.
		\begin{lemma}\label{mcbdlemti}
			For any $f:[0,T]\times \mathbb{R}_+^d\to \mathbb{R}$. There holds
			\begin{align*}
				\sup_{0<\tau<t<T}\frac{\tau^{\frac{\kappa}{2}}\|\nabla f(t)-\nabla f(\tau)\|_{L^\infty(\mathbb{R}^d_+)}}{(t-\tau)^{\frac{\kappa}{2}}}\lesssim \sup_{t\in[0,T]}\left(\|\nabla f(t)\|_{L^\infty(\mathbb{R}^d_+)}+t^{\frac{1}{2}}\|(\partial_t,\nabla^2)f(t)\|_{L^\infty(\mathbb{R}^d_+)}\right).
			\end{align*}
		\end{lemma}
		\begin{proof}
			For $g:[0,T]\times\mathbb{R}^d_+\to \mathbb{R}$, we define the mollification of $g$ on the half space
			\begin{equation*}
				g_\eps(t,x)=\int_{\mathbb{R}_+^d}g(t,y)\tilde\rho_\eps(x'-y',x_d,y_d)dy,
			\end{equation*}
			with 
			\begin{equation*}
				\begin{aligned}
					\tilde \rho_\eps(x'-y',x_d,y_d)=\rho_\eps(x'-y',x_d+y_d)+\rho_\eps(x'-y',x_d-y_d),
				\end{aligned}
			\end{equation*}
			and $\rho_\eps$ is the standard mollification multiplier. 
			For $0<\tau<t<T$, we split $\nabla f(t)-\nabla f(\tau )$ to the following three terms.
			\begin{equation*}
				\begin{aligned}
					|\nabla f(t,x)-\nabla f(\tau,x )|&\lesssim |(\nabla f)_\eps(t,x)-(\nabla f )_\eps(\tau,x)|+|\nabla f (t,x)-(\nabla f)_\eps(t,x)|+|\nabla f(\tau,x )-(\nabla f)_\eps(\tau,x )|\\
					&:=I_1(t,\tau,x)+I_2(t,x)+I_2(\tau,x),
				\end{aligned}
			\end{equation*}
			For $I_1$, using integration by parts we have
			\begin{equation*}
				\begin{aligned}
					&I_1(t,\tau,x)\lesssim\left|\int_\tau ^t\partial_{\tau'}(\nabla f)_\eps(\tau',x)d\tau'\right|\\
					&\lesssim\int_\tau ^t\int_{\mathbb{R}_+^d}|\partial_{\tau'}f(\tau',y)||\nabla_{y}(\tilde \rho_\eps(x'-y',x_d,y_d))|dyd\tau'+\int_\tau^t\int_{\mathbb{R}^{d-1}}|\partial_{\tau'}f(\tau',y',0)||\rho_\eps(x'-y',x_d)|dy'd\tau'\\
					&\lesssim \eps^{-1}|t-\tau |\|\partial_t f\|_{L^\infty([\tau ,t], L^\infty)}\lesssim  \eps^{-1}|t-\tau | \tau ^{-\frac{1}{2}}\sup_{\tau'\in[0,T]}\tau'^{\frac{1}{2}}\|\partial_{\tau'}f(\tau')\|_{L^\infty}.
				\end{aligned}
			\end{equation*}
			For $I_2$, we have
			\begin{equation*}
				\begin{aligned}
					I_2(t,x)\lesssim& \left|\int_{\mathbb{R}_+^d} (\nabla f(t,x)-\nabla f(t,y))\tilde \rho_\eps(x'-y',x_d,y_d)dy \right|\lesssim \eps t^{-\frac{1}{2}}\sup_{t\in[0,T]}\tau'^{\frac{1}{2}}\|\nabla^2f(\tau')\|_{L^\infty}.
				\end{aligned}
			\end{equation*}
			This also implies that  $   I_2(\tau,x)\lesssim \eps \tau^{-\frac{1}{2}}\sup_{t\in[0,T]}\tau'^{\frac{1}{2}}\|\nabla^2f(\tau')\|_{L^\infty}$.
			Take $\eps=|t-\tau |^{\frac{1}{2}}$ to get the following estimate
			\begin{equation*}
				\|\nabla f(t)-\nabla f(\tau )\|_{L^\infty}\lesssim \|I_1(t,\tau)\|_{L^\infty}+\|I_2(t)\|_{L^\infty}+\|I_2(\tau)\|_{L^\infty}\lesssim\tau ^{-\frac{1}{2}}(t-\tau )^{\frac{1}{2}}\sup_{\tau'\in[0,T]}\tau'^{\frac{1}{2}}\|(\partial_{\tau'},\nabla^2)f(\tau')\|_{L^\infty}.
			\end{equation*}
			Since obviously $|\nabla f(t)-\nabla f(\tau )|\lesssim\sup_{\tau'\in[0,T]}\|\nabla f(\tau')\|_{L^\infty}$, by interpolation we have
			\begin{equation*}
				\|\nabla f(t)-\nabla f(\tau )\|_{L^\infty}\lesssim \tau ^{-\frac{\kappa}{2}}(t-\tau )^{\frac{\kappa}{2}}\sup_{t\in[0,T]}\left(\|\nabla f(t)\|_{L^\infty}+t^{\frac{1}{2}}\|(\partial_t,\nabla^2)f(t)\|_{L^\infty}\right).
			\end{equation*}
			This completes the proof.
		\end{proof}\\
		We will need the following two boundary condition lemmas in the case of half space. We postpone the proof of the two lemmas to the end of this subsection.
		\begin{lemma}\label{bdcon}
			Assume $f,g\in C^{2k+2}(\bar{\mathbb{R}}^d_+)$, and 
			\begin{align*}
				\Delta^j g=\Delta^j f=0,\ \ \ \ j=0,1,\cdots,k, \ \ \ x\in\partial\mathbb{R}^d_+.
			\end{align*}
			Then
			\begin{align}\label{Acommu}
				\Delta^k\left(\A[\nabla g]:\nabla^2 f\right)- \A[\nabla g]:\nabla^2\Delta^k f=0,\ \ \ x\in\partial\mathbb{R}^d_+.
			\end{align}
		\end{lemma}
		\begin{lemma}\label{lembdc}
			Assume $g\in C^{2}(\bar{\mathbb{R}}^d_+)$ with $g(x)=0$ for any $x\in\partial\mathbb{R}_+^d$. Let $u\in L^\infty(0,T;C^{2}(\bar{\mathbb{R}}_+^d))$ be a solution to 
			\begin{equation}\label{bddmcbdyeq}
				\begin{aligned}
					&\partial_tu(t,x)-\A[\nabla g](t,x):\nabla^2u(t,x)=F(t,x),\quad \text{in}\ [0,T]\times\mathbb{R}^d_+,\\
					&u(0,x)=u_0(x),\quad \text{in}\ \mathbb{R}^d_+,\\
					&u(t,x)=0,\quad \text{on}\ [0,T]\times\partial\mathbb{R}^d_+,
				\end{aligned}
			\end{equation}
			with
			\begin{equation*}
				F(t,x)=0,\quad \text{on}\ (0,T]\times \partial\mathbb{R}^d_+.
			\end{equation*}
			We have 
			\begin{equation*}
				\Delta u(t,x)=0,\quad\text{on}\ (0,T]\times\partial\mathbb{R}^d_+.
			\end{equation*}
		\end{lemma}
		The following is a key lemma to prove Theorem . 
		\begin{lemma}\label{lemmbd}
			Consider $F=F^1+F^2\in C^\infty([0,T]\times \mathbb{R}_+^d)$ and $F^2|_{[0,T]\times \partial\mathbb{R}_+^d}=0$, $\A\in C^\infty(\mathbb{R}^d;\mathbb{R}^{d\times d})$ is a uniformly elliptic matrix with all derivatives bounded, $\phi\in C^\infty(\mathbb{R}_+^d)$. If $u$ is a solution to the system 
			\begin{equation}\label{mcbd}
				\begin{aligned}
					&\partial_tu(t,x)-\A[\nabla\phi](t,x):\nabla^2u(t,x)=F(t,x)=(F^1+F^2)(t,x),\quad \text{in}\ [0,T]\times\mathbb{R}_+^d,\\
					&u(0,x)=u_0(x),\quad \text{in}\ \mathbb{R}_+^d,\\
					&u(t,x)=0,\quad \text{on}\ [0,T]\times\partial\mathbb{R}_+^d.
				\end{aligned}
			\end{equation}
			Then we have
			\begin{equation}\label{uxha}
				\begin{aligned}
					&\|u\|_{X_T}\lesssim \|u_0\|_{\dot{W}^{1,\infty}}+T^{\frac{\kappa}{2}}(1+T)\|\phi\|_{C^3}\sup_{t\in[0,T]}\left(\|\nabla u(t)\|_{L^\infty}+t^{\frac{1}{2}}\|(\partial_t,\nabla^2)u(t)\|_{L^\infty}\right)\\
					&\quad+\sup_{t\in[0,T]}(t^{\frac{1}{2}}\|F^1(t)\|_{L^\infty}+t^{\frac{1+\kappa}{2}}\|F^1(t)\|_{\dot {C}^\kappa}+t^{\frac{1+\kappa}{2}}\|F^2(t)\|_{\dot{C}^\kappa})+\sup_{0<\tau<t<T}\tau^{\frac{1+\kappa}{2}}\frac{\|F^1(t)-F^1(\tau)\|_{L^\infty}}{(t-\tau)^{\frac{\kappa}{2}}},
				\end{aligned}
			\end{equation}
			Moreover, $\forall \gamma\geq 0$
			\begin{equation}\label{uzha}
				\begin{aligned}
					&\sup_{t\in[0,T]}t^{\gamma+\frac{1+\kappa}{2}}(\|\nabla^2u(t)\|_{\dot C^{\kappa}}+\|\partial_tu(t)\|_{\dot C^{\kappa}})\\
					&\lesssim \sup_{t\in[0,T]}t^{\gamma}\|\nabla u(t)\|_{L^\infty}+T^{\frac{1}{2}}(1+T)\|\phi\|_{C^3}\sup_{t\in[0,T]}t^{\gamma}\left(\|\nabla u(t)\|_{L^\infty}+t^{\frac{1}{2}}\|(\partial_t,\nabla^2)u(t)\|_{L^\infty}\right)\\
					&+\sup_{t\in[0,T]}(t^{\gamma+\frac{1}{2}}\|F^1(t)\|_{L^\infty}+t^{\gamma+\frac{1+\kappa}{2}}\|F^1(t)\|_{\dot C^\kappa}+t^{\gamma+\frac{1+\kappa}{2}}\|F^2(t)\|_{\dot{C}^{\kappa}})+\sup_{0<\tau<t<T}\tau^{\gamma+\frac{1+\kappa}{2}}\frac{\|F^1(t)-F^1(\tau)\|_{L^\infty}}{(t-\tau)^{\frac{\kappa}{2}}}.
				\end{aligned}
			\end{equation}
		\end{lemma}
		\begin{proof}
			Applying Lemma \ref{lemheat} with $\B=\A[\nabla\phi]$, there exists $T>0$ such that 
			\begin{equation}\label{wha}
				\begin{aligned}
					&\sup_{t\in[0,T]}(\|\nabla u(t)\|_{L^\infty(\mathbb{R}_+^d)}+t^{\frac{1+\kappa}{2}}\|\nabla^2u(t)\|_{\dot C^{\kappa}(\mathbb{R}_+^d)})\lesssim \|u_0\|_{\dot W^{1,\infty}(\mathbb{R}_+^d)}+\sup_{0<\tau<t<T}\tau^{\frac{1+\kappa}{2}}\frac{\|F^1(t)-F^1(\tau)\|_{L^\infty(\mathbb{R}_+^d)}}{(t-\tau)^{\frac{\kappa}{2}}}\\
					&\quad+\sup_{t\in[0,T]}\left(t^{\frac{1}{2}}\|F^1(t)\|_{L^\infty(\mathbb{R}^d_+)}+t^{\frac{1+\kappa}{2}}\|F^1(t)\|{\dot C^{\kappa}(\mathbb{R}_+^d)}+t^{\frac{1+\kappa}{2}}\|F^2(t)\|_{\dot C^{\kappa}(\mathbb{R}_+^d)}\right)\\
					&\quad+T^\frac{1}{2}\|\phi\|_{C^3}\sup_{t\in[0,T]}t^\frac{1}{2}\|\nabla ^2 u(t)\|_{L^\infty(\mathbb{R}^d_+)}.
				\end{aligned}
			\end{equation}
			Furthermore, by Lemma \ref{Hol} we can see that
			\begin{equation}\label{ux2ha}
				\begin{aligned}
					&\sup_{t\in[0,T]}t^{\frac{1+\kappa}{2}}\|\partial_tu(t)\|_{\dot C^{\kappa}}\lesssim \sup_{t\in[0,T]}t^{\frac{1+\kappa}{2}}\left(\|F(t)\|_{\dot C^{\kappa}}+\|\A[\nabla \phi]\nabla^2u(t)\|_{\dot C^{\kappa}}\right)\\
					&\lesssim \sup_{t\in[0,T]}t^{\frac{1+\kappa}{2}}\|F(t)\|_{\dot C^{\kappa}}+ \sup_{t\in[0,T]}t^{\frac{1+\kappa}{2}}\|\nabla^2u(t)\|_{\dot C^\kappa}+\sup_{t\in[0,T]}t^{\frac{1+\kappa}{2}}\big(\|\nabla^2u(t)\|_{L^\infty}\|\nabla \phi\|_{\dot C^\kappa}\big) .
				\end{aligned}
			\end{equation}
			Then \eqref{uxha} follows from \eqref{wha} and \eqref{ux2ha}.\\
			It remains to prove \eqref{uzha}. For any $\tau\in(0,T)$, one has 
			\begin{equation}
				\begin{aligned}
					&\sup_{t\in[\tau,T]}(\|\nabla u(t)\|_{L^\infty}+(t-\tau)^{\frac{1+\kappa}{2}}\|(\nabla^2,\partial_t)u(t)\|_{\dot C^{\kappa}})\\
					&\lesssim \|\nabla u(\tau)\|_{L^\infty}+T^{\frac{\kappa}{2}}(1+T)\|\phi\|_{C^3}\sup_{t\in[\tau,T]}(\|\nabla u(t)\|_{L^\infty}+t^{\frac{1+\kappa}{2}}\|(\partial_t,\nabla^2)u(t)\|_{L^\infty})\\
					&\quad+\sup_{t\in[\tau,T]}\big((t-\tau)^{\frac{1}{2}}\|F^1(t)\|_{L^\infty}+(t-\tau)^{\frac{1+\kappa}{2}}\|F^1(t)\|_{\dot{C}^\kappa}\big)+\sup_{\tau<\tau_1<\tau_2<T}(\tau_1-\tau)^{\frac{1+\kappa}{2}}\frac{\|F^1(\tau_1)-F^1(\tau_2)\|_{L^\infty}}{(\tau_1-\tau_2)^{\frac{\kappa}{2}}}\\
					&\quad+\sup_{t\in[\tau,T]}(t-\tau)^{\frac{1+\kappa}{2}}\|F^2(t)\|_{\dot{C}^\kappa}\\
					&\lesssim \tau^{-\gamma}\left(\sup_{t\in[0,T]}t^{\gamma}\| \nabla u(t)\|_{L^\infty}+T^{\frac{\kappa}{2}}(1+T)\|\phi\|_{C^3}\sup_{t\in[0,T]}(t^{\gamma}\|\nabla u(t)\|_{L^\infty}+t^{\gamma+\frac{1+\kappa}{2}}\|(\nabla^2,\partial_t)u(t)\|_{\dot C^{\kappa}})\right.\\
					&\quad+\sup_{t\in[0,T]}\big(t^{\gamma+\frac{1}{2}}\|F^1(t)\|_{L^\infty}+t^{\gamma+\frac{1+\kappa}{2}}\|F^1(t)\|_{ \dot{C}^\kappa}\big)+\sup_{0<\tau_1<\tau_2<T}\tau_1^{\gamma+\frac{1+\kappa}{2}}\frac{\|F^1(\tau_1)-F^1(\tau_2)\|_{L^\infty}}{(\tau_1-\tau_2)^{\frac{\kappa}{2}}}\\
					&\left.\quad+\sup_{t\in[0,T]}t^{\gamma+\frac{1+\kappa}{2}}\|F^2(t)\|_{\dot{C}^\kappa}\right),
				\end{aligned}
			\end{equation}
			which combines the fact that
			\begin{equation*}
				\sup_{t\in[0,T]}t^{\gamma+\frac{1+\kappa}{2}}\|(\nabla^2,\partial_t)u(t)\|_{\dot C^\kappa}\lesssim\sup_{\tau\in[0,\frac{T}{4}]}\tau^\gamma\sup_{t\in[\tau,T]}(t-\tau)^{\frac{1+\kappa}{2}}\|(\nabla^2,\partial_t)u(t)\|_{\dot C^\kappa},
			\end{equation*}
			gives \eqref{uzha}. Then we finish the proof.
		\end{proof}\vspace{0.3cm}\\
		\begin{proof}[Proof of Theorem \ref{eqmcls}]
			We will construct a solution to the system \eqref{eqmcbdd} by a fixed point argument.
			Let $g\in \mathcal{Y}^\sigma_{T,\phi}$, where $\sigma, T,\phi$ will be fixed later. We define a map $\mathcal{S}$ by $\mathcal{S}g=f$, where $f$ is the weak solution to the following system
			\begin{equation}\label{eqmcbd}
				\begin{aligned}
					&\partial_t f=\A[\nabla g]:\nabla^2 f,\ \ \ \text{in}\ [0,T]\times \mathbb{R}_+^d,\\
					&f|_{t=0}=f_0,\ \ \ \text{in}\ \mathbb{R}_+^d,\\
					&f=0,\ \ \text{on}\ [0,T]\times\partial \mathbb{R}_+^d.
				\end{aligned}
			\end{equation}
			We want to show that $\mathcal{S}$ is a contraction map in $\mathcal{Y}_{T,\phi}^\sigma$ for some suitable $T$, $\phi$ and $\sigma$. By Lemma \ref{lembdc}, we have
			\begin{equation*}
				\Delta f|_{\partial\mathbb{R}_+^d}(t)=0,\quad\forall t\in(0,T].
			\end{equation*}
			First consider the equation of $f-\phi$, with $\phi|_{\partial\mathbb{R}_+^d}=0$ and $\|f_0-\phi\|_{\dot W^{1,\infty}}\leq \sigma$ for some $\sigma$ to be fixed later, then we write the equation of $f-\phi$ as
			\begin{equation}
				\begin{aligned}
					&\partial_t (f-\phi)-\A[\nabla \phi]:\nabla^2 (f-\phi)=\A[\nabla \phi]:\nabla^2 \phi+(A[\nabla g])-A[\nabla\phi]):\nabla^2f,\ \ \ \text{in}\ [0,T]\times \mathbb{R}_+^d,\\
					&(f-\phi)|_{t=0}=f_0-\phi ,\ \ \ \text{in}\ \mathbb{R}_+^d,\\
					&f-\phi =0,\ \ \text{on}\ [0,T]\times\partial \mathbb{R}_+^d.
				\end{aligned}
			\end{equation}
			We apply Lemma \ref{lemmbd} with $F^1=\A[\nabla g]:\nabla^2 \phi$, $F_2=(\A[\nabla g])-\A[\nabla\phi]):\nabla^2f$. By the fact that $\Delta f|_{\partial\mathbb{R}_+^d}=0$ and the operator $(\A[\nabla g]-\A[\nabla \phi]):\nabla^2$ keeps invariant under rotation, we have $F_2|_{\partial\mathbb{R}_+^d}=0$. Since 
			\begin{equation*}
				\begin{aligned}
					&\sup_{t\in[0,T]}\|F^1(t)\|_{L^\infty}\lesssim \|\phi\|_{C^2},\\
					&\sup_{t\in[0,T]}t^{\frac{\kappa}{2}}\|F^1(t)\|_{\dot {C}^\kappa}\lesssim (1+\|g-\phi\|_{X_T}+T^{\frac{\kappa}{2}}\|\phi\|_{C^3})\|\phi\|_{C^3},\\
					&\sup_{0<\tau<t<T}\tau^{\frac{1+\kappa}{2}}\frac{\|F^1(t)-F^1(\tau)\|_{L^\infty}}{(t-\tau)^{\frac{\kappa}{2}}}\lesssim T^{\frac{1}{2}}\|\phi\|_{C^3}\|g-\phi\|_{X_T},\\
					&\sup_{t\in[0,T]}t^{\frac{1}{2}}\|F^2(t)\|_{L^\infty}\lesssim \|g-\phi\|_{X_T}(\|f-\phi\|_{X_T}+T^{\frac{1}{2}}\|\phi\|_{C^3}),\\
					&\sup_{t\in[0,T]}t^{\frac{1+\kappa}{2}}\|F^2(t)\|_{\dot C^\kappa}\lesssim (\|g-\phi\|_{X_T}+T^{\frac{\kappa}{2}}(1+T)\|\phi\|_{C^3})(\|f-\phi\|_{X_T}+T^{\frac{1}{2}}\|\phi\|_{C^3}),\\
				\end{aligned}
			\end{equation*}
			one has
			\begin{equation}\label{mchsprxt}
				\begin{aligned}
					\|f-\phi\|_{X_T}\leq& C_1\|f_0-\phi\|_{\dot{W}^{1,\infty}}+C_1(\|g-\phi\|_{X_T}+T^\frac{\kappa}{2}(1+T)\|\phi\|_{C^3})(\|f-\phi\|_{X_T}+T^{\frac{1}{2}}\|\phi\|_{C^3}).
				\end{aligned}
			\end{equation}
			Since $g\in \mathcal{Y}_{T,\phi}^\sigma$, we take $\sigma$, $\phi$ and $T$ such that
			\begin{align*}
				&\sigma\leq \frac{1}{200C_1},\quad \eps_0\leq \frac{\sigma}{200C_1}, \ \ \ T\leq 
				\left(\frac{\eps_0}{200C_1(2+\|\phi\|_{C^{2m+7}})^2}\right)^{10}.
			\end{align*}
			Then
			\begin{equation}\label{bdmcxts}
				\|f-\phi\|_{X_T}\leq 2C_1\eps_0.
			\end{equation}
			Then we estimate the  higher order norm $\|f-\phi\|_{Y_T}$. Denote $f^m=\Delta^mf$. Then one can write the equation of $f^m$ as 
			\begin{equation*}
				\begin{aligned}
					&\partial_t f^m-\A[\nabla \phi]:\nabla^2 f^m=F^m,\ \ \ \text{in}\ (0,T]\times \mathbb{R}_+^d,\\
					&f^m(t,x)=0,\ \ \text{on}\ (0,T]\times\partial \mathbb{R}_+^d,
				\end{aligned}
			\end{equation*}
			with 
			\begin{equation*}
				F^m=\Delta^m(\A[\nabla g]:\nabla^2f)-\A[\nabla g]:\nabla^2\Delta^mf+\left(\A[\nabla g]-\A[\nabla\phi]\right):\nabla^2f^m.
			\end{equation*}
			By Lemma \ref{bdcon}, we have $F^m|_{\partial\mathbb{R}_+^d}=0$. Apply Lemma \ref{lembdc}, we can get $\Delta^{n}f|_{\partial\mathbb{R}_+^d}=0$ with $n\leq m+1$. Apply \eqref{uzha} to $f^m$ with  $\gamma=m$, $F^1=0$, $F^2=F^m$, one has
			\begin{equation}\label{ukha}
				\begin{aligned}
					&\sup_{t\in[0,T]}t^{m+\frac{1+\kappa}{2}}\|(\nabla^2,\partial_t)f^m(t)\|_{\dot C^{\kappa}}\\
					&\quad\quad\quad\leq C_3\sup_{t\in[0,T]}t^{m}\| \nabla f^m(t)\|_{L^\infty}+C_3T^\frac{\kappa}{2}(1+T)\|\phi\|_{C^{3}}\sup_{t\in[0,T]}t^{m}\|\nabla f^m\|_{L^\infty}\\
					&\quad\quad\quad\quad\quad\quad+C_3\sup_{t\in[0,T]}t^{m+\frac{1+\kappa}{2}}\|F^m(t)\|_{\dot{C}^\kappa}.
				\end{aligned}
			\end{equation}
			By standard Schauder estimates, and by the fact that $\Delta^kf|_{\partial\mathbb{R}_+^d}=0$ for any $k\leq m$, we have
			\begin{equation*}
				\|\nabla^{2m}f\|_{\dot{C}^{2+\kappa}}\lesssim\|f^m\|_{\dot{C}^{\kappa}},
			\end{equation*}
			and by interpolation Lemma \ref{maininterpo},
			\begin{equation*}
				\|u\|_{\dot{C}^{m}}\lesssim\|u\|_{\dot{C}^{m_1}}^{\frac{m_2-m}{m_2-m_1}}\|u\|_{\dot{C}^{m_2}}^{\frac{m-m_1}{m_2-m_1}}.
			\end{equation*}
			By Lemma \ref{Hol}, notice that $f$ and $g$ has at least two order derivatives in $F^m$, one has
			\begin{equation}\label{mcbdhsesf}
				\begin{aligned}
					&\sup_{t\in[0,T]}t^{m+\frac{1+\kappa}{2}}\|F^m(t)\|_{\dot{C}^\kappa}\\
					&\quad\leq C_5M_{g,\phi,T}^{2m+3}(\|g-\phi\|_{Y_T}+T^{\frac{\kappa}{2}}(1+T)^{m+2}\|\phi\|_{C^{2k+3}})(\|f-\phi\|_{Y_T}+T^{\frac{\kappa}{2}}(1+T)^{m+2}\|\phi\|_{C^{2m+3}}),
				\end{aligned}
			\end{equation}
			with
			\begin{equation*}
				M_{g,\phi,T}=1+\|g-\phi\|_{Y_T}+T^{\frac{\kappa}{2}}(1+T)^m\|\phi\|_{C^{2m+5}}.
			\end{equation*}
			So by the interpolation Lemma \ref{maininterpo}, we have 
			\begin{equation*}
				\begin{aligned}
					&\sup_{t\in[0,T]}t^{m+\frac{1+\kappa}{2}}\|(\nabla^2,\partial_t)\Delta^m(f-\phi)(t)\|_{\dot C^\kappa}\\
					&\lesssim T^m(1+T)\|\phi\|_{C^{2m+3}}+\sup_{t\in[0,T]}t^{m}\| \nabla \Delta^m (f-\phi)(t)\|_{L^\infty}+M_{g,\phi,T}^{2m+3}(\|g-\phi\|_{X_T}+T^\frac{\kappa}{2}(1+T)^{m+2}\|\phi\|_{C^{3}})\\
					&\quad\quad\quad\quad\times(\sup_{t\in[0,T]}t^{m}\|\nabla \Delta^m(f-\phi)\|_{L^\infty}+\sup_{t\in[0,T]}t^{m+\frac{1+\kappa}{2}}\|(\nabla^2,\partial_t)\Delta^m(f-\phi)\|_{\dot C^{\kappa}}+T^{m+1}\|\phi\|_{C^{2m+3}}).
				\end{aligned}
			\end{equation*}
			Apply interpolation again, together with \eqref{mcbdhsesf}, we can get
			\begin{equation}\label{mchspryt}
				\begin{aligned}
					&\sup_{t\in[0,T]}t^{m+\frac{1+\kappa}{2}}\|(\nabla^2,\partial_t)\Delta^m(f-\phi)(t)\|_{\dot C^\kappa}\leq C_4 T^{\frac{\kappa}{2}}(1+T)^{m+2}\|\phi\|_{C^{2m+3}}\\
					&\quad+C_4M_{g,\phi,T}^{2m+3}(\|g-\phi\|_{X_T}+T^{\frac{1}{2}}(1+T)^{m+3}\|\phi\|_{C^{2k+3}})\sup_{t\in[0,T]}\|\nabla(f-\phi)(t)\|_{L^\infty}.
				\end{aligned}
			\end{equation}
			Combining the estimates above with \eqref{bdmcxts}, and take 
			\begin{align*}
				&\sigma\leq \frac{2^{-m-2}}{500C_1C_4},\quad\eps_0\leq\frac{2^{-2m-3}\sigma}{500C_1C_4},\ \ \ T\leq \left(\frac{2^{-m-2}C_1\eps_0}{200C_3C_4(2+\|\phi\|_{C^{2m+3}})^{2m+10}}\right)^{10},
			\end{align*}
			to ensure
			\begin{equation}
				\begin{aligned}
					\sup_{t\in[0,T]}t^{m+\frac{1+\kappa}{2}}\|(\nabla^2,\partial_t)(f^m-\Delta^m\phi)(t)\|_{\dot C^{\kappa}}\leq  2^{2m+4}C_1C_4\eps_0.
				\end{aligned}
			\end{equation}
			By \eqref{bdmcxts}, \eqref{mchspryt} and the choice of $\eps_0$, we have
			$$
			\|f-\phi\|_{Y_T}\leq \|f-\phi\|_{X_T}+\sup_{t\in[0,T]}t^{m+\frac{1+\kappa}{2}}\|(\nabla^2,\partial_t)\Delta^m(f-\phi)(t)\|_{\dot C^\kappa}\leq \sigma.
			$$
			Hence
			\begin{equation*}
				\mathcal{S}:\mathcal{Y}_{T,\phi}^\sigma\rightarrow\mathcal{Y}_{T,\phi}^\sigma.
			\end{equation*}
			It remains to prove that $\mathcal{S}$ is a contraction map.
			To prove this, we denote $f_1=\mathcal{S}g_1$, $f_2=\mathcal{S}g_2$, $f_1,f_2,g_1,g_2\in\mathcal{Y}_{T,\phi}^\sigma$, and $\mathbf{f}=f_1-f_2$, $\mathbf{g}=g_1-g_2$. We can write the equation of $\mathbf{f}$ as
			\begin{equation}
				\begin{aligned}
					&\partial_t\mathbf f+\A[\nabla g_1]:\nabla^2\mathbf f=(\A[\nabla g_1]-\A[\nabla g_2]):\nabla^2f_2=:\mathbf{F}_1,\quad \text{in}\ [0,T]\times\mathbb{R}_+^d,\\
					&\mathbf{f}(0,x)=0,\quad \text{in}\ \mathbb{R}_+^d,\\
					&\mathbf f(t,x)=0,\quad\text{on}\ [0,T]\times\partial\mathbb{R}_+^d.
				\end{aligned}
			\end{equation}
			Then apply Lemma \ref{lemmbd} for $\mathbf{f}$ with $F^1=0$, $F^2=\mathbf{F}_1$. By the fact that $\Delta f_2|_{\partial\mathbb{R}_+^d}=0$ and the operator $(\A[\nabla g_1]-\A[\nabla g_2])\nabla^2$ keeps invariant under rotation, $F^2|_{\partial\mathbb{R}_+^d}=0$. One has
			\begin{equation}\label{mchsfpxt}
				\begin{aligned}
					\|\mathbf{f}\|_{X_T}\leq C_1 (\|g_1-\phi\|_{X_T}+\|g_2-\phi\|_{X_T}+T^{\frac{\kappa}{2}}(1+T)\|\phi\|_{C^3})\|\mathbf{f}\|_{X_T}+C_1\sup_{t\in[0,T]}t^{\frac{1+\kappa}{2}}\|\mathbf{F}_1(t)\|_{C^\kappa}.
				\end{aligned}
			\end{equation}
			For $\mathbf{F}_1$, since $\mathbf{F}_1|_{\partial\mathbb{R}_+^d}=0$, by Lemma \ref{Hol}, one has
			\begin{equation}\label{mchsfpf1}
				\begin{aligned}
					&\sup_{t\in[0,T]}t^{\frac{1}{2}}\|\mathbf{F}_1(t)\|_{L^\infty}\leq C_5 (1+\sigma+T^{\frac{\kappa}{2}}(1+T)^m\|\phi\|_{C^3})^6\|\mathbf{g}\|_{X_T}(\|f_2-\phi\|_{X_T}+T^{\frac{1}{2}}\|\phi\|_{C^2}),\\
					&\sup_{t\in[0,T]}t^{\frac{1+\kappa}{2}}\|\mathbf{F}_1(t)\|_{\dot C^\kappa}\leq C_5(1+\sigma+T^{\frac{\kappa}{2}}(1+T)\|\phi\|_{C^3})^6\|\mathbf{g}\|_{X_T}(\|f_2-\phi\|_{X_T}+T^{\frac{1}{2}}\|\phi\|_{C^3}).
				\end{aligned}
			\end{equation}
			Since $f_i,g_i\in\mathcal{Y}_{T,\phi}^\sigma$, we can take $\sigma\leq\frac{1}{200C_1C_5}$ and $T\leq (\frac{C_0'}{800C_1(1+C_5)(2+\|\phi\|_{C^3})^7})^{\frac{2}{\kappa}}$ to get
			\begin{equation*}
				\|\mathbf{f}\|_{X_T}\leq C_0'\|\mathbf{g}\|_{X_T},
			\end{equation*}
			for some $C_0'<1$ to be decided later. For $\mathbf{f}^m=\Delta^m\mathbf{f}=\Delta^k(f_1-f_2)$, we can write the equation
			\begin{equation}
				\begin{aligned}
					&\partial_t\mathbf{f}^m+\A[\nabla \phi]:\nabla^2\mathbf{f}^m=\mathbf{F}_1^m,\quad [0,T]\times\mathbb{R}_+^d,\\
					&\mathbf{f}^m(t,x)=0,\quad\text{on}\ (0,T]\times\partial\mathbb{R}_+^d,
				\end{aligned}
			\end{equation}
			with 
			\begin{equation*}
				\mathbf{F}_1^m=\Delta^m((\A[\nabla g_1]-\A[\nabla g_2]):\nabla^2f_2)+\Delta^m(\A[\nabla g_1]:\nabla^2\mathbf{f})-\A[\nabla g_1]:\nabla^2\Delta^m\mathbf{f},\quad\mathbf{F}^m_1|_{\partial\mathbb{R}_+^d}=0.
			\end{equation*}
			By standard Schauder estimate, since $\mathbf{f}^k|_{\partial\mathbb{R}_+^d}=0$ for any $k\leq m$, we have
			\begin{equation*}
				\|\nabla^{2m}\mathbf{f}\|_{C^{\kappa}}\lesssim\|\Delta^m\mathbf{f}\|_{C^{\kappa}}.
			\end{equation*}
			Apply Lemma \ref{lemmbd} for $\mathbf{f}^m$ with $F_1=0$ and $F_2=\mathbf{F}^m_1$, one has \begin{equation*}
				\begin{aligned}
					&\sup_{t\in[0,T]}t^{m+\frac{1+\kappa}{2}}\|(\nabla^2,\partial_t)\mathbf{f}^m(t)\|_{\dot C^\kappa}\leq C_1\sup_{t\in[0,T]}t^{m-\frac{1}{2}}\|\mathbf{f}^m(t)\|_{L^\infty}+C_1\sup_{t\in[0,T]}t^{m+\frac{1+\kappa}{2}}\|\mathbf{F}^m_1(t)\|_{\dot C^\kappa}\\
					&\quad+C_1M_{g,\phi,T}^{2m+8}(\|g-\phi\|_{X_T}+T^\kappa(1+T)^6\|\phi\|_{C^3})\\
					&\quad\quad\quad\quad\times(\sup_{t\in[0,T]}t^{m}\|\nabla \mathbf f^m(t)\|_{L^\infty}+\sup_{t\in[0,T]}t^{m+\frac{1+\kappa}{2}}\|(\partial_t,\nabla^2)\mathbf{f}^m(t)\|_{\dot C^{\kappa}}).
				\end{aligned}
			\end{equation*}
			By interpolation,
			\begin{equation}\label{mchsfpyi}
				\sup_{t\in[0,T]}t^{m-\frac{1}{2}}\|\mathbf{f}^m(t)\|_{L^\infty}\leq \frac{1}{100C_1}\sup_{t\in[0,T]}t^{m+\frac{1+\kappa}{2·}}\|\nabla^2\mathbf{f}^m(t)\|_{\dot C^\kappa}+C_6\sup_{t\in[0,T]}\|\nabla\mathbf{f}(t)\|_{L^\infty}.
			\end{equation}
			For $\mathbf{F}_1^m$, by Lemma \ref{Hol} 
			\begin{equation}\label{mchsfpf2}
				\begin{aligned}
					&\sup_{t\in[0,T]}t^{m+\frac{1+\kappa}{2}}\|\mathbf{F}_1^m(t)\|_{C^\kappa}\leq C_7 (1+\eta+\|\phi\|_{C^{2m+3}})^{2m+3}\\
					&\quad\quad\times\big(\|\mathbf{g}\|_{Y_T}(\|f-\phi\|_{Y_T}+T^{\frac{\kappa}{2}}(1+T)^{m+2}\|\phi\|_{C^{2m+3}})+\|\mathbf{f}\|_{Y_T}(\|g-\phi\|_{X_T}+T^{\frac{\kappa}{2}}(1+T)^{m+2}\|\phi\|_{C^{2m+3}})\big).
				\end{aligned}
			\end{equation}
			Since $g_i,f_i\in\mathcal{Z}_{T,\phi}^\eta$, we can take 
			\begin{align*}
				\sigma<\frac{2^{-m-7}}{100C_1C_7},\quad T<\min\left\{\left(\frac{2^{-m-7}}{100C_1C_7(2+\|\phi\|_{C^{2m+3}})^{2m+9}}\right)^{\frac{2}{\kappa}},\left(\frac{2^{-m-7}}{100C_1C_7(2+\|\phi\|_{C^{2m+3}})^{2m+9}}\right)^2\right\},
			\end{align*}
			and $C_0'$ satisfies $C_0'C_7\leq \frac{1}{100}$ to get
			\begin{equation}
				\|\mathbf{f}\|_{Y_T}\leq\frac{1}{2}\|\mathbf{g}\|_{Y_T},
			\end{equation}
			and by contraction mapping theorem we get the existence and uniqueness of the solution. So we have proved the first part of Theorem \ref{eqmcls}.\\
			For the second part, we will apply $\phi\equiv 0$ and use similar methods. We shortly denote $\mathcal{X}^\sigma=\mathcal{X}_{\infty,0}^\sigma$, $\|\cdot\|_{X}=\|\cdot\|_{X_\infty}$, $\mathcal{Y}^\sigma=\mathcal{Y}_{\infty,0}^\sigma$ and $\|\cdot\|_{Y}=\|\cdot\|_{Y_\infty}$. Note that by taking $\phi\equiv 0$ in \eqref{mchsprxt}, we have
			\begin{equation}\label{mchsprxi}
				\|f\|_{X}\leq C_1\|f_0\|_{\dot{W}^{1,\infty}}+C_1(1+\|g\|_{X})\|g\|_{X}\|f\|_{X}.
			\end{equation}
			By taking $\sigma$ and $\eps_0$ small enough such that $C_1(1+\sigma)\sigma\leq\frac{1}{100C_1}$, one can prove 
			\begin{equation*}
				\|f\|_{X}\leq 2C_1\eps_0.
			\end{equation*}
			Furthermore, by taking $\phi\equiv 0$ in \eqref{mchspryt}, one can prove
			\begin{equation}\label{mchspryi}
				\sup_{t\in(0,T]}t^{m+\frac{1+\kappa}{2}}\|(\nabla^2,\partial_t)\Delta^mf(t)\|_{\dot{C}^\kappa}\leq C_4(1+\|g\|_{Y})^{2m+3}\|g\|_{Y}\sup_{t\in(0,T]}\|\nabla f(t)\|_{L^\infty}.
			\end{equation}
			So by taking $\sigma$ and $\eps_0$ small enough such that $C_1C_4\eps_0(1+\sigma)^{2m+3}\leq\frac{\sigma}{100}$, then one has 
			\begin{equation*}
				\mathcal{S}:\mathcal{Y}^\sigma\rightarrow\mathcal{Y}^\sigma.
			\end{equation*}
			For compression argument, by taking $\phi=0$ in \eqref{mchsfpxt} and \eqref{mchsfpf1}, one has
			\begin{equation*}
				\|\mathbf{f}\|_{X}\leq C_1(\|g_1\|_X+\|g_2\|_X)\|\mathbf{f}\|_X+C_1C_5(1+\sigma)^6\|\mathbf{g}\|_X\|f_2\|_X.
			\end{equation*}
			By taking $\sigma\leq \frac{C'_0}{500C_1C_5}$, one has
			\begin{equation}\label{mchsfpfx}
				\|\mathbf{f}\|_X\leq C'_0\|\mathbf{g}\|_X,
			\end{equation}
			with $C'_0$ to be decided later. Finally by taking $\phi=0$ in \eqref{mchsfpyi} and \eqref{mchsfpf2} and by interpolation, one has 
			\begin{equation*}
				\begin{aligned}
					&\sup_{t\in(0,T]}t^{m+\frac{1+\kappa}{2}}\|(\nabla^2,\partial_t)\mathbf{f}^m(t)\|_{\dot {C}^\kappa}\\
					&\quad\leq C'_1(1+\|g\|_Y)^{2m+3}\left(\|g\|_Y\sup_{t\in(0,T]}\|\mathbf{f}(t)\|_{L^\infty}+\|\mathbf{g}\|_Y(f_1\|_Y+\|f_2\|_Y)+\|\mathbf{f}\|_Y(\|g_1\|_Y+\|g_2\|_Y)\right).
				\end{aligned}
			\end{equation*}
			The estimate above combined with \eqref{mchsfpfx}, and take $C'_0$, $\sigma$ small enough, then we can get
			\begin{equation*}
				\|\mathbf{f}\|_Y\leq\frac{1}{2}\|\mathbf{g}\|_Y,
			\end{equation*}
			which completes the proof by the compression mapping theorem.
		\end{proof}\vspace{0.3cm}\\
		In the following, we prove Lemma \ref{bdcon} and \ref{lembdc}.\\
		\begin{proof}[Proof of Lemma \ref{bdcon}]
			By boundary condition, one has 
			\begin{align}\label{dex0}
				( \partial_1^{n_1}\cdots\partial_{d-1}^{n_{d-1}} \partial_d^{2j})f(x_0)=0,\ \ \ \forall n_1,\cdots,n_{d-1}\in\mathbb{N},\ \text{and}\ 0\leq j\leq k, 
			\end{align}
			and 
			\begin{align}\label{dex0g}
				( \partial_1^{n_1}\cdots\partial_{d-1}^{n_{d-1}} \partial_d^{2j})g(x_0)=0,\ \ \ \forall n_1,\cdots,n_{d-1}\in\mathbb{N},\ \text{and}\ 0\leq j\leq k. 
			\end{align}
			We can write 
			\begin{equation}\label{deldif}
				\begin{aligned}
					&|\Delta^k\left(\A[\nabla g]:\nabla^2 f\right)- \A[\nabla g]:\nabla^2\Delta^k f|\\
					&\quad\quad\quad\quad\leq \sum_{j,l} |\partial_1^{j_1}\partial_2^{j_2}\cdots\partial_d^{j_d}\A[\nabla g]:\nabla^2\partial_1^{l_1}\partial_2^{l_2}\cdots\partial_d^{l_d}f|,
				\end{aligned}
			\end{equation}
			here the summation is taken for some $j=(j_1,\cdots,j_d)$, $l=(l_1,\cdots,l_d)$ satisfying
			\begin{align}\label{indjl}
				\sum_{i=1}^d(j_i+l_i)=2k,\ \ \sum_{i=1}^dj_i> 0,\ \text{and}\  \{j_i+l_i\}_{i=1}^d\subset 2\mathbb{N},
			\end{align}
			where we denote $2\mathbb{N}=\{2n:n\in\mathbb{N}\}$ the set of even integers.\\
			In the following, we prove that the Frobenius inner product in the right hand side of \eqref{deldif} is $0$ for any $j,l$ satisfying \eqref{indjl}. We remark that, no cancellation property will be used to prove this. We only need to verify whether each element of the  matrices $\partial_1^{j_1}\partial_2^{j_2}\cdots\partial_d^{j_d}\A[\nabla g](x_0)$ and $\nabla^2\partial_1^{l_1}\partial_2^{l_2}\cdots\partial_d^{l_d}f(x_0)$ is $0$ or not. 
			
			For simplicity, for any $M\in\mathbb{R}^{d\times d}$, we define its indicator matrix $\mathbf{I}(M)$ as follows.
			\begin{align}
				\mathbf{I}_{ij}(M)=\begin{cases}
					1,\ \ \ \ \text{if}\ M_{ij}\neq 0,\\
					0,\ \ \ \ \text{if}\ M_{ij}= 0.\\
				\end{cases}
			\end{align}
			Note that any derivatives in $(x_1,\cdots,x_{d-1})$ do not change the indicator matrix. Hence it suffices to consider 
			$\mathbf{I}(\partial_d^{j}\A[\nabla g](x_0))$  and $\mathbf{I}(\nabla^2\partial_d^{l}f(x_0))$ with $j\leq 2k, l<2k$.
			We claim that, \\
			(i). If $l,j\in 2\mathbb{N}$, $j\leq 2k, l<2k$, then 
			\begin{align}
				&\mathbf{I}(\partial_d^{j}\A[\nabla g](x_0))\leq\left(\begin{array}{cc}
					1_{(d-1)\times(d-1) }  & 0_{d\times 1}  \\
					0_{1\times d}  & 1
				\end{array}\right),\label{eveA}\\    &\mathbf{I}(\nabla^2\partial_d^{l}f(x_0))\leq\left(\begin{array}{cc}
					0_{(d-1)\times(d-1) }  & 1_{d\times 1}  \\
					1_{1\times d}  & 0
				\end{array}\right).\label{evef}
			\end{align}
			Here $\leq $ is defined element-wise.\\ 
			(ii). If $l,j\in \mathbb{N}\backslash 2\mathbb{N}$, $j\leq 2k, l<2k$, then 
			\begin{align}
				&\mathbf{I}(\partial_d^{j}\A[\nabla g](x_0))\leq\left(\begin{array}{cc}
					0_{(d-1)\times(d-1) }  & 1_{d\times 1}  \\
					1_{1\times d}  & 0
				\end{array}\right),\label{oddA}\\    &\mathbf{I}(\nabla^2\partial_d^{l}f(x_0))\leq\left(\begin{array}{cc}
					1_{(d-1)\times(d-1) }  & 0_{d\times 1}  \\
					0_{1\times d}  & 1
				\end{array}\right).\label{oddf}
			\end{align}
			Then it follows from (i), (ii) and \eqref{deldif} that 
			\begin{align*}
				\Big(\Delta^k\left(\A[\nabla g]:\nabla^2 f\right)- \A[\nabla g]:\nabla^2\Delta^k f\big)(x_0)=0.
			\end{align*}
			This holds for any $x_0\in\partial\mathbb{R}_+^d$. This completes the proof of \eqref{Acommu}.
			
			Now we prove (i) and (ii). By \eqref{dex0} we obtain \eqref{evef} and \eqref{oddf}. It remains to verify \eqref{eveA} and \eqref{oddA}. 
			
			We first write $A[\nabla g]$ in terms of  Neumann series 
			\begin{align}\label{series}
				\A[\nabla g]=\left[\mathrm{Id}+\sum_{i=1}^n\nabla g^i\otimes \nabla g^i\right]^{-1}=\sum_{m=0}^\infty(-1)^m(\mathsf{P}^{-1}\B)^m\mathsf{P}^{-1},
			\end{align}
			where\begin{align*}
				\mathsf{P}=\mathrm{Id}+|\partial_d g|^2 e_d\otimes e_d,\ \ \ \B=\sum_{i=1}^n\nabla g^i\otimes \nabla g^i-|\partial_d g|^2 e_d\otimes e_d.
			\end{align*}
			We remark that the series in the right hand side of \eqref{series} is well-defined in a small 
			neighbourhood of $x_0$ because $|\mathsf{P}^{-1}|\leq 1$ and $\B(x)\lesssim O(|x-x_0|)$ for any $|x-x_0|\ll 1$.
			Moreover,
			it suffices to consider finite terms in the series \eqref{series}. In fact, one can check that all the elements in the matrix
			$
			(\mathsf{P}^{-1}\B)^2$ contain a factor in  the form $\partial_{i_1}g\partial_{i_2}g$, where $i_1,i_2<d$. Hence, all the elements in $(\mathsf{P}^{-1}\B)^m\mathsf{P}^{-1}$ contain a factor in the form $\prod_{k=1}^{m-1}\partial_{i_k}g$, with $i_k<d$. Then by \eqref{dex0g}, we have 
			\begin{align*}
				\left(\partial_d^{j}\left((\mathsf{P}^{-1}\B)^m\mathsf{P}^{-1}\right)\right)(x_0)=0,\ \ \ \forall m>j+1.
			\end{align*}
			Hence we only need to consider $\mathsf{D}_m=(\mathsf{P}^{-1}\B)^m\mathsf{P}^{-1}$ with $0\leq m\leq j+1$. To prove \eqref{eveA} and \eqref{oddA}, it suffices to prove  for any $0\leq m\leq j+1$, 
			\begin{equation}\label{indum}
				\begin{aligned}
					&\mathbf{I}(\partial_d^{j}\mathsf{D}_m(x_0))\leq
					\left(\begin{array}{cc}
						1_{(d-1)\times(d-1) }  & 0_{d\times 1}  \\
						0_{1\times d}  & 1
					\end{array}\right),\ \forall j\in 2\mathbb{N},j\leq 2k,\\
					&	    \mathbf{I}(\partial_d^{j}\mathsf{D}_m(x_0))\leq \left(\begin{array}{cc}
						0_{(d-1)\times(d-1) }  & 1_{d\times 1}  \\
						1_{1\times d}  & 0
					\end{array}\right),\ \forall j\in \mathbb{N}\backslash 2\mathbb{N},j\leq 2k.
				\end{aligned}
			\end{equation}
			When $m=0$, $\mathsf{D}_0=\mathsf{P}^{-1}$. Then by \eqref{dex0g} it is easy to check that \eqref{indum} holds for $m=0$.
			By induction, assume \eqref{indum} holds for $m\leq j$, we prove that it also holds for $m+1$. We can write 
			\begin{align*}
				\mathsf{D}_{m+1}=(\mathsf{P}^{-1}\B)^{m+1}\mathsf{P}^{-1}=\mathsf{D}_m \B\mathsf{P}^{-1}.
			\end{align*}
			And 
			\begin{align*}
				\partial_d^j\mathsf{D}_{m+1}=\sum_{j'=0}^j\partial_d^{j'}\mathsf{D}_m\ \partial_d^{j-j'}(\B\mathsf{P}^{-1}).
			\end{align*}
			Note that 
			\begin{align*}
				\B\mathsf{P}^{-1}=\left(\begin{array}{cc}
					\nabla'g\otimes \nabla'g & \frac{\partial_d g (\nabla'g)^\top  }{1+|\partial_d g|^2} \\
					\partial_d g \nabla'g & 0
				\end{array}\right).
			\end{align*}
			By \eqref{dex0g}, we have 
			\begin{equation}\label{indi}
				\begin{aligned}
					&\mathbf{I}(\partial_d^{j-j'}(\B\mathsf{P}^{-1})(x_0))\leq\left(\begin{array}{cc}
						0_{(d-1)\times (d-1)} & 1_{d\times 1}  \\
						1_{1\times d} & 0.
					\end{array}\right),\ \ \ \ j-j'\in \mathbb{N}\backslash 2\mathbb{N},\\
					&\mathbf{I}(\partial_d^{j-j'}(\B\mathsf{P}^{-1})(x_0))\leq\left(\begin{array}{cc}
						1_{(d-1)\times (d-1)} & 0_{d\times 1}  \\
						0_{1\times d} & 0
					\end{array}\right),\ \ \ \ j-j'\in 2\mathbb{N}.
				\end{aligned}
			\end{equation}
			If $j\in 2\mathbb{N}$, then we have $j',j-j'\in 2\mathbb{N}$ or $j',j-j'\in \mathbb{N}\backslash 2\mathbb{N}$, then by \eqref{indum} and \eqref{indi} we have
			\begin{align*}
				\mathbf{I}\left(\partial_d^{j'}\mathsf{D}_m\ \partial_d^{j-j'}(\B\mathsf{P}^{-1})\right)&= \mathbf{I}\left(  \left(\begin{array}{cc}
					1_{(d-1)\times(d-1) }  & 0_{d\times 1}  \\
					0_{1\times d}  & 1
				\end{array}\right)\left(\begin{array}{cc}
					1_{(d-1)\times (d-1)} & 0_{d\times 1}  \\
					0_{1\times d} & 0
				\end{array}\right)\right)\\
				&\leq\left(\begin{array}{cc}
					1_{(d-1)\times (d-1)} & 0_{d\times 1}  \\
					0_{1\times d} & 0
				\end{array}\right),\ \ \ j',j-j'\in 2\mathbb{N},
			\end{align*}
			and 
			\begin{align*}
				\mathbf{I}\left(\partial_d^{j'}\mathsf{D}_m\ \partial_d^{j-j'}(\B\mathsf{P}^{-1})\right)&= \mathbf{I}\left(  \left(\begin{array}{cc}
					0_{(d-1)\times(d-1) }  & 1_{d\times 1}  \\
					1_{1\times d}  & 0
				\end{array}\right)\left(\begin{array}{cc}
					0_{(d-1)\times (d-1)} & 1_{d\times 1}  \\
					1_{1\times d} & 0
				\end{array}\right)\right)\\
				&\leq\left(\begin{array}{cc}
					1_{(d-1)\times (d-1)} & 0_{d\times 1}  \\
					0_{1\times d} & 1
				\end{array}\right),\ \ \ j',j-j'\in \mathbb{N}\backslash 2\mathbb{N}.
			\end{align*}
			Hence we obtain 
			\begin{align}\label{evre}
				\mathbf{I}(\partial_d^j\mathsf{D}_{m+1}(x_0))\leq\left(\begin{array}{cc}
					1_{(d-1)\times (d-1)} & 0_{d\times 1}  \\
					0_{1\times d} & 1
				\end{array}\right),\ \ \ j\in2\mathbb{N}, j\leq 2k.
			\end{align}
			Similarly, if $j\in\mathbb{N}\backslash 2\mathbb{N}$, then we have $j'\in 2\mathbb{N}$, $j-j'\in \mathbb{N}\backslash 2\mathbb{N}$ or  $j'\in \mathbb{N}\backslash 2\mathbb{N}$, $j-j'\in 2\mathbb{N}$, then by \eqref{indum} and \eqref{indi} we have
			\begin{align*}
				\mathbf{I}\left(\partial_d^{j'}\mathsf{D}_m\ \partial_d^{j-j'}(\B\mathsf{P}^{-1})\right)&= \mathbf{I}\left(  \left(\begin{array}{cc}
					1_{(d-1)\times(d-1) }  & 0_{d\times 1}  \\
					0_{1\times d}  & 1
				\end{array}\right)\left(\begin{array}{cc}
					0_{(d-1)\times (d-1)} & 1_{d\times 1}  \\
					1_{1\times d} & 0
				\end{array}\right)\right)\\
				&\leq\left(\begin{array}{cc}
					0_{(d-1)\times (d-1)} & 1_{d\times 1}  \\
					1_{1\times d} & 0
				\end{array}\right),\ \ \ j'\in 2\mathbb{N}, j-j'\in \mathbb{N}\backslash 2\mathbb{N},
			\end{align*}
			and 
			\begin{align*}
				\mathbf{I}\left(\partial_d^{j'}\mathsf{D}_m\ \partial_d^{j-j'}(\B\mathsf{P}^{-1})\right)&= \mathbf{I}\left(  \left(\begin{array}{cc}
					0_{(d-1)\times(d-1) }  & 1_{d\times 1}  \\
					1_{1\times d}  & 0
				\end{array}\right)\left(\begin{array}{cc}
					1_{(d-1)\times (d-1)} & 0_{d\times 1}  \\
					0_{1\times d} & 0
				\end{array}\right)\right)\\
				&\leq\left(\begin{array}{cc}
					0_{(d-1)\times (d-1)} & 0_{d\times 1}  \\
					1_{1\times d} & 0
				\end{array}\right),\ \ \ j',j-j'\in \mathbb{N}\backslash 2\mathbb{N}.
			\end{align*}
			This implies that 
			\begin{align}\label{odre}
				\mathbf{I}(\partial_d^j\mathsf{D}_{m+1}(x_0))\leq\left(\begin{array}{cc}
					0_{(d-1)\times (d-1)} & 1_{d\times 1}  \\
					1_{1\times d} & 0
				\end{array}\right),\ \ \ j\in\mathbb{N}\backslash 2\mathbb{N}, j\leq 2k.
			\end{align}
			By \eqref{evre} and \eqref{odre}, we obtain that \eqref{indum} holds for $m+1$. This yields \eqref{eveA} and \eqref{oddA} and completes the proof.
		\end{proof}\\
		\begin{proof}[Proof of Lemma \ref{lembdc}]
			Note that at $x_0\in\partial\mathbb{R}_+^d$, by the condition $g|_{\partial\mathbb{R}_+^d}=0$, we know that $\nabla g^k(x_0)$ is parallel to $\mathbf{n}(x_0), \ k=1,2,\cdots,n$. Hence $\nabla g^k(x_0)$ is parallel to $e_d$ for any $k=1,2,\cdots,n$. This implies that $\A[\nabla g( x_0)]$ is a diagonal matrix. More precisely, we have
			\begin{equation}\label{mcbdfma}
				(\A[\nabla  g]:\nabla^2 u)( x_0)=\begin{pmatrix}
					1&0&\cdots&0&0\\
					0&1&\cdots&0&0\\
					\vdots&\vdots&\ddots&\vdots&\vdots\\
					0&0&\cdots&1&0\\
					0&0&\cdots&0&\frac{1}{1+|\partial_{x_d} g|^2}
				\end{pmatrix}:\nabla^2 u=\left(\Delta_{x'} u+\frac{\partial_{x_d}^2 u}{1+|\partial_{x_d} g|^2}\right)( x_0).
			\end{equation}
			Since $u|_{\partial\mathbb{R}_+^d}=0$, one has $\partial_t u|_{\partial\mathbb{R}_+^d}=0$. Moreover, we have $F|_{\partial\mathbb{R}_+^d}=0$, so by \eqref{bddmcbdyeq}, one has
			\begin{equation}\label{mcbdbdyv}
				\A[\nabla g]:\nabla^2u(x_0)=0.
			\end{equation}
			Combining this with \eqref{mcbdfma}, and the fact that 
			$\Delta_{x'}u( x_0)=0$, one has
			\begin{equation*}
				\partial_{x_d}^2 u( x_0)=0.
			\end{equation*}
			So we have 
			\begin{equation*}
				\Delta u(t,x_0)=0,\ \ \ \forall t\in(0,T].
			\end{equation*}
			Since these procedures holds for any $x_0\in\partial\mathbb{R}_+^d$, one has 
			\begin{equation*}
				\Delta u(t,x)=0,\quad\text {on}\ (0,T]\times \partial\mathbb{R}_+^d.
			\end{equation*}
			This completes the proof.
		\end{proof}\\
		\begin{remark}\label{bdydiff}
			The case of bounded domain is different from the half space. The boundary of half space is a flat surface, then we can obtain that some higher order derivatives of the solution remain $0$ on the boundary (see Lemma \ref{lembdc}). However, for general bounded domain where the boundary is a curved surface, we don't have this good property. In fact, by straightening out the boundary locally, we can only obtain that the first order tangential derivatives on the boundary remain $0$, higher order one loses this property. We need to analyse the boundary terms in more detail.
		\end{remark}
		\subsection{Well-posedness of the Mean curvature flow on bounded domain}
		In this section, we will prove the existence and regularity to the local and global solution of mean curvature system \eqref{eqmcbdd} on a bounded domain. The main result is stated in Theorem \ref{mcbdglo}. \vspace{0.1cm}\\
		As preparation, we first introduce some notation. For a bounded domain $\Omega\subset \mathbb{R}^d$, we take an open cover $\Omega_0\cup_{i=1}^K\Omega_i\supset \Omega$ such that $\bar{\Omega}_0\Subset\Omega$, $\Omega_i\cap\partial\Omega\neq\phi$. The way to choose this cover will be shown later. Then we define two families of smooth cut-off functions $\{\chi_i\}_{i=0}^K$ and $\{\tilde{\chi}_i\}_{i=0}^K$, such that $\Omega_i\subset\operatorname{Supp}(\chi_i)$, and $\operatorname{Supp}(\chi_i)\subset\tilde{\chi}^{-1}_i(1)$. Then, for any $f\in C^\infty(\Omega)$, consider $\{\tilde{\chi}_if\}_{i=1}^K$, there exists flatten transforms $\{\Phi_i\}_{i=1}^K$, such that $\{(\tilde{\chi}_if)\circ\Phi_i\}_{i=1}^k:\mathbb{R}_+^d\rightarrow\mathbb{R}$. Precisely, with a rotation, one can take $\partial\Omega\cap\operatorname{Supp}(\tilde{\chi}_i)$ to be the graph of a map $\varphi_i:\mathbb{R}^{d-1}\to\mathbb{R}$, which means that $\partial\Omega\cap\operatorname{Supp}(\tilde{\chi}_i)\subset\{(x',x_d)\in B(0,100)|
		x_d=\varphi(x')\}$, with $0\in\partial\Omega\cap\operatorname{Supp}(\chi_i)$ and $\nabla\varphi_i(0)=0$.
		\\
		The heat semigroup with Dirichlet boundary condition on bounded domain is defined as the solution $v=e^{t\Delta}u$ to the following system,
		\begin{align*}
			&\partial_tv-\Delta v=0,\quad\text{in}\ [0,\infty)\times\Omega,\\
			&v(0,x)=u(x),\quad\text{in}\ \Omega,\\
			&v(t,x)=0,\quad\text{on}\ [0,\infty)\times\partial\Omega.
		\end{align*}
		We denote  \begin{align*}
			\phi=e^{\eps_1\Delta} f_0.
		\end{align*}	
		Furthermore, since the operator $\A[\nabla g]:\nabla^2$ keeps invariant with rotation, we can make a rotation such that the normal vector at $0$ satisfies $\mathbf{n}=e_d$. By similar argument with \eqref{mcbdfma}, we have 
		\begin{equation*}
			\A[\nabla (\tilde{\chi}_i\phi)](0)=\operatorname{diag}(1,...1,(1+|\nabla \phi|^2(0))^{-1}).
		\end{equation*}
		Without loss of generality, we can take the cut off functions $\{\tilde{\chi}_i\}_{i=1}^K$ such that for any $i$, $|\A[\nabla (\tilde{\chi}_i\phi)](x)-\A[\nabla(\tilde{\chi}_i\phi)(0)]|\lesssim \eps$ for any $x\in\operatorname{Supp}(\tilde{\chi}_i)$ with $\eps$ to be fixed later. By the definition of $\phi$, we may see that the scale of $\operatorname{Supp}(\tilde{\chi}_i)$ is nearly $\eps_1^{\frac{1}{2}}\eps$.\\
		Denote $\nabla=(\nabla',\partial_d)$ with $\nabla'=(\partial_1,\cdots,\partial_{d-1})$, and 
		\begin{equation*}
			\nabla_{x'}^\alpha f=(\partial_{x_1}^{\alpha_1}\cdots\partial_{x_{d-1}}^{\alpha_{d-1}}f),\quad \alpha=(\alpha_1,\cdots,\alpha_{d-1}).
		\end{equation*}
		For any norm $\|\cdot\|_Z$, we write
		\begin{equation*}
			\|\nabla_{x'}^n f\|_Z=\sum_{|\alpha|=n}\|\partial_{x_1}^{\alpha_1}\cdots\partial_{x_{d-1}}^{\alpha_{d-1}}f\|_{Z},\quad|\alpha|=\sum_{i=1}^{d-1}\alpha_i,
		\end{equation*}
		To prove Theorem \ref{mcbdglo}, we will do the interior estimate and boundary estimate seperately. For the sake of brevity, we define the norms $\|\cdot\|_{Z_T^0,1}$, $\|\cdot\|_{Z_T,1}$ for the interior estimate, and $\|\cdot\|_{Z_T^0,2}$, $\|\cdot\|_{Z_T,2}$ for the boundary estimate. \\
		For $f:[0,T]\times\mathbb{R}^d\to  \mathbb{R}^N$,  define 
		\begin{align*}
			&\|f\|_{Z_T^0,1}=\sup_{t\in(0,T]}\left(\| f(t)\|_{C^1}+t^{\frac{1+\kappa}{2}}\|\nabla^{2}f(t)\|_{C^{\kappa}}\right),\\
			&\|f\|_{Z_T,1}=\sup_{t\in(0,T]}\left(\| f(t)\|_{C^1}+t^{\frac{2m+1+\kappa}{2}}\|\nabla^{2m+2}f(t)\|_{C^{\kappa}}\right). 
		\end{align*}
		For $f:[0,T]\times\mathbb{R}^d_{+}\to  \mathbb{R}^N$,  define 
		\begin{equation}\label{defnormZ}
			\begin{aligned}
				&\|f\|_{Z_T^0,2}=\sup_{t\in[0,T]}\left(\|f(t)\|_{C^1}+t^{\frac{1+\kappa}{2}}\|\nabla_{x'}\nabla f(t)\|_{C^{\kappa}}\right),\\
				&\|f\|_{Z_T,2}=\sup_{t\in[0,T]}\left(\|f(t)\|_{C^1}+\sum_{n+2k+2l\leq 2m}t^{\frac{n+2k+2l+1+\kappa}{2}}\|\nabla^{n+1}_{x'}\partial_{d}^{2k}\nabla\partial_t^lf(t)\|_{C^{\kappa}}\right).
			\end{aligned}
		\end{equation}
		We remark that $\|f\|_{Z_T^0,2}$ and $\|f\|_{Z_T,2}$ are anisotropic norms.
		Let $f:[0,T]\times\Omega\to \mathbb{R}^N$. We define 
		\begin{equation}\label{mcbdnmzt}
			\begin{aligned}
				\|f\|_{Z^0_T}:=&\|\chi_0f\|_{Z_T^0,1}+\sum_{i=1}^K\|(\chi_if)\circ\Phi_i\|_{Z_T^0,2},
			\end{aligned}
		\end{equation}
		and the corresponding higher order form
		\begin{equation*}
			\begin{aligned}
				&\|f\|_{Z_T}:=\|\chi_0f\|_{Z_T,1}+\sum_{i=1}^K\|(\chi_if)\circ\Phi_i\|_{Z_T,2}.
			\end{aligned}
		\end{equation*}
		We define the set
		\begin{equation*}
			\mathcal{Z}_{T,\phi}^\sigma=\{f\in L^\infty_TC^{2m+2}(\bar{\Omega}):f|_{(0,T]\times\partial\Omega}=0,\|f-\phi\|_{Z_T}\leq \sigma\}.
		\end{equation*}
		Since $\{\Phi_i\}_{i=1}^K$ are uniformly bounded maps, we can see that 
		\begin{equation*}
			\sup_{t\in[0,T]}\left(\|f\|_{C^{1}}+t^{\frac{2m+\kappa}{2}}\|f\|_{C^{2m+1+\kappa}}\right)\lesssim\|f\|_{Z_T}.
		\end{equation*}\vspace{0.1cm}\\
		Analogous to Lemma \ref{lembdc} in the previous section for the half space, we have the following lemma considering the behaviour of solution at the boundary. As we mentioned in Remark \ref{bdydiff}, unlike the half-space case, $\partial_d^2 f$ is not zero at boundary when considering general bounded domains. Thanks to the following lemma, we proved that higher order normal derivatives $\partial_d^{2k}f$ can be decomposed into a sum of terms with lower order normal derivatives on the boundary.
		\begin{lemma}\label{mcbdlembdy}
			If $f$ is the solution of the following system, 
			\begin{equation}\label{mcbdbdeq}
				\begin{aligned}
					&\partial_t f-\Delta f=F,\ \ \ \text{in}\ [0,T]\times \mathbb{R}_+^d,\\
					&f|_{t=0}=f_0,\ \ \ \text{in}\ \mathbb{R}_+^d,\\
					&f=0,\ \ \text{on}\ [0,T]\times\partial 
					\mathbb{R}_+^d,
				\end{aligned}
			\end{equation}
			then for any $k\in\mathbb{N}$, we have the boundary condition
			\begin{equation}\label{mcbdbdyfr}
				\begin{aligned}
					&\partial_d^{2k}f=F^{k-1}+R_{2k},\quad x\in \partial\mathbb{R}_+^d,\\
					&F^{k-1}=\sum_{2i+|\alpha|=2k-2}C^1_{i,\alpha}\partial_t^i\nabla^{\alpha}F,\quad R_{2k}=\sum_{\substack{l+|\beta|= 2k\\l<2k}}C^2_{l,\beta}\partial_{d}^{l}\nabla_{x'}^{\beta}f,
				\end{aligned}
			\end{equation}
			Furthermore, we can do similar decomposition for $\partial_d^{2k}\nabla_{x'}^\beta f$ such that
			\begin{equation}\label{mcbdbdyfrm}
				\begin{aligned}
					&\partial_d^{2k}\nabla_{x'}^\beta f=F^{k-1}_\beta+R^{2k}_\beta,\quad x\in(0,T]\times \partial\mathbb{R}_+^d\\
					&F^{k-1}_\beta=\sum_{2i+|\alpha|=2k-2}C^{1,\beta}_{i,\alpha}\partial_t^{i}\nabla^{\alpha}\nabla_{x'}^\beta F,\quad R^{2k}_n=\sum_{\substack{l+|\zeta|= 2k\\l<2k}}C_{l,\zeta}^{2,\beta}\partial_{d}^{l}\nabla_{x'}^{\zeta+\beta}f.
				\end{aligned}
			\end{equation}
		\end{lemma}
		\begin{proof}
			We prove by induction. First, for $k=1$ and $x_0\in\partial 
			\mathbb{R}_+^d$, since $f|_{\partial 
				\mathbb{R}_+^d}=0$, we have $\partial_tf|_{\partial 
				\mathbb{R}_+^d}=0$, then by the equation \eqref{mcbdbdeq} we have
			\begin{equation*}
				\begin{aligned}
					\partial_{dd}f=-F-\Delta_{x'}f:=-F-R_2,\quad\text{on}\ (0,T]\times \partial 
					\mathbb{R}_+^d.
				\end{aligned}
			\end{equation*}
			So we have proved the case of $k=1$. For $k>1$ By induction, assume \eqref{mcbdbdyfr} holds for any $l\leq k$. Taking $\partial_d^{2k}$ on both sides of \eqref{mcbdbdeq}, we have
			\begin{equation*}
				\begin{aligned}
					\partial_d^{2k+2}f&=\partial_t\partial_d^{2k}f-\Delta_{x'}\partial_d^{2k}f-\partial_d^{2k}F\\
					&=-\Delta_{x'}\partial_d^{2k}f-\partial_d^{2k}F+\partial_t(F^{k-1}+R_{2k})\\
					&=-\Delta_{x'}\partial_d^{2k}f-\partial_d^{2k}F+\sum_{2i+|\alpha|=2k-2}C^{1,k-1}_{i,\alpha}\partial_t^{i+1}\nabla^{\alpha}F+\sum_{\substack{l+|\beta|\leq 2k\\l<2k}}C_{l,\beta}^{2,k}\partial_{d}^{l}\nabla_{x'}^{\beta}(\Delta f+F)\\
					&=\sum_{2i+|\alpha|=2k}C^{1,k+1}_{i,\alpha}\partial_t^{i}\nabla^{\alpha}F+\sum_{\substack{l+|\beta|\leq 2k+2\\l<2k+2}}C_{l,\beta}^{2,k+2}\partial_{d}^{l}\nabla_{x'}^{\beta}f,\quad \text{on}\ (0,T]\times\partial\mathbb{R}_+^d,
				\end{aligned}
			\end{equation*}
			which implies that \eqref{mcbdbdyfr} holds for $k+1$, so by induction we have proved \eqref{mcbdbdyfr} for any $k\in\mathbb{N}$. The estimate \eqref{mcbdbdyfrm} follows directly by taking $\nabla^\beta_{x'}$ to the equation \eqref{mcbdlembdy}. This completes the proof.
		\end{proof}\vspace{0.3cm}\\
		Recalling the definition of $\|\cdot\|_{Z_T^0,2}$ and $\|\cdot\|_{Z_T,2}$ in \eqref{defnormZ}, the following is a key lemma for the boundary estimates.
		\begin{lemma}\label{mcbdlemby}
			Consider $F^1,F^2\in C^\infty([0,T]\times \mathbb{R}_+^d)$. If $u$ is a solution to the system 
			\begin{equation}\label{mcbd}
				\begin{aligned}
					&\partial_tu(t,x)-\Delta u(t,x)=F^1(t,x)+\partial_dF^2(t,x),\quad \text{in}\ [0,T]\times \mathbb{R}_+^d,\\
					&u(0,x)=u_0(x),\quad \text{in}\ \mathbb{R}_+^d,\\
					&u(t,x)=0,\quad \text{on}\ [0,T]\times\partial\mathbb{R}_+^d,
				\end{aligned}
			\end{equation}
			then we have
			\begin{equation}\label{mcbddlow}
				\|u\|_{Z_T^0,2}\lesssim \|u_0\|_{W^{1,\infty}}+\mathbf{C}_{0}(F^1,F^2)(T).
			\end{equation}
			and for higher order estimates, we have
			\begin{equation}\label{mcbddmrs}
				\begin{aligned}
					&\|u\|_{Z_T,2}\lesssim   \|u\|_{Z_T^0,2}+\mathbf{C}_{m}(F^1,F^2)(T),
				\end{aligned}
			\end{equation}
			where we denote 
			$$
			\mathbf{C}_{n}(F^1,F^2)(T)=\sup_{t\in[0,T]}\sum_{i+j\leq n}\left(t^{i+j+\frac{1+\kappa}{2}}\|\partial_t^i\nabla^{2j}(F^{1},\nabla_{x'}F^2)(t)\|_{C^\kappa}+t^{i+j+\frac{\kappa}{2}}\|\partial_t^i\nabla^{2j}F^2(t)\|_{C^\kappa}\right).
			$$
		\end{lemma}
		\begin{proof}
			We first prove  \eqref{mcbddlow}.  With the Green function $\H$ defined in \eqref{defH}, we transform \eqref{mcbd} to the following integral equation:
			\begin{equation*}
				\begin{aligned}
					u(t,x)&=\int_{\mathbb{R}_+^d} \H(t,x'-y',x_d,y_d)u_{0}(y)dy+\int_0^t\int_{\mathbb{R}^d_+} \H(t-\tau,x'-y',x_d,y_d) (F^1+\partial_dF^2)(\tau,y)dyd\tau\\
					&:=u_{L}(t,x)+\sum_{k=1}^2u^k_{N}(t,x).
				\end{aligned}
			\end{equation*}
			By \eqref{llha} we can prove that
			\begin{equation}\label{ulhalf}
				\begin{aligned}
					&\|\nabla u_{L}(t)\|_{L^\infty(\mathbb{R}^d_+)}\lesssim \|\nabla u_{0}\|_{L^\infty(\mathbb{R}^d_+)},
					\\
					&\| \delta_\alpha\nabla_{x'}\nabla u_{L}(t)\|_{L^\infty(\mathbb{R}^d_+)}\lesssim \int_{\mathbb{R}^d_+} |\delta_\alpha\nabla_{x}\H(t,x'-y',x_d,y_d)|dy\|\nabla_{x'} u_{0}\|_{L^\infty(\mathbb{R}^d_+)}\\
					&\quad\quad\quad\quad\quad\quad\quad\quad\quad\quad\quad\quad\lesssim \min\left\{1,\frac{|\alpha|}{t^{\frac{1}{2}}}\right\}t^{-\frac{1}{2}}\|\nabla u_0\|_{L^\infty(\mathbb{R}^d_+)}.
				\end{aligned}
			\end{equation}
			For $u_N^1$, using integration by parts to get
			\begin{equation}\label{un1half}
				\begin{aligned}
					&\|\nabla u_{N}^1(t)\|_{L^\infty(\mathbb{R}^d_+)}\lesssim \int_0^t\int_{\mathbb{R}^d_+} |\nabla_x\H(t-\tau,x'-y',x_d,y_d)||F^1(\tau,y)-F^1(\tau,x)|dyd\tau\\
					&\quad\quad\quad\quad\quad\quad\quad\quad\quad\quad\quad\quad\lesssim \sup_{t\in[0,T]}t^{\frac{1+\kappa}{2}}\|F^1(t)\|_{\dot C^\kappa},\\
					&\| \delta_\alpha\nabla_{x'}\nabla u_{N}^1(t)\|_{L^\infty(\mathbb{R}^d_+)}\lesssim \int_0^t\int_{\mathbb{R}^d_+} |\delta_\alpha\nabla_x\nabla_{x'}\H(t-\tau,x'-y',x_d,y_d)||F^1(\tau,y)-F^1(\tau,x)|dyd\tau\\
					&\quad\quad\quad\quad\quad\quad\quad\quad\quad\quad\quad\quad\lesssim \min\left\{1,\frac{|\alpha|}{t^{\frac{1}{2}}}\right\}t^{-\frac{1}{2}}\sup_{t\in[0,T]}t^{\frac{1+\kappa}{2}}\|F^1(t)\|_{\dot C^\kappa}.
				\end{aligned}
			\end{equation}
			For $u_N^2$, with the boundary condition $\partial_d\nabla_{x'}\H(t,x'-y',x_d,y_d)|_{y_d=0}=0$, using integration by parts  we can obtain that
			\begin{equation}\label{un2half}
				\begin{aligned}
					&\|\nabla u_{N}^2(t)\|_{L^\infty(\mathbb{R}^d_+)}\lesssim \int_0^t\int_{\mathbb{R}^d_+} |\nabla\partial_{x_d}\tilde{\H}(t-\tau,x'-y',x_d,y_d)||F^2(\tau,y)-F^2(\tau,x)|dyd\tau\\
					&\quad\quad\quad\quad\quad\quad\quad\quad\quad\quad\quad\quad\lesssim \sup_{t\in[0,T]}t^{\frac{\kappa}{2}}\|F^2(t)\|_{\dot C^\kappa},\\
					&\| \delta_\alpha\nabla_{x'}\nabla u_{N}^2(t)\|_{L^\infty(\mathbb{R}^d_+)}\lesssim \int_0^t\int_{\mathbb{R}^d_+} 
					|\delta_\alpha\nabla_x\partial_d\tilde{\H}(t-\tau,x'-y',x_d,y_d)||\nabla_{y'}F^2(\tau,y)-\nabla_{y'}F^2(\tau,x)|dyd\tau\\
					&\quad\quad\quad\quad\quad\quad\quad\quad\quad\quad\quad\quad\lesssim \min\left\{1,\frac{|\alpha|}{t^{\frac{1}{2}}}\right\}t^{-\frac{1}{2}}\sup_{t\in[0,T]}t^{\frac{1+\kappa}{2}}\|\nabla_{x'}F^2(t)\|_{\dot C^\kappa},
				\end{aligned}
			\end{equation}
			where $\tilde \H$ is defined in \eqref{deftilH}.
			By \eqref{ulhalf}, \eqref{un1half} and \eqref{un2half}, we obtain \eqref{mcbddlow}. 
			
			Moreover, following the proof of  \eqref{uzha}, we obtain 
			from  \eqref{mcbddlow} that for any $\gamma\geq 0$,
			\begin{equation}\label{mcbddmlemh}
				\begin{aligned}
					&\sup_{t\in(0,T]}t^{\gamma+\frac{1+\kappa}{2}}\|\nabla_{x'}\nabla u(t)\|_{C^\kappa}\\
					&\quad\lesssim \sup_{t\in(0,T]}t^{\gamma}\|\nabla u(t)\|_{L^\infty}+ \sup_{t\in[0,T]}t^{\gamma+\frac{1+\kappa}{2}}\|(F^1,\nabla_{x'}F^2)(t)\|_{C^\kappa}+ \sup_{t\in[0,T]}t^{\gamma+\frac{\kappa}{2}}\|F^2(t)\|_{C^\kappa}.
				\end{aligned}
			\end{equation}
			Now we prove \eqref{mcbddmrs}.
			It suffices to consider highest order spatial derivatives, since the time derivatives can be transformed into spatial derivatives by the equation of $u$. Hence, we only prove that 
			\begin{equation}\label{highes}
				\sup_{t\in[0,T]}\sum_{n+2k= 2m}t^{\frac{2m+1+\kappa}{2}}\|\nabla_{x'}\nabla(\partial_{d}^{2k}\nabla_{x'}^{n}u)(t)\|_{C^\kappa}\lesssim \|u\|_{Z_T^0,2}+\mathbf{C}_{m}(F^1,F^2)(T),
			\end{equation}
			For the case $k=0$, we consider the equation of  $\nabla_{x'}^{2m}u$, then the desired estimate follows \eqref{mcbddmlemh} by taking $\gamma=m$. It remains to consider $1\leq k\leq m$. First we denote $u_{n}=\nabla_{x'}^nu$, then we can write the equation of $u_n$ as
			\begin{equation}\label{equn}
				\begin{aligned}
					&\partial_tu_n(t,x)-\Delta u_n(t,x)=F_n^1(t,x)+\partial_dF_n^2(t,x),\quad \text{in}\ (0,T]\times\mathbb{R}^d_+,\\
					&u_n(t,x)=0,\ \text{on}\ (0,T]\times\partial \mathbb{R}^d_+,
				\end{aligned}
			\end{equation}
			with
			\begin{equation*}
				F_n^i=\nabla_{x'}^nF^i.
			\end{equation*}
			By Lemma \ref{mcbdlembdy},  we have the boundary condition
			\begin{equation*}
				\partial_{d}^{2k}u_n=F_{n}^{1,k-1}+\partial_dF_n^{2,k-1}+R^{2k}_n,\quad \text{on}\ \partial\mathbb{R}^d_+,
			\end{equation*}
			where 
			\begin{align*}
				F_{n}^{i,k-1}=\sum_{2j+|\alpha|=2k-2}C_{j,\alpha}^{1,n}\partial_t^j\nabla^\alpha F_n^i,\ \ i=1,2,\ \ R^{2k}_n=\sum_{\substack{l+|\beta|= 2k\\l<2k}}C_{l,\beta}^{2,n}\partial_{d}^{l}\nabla_{x'}^{\beta}u_n.
			\end{align*}
			Furthermore, denote $u_{n,2k}=\partial_{d}^{2k}u_n-(F_{n}^{1,k-1}+\partial_dF_{n}^{2,k-1}+R_{2k}^n)$, then we will get the equation of $u_{n,2k}$ as
			\begin{equation}\label{un2k}
				\begin{aligned}
					&\partial_tu_{n,2k}(t,x)-\Delta u_{n,2k}(t,x)=\tilde{F}_{n,2k}^1(t,x)+\partial_d\tilde{F}_{n,2k}^2(t,x)+\tilde{R}_n^{2k},\quad \text{in}\ (0,T]\times\mathbb{R}^d_+,\\
					&u_{n,2k}(t,x)=0,\ \text{on}\ (0,T]\times\partial \mathbb{R}^d_+,
				\end{aligned}
			\end{equation}
			with
			\begin{equation*}
				\tilde{F}_{n,2k}^{i}=(\partial_t-\Delta)F_n^{i,k-1}+\partial_d^{2k}F_n^i,\quad\tilde{R}^{2k}_n=(\partial_t-\Delta)R_n^{2k}.
			\end{equation*}
			By the form of $R_n^{2k}$ and equation of $u_n$ in \eqref{equn}, we have 
			\begin{equation*}
				\begin{aligned}
					\tilde{R}^{2k}_n&=  (\partial_t-\Delta) R^{2k}_n=\sum_{\substack{l+|\beta|= 2k\\l<2k}}C_{l,\beta}^{2,n}\partial_{d}^{l}\nabla_{x'}^{\beta}(\partial_t-\Delta)u_n\\
					&=\sum_{\substack{l+|\beta|= 2k\\l<2k}}C_{l,\beta}^{2,n}\partial_{d}^{l}\nabla_{x'}^{\beta}F^1_{n}+\sum_{\substack{l+|\beta|= 2k\\l<2k}}C_{l,\beta}^{2,n}\partial_{d}^{l+1}\nabla_{x'}^{\beta}F^2_{n}=R^{1,2k}_n+\partial_dR_{n}^{2,2k}.
				\end{aligned}
			\end{equation*}
			Applying the estimate \eqref{mcbddmlemh} to the equation \eqref{un2k} with $\gamma=\frac{n+2k}{2}$, we will get
			\begin{equation*}
				\begin{aligned}
					&\sup_{t\in[0,T]}t^{\frac{n+2k+1+\kappa}{2}}\|\nabla_{x'}\nabla u_{n,2k}(t)\|_{C^\kappa}\\
					&\quad\quad\lesssim  \sup_{t\in[0,T]}t^{\frac{n+2k}{2}}\|\nabla u_{n,2k}(t)\|_{L^\infty}+\sup_{t\in[0,T]}t^{\frac{n+2k+\kappa}{2}}\left(\|\tilde{F}_{n,2k}^{2}(t)\|_{C^\kappa}+\|R_n^{2,2k}(t)\|_{C^\kappa}\right)\\
					&\quad\quad\quad\quad+\sup_{t\in(0,T]}t^{\frac{n+2k+1+\kappa}{2}}\left(\|(\tilde{F}_{n,2k}^{1},\nabla_{x'}\tilde{F}_{n,2k}^{2})(t)\|_{C^\kappa}+\|(R_n^{1,2k},\nabla_{x'}R_n^{2,2k})(t)\|_{C^\kappa}\right).
				\end{aligned}
			\end{equation*}
			By the definition of $\tilde{F}_{n,2k}^{1},\tilde{F}_{n,2k}^{2},R_n^{1,2k},R_n^{2,2k}$, one has 
			\begin{align*}
				&\sup_{t\in[0,T]}t^{\frac{n+2k+\kappa}{2}}\left(\|\tilde{F}_{n,2k}^{2}(t)\|_{C^\kappa}+\|R_n^{2,2k}(t)\|_{C^\kappa}\right)\\
				&\ \ \ +\sup_{t\in[0,T]}t^{\frac{n+2k+1+\kappa}{2}}\left(\|(\tilde{F}_{n,2k}^{1},\nabla_{x'}\tilde{F}_{n,2k}^{2})(t)\|_{C^\kappa}+\|(R_n^{1,2k},\nabla_{x'}R_n^{2,2k})(t)\|_{C^\kappa}\right)\lesssim \mathbf{C}_{m}(F^1,F^2)(T).
			\end{align*}
			On the other hand, we have the interpolation inequality
			\begin{align*}
				\sup_{t\in[0,T]}t^{\frac{n+2k}{2}}\|\nabla u_{n,2k}(t)\|_{L^\infty}&\lesssim C(\eps)\sup_{t\in[0,T]}\|\nabla u(t)\|_{L^\infty}+\eps\sup_{t\in[0,T]}t^{\frac{n+2k+1+\kappa}{2}}\|\nabla_{x'}\nabla u_{n,2k}(t)\|_{C^\kappa}\\
				&\lesssim C(\eps)\|u\|_{Z^0_T,2}+\eps\sup_{t\in[0,T]}t^{\frac{n+2k+1+\kappa}{2}}\|\nabla_{x'}\nabla u_{n,2k}(t)\|_{C^\kappa},
			\end{align*}
			holds for any $\eps\in(0,1)$. Then we obtain \eqref{highes} by taking $\eps$ small enough such that the remainder term in the right hand side can be absorbed by the left hand side. This completes the proof.
		\end{proof}\vspace{0.3cm}\\
		\begin{proof}[Proof of Theorem \ref{mcbdglo}, i)] We will construct a solution to the system \eqref{eqmcbdd} by a fixed point argument. Let $g\in\mathcal{Z}_{T,\phi}^\sigma$ with $\phi=e^{\eps_1\Delta}f_0$, and  \begin{align}\label{inismal}
				\|f_0-\phi\|_{W^{1,\infty}}\leq\eps_0.
			\end{align} 
			We define a map $\mathcal{S}$ by $\mathcal{S}g=f$, where $f$ is the weak solution to the following system
			\begin{equation}\label{eqmcbd}
				\begin{aligned}
					&\partial_t f=\A[\nabla g]:\nabla^2 f,\ \ \ \text{in}\ [0,T]\times \Omega,\\
					&f|_{t=0}=f_0,\ \ \ \text{in}\ \Omega,\\
					&f=0,\ \ \text{on}\ [0,T]\times\partial \Omega.
				\end{aligned}
			\end{equation}
			We want to show that $\mathcal{S}$ is a contraction map in $\mathcal{Z}_{T,\phi}^\sigma$ for some suitable $T$ and $\sigma$. \\
			First we want to show that if $g\in\mathcal{Z}_{T,\phi}^\sigma$, then $f=\mathcal{S}g\in \mathcal{Z}_{T,\phi}^\sigma$. 
			Denote $f_{in}=\chi_0 (f-\phi)$. Then one has 
			\begin{equation}\label{eqinter}
				\begin{aligned}
					&\partial_tf_{in}-\A[\nabla \phi]:\nabla^2f_{in}=F^0,\quad \text{in}\ [0,T]\times\mathbb{R}^d,\\
					&f_{in}|_{t=0}=\chi_0(f_0-\phi),\quad \text{in} \ \mathbb{R}^d,
				\end{aligned}
			\end{equation}
			with
			\begin{equation*}
				F^0 =-2\A[\nabla g]:(\nabla f\otimes\nabla\chi_0)-\A[\nabla g]:(f\nabla^2\chi_0)+\A[\nabla g]:\nabla^2\phi+(\A[\nabla g]-\A[\nabla\phi]):\nabla^2 f_{in}.
			\end{equation*}
			By Theorem \ref{thmexis}, since $g\in\mathcal{Z}_{T,\phi}$ with small $\sigma$, there exists $T=T(\|\phi\|_{C^{n+5}})>0$ such that 
			\begin{equation}\label{ee1ho}
				\begin{aligned}
					& \|f_{in}\|_{Z_T^0,1}\lesssim \|\chi_0(f_0-\phi)\|_{ W^{1,\infty}}+\sup_{t\in[0,T]}t^\frac{1+\kappa}{2}\|F^0(t)\|_{ C^\kappa},\\
					& \|f_{in}\|_{Z_T,1}\lesssim \|\chi_0(f_0-\phi)\|_{ W^{1,\infty}}+\sup_{t\in[0,T]}\left(t^\frac{1+\kappa}{2}\|F^0(t)\|_{ C^\kappa}+t^\frac{n+1+\kappa}{2}\|F^0(t)\|_{ C^{n+\kappa}}\right).
				\end{aligned}
			\end{equation}
			For data term $\|\chi_0(f_0-\phi)\|_{W^{1,\infty}}$, note that since $(f_0-\phi)|_{\partial\Omega}=0$ and our definition of $\phi$ as in Theorem \ref{mcbdglo}, we take $\chi_0$ such that the size of $\operatorname{Supp}(\chi_0)\sim\eps_1^{\frac{1}{2}}$, we can see that $|\chi_0(f_0-\phi)|\lesssim \eps_1^{\frac{1}{2}} \|f_0-\phi\|_{W^{1,\infty}}$. Together with the scale of $|\nabla\chi_0|\sim\eps_1^{-\frac{1}{2}}$, we will have
			\begin{equation}\label{mcbdbdysc}
				\|\chi_0(f_0-\phi)\|_{W^{1,\infty}}\leq C'\|f_0-\phi\|_{W^{1,\infty}}.
			\end{equation}
			Then by interpolation we obtain 
			\begin{equation}\label{bddmcior}
				\begin{aligned}
					&\|f_{in}\|_{Z_T^0,1}\lesssim C_\Omega\left(\|  f_0-\phi\|_{W^{1,\infty}}+(\|g-\phi\|_{Z_T}+T^{\frac{\kappa}{2}}\|\phi\|_{C^{5}})(\|f_{in}\|_{Z_T^0,1}+T^{\frac{1}{2}}\|\phi\|_{C^{5}})\right.\\
					&\quad\left.+(1+T)^{4}(\|g-\phi\|_{Z_T}+T^{\frac{\kappa}{2}}\|\phi\|_{C^{5}})(1+\|g-\phi\|_{Z_T}+\|\phi\|_{C^{5}})^{5}\sup_{t\in(0,T]}t^{\frac{1}{2}}\|f_{in}(t)\|_{C^{2}}\right).
				\end{aligned}
			\end{equation}
			Furthermore, we will have
			\begin{equation*}
				\begin{aligned}
					\|f_{in}\|_{Z_T^0,1}\leq C\|f_0-\phi\|_{W^{1,\infty}}+C(\sigma+T^{\frac{\kappa}{2}}\|\phi\|_{C^5})\|f_{in}\|_{Z_T^0,1}+CT^{\frac{\kappa}{2}}(1+\sigma)\|\phi\|_{C^5}^2,
				\end{aligned}
			\end{equation*}
			for some $T<1$, and we can take
			\begin{equation*}
				\sigma\leq\frac{1}{500C},\quad T\leq (\frac{\sigma}{500C\|\phi\|_{C^5}})^{\frac{2}{\kappa}},
			\end{equation*}
			then we can get
			\begin{equation}\label{bdmczt01}
				\|f_{in}\|_{Z_T^0,1}\leq C_1\eps_0.
			\end{equation}
			Furthermore, by same methods, we will have the higher order estimates
			\begin{equation}\label{bddmciorh}
				\begin{aligned}
					&\|f_{in}\|_{Z_T,1}\lesssim C_\Omega\left(\sup_{t\in(0,T]}t^{m}\|\nabla^{2m+1}f_{in}(t)\|_{L^\infty}+(\|g-\phi\|_{Z_T}+T^{\frac{\kappa}{2}}\|\phi\|_{C^{5}})(\|f_{in}\|_{Z_T,1}+(1+T)^{m+4}T^{\frac{1}{2}}\|\phi\|_{C^{5}})\right.\\
					&\quad+(1+T)^{m+4}(\|g-\phi\|_{Z_T}+T^{\frac{\kappa}{2}}\|\phi\|_{C^{5}})(1+\|g-\phi\|_{Z_T}+\|\phi\|_{C^{5}})^{m+5}\\
					&\quad\quad\quad\quad\quad\times\left.\left(\sup_{t\in(0,T]}t^{m+\frac{1}{2}}\|\nabla^{2m+1}f_{in}(t)\|_{L^\infty}+T^{\frac{1}{2}}\|\phi\|_{C^{2m+5}}\right)\right).
				\end{aligned}
			\end{equation}
			Then by interpolation, \eqref{bddmciorh} and \eqref{bdmczt01}, we can also take $\|f_0-\phi\|_{C^1}$, $\sigma$, $T$ small enough to get
			\begin{equation}\label{bdmczt1}
				\|f_{in}\|_{Z_T,1}\leq\frac{\sigma}{10}.
			\end{equation}
			For boundary estimate, we multiply \eqref{mcbd} by $\chi_j$ (with a little abuse of notation, we will omit $j$ in the following proof), with the rotation defined as before such that $0\in\operatorname{Supp}(\chi)$ and $\nabla\varphi(0)=0$, then
			\begin{equation}\label{eqbd}
				\begin{aligned}
					&\partial_t(\chi (f-\phi))-\A[\nabla g]:\nabla^2(\chi (f-\phi))=G,\quad \text{in}\ [0,T]\times\Omega,\\
					&\chi (f-\phi)|_{t=0}=\chi (f_0-\phi),\quad\text{in}\ \Omega,\\
					&\chi (f-\phi)(t,x)=0,\quad\text{in}\ [0,T]\times\partial \Omega,
				\end{aligned}
			\end{equation}
			with
			\begin{equation*}
				G=-2\A[\nabla g]:(\nabla\chi\otimes\nabla f)-\A[\nabla g]:(f\nabla^2\chi)+\A[\nabla g]:\nabla^2\phi.
			\end{equation*}
			Assume $\partial\Omega\cap\mathrm{Supp}(\tilde\chi)$ can be written as $x_d=\varphi(x')$, and denote $\Phi$ the flatten map, which is $\Phi(x)=(x',x_d-\varphi(x'))$, and denote the inverse map $\Psi=\Phi^{-1}$. We notice that after this transform, the coefficient matrix $\A[\nabla g]$ will be replaced by $\tilde{\A}[\nabla((\tilde{\chi}g)\circ\Phi)]$, with
			\begin{equation}\label{mcbdnewc}
				\begin{aligned}
					&\left(\A[\nabla (\tilde{\chi}g)]:\nabla^2f\right)\circ\Phi=\tilde{\A}[\nabla ((\tilde{\chi}g)\circ\Phi)]:\nabla^2(f\circ\Phi)+R,\\
					&\tilde{\A}[\nabla ((\tilde{\chi}g)\circ\Phi)]=(\nabla\Psi)^T\A[\nabla\Psi\nabla((\tilde{\chi}g)\circ\Phi)]\nabla\Psi,\\
					&R=\A_{ij}[\nabla\Psi\nabla((\tilde{\chi}g)\circ\Phi)]\partial_{ij}\Psi^k(\partial_kf)\circ\Phi.
				\end{aligned}
			\end{equation}
			Then we transform \eqref{eqbd} to the following system on half space, and freeze coefficient at $0$ and apply the rescale $\Theta$, such that $\Theta(x)=(x',(1+|\nabla\phi(0)|^2)^{\frac{1}{2}}x_d)$, and consider the equation of $\tilde{f}(t,x)=\big((\chi (f-\phi)(t,\cdot))\circ\Phi\circ\Theta\big)(x)$ as
			\begin{equation}
				\begin{aligned}
					&\partial_t\tilde{f}(t,x)-\Delta\tilde{f} (t,x)=\tilde{F}^1(t,x)+\partial_d\tilde{F}^2(t,x),\quad \text{in}\ [0,T]\times\mathbb{R}^d_+,\\
					&\tilde{f}|_{t=0}=\tilde{f}_0=(\chi (f_0-\phi))\circ\Phi\circ\Theta,\quad \text{in}\ \mathbb{R}^d_+,\\
					&\tilde{f}(t,x)=0,\ \text{on}\ [0,T]\times\partial \mathbb{R}^d_+,
				\end{aligned}
			\end{equation}
			with
			\begin{equation*}
				\begin{aligned}
					&\tilde{F}^1=\left(G\circ\Phi+R+\sum_{i\neq d\ or\ j\neq d}(\tilde{\A}[\nabla\tilde{g}]-\tilde{\A}_0)_{ij}\partial_{ij}\tilde{f}+\left(\partial_d(\tilde{\A}[\nabla\tilde{g}]-\tilde{\A}_0)_{dd}\right)\partial_d\tilde{f}\right)\circ\Theta,\\
					&\tilde{F}^2=\left((\tilde{\A}[\nabla\tilde{g}]-\tilde{\A}_0)_{dd}\partial_df\right)\circ\Theta,
				\end{aligned}
			\end{equation*}
			where $\tilde{\A}[\nabla\tilde{g}]=\tilde{\A}[\nabla((\tilde{\chi}g)\circ\Phi)]=\tilde{\A}[\nabla \tilde{g}]$ is the new coefficient defined in \eqref{mcbdnewc}, $\tilde{\A}_0=\tilde{\A}[\nabla((\tilde{\chi}\phi)\circ\Phi)](0)=\tilde{\A}[\nabla\tilde{\phi}](0)$. We emphasize that we will take $\tilde\chi$ and $\chi$ such that $|\tilde{\A}[\tilde{\phi}](x)-\tilde{\A}_0|\leq\eps$ with $\phi=e^{\eps_1\Delta}f_0$, which gives the scale of $|\nabla\chi|$ is nearly $\eps_1^{-\frac{1}{2}}$. Then we apply \eqref{mcbddlow} to get 
			\begin{equation}\label{mcbddzt02}
				\begin{aligned}
					&\|\tilde{f}\|_{Z_T^0,2}\lesssim\|\tilde{f}_0\|_{W^{1,\infty}}+\mathbf{C}_0(\tilde{F}^1,\tilde{F}^2)(T)\\
					&\quad\lesssim \|\tilde{f}_0\|_{W^{1,\infty}}+(\eps+\|g-\phi\|_{Z_T}+T^{\frac{\kappa}{2}}\|\phi\|_{C^{5}})(\|\tilde{f}\|_{Z_T^0,2}+T^{\frac{\kappa}{2}}\|\phi\|_{C^5})\\
					&\quad\quad+(1+T)^4(\|g-\phi\|_{Z_T}+T^{\frac{\kappa}{2}}\|\phi\|_{C^{5}})(1+\|g-\phi\|_{Z_T}+\|\phi\|_{C^{5}})\sup_{t\in(0,T]}t^{\frac{1}{2}}\|\nabla_{x'}\tilde{f}(t)\|_{C^1}.
				\end{aligned}
			\end{equation}
			By interpolation, and since $\Phi$, $\Theta$ are uniformly bounded maps, we have
			\begin{equation*}
				\|\tilde{f}\|_{Z_T^0,2}\leq C\|\tilde{f}_0\|_{W^{1,\infty}}+C(\eps+\|g-\phi\|_{Z_T}+T^{\frac{\kappa}{2}}\|\phi\|_{C^{5}})(\|\tilde{f}\|_{Z_T^0,2}+T^{\frac{\kappa}{2}}\|\phi\|_{C^5}),
			\end{equation*}
			for some $T<1$. For data term $\|\tilde{f}_0\|_{W^{1,\infty}}$, similar with \eqref{mcbdbdysc}, we will have
			\begin{equation*}
				\|\tilde{f}_0\|_{W^{1,\infty}}\leq C'\|f_0-\phi\|_{W^{1,\infty}}.
			\end{equation*}
			We take
			\begin{equation*}
				\sigma\leq\frac{1}{500C},\quad\eps\leq\frac{1}{500C},\quad \|f_0-\phi\|_{W^{1,\infty}}\leq\frac{\sigma}{500CC'},\quad T\leq(\frac{\sigma}{500C\|\phi\|_{C^5}})^{\frac{2}{\kappa}},
			\end{equation*}
			then we can prove that
			\begin{equation}\label{mcbddestz02}
				\|\tilde{f}\|_{Z_T^0,2}\leq\frac{\sigma}{100}.
			\end{equation}
			Now we only need to do estimates for highest order spacial derivatives, since the time derivatives can be transformed into spacial derivatives by the equation of $f$. In general, we will consider
			\begin{equation*}
				\sup_{t\in(0,T]}\sum_{n+2k\leq 2m}t^{n+2k+\frac{1+\kappa}{2}}\|\partial_{d}^{2k+1}\nabla_{x'}^{n+1}\tilde{f}(t)\|_{C^\kappa}.
			\end{equation*}
			By Lemma \ref{mcbdlemby} and interpolation, and the fact that $|\tilde{\A}[\nabla\tilde{\phi}](t,x)-\tilde{\A}_0(t,x)|\lesssim\eps$ when $x\in\operatorname{Supp}(\chi)$, we will have the estimate for force term $\tilde{F}$ that
			\begin{equation*}
				\begin{aligned}
					&\mathbf{C}_m(\tilde{F}^1,\tilde{F}^2)=\sum_{i+j=m}\left(t^{\frac{2m+1+\kappa}{2}}\|\partial_t^i\nabla^{2j}(\tilde{F}^1,\nabla_{x'}\tilde{F}^2)(t)\|_{C^\kappa}+t^{\frac{2m+\kappa}{2}}\|\partial_t^i\nabla^{2j}\tilde{F}^2(t)\|_{C^\kappa}\right)\\
					&\leq C_2(\eps+\|g-\phi\|_{Z_T}+T^{\frac{\kappa}{2}}\|\phi\|_{C^{2m+5}})(\|f-\phi\|_{Z_T}+T^{\frac{\kappa}{2}}\|\phi\|_{C^{2m+5}})\\
					&+(1+T)^{n+4}(\|g-\phi\|_{Z_T}+T^{\frac{\kappa}{2}}\|\phi\|_{C^{2m+5}})(1+\|g-\phi\|_{Z_T}+\|\phi\|_{C^{2m+5}})^{2m+5}\\
					&\quad\quad\quad\quad\quad\times\sup_{t\in(0,T]}\left(\|\nabla f(t)\|_{L^\infty}+\sum_{n+2k\leq 2m}t^{\frac{n+2k+1}{2}}(\|\nabla_{x'}^{n+1}\partial_d^{2k+1}f(t)\|_{L^\infty}+\|\nabla_{x'}^{n+2}\partial_d^{2k}f(t)\|_{L^\infty})\right).
				\end{aligned}
			\end{equation*}
			So we combine the estimate of $\tilde{F}$ and \eqref{mcbddmrs} to get
			\begin{equation}\label{bdmcestbd}
				\begin{aligned}
					&\|\tilde{f}\|_{Z_T,2}\leq C_1\left(\|\tilde{f}\|_{Z_T^0,2}+\mathbf{C}_m(\tilde{F}^1,\tilde{F}^2)(T)\right)\\
					&\quad\leq C_1\|\tilde{f}\|_{Z_T^0,2}+C_2(\eps+\|g-\phi\|_{Z_T}+T^{\frac{\kappa}{2}}\|\phi\|_{C^{2m+5}})(\|f-\phi\|_{Z_T}+T^{\frac{\kappa}{2}}\|\phi\|_{C^{2m+5}})\\
					&\quad\quad+(1+T)^{m+4}(\|g-\phi\|_{Z_T}+T^{\frac{\kappa}{2}}\|\phi\|_{C^{2m+5}})(1+\|g-\phi\|_{Z_T}+\|\phi\|_{C^{2m+5}})^{2m+5}\\
					&\quad\quad\quad\quad\quad\times\sup_{t\in(0,T]}\left(\|\nabla f(t)\|_{L^\infty}+\sum_{n+2k\leq 2m}t^{\frac{n+2k+1}{2}}(\|\nabla_{x'}^{n+1}\partial_d^{2k+1}f(t)\|_{L^\infty}+\|\nabla_{x'}^{n+2}\partial_d^{2k}f(t)\|_{L^\infty})\right)\\
					&\leq \frac{1}{2}\|\tilde{f}\|_{Z_T,2}+ C_2\|\tilde{f}\|_{Z_T^0,2}\\
					&\quad\quad+C_3(\eps+\|g-\phi\|_{Z_T}+T^{\frac{\kappa}{2}}\|\phi\|_{C^{2m+5}})(1+\|g-\phi\|_{Z_T}+\|\phi\|_{C^{2m+5}})^{2m+5}(\|f-\phi\|_{Z_T}+T^{\frac{\kappa}{2}}\|\phi\|_{C^{2m+5}}),
				\end{aligned}
			\end{equation}
			where the last inequality we use interpolation. Combine \eqref{bdmcestbd} with $\eqref{mcbddestz02}$, and take $\eps_0$, $\eps$, $\sigma$ and $T$ small enough, we can get 
			\begin{equation*}
				\|\tilde{f}\|_{Z_T,2}\leq \frac{\sigma}{100}.
			\end{equation*}
			Combine the result with \eqref{bdmczt1}, and we can prove that
			\begin{equation}\label{ztpriest}
				\|f-\phi\|_{Z_T}\leq \frac{\sigma}{2}.
			\end{equation}  
			Now we prove the contraction property. Let $g_1,g_2\in\mathcal{Z}_{T,\phi}^\sigma$, and $f_1=\mathcal{S}g_1$, $f_2=\mathcal{S}g_2$, and denote $\mathbf{f}=f_1-f_2$, $\mathbf{g}=g_1-g_2$, $\mathbf{A}=\A[\nabla g_1]-\A[\nabla g_2]$, then we can write the equation of $\mathbf{f}$ as
			\begin{equation*}
				\begin{aligned}
					&\partial_t\mathbf{f}-\A[\nabla g_1]:\nabla^2\mathbf{f}=\mathbf{A}:\nabla^2f_2,\ \ \ \text{in}\ [0,T]\times \Omega,\\
					&\mathbf{f}|_{t=0}=0,\ \ \ \text{in}\ \Omega,\\
					&\mathbf{f}=0,\ \ \text{on}\ [0,T]\times\partial \Omega.
				\end{aligned}
			\end{equation*}
			For interior estimates we can just denote $\mathbf{f}_{in}=\chi_0\mathbf{f}$ and we can write the equation of $\mathbf{f}_{in}$ as
			\begin{equation*}
				\begin{aligned}
					&\partial_t\mathbf{f}_{in}-\A[\nabla g_1]:\nabla^2\mathbf{f}_{in}=\mathbf{A}:\nabla^2f_{2,in}+\chi_0\A[\nabla g_1]:\nabla^2\mathbf{f}-\A[\nabla g_1]:\nabla^2\mathbf{f}_{in},\ \ \ \text{in}\ [0,T]\times \mathbb{R}^d,\\
					&\mathbf{f}_{in}|_{t=0}=0,\ \ \ \text{in}\ \mathbb{R}^d.
				\end{aligned}
			\end{equation*}
			Apply Theorem \ref{lemmain} to obtain the interior estimate
			\begin{equation*}
				\begin{aligned}
					&\|\mathbf{f}_{in}\|_{Z_T,1}\lesssim (\|\mathbf{f}\|_{Z_T}+\|\mathbf{g}\|_{Z_T})(\|f_1-\phi\|_{Z_T}+\|f_2-\phi\|_{Z_T}+\|g_1-\phi\|_{Z_T}+\|g_2-\phi\|_{Z_T}+T^{\frac{\kappa}{2}}\|\phi\|_{C^{2m+5}})\\
					&\quad\quad\times(1+T)^m (1+\|f_1-\phi\|_{Z_T}+\|f_2-\phi\|_{Z_T}+\|g_1-\phi\|_{Z_T}+\|g_2-\phi\|_{Z_T}+T^{\frac{\kappa}{2}}\|\phi\|_{C^{2m+5}})^{2m+4}\\
					&\quad\lesssim (1+T)^{m}(\sigma+T^{\frac{\kappa}{2}}\|\phi\|_{C^{2m+5}})(1+\sigma+T^{\frac{\kappa}{2}}\|\phi\|_{C^{2m+5}})^{2m+4}(\|\mathbf{f}\|_{Z_T}+\|\mathbf{g}\|_{Z_T}).
				\end{aligned}
			\end{equation*}
			We can take $\sigma$ and $T$ small enough to ensure
			\begin{equation*}
				\|\mathbf{f}_{in}\|_{Z_T,1}\leq\frac{1}{100}\|\mathbf{g}\|_{Z_T,1}.
			\end{equation*}
			Now we only need to consider the boundary case. Denote $h_j^i=(\chi_if_j)\circ\Phi_j\circ\Theta_j$, $\mathbf{h}^i=(\chi_i\mathbf{f})\circ\Phi_i\circ\Theta_i$ for $\Phi_i$ and $\Theta_i$ defined as before, such that we can take $\A[\nabla \tilde{\phi}](0)=\operatorname{diag}(1,...1,(1+|\nabla \tilde{\phi}|^2)^{-1})$, and we take suitable $\chi$ such that $|\tilde{\A}[\nabla 
			\tilde{\phi}]-\tilde{\A}_0|\leq \eps$ in $\operatorname{Supp}(\chi)$ (with a little abuse of notation, we will omit $i$ in the following proof). Denote $\tilde{\A}^j=\tilde{\A}[\nabla ((\tilde{\chi}g_j)\circ\Phi)]$ as the new coefficient matrix, and $\mathbf{R}$ corresponding remainder terms as \eqref{mcbdnewc}, and $\mathbf{A}=\tilde{\A}^1-\tilde{\A}^2$. Then we can write the equation of $\mathbf{h}$ as 
			\begin{equation*}
				\begin{aligned}
					&\partial_t\mathbf{h}-\Delta\mathbf{h}=\mathbf{F}^1+\partial_d\mathbf{F}^2,\ \ \ \text{in}\ [0,T]\times \mathbb{R}_+^d,\\
					&\mathbf{h}|_{t=0}=0,\ \ \ \text{in}\ \mathbb{R}_+^d,\\
					&\mathbf{h}=0,\ \ \text{on}\ [0,T]\times\partial \mathbb{R}_+^d,
				\end{aligned}
			\end{equation*}
			with 
			\begin{equation*}
				\begin{aligned}
					&\mathbf{F}^1=\left(\tilde{\A}^1:(\nabla\mathbf{f}\otimes\nabla\chi+\mathbf{f}\nabla^2\chi)\right)\circ\Phi\circ\Theta+\mathbf{R}\circ\Theta+\mathbf{A}:(\nabla h_2\otimes\nabla\psi)+\sum_{i,j\neq d}\mathbf{A}_{ij}\partial_{ij}h_2\\
					&\quad+(\partial_{d}\mathbf{A}_{dd})\partial_dh_2+\sum_{i,j\neq d}(\tilde{\A}-\tilde{\A}_0)_{ij}\partial_{ij}\mathbf{h}+\partial_{d}(\tilde{\A}-\tilde{\A}_0)_{dd}\partial_d\mathbf{h},\\
					&\mathbf{F}^2=\mathbf{A}_{dd}\partial_dh_2+(\tilde{\A}^1-\tilde{\A}^1_0)\partial_d\mathbf{h}.
				\end{aligned}
			\end{equation*}
			Notice that by elementary calculus we have
			\begin{equation*}
				\begin{aligned}
					&\mathbf{C}_m(\mathbf{F}^1,\mathbf{F}^2)(T)=\sum_{i+j=m}\left(t^{\frac{2m+1+\kappa}{2}}\|\partial_t^i\nabla^{2j}(\mathbf{F}^1,\nabla_{x'}\mathbf{F}^2)(t)\|_{C^\kappa}+t^{\frac{2m+\kappa}{2}}\|\partial_t^i\nabla^{2j}\mathbf{F}^2(t)\|_{C^\kappa}\right)\\
					&\leq C'_2(\eps+\|g-\phi\|_{Z_T}+\|f-\phi\|_{Z_T}+T^{\frac{\kappa}{2}}\|\phi\|_{C^{2m+5}})(1+\|g-\phi\|_{Z_T}+\|f-\phi\|_{Z_T}+\|\phi\|_{C^{2m+5}})^{2m+5}\\
					&\quad\quad\quad\times\left(\|\mathbf{f}\|_{Z_T}+\|\mathbf{g}\|_{Z_T}\right).
				\end{aligned}
			\end{equation*}
			So apply \eqref{mcbddmrs} with the estimates of $\mathbf{F}$ above, we have
			\begin{equation*}
				\begin{aligned}
					&\|\mathbf{h}\|_{Z_T,2}\leq C'_1\mathbf{C}_m(\mathbf{F}^1,\mathbf{F}^2)(T)\\
					&\quad\leq C'_2(\eps+\|g-\phi\|_{Z_T}+\|f-\phi\|_{Z_T}+T^{\frac{\kappa}{2}}\|\phi\|_{C^{2m+5}})(1+\|g-\phi\|_{Z_T}+\|f-\phi\|_{Z_T}+T^{\frac{\kappa}{2}}\|\phi\|_{C^{2m+5}})^{2n+5}\\
					&\quad\quad\quad\quad\times(1+T)^m\left(\|\mathbf{f}\|_{Z_T}+\|\mathbf{g}\|_{Z_T}\right).
				\end{aligned}
			\end{equation*}
			Since $f_i,g_i\in\mathcal{Z}_{T,\phi}^\sigma$, by taking $\eps$, $\sigma$ and $T$ small enough, we have proved 
			\begin{equation*}
				\|\mathbf{f}\|_{Z_T}\leq \frac{1}{2} \|\mathbf{g}\|_{Z_T},
			\end{equation*}
			and the whole proof of Theorem \ref{mcbdglo} \textit{i)} is done.
		\end{proof}
		\begin{remark}
			The compatibility condition $u_0|_{\partial \Omega}=0$ is necessary. For a counterexample,
			\begin{equation*}
				\begin{aligned}
					&    \partial_tu-\partial_{x}^2 u=0,\quad \text{in} \ [0,T]\times[0,\infty),\\
					&u(0,x)=1, \quad \text{in}\ [0,T]\times[0,\infty),\\
					&  u(t,0)=0,\quad\forall t\geq 0,
				\end{aligned}
			\end{equation*}
			which is a smooth data. Then by Duhamel formula, we have
			\begin{equation*}
				u(t,x)=\int_{0}^\infty (K(t,x-y)-K(t,x+y))dy=2\int_0^xK(t,z)dz,
			\end{equation*}
			with
			\begin{equation*}
				K(t,x)=\frac{1}{\sqrt{4\pi t}}e^{-\frac{|x|^2}{4t}}.
			\end{equation*}
			So we take derivative to obtain
			\begin{equation*}
				u'(t,x)=2\frac{1}{\sqrt{4\pi t}}e^{-\frac{x^2}{4t}},
			\end{equation*}
			which infers
			\begin{equation*}
				\|u\|_{L_T^\infty(\dot C^1)}\lesssim \|\frac{u(0,x)}{x}\|_{L^\infty}.
			\end{equation*}
			In fact one can see that $\partial_x^{2n}u(t,x)|_{x=0}=0$ and $\partial_x^{2n}u(t,x)\sim cx$ near boundary for any $n\in\mathbb{N}_+$ and $t>0$. But $\|\frac{u_0(x)}{x}\|_{L^\infty}$ cannot be controlled by $\|u_0\|_{\dot C^1}$ unless $u(0,0)=0$.
		\end{remark}
		\begin{remark}\label{mcbdrmklt}
			For small initial data $f_0$ satisfying $\|f_0\|_{\dot W^{1,\infty}(\Omega)}=\eps_0\eta\leq \eps_0\eta_0$, then by Theorem \ref{eqmcls} one can take $\phi\equiv 0$, and there exists $T_0$ depending only on $\eta_0$ such that \eqref{eqmcbdd} has a unique solution $f\in \mathcal{Z}_{T_0,0}^\eta$.
		\end{remark}
		\begin{proof}[Proof of Theorem \ref{mcbdglo}, ii)]\\
			\textbf{Step 1: Large time existence}\\
			In this step, we assume $T\gg 1$, and we will give the estimates of solution $f$ in $[0,T]$ for sufficiently small initial data. By Remark \ref{mcbdrmklt}, if $\|f_0\|_{ W^{1,\infty}}=\eta\leq C^{-1}\sigma$, we can take $\phi=0$ to prove the local existence with initial data $f_0$ in the time interval $[0,T_0]$, and the solution $f$ satisfying 
			\begin{equation*}
				\|f\|_{Z_{T_0}}\leq C\eta\leq\sigma.
			\end{equation*}
			Next if we take $\eta$ small such that $2^{\frac{1+\kappa}{2}}C\eta\leq C^{-1}\sigma$, then by Remark \ref{mcbdrmklt} again and uniqueness, we can take $f(\frac{T_0}{2})$ as the new initial data and give the expansion of the solution $f$ onto $[0,\frac{3T_0}{2}]$, such that 
			\begin{equation*}
				\|f\|_{Z_{\frac{3T_0}{2}}}\leq 2^{\frac{1+\kappa}{2}}C^{2}\eta\leq\sigma.
			\end{equation*}
			We set $n=[\frac{T}{T_0}]+1$, and if $\sigma$ is small enough such that $((2n-1)!)^{\frac{1+\kappa}{2}}C^{2n-1}\eta\leq C^{-1}\sigma$, then by Remark \ref{mcbdrmklt} we can repeat this procedure $2n-1$ times to get
			\begin{equation*}
				\|f\|_{Z_{nT_0}}\leq ((2n-1)!)^{\frac{1+\kappa}{2}}C^{2n-1}\eta.
			\end{equation*}
			By Stirling formula, the condition $((2n-1)!)^{\frac{1+\kappa}{2}}C^{2n-1}\sigma\leq C^{-1}\sigma$ leads to
			\begin{equation}\label{z1}
				e^{C_0T\log T}\eta\leq \sigma,
			\end{equation}
			for some $C_0$ depends on $C$. \\
			Then \eqref{eqmcbdd} admits a unique solution in $[0,T]$ such that
			\begin{equation*}
				\|f\|_{Z_T}\leq e^{C_0T\log T}\sigma
			\end{equation*}
			provided that \eqref{z1} holds.\\
			Let $T_*\geq 2$ satisfy
			\begin{equation}\label{mcbdgi}
				C_1e^{-\frac{1}{d+5}c_0t}(1+t)^{100}\leq\frac{\iota}{10},\quad \forall t\geq T_*, 
			\end{equation}
			where $\iota\leq\sigma$, $C_1$ and $c_0$ will be decided later. Take $\eta$ small enough such that
			\begin{equation}\label{mcbddefs}
				e^{C_0T_*\log T_*}\eta\leq \frac{\iota}{100},
			\end{equation}
			Moreover,  \eqref{eqmcbdd} admits a unique solution on $[0,T_*]$ satisfying $\|f(T_*)\|_{C^2}\leq \frac{\iota}{50}$. So, we can  define
			\begin{equation}\label{bddmct*}
				T^*=\sup\left\{t>T_*:\sup_{t'\in[T_*,t]}\|f(t')\|_{C^2}\leq\iota\right\}.
			\end{equation}
			We want to prove by contradiction that $T^*=\infty$.
			\\
			\textbf{Step 2: Exponential decay}\\
			We first prove the exponential decay of $L^2$ norm, which is
			\begin{equation}\label{mcbdexpd}
				\|f(t)\|_{L^2}\leq e^{-c_0(t-1)}\|f(1)\|_{L^2}, \quad T^*\geq t\geq 1,
			\end{equation}
			for some $c_0>0$. To prove this, multiply $f$ as a test function to both sides of the system \eqref{eqmcbdd}. By the boundary condition $f|_{\partial\Omega}=0$ and integration by parts, we have
			\begin{equation*}
				\begin{aligned}
					\frac{1}{2}\frac{d}{dt}\|f(t)\|_{L^2}^2+\int_{\Omega} \A_{ij}(\nabla f)\partial_if\partial_jfdx&\leq \|\nabla \A(\nabla f)(t)\|_{L^\infty}\|\nabla^2f(t)\|_{L^\infty}\|\nabla f(t)\|_{L^2}^2\|f(t)\|_{L^\infty}\\
					&\leq C_2\text{diam}(\Omega) \iota^2\|\nabla f(t)\|_{L^2}^2,\quad T^*>t>1.
				\end{aligned}
			\end{equation*}
			Then by the positive definite property of $\A$, we have $\int_{\Omega} \A_{ij}(\nabla f)\partial_if\partial_jfdx\geq C_3\|\nabla f\|_{L^2}^2$. For small $\iota$ such that 
			\begin{equation*}
				C_2\text{diam}(\Omega)\iota^2\leq\frac{C_3}{2},
			\end{equation*}
			we have
			\begin{equation*}
				\partial_t\|f(t)\|_{L^2}^2+C_3\|\nabla f(t)\|_{L^2}^2\leq 0.
			\end{equation*}
			So by Poincare inequality and Gronwall lemma we get \eqref{mcbdexpd} with $c_0=\frac{1}{2}C_3\text{diam}(\Omega)$.\\
			\textbf{Step 3: Long time behavior}\\
			We will take an exact index, say that
			\begin{equation}\label{bdmcltu}
				\|\nabla^2f(T^*)\|_{\dot C^{\frac{3}{4}}}\leq C_4(1+T^*)^{100},
			\end{equation}
			which is a direct corollary if we can prove the same result for $Z_T$ norm. We have 
			\begin{equation*}
				\begin{aligned}
					&\partial_tf(t,x) -\Delta f(t,x) =F(t,x) ,\quad\text{in}\ [0,T^*]\times\Omega,\\
					&f(t,x)=0,\quad\text{on}\ [0,T^*]\times\partial\Omega,\\
					&f(0,x)=f_0(x),
				\end{aligned}
			\end{equation*}
			with
			\begin{equation*}
				F = (\A[\nabla f]:\nabla^2f-\Delta f)= ((\A[\nabla f]-\A[0]):\nabla^2f),\quad F |_{\partial\Omega}=0.
			\end{equation*}
			Apply \eqref{bddmciorh} and $\eqref{bdmcestbd}$ with $\phi,g\equiv 0$ and $\eps$ small enough, together with interpolation, we get the estimate
			\begin{equation*}
				\begin{aligned}
					\|f\|_{Z_T}\leq C_5\left(\|f_0\|_{W^{1,\infty}}+T^{\frac{1}{2}}(1+T)^{n+4}\sup_{t\in[0,T]}\left(\|\nabla f(t)\|_{L^\infty}+t^{\frac{1}{2}}\|\nabla^2f(t)\|_{L^\infty}\right)\right).
				\end{aligned}
			\end{equation*}
			So we have proved the estimate \eqref{bdmcltu} for some suitable $n$ and $C_4=100C_5$ by taking $\iota<\frac{1}{20C_5}$ and take valut at $t=T^*$. \\
			\textbf{Step 4: Global existence}\\
			We prove that $T^*=\infty$. If this is not true, \textit{i.e.} $T^*<\infty$, we denote 
			\begin{align*}
				&\|f(T^*)\|_{C^2}=||\nabla f(T^*)||_{L^\infty}+||\nabla^2 f(T^*)||_{L^\infty}:=X_1+X_2.
			\end{align*}
			Fix 
			\begin{equation*}
				\iota=\frac{1}{10}\min\left\{\eta_0,\frac{1}{20C_5},\left(\frac{C_3}{2C_2\text{diam}(\Omega)}\right)^{-\frac{1}{2}}\right\},
			\end{equation*} 
			such that all the estimates in Step 2 and Step 3 hold. \\
			Using the exponential decay \eqref{mcbdexpd}, long time estimate \eqref{bdmcltu},  assumption of $T_*$ \eqref{mcbdgi}  and the following  interpolation
			\begin{equation}
				||\nabla^m f||_{L^\infty} \leq C_7||f||_{L^2}^{a_m} \|\nabla^2f\|_{\dot C^{\frac{3}{4}}}^{1-a_m}, \quad m=1,2,\quad a_m=\frac{11-4m}{2d+11},
			\end{equation}
			we get 
			\begin{align*}
				&X_1\leq C_7\|f(T^*)\|_{L^2}^{\frac{7}{2d+11}}\|\nabla^2f(T^*)\|_{\dot C^{\frac{3}{4}}}^{\frac{2d+4}{2d+11}}\leq(1+C_4)C_7e^{-\frac{7}{2d+11}c_0T^*}(1+T^*)^{100},\\
				&X_2\leq  C_7\|f(T^*)\|_{L^2}^{\frac{3}{2d+11}}\|\nabla^2f(t)\|_{\dot C^{\frac{3}{4}}}^{\frac{2d+8}{2d+11}}\leq (1+C_4)C_7e^{-\frac{3}{2d+11}c_0T^*}(1+T^*)^{100}.
			\end{align*}
			We take $\iota$ satisfying \eqref{mcbdgi}, and let $C_1=(1+C_4)C_7$, then 
			\begin{equation}
				X_1+X_2\leq \frac{\iota}{2},
			\end{equation}
			which contradicts the definition of $T^*$, and gives $T^*=\infty$. \\
			\textbf{Step 5: Higher regularity} \\
			We have proved 
			\begin{align*}
				\|f(t)\|_{C^2}\leq Ce^{-ct}\|f_0\|_{C^2},\quad t\geq 1.
			\end{align*}
			Then we repeat the proof for Step 1 below again for any $t_0$ with $\|f(t_0)\|_{C^2}\leq \sigma$, to get 
			\begin{equation*}
				\sup_{t\in[t_0,t_0+T_0]}(t-t_0)^{m+\frac{1+\kappa}{2}}\|f(t)\|_{\dot C^{2m+2+\kappa}}\leq \iota,
			\end{equation*}
			then by interpolation, we can get \eqref{bddmcgloe}. This completes the proof of Theorem \ref{mcbdglo} \textit{ii)}.
			
		\end{proof}
		
		\section{Comments}
		
		1) For the Muskat equation without surface tension,
		\begin{align*}
			\partial_t f+\frac{\Lambda f}{\langle\partial_x f\rangle^2}=\mathcal{N}[f].
		\end{align*}
		Denote $F=\partial_x f$, one has 
		\begin{align}\label{musnost}
			\partial_t F+\frac{\Lambda F}{\langle F\rangle^2}-\frac{F\mathcal{H}F} {\langle F\rangle^2}\partial_x F=\partial_x N[f].
		\end{align}
		In this case, we can not apply our methods directly to obtain well-posedness results in critical spaces because the highest order terms include not only a dissipative term but also an extra drift term $\frac{F\mathcal{H}F} {\langle F\rangle^2}\partial_x F$. However, if we have $\mathcal{H}F$ bounded, our methods still work. To see this, we consider the  evolution equation 
		\begin{align*}
			&\partial_t f+\mathcal{L} f=0,\\
			&f|_{t=0}=f_0
		\end{align*}
		Suppose $K(t,x,y)$ is the fundamental solution to the above Cauchy problem. 
		Then, when taking a drift term into account,
		\begin{align}\label{drift}
			\partial_t f+\mathcal{L} f+A(t)\partial_xf=0,
		\end{align}
		the fundamental solution can be written as 
		$\tilde K(t,s,x,y)=K(t-s,x,y-\int_s^t A(\tau)d\tau)$. The Schauder-type estimates still hold for \eqref{drift} provided the boundedness of $A(t)$. 
		For the Muskat problem \eqref{musnost}, to ensure the coefficient of the drift term bounded, it requires the initial data $f_0\in \dot C^{1,\log^\varkappa  }$ with $\varkappa>1$.  Hence, our method can be applied to obtain the well-posedness results for the Muskat equation without surface tension with initial data $f_0\in \dot C^{1,\log^\varkappa  }$, $\varkappa>1$.\\
		2) The ideas in this paper can also be applied to obtain well-posedness results of equations in kinetic theory: Landau/Boltzmann equations with critical data.
		The well-posedness of the Landau/Boltzmann equations appears to be a challenging and important topic in kinetic theory. We expect that the ideas in this paper can also be applied to solve well-posedness problems of Landau/Boltzmann equations with critical data. This will be further explored in our future work.
		\section{Appendix}
		\subsection{Besov space and some interpolation inequalities}
		We first recall the definition of homogeneous Besov spaces.
		\begin{definition}
			Let $s$ be a real number. The homogeneous Besov space $\dot B^{s}_{\infty,\infty}$ consists of those distributions $u$  such that 
			\begin{align*}
				\|u\|_{\dot B^s_{\infty,\infty}}:=\sup_{j\in\mathbb{Z}}2^{js}\|\dot \Delta_ju\|_{L^\infty}<\infty.
			\end{align*}
			Here $\{\dot \Delta_j\}_{j\in\mathbb{Z}}$ are the standard Littlewood-Paley decomposition blocks: $\dot \Delta_j f=\mathcal{F}^{-1}(\phi_j(\xi)\hat f(\xi))$, where \begin{align}\label{phidecom}
				\phi_j(\xi)=\phi(2^{-j}\xi),\ \  \operatorname{supp}(\phi)\subset \{\xi: \frac{3}{4}\leq|\xi|\leq \frac{8}{3}\}, \ \ \sum_{j\in\mathbb{Z}}\phi(2^{-j}\xi)\equiv 1.
			\end{align}
		\end{definition}
		\begin{remark}
			For $s\in\mathbb{R}^+\backslash \mathbb{N}$, the Besov space $\dot B^s_{\infty,\infty}$ coincides the homogeneous H\"{o}lder space $\dot C^s$. This is not true for $s\in \mathbb{N}$. The $\dot B^1_{\infty,\infty}$ norm has the following equivalent definition (see \cite{Triebel}).
			$$\|u\|_{\dot B^1_{\infty,\infty}}\sim \sup_{x,y}\frac{|2u(x)-u(x+y)-u(x-y)|}{|y|}.
			$$
			Since $\|\nabla u\|_{\dot B^{s}_{p,r}}\sim \|u\|_{\dot B^{s+1}_{p,r}}$, we also have
			$$\|u\|_{\dot B^n_{\infty,\infty}}\sim \sup_{x,y}\frac{|2\nabla^{n-1}u(x)-\nabla^{n-1}u(x+y)-\nabla^{n-1}u(x-y)|}{|y|},
			$$
			for any $n\in\mathbb{N}$.
		\end{remark}
		\begin{proposition}\cite{Triebel}
			Let $1\leq p,q\leq \infty$ and $s,\alpha \in\mathbb{R}$. Then $(-\Delta)^\frac{\alpha}{2}:\dot B^s_{p,q}\to \dot B^{s-\alpha}_{p,q}$ is an
			isomorphism.
		\end{proposition}
		In the following proposition, we present a characterization of Besov space $\dot B^0_{\infty,\infty}$ that relies on the the fundamental solution defined by \eqref{defk0}.
		\begin{proposition}\label{normequ}
			For any $a>0$, $x_0\in\mathbb{R}^d$,  denote $\mathcal{K}(f)(t,x)=\K_{x_0}(t,0)\ast f$, where $\K_{x_0}(t,0,x)$ is defined by \eqref{defk0}.
			then there holds 
			\begin{align}\label{rrr}
				\sup_{t>0} t^\frac{a}{s}\|\mathcal{K}(f)\|_{\dot C^a}(t)\lesssim 	\|f\|_{\dot B^0_{\infty,\infty}}.
			\end{align}
		\end{proposition}
		\begin{proof}
			By the homogeneous Littlewood-Paley decomposition $\mathrm{Id}=\sum_{j\in\mathbb{Z}}\dot \Delta_j$, we can write 
			\begin{align}\label{e1}
				\|\mathcal{K}(f)\|_{\dot C^\kappa}\lesssim \sum_{j\in\mathbb{Z}}\|\dot\Delta_j\mathcal{K}(f)\|_{\dot C^\kappa}.
			\end{align}
			For any $j\in\mathbb{Z}$, we can write 
			\begin{align*}
				\dot\Delta_j\mathcal{K}(f)=G_j(t,\cdot) \ast \dot\Delta_jf,
			\end{align*}
			where $G_j(t,x)=\left(\mathcal{F}^{-1}(\hat{\K}(t,0,\xi)\sum_{i=-2}^2\phi_{j+i}(\xi))\right)(x)$ with $\{\phi_k\}_{k\in\mathbb{Z}}$ defined in \eqref{phidecom}. 
			Then 
			\begin{align}\label{e2}
				\|	\delta_\alpha\nabla^{[\kappa]}\dot\Delta_j\mathcal{K}(f)\|_{L^\infty}(t)\lesssim \|\delta_{\alpha}\nabla^{[\kappa]}G_j(t,\cdot)\|_{L^1}\|\dot\Delta_jf\|_{L^\infty},
			\end{align}
			here $[\kappa]$ denotes the integer part of $\kappa$.
			It is easy to check that 
			\begin{align*}
				\|\nabla^m G_j(t,\cdot )\|_{L^1}\lesssim e^{-ct2^{sj}}2^{jm},\ \ m\in\mathbb{N},
			\end{align*}
			for some $c>0$.
			This implies that 
			\begin{align*}
				\sup_\alpha\frac{\|\delta_\alpha \nabla^{[\kappa]} G_j(t,\cdot)\|_{L^1}}{|\alpha|^{\kappa-[\kappa]}}\lesssim e^{-ct2^{sj}}2^{j\kappa}.
			\end{align*}
			Combining this with \eqref{e1} and \eqref{e2} to obtain 
			\begin{align*}
				\|\mathcal{K}(f)\|_{\dot C^\kappa}(t)\lesssim \sum_{j\in\mathbb{Z}}\left( e^{-ct2^{sj}}2^{j\kappa} \right)\sup_{l}\|\dot \Delta_l f\|_{L^\infty}\lesssim t^{-\frac{\kappa}{s}}\sup_{l}\|\dot \Delta_l f\|_{L^\infty}.
			\end{align*}
			We obtain 
			\begin{align*}
				\sup_{t>0}t^\frac{\kappa}{s}\|\mathcal{K}(f)\|_{\dot C^\kappa}(t)\lesssim \|f\|_{\dot B^0_{\infty,\infty}}.
			\end{align*}
			This complete the proof of \eqref{rrr}.
		\end{proof}
		\begin{remark}
			The proof for the special case when $s=2$, $\A(\xi)=|\xi|^2$ and negative Besov index can be found in \cite{Fourierbook}.
		\end{remark}
		We list some elementary inequalities that used frequently throughout the paper. 
		\begin{lemma}~~\\
			1. For any $\gamma_1,\gamma_2\in(-\infty,1)$, there holds
			\begin{align}
				\int_0^t (t-\tau)^{-\gamma_1}\tau^{-\gamma_2}d\tau\sim t^{1-\gamma_1-\gamma_2}.
			\end{align}
			2. For any $0<\eps<1$, $a_1,a_2\in\mathbb{R}^+$, there holds
			\begin{align*}
				\int_{\mathbb{R}^d}\min\left\{\frac{a_1}{|x|^{d+\eps}},\frac{a_2}{|x|^{d-\eps}}\right\}dx\lesssim a_1^\frac{1}{2}a_2^\frac{1}{2}.
			\end{align*}
		\end{lemma}
		\begin{lemma}\label{maininterpo}~~\\
			1. For any $0<\gamma_1<\gamma<\gamma_2$, there holds
			\begin{align*}
				\|f\|_{\dot C^\gamma}\lesssim \|f\|_{\dot C^{\gamma_1}}^\frac{\gamma_2-\gamma}{\gamma_2-\gamma_1} \|f\|_{\dot C^{\gamma_2}}^\frac{\gamma-\gamma_1}{\gamma_2-\gamma_1}.
			\end{align*}
			2. For any $j,j_1,j_2\in\mathbb{N}$ with $0\leq j_1\leq j\leq j_2$, $p\in[1,\infty]$, there holds
			\begin{align*}
				\|\nabla^jf\|_{L^p}\lesssim\|\nabla^{j_1}f\|_{L^p}^\frac{j_2-j}{j_2-j_1}\|\nabla^{j_2}f\|_{L^p}^{\frac{j-j_1}{j_2-j_1}}.
			\end{align*}
		\end{lemma}
		\begin{lemma}\label{douint}
			Let $\kappa_1,\kappa_2\in (0,1)$, for any function $f:\mathbb{R}^d\to \mathbb{R}^N$, there holds
			\begin{equation*}
				\sup_{\alpha\in\mathbb{R}^d}  \frac{\|\delta_\alpha f\|_{L^\infty}}{|\alpha|^{\kappa_1}}\lesssim\|f\|_{\dot C^{\kappa_1}},\quad\quad\quad\quad  \sup_{\alpha,\beta\in\mathbb{R}^d}\frac{\|\delta_\alpha\delta_\beta f\|_{L^\infty}}{|\alpha|^{\kappa_1}|\beta|^{\kappa_2}}\lesssim\|f\|_{\dot C^{\kappa_1+\kappa_2}}.
			\end{equation*}
			Moreover, for $\E^\alpha, \E_\alpha$ defined in \eqref{defEal}, 
			\begin{align*}
				\sup_{\alpha\in\mathbb{R}^d}\frac{\|\E^\alpha f\|_{L^\infty}+\|\E_\alpha f\|_{L^\infty}}{|\alpha|^{\kappa_1}}\lesssim \|f\|_{\dot C^{1+\kappa_1}}.
			\end{align*}
			For  $\tilde{\E}^\alpha$  defined in \eqref{notapes}, and $g:\mathbb{S}\to \mathbb{R}^N$,
			\begin{align*}
				\sup_{\alpha\in\mathbb{S}} \frac{  \|\tilde\E^\alpha g\|_{L^\infty}}{|\alpha|^{\kappa_1}}\lesssim \|g\|_{\dot C^{1+\kappa_1}}.
			\end{align*}
		\end{lemma}
		\begin{lemma}
			\label{lemoalin}
			Denote $\mathcal{O}_\alpha g(x)=\frac{\delta_\alpha\delta_{-\alpha}g(x)}{|\alpha|}$. Then  for any $\nu\in(0,1)$, $\varkappa>1$,
			\begin{align*}
				&\int_{\mathbb{R}^d}\|\mathcal{O}_\alpha g\|_{L^\infty}\frac{d\alpha}{|\alpha|^{d}}\lesssim \min\left\{\|g\|_{\dot C^{1+\nu}}^\frac{1}{2} \|g\|_{\dot C^{1-\nu}}^\frac{1}{2}, \|g\|_{ C^{1,\log^\varkappa}}\right\}.
			\end{align*}
			Moreover, for any $\beta\in\mathbb{R}^d$, and $a\in(0,1)$,
			\begin{align*}
				&\int_{\mathbb{R}^d}\|\delta_\beta\mathcal{O}_\alpha g\|_{L^\infty}\frac{d\alpha}{|\alpha|^{d}}\lesssim|\beta|^{a} \|g\|_{\dot C^{1+a}}.
			\end{align*}
		\end{lemma}
		\begin{lemma}\label{Hol}
			Let $a\in(0,1)$, for any $f,g\in C^a$, there holds
			\begin{align*}
				\|fg\|_{\dot C^a}\lesssim \|f\|_{L^\infty}\|g\|_{\dot C^a}+\|f\|_{\dot C^a}\|g\|_{L^\infty}.
			\end{align*}
		\end{lemma}
		\begin{lemma}(\cite[Lemma 2.2]{KN}) Let $\theta_1\in(0,1)$. For any function $f$ and any $0<\varepsilon_0<\frac{1}{2}\min\{\theta_1,1-\theta_1\}$, there hold
			\begin{equation}\label{interpfrac}
				\|\Lambda^{\theta_1}f\|_{L^\infty}\lesssim (\|f\|_{\dot C^{\theta_1-\varepsilon_0}}\|f\|_{\dot C^{\theta_1+\varepsilon_0}})^{\frac{1}{2}},
			\end{equation}	
			\begin{equation}\label{interpfrac2}
				\|f\|_{\dot C^{\theta_1}}\lesssim \|\Lambda^{\theta_1}f\|_{L^\infty}.
			\end{equation}	
			
		\end{lemma}
		\begin{lemma}\label{lemcom}
			Let $m\in\mathbb{N}$, $\alpha\in(0,1)$. Consider $g,g_1,g_2:\mathbb{R}^d\to\mathbb{R}$, and $f:\mathbb{R}\to \mathbb{R}$ satisfying 
			\begin{align*}
				\sum_{k=1}^{m+2}\|f^{(k)}\|_{L^\infty}\lesssim 1.
			\end{align*}
			Then 
			\begin{equation}\label{mainint11}
				\|\nabla^m(f\circ g)\|_{L^\infty}\lesssim  \|g\|_{\dot C^1}^m+\| g\|_{\dot C^m},
			\end{equation}
			\begin{equation}\label{mainint12}
				\|\nabla^m(f\circ g)\|_{\dot C^\alpha}\lesssim  \| g\|_{\dot C^\alpha}^\frac{m+\alpha}{\alpha}+\|g\|_{\dot C^{m+\alpha}},
			\end{equation}
			\begin{equation}\label{mainint21}
				\|\nabla^m(f\circ g_1-f\circ g_2)\|_{L^\infty}\lesssim\sum_{n=0}^m \|g_1-g_2\|_{\dot C^n}(\|(g_1,g_2)\|_{\dot C^1}^{m-n}+\|(g_1,g_2)\|_{\dot C^{m-n}}),
			\end{equation}
			\begin{equation}\label{mainint22}
				\begin{aligned}
					\|\nabla^m(f\circ g_1-f\circ g_2)\|_{\dot C^\alpha}\lesssim &
					\sum_{n=0}^m\left\{\|g_1-g_2\|_{\dot C^{n+\alpha}}(\|(g_1,g_2)\|_{\dot C^1}^{m-n}+\|(g_1,g_2)\|_{\dot C^{m-n}})\right.\\
					&\left.\quad\quad+\|g_1-g_2\|_{\dot C^n}(\|(g_1,g_2)\|_{\dot C^\alpha}^\frac{m-n+\alpha}{\alpha}+\|(g_1,g_2)\|_{\dot C^{m-n+\alpha}})\right\}.
				\end{aligned}
			\end{equation}
			\begin{equation}\label{mainint31}
				\begin{aligned}
					&\left\|\nabla^m\big((f\circ g_1-f\circ g_2)-(f\circ g_3-f\circ g_4)\big)\right\|_{L^\infty}\\
					&\lesssim \sum_{n=0}^m\|(g_1-g_2)-(g_3-g_4)\|_{\dot C^n}(\|(g_1,g_2)\|_{\dot C^1}^{m-n}+\|(g_1,g_2)\|_{\dot C^{m-n}})\\
					&\quad\quad\quad+\sum_{n_1+n_2+n_3=n}\|g_3-g_4\|_{\dot C^{n_1}}\|(g_1-g_3,g_2-g_4)\|_{\dot C^{n_2}}\left(\sum_{k=1}^4(\|g_k\|_{\dot C^1}^{n_3}+\|g_k\|_{\dot C^{n_3}})\right),
				\end{aligned}
			\end{equation}
			
		\end{lemma}
		\begin{proof}
			The first two inequalities can be obtained by Lemma \ref{maininterpo} and Lemma \ref{Hol}, and we omit the proof.\\
			For \eqref{mainint21} and \eqref{mainint22}, since
			\begin{equation*}
				f\circ g_1-f\circ g_2=\int_0^1(g_1-g_2)\cdot\nabla f(g_1-\lambda(g_1-g_2))d\lambda.
			\end{equation*}
			We can get the result by Lemma \ref{maininterpo}, Lemma \ref{Hol}, \eqref{mainint11} and \eqref{mainint12}.\\
			For the last two inequalities, similarly,
			\begin{align*}
				&(f\circ g_1-f\circ g_2)-(f\circ g_3-f\circ g_4)\\
				&=\int_0^1(g_1-g_2)\cdot\nabla f(g_1-\lambda(g_1-g_2))d\lambda-\int_0^1(g_3-g_4)\cdot\nabla f(g_3-\lambda(g_3-g_4))d\lambda\\
				&=\int_0^1\int_0^1(g_3-g_4)(g_3-\lambda(g_3-g_4))\nabla^2f\left(g_3-\lambda(g_3-g_4)-\mu\big((g_1-g_3)-\lambda((g_1-g_3)-(g_2-g_4))\big)\right)\\
				&\quad+\int_0^1\left((g_1-g_2)-(g_3-g_4)\right)\cdot\nabla f(g_1-\lambda(g_1-g_2))d\lambda.
			\end{align*}
			Then we get \eqref{mainint31} by Lemma \ref{maininterpo}, Lemma \ref{Hol}, \eqref{mainint11} and \eqref{mainint12}.
		\end{proof}
		\subsection{Estimates of nonlinear terms in the Muskat equation}
		Let $\kappa,m$ be as fixed in Theorem \ref{thmmst}, the norm $\|\cdot\|_T$, $\|\cdot \|_{X_T}$ be defined in \eqref{normmst} and \eqref{musnol}. 
		For brevity of statement, we also define the following non-endpoint norm
		\begin{align}\label{defnn}
			&\|h\|_{T,*}:=\sup_{t\in[0,T]}\sum_{\iota=\frac{1}{5},m+\kappa}t^\frac{1}{3\iota}\|\partial_xh(t)\|_{\dot C^\iota}.
		\end{align}
		The non-endpoint norm provides desired smallness to proceed the fixed point argument and get local solution under the condition \eqref{stcon}. More precisely, as we can see in the proof of the local well-posedness result of Theorem \ref{thmmst},  the contraction mapping theorem is performed in the set with center $\phi=f_0\ast \rho_{\eps_1}$. In this case, $\|f\|_T$ is merely bounded, but the non-endpoint norm
		\begin{align*}
			\|f\|_{T,*}\leq \|f-\phi\|_{T,*}+\|\phi\|_{T,*}\lesssim \|f-\phi\|_{T,*}+T^\frac{1}{15}\|
			\phi\|_{C^{m+3}}
		\end{align*}
		is small because of the smallness of $\|f-\phi\|_{T,*}$ and $T$. Hence, in the following estimates of nonlinear terms, we will keep a $\|\cdot\|_{T,*}$ norm in the result.
		\begin{lemma}\label{nonst}
			Let $\N[f](t,x)$  be as defined in \eqref{defNst} with the constant $\varrho_0$ defined in \eqref{defvr0},  then for any $T>0$,
			\begin{equation}\label{NNN}
				\begin{aligned}
					&		\sum_{j=0}^m\sup_{t\in[0,T]}t^\frac{j+\kappa}{3}\| \N[f](t)\|_{\dot C^{j+\kappa-2}}\lesssim \|f\|_{T,*}^2(1+\|f\|_T)^{2m+5}+|\varrho_0|(T^\frac{2}{3}+T^\frac{1}{10})\|f\|_{X_T}(1+\|f\|_{X_T})^{m+5},
				\end{aligned}
			\end{equation}
			and
			\begin{equation}\label{Ndi}
				\begin{aligned}
					&		\sum_{j=0}^m\sup_{t\in[0,T]}t^\frac{j+\kappa}{3}\| \N[f_1](t)- \N[f_2](t)\|_{\dot C^{j+\kappa-2}})\\
					&\lesssim \|f_1-f_2\|_T\|(f_1,f_2)\|_{T,*}(1+\|(f_1,f_2)\|_T)^{2m+5}+|\varrho_0|(T^\frac{2}{3}+T^\frac{1}{10})||f_1-f_2||_{X_T}(1+||(f_1,f_2)||_{X_T})^{m+4}.
				\end{aligned}
			\end{equation}	
		\end{lemma}
		\begin{proof}
			For simplicity, we denote the quantity
			\begin{align*}
				\mathbf{C}_n(f_1,f_2):=\left(\|f_1-f_2\|_{T,*}\|(f_1,f_2)\|_{T,*}+\|f_1-f_2\|_{T}\|(f_1,f_2)\|_{T,*}^2\right)(1+\|(f_1,f_2)\|_T)^{2n+5}.
			\end{align*}
			To prove \eqref{NNN} and \eqref{Ndi}, it suffices to show that, for any $n\in\mathbb{N}$, $n\leq m+1$,
			\begin{equation}\label{eee1}
				\begin{aligned}
					&\sup_{t\in[0,T]}t^{\frac{n+2}{3}}\|\N[f_1](t)-\N[f_2](t)\|_{\dot C^{n}}\\
					&\ \ \quad\quad\quad\quad\lesssim   \mathbf{C}_n(f_1,f_2)+|\varrho_0|(T^\frac{2}{3}+T^\frac{1}{10})||f_1-f_2)||_{X_T}(1+||(f_1,f_2)||_{X_T})^{n+4}.
				\end{aligned}
			\end{equation}
			Note that $\kappa-2\in(0,1)$, then by Lemma \ref{maininterpo}, the desired estimate \eqref{Ndi} follows from \eqref{eee1}, and \eqref{NNN} follows by taking $f_1=f$, $f_2\equiv 0$ in \eqref{eee1} since $\N[0]\equiv 0$.\\
			For simplicity, we fix $t\in[0,T]$ and omit the time variable in our proof. Denote 
			\begin{equation}\label{defEal}
				\begin{aligned}
					&	\mathsf E_\alpha f(x)=\partial_x f(x)-\Delta_\alpha f(x),\ \ \ \ \ \ \E^\alpha f(x)=\partial_x f(x-\alpha)-\Delta_\alpha f(x),\\
					&\B[f](x,\alpha)=\frac{ \Delta_{\alpha} f(x)\E_\alpha f(x)}{\alpha\left\langle\Delta_{\alpha} f(x)\right\rangle^{2}},\ \ \quad\quad \mathsf{M}[f](x)=\frac{\partial_{x}^{2} f(x)}{\left\langle\partial_{x} f(x)\right\rangle^{3}}.
				\end{aligned}
			\end{equation}
			Then we have $\mathsf{N}[f]=\N_1[f]+\N_2[f]+\varrho_0\N_3[f]$ with 
			\begin{align*}
				&\N_{1}[f](x)=\frac{1}{\pi} \int_{\mathbb{R}}\B[f](\alpha,x)\partial_{x}\M[f](x-\alpha) {d \alpha},\\
				&\N_2[f](x)=-\frac{1}{\pi} \int_{\mathbb{R}}\partial_x^2f(x-\alpha)\left(\frac{1}{\left\langle\partial_{x} f(x-\alpha)\right\rangle^{3}}-\frac{1}{\left\langle\partial_{x} f(x)\right\rangle^{3}}\right)\frac{d\alpha}{\alpha^2},\\
				&\N_3[f](x)=-\frac{1}{\pi}\int_{\mathbb{R}}\B[f](\alpha,x)\partial_x f(x-\alpha){d\alpha}-\Lambda f(x).
			\end{align*}
			We first consider $\N_1[f_1]-\N_1[f_2]$, we have 
			\begin{align*}
				\N_1[f_1](x)-
				\N_1[f_2](x)=&\frac{1}{\pi} \int_{\mathbb{R}}\B[f_1](\alpha,x) \partial_{x}\left(\M[f_1]-\M[f_2]\right)(x-\alpha) d \alpha\\
				&+\frac{1}{\pi} \int_{\mathbb{R}}(\B[f_1]-\B[f_2])(\alpha,x) \partial_{x}\M[f_2](x-\alpha) d \alpha\\
				:=&\mathrm{I}_1+\mathrm{I}_2.
			\end{align*}
			Integrating  by parts to obtain 
			\begin{align*}
				\mathrm{I}_1&=-\frac{1}{\pi} \int_{\mathbb{R}}\B[f_1](\alpha,x) \partial_{\alpha}\left(\M[f_1]-\M[f_2]\right)(x-\alpha)d \alpha\\
				&=-\frac{1}{\pi} \int_{\mathbb{R}}\partial_{\alpha}\B[f_1](\alpha,x)\delta_\alpha\left(\M[f_1]-\M[f_2]\right)(x) {d \alpha}.
			\end{align*}
			Then, for $n\in\mathbb{N}$, $n\leq m+1$,
			\begin{align}\label{fori1}
				\partial_x^{n}\mathrm{I}_1=-\frac{1}{\pi}\sum_{k+l=n}\int_{\mathbb{R}} \partial_x^k\partial_{\alpha}\B[f_1](\alpha,x)\delta_\alpha\partial_x^l\left(\M[f_1]-\M[f_2]\right)(x)d\alpha.
			\end{align}
			Denote $\psi(r)=\frac{r}{\langle r\rangle^2}$, then $\B[f](x,\alpha)=\frac{\psi(\Delta_\alpha f(x))\E_\alpha f(x)}{\alpha}$. 	Note that $-\partial_\alpha(\E_\alpha f(x))=\partial_\alpha (\Delta_\alpha f(x))=\frac{\E^\alpha f(x)}{\alpha}$. Hence,  \begin{align*}
				\partial_\alpha \B[f_1](\alpha,x)=\frac{1}{\alpha^2}\left(\psi'(\Delta_\alpha f(x))\E^\alpha f(x)\E_\alpha f(x)-\psi(\Delta_\alpha f(x))(\E^\alpha f(x)+\E_\alpha f(x))\right).
			\end{align*}
			Then we get 
			\begin{align*}
				| \partial_x^k  \partial_\alpha \B[f](\alpha,x)|\lesssim& \frac{1}{\alpha^2}\left(\sum_{k_1+k_2+k_3=k}\|\partial_x^{k_1}\psi'(\Delta_\alpha f(\cdot))\|_{L^\infty}\|\partial_x^{k_2}\E^\alpha f\|_{L^\infty}\|\partial_x^{k_3}\E_\alpha f\|_{L^\infty}\right.\\
				&\ \left.+\sum_{k_4+k_5=k}\|\partial_x^{k_4}\left(\psi(\Delta_\alpha f(\cdot))\right)\|_{L^\infty}(\|\partial_x^{k_5}\E^\alpha f\|_{L^\infty}+\|\partial_x^{k_5}\E_\alpha f\|_{L^\infty})\right).
			\end{align*}
			Recall the definition of $\|\cdot\|_T$ in \eqref{normmst}. By \eqref{mainint11} in Lemma \ref{lemcom}, we have
			\begin{align*}
				\|\partial_x^{k_1}\left(\psi(\Delta_\alpha f(\cdot))\right)\|_{L^\infty}&+\|\partial_x^{k_1}\left(\psi'(\Delta_\alpha f(\cdot))\right)\|_{L^\infty}\lesssim \|\Delta_\alpha f\|_{\dot C^1}^{k_1}+\|\Delta_\alpha f\|_{\dot C^{k_1}}\\
				&\lesssim \|\partial_x f\|_{\dot C^{k_1}}+\|\partial_x f\|_{\dot C^1}^{k_1}\lesssim t^{-\frac{k_1}{3}}\|f\|_T(1+\|f\|_T)^{k_1}.
			\end{align*}
			Moreover, we have 
			\begin{align*}
				\|\partial_x^{k_3}\E_\alpha f\|_{L^\infty}\lesssim \|\partial_x f\|_{\dot C^{k_3}}\lesssim t^{-\frac{k_3}{3}}\|f\|_T.
			\end{align*}
			This yields that 
			\begin{equation}\label{b1}
				\begin{aligned}
					&|\partial_{\alpha}\partial_x^k\B[f](\alpha,x)|
					\lesssim \frac{1}{\alpha^2}\|f\|_T(1+\|f\|_T)^{k+1}\sum_{k_1+k_2=k} t^{-\frac{k_1}{3}}(\|\E^\alpha \partial_x^{k_2}f\|_{L^\infty}+\|\E_\alpha \partial_x^{k_2}f\|_{L^\infty}).
				\end{aligned}
			\end{equation}
			Moreover, by Lemma \ref{douint}, one has 
			\begin{align*}
				\|\partial_x^{k_2}\E^\alpha f\|_{L^\infty}+\|\partial_x^{k_2}\E_\alpha f\|_{L^\infty}&\lesssim \min\{|\alpha|^\frac{3}{4}\|\partial_xf\|_{\dot C^{k_2+\frac{3}{4}}},|\alpha|^\frac{3}{4}\|\partial_xf\|_{\dot C^{k_2+\frac{1}{4}}}\}\\
				&\lesssim t^{-\frac{k_2}{3}}\|f\|_{T,*}\min\{(|\alpha|t^{-\frac{1}{3}})^\frac{3}{4},(|\alpha|t^{-\frac{1}{3}})^\frac{1}{4}\}.
			\end{align*}
			Combining this with \eqref{b1}, we have
			\begin{equation}\label{i1esb}
				\begin{aligned}
					&|\partial_{\alpha}\partial_x^k\B[f](\alpha,x)|\lesssim \frac{1}{\alpha^2}\|f\|_{T,*}\| f\|_T (1+\| f\|_T )^{k+1}t^{-k/3}\min\{(|\alpha|t^{-\frac{1}{3}})^\frac{3}{4},(|\alpha|t^{-\frac{1}{3}})^\frac{1}{4}\}.
				\end{aligned}
			\end{equation}
			Here and in the following we will use Lemma \ref{lemcom} without claim.
			One has 
			\begin{equation*}
				\begin{aligned}
					&\left|\delta_\alpha\partial_x^l\left(\M[f_1]-\M[f_2]\right)(x)\right|\\
					&\quad\lesssim \left|\delta_\alpha\partial_x^l\left( \frac{\partial_{x}^{2} (f_1-f_2)(\cdot )}{\left\langle\partial_{x} f_2(\cdot )\right\rangle^{3}}\right)(x)\right|+\left|\delta_\alpha \partial_x^l\left(\partial_x^2f_1\left(\frac{1}{\left\langle\partial_{x} f_1(\cdot )\right\rangle^{3}}-\frac{1}{\left\langle\partial_{x} f_2(\cdot )\right\rangle^{3}}\right)\right)\right|\\
					&\quad\lesssim \sum_{l_1+l_2=l}\Big(\|\delta_\alpha \partial_x^{2+l_1}(f_1-f_2)\|_{L^\infty}(\|\partial_x^2f_2\|_{L^\infty}^{l_2}+\|\partial_x^{1+l_2}f_2\|_{L^\infty})\\
					&\quad\quad\quad+ \|\partial_x^{2+l_1}(f_1-f_2)\|_{L^\infty}\sum_{l_3+l_4=l_2}\|\delta_\alpha \partial_x^{l_3}(f_1,f_2)\|_{L^\infty}(\|\partial_x^2(f_1,f_2)\|_{L^\infty}^{l_4}+\|\partial_x^{1+l_4}(f_1,f_2)\|_{L^\infty})\\
					&\quad\quad\quad+\|\delta_\alpha \partial_x^{1+l_1}(f_1-f_2)\|_{L^\infty}(\|\partial_x^2(f_1,f_2)\|_{L^\infty}^{1+l_2}+\|\partial_x^{2+l_2}(f_1,f_2)\|_{L^\infty})\\
					&\quad\quad\quad +\|\partial_x^{1+l_1}(f_1-f_2)\|_{L^\infty}\sum_{l_3+l_4=l_2}\|\delta_\alpha \partial_x^{2+l_3}(f_1,f_2)\|_{L^\infty}(\|\partial_x^2(f_1,f_2)\|_{L^\infty}^{l_4}+\|\partial_x^{1+l_4}(f_1,f_2)\|_{L^\infty})\Big)\\
					&\quad\lesssim t^{-\frac{l+1}{3}}(1+\|(f_1,f_2)\|_{T})^{l+1}(\|f_1-f_2\|_T\|(f_1,f_2)\|_{T,*}+\|f_1-f_2\|_{T,*})\min\{(|\alpha|t^{-\frac{1}{3}})^\frac{3}{4},(|\alpha|t^{-\frac{1}{3}})^\frac{1}{4}\}.
				\end{aligned}
			\end{equation*}
			Combining this with \eqref{i1esb} and \eqref{fori1},  we obtain 
			\begin{equation}\label{i1}
				\begin{aligned}
					\|\partial_x^n\mathrm{I}_1\|_{L^\infty}&\lesssim t^{-\frac{n+1}{3}}\mathbf{C}_n(f_1,f_2)\int_{\mathbb{R}}\min\{(|\alpha|t^{-\frac{1}{3}})^\frac{3}{2},(|\alpha|t^{-\frac{1}{3}})^\frac{1}{2}\}\frac{d\alpha}{\alpha^2}\\
					&\lesssim t^{-\frac{n+2}{3}}\mathbf{C}_n(f_1,f_2).
				\end{aligned}
			\end{equation}
			Similarly, integrating by parts one has 
			\begin{align*}
				\mathrm{I}_2=-\frac{1}{\pi} \int_{\mathbb{R}}\partial_\alpha (\B[f_1]-\B[f_2])(\alpha,x) \delta_\alpha \M[f_2](x) {d \alpha},
			\end{align*}
			and
			\begin{align*}
				\partial_x^n\mathrm{I}_2=-\frac{1}{\pi} \sum_{k+l=n}\int_{\mathbb{R}}\partial_\alpha \partial_x^k(\B[f_1]-\B[f_2])(\alpha,x) \delta_\alpha \partial_x^l\M[f_2](x) {d \alpha}.
			\end{align*}
			We can see that
			\begin{equation*}
				(\B[f_1]-\B[f_2])(\alpha,x)= \frac{\Delta_\alpha(f_1-f_2)\E_\alpha f_1}{\alpha\langle \Delta_\alpha f\rangle^2}+\frac{\Delta_\alpha f_2\E_\alpha (f_1-f_2)}{\alpha\langle \Delta_\alpha f\rangle^2}+\frac{\Delta_\alpha f_2\E_\alpha f_2}{\alpha}\left(\frac{1}{\langle\Delta_\alpha f_1\rangle^2}-\frac{1}{\langle\Delta_\alpha f_2\rangle^2}\right).
			\end{equation*}
			Then we can use similar estimates as \eqref{b1} to obtain
			\begin{equation}\label{b2}
				\begin{aligned}
					&\left|\partial_\alpha \partial_x^k(\B[f_1]-\B[f_2])(\alpha,x) \right|\\
					&\lesssim \frac{1}{\alpha^2}\left(\sum_{k_1+k_2+k_3=k}(\|(f_1,f_2)\|_{\dot C^{k_1+1}}+\|(f_1,f_2)\|_{\dot C^2}^{k_1})(\|(\E_{\alpha}f_1,\E_{\alpha}f_2)\|_{\dot C^{k_2}}+\|(\E_{\alpha}f_1,\E_{\alpha}f_2)\|_{\dot C^1}^{k_2})\right.\\
					&\quad\quad\quad\quad\quad\quad\times(\|f_1-f_2\|_{\dot C^{k_3+1}}+\|f_1-f_2\|_{\dot C^2}^{k_3})\\
					&\left.\quad\quad+\sum_{k_1+k_2=k}(\|(f_1,f_2)\|_{\dot C^{k_1+1}}+\|(f_1,f_2)\|_{\dot C^2}^{k_1})(\|(\E_{\alpha}(f_1-f_2)\|_{\dot C^{k_2}}+\|(\E_{\alpha}(f_1-f_2)\|_{\dot C^1}^{k_2})\right)
				\end{aligned}
			\end{equation}
			and by definition of $\|\cdot\|_T$, we have
			\begin{equation*}
				\begin{aligned}
					&\left|\partial_\alpha \partial_x^k(\B[f_1]-\B[f_2])(\alpha,x) \right|\\
					&\lesssim \frac{1}{t^{\frac{k}{3}}\alpha^2}(1+\|(f_1,f_2)\|_T)^{k+3}(\|f_1-f_2\|_{T}\|(f_1,f_2)\|_{T,*}+\|f_1-f_2\|_{T,*})\min\{(|\alpha|t^{-\frac{1}{3}})^\frac{3}{4},(|\alpha|t^{-\frac{1}{3}})^\frac{1}{4}\},
				\end{aligned}
			\end{equation*}
			By Lemma \ref{lemcom},
			\begin{align*}
				\left|\delta_\alpha \partial_x^l\M[f_2](x)\right|\lesssim &\sum_{l_1+l_2=l}\|\delta_\alpha f_2\|_{\dot C^{l_1+2}}(\|\partial_x f_2\|_{\dot C^{l_2}}+\|\partial_x f_2\|_{\dot C^{1}}^{l_2})\\
				&\quad+\sum_{l_1+l_2=l}\|\partial_x^2f_2\|_{\dot C^{l_1}}(\|\delta_\alpha\partial_x f_2\|_{\dot C^{l_2}}+\|\delta_\alpha\partial_x f_2\|_{\dot C^{1}}^{l_2})
			\end{align*}
			and by definition of $\|\cdot\|_T$, we have
			\begin{equation*}
				\left|\delta_\alpha \partial_x^l\M[f_2](x)\right|\lesssim t^{-\frac{l+1}{3}}\|f_2\|_{T,*}(1+\| f_2\|_T)^l\min\{(|\alpha|t^{-\frac{1}{3}})^\frac{3}{4},(|\alpha|t^{-\frac{1}{3}})^\frac{1}{4}\}.
			\end{equation*}
			Hence one has
			\begin{equation}\label{i2}
				\begin{aligned}
					\|\partial_x^n\mathrm{I}_2\|_{L^\infty}&\lesssim t^{-\frac{n+1}{3}}\mathbf{C}_n(f_1,f_2)\int_{\mathbb{R}}\min\{(|\alpha|t^{-\frac{1}{3}})^\frac{3}{2},(|\alpha|t^{-\frac{1}{3}})^\frac{1}{2}\}\frac{d\alpha}{|\alpha|^2}\\
					&\lesssim  t^{-\frac{n+2}{3}}\mathbf{C}_n(f_1,f_2).
				\end{aligned}
			\end{equation}
			We conclude from \eqref{i1} and \eqref{i2} that 
			\begin{equation}\label{N1}
				\begin{aligned}
					&\sup_{t\in[0,T]}t^\frac{n+2}{3}\|\partial_x^n(N_1[f_1]-N_1[f_2])(t)\|_{L^\infty}\lesssim \mathbf{C}_n(f_1,f_2),
				\end{aligned}
			\end{equation}
			for any $0\leq n\leq m+1$.\\
			Then we consider $N_2[f_1]-N_2[f_2]$, note that 
			\begin{align*}
				N_2[f](x)=-\frac{1}{\pi} \int_{\mathbb{R}}\delta_\alpha \partial_x^2f(x)\delta_\alpha\left(\frac{1}{\left\langle\partial_{x} f(\cdot)\right\rangle^{3}}\right)(x)\frac{d\alpha}{\alpha^2}+\frac{1}{2\pi}\partial_x^2f(x)\int _{\mathbb{R}}\mathcal{O}_\alpha \left(\frac{1}{\left\langle\partial_{x} f(\cdot)\right\rangle^{3}}\right)(x)\frac{d\alpha}{|\alpha|},
			\end{align*}
			where the $\mathcal{O}_\alpha$ is defined in Lemma \ref{lemoalin}.
			We split $N_2[f_1]-N_2[f_2]$ into 
			\begin{equation}\label{splitN2}
				\begin{aligned}
					&(N_2[f_1]-N_2[f_2])(x)\\
					&\quad\quad=-\frac{1}{\pi} \int_{\mathbb{R}}\delta_\alpha \partial_x^2(f_1-f_2)(x)\delta_\alpha\left(\frac{1}{\left\langle\partial_{x} f_1(\cdot)\right\rangle^{3}}\right)(x)\frac{d\alpha}{\alpha^2}\\
					&\quad\quad\quad\quad\ -\frac{1}{\pi} \int_{\mathbb{R}}\delta_\alpha \partial_x^2f_2(x)\delta_\alpha\left(\frac{1}{\left\langle\partial_{x} f_1(\cdot)\right\rangle^{3}}-\frac{1}{\left\langle\partial_{x} f_2(\cdot)\right\rangle^{3}}\right)(x)\frac{d\alpha}{\alpha^2}\\
					&\quad\quad\quad\quad\ +\frac{1}{2\pi}\partial_x^2(f_1-f_2)(x)\int_{\mathbb{R}} \mathcal{O}_\alpha \left(\frac{1}{\left\langle\partial_{x} f_1(\cdot)\right\rangle^{3}}\right)(x)\frac{d\alpha}{|\alpha|}\\
					&\quad\quad\quad\quad\ +\frac{1}{2\pi}\partial_x^2f_2(x)\int_{\mathbb{R}} \mathcal{O}_\alpha \left(\frac{1}{\left\langle\partial_{x} f_1(\cdot)\right\rangle^{3}}-\frac{1}{\left\langle\partial_{x} f_2(\cdot)\right\rangle^{3}}\right)(x)\frac{d\alpha}{|\alpha|}\\
					&\quad\quad:=\sum_{j=1,2}\int_{\mathbb{R}} \mathrm{II}_j(x,\alpha)\frac{d\alpha}{\alpha^2}+\partial_x^2(f_1-f_2)(x) \mathrm{II}_3(x)+\partial_{x}^2f_2(x) \mathrm{II}_4(x).
				\end{aligned}  
			\end{equation}
			It is easy to check that 
			\begin{align*}
				&\|\partial_x^n\mathrm{II}_1(\cdot,\alpha)\|_{L^\infty}\lesssim\left\|\partial_x^n\left(\delta_\alpha \partial_x^2(f_1-f_2)\delta_\alpha\left(\frac{1}{\left\langle\partial_{x} f_1(\cdot)\right\rangle^{3}}\right)\right)\right\|_{L^\infty}\\
				&\lesssim \sum_{l+k=n}\|\delta_\alpha \partial_{x}^{l+2}(f_1-f_2)\|_{L^\infty}\left(\sum_{k_1+k_2=k}\|\delta_\alpha \partial_x^{1+k_1}f_1\|_{L^\infty}(\|\partial_x^{1+k_2}f_1\|_{L^\infty}+\|\partial_x^2f_1\|_{L^\infty}^{k_2})\right)\\
				&\lesssim \sum_{l+k=n}\min\{\| \partial_{x}^{l+2}(f_1-f_2)\|_{\dot C^\frac{3}{4}}|\alpha|^\frac{3}{4},\| \partial_{x}^{l+2}(f_1-f_2)\|_{\dot C^\frac{1}{4}}|\alpha|^\frac{1}{4}\}\\
				&\quad\quad\quad\times\left(\sum_{k_1+k_2=k}\min\{\|\partial_x^{1+k_1}f_1\|_{\dot C^\frac{3}{4}}|\alpha|^\frac{3}{4},\|\partial_x^{1+k_1}f_1\|_{\dot C^\frac{1}{4}}|\alpha|^\frac{1}{4}\}(\|\partial_x^{1+k_2}f_1\|_{L^\infty}+\|\partial_x^2f_1\|_{L^\infty}^{k_2})\right).
			\end{align*}
			Then by the definition of $\|\cdot\|_T$ and $\|\cdot\|_{T,*}$ in \eqref{normmst} and \eqref{defnn}, we obtain 
			\begin{equation}\label{ii1}
				\begin{aligned}
					\|\partial_x^n\mathrm{II}_1(\cdot,\alpha)\|_{L^\infty}\lesssim& t^{-\frac{n+1}{3}}\mathbf{C}_n(f_1,f_2) \min\{(|\alpha|t^{-\frac{1}{3}})^\frac{3}{4},(|\alpha|t^{-\frac{1}{3}})^\frac{1}{4}\}^2.
				\end{aligned}
			\end{equation}
			For $\mathrm{II}_2$, we have 
			\begin{align*}
				&\|\partial_x^n\mathrm{II}_2(\cdot,\alpha)\|_{L^\infty}\lesssim\left\|\partial_x^n\left(\delta_\alpha \partial_x^2f_2\delta_\alpha\left(\frac{1}{\left\langle\partial_{x} f_1(\cdot)\right\rangle^{3}}-\frac{1}{\left\langle\partial_{x} f_2(\cdot)\right\rangle^{3}}\right)\right)\right\|_{L^\infty}\\
				&\lesssim \sum_{k+l=n}\|\delta_\alpha\partial_x^{2+l}f_2\|_{L^\infty}\sum_{k_1+k_2=k}\left\{\|\delta_\alpha \partial_x^{1+k_1}(f_1-f_2)\|_{L^\infty}(\|\partial_x^{1+k_2}(f_1,f_2)\|_{L^\infty}+\|\partial_x^{2}(f_1,f_2)\|_{L^\infty}^{k_2})\right.\\
				&\quad\quad+\left.\|\partial_x^{1+k_1}(f_1-f_2)\|_{L^\infty}\sum_{k_3+k_4=k_2}\Bigg(\Big(\|\delta_\alpha \partial_x^{1+k_3}(f_1,f_2)\|_{L^\infty}\Big)(\|\partial_x^{1+k_4}(f_1,f_2)\|_{L^\infty}+\|\partial_x^{2}(f_1,f_2)\|_{L^\infty}^{k_4})\Bigg)\right\}.
			\end{align*}
			This implies that 
			\begin{equation}\label{ii2}
				\begin{aligned}
					\|\partial_x^n\mathrm{II}_2(\cdot,\alpha)\|_{L^\infty}&\lesssim t^{-\frac{n+1}{3}}\mathbf{C}_n(f_1,f_2)\min\{(|\alpha|t^{-\frac{1}{3}})^\frac{3}{4},(|\alpha|t^{-\frac{1}{3}})^\frac{1}{4}\}^2.
				\end{aligned}
			\end{equation}
			By \eqref{ii1} and \eqref{ii2}, we obtain 
			\begin{equation}\label{ii12}
				\begin{aligned}
					\sum_{i=1,2}\int_{\mathbb{R}} \|\partial_x^n\mathrm{II}_i(\cdot,\alpha)\|_{L^\infty}\frac{d\alpha}{\alpha^2}&\lesssim t^{-\frac{n+1}{3}}\mathbf{C}_n(f_1,f_2) \int_{\mathbb{R}}\min\{(|\alpha|t^{-\frac{1}{3}})^\frac{3}{4},(|\alpha|t^{-\frac{1}{3}})^\frac{1}{4}\}^2\frac{d\alpha}{\alpha^2}\\
					&\lesssim t^{-\frac{n+2}{3}}\mathbf{C}_n(f_1,f_2).
				\end{aligned}
			\end{equation}
			For $\mathrm{II}_3$ and $\mathrm{II}_4$, we apply Lemma \ref{lemoalin} to obtain  
			\begin{equation}\label{ii3}
				\begin{aligned}
					\|\partial_x^n\mathrm{II}_3\|_{L^\infty}&\lesssim\left\|\partial_x^n\left(\frac{1}{\left\langle\partial_{x} f_1(\cdot)\right\rangle^{3}}\right)\right\|_{\dot C^{\frac{1}{2}}}^\frac{1}{2}\left\|\partial_x^n\left(\frac{1}{\left\langle\partial_{x} f_1(\cdot)\right\rangle^{3}}\right)\right\|_{\dot C^{\frac{3}{2}}}^\frac{1}{2}\\
					&\lesssim \|\partial_x f_1\|_{\dot C^\frac{1}{2}}^{2(n+1)}+\|\partial_x f_1\|_{\dot C^{n+3/2}}^\frac{n+1}{n+\frac{3}{2}}\\
					&\lesssim t^{-\frac{n+1}{3}}\|f_1\|_{T,*}(1+\|f_1\|_T)^{2n+2},
				\end{aligned}
			\end{equation}
			and
			\begin{equation}\label{ii4}
				\begin{aligned}
					&\|\partial_x^n\mathrm{II}_4\|_{L^\infty}\lesssim \left\|\partial_x^n\left(\frac{1}{\left\langle\partial_{x} f_1(\cdot)\right\rangle^{3}}-\frac{1}{\left\langle\partial_{x} f_2(\cdot)\right\rangle^{3}}\right)\right\|_{\dot C^{\frac{1}{2}}}^\frac{1}{2}\left\|\partial_x^n\left(\frac{1}{\left\langle\partial_{x} f_1(\cdot)\right\rangle^{3}}-\frac{1}{\left\langle\partial_{x} f_2(\cdot)\right\rangle^{3}}\right)\right\|_{\dot C^{\frac{3}{2}}}^\frac{1}{2}\\
					&\quad\quad\quad\lesssim t^{-\frac{n+1}{3}}(\|f_1-f_2\|_T\|(f_1,f_2)\|_{T,*}+\|f_1-f_2\|_{T,*})(1+\|(f_1,f_2)\|_T)^n.
				\end{aligned}
			\end{equation}
			Combining \eqref{ii3} and \eqref{ii4} with the estimate
			\begin{align*}
				\|\partial_x^{2+n}f\|_{L^\infty}\lesssim t^{-\frac{n+1}{3}}\|f\|_{T,*},
			\end{align*}
			we obtain that 
			\begin{align}\label{ii34}
				\|\partial_x^n(\partial_x^2(f_1-f_2)\mathrm{II}_3)\|_{L^\infty}+\|\partial_x^n(\partial_x^2f_2\mathrm{II}_4)\|_{L^\infty}\lesssim t^{-\frac{n+2}{3}}\mathbf{C}_n(f_1,f_2).
			\end{align}
			We conclude from \eqref{splitN2}, \eqref{ii12} and  \eqref{ii34}  that 
			\begin{equation}\label{N2}
				\begin{aligned}
					&\sup_{t\in[0,T]}t^\frac{n+2}{3}\|\partial_x^n(N_2[f_1]-N_2[f_2])(t)\|_{L^\infty}\lesssim \mathbf{C}_n(f_1,f_2).
				\end{aligned}
			\end{equation}
			Finally, we estimate lower order terms
			\begin{align*}
				N_3[f_1]-N_3[f_2]=&-\frac{1}{\pi}\int_{\mathbb{R}}B[f_1](\alpha,x)\partial_{x}(f_1-f_2)(x-\alpha){d \alpha}-\frac{1}{\pi}\int_{\mathbb{R}}(B[f_1]-B[f_2])(\alpha,x)\partial_{x}f_2(x-\alpha){d \alpha}\\
				&-\Lambda (f_1-f_2)(x)\\
				:=&\mathrm{III}_1+\mathrm{III}_2+\mathrm{III}_3.
			\end{align*}
			Using integration by parts to get 
			\begin{align*}
				&|\partial_x^n\mathrm{III}_1|\lesssim \sum_{k+l=n}\int_{\mathbb{R}}|\partial_\alpha\partial_x^kB[f_1](\alpha,x)||\delta_\alpha\partial_x^l (f_1-f_2)(x)|d\alpha,\\
				&|\partial_x^n\mathrm{III}_2|\lesssim \sum_{k+l=n}\int_{\mathbb{R}}|\partial_\alpha\partial_x^k(B[f_1]-B[f_2])(\alpha,x)||\delta_\alpha \partial_x^lf_2(x)|d\alpha.
			\end{align*}
			Then we use \eqref{b1}, \eqref{b2}, together with the facts that
			\begin{align*}
				|\delta_\alpha\partial_x^l f(x)|\lesssim |\alpha|^\frac{1}{2}\|f\|_{\dot C^{l+\frac{1}{2}}}\lesssim |\alpha|^\frac{1}{2}(\mathbf{1}_{l\geq 1}t^{-\frac{1}{3}(l-\frac{1}{2})}\|f\|_T+\mathbf{1}_{l=0}||f||_{X_T}),
			\end{align*}
			to obtain
			\begin{align*}
				|\partial_x^n\mathrm{III}_1|&\lesssim \sum_{k+l=n}\|f_1\|_{T,*}\| f_1\|_T ^2(1+\| f_1\|_T )^kt^{-\frac{k}{3}}\int_{\mathbb{R}}\min\{|\alpha|^\frac{3}{4}t^{-\frac{1}{4}},|\alpha|^\frac{1}{4}t^{-\frac{1}{12}}\}|\delta_\alpha\partial_x^l (f_1-f_2)(x)|\frac{d\alpha}{\alpha^2}\\
				&\lesssim \|f_1\|_{T,*}\| f_1\|_T ^2(1+\| f_1\|_T )^nt^{-\frac{n}{3}}(\|f_1-f_2\|_T+t^{-\frac{1}{6}}||f_1-f_2||_{X_T}).\\
				|\partial_x^n\mathrm{III}_2|&\lesssim \sum_{k+l=n}(1+\|(f_1,f_2)\|_T)^{k+4}\|f_1-f_2\|_{T}t^{-\frac{k}{3}}\int_{\mathbb{R}}\min\{|\alpha|^\frac{3}{4}t^{-\frac{1}{4}},|\alpha|^\frac{1}{4}t^{-\frac{1}{12}}\}|\delta_\alpha \partial_x^lf_2(x)|\frac{d\alpha}{\alpha^2}\\
				&\lesssim (1+\|(f_1,f_2)\|_T)^{n+4}\|f_1-f_2\|_{T}t^{-\frac{n}{3}}(\|f_2\|_T+t^{-\frac{1}{6}}||f_2||_{X_T}).
			\end{align*}
			For $\mathrm{III}_3$, there holds
			\begin{align*}
				\|\partial_x^n\mathrm{III}_3\|_{L^\infty}\lesssim t^{-\frac{n}{3}}||f_1-f_2||_{X_T}.
			\end{align*}
			We conclude that 
			\begin{align*}
				\sup_{t\in[0,T]}t^\frac{n}{3}\|\partial_x^n (N_3[f_1]-N_3[f_2])(t)\|_{L^\infty}\lesssim (T^\frac{2}{3}+T^\frac{1}{10})||f_1-f_2)||_{X_T}(1+||(f_1,f_2)||_{X_T})^{n+5}.
			\end{align*}
			Combining this with \eqref{N1} and \eqref{N2}, we get \eqref{eee1}. This completes the proof.
		\end{proof}
		\begin{lemma}\label{GGG}
			Let $\G[f,\phi]$ be as defined in \eqref{defffg}, then for any $n\in\mathbb{N}$, $n\leq m$, and any $T>0$,
			\begin{align*}
				&\sup_{t\in[0,T]}t^{\frac{n+\kappa}{3}}\|(\G[f_1,\phi]-\G[f_2,\phi])(t)\|_{\dot C^{n+\kappa-2}}\\
				&\quad\quad\lesssim \|f_1-f_2\|_T(\|f_1\|_{T,*}+\|f_2-\phi\|_T)(1+\|(f_1,f_2,\phi)\|_T)^{m+5}.
			\end{align*}
		\end{lemma}
		\begin{proof}
			We can write 
			\begin{align*}
				\G[f_1,\phi]-\G[f_2,\phi]=\left(\frac{1}{\langle\partial_xf_2\rangle^3}-\frac{1}{\langle\partial_xf_1\rangle^3}\right)\Lambda^3 f_1+\left(\frac{1}{\langle\partial_x\phi\rangle^3}-\frac{1}{\langle\partial_xf_2\rangle^3}\right)\Lambda^3 (f_1-f_2).
			\end{align*}
			By Lemma \ref{lemcom}, for any $n\in\mathbb{N}$, $n\leq m$,
			
			\begin{align*}
				\left\|\partial_x^n\left(\frac{1}{\langle\partial_xf_2\rangle^3}-\frac{1}{\langle\partial_xf_1\rangle^3}\right)\right\|_{L^\infty}
				&\lesssim \sum_{l=0}^n\|\partial_x(f_1-f_2)\|_{\dot C^{l}}(\|\partial_x(f_1,f_2)\|_{\dot C^1}^{n-l}+\|\partial_x (f_1,f_2)\|_{\dot C^{n-l}}),\\
				\left\|\partial_x^n\left(\frac{1}{\langle\partial_xf_2\rangle^3}-\frac{1}{\langle\partial_xf_1\rangle^3}\right)\right\|_{\dot C^{\kappa-2}}
				&\lesssim \sum_{l=0}^n\left\{\|\partial_x(f_1-f_2)\|_{\dot C^{l+\kappa-2}}(\|\partial_x(f_1,f_2)\|_{\dot C^1}^{n-l}+\|\partial_x (f_1,f_2)\|_{\dot C^{n-l}})\right.\\
				&\quad\quad\quad\left.+\|\partial_x(f_1-f_2)\|_{\dot C^l}(\|\partial_x(f_1,f_2)\|_{\dot C^{\kappa-2}}^\frac{n-l}{\kappa-2}+\|\partial_x (f_1,f_2)\|_{\dot C^{n-l+\kappa-2}})\right\}.
			\end{align*}
			Then by the definition of $\|\cdot\|_T$, we have
			\begin{align*}
				&\left\|\partial_x^n\left(\frac{1}{\langle\partial_xf_2\rangle^3}-\frac{1}{\langle\partial_xf_1\rangle^3}\right)\right\|_{L^\infty}\lesssim t^{-\frac{n}{3}}\|f_1-f_2\|_T(1+\|(f_1,f_2)\|_T)^{m+5},\\
				&\left\|\partial_x^n\left(\frac{1}{\langle\partial_xf_2\rangle^3}-\frac{1}{\langle\partial_xf_1\rangle^3}\right)\right\|_{\dot C^{\kappa-2}}\lesssim t^{-\frac{n+\kappa-2}{3}}\|f_1-f_2\|_T(1+\|(f_1,f_2)\|_T)^{m+5}.
			\end{align*}
			Here we assume without loss of generality that $\kappa-2>\frac{m}{m+5}$ in the last inequality. \\
			On the other hand, we have 
			\begin{align*}
				\|\partial_x^n \Lambda^3 f_1\|_{L^\infty}\lesssim t^{-\frac{n+2}{3}}\|f_1\|_{T,*},\ \ \  \|\partial_x^n \Lambda^3 f_1\|_{\dot C^{\kappa-2}}\lesssim t^{-\frac{n+\kappa}{3}}\|f_1\|_{T,*}.
			\end{align*}
			Combining the above estimates to get 
			\begin{align*}
				\left\|\partial_x^n \left(\left(\frac{1}{\langle\partial_xf_2\rangle^3}-\frac{1}{\langle\partial_xf_1\rangle^3}\right)\Lambda^3 f_1\right)\right\|_{\dot C^{\kappa-2}}\lesssim t^{-\frac{n+\kappa}{3}}\|f_1-f_2\|_T\|
				f_1\|_{T,*}(1+\|(f_1,f_2)\|_T)^{m+5}.
			\end{align*}
			Similarly, we have 
			\begin{align*}
				\left\|\partial_x^n \left(\left(\frac{1}{\langle\partial_x\phi\rangle^3}-\frac{1}{\langle\partial_xf_2\rangle^3}\right)\Lambda^3 (f_1-f_2)\right)\right\|_{\dot C^{\kappa-2}}\lesssim t^{-\frac{n+\kappa}{3}}\|f_2-\phi\|_T\|f_1-f_2\|_T(1+\|(f_2,\phi)\|_T)^{m+5}.
			\end{align*}
			Then we complete the proof.	
		\end{proof}~\vspace{0,3cm}\\
		\subsection{Estimates of nonlinear terms in 2D Peskin equation}
		The following part is devoted to estimate the nonlinear terms of the Peskin equation. In this subsection, we fix $m\in\mathbb{N}$, $\kappa\in(0,1)$ and $1-\kappa\ll 1$. We define 
		\begin{equation}\label{normstar}
			\begin{aligned}
				&\|h\|_{T}=\sup_{t\in[0,T]}(\|h(t)\|_{L^\infty}+t^{m+\kappa}\|h(t)\|_{\dot C^{m+\kappa}}),\\
				&	\|h\|_{T,*}=\sup_{t\in[0,T]}(t^\frac{1}{10}\|h(t)\|_{\dot C^\frac{1}{10}}+t^{m+\kappa}\|h(t)\|_{\dot C^{m+\kappa}}).
			\end{aligned}
		\end{equation}
		Recall the definition of the set $\mathcal{X}^\sigma_{T,\Phi}$ in \eqref{defset}, and the constant $C_{\mathcal{T},n}$ in \eqref{consct}. We have the following results.
		\begin{lemma}
			\label{lemnonpes}Let $\mathcal{N}(X)$ be as defined in \eqref{defnonpes}. For any $T>0$, and any $X,Y,Z\in \mathcal{X}^\sigma_{T,\Phi}$, \begin{align}\label{resn}
				&\sum_{j=0,m}\sup_{t\in[0,T]}t^{j+\kappa}	\|\mathcal{N}(X)(t)\|_{\dot C^{j+\kappa}}
				\lesssim C_{\mathcal{T},m}(1+\mathbf{\Theta}_X(T))^{m+5}\|\partial_x X\|_{T,*}^2(1+\|\partial_x X\|_{T})^{m+5},
			\end{align}
			and 
			\begin{equation}\label{rendiff}
				\begin{aligned}
					&\sum_{j=0,m}\sup_{t\in[0,T]} t^{j+\kappa}	\|\mathcal{N}(Y)(t)-\mathcal{N}(Z)(t)\|_{\dot C^{j+\kappa}}\\
					&\quad\lesssim C_{\mathcal{T},m}(1+\mathbf{\Theta}_{(Y,Z)}(T))^{m+5} \|\partial_x(Y-Z)\|_{T}(\|\partial_x( Y,Z)\|_{T,*}+T^{\frac{1}{2}}\|\partial_x( Y,Z)\|_{T})(1+\|\partial_x (Y,Z)\|_T)^{5(m+1)}.
				\end{aligned}
			\end{equation}
		\end{lemma}
		\begin{lemma}\label{lempesR}
			Let  $\mathcal{M}(\partial_xX)$ be as defined in \eqref{defpeR} with $a\in(0,1)$ and $1-a\ll 1$. Then for any $T>0$, any $Y,Z\in C((0,T);C^{1+\kappa})$, there holds 
			\begin{align*}
				&\sum_{j=0,m}\sup_{t\in[0,T]} t^{j+\kappa}\|(\mathcal{M}(\partial_xY)-\mathcal{M}(\partial_xZ))(t)\|_{\dot C^{j+\kappa-a}}\\
				&\quad\quad\quad\quad\lesssim C_{\mathcal{T},m}\|\partial_x(Y-Z)\|_T\|\partial_x(Y, Z)\|_{T,*}(1+\|\partial_x(Y,Z)\|_T)^{m+5}.
			\end{align*}
		\end{lemma}
		Recalling the definition \eqref{notapes},	we first prove the following lemmas. 
		\begin{lemma}\cite[Lemma 2.4]{KN}\label{lempesibp}
			For any function $f,g:{\mathbb{R}}\rightarrow\mathbb{R}$, denote $\tilde g(\alpha)={\Delta}_\alpha g(0)$. Then for any $\sigma\in(0,1)$ and $0<\varepsilon< 10^{-3}\min\{1-\sigma,\sigma\}$, there holds
			$$
			\left|\int_{\mathbb{R}} f(\alpha)(\partial_\alpha \tilde g)(\alpha)d\alpha\right|\lesssim \| f\|_{\dot C^{\sigma+\varepsilon}}^\frac{1}{2} \| f\|_{\dot C^{\sigma-\varepsilon}}^\frac{1}{2}\|g\|_{\dot C^{1-\sigma+\varepsilon}}^\frac{1}{2}\|g\|_{\dot C^{1-\sigma-\varepsilon}}^\frac{1}{2}.$$
		\end{lemma}
		\begin{lemma}\label{dede}
			For any function $f:\mathbb{S}\to\mathbb{R}$, then for any $0<\gamma_1\leq\gamma<1$,
			\begin{align}
				\label{es2}
				&	
				\sup_x\sup_{\alpha\neq h} \frac{|\tilde \Delta_\alpha f(x)-\tilde \Delta_h f(x)||\alpha|^{\gamma-\gamma_1}}{|\alpha-h|^\gamma}\lesssim {\|f'\|_{\dot C^{\gamma_1}}}.
			\end{align}
		\end{lemma}
		\begin{proof} 
			Consider $x=0$ without loss of generality, and denote $\tilde f (\alpha)=\tilde \Delta_\alpha f(0)$.
			Thanks to the periodicity, it is enough to consider $\alpha\in (0,\pi)$ and $h\in[\alpha/2,\alpha/2+2\pi]$. 
			We have 
			\begin{align*}
				\tilde f(\alpha)-\tilde f(h)=&\left(\frac{1}{2}\cot\left(\frac{h}{2}\right)-\frac{1}{h}\right)\int_0^{-h}f'(\omega)d\omega-\left(\frac{1}{2}\cot\left(\frac{\alpha}{2}\right)-\frac{1}{\alpha}\right)\int_0^{-\alpha}f'(\omega)d\omega\\
				&+\frac{1}{h}\int_0^{-h}f'(\omega)d\omega-\frac{1}{\alpha}\int_0^{-\alpha}f'(\omega)d\omega\\
				=&\left(\frac{1}{2}\cot\left(\frac{\alpha}{2}\right)-\frac{1}{\alpha}\right)\int_{-\alpha}^{-h}f'(\omega)d\omega+\left(\frac{1}{2}\cot\left(\frac{h}{2}\right)-\frac{1}{h}-\left(\frac{1}{2}\cot\left(\frac{\alpha}{2}\right)-\frac{1}{\alpha}\right)\right)\int_0^{-h}f'(\omega)d\omega\\
				&+\frac{1}{h}\int_{-\alpha}^{-h}(f'(\omega)-f'(-\alpha))d\omega+\frac{\alpha-h}{\alpha h}\int_0^{-\alpha}(f'(\omega)-f'(-\alpha))d\omega.
			\end{align*}
			Here, with a slight abuse of notation, we denote $\int_x^yg(\omega)d\omega=-\int_y^xg(\omega)d\omega$ if $y<x$.
			Note that $$\left|\frac{1}{2}\cot\left(\frac{\alpha}{2}\right)-\frac{1}{\alpha}\right|\lesssim |\alpha|,\ \ \ \ \  \left|\frac{1}{2}\cot\left(\frac{\alpha}{2}\right)-\frac{1}{\alpha}-\frac{1}{2}\cot\left(\frac{h}{2}\right)+\frac{1}{h}\right|\lesssim |\alpha-h|.$$ Hence 
			\begin{align*}
				|\tilde f(\alpha)-\tilde f(h)|\lesssim \left(|h|+|\alpha|+\frac{|h-\alpha|^{\gamma_1}}{|h|}+\frac{|\alpha|^{\gamma_1}}{|h|}\right) |h-\alpha|\|f'\|_{\dot C^{\gamma_1}}.
			\end{align*}
			We obtain 
			\begin{align*}
				\frac{|\tilde \Delta_\alpha f-\tilde \Delta_h f|}{|\alpha-h|^\gamma}\lesssim \left(|h|+|\alpha|+\frac{|h-\alpha|^{\gamma_1}}{|h|}+\frac{|\alpha|^{\gamma_1}}{|h|}\right) |h-\alpha|^{1-\gamma}\|f'\|_{\dot C^{\gamma_1}}\lesssim \frac{\|f'\|_{\dot C^{\gamma_1}}}{|h|^{\gamma-\gamma_1}}\lesssim \frac{\|f'\|_{\dot C^{\gamma_1}}}{|\alpha|^{\gamma-\gamma_1}}.
			\end{align*}
			This completes the proof.

		\end{proof}

		For $f(t,x,\alpha): [0,T]\times \mathbb{S}\times \mathbb{S}\to \mathbb{R}$, we define 
		\begin{equation}\label{defTn}
			\begin{aligned}
				&	[f_1]_{T,*}:=	\sup_{\vartheta\in[\frac{1}{2},\kappa]}\sup_{t\in[0,T]}(t^{\vartheta-\frac{1}{4}}\|f_1(t)\|_{\vartheta,\frac{1}{4}}+t^{m+\vartheta-\frac{1}{4}}\|\partial_x^mf_1(t)\|_{\vartheta,\frac{1}{4}})\\
				&\quad\quad\quad\quad\quad+\sup_{t\in[0,T]}\sup_x(t^\frac{1}{10}\|f_1(t,x,\cdot)\|_{\dot C^\frac{1}{10}}+t^{m+\kappa}\|\partial_x^{m}f_1(t,x,\cdot)\|_{\dot C^\kappa}),\\
				&[f_1]_{T}:=\sup_{t\in[0,T]}\sup_\alpha(\|f_1(t,\cdot,\alpha)\|_{L^\infty}+t^{m+\kappa}\|f_1(t,\cdot,\alpha)\|_{\dot C^{m+\kappa}})+	[f_1]_{T,*},
			\end{aligned}
		\end{equation}
		where 
		\begin{align}
			\label{defnorab}
			\|f_1(t)\|_{a,b}:=	\sup_{\alpha,z\in\mathbb{S}}|\alpha|^b\left(\|f_1(t,\cdot,\alpha)\|_{\dot C^a}+\frac{\|f_1(t,\cdot,\alpha )-f_1(t,\cdot,\alpha -z)\|_{L^\infty}}{|z|^a}\right),\ \ \ \text{for}\ a,b\in[0,1].
		\end{align}
		We have the following lemma.
		\begin{lemma}\label{lemeee} Fix $m\in\mathbb{N}$, $\kappa\in(0,1)$ and $1-\kappa\ll 1$. 
			Let $f_1:[0,T]\times \mathbb{S}\times \mathbb{S}\to \mathbb{R}$, and  $f_2,f_3: [0,T]\times  \mathbb{S}\to \mathbb{R} $ satisfy
			\begin{equation*}
				\begin{aligned}
					& [f_1]_{T}+\|\partial_xf_2\|_{T}+ \|f_3\|_{T} <\infty
				\end{aligned}
			\end{equation*}
			with $0<T<1$, where  $\|\cdot \|_T$ and $\|\cdot\|_{T,*}$ are defined in \eqref{normstar}, and $[\cdot]_T$ is defined in \eqref{defTn}.
			Then for
			\begin{align*}
				f(t,x)=\int_{\mathbb{S}} f_1(t,x,\alpha)\tilde{\E}^\alpha f_2(t,x) f_3(t,x-\alpha)\frac{d\alpha}{\tilde \alpha},
			\end{align*}
			where $\tilde \alpha$ and $\tilde{\E}^\alpha$ are  defined in \eqref{notapes},
			there holds
			\begin{equation}\label{res1}
				\begin{aligned}
					\sup_{t\in[0,T]}t^{m+\kappa }\|f(t,\cdot)\|_{\dot C^{m+\kappa }}
					\lesssim \|\partial_x f_2\|_{T,*}([f_1]_{T,*}\|f_3\|_{T}+[f_1]_T\|f_3\|_{T,*}+T^\frac{1}{2}[f_1]_T\|f_3\|_{T}).
				\end{aligned}
			\end{equation}
			If we further have $\int_{\mathbb{S}} f_1(t,x,\alpha)\tilde{\E}^\alpha f_2(t,x)\frac{d\alpha}{\tilde \alpha}=0$, then we can only keep the second term, which is
			\begin{align}\label{res2}
				\sup_{t\in[0,T]}t^{m+\kappa }\|f(t,\cdot)\|_{\dot C^{m+\kappa }}\lesssim [f_1]_T\|\partial_x f_2\|_{T,*}\|f_3\|_{T,*}.
			\end{align}
		\end{lemma}
		\begin{proof}
			We first prove \eqref{res1}. 
			For simplicity, we drop the time variable $t$ in this proof, and denote $f_i^{m_i}=\partial_x ^{m_i}f_i$ for $i=1,2,3$, $m_i\in\mathbb{N}$. We have 
			\begin{align*}
				\partial_x^m f(x)=\sum_{m_1+m_2+m_3=m}\int_{\mathbb{S}}f_1^{m_1}(x,\alpha)\tilde{\E}^\alpha f_2^{m_2}(x) f_3^{m_3}(x-\alpha)\frac{d\alpha}{\tilde \alpha}.
			\end{align*}
			In the following,  we drop the summation for $m_1+m_2+m_3=m$ with a slight abuse of notation. \\
			For any $\beta \neq 0$, we can write 
			\begin{align*}
				\delta_\beta\partial_x^m f(x)=&\int_{\mathbb{S}}\delta_\beta f_1^{m_1}(x,\alpha)\tilde{\E}^\alpha f_2^{m_2}(x) f_3^{m_3}(x-\alpha)\frac{d\alpha}{\tilde \alpha}\\
				&+\int_{\mathbb{S}} f_1^{m_1}(x-\beta,\alpha)\delta_\beta\tilde{\E}^\alpha f_2^{m_2}(x) f_3^{m_3}(x-\alpha)\frac{d\alpha}{\tilde \alpha}\\
				&+\int_{\mathbb{S}} f_1^{m_1}(x-\beta,\alpha) \tilde{\E}^\alpha f_2^{m_2}(x-\beta)\delta_\beta f_3^{m_3}(x-\alpha)\frac{d\alpha}{\tilde \alpha}\\
				:=&P_1+P_2+P_3.
			\end{align*}
			By the definition of $\|\cdot\|_{a,b}$ in \eqref{defnorab}, we have
			\begin{align*}
				&|\alpha|^\frac{1}{4}|\delta_\beta f_1^{m_1}(x,\alpha)|\lesssim |\beta|^\kappa\|f_1^{m_1}\|_{\kappa,\frac{1}{4}}.
			\end{align*}
			Hence it is easy to check that 
			\begin{equation}\label{1P1}
				\begin{aligned}
					|P_1|&\lesssim |\beta|^\kappa t^{-(m_1+\kappa-\frac{1}{4})}[f_1]_{T,*}\|f_3^{m_3}\|_{L^\infty} \int_{\mathbb{S}}|\tilde{\E}^\alpha f_2^{m_2}|\frac{d\alpha}{|\alpha|^\frac{5}{4}}\\
					&\lesssim |\beta|^\kappa t^{-(m+\kappa)}[f_1]_{T,*}\|\partial_xf_2\|_{T,*}\|f_3\|_{T}.
				\end{aligned}
			\end{equation}
			Then we deal with $P_{2}$. Observe that 
			\begin{align}\label{ealp}
				\tilde{\E}^\alpha f_2^{m_2}(x)=\alpha\partial_\alpha (\Delta_\alpha f_2^{m_2})(x)+\delta_\alpha f_2^{m_2}(x)\left(\frac{1}{\alpha}-\frac{1}{\tilde \alpha}\right).
			\end{align}
			We observe that, for any $h(\alpha):\mathbb{S}\to\mathbb{R}$, there holds
			\begin{align}\label{1perio}
				\int_{\mathbb{R}}h(\alpha)\frac{d\alpha}{\alpha}=\int_{-\pi}^\pi h(\alpha)\frac{d\alpha}{\tilde\alpha}.
			\end{align}
			Hence we can write 
			\begin{align*}
				P_{2}=
				&\int_{\mathbb{R}} f_1^{m_1}(x-\beta,\alpha)\delta_\beta\tilde{\E}^\alpha f_2^{m_2}(x) f_3^{m_3}(x-\alpha)\frac{d\alpha}{ \alpha}\\	=&\int_{\mathbb{R}} f_1^{m_1}(x-\beta,\alpha)\partial_\alpha (\Delta_\alpha \delta_\beta f_2^{m_2})(x) f_3^{m_3}(x-\alpha){d\alpha}\\
				&+\int_{\mathbb{R}} f_1^{m_1}(x-\beta,\alpha)(\delta_\alpha\delta_\beta f_2^{m_2})(x)  f_3^{m_3}(x-\alpha)\left(\frac{1}{\alpha}-\frac{1}{\tilde\alpha}\right)\frac{d\alpha}{\alpha}\\
				:=&P_{2,1}+P_{2,2}.
			\end{align*}
			Since $|\frac{1}{\alpha}(\frac{1}{\alpha}-\frac{1}{\tilde\alpha})|\lesssim\min\{1,|\alpha|^{-2}\}$, we remark that $P_{2,1}$ is the main term, and $P_{2,2}$ is the remainder error term.\\
			Applying Lemma \ref{lempesibp} with $f(\alpha)=f_1^{m_1}(x,\alpha) f_3^{m_3}(x-\alpha)$, $g=\delta_\beta f_2^{m_2}$, one has 
			\begin{align*}
				|P_{2,1}|\lesssim\prod_{+,-} \left\{ \|\delta_\beta f_2^{m_2}\|_{\dot C^{\frac{1}{4}\pm\varepsilon}}^\frac{1}{2}\sup_x\|f_1^{m_1}(x,\alpha) f_3^{m_3}(x-\alpha)\|_{\dot C_\alpha^{\frac{3}{4}\pm\varepsilon}}^\frac{1}{2}\right\}.
			\end{align*}
			Here we denote $\dot C^{b}_\alpha$ the H\"{o}lder semi-norm in variable $\alpha$. By Lemma \ref{Hol}, it is easy to check  that for $\gamma=\frac{3}{4}\pm\varepsilon$, 
			\begin{align*}
				\|f_1^{m_1}(x,\alpha) f_3^{m_3}(x-\alpha)\|_{L^\infty_x\dot C_\alpha^{\gamma}}
				&\lesssim \|f_1^{m_1}\|_{L^\infty_x\dot C^\gamma_\alpha}\|f_3^{m_3}\|_{L^\infty}+\|f_1^{m_1}\|_{L^\infty_{x,\alpha}}\|f_3^{m_3}\|_{\dot C^\gamma}\\
				&\lesssim t^{-(m_1+m_3)-\gamma}([f_1]_{T,*}\|f_3\|_{T}+[f_1]_T\|f_3\|_{T,*}).
			\end{align*}
			Hence we obtain that 
			\begin{align}
				|P_{2,1}|&\lesssim |\beta|^\kappa  \prod_{+,-} \left\{t^{-(m_1+m_3)-\frac{3}{4}\mp\varepsilon}\|\partial_x f_2^{m_2}\|_{\dot C^{\kappa-\frac{3}{4}\pm\varepsilon}} \right\}^\frac{1}{2}([f_1]_{T,*}\|f_3\|_{T}+[f_1]_T\|f_3\|_{T,*})\nonumber\\
				&\lesssim t^{-\kappa-m} |\beta|^\kappa \|\partial_x f_2\|_{T,*}([f_1]_{T,*}\|f_3\|_{T}+[f_1]_T\|f_3\|_{T,*}),\label{1lP21}
			\end{align}
			where $m=m_1+m_2+m_3$. 
			For $P_{2,2}$, by \eqref{1perio} we have 
			\begin{equation}\label{aaaaa}
				\begin{aligned}
					\left|	\int_{\mathbb{R}}h(\alpha)\left(\frac{1}{\alpha}-\frac{1}{\tilde\alpha}\right)\frac{d\alpha}{\alpha}\right|&=\left|\int_{\mathbb{R}} h(\alpha)\left(\frac{1}{\alpha^2}-\frac{\mathbf{1}_{|\alpha|\leq \pi}}{\tilde\alpha^2}\right)d\alpha\right|\\
					&\lesssim \int_{\mathbb{R}}\left| h(\alpha)\right|\min\{1,|\alpha|^{-2}\}d\alpha.
				\end{aligned}
			\end{equation}
			Then we obtain that 
			\begin{align*}
				|P_{2,2}|&\lesssim \|f_1^{m_1}\|_{L^\infty}\| f_3^{m_3}\|_{L^\infty}\int _{\mathbb{R}} \|\delta_\alpha\delta_\beta f_2^{m_2}\|_{L^\infty}\min\{1,|\alpha|^{-2}\}d\alpha \\
				&\lesssim |\beta|^\kappa\|f_1^{m_1}\|_{L^\infty}\| f_3^{m_3}\|_{L^\infty}\|\partial_xf_2^{m_2}\|_{\dot C^{\kappa-\frac{1}{2}}}\\
				&\lesssim t^{-\kappa-m+\frac{1}{2}} |\beta|^\kappa [f_1]_T\|\partial_x f_2\|_{T,*}\|f_3\|_{T}.
			\end{align*}
			Combining this with \eqref{1lP21} to obtain that for any $t\in[0,T]$,
			\begin{align}\label{1P2}
				t^{m+\kappa } |\beta|^{-\kappa} |P_{2}|\lesssim \|\partial_x f_2\|_{T,*}([f_1]_{T,*}\|f_3\|_{T}+[f_1]_T\|f_3\|_{T,*}+T^\frac{1}{2}[f_1]_T\|f_3\|_{T}).
			\end{align}
			Finally, we deal with $P_3$. If $m_1+m_2\neq 0$, \textit{i.e.} $m_3<m$, we directly have
			\begin{align*}
				|P_3|&\lesssim \|f_1^{m_1}\|_{L^\infty}\|\partial_x f_2\|_{\dot C^{m_2+\frac{3}{4}}}\|\delta_\beta f_3\|_{\dot C^{m_3}}\\
				&\lesssim |\beta|^\kappa t^{-(m+\kappa)}[f_1]_T\|\partial_xf_2\|_{T,*}\|f_3\|_{T,*}.
			\end{align*}
			If $m_1+m_2=0$, we observe that 
			for any function $h_1,h_2$, by a change of variable, 
			\begin{align*}
				\int_{\mathbb{S}} h_1(\alpha)\delta_\beta h_2(x-\alpha)d\alpha&=\int_{\mathbb{S}} \left(h_1(\alpha)-h_1(\alpha-\beta)\right) h_2(x-\alpha)d\alpha\\
				&=-\int_{\mathbb{S}} \left(h_1(\alpha)-h_1(\alpha-\beta)\right) \delta_\alpha h_2(x)d\alpha.
			\end{align*}
			Applying this to $P_3$,
			\begin{align*}
				P_3&=\int_{\mathbb{S}}\left( \frac{f_1(x-\beta,\alpha)\tilde{\E}^\alpha f_2(x-\beta)}{\alpha}-\frac{f_1(x-\beta,\alpha-\beta)\tilde{\E}^{\alpha-\beta} f_2(x-\beta)}{\alpha-\beta}\right)\delta_\alpha f_3^m(x)d\alpha\\
				&=\int_{\mathbb{S}}{ (f_1(x-\beta,\alpha)-f_1(x-\beta,\alpha-\beta))\tilde{\E}^\alpha f_2(x-\beta)}\delta_\alpha f_3^m(x) \frac{d\alpha}{\alpha}\\
				&\ \ \ \ +\int_{\mathbb{S}} f_1(x-\beta,\alpha-\beta)(\tilde{\E}^\alpha f_2-\tilde{\E}^{\alpha-\beta}f_2)(x-\beta)\delta_\alpha f_3^m(x)\frac{d\alpha}{\alpha}\\
				&\ \ \ \ +\int_{\mathbb{S}} f_1(x-\beta,\alpha-\beta)\tilde{\E}^{\alpha-\beta}f_2(x-\beta)\delta_\alpha f_3^m(x)\left(\frac{1 }{\alpha}-\frac{1}{\alpha-\beta}\right)d\alpha \\
				&:=P_{3,1}+P_{3,2}+P_{3,3}.
			\end{align*}
			For $P_{3,1}$, by \eqref{defnorab} and  Lemma \ref{douint} we have 
			we have 
			\begin{align}
				&	|f_1(x-\beta,\alpha)-f_1(x-\beta,\alpha-\beta)|\lesssim|\beta|^\kappa  \min\left\{\sup_x\|f_1(x,\cdot)\|_{\dot C^\kappa},|\alpha|^{-\frac{1}{4}}\|f_1\|_{\kappa,\frac{1}{4}}\right\},\nonumber\\
				&	|\tilde{\E}^\alpha f_2(x,\beta)|\lesssim |\alpha|^a \|\partial_x f_2\|_{\dot C^{a}},\quad\quad |\delta_\alpha f_3^m(x)|\lesssim |\alpha |^a \|f_3\|_{\dot C^{m+a}}, \ \ \forall a\in(0,1).\label{eal111}
			\end{align}
			Hence  we get
			\begin{equation}\label{1lP31}
				\begin{aligned}
					|P_{3,1}|	&\lesssim 	\int_{\mathbb{S}} | (f_1(x-\beta,\alpha)-f_1(x-\beta,\alpha-\beta))\tilde{\E}^\alpha f_2(x-\beta)\delta_\alpha f_3^m(x) |\frac{d\alpha}{|\alpha|}\\
					&\lesssim |\beta|^\kappa\int_{\mathbb{S}} \min\left\{\sup_x\|f_1(x,\cdot)\|_{\dot C^\kappa},|\alpha|^{-\frac{1}{4}}\|f_1\|_{\kappa,\frac{1}{4}}\right\}\frac{d\alpha}{|\alpha|^\frac{4}{5}}  \|\partial_x f_2\|_{\dot C^{\frac{1}{10}}}\|f_3\|_{\dot C^{\frac{1}{10}+m}}\\
					&\lesssim |\beta|^\kappa t^{-(m+\kappa)}[f_1]_{T,*}\|\partial_xf_2\|_{T,*}\|f_3\|_{T,*}.
				\end{aligned}
			\end{equation}
			For $P_{3,2}$, observe that
			\begin{align*}
				({\tilde{\E}^\alpha  f_2}-{\tilde{\E}^{\alpha-\beta}  f_2})(x)&=\partial_x(\delta_{-\beta}f_2)(x-\alpha)-(\tilde \Delta_\alpha f_2(x)-\tilde \Delta_{\alpha-\beta} f_2(x))\\
				&=\partial_\alpha (\alpha \Delta_\alpha \delta_{-\beta}f_2)(x)-(\tilde \Delta_\alpha f_2(x)-\tilde \Delta_{\alpha-\beta} f_2(x))\\
				&=\alpha\partial_\alpha ( \Delta_\alpha \delta_{-\beta}f_2)(x)+\Delta_\alpha \delta_{-\beta}f_2(x)-(\tilde \Delta_\alpha f_2(x)-\tilde \Delta_{\alpha-\beta} f_2(x)).
			\end{align*}
			We further decompose $P_{3,2}$ into 
			\begin{align*}
				P_{3,2}=&\int_{\mathbb{S}}  f_1(x-\beta,\alpha-\beta)\partial_\alpha ( \Delta_\alpha \delta_{-\beta}f_2)(x) \delta_\alpha f_3^m(x){d\alpha}\\
				&+\int_{\mathbb{S}}  f_1(x-\beta,\alpha-\beta)( \Delta_\alpha \delta_{-\beta}f_2)(x) \delta_\alpha f_3^m(x)\frac{d\alpha}{\alpha}\\
				&+\int_{\mathbb{S}}  f_1(x-\beta,\alpha-\beta)(\tilde \Delta_\alpha f_2-\tilde \Delta_{\alpha-\beta} f_2)(x) \delta_\alpha f_3^m(x)\frac{d\alpha}{\alpha}\\
				=&P_{3,2,1}+P_{3,2,2}+P_{3,2,3}.
			\end{align*}
			We estimate $P_{3,2,1}$ similarly as $P_{2,1}$, which yields
			\begin{equation}\label{1lP321}
				\begin{aligned}
					&|P_{3,2,1}|
					\lesssim  t^{-\kappa-m} |\beta|^\kappa \|\partial_x f_2\|_{T,*}([f_1]_{T,*}\|f_3\|_{T}+[f_1]_{T}\|f_3\|_{T,*}),.
				\end{aligned}
			\end{equation}	
			Moreover, we have 
			\begin{align*}
				|P_{3,2,2}|&\lesssim [f_1]_{T}\int_{\mathbb{S}} |\delta_\alpha\delta_{-\beta} f_2(x)\delta_\alpha f_3 ^m(x)|\frac{d\alpha}{|\alpha|^2}.
			\end{align*}
			By Lemma \ref{douint}, it is easy to check that 
			\begin{align*}
				&\|\delta_\alpha\delta_{-\beta} f_2\|_{L^\infty}\lesssim |\beta|^\kappa \min\left\{|\alpha|^{\frac{1}{4}}\|\partial_x f_2\|_{\dot C^{\kappa-\frac{3}{4}}},|\alpha|^\frac{3}{4}\|\partial_x f_2\|_{\dot C^{\kappa-\frac{1}{4}}}\right\},\\
				&\|\delta_\alpha f_3^m\|_{L^\infty}\lesssim \min\left\{|\alpha|^\frac{1}{4}\|f_3\|_{\dot C^{\frac{1}{4}+m}}, |\alpha|^\frac{3}{4}\|f_3\|_{\dot C^{\frac{3}{4}+m}}\right\}.
			\end{align*}
			Hence we obtain 
			\begin{equation}\label{1lP322}
				\begin{aligned}
					&|P_{3,2,2}|\\
					&\lesssim |\beta|^\kappa \|f_1\|_{L^\infty} \int_{\mathbb{S}}  \min\left\{\|\partial_x f_2\|_{\dot C^{\kappa-\frac{3}{4}}},|\alpha|^{\frac{1}{2}}\|\partial_x f_2\|_{\dot C^{\kappa-\frac{1}{4}}}\right\}\min\{\|f_3\|_{\dot C^{\frac{1}{4}+m}}, |\alpha|^\frac{1}{2}\|f_3\|_{\dot C^{\frac{3}{4}+m}}\}\frac{d\alpha}{|\alpha|^\frac{3}{2}}\\
					&\lesssim |\beta|^\kappa \|f_1\|_{L^\infty} (\|\partial_x f_2\|_{\dot C^{\kappa-\frac{3}{4}}}\|\partial_x f_2\|_{\dot C^{\kappa-\frac{1}{4}}}\|f_3\|_{\dot C^{\frac{1}{4}+m}}\|f_3\|_{\dot C^{\frac{3}{4}+m}})^\frac{1}{2}.
				\end{aligned}
			\end{equation}
			Finally, similar to $P_{3,1}$, by Lemma \ref{dede},
			\begin{align*}
				|\alpha|^\frac{1}{4}\frac{\|\tilde \Delta_\alpha f_2-\tilde \Delta_{\alpha-\beta} f_2\|_{L^\infty}}{|\beta|^\kappa}\lesssim  \|\partial_x f_2\|_{\dot C^{\kappa-\frac{1}{4}}}.
			\end{align*}
			Hence 
			\begin{align*}
				|P_{3,2,3}|&\lesssim \int_{\mathbb{S}}| f_1(x-\beta,\alpha-\beta)(\tilde \Delta_\alpha f_2-\tilde \Delta_{\alpha-\beta} f_2)(x) \delta_\alpha f_3^m(x)|\frac{d\alpha}{|\tilde \alpha|}\\&\lesssim |\beta|^\kappa \|f_1\|_{L^\infty_{x,\alpha}} \|\partial_x f_2\|_{\dot C^{\kappa-\frac{1}{4}}}\int_{\mathbb{S}}|\delta_\alpha f_3^m|\frac{d\alpha}{|\tilde\alpha|^\frac{5}{4}}\\
				&\lesssim t^{-{m+\kappa }} |\beta|^\kappa [f_1]_T\|\partial_x f_2\|_{T,*}\|f_3\|_{T,*}.
			\end{align*}
			Combining this with \eqref{1lP321} and \eqref{1lP322} to deduce that 
			\begin{align}\label{1lP32}
				t^{m+\kappa}|\beta|^{-\kappa}|P_{3,2}|\lesssim [f_1]_T\|\partial_x f_2\|_{T,*}\|f_3\|_{T,*}.
			\end{align}
			Then we estimate $$P_{3,3}=-\int_{\mathbb{S}}  f_1(x-\beta,\alpha-\beta)\tilde\E^{\alpha-\beta}f_2(x-\beta)\delta_\alpha f_3^m(x)\left(\frac{1 }{\alpha}-\frac{1}{\alpha-\beta}\right)d\alpha. $$
			By \eqref{eal111}, one has
			\begin{align*}
				|P_{3,3}|\lesssim \|f_1\|_{L^\infty_{x,\alpha}}\|\partial_x f_2\|_{\dot C^\frac{\kappa}{2}} \|f_3\|_{\dot C^{m+\frac{\kappa}{2}}}\int_{\mathbb{R}}|\alpha-\beta|^\frac{\kappa}{2}|\alpha|^\frac{\kappa}{2}\left|\frac{1 }{\alpha}-\frac{1}{\alpha-\beta}\right| d\alpha.
			\end{align*}
			Note that 
			\begin{align*}
				\int_{\mathbb{R}}|\alpha-\beta|^\frac{\kappa}{2}|\alpha|^\frac{\kappa}{2}\left|\frac{1 }{\alpha}-\frac{1}{\alpha-\beta}\right| d\alpha\lesssim |\beta|^\kappa\int_{\mathbb{R}}|\alpha-1|^\frac{\kappa}{2}|\alpha|^\frac{\kappa}{2}\left|\frac{1 }{\alpha}-\frac{1}{\alpha-1}\right| d\alpha\lesssim |\beta|^\kappa.
			\end{align*}
			This yields that 
			\begin{align}\label{1lP33}
				|P_{3,3}|\lesssim |\beta|^\kappa \|f_1\|_{L^\infty_{x,\alpha}}\|\partial_x f_2\|_{\dot C^\frac{\kappa}{2}} \|f_3\|_{\dot C^{m+\frac{\kappa}{2}}}\lesssim|\beta|^\kappa t^{-(m+\kappa)}[f_1]_T\|\partial_x f_2\|_{T,*}\|f_3\|_{T,*}.
			\end{align}
			Combining \eqref{1lP31},  \eqref{1lP32}, and \eqref{1lP33} to obtain 
			\begin{align}\label{1lP3}
				t^{m+\kappa}|\beta |^{-\kappa}|P_{3}|\lesssim [f_1]_T\|\partial_x f_2\|_{T,*}\|f_3\|_{T,*}.
			\end{align}
			Hence we obtain from \eqref{1P1}, \eqref{1P2} and \eqref{1lP3} that 
			\begin{align*}
				&	\|\delta_\beta \partial_x f^m\|_{L^\infty}\lesssim 	|P_1|+|P_2|+|P_3|\\
				&\quad\quad\quad\lesssim |\beta|^\kappa t^{-(m+\kappa)}\|\partial_x f_2\|_{T,*}([f_1]_{T,*}\|f_3\|_{T}+[f_1]_{T}\|f_3\|_{T,*}+T^\frac{1}{2}[f_1]_T\|f_3\|_{T}).
			\end{align*}
			This completes the proof of \eqref{res1}.\vspace{0.3cm}\\
			Then we prove \eqref{res2}. If $$
			\int_{\mathbb{S}} f_1(t,x,\alpha)\tilde{\E}^\alpha f_2(t,x)\frac{d\alpha}{\tilde \alpha}=0
			$$
			holds, then 
			\begin{align*}
				f(t,x)=-\int_{\mathbb{S}}f_1(t,x,\alpha)\tilde{\E}^\alpha f_2(t,x)\delta_\alpha f_3(t,x)\frac{d\alpha}{\tilde \alpha}.
			\end{align*}
			We note the finite difference operator in $f_3$ allows $f_3$ to share derivative, hence we avoid the endpoint norm $\|f_3\|_{L^\infty}$ in the estimate \eqref{res2}. The ideas of the proof is similar to that of \eqref{res1}. For $m\in\mathbb{N}$, we write 
			\begin{align*}
				\partial_x^m f(x)=-\sum_{m_1+m_2+m_3=m}\int_{\mathbb{S}}f_1^{m_1}(x,\alpha)\tilde{\E}^\alpha f_2^{m_2}(x)\delta_\alpha f_3^{m_3}(x)\frac{d\alpha}{\tilde \alpha}.
			\end{align*}
			We drop the summation for $m_1+m_2+m_3=m$ with a slight abuse of notation. \\
			For any $\beta \neq 0$, we can write 
			\begin{align*}
				\delta_\beta\partial_x^m f(x)=&\int_{\mathbb{S}}\delta_\beta f_1^{m_1}(x,\alpha)\tilde{\E}^\alpha f_2^{m_2}(x)\delta_\alpha f_3^{m_3}(x)\frac{d\alpha}{\tilde \alpha}\\
				&+\int_{\mathbb{S}} f_1^{m_1}(x-\beta,\alpha)\delta_\beta \tilde{\E}^\alpha f_2^{m_2}(x)\delta_\alpha f_3^{m_3}(x)\frac{d\alpha}{\tilde \alpha}\\
				&+\int_{\mathbb{S}} f_1^{m_1}(x-\beta,\alpha) \tilde{\E}^\alpha f_2^{m_2}(x-\beta)\delta_\beta\delta_\alpha f_3^{m_3}(x)\frac{d\alpha}{\tilde \alpha}\\
				:=&\tilde P_1+\tilde P_2+\tilde P_3.
			\end{align*}
			Note that 
			\begin{align*}
				&|\alpha|^\frac{1}{4}|\delta_\beta f_1^{m_1}(x,\alpha)|\lesssim |\beta|^\kappa\|f_1^{m_1}\|_{\kappa,\frac{1}{4}}.
			\end{align*}
			Hence it is easy to check that 
			\begin{equation}\label{P1}
				\begin{aligned}
					|\tilde P_1|&\lesssim |\beta|^\kappa t^{-(m_1+\kappa-\frac{1}{4})}[f_1]_{T,*} \|f_2^{m_2+1}\|_{\dot C^\frac{1}{8}}\int_{\mathbb{S}}|\delta_\alpha f_3^{m_3}|\frac{d\alpha}{|\alpha|^\frac{9}{8}}\\
					&\lesssim |\beta|^\kappa t^{-(m+\kappa)}[f_1]_{T,*}\|\partial_xf_2\|_{T,*}\|f_3\|_{T,*}.
				\end{aligned}
			\end{equation}
			Then we deal with $P_{2}$. By \eqref{ealp}, we can write 
			\begin{align*}
				\tilde P_{2}=&\int f_1^{m_1}(x-\beta,\alpha)\partial_\alpha (\Delta_\alpha \delta_\beta f_2^{m_2})(x)\delta_\alpha f_3^{m_3}(x){d\alpha}\\
				&+\int f_1^{m_1}(x-\beta,\alpha)(\delta_\alpha\delta_\beta f_2^{m_2})(x) \delta_\alpha f_3^{m_3}(x)\left(\frac{1}{\alpha}-\frac{1}{\tilde\alpha}\right)\frac{d\alpha}{\alpha}\\
				:=&\tilde P_{2,1}+\tilde P_{2,2}.
			\end{align*}
			By Lemma \ref{lempesibp}, one has 
			\begin{align*}
				|\tilde P_{2,1}|\lesssim\prod_{+,-} \left\{ \|\delta_\beta f_2^{m_2}\|_{\dot C^{\frac{1}{4}\pm\varepsilon}}^\frac{1}{2}\sup_x\|f_1^{m_1}(x,\alpha)\delta_{\alpha} f_3^{m_3}(\tau)\|_{\dot C_\alpha^{\frac{3}{4}\pm\varepsilon}}^\frac{1}{2}\right\}.
			\end{align*}
			Here we use $\dot C^{a}_\alpha$ denote the $\dot C^a$ H\"{o}lder semi-norm in variable $\alpha$. By Lemma \ref{Hol} we obtain that for $\gamma=\frac{3}{4}\pm\varepsilon$, 
			\begin{align*}
				&\frac{|	f_1^{m_1}(x,\alpha)\delta_{\alpha} f_3^{m_3}(x)-f_1^{m_1}(x,\alpha-z)\delta_{\alpha-z} f_3^{m_3}(x)|}{|z|^\gamma}
				\lesssim \|f_1^{m_1}\|_{\gamma,\frac{1}{4}}\|f_3^{m_3}\|_{\dot C^\frac{1}{4}}+\|f_1^{m_1}\|_{L^\infty}\|f_3^{m_3}\|_{\dot C^\gamma}.
			\end{align*}
			Hence we obtain that 
			\begin{align}
				|\tilde P_{2,1}|&\lesssim |\beta|^\kappa  \prod_{+,-} \left\{\|\partial_x f_2^{m_2}\|_{\dot C^{\kappa-\frac{3}{4}\pm\varepsilon}}\left(\|f_1^{m_1}\|_{\frac{3}{4}\pm\varepsilon,\frac{1}{4}}\|f_3^{m_3}\|_{\dot C^\frac{1}{4}}+\|f_1^{m_1}\|_{L^\infty}\|f_3^{m_3}\|_{\dot C^{\frac{3}{4}\pm\varepsilon}} \right) \right\}\nonumber\\
				&\lesssim t^{-\kappa-m} |\beta|^\kappa [f_1]_T\|\partial_x f_2\|_{T,*}\|f_3\|_{T,*}.\label{lP21}
			\end{align}
			For $\tilde P_{2,2}$, by \eqref{aaaaa},
			\begin{align*}
				|\tilde P_{2,2}|&\lesssim \|f_1^{m_1}\|_{L^\infty}\| f_3^{m_3}\|_{\dot C^\frac{1}{4}}\int  \|\delta_\alpha\delta_\beta f_2^{m_2}\|_{L^\infty}|\alpha|^\frac{1}{4}\min\{1,|\alpha|^{-2}\}d\alpha \\
				&\lesssim |\beta|^\kappa\|f_1^{m_1}\|_{L^\infty}\| f_3^{m_3}\|_{\dot C^\frac{1}{4}}\|f_2^{m_2}\|_{\dot C^{\kappa+\frac{1}{4}}}\int \min\{|\alpha|^\frac{1}{2},|\alpha|^{-\frac{3}{2}}\}d\alpha\\
				&\lesssim |\beta|^\kappa\|f_1^{m_1}\|_{L^\infty}\| f_3^{m_3}\|_{\dot C^\frac{1}{4}}\|f_2^{m_2}\|_{\dot C^{\kappa+\frac{1}{4}}}\\
				&\lesssim t^{-\kappa-m} |\beta|^\kappa [f_1]_T\|\partial_x f_2\|_{T,*}\|f_3\|_{T,*}.
			\end{align*}
			Combining this with \eqref{lP21} to obtain that for any $t\in[0,T]$,
			\begin{align}\label{lP2}
				t^{m+\kappa } |\beta|^{-\kappa} |\tilde P_{2}|\lesssim [f_1]_T\|\partial_x f_2\|_{T,*}\|f_3\|_{T,*}.
			\end{align}
			Finally, for $\tilde P_3$, we consider two cases.
			If $m_1+m_2\neq 0$, we directly have
			\begin{align*}
				|\tilde P_3|&\lesssim \|f_1^{m_1}\|_{L^\infty}\|\partial_x f_2\|_{\dot C^{m_2+\frac{3}{4}}}\|f_3\|_{\dot C^{m_3+\kappa+\frac{1}{4}}}\lesssim |\beta|^\kappa t^{-(m+\kappa)}[f_1]_T\|\partial_xf_2\|_{T,*}\|f_3\|_{T,*}.
			\end{align*}
			When $m_1=m_2=0$, by the condition $\int f_1(x,\alpha)E^\alpha f_2(x)\frac{d\alpha}{\alpha}=0$, the estimate is the same with $P_3$. Hence
			\begin{align}\label{lP3}
				t^{m+\kappa}|\beta |^{-\kappa}|\tilde P_{3}|\lesssim [f_1]_T\|\partial_x f_2\|_{T,*}\|f_3\|_{T,*}.
			\end{align}
			Hence we obtain from \eqref{P1}, \eqref{lP2} and \eqref{lP3} that 
			\begin{align*}
				\|\delta_\beta \partial_x f^m\|_{L^\infty}\lesssim 	|\tilde P_1|+|\tilde P_2|+|\tilde P_3|\lesssim |\beta|^\kappa t^{-(m+\kappa)}[f_1]_T\|\partial_xf_2\|_{T,*}\|f_3\|_{T,*}.
			\end{align*}
			This holds for any $\beta \neq 0$. Then we  complete the proof of \eqref{res2}. \vspace{0.5cm}
		\end{proof}\\
		\begin{proof}[Proof of Lemma \ref{lemnonpes}]
			In this proof, we consider $t\in[0,T]$. We first prove \eqref{resn}.	We rewrite  \eqref{defnonpes} to the following form
			\begin{align*}
				\mathcal{N}(X)(x)=&-\sum \int _{\mathbb{R}}H(\tilde \Delta_\alpha X(x))\tilde{\E}^\alpha X_i(x) (\mathbf{T}(|\partial_x X|)\partial_xX_j)(x-\alpha)\frac{d\alpha}{\alpha}.
			\end{align*}
			where the sum is for some $i,j=1,2$ and $H(x)=C_{i_1,i_2,i_3}\frac{x_{i_1}x_{i_2}x_{i_3}}{|x|^4}$, $i_1,i_2,i_3=1,2$.  With a slight abuse of notation, we omit the summation symbol and subscript $i,j$ in our proof. And we do not distinguish $\partial_x f$ and $f'$.
			Note that $H(\tilde \Delta_\alpha X(x))\tilde{\E}^\alpha X_i(x)=\partial_\alpha G(\delta_\alpha X(x))-\frac{1}{2}\cot \frac{\alpha}{2}$, hence
			\begin{align*}
				\int_{\mathbb{S}}H(\tilde {\Delta}_\alpha X(t,x))\tilde \E^\alpha X_i \frac{d\alpha}{\tilde \alpha}=0.
			\end{align*}
			By \eqref{eqpesk}, we can applying \eqref{res2} with 
			\begin{align*}
				f_1(t,x,\alpha)=H(\tilde {\Delta}_\alpha X(t,x)),\ \ \ \ f_2(t,x)=X_i(t,x),\ \ \ \ f_3(t,x)=\big(\mathbf{T}(|\partial_x X|)\partial_xX_j\big)(t,x),
			\end{align*}
			one has 
			\begin{align*}
				\sup_{t\in[0,T]}t^{m+\kappa }\|\mathcal{N}(X)(t)\|_{\dot C^{m+\kappa }}\lesssim [f_1]_T\|\partial_x f_2\|_{T,*}\|f_3\|_{T,*}.
			\end{align*}
			We denote $\tilde{H}[X](t,x,\alpha)=H(\tilde \Delta_\alpha X)(t,x)$. We first consider 
			\begin{align*}
				[f_1]_T=[\tilde{H}[X](t)]_T,
			\end{align*}
			with $[\cdot]_T$ defined in \eqref{defTn}. For any $k\in\mathbb{N}$, $\gamma\in[\frac{1}{2},\kappa]$,
			\begin{align}\label{Tn1}
				\|\partial_x^k\tilde{H}[X](t))\|_{\gamma,\frac{1}{4}}=\sup_{\alpha,z\in\mathbb{S}}|\alpha|^\frac{1}{4}\left(\|\partial_x^kH(\tilde {\Delta}_\alpha X(t,\cdot))\|_{\dot C^\gamma}+\frac{\|\partial_x^kH(\tilde {\Delta}_\alpha X(t,.))-\partial_x^kH(\tilde {\Delta}_{\alpha-z} X(t,.))\|_{L^\infty}}{|z|^\gamma}\right).
			\end{align}
			By definition, it is easy to check that for any $0\leq k\leq m$, $0\leq t\leq T$,
			\begin{equation}\label{Tn2}
				\sup_{\alpha\in\mathbb{S}}|\alpha|^\frac{1}{4}\|\partial_x^kH(\tilde {\Delta}_\alpha X(t,\cdot))\|_{\dot C^\gamma}\lesssim  t^{-(k+\gamma-\frac{1}{4})}(1+\mathbf{\Theta}(X))^{m+2}\|\partial_x X\|_{T,*}(1+\|\partial_x X\|_T)^{m+5}.
			\end{equation}
			Moreover, by Lemma \ref{dede}, 
			\begin{equation}\label{Tn3}
				\begin{aligned}
					&\sup_{\alpha,z\in\mathbb{S}}|\alpha|^\frac{1}{4}\frac{\|\partial_x^kH(\tilde {\Delta}_\alpha X(t,x))-\partial_x^kH(\tilde {\Delta}_{\alpha-z} X(t,x))\|_{L^\infty}}{|z|^\gamma}\\
					&\quad\lesssim t^{-(k+\gamma-\frac{1}{4})}(1+\mathbf{\Theta}(X))^{1+m}\|\partial_x X\|_{T,*}(1+\|\partial_x X\|_T)^{m+5}.
				\end{aligned}
			\end{equation}
			Combine \eqref{Tn1}, \eqref{Tn2} with \eqref{Tn3}, we obtain 
			\begin{align}\label{Hs}
				[\tilde{H}[X](t)]_{T,*}\lesssim (1+\mathbf{\Theta}(X))^{1+m}\|\partial_x X\|_{T,*}(1+\|\partial_x X\|_T)^{m+5}.
			\end{align}
			Note that $\mathbf{T}(|\partial_x X|)\partial_xX_j=\mathcal{T}(|\partial_x X|)\frac{\partial_xX_j}{|\partial_x X|}$, then for any $k\in\mathbb{N}$, $a\in(0,1)$,
			\begin{align*}
				\|\mathbf{T}(|\partial_x X|)\partial_xX_j\|_{\dot C^{k+a}}\lesssim C_{\mathcal{T},k} t^{-(k+a)} \mathbf{\Theta}_X(T) \|\partial_x X\|_{T,*}(1+\|\partial_x X\|_{T})^{1+k}.
			\end{align*}
			This yields 
			\begin{align}\label{Ts}
				\|\mathbf{T}(|\partial_x X|)\partial_xX_j\|_{T,*}\lesssim C_{\mathcal{T},m}\mathbf{\Theta}_X(T) \|\partial_x X\|_{T,*}(1+\|\partial_x X\|_{T})^{1+m}.
			\end{align}
			Then we obtain 
			\begin{align*}
				\sup_{t\in[0,T]}(t^{\kappa}\|\mathcal{N}(X)(t)\|_{\dot C^{\kappa}}+	t^{m+\kappa }\|\mathcal{N}(X)(t)\|_{\dot C^{m+\kappa }})\lesssim C_{\mathcal{T},m}(1+\mathbf{\Theta}_X(T))^{m+5}\|\partial_x X\|_{T,*}^2(1+\|\partial_x X\|_T)^{5(m+1)}.
			\end{align*}
			This completes the proof of \eqref{resn}. \\
			Then we deal with \eqref{rendiff}.
			Note that 
			\begin{align*}
				(\mathcal{N}(Z)-\mathcal{N}(Y))(x)=& \int_{\mathbb{R}} (H(\tilde \Delta_\alpha Y)-H(\tilde \Delta_\alpha Z))(x)\tilde{\E}^\alpha Y(x) (\mathbf{T}(|\partial_x Y|)\partial_xY)(x-\alpha)\frac{d\alpha}{\alpha}\\
				&+ \int_{\mathbb{R}} H(\tilde \Delta_\alpha Z)(x)\tilde{\E}^\alpha (Y-Z)(x)  (\mathbf{T}(|\partial_x Y|)\partial_xY)(x-\alpha)\frac{d\alpha}{\alpha}\\
				&+\int_{\mathbb{R}} H(\tilde \Delta_\alpha Z)(x)\tilde{\E}^\alpha Z(x) (\mathbf{T}(|\partial_x Y|)\partial_xY_j-\mathbf{T}(|\partial_x Z|)\partial_xZ_j)(x-\alpha)\frac{d\alpha}{\alpha}\\
				:=&B_1+B_2+B_3.
			\end{align*}
			We apply \eqref{res1} with $(f_1,f_2,f_3)$ equals $\Big((H(\tilde \Delta_\alpha Y)-H(\tilde \Delta_\alpha Z))(x),Y(x),\mathbf{T}(|\partial_x Y|)\partial_xY(x)\Big)$, $\Big(H(\tilde \Delta_\alpha Z)(x),(Y-Z)(x),\mathbf{T}(|\partial_x Y|)\partial_xY(x)\Big)$ and $\Big(H(\tilde \Delta_\alpha Z)(x),Z(x),(\mathbf{T}(|\partial_x Y|)\partial_xY-\mathbf{T}(|\partial_x Z|)\partial_xZ)(x)\Big)$ respectively. Then we can obtain 
			\begin{align*}
				&\sup_{t\in[0,T]}t^{m+\kappa}\|\partial_x^m B_1(t)\|_{\dot C^\kappa}\\
				&\quad\lesssim \|\partial_x Y\|_{T,*}\left([(\tilde{H}[Y]-\tilde{H}[Z])(t)]_{T,*}\|\mathbf{T}(|\partial_x Y|)\partial_xY\|_T+[(\tilde{H}[Y]-\tilde{H}[Z])(t)]_T\|\mathbf{T}(|\partial_x Y|)\partial_xY\|_{T,*}\right.\\
				&\quad\quad\quad\quad\quad\quad\quad\quad\left.+T^{\frac{1}{2}}[(\tilde{H}[Y]-\tilde{H}[Z])(t)]_T\|\mathbf{T}(|\partial_x Y|)\partial_xY\|_T\right),\\
				&\sup_{t\in[0,T]}t^{m+\kappa}	\|\partial_x^m B_2(t)\|_{\dot C^\kappa}\\
				&\quad\lesssim \|\partial_x (Y-Z)\|_{T,*}\left([\tilde{H}[Z]_{T,*}\|\mathbf{T}(|\partial_x Y|)\partial_xY\|_T+[\tilde{H}[Z]_T\|\mathbf{T}(|\partial_x Y|)\partial_xY\|_{T,*}\right.\\
				&\quad\quad\quad\quad\quad\quad\quad\quad\left.+T^{\frac{1}{2}}[\tilde{H}[Z]_T\|\mathbf{T}(|\partial_x Y|)\partial_xY\|_T\right),\\
				&\sup_{t\in[0,T]}t^{m+\kappa}	\|\partial_x^m B_3(t)\|_{\dot C^\kappa}\\
				&\quad\lesssim  \|\partial_x Z\|_{T,*}\left([\tilde{H}[Z]] ]_{T,*}\|(\mathbf{T}(|\partial_x Y|)\partial_xY-\mathbf{T}(|\partial_x Z|)\partial_xZ)\|_T\right.\\
				&\quad\quad\quad\quad\quad\quad\quad\quad+[\tilde{H}[Z]] ]_T\|(\mathbf{T}(|\partial_x Y|)\partial_xY-\mathbf{T}(|\partial_x Z|)\partial_xZ)\|_{T,*}\\
				&\quad\quad\quad\quad\quad\quad\quad\quad\left.+T^{\frac{1}{2}}[\tilde{H}[Z]] ]_T\|\partial_x Z\|_{T,*}\|(\mathbf{T}(|\partial_x Y|)\partial_xY-\mathbf{T}(|\partial_x Z|)\partial_xZ)\|_T\right).
			\end{align*}
			We will shortly denote $W=Y-Z$, $\vec{V}=(Y,Z)$ in the following proof. We can prove similarly as \eqref{Hs} that 
			\begin{align*}
				[(\tilde{H}[Y]-\tilde{H}[Z]) ]_{T,*}\lesssim (1+\mathbf{\Theta}(\vec{V}))^{1+m} (\|\partial_xW\|_{T,*}+\|\partial_xW\|_{T}\|\partial_xW\|_{T,*})(1+\|\partial_x \vec{V} \|_T)^{m+5},
			\end{align*}
			and 
			\begin{align*}
				\|\mathbf{T}(|\partial_x Y|)\partial_xY-\mathbf{T}(|\partial_x Z|)\partial_xZ\|_{T,*}\lesssim C_{\mathcal{T},m}(1+\mathbf{\Theta}_{\vec{V}}(T))(\|\partial_x W\|_{T,*}+\|\partial_xW\|_T\|\partial_x\vec{V}\|_{T,*})(1+\|\partial_x \vec{V}\|_T)^{m+5}.
			\end{align*}
			Combining the above estimates to obtain
			\begin{align*}
				&\sup_{t\in[0,T]}t^{m+\kappa}(\|\partial_x^m B_1(t)\|_{\dot C^\kappa}+\|\partial_x^m B_2(t)\|_{\dot C^\kappa}+\|\partial_x^m B_3(t)\|_{\dot C^\kappa})\\
				&\lesssim C_{\mathcal{T},m}(1+\mathbf{\Theta}_{\vec{V}}(T))^{5(m+1)}\|\partial_x W\|_{T}(\|\partial_x \vec{V}\|_{T,*}+T^{\frac{1}{2}}\|\partial_x \vec{V}\|_{T})(1+\|\partial_x \vec{V}\|_T)^{5(m+1)}.
			\end{align*} 
			This completes the proof of \eqref{rendiff}.
		\end{proof}
		~~\vspace{0.3cm}\\	
		\begin{proof}[Proof of Lemma \ref{lempesR}] For simplicity, we only estimate the norm $\sup_{t\in[0,T]}t^\kappa\|(\mathcal{M}(\partial_x Y)-\mathcal{M}(\partial_x Z))(t)\|_{\dot C^{\kappa-a}}$. The estimates for higher order derivatives are parallel. We denote $W=Y-Z$ and $\vec{V}=(Y,Z)$ in this proof. Note that
			\begin{align*}
				\mathcal{M}(\partial_xY)(x)-\mathcal{M}(\partial_xZ)(x)=-\frac{1}{\pi}\int_{\mathbb{R}} \mathfrak{R}(\alpha,s)\frac{d\alpha}{|\alpha|^{1+a}}.
			\end{align*}
			Here we denote 
			\begin{align*}
				&\mathfrak{R}(\alpha,s)= \delta_\alpha (\mathbf{T}(|\partial_x Y|)-\mathbf{T}(|\partial_x Z|))(x)\delta_\alpha \partial_x Y(x)+\delta_\alpha (\mathbf{T}(|\partial_x Z|))(x)\delta_\alpha \partial_x W(x)\\
				&\ \ \quad\quad\quad+\partial_xW(x)D_{\mathbf{T},\alpha}(\partial_x Y)(x)+\partial_x Z(x)(D_{\mathbf{T},\alpha}(\partial_x Y)-D_{\mathbf{T},\alpha}(\partial_x Z))(x)\end{align*}
			with 
			\begin{align*}
				&D_{\mathbf{T},\alpha}(A)=\delta_\alpha (\mathbf{T}(|A|))-\delta_\alpha A\cdot \nabla \big(\mathbf{T}(|A|)\big),\ \ \ \forall \ A:\mathbb{R}\to \mathbb{R}^2.
			\end{align*}
			By Lemma \ref{Hol} we have 
			\begin{equation}\label{eR1}
				\begin{aligned}
					&\|\mathfrak{R}(\alpha,\cdot)\|_{\dot C^{\kappa-a}}\lesssim
					\|\mathbf{T}(|\partial_x Y|)-\mathbf{T}(|\partial_x Z|)\|_{\dot C^{\kappa-a}}\|\delta_\alpha \partial_x Y\|_{L^\infty}+\|\delta_\alpha (\mathbf{T}(|\partial_x Y|)-\mathbf{T}(|\partial_x Z|))\|_{L^\infty}\|\partial_x Y\|_{\dot C^{\kappa-a}}\\
					&\quad\quad+\| \mathbf{T}(|\partial_x Z|)\|_{\dot C^{\kappa-a}}\|\delta_\alpha \partial_x W\|_{L^\infty}+\|\delta_\alpha (\mathbf{T}(|\partial_x Z|))\|_{L^\infty}\|\partial_x W\|_{\dot C^{\kappa-a}}\\
					&\quad\quad+\|\partial_xW\|_{\dot C^{\kappa-a}}\|D_{\mathbf{T},\alpha}(\partial_x Y)\|_{L^\infty}+\|\partial_xW\|_{L^\infty}\|D_{\mathbf{T},\alpha}(\partial_x Y)\|_{\dot C^{\kappa-a}}\\
					&\quad\quad+\|\partial_x Z\|_{\dot C^{\kappa-a}}\|(D_{\mathbf{T},\alpha}(\partial_x Y)-D_{\mathbf{T},\alpha}(\partial_x Z))\|_{L^\infty}+\|\partial_x Z\|_{L^\infty}\|(D_{\mathbf{T},\alpha}(\partial_x Y)-D_{\mathbf{T},\alpha}(\partial_x Z))\|_{\dot C^{\kappa-a}},
				\end{aligned}
			\end{equation}
			where all $\dot C^{\kappa-a}$ denote H\"{o}lder semi-norm of $x$ variable. We remark that, all the estimates related to $\mathbf{T}$ will result in a constant $C_{\mathcal{T},m}$ and $(1+\mathbf{\Theta}_Y(T)+\mathbf{\Theta}_Z(T))^m$. We omit these constants in the following to avoid redundancy. 
			Lemma \ref{lemcom} implies that 
			\begin{align*}
				&\|\mathbf{T}(|\partial_x Y|)-\mathbf{T}(|\partial_x Z|)\|_{\dot C^{\kappa-a}}\lesssim \|W'\|_{\dot C^{\kappa-a}}+\| W'\|_{L^\infty}(\|\vec{V}'\|_{\dot C^{\kappa-a}}),\\
				&\| \mathbf{T}(|\partial_x Z|)\|_{\dot C^{\kappa-a}}\lesssim \| Z'\|_{\dot C^{\kappa-a}},\ \ \ \ \ \ \ \ \|D_{\mathbf{T},\alpha}(\partial_x Y)\|_{L^\infty}\lesssim \|\delta_\alpha Y'\|_{L^\infty}^2,\\
				&\|D_{\mathbf{T},\alpha}(\partial_x Y)\|_{\dot C^{\kappa-a}}\lesssim \|\delta_\alpha Y'\|_{\dot C^{\kappa-a}}\|\delta_\alpha Y'\|_{L^\infty}+\| Y'\|_{\dot C^{\kappa-a}}\|\delta_\alpha Y'\|_{L^\infty}^2,\\
				&\|(D_{\mathbf{T},\alpha}(\partial_x Y)-D_{\mathbf{T},\alpha}(\partial_x Z))\|_{L^\infty}\lesssim \|\delta_\alpha W'\|_{L^\infty}\|\delta_\alpha \vec{V}'\|_{L^\infty}+\|\delta_\alpha \vec{V}'\|_{L^\infty}^2\|W'\|_{L^\infty},\\
				&\|(D_{\mathbf{T},\alpha}(\partial_x Y)-D_{\mathbf{T},\alpha}(\partial_x Z))\|_{\dot C^{\kappa-a}}\\
				&\quad\lesssim \|\delta_\alpha W'\|_{\dot C^{\kappa-a}}\|\delta_\alpha \vec{V}'\|_{L^\infty}+\|\delta_\alpha W'\|_{L^\infty}\|\delta_\alpha \vec{V}'\|_{\dot C^{\kappa-a}}+\|\delta_\alpha W'\|_{L^\infty}\|\delta_\alpha \vec{V}'\|_{L^\infty}\|\vec{V}'\|_{\dot C^{\kappa-a}}\\
				&\quad\quad+\|\delta_\alpha \vec{V}'\|_{\dot C^{\kappa-a}}\|\delta_\alpha \vec{V}'\|_{L^\infty}\| W'\|_{L^\infty}+\|\delta_\alpha \vec{V}'\|_{L^\infty}^2(\|W'\|_{\dot C^{\kappa-a}}+\| W'\|_{L^\infty}\|\vec{V}'\|_{\dot C^{\kappa-a}}),
			\end{align*}
			where all norms are related to $x$ variable. Hence 
			\begin{equation}
				\label{eR2}
				\begin{aligned}
					&\|D_{\mathbf{T},\alpha}(\partial_x Y)\|_{L^\infty}\lesssim \min\{|\alpha|^\frac{2}{5}\|Y'\|_{\dot C^\frac{1}{5}}^2,|\alpha|^{2a}\|Y'\|_{\dot C^a}^2\} \lesssim \min\{|\alpha|^\frac{2}{5}t^{-\frac{2}{5}},|\alpha|^{2a}t^{-2a}\}\|Y'\|_{T,*}\|Y'\|_T.\\
					&\|D_{\mathbf{T},\alpha}(\partial_x Y)\|_{\dot C^{\kappa-a}}
					\lesssim  \min\{\|Y'\|_{\dot C^{\kappa-a}}\| Y'\|_{L^\infty}(1+\| Y'\|_{L^\infty}),|\alpha|^{2a}(\|Y'\|_{\dot C^{\kappa}}\|Y'\|_{\dot C^{a}}+\| Y'\|_{\dot C^{\kappa-a}}\|Y'\|_{\dot C^{a}}^2)\}\\
					&\quad\quad\quad\quad \quad\quad\quad\quad \lesssim 	t^{-(\kappa-a)}\min\{1,|\alpha|^{2a}t^{-2a}\}\|Y'\|_{T,*}(1+\|Y'\|_T)^2,\\
					&\|(D_{\mathbf{T},\alpha}(\partial_x Y)-D_{\mathbf{T},\alpha}(\partial_x Z))\|_{L^\infty}\\
					&\quad\quad \lesssim\min\left\{|\alpha|^\frac{1}{5}\|W'\|_{L^\infty}\|\vec{V}'\|_{\dot C ^\frac{1}{5}}(1+\|\vec{V}'\|_{L^\infty}), |\alpha|^{2a}(\|W'\|_{\dot C^a}\|\vec{V}'\|_{\dot C^a}+\|\vec{V}'\|_{\dot C^a}^2\|W'\|_{L^\infty})\right\} \\
					&\quad\quad\lesssim\min\{|\alpha|^\frac{1}{5}t^{-\frac{1}{5}},|\alpha|^{2a}t^{-2a}\}\|W'\|_T\|\vec{V}'\|_{T,*}(1+\|\vec{V}'\|_T).
				\end{aligned}
			\end{equation}
			Similarly, we have both
			{\small\begin{align*}
					&\|(D_{\mathbf{T},\alpha}(\partial_x Y)-D_{\mathbf{T},\alpha}(\partial_x Z))\|_{\dot C^{\kappa-a}}\\
					&\quad\lesssim |\alpha|^\frac{1}{5}\big(\|W'\|_{\dot C^{\kappa-a}}\|\vec{V}'\|_{\dot C^\frac{1}{5}}+\|W'\|_{\dot C^\frac{1}{5}}\|\vec{V}'\|_{\dot C^{\kappa-a}}\big)(1+\|\vec{V}'\|_{L^\infty})+|\alpha|^\frac{1}{5}\|W'\|_{L^\infty}\|\vec{V}'\|_{\dot C^{\kappa-a}}\|\vec{V}'\|_{\dot C^\frac{1}{5}}\\
					&\quad\lesssim|\alpha|^\frac{1}{5} t^{-(\kappa-a+\frac{1}{5})}\|W'\|_T\|\vec{V}'\|_{T,*}(1+\|\vec{V}'\|_T),
			\end{align*}}
			and 
			\begin{align*}
				&\|(D_{\mathbf{T},\alpha}(\partial_x Y)-D_{\mathbf{T},\alpha}(\partial_x Z))\|_{\dot C^{\kappa-a}}\\
				&\quad\lesssim |\alpha|^{2a}\left\{\|W'\|_{\dot C^{\kappa}}\|\vec{V}'\|_{\dot C^{a}}+\|W'\|_{\dot C^{a}}\|\vec{V}'\|_{\dot C^{\kappa}}(1+\|\vec{V}'\|_{L^\infty})\right.\\
				&\quad\quad\quad\left.+\|W'\|_{L^\infty}\|\vec{V}'\|_{\dot C^{\kappa}}\|\vec{V}'\|_{\dot C^{a}}+\|\vec{V}'\|_{\dot C^a}^2(\|W'\|_{\dot C^{\kappa-a}}+\| W'\|_{L^\infty}\|\vec{V}'\|_{\dot C^{\kappa-a}})\right\}\\
				&\quad\lesssim |\alpha|^{2a}t^{-(\kappa+a)}\|W'\|_T\|\vec{V}'\|_{T,*}(1+\|\vec{V}'\|_T)^2.
			\end{align*}
			Hence 
			\begin{align}\label{eR3}
				&\|(D_{\mathbf{T},\alpha}(\partial_x Y)-D_{\mathbf{T},\alpha}(\partial_x Z))\|_{\dot C^{\kappa-a}}\lesssim \frac{\min\{1,|\alpha|^{2a}t^{-2a}\}}{t^{\kappa-a}}\|W'\|_T\|\vec{V}'\|_{T,*}(1+\|\vec{V}'\|_T)^2.
			\end{align}
			Combining \eqref{eR1}-\eqref{eR3}, we obtain 
			\begin{align*}
				\|\mathfrak{R}(\alpha,\cdot)\|_{\dot C^{\kappa-a}}\lesssim \frac{\min\{1,|\alpha|t^{-1}\}^{2a}}{t^{\kappa-a}}\|W'\|_T\|\vec{V}'\|_{T,*}(1+\|\vec{V}'\|_T)^2. 
			\end{align*}
			Hence 
			\begin{align*}
				&\|\mathcal{M}(\partial_xY)-\mathcal{M}(\partial_xZ)\|_{\dot C^{\kappa-a}}\lesssim\int 	\|\mathfrak{R}(\alpha,\cdot)\|_{\dot C^{\kappa-a}}\frac{d\alpha}{|\alpha|^{1+a}}\\
				&\lesssim t^{-(\kappa-a)}\int \min\{1,|\alpha|^{2a}t^{-2a}\}\frac{d\alpha}{|\alpha|^{1+a}}\|W'\|_T\|\vec{V}'\|_{T,*}(1+\|\vec{V}'\|_T)^2\\
				&\lesssim t^{-\kappa }\|W'\|_T\|\vec{V}'\|_{T,*}(1+\|\vec{V}'\|_T)^2.
			\end{align*}
			This completes the proof.
		\end{proof}\\
		\subsection{Estimates of nonlinear terms in 3D Peskin equation}
		The following lemma will be useful in the estimate of 3D Peskin equation.
		\begin{lemma}
			\label{lemG11}
			Suppose there exists $C_{\mathbf{G}}>0$ such that  \begin{align}\label{conG}|\mathbf{G}(\theta,\eta)|+|\theta-\eta||\nabla _{\theta,\eta} \mathbf{G}(\theta,\eta)|+|\theta-\eta|^2|\nabla _\theta\nabla _\eta \mathbf{G}(\theta,\eta)|\leq \frac{C_{\mathbf{G}}}{|\theta-\eta|},\ \ \ \forall \theta,\eta\in\mathbb{R}^2,\ \theta\neq \eta.\end{align}   Define 
			\begin{align*}
				\mathcal{M}f(\theta)=\int_{\mathbb{R}^2}\mathbf{G}(\theta,\eta)\nabla \cdot f(\eta)d\eta.
			\end{align*}
			Then there holds 
			\begin{align*}
				\|\mathcal{M}f\|_{\dot C^\kappa}\lesssim C_{\mathbf{G}}\|f\|_{\dot C^\kappa}.
			\end{align*}
		\end{lemma}
		\begin{proof}
			For $\theta\neq \theta'$, denote $\delta=|\theta-\theta'|$. We further decompose $\mathcal{M} f$ into high and low frequency parts:
			\begin{align*}
				\mathcal{M}f(\theta)&=\int_{\mathbb{R}^2}\chi_{\delta}(|\theta-\eta|)\mathbf{G}(\theta,\eta)\nabla \cdot f(\eta)d\eta+\int_{\mathbb{R}^2}\chi_{\delta}^c(|\theta-\eta|)\mathbf{G}(\theta,\eta)\nabla \cdot f(\eta)d\eta\\
				:&=\mathcal{M}_1f(\theta)+\mathcal{M}_2f(\theta),
			\end{align*}
			where we denote $\chi^c_\delta=1-\chi_\delta$ for simplicity.
			One has 
			\begin{align*}
				| \mathcal{M}f(\theta)- \mathcal{M}f(\theta')|\leq |\mathcal{M}_1f(\theta)|+| \mathcal{M}_1f(\theta')|+| \mathcal{M}_2f(\theta)- \mathcal{M}_2f(\theta')|.
			\end{align*}
			Integrate by parts to obtain 
			\begin{equation}\label{hf2}
				\begin{aligned}
					|\mathcal{M}_1f(\theta)|&\lesssim \int_{\mathbb{R}^2}|\nabla _\eta \left(\chi_{\delta}(|\theta-\eta|)\mathbf{G}(\theta,\eta)\right)|| f(\eta)-f(\theta)|d\eta\\
					&\lesssim \left(\int_{|\theta-\eta|\leq 2\delta}|\nabla _\eta\mathbf{G}(\theta,\eta)||\theta-\eta|^\kappa d\eta +\delta^{-1}\int_{\delta\leq|\theta-\eta|\leq 2\delta}|\mathbf{G}(\theta,\eta)||\theta-\eta|^\kappa d\eta\right)\|f\|_{\dot C^\kappa}\\
					&\lesssim C_{\mathbf{G}}\delta^\kappa \|f\|_{\dot C^\kappa}.
				\end{aligned}
			\end{equation}
			Similarly, there holds 
			\begin{align}\label{hf1}
				|\mathcal{M}_1f(\theta')|&\lesssim C_{\mathbf{G}}\delta^\kappa \|f\|_{\dot C^\kappa}.
			\end{align}
			It remains to deal with $| \mathcal{M}_2f(\theta)- \mathcal{M}_2f(\theta')|$. Observe that 
			\begin{align*}
				&\chi_{\delta}^c(|\theta-\eta|)\mathbf{G}(\theta,\eta)-\chi_{\delta}^c(|\theta'-\eta|)\mathbf{G}(\theta',\eta)\\
				&\quad\quad=\int_0^1 (\theta-\theta')\cdot \nabla _\theta \left(\mathbf{G}(d_{\theta,\theta'}(\tau),\eta)\chi^c_\delta(|d_{\theta,\theta'}(\tau)-\eta|)\right)d\tau ,
			\end{align*}
			where we denote $d_{\theta,\theta'}(\tau)=\tau\theta+(1-\tau)\theta'$.
			Hence we  can write 
			\begin{align*}
				\mathcal{M}_2f(\theta)- \mathcal{M}_2f(\theta')&=\int_{\mathbb{R}^2}\big(\chi_{\delta}^c(|\theta-\eta|)\mathbf{G}(\theta,\eta)-\chi_{\delta}^c(|\theta'-\eta|)\mathbf{G}(\theta',\eta)\big)\nabla \cdot f(\eta)d\eta\\
				&=\int_{\mathbb{R}^2}\int_0^1 (\theta-\theta')\cdot \nabla _\theta \left(\mathbf{G}(d_{\theta,\theta'}(\tau),\eta)\chi^c_\delta(|d_{\theta,\theta'}(\tau)-\eta|)\right)d\tau \nabla\cdot f(\eta)d\eta\\
				&=\int_0^1\int_{\mathbb{R}^2} (\theta-\theta')\cdot \nabla _\theta \nabla _\eta\left(\mathbf{G}(d_{\theta,\theta'}(\tau),\eta)\chi^c_\delta(|d_{\theta,\theta'}(\tau)-\eta|)\right)\cdot (f(d_{\theta,\theta'}(\tau))-f(\eta)) d\eta d\tau.
			\end{align*}
			Note that 
			\begin{align*}
				&|\nabla _\theta \nabla _\eta\left(\mathbf{G}(d_{\theta,\theta'}(\tau),\eta)\chi^c_\delta(|d_{\theta,\theta'}(\tau)-\eta|)\right)|\\
				&\lesssim C_{\mathbf{G}}\left(\mathbf{1}_{|\eta-d_{\theta,\theta'}(\tau)|\geq \delta}|\eta-d_{\theta,\theta'}(\tau)|^{-3}+\mathbf{1}_{\delta\leq|\eta-d_{\theta,\theta'}(\tau)|\leq 2\delta}(\delta^{-1}|\eta-d_{\theta,\theta'}(\tau)|^{-2}+\delta^{-2}|\eta-d_{\theta,\theta'}(\tau)|^{-1})\right).
			\end{align*}
			Hence 
			\begin{equation}\label{lf}
				\begin{aligned}
					|\mathcal{M}_2f(\theta)- \mathcal{M}_2f(\theta')|&\lesssim C_{\mathbf{G}} \delta\|f\|_{\dot C^\kappa} \int_{\mathbb{R}^2}\mathbf{1}_{|\eta|\geq \delta} |\eta|^{-3+\kappa} +\mathbf{1}_{\delta\leq |\eta|\leq 2\delta}\delta^{-1}|\eta|^{-2+\kappa}d\eta \\
					&\lesssim C_{\mathbf{G}}\delta^\kappa \|f\|_{\dot C^\kappa}.
				\end{aligned}
			\end{equation}
			We conclude from \eqref{hf2}, \eqref{hf1} and \eqref{lf} that
			\begin{align*}
				|\mathcal{M}f(\theta)-\mathcal{M}f(\theta')|\lesssim C_{\mathbf{G}}\delta^\kappa\|f\|_{\dot C^\kappa}.
			\end{align*}
			This completes the proof.
		\end{proof}\\
		\begin{remark}\label{rmk3dpesk}
			The result of Lemma \ref{lemG11} on $\mathbb{R}^2$ can be extended parallel to $\mathbb{S}^2$, which means, if we denote
			\begin{equation*}
				\tilde{\mathcal{M}} F(\widehat{\boldsymbol{x}})=\int_{\mathbb{S}^2}\mathbf{G}(\widehat{\boldsymbol{x}},\widehat{\boldsymbol{y}})\widetilde\nabla\cdot F(\widehat{\boldsymbol{y}})d \mu_{\mathbb{S}^2}(\widehat{\boldsymbol{y}}),
			\end{equation*}
			and 
			\begin{equation*}
				|\mathbf{G}(\widehat{\boldsymbol{x}},\widehat{\boldsymbol{y}})|+|\widehat{\boldsymbol{x}}-\widehat{\boldsymbol{y}}||\widetilde\nabla_{ \widehat{\boldsymbol{x}},\widehat{\boldsymbol{y}}}\mathbf{G}(\widehat{\boldsymbol{x}},\widehat{\boldsymbol{y}})|+|\widehat{\boldsymbol{x}}-\widehat{\boldsymbol{y}}|^2|\widetilde\nabla_{ \widehat{\boldsymbol{x}}}\widetilde\nabla_{ \widehat{\boldsymbol{y}}}\mathbf{G}(\widehat{\boldsymbol{x}},\widehat{\boldsymbol{y}})|\leq \frac{C_{\mathbf G}}{|\widehat{\boldsymbol{x}}-\widehat{\boldsymbol{y}}|},
			\end{equation*}
			then
			\begin{equation*}
				\|\tilde{M}F\|_{\dot C^{\kappa}(\mathbb{S}^2)}\lesssim C_{\mathbf{G}}\|F\|_{\dot C^{\kappa}(\mathbb{S}^2)}.
			\end{equation*}
			We can just copy the proof of $\mathbb{R}^2$ or make localization on $\mathbb{S}^2$, we omit the details of proof here.
		\end{remark}
		\begin{lemma}\label{3dpeslemn1}
			Let $N_1(F,\Phi)$ be defined in \eqref{3dpesnlt}. For any $T>0$, and any $X,Y,Z\in\mathcal{Z}_{T,\Phi}^\sigma$, 
			\begin{equation*}
				\begin{aligned}
					&\|N_1(X,\Phi)\|_{L_T^1L^\infty}+\sup_{t\in[0,T]}(t^{\kappa}\|N_1(X,\Phi)\|_{\dot C^\kappa}+t^{m+\kappa}\|N_1(X,\Phi)\|_{\dot C^{m+\kappa}})\\
					&\quad\lesssim \|\widetilde\nabla(X-\Phi)\|_{T}(\|\widetilde\nabla(X-\Phi)\|_{T}+T^{\kappa}\|\Phi\|_{C^{m+3}})^{m+2},
				\end{aligned}
			\end{equation*}
			and
			\begin{equation*}
				\begin{aligned}
					&\|N_1(Y,\Phi)-N_1(Z,\Phi)\|_{L_T^1L^\infty}+\sup_{t\in[0,T]}(t^{\kappa}\|N_1(Y,\Phi)-N_1(Z,\Phi)\|_{\dot C^\kappa}+t^{m+\kappa}\|N_1(Y,\Phi)-N_1(Z,\Phi)\|_{\dot C^{m+\kappa}})\\
					&\quad\quad\lesssim \|\widetilde\nabla(Y-Z)\|_{T}(\|\widetilde\nabla(Y-\Phi,Z-\Phi)\|_{T}+T^{\kappa}\|\Phi\|_{C^{m+3}})(1+\|\widetilde\nabla(Y-\Phi,Z-\Phi)\|_{T}+T\|\Phi\|_{C^{m+3}})^{m+2}.
				\end{aligned}
			\end{equation*}
		\end{lemma}
		\begin{proof}
			For $L^\infty$ norm, we can see that
			\begin{equation}\label{3dpeslemn10}
				\begin{aligned}
					&\|N_1(X,\Phi)\|_{L^\infty}\\
					&\lesssim\left\|\int \left(G(X(\widehat{\boldsymbol{x}})-X(\widehat{\boldsymbol{y}}))-G(\Phi(\widehat{\boldsymbol{x}})-\Phi(\widehat{\boldsymbol{y}}))\right)\widetilde\nabla\cdot ((\mathbf T(|\widetilde\nabla X|)\widetilde\nabla X)(\widehat{\boldsymbol{y}})-(\mathbf T(|\widetilde\nabla X|)\widetilde\nabla X)(\widehat{\boldsymbol{y}}))d\mu_{\mathbb{S}^2}(\widehat{\boldsymbol{x}})\right\|_{L^\infty}\\
					&\lesssim\|X-\Phi\|_{Z_T}\|\mathbf T(|\widetilde\nabla X|)\widetilde\nabla X\|_{\dot C^{\eps}}\lesssim t^{-\eps}\|X-\Phi\|_{Z_T}(\|X-\Phi\|_{Z_T}+T^{\eps}\|\Phi\|_{C^3}).
				\end{aligned}
			\end{equation}
			and similarly, by Lemma \ref{lemcom},
			\begin{equation}\label{3dpeslemn11}
				\begin{aligned}
					&\|N_1(Y,\Phi)-N_1(Z,\Phi)\|_{L^\infty}\\
					&\lesssim\left\|\int \left(G(Y(\widehat{\boldsymbol{x}})-Y(\widehat{\boldsymbol{y}}))-G(Z(\widehat{\boldsymbol{x}})-Z(\widehat{\boldsymbol{y}}))\right)\widetilde\nabla\cdot ((\mathbf T(|\widetilde\nabla Y|)\widetilde\nabla Y)(\widehat{\boldsymbol{y}})-(\mathbf T(|\widetilde\nabla Y|)\widetilde\nabla Y)(\widehat{\boldsymbol{y}}))d\mu_{\mathbb{S}^2}(\widehat{\boldsymbol{x}})\right\|_{L^\infty}\\
					&+\left\|\int \left(G(Z(\widehat{\boldsymbol{x}})-Z(\widehat{\boldsymbol{y}}))-G(\Phi(\widehat{\boldsymbol{x}})-\Phi(\widehat{\boldsymbol{y}}))\right)\widetilde\nabla\cdot \big(((\mathbf T(|\widetilde\nabla Y|)\widetilde\nabla Y)(\widehat{\boldsymbol{y}})-(\mathbf T(|\widetilde\nabla Y|)\widetilde\nabla Y)(\widehat{\boldsymbol{x}})\right.\\
					&\quad\quad\quad\left.-(\mathbf T(|\widetilde\nabla Z|)\widetilde\nabla Z)(\widehat{\boldsymbol{y}})-(\mathbf T(|\widetilde\nabla Z|)\widetilde\nabla Z)(\widehat{\boldsymbol{x}}))d\mu_{\mathbb{S}^2}(\widehat{\boldsymbol{y}})\big)\right\|_{L^\infty}\\
					&\lesssim t^{-\eps}\|Y-Z\|_{Z_T}(\|Y-\Phi\|_{Z_T}+\|Z-\Phi\|_{Z_T}+T^{\eps}\|\Phi\|_{C^3}).
				\end{aligned}
			\end{equation}
			For H\"{o}lder norm, let $\mathbf{G}(\widehat{\boldsymbol{x}},\widehat{\boldsymbol{y}})=G(X(\widehat{\boldsymbol{x}})-X(\widehat{\boldsymbol{y}}))-G(\Phi(\widehat{\boldsymbol{x}})-\Phi(\widehat{\boldsymbol{y}}))$ and $C_{\mathbf{G}}=C\|\widetilde\nabla(X-\Phi)\|_{L^\infty}(1+\|\widetilde\nabla(X-\Phi)\|_{L^\infty})$. By \eqref{3dpesnc} we can see that the condition \eqref{conG} is satisfied. Then apply Remark \ref{rmk3dpesk} we have
			\begin{equation}\label{3dpesn1}
				\begin{aligned}
					\|N_1\|_{\dot C^{\kappa}}&\lesssim \|\widetilde\nabla(X-\Phi)\|_{L^\infty}(1+\|\widetilde\nabla(X-\Phi)\|_{L^\infty})\|\mathbf{T}(|\widetilde\nabla X|)\widetilde\nabla X\|_{\dot C^\kappa(\mathbb{S}^2)}\\
					&\lesssim t^{-\kappa}\|X-\Phi\|_{Z_T}(1+\|X-\Phi\|_{Z_T})(\|X-\Phi\|_{Z_T}+T^{\kappa}\|\Phi\|_{C^3(\mathbb{S}^2)}).
				\end{aligned}
			\end{equation}
			For higher order derivatives of $N_1(X,\Phi)$, we apply \eqref{3dpeskerd}, integral by parts, and define the kernel
			\begin{equation*}
				G_l(\widehat{\boldsymbol{x}},\widehat{\boldsymbol{y}})=(\widetilde\nabla_{ \widehat{\boldsymbol{x}}}+\widetilde\nabla_{ \widehat{\boldsymbol{y}}})^l\left(G(X(\widehat{\boldsymbol{x}})-X(\widehat{\boldsymbol{y}}))-G(\Phi(\widehat{\boldsymbol{x}})-\Phi(\widehat{\boldsymbol{y}}))\right).
			\end{equation*}
			Use the fact that
			\begin{equation*}
				\begin{aligned}
					&|G_l(\widehat{\boldsymbol{x}},\widehat{\boldsymbol{y}})|+|\widehat{\boldsymbol{x}}-\widehat{\boldsymbol{y}}||\widetilde\nabla_{ \widehat{\boldsymbol{x}},\widehat{\boldsymbol{y}}}G_l(\widehat{\boldsymbol{x}},\widehat{\boldsymbol{y}})|+|\widehat{\boldsymbol{x}}-\widehat{\boldsymbol{y}}|^2|\widetilde\nabla_{ \widehat{\boldsymbol{x}}}\widetilde\nabla_{ \widehat{\boldsymbol{y}}}G_l(\widehat{\boldsymbol{x}},\widehat{\boldsymbol{y}})|\\
					&\quad\lesssim t^{-l}\frac{\|\widetilde\nabla(X-\Phi)\|_{T}(1+\|\widetilde\nabla(X-\Phi)\|_{T}+T\|\Phi\|_{C^{m+3}})^{m}}{|\widehat{\boldsymbol{x}}-\widehat{\boldsymbol{y}}|}.
				\end{aligned}
			\end{equation*}
			Then also by Remark \ref{rmk3dpesk} with corresponding $G_l$ and $C_{G_l}$, we obtain
			\begin{equation}\label{3dpesn1h}
				\begin{aligned}
					\|N_1\|_{\dot C^{m+\kappa}}\lesssim t^{-m-\kappa}\|X-\Phi\|_{Z_T}(1+\|X-\Phi\|_{Z_T}+T\|\Phi\|_{C^{m+3}})^{m}(\|X-\Phi\|_{Z_T}+T^{\kappa}\|\Phi\|_{C^3(\mathbb{S}^2)})^m.
				\end{aligned}
			\end{equation}
			For $N_1(Y,\Phi)-N_1(Z,\Phi)$ term, we will shortly denote $W=Y-Z$, $\vec{V}=(Y,Z)$ and $\Vec{V}-\Phi=(Y-\Phi,Z-\Phi)$ in the following proof. We have
			\begin{equation*}
				\begin{aligned}
					&N_1(Y,\Phi)-N_1(Z,\Phi)=\int_{\mathbb{S}^2}\left(G(Y(\widehat{\boldsymbol{x}})-Y(\widehat{\boldsymbol{y}}))-G(Z(\widehat{\boldsymbol{x}})-Z(\widehat{\boldsymbol{y}}))\right)\widetilde\nabla\cdot(\mathbf{T}(|\widetilde\nabla Y|)\widetilde\nabla Y)(\widehat{\boldsymbol{y}})d \mu_{\mathbb{S}^2}(\widehat{\boldsymbol{y}})\\
					&\quad+\int_{\mathbb{S}^2}\left(G(Z(\widehat{\boldsymbol{x}})-Z(\widehat{\boldsymbol{y}}))-G(\Phi(\widehat{\boldsymbol{x}})-\Phi(\widehat{\boldsymbol{y}}))\right)\widetilde\nabla\cdot(\mathbf{T}(|\widetilde\nabla Y|)\widetilde\nabla Y-\mathbf{T}(|\widetilde\nabla Z|)\widetilde\nabla Z)(\widehat{\boldsymbol{y}})d \mu_{\mathbb{S}^2}(\widehat{\boldsymbol{y}}).
				\end{aligned}
			\end{equation*}
			Separately take $\mathbf{G}=G(Y(\widehat{\boldsymbol{x}})-Y(\widehat{\boldsymbol{y}}))-G(Z(\widehat{\boldsymbol{x}})-Z(\widehat{\boldsymbol{y}})$ and $\mathbf{G}=G(Z(\widehat{\boldsymbol{x}})-Z(\widehat{\boldsymbol{y}}))-G(\Phi(\widehat{\boldsymbol{x}})-\Phi(\widehat{\boldsymbol{y}}))$, by Remark \ref{rmk3dpesk}, and by similar methods of \eqref{3dpesn1} and \eqref{3dpesn1h}, one can obtain
			\begin{equation}\label{pesN1}
				\begin{aligned}
					&\|N_1(Y,\Phi)-N_1(Z,\Phi)\|_{\dot C^{m+\kappa}}\\
					&\lesssim  t^{-m+\kappa}\|\widetilde\nabla W\|_{T}(\|\vec{V}-\Phi\|_{Z_T}+T^{\kappa}\|\Phi\|_{C^3(\mathbb{S}^2)})(1+\|\vec{V}-\Phi\|_{Z_T}+T\|\Phi\|_{C^3(\mathbb{S}^2)})^m.
				\end{aligned}
			\end{equation}
			This completes the proof.
		\end{proof}
		\begin{lemma}\label{3dpeslemn2}
			Let $N_2(F,\Phi)$ be defined in \eqref{3dpesnlt}. For any $T>0$, and any $X,Y,Z\in\mathcal{Z}_{T,\Phi}^\sigma$, 
			\begin{equation*}
				\begin{aligned}
					&\|N_2(X,\Phi)\|_{L_T^1L^\infty}+\sup_{t\in[0,T]}(t^{\kappa}\|N_2(X,\Phi)\|_{\dot C^\kappa}+t^{m+\kappa}\|N_2(X,\Phi)\|_{\dot C^{m+\kappa}})\\
					&\lesssim \|\Phi\|_{C^{1}}\|\widetilde\nabla(X-\Phi)\|_{T}(\|\widetilde\nabla(X-\Phi)\|_{T}+T^{\kappa}\|\Phi\|_{C^{m+3}})^{m+2},
				\end{aligned}
			\end{equation*}
			and
			\begin{equation*}
				\begin{aligned}
					&\|N_2(Y,\Phi)-N_2(Z,\Phi)\|_{L_T^1L^\infty}+\sup_{t\in[0,T]}(t^{\kappa}\|N_2(Y,\Phi)-N_2(Z,\Phi)\|_{\dot C^\kappa}+t^{m+\kappa}\|N_2(Y,\Phi)-N_2(Z,\Phi)\|_{\dot C^{m+\kappa}})\\
					&\lesssim \|\Phi\|_{C^{1}}\|\widetilde\nabla(Y-Z)\|_{T}(\|\widetilde\nabla(Y-\Phi,Z-\Phi)\|_{T}+T^{\kappa}\|\Phi\|_{C^{m+3}})(1+\|\widetilde\nabla(Y-\Phi,Z-\Phi)\|_{T}+T\|\Phi\|_{C^{m+3}})^{m+2}.
				\end{aligned}
			\end{equation*} 
		\end{lemma}
		\begin{proof}
			The proof is just similar with Lemma \ref{lempesR}. We denote 
			\begin{equation*}
				R(F,\Phi)=\mathbf{T}(|\widetilde\nabla F|) \widetilde\nabla F-\mathbf{T}(|\widetilde\nabla\Phi|) \widetilde\nabla\Phi-J(\widetilde\nabla\Phi) \widetilde\nabla(F-\Phi).
			\end{equation*}
			For $R$, notice that
			\begin{equation*}
				\begin{aligned}
					&\mathbf{T}(|\widetilde\nabla X|)\widetilde\nabla X-\mathbf{T}(|\widetilde\nabla\Phi|)\widetilde\nabla\Phi\\&\quad\quad\quad\quad\quad
					=\mathbf{T}(|\widetilde\nabla\Phi|)\widetilde\nabla(X-\Phi)+(\mathbf{T}(|\widetilde\nabla X|)-\mathbf{T}(|\widetilde\nabla\Phi|))\widetilde\nabla(X-\Phi)+(\mathbf{T}(|\widetilde\nabla X|)-\mathbf{T}(|\widetilde\nabla\Phi|))\widetilde\nabla\Phi\\
					&\quad\quad\quad\quad\quad
					=\mathbf{T}(|\widetilde\nabla\Phi|)\widetilde\nabla(X-\Phi)+\mathbf{T}'(|\widetilde\nabla\Phi|)(|\widetilde\nabla X|-|\widetilde\nabla\Phi|)\widetilde\nabla\Phi+R_1+R_2\\
					&\quad\quad\quad\quad\quad
					=\mathbf{T}(|\widetilde\nabla\Phi|)\widetilde\nabla(X-\Phi)+J_2(\widetilde\nabla\Phi)\widetilde\nabla(X-\Phi)+R_1+R_2+R_3\\
					&\quad\quad\quad\quad\quad=J(\widetilde\nabla\Phi)\widetilde\nabla(X-\Phi)+R_1+R_2+R_3,\\
				\end{aligned}
			\end{equation*}
			where
			\begin{align*}
				&R_1=(\mathbf{T}(|\widetilde\nabla X|)-\mathbf{T}(|\widetilde\nabla\Phi|))\widetilde\nabla\Phi-\mathbf{T}'(|\widetilde\nabla\Phi|)(|\widetilde\nabla X|-|\widetilde\nabla\Phi|)\widetilde\nabla\Phi,\\
				&R_2=(\mathbf{T}(|\widetilde\nabla X|)-\mathbf{T}(|\widetilde\nabla\Phi|))\widetilde\nabla(X-\Phi),\\
				&R_3=\mathbf{T}'(|\widetilde\nabla\Phi|)(|\widetilde\nabla X|-|\widetilde\nabla\Phi|)\widetilde\nabla\Phi-J_2(\widetilde\nabla\Phi)\widetilde\nabla(X-\Phi),
			\end{align*}
			with $J_2$ defined as \eqref{3dpesnlt}. For $L^\infty$ norm, similar with \eqref{3dpeslemn10} and \eqref{3dpeslemn11},
			\begin{equation*}
				\begin{aligned}
					&\|N_2(X,\Phi)\|_{L^\infty}\lesssim \|\Phi\|_{C^3}(1+\|\Phi\|_{C^3})\|R(X,\Phi)\|_{\dot C^{\eps}}\\
					&\quad\quad\quad\lesssim t^{-\eps}\|\Phi\|_{C^3}(1+\|\Phi\|_{C^3})\|X-\Phi\|_{Z_T}(\|X-\Phi\|_{Z_T}+T^{\eps}\|\Phi\|_{C^3}),\\
					&\|N_2(Y,\Phi)-N_2(Z,\Phi)\|_{L^\infty}\lesssim t^{-\eps}\|\Phi\|_{C^3}(1+\|\Phi\|_{C^3})\|Y-Z\|_{Z_T}(\|Y-\Phi\|_{Z_T}+\|Z-\Phi\|_{Z_T}+T^{\eps}\|\Phi\|_{C^3}).
				\end{aligned}
			\end{equation*}
			For H\"{o}lder norm, let $\mathbf{G}=G(\Phi(\widehat{\boldsymbol{x}})-\Phi(\widehat{\boldsymbol{y}}))$ and $C_{\mathbf{G}}=\|\Phi\|_{C^1}(1+\|\Phi\|_{C^1})$. By \eqref{3dpesnc}, the condition \eqref{conG} is satisfied. Then we apply Remark \ref{rmk3dpesk} to obtain
			\begin{equation*}
				\|N_2(X,\Phi)\|_{\dot C^\kappa}\lesssim \|\Phi\|_{C^1}(1+\|\Phi\|_{C^1})\|R(X,\Phi)\|_{\dot C^{\kappa}}.
			\end{equation*}
			By Newton-Leibniz formula, we can rewrite the remainder terms as
			\begin{equation*}
				\begin{aligned}
					&R_1=\int_0^1\int_0^1\mathbf{T}''(|\nabla_{\mathbb{S}^2}\Phi|+\lambda\mu(|\nabla_{\mathbb{S}^2}X|-|\nabla_{\mathbb{S}^2}\Phi|))\lambda (|\nabla_{\mathbb{S}^2}X|-|\nabla_{\mathbb{S}^2}\Phi|)^2\nabla_{\mathbb{S}^2}\Phi d\mu d\lambda,\\
					&R_2=\int_0^1\mathbf{T}'(|\nabla_{\mathbb{S}^2}\Phi|+\lambda(|\nabla_{\mathbb{S}^2}X|-|\nabla_{\mathbb{S}^2}\Phi|))(|\nabla_{\mathbb{S}^2}X|-|\nabla_{\mathbb{S}^2}\Phi|)\nabla_{\mathbb{S}^2}(X-\Phi)d\lambda,\\
					&R_3=\mathbf{T}'(|\nabla_{\mathbb{S}^2}\Phi|)\left(\frac{\nabla_{\mathbb{S}^2}(X-\Phi)\otimes\nabla_{\mathbb{S}^2}\Phi}{|\nabla_{\mathbb{S}^2}X|+|\nabla_{\mathbb{S}^2}\Phi|}+\left(\frac{2\nabla_{\mathbb{S}^2}\Phi\otimes\nabla_{\mathbb{S}^2}\Phi}{|\nabla_{\mathbb{S}^2}X|+|\nabla_{\mathbb{S}^2}\Phi|}-\frac{\nabla_{\mathbb{S}^2}\Phi\otimes\nabla_{\mathbb{S}^2}\Phi}{|\nabla_{\mathbb{S}^2}\Phi|}\right)\right)\nabla_{\mathbb{S}^2}(X-\Phi).
				\end{aligned}
			\end{equation*}
			By Lemma \ref{Hol} we can see that
			\begin{equation*}
				\|R(X,\Phi)\|_{\dot C^{\kappa}}\lesssim t^{-\kappa}\|\nabla_{\mathbb{S}^2}(X-\Phi)\|_{T}(\|\nabla_{\mathbb{S}^2}(X-\Phi)\|_{T}+T^{\kappa}\|\Phi\|_{C^3})(1+\|\nabla_{\mathbb{S}^2}(X-\Phi)\|_{T}+\|\Phi\|_{C^3}),
			\end{equation*}
			which infers
			\begin{equation*}
				\|N_2(X,\Phi)\|_{\dot C^\kappa}\lesssim t^{-\kappa}(1+\|\Phi\|_{C^3})^2\|\nabla_{\mathbb{S}^2}(X-\Phi)\|_{T}(\|\nabla_{\mathbb{S}^2}(X-\Phi)\|_{T}+T^{\kappa}\|\Phi\|_{C^3})(1+\|\nabla_{\mathbb{S}^2}(X-\Phi)\|_{T}+\|\Phi\|_{C^3}).
			\end{equation*}
			For higher order derivatives of $N_2(X,\Phi)$, we only need to use \eqref{3dpeskerd}, integral by parts and take derivatives on $R$, then by similar way of \eqref{3dpesn1h} we can obtain
			\begin{equation*}
				\begin{aligned}
					&\|N_2(X,\Phi)\|_{\dot C^{m+\kappa}}\\
					&\lesssim t^{-m-\kappa}\|\Phi\|_{C^{1}}\|X-\Phi\|_{Z_T}(\|X-\Phi\|_{Z_T}+T^{\kappa}\|\Phi\|_{C^{m+3}})(1+\|X-\Phi\|_{Z_T}+T\|\Phi\|_{C^{m+3}}).
				\end{aligned}
			\end{equation*}
			For $N_2(Y,\Phi)-N_2(Z,\Phi)$, we will shortly denote $W=Y-Z$, $\vec{V}=(Y,Z)$ and $\vec{V}-\Phi=(Y-\Phi,Z-\Phi)$ in the following proof. We have
			\begin{equation*}
				\begin{aligned}
					\mathbf{T}(|\widetilde\nabla Y|)\widetilde\nabla Y-\mathbf{T}(|\widetilde\nabla Z|)\widetilde\nabla Z-J(\widetilde\nabla \Phi)\widetilde\nabla(Y-Z)=\mathbf{R}=\mathbf{R}_1+\mathbf{R}_2+\mathbf{R}_3+\mathbf{R}_4,
				\end{aligned}
			\end{equation*}
			with
			\begin{align*}
				&\mathbf{R}_1=(\mathbf{T}(|\widetilde\nabla Y|)-\mathbf{T}(\widetilde\nabla\Phi))\widetilde\nabla(Y-Z)+(\mathbf{T}(|\widetilde\nabla Y|)-\mathbf{T}(|\widetilde\nabla Z|))\widetilde\nabla(Z-\Phi),\\
				&\mathbf{R}_2=\big(\mathbf{T}(|\widetilde\nabla Y|)-\mathbf{T}(|\widetilde\nabla Z|)-\mathbf{T}'(|\widetilde\nabla\Phi|)(|\widetilde\nabla Y|-|\widetilde\nabla Z|)\big)\widetilde\nabla\Phi,\\
				&\mathbf{R}_3=(\mathbf{T}'(|\widetilde\nabla Y|-\mathbf{T}'(|\widetilde\nabla\Phi|))(|\widetilde\nabla Y|-|\widetilde\nabla Z|),\\
				&\mathbf{R}_4=\big(\mathbf{T}'(|\widetilde\nabla\Phi|)(|\widetilde\nabla Y|-|\widetilde\nabla Z|)\big)\widetilde\nabla\Phi-J(\widetilde\nabla\Phi)\widetilde\nabla(Y-Z),
			\end{align*}
			and by taking $\mathbf{G}=G(\Phi(\widehat{\boldsymbol{x}})-\Phi(\widehat{\boldsymbol{y}}))$, applying Remark \ref{rmk3dpesk}, we have
			\begin{equation*}
				\|N_2(Y,\Phi)-N_2(Z,\Phi)\|_{\dot C^\kappa}\lesssim \|\Phi\|_{C^3}(1+\|\Phi\|_{C^3})\|\mathbf{R}\|_{\dot C^{\kappa}}.
			\end{equation*}
			By Lemma \ref{Hol} we have
			\begin{equation*}
				\|\mathbf{R}\|_{\dot C^{\kappa}}\lesssim t^{-\kappa}\|\widetilde\nabla W\|_{T}(\|\widetilde\nabla\vec{V}-\Phi\|_{T}+T^\kappa\|\Phi\|_{C^3})(1+(\|\widetilde\nabla(\vec{V}-\Phi)\|_{T}+\|\Phi\|_{C^3})).
			\end{equation*}
			So we have
			\begin{equation*}
				\begin{aligned}
					&\|N_2(Y,\Phi)-N_2(Z,\Phi)\|_{\dot C^\kappa}\\
					&\lesssim t^{-\kappa}(1+\|\Phi\|_{C^3})^2\|W\|_{Z_T}(\|\vec{V}-\Phi\|_{Z_T}+T^\kappa\|\Phi\|_{C^3})(1+(\|\vec{V}-\Phi\|_{Z_T}+\|\Phi\|_{C^3})).
				\end{aligned}
			\end{equation*}
			For higher order derivatives, we use \eqref{3dpeskerd}, integrating by parts and taking derivatives on $\mathbf{R}$ to obtain
			\begin{equation*}
				\begin{aligned}
					&\|N_2(Y,\Phi)-N_2(Z,\Phi)\|_{\dot C^{m+\kappa}}\\
					&\quad\lesssim t^{-m-\kappa}\|\Phi\|_{C^{1}}\|W\|_{Z_T}(\|\vec{V}-\Phi\|_{Z_T}+T^\kappa\|\Phi\|_{C^{m+3}})(1+(\|\vec{V}-\Phi\|_{Z_T}+T\|\Phi\|_{C^{m+3}}))^{m+2}.
				\end{aligned}
			\end{equation*}
		\end{proof}\\
		Finally for $N_3(\Phi)$, since $\Phi$ is smooth, it is easy to obtain that
		\begin{lemma}\label{3dpeslemn3}
			Let $N_3(\Phi)$ be defined in \eqref{3dpesnlt}. For any $T>0$, we have
			\begin{equation*}
				\begin{aligned}
					\|N_3(\Phi)\|_{C^{m+\kappa}}\lesssim \|\Phi\|_{C^{m+3}}(1+\|\Phi\|_{C^{m+3}})^{m+2}.
				\end{aligned}
			\end{equation*}
		\end{lemma}
		\subsection{Estimates of nonlinear terms in the general Muskat equation}
		This subsection is devoted to estimate force terms of \eqref{infaeq11}.
		For a family of functions $\{h_\b:\mathbb{K}\to \mathbb{R}\}_{\b\in \mathbb{R}^d}$, where $\mathbb{K}$ is a subset of $\mathbb{R}^d$ or $\mathbb{R}^{d+1}$,  denote 
		\begin{equation}\label{defnormbb}
			\begin{aligned}
				&	|||h_\b|||_{a,\mathbb{K}}:=\sup_{x,y\in\mathbb{K}}\frac{|(h_\b(x)- h_\b(y))|_{\b=\nabla\phi(x)}|}{|x-y|^{a}},\ \ \  \text{for} \ a\in(0,1).
			\end{aligned}
		\end{equation}
		where $\phi$ is a smooth function defined in \eqref{defphigm}. And we will drop the domain if $\mathbb{K}$ is the whole space $\mathbb{R}^d$ or $\mathbb{R}^{d+1}$.
		We remark that in this section, we write $\sup_\b$ for short to denote $\sup_{\b\in\operatorname{Ran}(\nabla\phi)}$.
		
		We first estimate the explicit term $\H^{ex}_\b[f]$. 
		We have the following result.
		\begin{lemma}\label{lemH1}
			Let $\H^{ex}_\b[f]$ be as defined in \eqref{defH}, the norm $\|\cdot\|_{T}$ be as defined in \eqref{defnorgm} with constants $\varkappa,m,\kappa$ fixed in \eqref{consgm}. Then for any $T\in(0,\frac{1}{2})$, there hold
			\begin{align}
				&\sup_{t\in[0,T]}|\log t|^\varkappa  (t^\frac{\kappa}{3}|||(\H^{ex}_\b[f])(t)|||_{\kappa-2}+t^\frac{m+\kappa}{3}|||(\nabla^m\H^{ex}_\b[f])(t)|||_{\kappa-2})\nonumber\\
				&\quad\quad\lesssim\|f-\phi\|_T(\|f-\phi\|_T+T^\frac{1}{20}\|\nabla\phi\|_{C^{m+3}})(1+\|(f,\phi)\|_T)^{2m+5}+T^\frac{1}{5}(1+\|\nabla\phi\|_{C^{n+3}})^{m+3}, \label{m}\\
				&\sup_{t\in[0,T]}|\log t|^\varkappa  \left(t^\frac{\kappa}{3}||| (\H^{ex}_\b[f_1]-\H^{ex}_\b[f_2])(t)|||_ {\kappa-2}+t^\frac{m+\kappa}{3}||| \nabla^m(\H^{ex}_\b[f_1]-\H^{ex}_\b[f_2])(t)|||_{\kappa-2}\right)\nonumber\\
				&\quad\quad\lesssim \|f_1-f_2\|_{T}(\|(f_1-\phi,f_2-\phi)\|_{T}+T^\frac{1}{20}\|\nabla \phi\|_{C^{m+3}})(1+\|( f_1,f_2,\phi)\|_{T})^{2m+5}.\label{dif}
			\end{align}
		\end{lemma}
		\begin{proof}
			For brevity, we drop the constants $\frac{\tilde c_d}{\mu_++\mu_-}$ and $\frac{c_-}{\mu_-}$  in \eqref{defH}. Then
			\begin{align*}
				\H^{ex}_\b[f](x)=	&\int _{\mathbb{R}^d}(\A(\b,\alpha)-\A(\nabla\phi(x),\alpha)):\delta_\alpha \nabla^2 f(x)\frac{d\alpha}{|\alpha|^{d+1}}\\
				&	+		\int_{\mathbb{R}^d}\frac{(\mathbf{B}(\nabla f(x))-\mathbf{B}(\mathbf{b})):\delta_\alpha \nabla^2f(x)}{\langle\hat \alpha \cdot \mathbf{b}\rangle^{d+1}}\frac{d\alpha}{|\alpha|^{d+1}}\\
				&	+{\langle\nabla f(x)\rangle^2 }\int_{\mathbb{R}^d}\frac{\delta_\alpha B(\nabla f): \nabla^2f(x-\alpha)}{\langle\hat \alpha \cdot \mathbf{b}\rangle^{d+1}}\frac{d\alpha}{|\alpha|^{d+1}}\\
				&	+{\langle\nabla f(x)\rangle^2 }\left(\frac{\nabla f(x)}{\langle\nabla f(x)\rangle^2 }-\frac{\mathbf{b}}{\langle \mathbf{b}\rangle^2}\right)\cdot \nabla\mathcal{G}[f](x)\\
				:=&\mathcal{P}_\b[f](x)+\mathcal{Q}_\b[f](x)+\mathcal{W}_\b[f](x)+\mathcal{U}_\b[f](x).
			\end{align*}
			where  $B(\cdot)$, $\mathbf{B}(\cdot)$, and  $\mathcal{G}[\cdot]$ are defined in  \eqref{ndefB}, \eqref{ndefbfb}, and  \eqref{ndefG}, respectively.\\
			Denote 
			\begin{align*}
				&\|h\|_{T,*}:=\sup_{t\in[0,T]}(t^\frac{1}{15}|\log t|^\varkappa  \|\nabla h(t)\|_{\dot C^{\frac{1}{5}}}+t^\frac{m+\kappa}{3}|\log t|^\varkappa  \|\nabla h(t)\|_{\dot C^{m+\kappa}}).
			\end{align*}
			We first prove \eqref{dif}. For simplicity, we denote $a=\kappa-2$ in this proof. It suffices to prove that for any $n\in\mathbb{N}$, $n\leq m$, 
			\begin{equation}\label{difne}
				\begin{aligned}
					&\sup_{t\in[0,T]}|\log t|^\varkappa  t^\frac{n+\kappa}{3}||| \nabla^n(\H^{ex}_\b[f_1]-\H^{ex}_\b[f_2])(t)|||_{a}\\\
					&\ \ \quad\quad\quad\quad\quad\quad\quad\quad\lesssim \|f_1-f_2\|_{T}(\|(f_1,f_2)\|_{T,*}+T^\frac{1}{6}\|\nabla\phi\|_{C^{m+3}})(1+\|(f_1,f_2,\phi)\|_{T})^{2m+5}.
				\end{aligned} 
			\end{equation}
			In fact, \eqref{dif} follows by \eqref{difne} and the fact that 
			\begin{align*}
				\|(f_1,f_2)\|_{T,*}\lesssim \|(f_1-\phi,f_2-\phi)\|_{T}+\|\phi\|_{T,*}\lesssim \|(f_1-\phi,f_2-\phi)\|_{T}+T^\frac{1}{20}\|\nabla\phi\|_{C^{m+3}}.
			\end{align*}
			Now we start to prove \eqref{difne}.\\  
			1) Estimate $\mathcal{P}_\b$. \\
			By the symmetry $\A(\cdot,\alpha)=\A(\cdot,-\alpha)$,  we have 
			\begin{align*}
				(	\mathcal{P}_\b[f_1]-\mathcal{P}_\b[f_2])(x)=&\int _{\mathbb{R}^d}(\A(\b,\alpha)-\A(\nabla\phi(x),\alpha)):\delta_\alpha \nabla^2 (f_1-f_2)(x)\frac{d\alpha}{|\alpha|^{d+1}}\\
				=&\frac{1}{2}\int _{\mathbb{R}^d}(\A(\b,\alpha)-\A(\nabla\phi(x),\alpha)):\mathcal{O}_\alpha \nabla^2 (f_1-f_2)(x)\frac{d\alpha}{|\alpha|^{d}},
			\end{align*}
			where  $\mathcal{O}_\alpha g(x)=\Delta_\alpha g(x)+\Delta_{-\alpha}g(x)$.
			Then for any $\beta\in\mathbb{R}^d$, and any $0\leq n\leq m$,
			\begin{align*}
				|\delta_\beta 	\nabla^n(	\mathcal{P}_\b[f_1]-\mathcal{P}_\b[f_2])(x)|\lesssim& \sum_{n_1=0}^n\int _{\mathbb{R}^d}\left|\nabla ^{n_1}(\A(\b,\alpha)-\A(\nabla\phi(x),\alpha)):\delta_\beta\mathcal{O}_\alpha  \nabla^{2+n-n_1} (f_1-f_2)(x)\right|\frac{d\alpha}{|\alpha|^{d}}\\ &+\sum_{n_1=0}^n\int _{\mathbb{R}^d}|\delta_\beta\nabla ^{n_1}\A(\nabla\phi(\cdot ),\alpha)(x)\mathcal{O}_\alpha \nabla^{2+n-n_1} (f_1-f_2)(x-\beta)|\frac{d\alpha}{|\alpha|^{d}}\\
				:=&\mathrm{I}_{1,\b}+\mathrm{I}_2.
			\end{align*}
			We first consider $I_{1,\b}$.
			Note that $I_{1,\b}|_{\b=\nabla\phi(x)}$ vanishes if $n=0$. Moreover, if $n>0$,  the contribution of $n_1=0$ turns out to be $0$ if we take $\b=\nabla\phi(x)$. Hence  it suffices to consider $1\leq n_1\leq n$. We have 
			\begin{align*}
				\|\I_{1,\b}|_{\b=\nabla\phi(x)}\|_{L^\infty}\lesssim \sum_{n_1=1}^n	\|\nabla ^{n_1}(\A(\b,\alpha)-\A(\nabla\phi(x),\alpha))\|_{L^\infty_{x,\alpha}}\int _{\mathbb{R}^d}\left\|\delta_\beta\mathcal{O}_\alpha  \nabla^{2+n-n_1} (f_1-f_2)(x)\right\|_{L^\infty}\frac{d\alpha}{|\alpha|^{d}}.
			\end{align*}
			Then by Lemma \ref{lemcom},
			\begin{align*}
				\|\nabla ^{n_1}(\A(\b,\alpha)-\A(\nabla\phi(x),\alpha))\|_{L^\infty_{x,\alpha}}=	\|\nabla ^{n_1}\A(\nabla\phi(x),\alpha)\|_{L^\infty_{x,\alpha}}\lesssim (1+\|\nabla\phi\|_{C^n})^{n+2}.
			\end{align*}
			Moreover, from Lemma \ref{lemoalin} we get 
			\begin{align*}
				\int_{\mathbb{R}^d}	\|\delta_\beta\mathcal{O}_\alpha  \nabla^{2+n-n_1} (f_1-f_2)\|_{L^\infty}\frac{d\alpha}{|\alpha|^d}\lesssim  |\beta|^a |\log t|^{-\varkappa}t^{-\frac{2+n-n_1+a}{3}}\|f_1-f_2\|_{T,*}.
			\end{align*}
			Hence we get 
			\begin{equation}\label{i1b}
				\begin{aligned}
					\|\I_{1,\b}|_{\b=\nabla\phi(x)}\|_{L^\infty}	\lesssim &|\beta|^a|\log t|^{-\varkappa}t^{-\frac{2+n-1+a}{3}}\|f_1-f_2\|_{T,*}(1+\|\nabla\phi\|_{C^n})^{n+2} \|f_1-f_2\|_{T,*}.
				\end{aligned}
			\end{equation}
			For $\mathrm{I_2}$, we have 
			\begin{align*}
				\|\delta_\beta\nabla ^{n_1}\A(\nabla\phi(\cdot ),\alpha)(x)\|_{L^\infty_{x,\alpha}}&\lesssim |\beta|^a\|\nabla\phi\|_{C^{n+1}},\end{align*}
			Apply Lemma \ref{lemoalin} again, we get 
			\begin{align*}
				\int_{\mathbb{R}^d}	\|\mathcal{O}_\alpha \nabla^{2+n-n_1} (f_1-f_2)\|_{L^\infty}\frac{d\alpha}{|\alpha|^d}&\lesssim \|\nabla(f_1-f_2)\|_{\dot C^{2+n-n_1+\frac{1}{2}}}^\frac{1}{2}\|\nabla(f_1-f_2)\|_{\dot C^{2+n-n_1-\frac{1}{2}}}^\frac{1}{2}\\
				&\lesssim |\log t|^{-\varkappa}t^{-\frac{2+n-n_1}{3}}\|f_1-f_2\|_{T,*}.
			\end{align*}
			Hence 
			\begin{align*}
				\|	\mathrm{I}_2\|_{L^\infty}&\lesssim |\beta|^a	\|\nabla ^{n_1}\A(\nabla\phi(\cdot ),\alpha)(x)\|_{L^\infty_{\alpha}\dot C^a_x}\sum_{n_1=0}^n \int_{\mathbb{R}^d}\|\mathcal{O}_\alpha \nabla^{2+n-n_1} (f_1-f_2)\|_{L^\infty}\frac{d\alpha}{|\alpha|^d}\\
				&\lesssim|\beta|^a\|\nabla\phi\|_{C^{n+1}}|\log t|^{-\varkappa}t^{-\frac{2+n-n_1}{3}}\|f_1-f_2\|_{T,*}.
			\end{align*}
			Combining this with \eqref{i1b}, we have
			\begin{align}\label{pbhol}
				\|(\delta_\beta 	\nabla^n(	\mathcal{P}_\b[f_1]-\mathcal{P}_\b[f_2]))|_{\b=\nabla\phi(x)}\|_{L^\infty}\lesssim |\beta|^a(1+\|\nabla\phi\|_{C^n})^{n+2}|\log t|^{-\varkappa}t^{-\frac{2+n}{3}}\|f_1-f_2\|_{T,*}.
			\end{align}
			2) Estimate $\mathcal{Q}_\b$. \\
			By symmetry, one has 
			\begin{align*}
				\mathcal{Q}_\b[f](x)=	\frac{1}{2}\int_{\mathbb{R}^d}\frac{(\mathbf{B}(\nabla f(x))-\mathbf{B}(\mathbf{b})):\mathcal{O}_\alpha \nabla^2f(x)}{\langle\hat \alpha \cdot \mathbf{b}\rangle^{d+1}}\frac{d\alpha}{|\alpha|^{d}}.
			\end{align*}
			Hence
			\begin{align*}
				\mathcal{Q}_\b[f_1]-\mathcal{Q}_\b[f_2]
				=&\frac{1}{2}\int_{\mathbb{R}^d}\frac{(\mathbf{B}(\nabla f_1(x))-\mathbf{B}(\nabla f_2(x))):\mathcal{O}_\alpha \nabla^2f_1(x)}{\langle\hat \alpha \cdot \mathbf{b}\rangle^{d+1}}\frac{d\alpha}{|\alpha|^{d}}\\
				&+\frac{1}{2}\int_{\mathbb{R}^d}\frac{(\mathbf{B}(\nabla f_2(x))-\mathbf{B}(\mathbf{b})):\mathcal{O}_\alpha \nabla^2(f_1-f_2)(x)}{\langle\hat \alpha \cdot \mathbf{b}\rangle^{d+1}}\frac{d\alpha}{|\alpha|^{d}}\\
				:=&\mathcal{Q}_{1,\b}+\mathcal{Q}_{2,\b}.
			\end{align*}
			Using the formula $\delta_\beta (h_1h_2)(x)=\delta_\beta h_1(x) h_2(x)+h_1(x-\beta)\delta_\beta h_2(x)$, we have 
			\begin{equation}\label{q1b}
				\begin{aligned}
					&\left\|(\delta_\beta\nabla^n\mathcal{Q}_{1,\b})|_{\b=\nabla\phi(x)}\right\|_{L^\infty}\\
					&\quad\quad \lesssim  \sum_{n_1+n_2=n}\left(|\beta|^a\left\|\nabla^{n_1} (  \mathbf{B}(\nabla f_1)- \mathbf{B}(\nabla f_2)) \right\|_{\dot C^a}\int _{\mathbb{R}^d}\left\|\mathcal{O}_\alpha \nabla^{2+n_2} f_1\right\|_{L^\infty}\frac{d\alpha}{|\alpha|^{d}}\right.\\
					&\quad\quad\quad\quad\quad\quad \quad\ \ \  \left.+\left\|\nabla^{n_1} (  \mathbf{B}(\nabla f_1)- \mathbf{B}(\nabla f_2))\right\|_{L^\infty}\int_{\mathbb{R}^d} \left\|\delta_\beta\mathcal{O}_\alpha \nabla^{2+n_2} f_1\right\|_{L^\infty}\frac{d\alpha}{|\alpha|^{d}}\right).
				\end{aligned}
			\end{equation}
			By \eqref{mainint21} and \eqref{mainint22} in Lemma \ref{lemcom},  one can check for any $k\in\mathbb{N}, k\leq m$,
			\begin{align*}
				\left\|\nabla^{k} (  \mathbf{B}(\nabla f_1)- \mathbf{B}(\nabla f_2)) \right\|_{L^\infty}&\lesssim \sum_{k_1=0}^k\|\nabla (f_1-f_2)\|_{\dot C^{k_1}}\left(\|(\nabla f_1,\nabla f_2)\|_{\dot C^1}^{k-k_1}+\|(\nabla f_1,\nabla f_2)\|_{\dot C^{k-k_1}}\right)\\
				&\lesssim t^{-\frac{k}{3}}\|f_1-f_2\|_T(1+\|(f_1,f_2)\|_T)^{k},\\
				\left\|\nabla^{k} (  \mathbf{B}(\nabla f_1)- \mathbf{B}(\nabla f_2)) \right\|_{\dot C^a}&\lesssim \sum_{k_1=0}^k\|\nabla (f_1-f_2)\|_{\dot C^{k_1+a}}\left(\|(\nabla f_1,\nabla f_2)\|_{\dot C^1}^{k-k_1}+\|(\nabla f_1,\nabla f_2)\|_{\dot C^{k-k_1}}\right)\\
				&\quad+\sum_{k_1=0}^k\|\nabla (f_1-f_2)\|_{\dot C^{k_1}}\left(\|(\nabla f_1,\nabla f_2)\|_{\dot C^a}^{\frac{k-k_1+a}{a}}+\|(\nabla f_1,\nabla f_2)\|_{\dot C^{k-k_1+a}}\right)\\
				&\lesssim t^{-\frac{k+a}{3}}\|f_1-f_2\|_T(1+\|(f_1,f_2)\|_T)^{k}.
			\end{align*}
			Moreover, applying Lemma \ref{lemoalin}, we obtain for any $k\in\mathbb{N}, k\leq m$,
			\begin{equation}\label{ososf1}
				\begin{aligned}
					&\int _{\mathbb{R}^d}\left\|\mathcal{O}_\alpha \nabla^{2+k} f_1\right\|_{L^\infty}\frac{d\alpha}{|\alpha|^{d}}\lesssim \|\nabla f_1\|_{\dot C^{2+k+\frac{1}{2}}}^\frac{1}{2}\|\nabla f_1\|_{\dot C^{2+k-\frac{1}{2}}}^\frac{1}{2}\lesssim  |\log t|^{-\varkappa}t^{-\frac{2+k}{3}}\|f_1\|_{T,*},\\
					&\int_{\mathbb{R}^d} \left\|\delta_\beta\mathcal{O}_\alpha \nabla^{2+k} f_1\right\|_{L^\infty}\frac{d\alpha}{|\alpha|^{d}}\lesssim |\beta|^a\|\nabla f_1\|_{\dot C^{2+k+a}}\lesssim |\beta|^a |\log t|^{-\varkappa}t^{-\frac{2+k+a}{3}}\|f_1\|_{T,*}.
				\end{aligned}
			\end{equation}
			Putting the above estimates in \eqref{q1b}, we get 
			\begin{align}\label{dq1b}
				\left\|(\delta_\beta\nabla^n\mathcal{Q}_{1,\b})|_{\b=\nabla\phi(x)}\right\|_{L^\infty}\lesssim |\beta|^a|\log t|^{-\varkappa}t^{-\frac{n+a}{3}}\|f_1-f_2\|_T\|f_1\|_{T,*}(1+\|(f_1,f_2)\|_T)^{n}.
			\end{align} 
			Then we deal with $\mathcal{Q}_{2,\b}$.
			\begin{align*}
				&\left\|(\delta_\beta\nabla^n\mathcal{Q}_{2,\b})|_{\b=\nabla\phi(x)}\right\|_{L^\infty}\\
				&\quad\quad \lesssim  \sum_{n_1=0}^n\left( \int_{\mathbb{R}^d}\left\|(\delta_\beta\nabla^{n_1}(  \mathbf{B}(\nabla f_2)- \mathbf{B}(\b)))|_{\b=\nabla\phi(x)}\right\|_{L^\infty}\left\|\mathcal{O}_\alpha \nabla^{2+n-n_1} (f_1-f_2)\right\|_{L^\infty}\frac{d\alpha}{|\alpha|^{d}}\right.\\
				&\quad\quad\quad\quad\quad\quad\quad\ \ \ +\left.\int_{\mathbb{R}^d}\left\|(\nabla^{n_1} (\mathbf{B}(\nabla f_2)- \mathbf{B}(\b)))|_{\b=\nabla\phi(x)}\right\|_{L^\infty}\left\|\delta_\beta\mathcal{O}_\alpha \nabla^{2+n-n_1} (f_1-f_2)\right\|_{L^\infty}\frac{d\alpha}{|\alpha|^{d}}\right)\\
				&\quad\quad \lesssim  |\beta|^a\sum_{n_1=0}^n\left(\left\|\nabla^{n_1} \mathbf{B}(\nabla f_2)\right\|_{\dot C^{a}} \int_{\mathbb{R}^d}\left\|\mathcal{O}_\alpha \nabla^{2+n-n_1} (f_1-f_2)\right\|_{L^\infty}\frac{d\alpha}{|\alpha|^{d}}\right)\\
				&\quad\quad\quad\quad\ \ \ +\left\|(\mathbf{B}(\nabla f_2)- \mathbf{B}(\nabla\phi))\right\|_{L^\infty}\int_{\mathbb{R}^d}\left\|\delta_\beta\mathcal{O}_\alpha \nabla^{2+n} (f_1-f_2)\right\|_{L^\infty}\frac{d\alpha}{|\alpha|^{d}}\\
				&\quad\quad\quad\quad\ \ \ +\sum_{n_1=1}^n\left(\left\|\nabla^{n_1}\mathbf{B}(\nabla f_2)\right\|_{L^\infty}\int_{\mathbb{R}^d}\left\|\delta_\beta\mathcal{O}_\alpha \nabla^{2+n-n_1} (f_1-f_2)\right\|_{L^\infty}\frac{d\alpha}{|\alpha|^{d}}\right).
			\end{align*}
			By Lemma \ref{lemcom}, we get
			\begin{align*}
				&\left\|\nabla^{n_1} \mathbf{B}(\nabla f_2)\right\|_{\dot C^{a}}\lesssim (\|\nabla f_2\|_{\dot C^{n_1+a}}+\|\nabla f_2\|_{\dot C^{a}}^\frac{n_1+a}{a})(1+\|\nabla f_2\|_{L^\infty})^2\lesssim t^{-\frac{n_1+a}{3}}\|f_2\|_{T,*},\\
				&\left\|(\mathbf{B}(\nabla f_2)- \mathbf{B}(\nabla\phi))\right\|_{L^\infty}\lesssim \|\nabla (f_2-\phi)\|_{l^\infty}(1+\|\nabla f_2\|_{L^\infty}+\|\nabla\phi\|_{L^\infty})^2\lesssim \|f_2-\phi\|_T(1+\|f_2\|_T+\|\nabla\phi\|_{L^\infty})^2,\\
				&\left\|\nabla^{n_1}\mathbf{B}(\nabla f_2)\right\|_{L^\infty}\lesssim (\|\nabla f_2\|_{\dot C^{n_1}}+\|\nabla f_2\|_{\dot C^{2}}^{n_1})(1+\|\nabla f_2\|_{L^\infty})^2\lesssim \|f_2\|_{T,*}(1+\|f_2\|_T)^{n_1+2}.
			\end{align*}
			Combining this with \eqref{ososf1}, we conclude that \begin{align}\label{dq2b}
				\left\|(\delta_\beta\nabla^n\mathcal{Q}_{2,\b})|_{\b=\nabla\phi(x)}\right\|_{L^\infty}\lesssim |\log t|^{-\varkappa} |\beta|^a(\|f_2\|_{T,*}+\|f_2-\phi\|_T)(1+\|f_2\|_T)^{n+2}.
			\end{align}
			It follows from \eqref{dq1b} and \eqref{dq2b} that 
			\begin{equation}\label{difP}
				\begin{aligned}
					&\|\delta_\beta\nabla^n (\mathcal{Q}_\b[f_1]-\mathcal{Q}_\b[f_2])\|_{L^\infty}\\
					&\lesssim |\beta|^{a}t^{-\frac{n+a}{3}}|\log t|^{-\varkappa } \|f_1-f_2\|_{T}(\|f_1\|_{T,*}+\|f_2-\phi\|_T)(1+\|(f_1,f_2,\phi)\|_T)^{n+4},\ \ \ \forall\ n\in\mathbb{N}, n\leq m.
				\end{aligned}
			\end{equation}
			3) Estimate $\mathcal{W}_\b$.\\
			We can write
			\begin{align*}
				\mathcal{W}_\b[f](x)&=	{\langle \nabla f(x)\rangle ^2}\int_{\mathbb{R}^d}\frac{\delta_\alpha B(\nabla f(x)):\nabla^2 f(x-\alpha)}{\langle\hat \alpha \cdot \b\rangle ^{d+1}} \frac{d\alpha}{|\alpha|^{d+1}}\\
				&=	{\langle \nabla f(x)\rangle ^2}\int_{\mathbb{R}^d}\frac{\delta_\alpha B(\nabla f(x)):\nabla^2 f(x)}{\langle\hat \alpha \cdot \b\rangle ^{d+1}} \frac{d\alpha}{|\alpha|^{d+1}}-	{\langle \nabla f(x)\rangle ^2}\int_{\mathbb{R}^d}\frac{\delta_\alpha B(\nabla f(x)):\delta_\alpha\nabla^2 f(x)}{\langle\hat \alpha \cdot \b\rangle ^{d+1}} \frac{d\alpha}{|\alpha|^{d+1}}\\
				&:=\mathcal{W}_{1,\b}[f](x)+\mathcal{W}_{2,\b}[f](x).
			\end{align*}
			First, we claim that for any $k\in\mathbb{N}$, $k\leq m+1$,
			\begin{equation}\label{esA}
				\begin{aligned}
					&\sup_{t\in[0,T]}|\log t|^\varkappa  t^{\frac{k+2}{3}}\|\nabla^k(\mathcal{W}_\b[f_1]-\mathcal{W}_\b[f_2])\|_{L^\infty})\\
					&\quad\quad\quad\quad\quad \lesssim \|f_1-f_2 \|_{T} \|(f_1,f_2)\|_{T,*}(1+\|(f_1,f_2,\phi)\|_T)^{k+4}.
				\end{aligned}
			\end{equation}
			By symmetry, we can write 
			\begin{align*}
				\mathcal{W}_{1,\b}[f](x)&=\frac{1}{2}	{\langle \nabla f(x)\rangle ^2}\int_{\mathbb{R}^d}\frac{\mathcal{O}_{\alpha} B(\nabla f(x)):\nabla^2 f(x)}{\langle\hat \alpha \cdot \b\rangle ^{d+1}} \frac{d\alpha}{|\alpha|^{d}}.
			\end{align*}
			Then 
			\begin{align*}
				\mathcal{W}_{1,\b}[f_1](x)-	\mathcal{W}_{1,\b}[f_2](x)=&\frac{1}{2}\langle \nabla f_1(x)\rangle ^2\int_{\mathbb{R}^d}\frac{ \big(\mathcal{O}_{\alpha} (B(\nabla f_1(x))-B(\nabla f_2(x))):\nabla^2 f_1(x)}{\langle\hat \alpha \cdot \b \rangle ^{d+1}}\frac{d\alpha}{|\alpha|^{d}}\\&
				+\frac{1}{2}\langle \nabla f_1(x)\rangle ^2\int_{\mathbb{R}^d}\frac{\mathcal{O}_{\alpha}B(\nabla f_2(x)):\nabla^2 (f_1-f_2)(x)}{\langle\hat \alpha \cdot \b\rangle ^{d+1}}\frac{d\alpha}{|\alpha|^{d}}\\
				&+\frac{1}{2}(\langle \nabla f_1(x)\rangle ^2-\langle \nabla f_2(x)\rangle ^2)\int_{\mathbb{R}^d}\frac{\mathcal{O}_{\alpha} B(\nabla f_2(x)):\nabla^2 f_2(x)}{\langle\hat \alpha \cdot \b\rangle ^{d+1}} \frac{d\alpha}{|\alpha|^{d}}\\
				:=&\textit{II}_1+\textit{II}_2+\textit{II}_3.
			\end{align*}
			For any $k\in\mathbb{N}$, $k\leq m+1$,	we have 
			\begin{align*}
				&\|\nabla^k\textit{II}_1\|_{L^\infty}\lesssim \sum_{l=0}^k\left( \int_{\mathbb{R}^d}\|\mathcal{O}_\alpha \nabla^l(B(\nabla f_1)-B(\nabla f_2))\|_{L^\infty}\frac{d\alpha}{|\alpha|^{d}} \left\|\nabla^{k-l}\left({\langle \nabla f_1\rangle ^2\nabla^2f_1}\right)\right\|_{L^\infty}\right),\\
				&\|\nabla^k\textit{II}_2\|_{L^\infty}\lesssim \sum_{l=0}^k\left( \int_{\mathbb{R}^d}\|\mathcal{O}_\alpha \nabla^l B(\nabla f_2)\|_{L^\infty}\frac{d\alpha}{|\alpha|^{d}} \left\|\nabla^{k-l}\left({\langle \nabla f_1\rangle ^2\nabla^2(f_1-f_2)}\right)\right\|_{L^\infty}\right),\\
				&\|\nabla^k\textit{II}_3\|_{L^\infty}\lesssim \sum_{l=0}^k \left(\int_{\mathbb{R}^d}\|\mathcal{O}_\alpha \nabla^l B(\nabla f_2)\|_{L^\infty}\frac{d\alpha}{|\alpha|^{d}} \left\|\nabla^{k-l}\left({(\langle \nabla f_1\rangle ^2-\langle \nabla f_2\rangle ^2)\nabla^2f_2}\right)\right\|_{L^\infty}\right).
			\end{align*}
			By Lemma \ref{lemoalin},
			\begin{align*}
				&\|\nabla^k\textit{II}_1\|_{L^\infty}\lesssim \sum_{l=0}^k\left( \prod_{j=l,l+1}\| B(\nabla f_1)-B(\nabla f_2)\|_{\dot C^{j+\frac{1}{2}}} ^\frac{1}{2} \left\|\nabla^{k-l}\left({\langle \nabla f_1\rangle ^2\nabla^2f_1}\right)\right\|_{L^\infty}\right),\\
				&\|\nabla^k\textit{II}_2\|_{L^\infty}\lesssim \sum_{l=0}^k\left(  \prod_{j=l,l+1}\| B(\nabla f_2)\|_{\dot C^{j+\frac{1}{2}}} ^\frac{1}{2} \left\|\nabla^{k-l}\left({\langle \nabla f_1\rangle ^2\nabla^2(f_1-f_2)}\right)\right\|_{L^\infty}\right),\\
				&\|\nabla^k\textit{II}_3\|_{L^\infty}\lesssim \sum_{l=0}^k\left(  \prod_{j=l,l+1}\| B(\nabla f_2)\|_{\dot C^{j+\frac{1}{2}}} ^\frac{1}{2} \left\|\nabla^{k-l}\left({(\langle \nabla f_1\rangle ^2-\langle \nabla f_2\rangle ^2)\nabla^2f_2}\right)\right\|_{L^\infty}\right).
			\end{align*}
			By Lemma \ref{lemcom}, we have
			\begin{align*}
				\| B(\nabla f_1)-B(\nabla f_2)\|_{\dot C^{j+\frac{1}{2}}}\lesssim& \sum_{j'=0}^j\left\{\|\nabla(f_1-f_2)\|_{\dot C^{j'+\frac{1}{2}}}(\|\nabla(f_1,f_2)\|_{\dot C^1}^{j-j'}+\|\nabla(f_1,f_2)\|_{\dot C^{j-j'}})\right.\\
				&\quad\quad\left. + \|\nabla(f_1-f_2)\|_{\dot C^{j'}}(\|\nabla(f_1,f_2)\|_{\dot C^\frac{1}{2}}^{2(j-j')+1}+\|\nabla(f_1,f_2)\|_{\dot C^{j-j'+\frac{1}{2}}}) \right\}\\
				\lesssim  &t^{-\frac{j}{3}-\frac{1}{6}}|\log t|^{-\varkappa } \|f_1-f_2\|_T(1+\|(f_1,f_2)\|_{T})^{2j+1},
			\end{align*}
			and 
			\begin{align*}
				\| B(\nabla f_2)\|_{\dot C^{j+\frac{1}{2}}} \lesssim \|\nabla f_2\|_{\dot C^\frac{1}{2}}^{2j+1}+\|\nabla f_2\|_{\dot C^{j+\frac{1}{2}}}\lesssim t^{-\frac{j}{3}-\frac{1}{6}}|\log t|^{-\varkappa }  \|f_2\|_{T,*}(1+\|f_2\|_T)^{2j+1}.
			\end{align*}
			Moreover, by Lemma \ref{lemcom}, we obtain that for any $k\in\mathbb{N}, k\leq m+1$,
			\begin{align*}
				\left\|\nabla^{k}\left(\langle \nabla f_1\rangle ^2\nabla^2f_1\right)\right\|_{L^\infty}
				&\lesssim (1+\|\nabla f_1\|_{L^\infty})^2(\|\nabla f_1\|_{\dot C^1}^{k+1}+\|\nabla f_1\|_{\dot C^{k+1}})\\
				&\lesssim t^{-\frac{k+1}{3}}|\log t|^{-\varkappa } \|f_1\|_{T,*}(1+\|f_1\|_T)^{k+3}.
			\end{align*}
			Similarly, we have 
			\begin{align*}
				&\left\|\nabla^{k}\left({\langle \nabla f_1\rangle ^2\nabla^2(f_1-f_2)}\right)\right\|_{L^\infty}+\left\|\nabla^{k}\left((\langle \nabla f_1\rangle ^2-\langle \nabla f_2\rangle ^2)\nabla^2f_2\right)\right\|_{L^\infty}\\
				&\quad\quad\quad\quad\quad\quad\quad\quad\quad\quad\lesssim t^{-\frac{k+1}{3}}|\log t|^{-\varkappa } \|f_1-f_2\|_{T}(1+\|(f_1,f_2)\|_T)^{k+3}.
			\end{align*}
			Collecting the estimates above, we obtain that 
			\begin{align*}
				\|\nabla^k(	\mathcal{W}_{\b}[f_1]-	\mathcal{W}_{\b}[f_2])(t)\|_{L^\infty}&\leq \|\nabla ^k \textit{II}_1\|_{L^\infty}+\|\nabla ^k \textit{II}_2\|_{L^\infty}+\|\nabla ^k \textit{II}_3\|_{L^\infty}\\
				&\lesssim t^{-\frac{k+2}{3}}|\log t|^{-\varkappa } \|f_1-f_2\|_T\|(f_1,f_2)\|_{T,*}(1+\|(f_1,f_2)\|_T)^{2k+5}.
			\end{align*}
			Then by interpolation, we get 
			\begin{equation}\label{difA}
				\begin{aligned}
					& \sup_{t\in[0,T]}t^{\frac{n+\kappa}{3}}|\log t|^\varkappa   \|(\delta_\beta\nabla^n(	\mathcal{W}_{\b}[f_1]-	\mathcal{W}_{\b}[f_2]))(t)|_{\b=\nabla\phi(x)}\|_{L^\infty}\\
					&\quad\quad\quad\quad\quad\lesssim |\beta|^a \|f_1-f_2\|_T\|(f_1,f_2)\|_{T,*}(1+\|(f_1,f_2)\|_T)^{2n+5}, \ \ \ \forall n\in\mathbb{N}, n\leq m. 
				\end{aligned}
			\end{equation}
			4) Estimate $\mathcal{U}_\b$.\\
			Note that 
			\begin{align*}
				\mathcal{U}_\b[f](x)=\left({\nabla f(x)}-\frac{\mathbf{b}\langle\nabla f(x)\rangle^2 }{\langle \mathbf{b}\rangle^2}\right)\cdot \nabla\mathcal{G}[f](x).
			\end{align*}
			Hence 
			\begin{align*}
				(\mathcal{U}_\b[f_1]-\mathcal{U}_\b[f_2])(x)=&\left(\mathrm{Id}-\frac{\mathbf{b}\otimes\nabla (f_1+f_2)(x)}{\langle \mathbf{b}\rangle^2}\right)\nabla(f_1-f_2)(x)\cdot \nabla\mathcal{G}[f_1](x)\\
				&+\left({\nabla f_2(x)}-\frac{\mathbf{b}\langle\nabla f_2(x)\rangle^2 }{\langle \mathbf{b}\rangle^2}\right)\cdot \nabla(\mathcal{G}[f_1]-\mathcal{G}[f_2])(x)\\
				:=&\tilde {\mathcal{U}}_{1,\b}(x)+\tilde {\mathcal{U}}_{2,\b}(x).
			\end{align*}
			Applying \eqref{mainint11} and \eqref{mainint12} in Lemma \ref{lemcom}, one can check that for any $n\in\mathbb{N}, n\leq m$,
			\begin{align*}
				&	\|\nabla \mathcal{G}[f]\|_{\dot C^{n}}\lesssim \|\nabla f\|_{\dot C^{2+n}}+\|\nabla f\|_{\dot C^1}^{n+2}\lesssim |\log t|^{-\varkappa}t^{-\frac{n+2}{3}}\|f\|_{T,*}(1+\|f\|_T)^{n+2},\\
				&	\|\nabla \mathcal{G}[f]\|_{\dot C^{n+a}}\lesssim \|\nabla f\|_{\dot C^{2+n+a}}+\|\nabla f\|_{\dot C^a}^\frac{n+2}{a}\lesssim |\log t|^{-\varkappa}t^{-\frac{n+2+a}{3}}\|f\|_{T,*}(1+\|f\|_T)^{n+10}.
			\end{align*}
			From this we get 
			\begin{align}\label{u1bhol}
				\|(\delta_\beta\nabla^n  \tilde {\mathcal{U}}_{1,\b})|_{\b=\nabla\phi(x)}\|_{L^\infty}\lesssim |\beta|^a |\log t|^{-\varkappa}t^{-\frac{n+2+a}{3}}\|f_1\|_{T,*}\|f_1-f_2\|_T(1+\|(f_1,f_2)\|_T)^{n+10}.
			\end{align}
			Then we consider ${\mathcal{U}}_{2,\b}$. Denote $F_{\b}={\nabla f_2(x)}-\frac{\mathbf{b}\langle\nabla f_2(x)\rangle^2 }{\langle \mathbf{b}\rangle^2}$. 
			We can get desired smallness directly if any derivative hits $F_{\b}$, namely,
			\begin{align*}
				&\|\nabla^n F_{\b}|_{\b=\nabla\phi(x)}\|_{L^\infty}\lesssim(1+\|\nabla f_2\|_{L^\infty})^2(\|\nabla f_2\|_{\dot C^1}^n+\|\nabla f_2\|_{\dot C^n})\lesssim   |\log t|^{-\varkappa}t^{-\frac{n}{3}} \|f_2\|_{T,*}(1+\|f_2\|_T)^{n},\\	&\|\delta_\beta \nabla^nF_{\b}|_{\b=\nabla\phi(x)}\|_{L^\infty}\lesssim |\beta|^a(1+\|\nabla f_2\|_{L^\infty})^2(\|\nabla f_2\|_{\dot C^{a}}^{\frac{n+a}{a}}+\|\nabla f_2\|_{\dot C^{n+a}})\\
				&\quad\quad\quad\quad\quad\quad\quad\quad\lesssim |\beta|^a  |\log t|^{-\varkappa}t^{-\frac{n+a}{3}} \|f_2\|_{T,*}(1+\|f_2\|_T)^{n}.
			\end{align*}
			Note that the above estimate does not rely on the exact value of $\b$, it holds for any $\b\in\mathbb{R}^d$. On the other hand, if no derivative hits $F_{\b}$, then $\b=\nabla\phi(x)$ is necessary to get smallness. Indeed, one has
			\begin{align*}
				F_{\b}(x)|_{\b=\nabla\phi(x)}=\langle \nabla f_2(x)\rangle^2\left(\frac{\nabla f_2(x)}{\langle\nabla f_2(x)\rangle^2}-\frac{\nabla \phi(x)}{\langle\nabla \phi(x)\rangle^2}\right),
			\end{align*}
			which yields 
			\begin{align*}
				\|	F_{\b}|_{\b=\nabla\phi(x)}\|_{L^\infty}\lesssim (1+\|\nabla f_2\|_{L^\infty})^2\|\nabla (f_2-\phi)\|_{L^\infty}.
			\end{align*}
			Moreover, apply \eqref{mainint21} and \eqref{mainint22} in Lemma \ref{lemcom}, we get 
			\begin{align*}
				&	\|\nabla(\mathcal{G}[f_1]-\mathcal{G}[f_2])\|_{\dot C^n}\\
				&	\lesssim\|\nabla(f_1-f_2)\|_{\dot C^{n+2}}+ \sum_{n_1=0}^{n+1}\|\nabla(f_1-f_2)\|_{\dot C^{n_1}}\left(\|(\nabla f_1,\nabla f_2)\|_{\dot C^1}^{n+2-n_1}+\|(\nabla f_1,\nabla f_2)\|_{\dot C^{n+2-n_1}}\right)\\
				&	\lesssim  |\log t|^{-\varkappa}t^{-\frac{n+2}{3}}\|f_1-f_2\|_T(1+\|(f_1,f_2)\|_{T})^{n+2},
			\end{align*}
			and 
			\begin{align*}
				&	\|\nabla(\mathcal{G}[f_1]-\mathcal{G}[f_2])\|_{\dot C^{n+a}}\\
				&\lesssim\|\nabla(f_1-f_2)\|_{\dot C^{n+2+a}}+ \sum_{n_1=0}^{n+1}\|\nabla(f_1-f_2)\|_{\dot C^{n_1+a}}\left(\|(\nabla g_1,\nabla g_2)\|_{\dot C^1}^{n+2-n_1}+\|(\nabla g_1,\nabla g_2)\|_{\dot C^{n+2-n_1}}\right)\\
				&\quad\quad+\sum_{n_1=0}^{n+2}\|\nabla(f_1-f_2)\|_{\dot C^{n_1}}\left(\|(\nabla f_1,\nabla f_2)\|_{\dot C^a}^{\frac{n+2-n_1+a}{a}}+\|(\nabla f_1,\nabla f_2)\|_{\dot C^{n+2-n_1+a}}\right)\\
				&\lesssim  |\log t|^{-\varkappa}t^{-\frac{n+2+a}{3}}\|f_1-f_2\|_T(1+\|(f_1,f_2)\|_{T})^{n+10}.
			\end{align*}
			Hence we get 
			\begin{align*}
				\|(\delta_\beta\nabla^n  \tilde {\mathcal{U}}_{2,\b})|_{\b=\nabla\phi(x)}\|_{L^\infty}\lesssim |\beta|^a |\log t|^{-\varkappa}t^{-\frac{n+2+a}{3}}\|f_1-f_2\|_T(\|f_2\|_{T,*}+\|f_2-\phi\|_T)(1+\|(f_1,f_2)\|_{T,*})^{n+10}.
			\end{align*}
			Combining this with \eqref{u1bhol}, we obtain 
			\begin{equation}\label{difw}
				\begin{aligned}
					&	\|(\delta_\beta\nabla^n   ({\mathcal{U}}_{\b}[f_1]-{\mathcal{U}}_{\b}[f_2]))|_{\b=\nabla\phi(x)}\|_{L^\infty}\\
					&\quad\quad\quad\quad\lesssim |\beta|^a |\log t|^{-\varkappa}t^{-\frac{n+2+a}{3}}\|f_1-f_2\|_T(\|f_1,f_2\|_{T,*}+\|f_2-\phi\|_T)(1+\|(f_1,f_2)\|_{T,*})^{n+10}.
				\end{aligned}
			\end{equation}
			By \eqref{pbhol}, \eqref{difP}, \eqref{difA}, \eqref{difw}, we obtain \eqref{difne}. This completes the proof of \eqref{dif}. \vspace{0.3cm}\\
			To get \eqref{m}, we note that 
			\begin{align*}
				&|||\nabla^n \H^{ex}_\b[f]|||_{a}\leq
				|||\nabla^n (\H^{ex}_\b[f]-\H^{ex}_\b[\phi])|||_{a}+|||\nabla^n \H^{ex}_\b[\phi]|||_{a}.
			\end{align*}
			Applying \eqref{dif} with $(f_1,f_2)=(f,\phi)$ to get 
			\begin{align*}
				\sup_{t\in[0,T]}|\log t|^\varkappa  t^\frac{n+\kappa}{3}	|||\nabla^n (\H^{ex}_\b[f]-
				\H^{ex}_\b[\phi])|||_{a}\lesssim \|f-\phi\|_T(\|f-\phi\|_T+T^\frac{1}{20}\|\nabla\phi\|_{C^{m+3}})(1+\|(f,\phi)\|_T)^{2m+5}.
			\end{align*}
			On the other hand,
			by the smoothness of $\phi$, one can check that 
			\begin{align*}
				|||\nabla^n \H^{ex}_\b[\phi]|||_{a}\lesssim (1+\|\nabla\phi\|_{C^{n+3}})^{n+3}. 
			\end{align*}
			Hence we obtain 
			\begin{align*}
				\sup_{t\in[0,T]}|\log t|^\varkappa  t^\frac{n+\kappa}{3}|||\nabla^n \H^{ex}_\b[f](t)|||_{a}\lesssim &  \|f-\phi\|_T(\|f-\phi\|_T+T^\frac{1}{20}\|\nabla\phi\|_{C^{m+3}})(1+\|(f,\phi)\|_T)^{2m+5}\\
				&\quad\quad\quad\quad\quad+T^\frac{1}{5}(1+\|\nabla\phi\|_{C^{m+3}})^{m+3},\quad\quad\forall n\in\mathbb{N}, \ n\leq m.
			\end{align*}
			This completes the proof.
		\end{proof}
		
		The following lemma gives the estimate for lower order term $R[g]=(R_i[g])_{i=1}^d$, where $R_i[g]$ was defined in \eqref{defR}.
		\begin{lemma}\label{lemR}
			Let $\|\cdot\|_{T}$ be as defined in \eqref{defnorgm} with constants $\varkappa,m,\kappa$ fixed in \eqref{consgm}. Then for any $T\in(0,\frac{1}{2})$, there holds
			\begin{align*}
				&\sup_{t\in[0,T]}\left(\int_0^t|\log (t-\tau)|^\varkappa\|R[g](\tau)\|_{L^\infty}d\tau+|\log t|^\varkappa  t^{\frac{m}{3}+1}\|R[g](t)\|_{\dot C^{m}}\right)		\\
				&\hspace{7cm}\quad\quad	\lesssim T^\frac{1}{6}(1+\|\nabla \phi\|_{C^{m+4}})^{m+3}(\|g-\phi\|_{T}+1),\\
				&\sup_{t\in[0,T]}\left(\int_0^t|\log (t-\tau)|^\varkappa{\|(R[g_1]-R[g_2])(\tau)\|_{L^\infty}}d\tau+|\log t|^\varkappa  t^{\frac{m}{3}+1}\|(R[g_1]-R[g_2])(t)\|_{\dot C^{m}}\right)\\
				&\hspace{7cm}\quad\quad	\lesssim T^\frac{1}{6}(1+\|\nabla \phi\|_{C^{m+4}})^{m+3}\|g_1-g_2\|_{T}.
			\end{align*}
		\end{lemma}
		\begin{proof}
			By symmetry,	we can write 
			\begin{align*}
				R[g](x)=&	\frac{1}{2}\int_{\mathbb{R}^d} \nabla_x\left( \A(\nabla \phi(x),\alpha)\right):\mathcal{O}_\alpha \nabla^2 (g-\phi)(x)\frac{d\alpha}{|\alpha|^{d}}+\frac{1}{2}\nabla\left(\int_{\mathbb{R}^d}\A(\nabla \phi(x),\alpha):\mathcal{O}_\alpha \nabla ^2\phi(x)\frac{d\alpha}{|\alpha|^{d}}\right).
			\end{align*}
			Then by Lemma \ref{lemoalin}, it is easy to check that 
			\begin{align*}
				\|R[g](t)\|_{L^\infty}&\lesssim (1+\|\nabla \phi\|_{C^4})^3t^{-\frac{2}{3}}(\|g-\phi\|_{T}+1),\\
				\|R[g](t)\|_{\dot C^m}&\lesssim (1+\|\nabla \phi\|_{C^{m+4}})^{m+3}t^{-\frac{m+2}{3}}(\|g-\phi\|_{T}+1).
			\end{align*}
			Hence 
			\begin{align*}
				\sup_{t\in[0,T]}\left(\int_0^t|\log (t-\tau)|^\varkappa\|R[g](\tau)\|_{L^\infty}d\tau+|\log t|^\varkappa  t^{\frac{m}{3}+1}\|R[g](t)\|_{\dot C^{m}}\right)&\\
				\quad\quad	\quad\quad\quad\quad\quad\quad\quad	\lesssim T^\frac{1}{6}(1+\|\nabla \phi\|_{C^{m+4}})^{m+3}&(\|g-\phi\|_{T}+1) .
			\end{align*}
			Then we obtain the first estimate. The second estimate follows by similar arguments. The proof is completed. 
		\end{proof}
		~\vspace{0.3cm}\\
		It remains to estimate the implicit force term $\H^{im}_\b[f]$ defined in \eqref{defHim}. This is the most difficult part because it involves
		the estimate of $\omega_\b^\pm$, which is a solution to the elliptic system \eqref{elliw} and for which we do not have an explicit formula. 
		Now we start to estimate the H\"{o}lder norm of solution to system \eqref{elliw}. 
		Recalling \eqref{defwb}, we have 
		\begin{align*}
			v_f^\pm=\omega_\b^\pm[f]+\tilde v_\b^\pm[f].
		\end{align*}
		The estimation of force terms in \eqref{elliw} requires estimation of $\tilde v_\b^\pm[f]$ which is defined in 
		\eqref{deftvb}:
		\begin{align}\label{fortvb}
			\tilde v_\b^\pm[f](x,z)=c^\pm\int_{\mathbb{R}^d}K_\b(x-y,z)\mathcal{G}[f](y)dy.
		\end{align}
		As preparation, we first give some estimates of $\tilde v_\b^\pm[f]$.
		
		The following lemma gives an explicit formula of the kernel $K_\b$ defined in \eqref{defKor}.
		\begin{lemma}\label{Z10} Let $\b\in \mathbb{R}^d$ and $	\lambda^\pm(\xi,\b)=\frac{i\b\cdot\xi\pm\sqrt{(1+|\b|^2)|\xi|^2-(\b\cdot \xi)^2}}{1+|\b|^2},$
			\begin{equation*}
				K_\b(x,z)=\begin{cases}
					\frac{1}{(2\pi)^\frac{d}{2}}\int_{\mathbb{R}^d}\exp\left(ix\cdot\xi+\lambda^-(\xi,\b)z\right)d\xi,& z\geq0,\\
					\frac{1}{(2\pi)^\frac{d}{2}}\int_{\mathbb{R}^d}\exp\left(ix\cdot\xi+\lambda^+(\xi,\b)z\right)d\xi,& z<0.
				\end{cases}	
			\end{equation*}
			Then
			\begin{equation}\label{defK}
				K_\b(x,z)=c_d\frac{|z|}{\left((x\cdot \b+z)^2+|x|^2\right)^{\frac{d+1}{2}}},\ \ \ \forall x\in\mathbb{R}^d,\ \ z\neq 0.
			\end{equation}
			where the constant $c_d=\frac{\Gamma(\frac{d+1}{2})}{\pi^\frac{d+1}{2}}$ ensures that $\int _{\mathbb{R}^d}K_\b(x,1)dx=1$. Moreover, when $z=0$, there holds $K_\b(x,0)=\delta(x)$.
		\end{lemma}
		\begin{proof}
			When $z=0$, the result is well-known. 
			For simplicity, we only consider $z<0$. And the case $z> 0$ can be done similarly. We first  recall the explicit formula of Poisson kernel (see \cite[Chapter \rm{III}, section 2.1]{Steinbook})
			\begin{equation}\label{Z1}
				\frac{1}{(2\pi)^\frac{d}{2}}\int_{\mathbb{R}^d}\exp\left(ix\cdot\xi+|\xi|z\right)d\xi=c_d\frac{-z}{(|x|^2+|z|^2)^{\frac{d+1}{2}}},~z<0.
			\end{equation}
			1) Case $d=1$, we have, 
			\begin{align*}
				K_\b(x,z)&=\frac{1}{2\pi}\int_{\mathbb{R}}\exp\left(ix\cdot\xi+\frac{i\b\xi+|\xi|}{\langle \b\rangle^2}z\right)d\xi\\&= \frac{1}{2\pi}\int_{\mathbb{R}}\exp\left(i\left(x+\frac{\b z}{\langle \b\rangle^2}\right)\cdot\xi+|\xi|\frac{z}{\langle \b\rangle^2}\right)d\xi.
			\end{align*}
			Applying \eqref{Z1} with $d=1$ to get 
			\begin{align*}
				K_\b(x,z)&=c_1\frac{-\frac{z}{\langle \b\rangle^2}}{\left(x+\frac{\b z}{\langle \b\rangle^2}\right)^2+\left(\frac{z}{\langle \b\rangle^2}\right)^2}=c_1\frac{-z }{\left( x\cdot \b+z\right)^2+|x|^2}.
			\end{align*}
			2) Case $d\geq 2$. First we have, for any $\theta_1,\theta_2\in \mathbb{S}^{d-1}$, and any function $G$,
			\begin{equation}\label{Z2}
				\int_{\mathbb{R}^d} G(\theta_1\cdot x,\theta_2\cdot x,|x|)dx=\int_{\mathbb{R}^d} G((1-(\theta_1\cdot\theta_2)^2)^{\frac{1}{2}}x_1+(\theta_1\cdot\theta_2)x_2,x_2,|x|)dx.
			\end{equation} 
			Indeed, the unit vectors 
			\begin{equation*}
				\theta_2, \tilde{\theta}_1= \frac{\theta_1-(\theta_1\cdot\theta_2)\theta_2}{(1-(\theta_1\cdot\theta_2)^2)^{\frac{1}{2}}}\in \mathbb{S}^{d-1},
			\end{equation*}
			are orthogonal. 
			Hence,  by rotation, 
			\begin{align*}
				\int_{\mathbb{R}^d} G(\theta_1\cdot x,\theta_2\cdot x,|x|)dx&=\int_{\mathbb{R}^d} G((\tilde{\theta}_1\cdot\theta_1)x_1+(\theta_1\cdot\theta_2)x_2,x_2,|x|)dx\\&=\int_{\mathbb{R}^d} G((1-(\theta_1\cdot\theta_2)^2)^{\frac{1}{2}}x_1+(\theta_1\cdot\theta_2)x_2,x_2,|x|)dx.
			\end{align*}
			Therefore, 
			\begin{align*}
				&	K_\b(x,z)=\frac{1}{(2\pi)^d}\int_{\mathbb{R}^d}\exp\left(ix\cdot\xi+\frac{i\b\cdot\xi+\sqrt{(1+|\b|^2)|\xi|^2-(
						\b\cdot \xi)^2}}{1+|
					\b|^2}z\right)d\xi\\&=\frac{1}{(2\pi)^d}\int_{\mathbb{R}^d}\exp\left(i|x|\left((1-(\frac{x\cdot \b}{|x||
					\b|})^2)^{\frac{1}{2}}\xi_1+\frac{x\cdot \b}{|x||\b|}\xi_2\right)+\frac{i|\b|\xi_2+\sqrt{(1+|\b|^2)|\xi|^2-|\b|^2\xi_2^2}}{1+|\b|^2}z\right)d\xi
				\\&=\frac{1}{(2\pi)^d}\int_{\mathbb{R}^d}\exp\left(i\frac{((|x||\b|)^2-(x\cdot \b)^2)^{\frac{1}{2}}}{|\b|}\xi_1+i\left(\frac{x\cdot \b}{|\b|}+\frac{|\b|z}{\langle \b\rangle^2}\right)\xi_2+\frac{\sqrt{(1+|\b|^2)|\xi|^2-|\b|^2\xi_2^2}}{1+|\b|^2}z\right)d\xi.
			\end{align*}
			By simple changes of variables, we obtain 
			\begin{align*}
				K_\b(x,z)=\frac{1}{(2\pi)^d}\int_{\mathbb{R}^d}\exp\left(i\frac{((|x||\b|)^2-(x\cdot \b)^2)^{\frac{1}{2}}}{|\b|\langle \b\rangle}\xi_1+i\left(\frac{x\cdot \b}{|\b|}+\frac{|\b|z}{\langle \b\rangle^2}\right)\xi_2+|\xi|\frac{z}{\langle \b\rangle^2}\right)\frac{d\xi}{\langle \b\rangle^{d-1}}.
			\end{align*}
			Finally, applying  \eqref{Z1} to get
			\begin{align*}
				K_\b(x,z)	&=\frac{c_d}{\langle \b\rangle^{d-1}}\frac{-\frac{z}{\langle \b\rangle^2}}{\left(\frac{(|x||\b|)^2-(x\cdot \b)^2}{|\b|^2\langle \b\rangle^2}+\left(\frac{x\cdot \b}{|\b|}+\frac{|\b|z}{\langle \b\rangle^2}\right)^2+\left(\frac{z}{\langle \b\rangle^2}\right)^2\right)^{\frac{d+1}{2}}}\\&=
				\frac{c_d}{\langle \b\rangle^{d-1}}\frac{-\frac{z}{\langle \b\rangle^2}}{\left(\frac{(x\cdot \b)^2+2zx\cdot \b +z^2}{\langle \b\rangle^2}+\frac{z^2}{\langle \b\rangle^2}\right)^{\frac{d+1}{2}}}
				=
				c_d\frac{-z}{\left((x\cdot \b+z)^2+|x|^2\right)^{\frac{d+1}{2}}}.
			\end{align*}
			This completes the proof.
		\end{proof}
		
		Based on the explicit formula of $K_\b$ in \eqref{defK}, we collect the following estimates for later use.
		\begin{lemma}\label{esK11}
			We have,\\
			(1) For any $n\in\mathbb{N}$,
			\begin{align}\label{Knxb}
				|\nabla_x^{n}K_\b(x,z)|\lesssim \frac{|z|\langle \b\rangle^{n}}{\left((x\cdot \b+z)^2+|x|^2\right)^{\frac{d+1+n}{2}}},\quad\forall x\in\mathbb{R}^d, z\neq 0.
			\end{align}
			(2)   For any $n_1,n_2\in\mathbb{N}$,
			\begin{align}\label{deriK}
				|\nabla_x^{n_1}\partial_z^{n_2}K_\b(x,z)|\lesssim {\langle \b\rangle^{1+n_1+n_2}}{\left(|x|^2+\frac{|z|^2}{\langle \b\rangle^2}\right)^{-\frac{d+n_1+n_2}{2}}},\quad \quad\forall x\in\mathbb{R}^d, z\neq 0.
			\end{align}
			(3) For any $z\neq 0$, $\b\in\mathbb{R}^d$,
			\begin{align}\label{remintk}
				\int_{\mathbb{R}^d}K_\b(x,z)dx=C,
			\end{align}
			where $C$ is a constant independent in $z$ and $\b$.\\
			(4)	For any $a\in(0,1)$,	
			\begin{equation}\label{esl1}
				\begin{aligned}
					&\int_{\mathbb{R}^d}|K_\b(y,z)-K_\b(y,\beta)||y|^ady \lesssim \langle \b\rangle^2|z-\beta|^a,\ \ \ \forall a\in(0,1).
				\end{aligned}
			\end{equation}
			
		\end{lemma}
		\begin{proof} The point-wise estimates \eqref{Knxb} and \eqref{deriK} follow directly from the formula \eqref{defK}.\\
			It remains to prove \eqref{remintk}. Observing the scaling invariant property
			$$
			K_\b(x,z)=|z|^{-d}K_\b(|z|^{-1}x,|z|^{-1}z),
			$$
			we obtain that for any $\b\in\mathbb{R}^d$, $z\neq 0$,
			\begin{align*}
				\int_{\mathbb{R}^d}K_\b(x,z)dx=|z|^{-d}\int_{\mathbb{R}^d} K_\b(|z|^{-1}x,|z|^{-1}z)dx=\int_{\mathbb{R}^d}K_\b(x,\pm 1)dx,\ \ \ z\in\mathbb{R}^\pm.
			\end{align*}
			By a rotation transformation and changes of variables, we obtain
			\begin{align*}
				\int_{\mathbb{R}^d}K_\b(x,\pm 1)dx&=c_d \int_{\mathbb{R}^d}\frac{1}{\left((x\cdot \b\pm 1)^2+|x|^2\right)^{\frac{d+1}{2}}}dx=\int_{\mathbb{R}^d}\frac{1}{\left((x_1|\b|\pm 1)^2+|x|^2\right)^{\frac{d+1}{2}}}dx\\
				&=\langle\b\rangle^{-1}\int_{\mathbb{R}^d}\frac{1}{\left(|x|^2+\frac{1}{\langle\b\rangle^{2}}\right)^\frac{d+1}{2}}dx=\int_{\mathbb{R}^d}\frac{1}{\langle|x|\rangle^{d+1}}dx=C.
			\end{align*}
			This completes the proof of \eqref{remintk}. Finally, for any $z,\beta>0$, $m\in\mathbb{N}$, by \eqref{deriK} we have
			\begin{align*}
				&|K_\b(y,z)-K_\b(y,\beta)|\lesssim \frac{\langle \b\rangle^{2}}{|
					y|^d} \min\left\{1,\frac{|z-\beta|}{|y|}\right\}.
			\end{align*}
			This implies \eqref{esl1}. Then we complete the proof.
		\end{proof}
		
		In view of the formula \eqref{fortvb}, when dealing with $\nabla_x^n \tilde v_\b^\pm[f]$, one can transfer derivative in $x$ freely between the kernel $K_\b$ and the function $\mathcal{G}[f]$. Hence Lemma \ref{esK11} is sufficient to estimate any order derivative of $\tilde v_\b^\pm[f]$ with respect to $x$. However, when estimating $\nabla_x^n \partial_z \tilde v_\b^\pm[f]$, we have to consider $\partial_z K_\b$ which appears to be singular near the interface $\{z=0\}$ because 
		\begin{align*}
			\int_{\mathbb{R}^d}|\partial_z K_\b(x,z)|dx\lesssim \langle\b\rangle^3|z|^{-1}.
		\end{align*}
		To solve this issue, we further analyze $\partial_z\tilde v_\b^\pm[f]$ by its Fourier transform 
		\begin{align}\label{pztilvb}
			\partial_z\tilde {v}_{\b}^{\pm}[f]=c^{\pm}\mathcal{F}^{-1}\left(\lambda^{\mp}(\xi,\b)e^{\lambda^{\mp}(\xi,\b)}\hat{\mathcal{G}}[f](\xi)\right).
		\end{align}
		More precisely, we give the following Lemma \ref{lemz0} and Lemma \ref{lemfdz0}. 
		\begin{lemma}\label{lemz0} Let $e\in \mathbb{S}^d$, $\b\in\mathbb{R}^d$. There holds, 
			\begin{align}\label{contc}
				\tilde c_d\int_{\mathbb{R}^d} \frac{1-\cos(e\cdot\alpha)}{\left((\alpha\cdot \b)^2+|\alpha|^2\right)^{\frac{d+1}{2}}}d\alpha=\frac{(\langle \b\rangle^2-(\b\cdot e)^2)^{\frac{1}{2}}}{\langle \b\rangle^2},
			\end{align}
			where $\tilde c_d$ is a positive constant depending only on the dimension $d$.
		\end{lemma}
		\begin{proof}
			Applying the rotation transform \eqref{Z2} to get
			\begin{align*}
				&\int_{\mathbb{R}^d} \frac{1-\cos(e\cdot \alpha)}{\left((\alpha\cdot \b)^2+|\alpha|^2\right)^{\frac{d+1}{2}}}d\alpha=\int_{\mathbb{R}^d} \frac{1-\cos(\alpha_1)}{\left(((\b\cdot e)\alpha_1+(|\b|^2-(\b\cdot e)^2)^{\frac{1}{2}}\alpha_2)^2+|\alpha|^2\right)^{\frac{d+1}{2}}}d\alpha\\&=\int_{\mathbb{R}^d} \frac{1-\cos(\alpha_1)}{\left((1+(\b\cdot e)^2)\alpha_1^2+(\langle \b\rangle^2-(\b\cdot e)^2)\alpha_2^2+2(\b\cdot e)(|\b|^2-(\b\cdot e)^2)^{\frac{1}{2}}\alpha_1\alpha_2+|\alpha'|^2\right)^{\frac{d+1}{2}}}d\alpha,
			\end{align*}
			where we denote $\alpha'=\alpha-\alpha_1e_1-\alpha_2e_2$.
			Since, 
			\begin{align*}
				&\left(1+(\b\cdot e)^2\right)\alpha_1^2+\left(\langle \b\rangle^2-(\b\cdot e)^2\right)\alpha_2^2+2(\b\cdot e)\left(|\b|^2-(\b\cdot e)^2\right)^{\frac{1}{2}}\alpha_1\alpha_2\\&=\frac{\langle \b\rangle^2}{\langle \b\rangle^2-(\b\cdot  e)^2}\alpha_1^2+\left((\langle \b\rangle^2-(\b\cdot e)^2)^{\frac{1}{2}}\alpha_2+\frac{(\b\cdot e)(|\b|^2-(\b\cdot e)^2)^{\frac{1}{2}}}{(\langle \b\rangle^2-(\b\cdot e)^2)^{\frac{1}{2}}}\alpha_1\right)^2.
			\end{align*}
			Using elementary linear transformations, we get 
			\begin{align*}
				\int_{\mathbb{R}^d} \frac{1-\cos(e\cdot \alpha)}{\left((\alpha\cdot \b)^2+|\alpha|^2\right)^{\frac{d+1}{2}}}d\alpha	&=\frac{1}{(\langle \b\rangle^2-(\b\cdot e)^2)^{\frac{1}{2}}}\int_{\mathbb{R}^d} \frac{1-\cos(\alpha_1)}{\left(\frac{\langle \b\rangle^2}{\langle \b\rangle^2-(\b\cdot e)^2}\alpha_1^2+\alpha_2^2+|\alpha'|^2\right)^{\frac{d+1}{2}}}d\alpha\\&
				=\frac{(\langle \b\rangle^2-(\b\cdot e)^2)^{\frac{1}{2}}}{\langle \b\rangle^2} \int_{\mathbb{R}^d} \frac{1-\cos(\alpha_1)}{|\alpha|^{d+1}}d\alpha=\tilde c_d^{-1}\frac{(\langle \b\rangle^2-(\b\cdot e)^2)^{\frac{1}{2}}}{\langle \b\rangle^2},
			\end{align*}
			where $\tilde c_d=\left(\int_{\mathbb{R}^d} \frac{1-\cos(\alpha_1)}{|\alpha|^{d+1}}d\alpha\right)^{-1}$.
			This completes the proof.
		\end{proof}
		\begin{lemma}\label{lemfdz0}
			Let $\b\in \mathbb{R}^d$ and $	\lambda^\pm(\xi,\b)=\frac{i\b\cdot\xi\pm\sqrt{(1+|\b|^2)|\xi|^2-(\b\cdot \xi)^2}}{1+|\b|^2}$. Then
			\begin{align}\label{opL}
				\mathcal{F}^{-1}\left(\lambda^\pm(\xi,\b)\hat f(\xi)\right)(x)=\frac{\b\cdot\nabla f(x)}{\langle \b\rangle^2}\pm\tilde c_d \mathrm{P.V.}\int_{\mathbb{R}^d} \frac{\delta_\alpha f(x)}{\langle \hat\alpha\cdot \b\rangle^{d+1}}\frac{d\alpha}{|\alpha|^{d+1}}:=\mathfrak{L}_{\pm,\b} f(x),
			\end{align}
			where $\tilde c_d $ is the positive constant in \eqref{contc}.
			Moreover, for any $a\in (0,1)$, $m\in\mathbb{N}$,
			\begin{align}\label{LHol}
				\sup_{\b\in\mathbb{R}^d}	\|\mathfrak{L}_{\pm,\b} f\|_{\dot C^{m+a}_x}\lesssim \|f\|_{\dot C^{m+1+a}},
			\end{align}
			and \begin{align}\label{LLLinf}
				\sup_{\b\in\mathbb{R}^d}	\|\mathfrak{L}_{\pm,\b} f\|_{L^\infty}\lesssim \|f\|_{ C^{1,\log^\varkappa  }},
			\end{align}
			where $\varkappa>1$.
		\end{lemma}
		\begin{proof}
			We have
			\begin{equation}\label{fopz}
				\begin{aligned}
					&	\mathcal{F}^{-1}\left(\lambda^\pm(\xi,\b)\hat f(\xi)\right)=\frac{1}{(2\pi)^\frac{d}{2}}\int_{\mathbb{R}^d}\lambda^+(\xi,\b)\hat f (\xi)\exp\left(ix\cdot\xi\right)d\xi \\
					&\ \ \ \ \ =\frac{\b\cdot\nabla f(x)}{\langle \b\rangle^2}+\frac{1}{(2\pi)^\frac{d}{2}}\int_{\mathbb{R}^d}\int_{\mathbb{R}^d}\frac{\sqrt{\langle \b\rangle^2|\xi|^2-(\b\cdot\xi)^2}}{\langle \b\rangle^2}\hat{f}(\xi)\exp\left(ix\cdot \xi\right)d\xi .
				\end{aligned}
			\end{equation}
			By Lemma \ref{lemz0}, we have 
			$$
			\frac{\sqrt{\langle \b\rangle^2|\xi|^2-(\b\cdot\xi)^2}}{\langle \b\rangle^2}=	\tilde c_d|\xi|\int_{\mathbb{R}^d}\frac{1-\exp{(i\hat \xi\cdot \alpha)}}{\langle\hat \alpha \cdot \b\rangle^{d+1}}\frac{d\alpha}{|\alpha|^{d+1}},
			$$
			where $\hat \xi=\frac{\xi}{|\xi|}$.
			This implies that 
			\begin{align*}
				\frac{1}{(2\pi)^\frac{d}{2}}\int_{\mathbb{R}^d}\frac{\sqrt{\langle \b\rangle^2|\xi|^2-(\b\cdot\xi)^2}}{\langle \b\rangle^2}\hat f(\xi)\exp(ix\cdot\xi)d\xi&=\tilde c_d\iint_{\mathbb{R}^{2d}}\frac{(1-\exp(i\xi\cdot\alpha))\hat{f}(\xi)}{\langle\hat{\alpha}\cdot b\rangle^{d+1}}\exp(ix\cdot\xi)d\xi\frac{d\alpha}{|\alpha|^{d+1}} \\
				&=\tilde c_d\int_{\mathbb{R}^d} \frac{\delta_\alpha f(x)}{\langle \hat\alpha\cdot \b\rangle^{d+1}}\frac{d\alpha}{|\alpha|^{d+1}}.
			\end{align*}
			Combining this with \eqref{fopz}, we get \eqref{opL}.\\
			Moreover, observe that 
			\begin{align*}
				\nabla_x^m\mathfrak{L}_{\pm,\b} f(x)= \mathfrak{L}_{\pm,\b}(\nabla_x^mf)(x).
			\end{align*}
			Hence, to prove \eqref{LHol} and \eqref{LLLinf}, it suffices to prove that 	
			\begin{align}\label{gggg1}
				\sup_{\b\in\mathbb{R}^d}\left\|\int_{\mathbb{R}^d} \frac{\delta_\alpha f(x)}{\langle \hat\alpha\cdot \b\rangle^{d+1}}\frac{d\alpha}{|\alpha|^{d+1}}\right\|_{\dot C^a_x}\lesssim \|f\|_{\dot C^{1+a}},\quad\quad\quad 	\sup_{\b\in\mathbb{R}^d}\left\|\int_{\mathbb{R}^d} \frac{\delta_\alpha f(x)}{\langle \hat\alpha\cdot \b\rangle^{d+1}}\frac{d\alpha}{|\alpha|^{d+1}}\right\|_{L^\infty_x}\lesssim \|f\|_{ C^{1,\log^\varkappa}}.
			\end{align}
			Indeed, by symmetry we have 
			\begin{align*}
				\int_{\mathbb{R}^d} \frac{\delta_\alpha f(x)}{\langle \hat\alpha\cdot \b\rangle^{d+1}}\frac{d\alpha}{|\alpha|^{d+1}}=\frac{1}{2}\int_{\mathbb{R}^d} \frac{\mathcal{O}_\alpha f(x)}{\langle \hat\alpha\cdot \b\rangle^{d+1}}\frac{d\alpha}{|\alpha|^{d}},
			\end{align*}
			where $\mathcal{O}_\alpha$ is the double finite difference operator.
			Then \eqref{gggg1} follows from Lemma \ref{lemoalin}.
			This completes the proof.
		\end{proof}
		\begin{lemma}\label{lemKf}
			Consider $K_\b(x,z)$ defined in \eqref{defK}.	For any fixed $\b\in\mathbb{R}^d$, $m_1\in\{0,1\}$, $a\in(0,1)$, there hold\footnote{{Note that $\|\cdot \|_{\dot C^a_{x,z}(\mathbb{R}^{d+1}_+\cup \mathbb{R}^{d+1}_-)}=\|\cdot \|_{\dot C^a_{x,z}(\mathbb{R}^{d+1}_+)}+\|\cdot \|_{\dot C^a_{x,z}(\mathbb{R}^{d+1}_-)}$ is different from  $\|\cdot\|_{\dot C^a_{x,z}(\mathbb{R}^{d+1})}$.}}
			\begin{align}
				&\left\|(\partial_z^{m_1}K_\b(\cdot,z)\ast f)(x)\right\|_{\dot C^a_{x,z}(\mathbb{R}^{d+1}_+\cup \mathbb{R}^{d+1}_-)}\lesssim \langle \b\rangle^2 \|f\|_{\dot C^{m_1+a}},\label{hhhkf}\\
				&\left\|(\partial_z^{m_1}K_\b(\cdot,z)\ast f)(x)\right\|_{L^\infty_{x,z}(\mathbb{R}^{d+1}_+\cup \mathbb{R}^{d+1}_-)}\lesssim  \|f\|_{ C^{m_1,\log^\varkappa  }},\label{lllkf}\ \ \varkappa>1.
			\end{align}
		\end{lemma}
		\begin{proof}
			Note that 
			\begin{align*}
				\mathcal{F}(\partial_z^{m_1}K_\b(\cdot,z))=\left(\lambda^\pm(\xi,\b)\right)^{m_1}\exp(\lambda^\pm(\xi,\b)z),\ \ \ z\in\mathbb{R}^\mp.
			\end{align*}
			Hence we can write 
			\begin{align*}
				\partial_z^{m_1}(K_\b(\cdot,z)\ast f)(x)=(K_\b(\cdot,z)\ast( \mathfrak{L}^{m_1}_{\pm,\b} f)(\cdot))(x),
			\end{align*}
			where the operator $\mathfrak{L}_{\pm,\b}$ is defined in \eqref{opL}.
			By \eqref{LLLinf} and \eqref{remintk}, we obtain \eqref{lllkf}.\\
			By \eqref{LHol}, we obtain 
			\begin{align}\label{holx}
				\|	\partial_z^{m_1}(K_\b(\cdot,z)\ast f)(x)	\|_{L^\infty_z(\mathbb{R})\dot C^a_{x}(\mathbb{R}^d)}\lesssim\|\mathfrak{L}^{m_1}_{\pm,\b} f\|_{\dot C^a}\lesssim \|f
				\|_{\dot C^{m_1+a}}.
			\end{align}
			Then we estimate the  H\"{o}lder norm in $z$-variable. For simplicity, we only consider $z\in\mathbb{R}^+$. Then for any  $\beta\in \mathbb{R}^+$,
			\begin{align*}
				&\partial_z^{m_1}(K_\b(\cdot,z)\ast f)(x)-	\partial_z^{m_1}(K_\b(\cdot,\beta)\ast f)(x)\\
				&\quad\quad=\int_{\mathbb{R}^{d}}(K_\b(y,z)-K_\b(y,\beta))( \mathfrak{L}_{-,\b}^{m_1}f)(x-y)dy\\
				&\quad\quad=\int_{\mathbb{R}^{d}}(K_\b(y,z)-K_\b(y,\beta))\left(( \mathfrak{L}_{-,\b}^{m_1}f)(x-y)-( \mathfrak{L}_{-,\b}^{m_1}f)(x)\right)dy.
			\end{align*}
			Here in the last inequality we applied \eqref{remintk} in Lemma \ref{esK11} to get
			\begin{align*}
				\int_{\mathbb{R}^{d}}\big(K_\b(y,z)-K_\b(y,\beta)\big)dy	=0.
			\end{align*}
			Then, by \eqref{esl1} and \eqref{LHol}, we obtain that 
			\begin{align*}
				&|((\partial_z^{m_1}K_\b)(\cdot,z)\ast f)(x)-	((\partial_z^{m_1}K_\b)(\cdot,\beta)\ast f)(x)|\\
				&\quad\quad\lesssim \int_{\mathbb{R}^{d}}|K_\b(y,z)-K_\b(y,\beta)||y|^ady\| \mathfrak{L}_{\pm,\b}^{m_1}f\|_{\dot C^a_x}\\
				&\quad\quad\lesssim \langle \b\rangle^2|z-\beta|^a \|f
				\|_{\dot C^{m_1+a}},
			\end{align*}
			which implies that 
			\begin{align*}
				\|\partial_z^{m_1}(K_\b(\cdot,z)\ast f)(x)\|_{L^\infty_x(\mathbb{R}^d)\dot C^a_z(\mathbb{R}^+)}\lesssim \langle \b\rangle^2\|f
				\|_{\dot C^{m_1+a}}.
			\end{align*}
			Combining this with \eqref{holx}, we get 
			\begin{align*}
				\|\partial_z^{m_1}(K_\b(\cdot,z)\ast f)(x)\|_{\dot C^a(\mathbb{R}^{d+1}_+)}\lesssim \langle \b\rangle^2\|f
				\|_{\dot C^{m_1+a}}.
			\end{align*}
			The estimate in $\mathbb{R}^{d+1}_-$ follows by the same argument.		This completes the proof.
		\end{proof}
		
		As a direct result of Lemma \ref{lemKf},		we have the following estimates of $\tilde v_\b^\pm[f]$.
		\begin{lemma}\label{lemtvbhol}
			For $\tilde v_\b^\pm[f]$ defined in \eqref{fortvb},  it holds
			\begin{align*}
				&\sup_{\b}\|\nabla_x^n\nabla_{x,z}\tilde v_\b^\pm[f]\|_{L^\infty(\mathbb{R}^{d+1}_\pm)}\lesssim \|\mathcal{G}[f]\|_{ C^{n+1,\log^\varkappa}},\ \ \ \varkappa>1, n\in\mathbb{N},\\
				&\sup_{\b}\|\nabla_x^n\nabla_{x,z}\tilde v_\b^\pm[f]\|_{\dot C^a(\mathbb{R}^{d+1}_\pm)}\lesssim (1+\|\nabla\phi\|_{L^\infty})^2\|\mathcal{G}[f]\|_{\dot C^{n+1+a}},\ \ \ a\in (0,1).
			\end{align*}
		\end{lemma}
		Here we recall that we denote $\sup_{\b\in\operatorname{Ran}(\nabla\phi)}$ by $\sup_\b$ 	 for short.\\
		Moreover, we have the following estimate in consideration of the convolution of $\nabla_{x,z}\tilde v_\b^\pm$   with some singular integral. 
		\begin{lemma}\label{lemsiG}
			Let $\mathbf{G}(x,z,y,w):\mathbb{R}^{d+1}_+\times \mathbb{R}^{d+1}_+\to \mathbb{R}$ with $(x,z),(y,w)\in\mathbb{R}^d\times \mathbb{R}^+$ be a singular kernel satisfying 
			\begin{align}\label{conbfG}
				\sup_{0<r_1<r_2<\infty}\left|\int_{\mathbb{R}^+}\int_{r_1\leq|x-y|\leq r_2}\mathbf{G}(x,z,y,w)dy dw\right|\leq c,\ \ \ \ \ \left|\int_{\mathbb{R}^+}\mathbf{G}(x,z,y,w)dw\right|\leq C|x-y|^{-d},
			\end{align}
			and for any $y_1,y_2\in\mathbb{R}^{d}$,
			\begin{align}\label{condiffy}
				\left|\int_{\mathbb{R}^+}\mathbf{G}(x,z,y_1,w)-\mathbf{G}(x,z,y_2,w)dw\right|\lesssim |y_1-y_2|\left(\frac{1}{|x-y_1|^{d+1}}+\frac{1}{|x-y_2|^{d+1}}\right),\ \ \forall (x,z)\in\mathbb{R}^{d+1}_+.
			\end{align}
			Then for 
			\begin{align*}
				\mathbf{V}_\b(x,z)=\int_{\mathbb{R}^{d+1}_+}\mathbf{G}(x,z,y,w)\nabla_{x,z}\tilde v_\b^+(y,w)dydw,
			\end{align*}
			there holds 
			\begin{align*}
				\sup_{\b}\|\nabla_x^n\mathbf{V}_\b^\pm\|_{L^\infty}\lesssim \|\mathcal{G}[f]\|_{\dot C^{n+1,\log^\varkappa}}+\|\mathcal{G}[f]\|_{ L^2},\ \ \ \varkappa>1, n\in\mathbb{N}.
			\end{align*}
		\end{lemma}
		\begin{proof}
			By \eqref{fortvb}, \eqref{pztilvb} and \eqref{opL}, we can write
			\begin{align*}
				\nabla_x \tilde v_\b^+(x,z)=c^+\left(K_\b(\cdot,z)\ast \nabla_x\mathcal{G}\right)(x),\ \ \ 
				\partial_z \tilde v_\b^+(x,z)=c^+\left(K_\b(\cdot,z)\ast \mathfrak{L}_{+,\b}\mathcal{G}\right)(x),
			\end{align*}
			where we write $\mathcal{G}=\mathcal{G}[f]$, which is defined in \eqref{ndefG}, for short. Then it is clear that 
			\begin{align*}
				\| \mathbf{V}_\b\|_{L^\infty(\mathbb{R}^{d+1}_+)}\lesssim \| \mathbf{V}_\b^1\|_{L^\infty(\mathbb{R}^{d+1}_+)}+\| \mathbf{V}_\b^2\|_{L^\infty(\mathbb{R}^{d+1}_+)},
			\end{align*}
			where 
			\begin{align*}
				&\mathbf{V}_\b^1(x,z)=\int_{\mathbb{R}^{d+1}_+}\mathbf{G}(x,z,y,w)\left(K_\b(\cdot,w)\ast \nabla\mathcal{G}\right)(y)dydw,\\
				&\mathbf{V}_\b^2(x,z)=\int_{\mathbb{R}^{d+1}_+}\mathbf{G}(x,z,y,w)\left(K_\b(\cdot,w)\ast \tilde {\mathfrak{L}}\mathcal{G}\right)(y)dydw,\ \ \text{with}\ \tilde{\mathfrak{L}}\mathcal{G}(x)=\int_{\mathbb{R}^d}\frac{\delta_\alpha \mathcal{G}(x)}{\langle\hat \alpha\cdot \b\rangle^{d+1}}\frac{d\alpha}{|\alpha|^{d+1}}.
			\end{align*}
			Using the fact that $\|K_\b\|_{L^1_xL^\infty_z}\lesssim 1$, we have 
			\begin{align*}
				|\mathbf{V}_\b^1(x,z)|\lesssim \sup_y\left|\int_{\mathbb{R}^{d+1}_+}\mathbf{G}(x,z,y',w)\nabla \mathcal{G}(y'-y)dy'dw\right|,\\
				|\mathbf{V}_\b^2(x,z)|\lesssim \sup_y\left|\int_{\mathbb{R}^{d+1}_+}\mathbf{G}(x,z,y',w)\tilde{\mathfrak{L}} \mathcal{G}(y'-y)dy'dw\right|.
			\end{align*}
			For simplicity, we assume the supremum is taken at $y = 0$. And denote $$
			\bar {\mathbf{G}}(x,z,y')=\int_0^\infty{\mathbf{G}}(x,z,y',w)dw.
			$$ By \eqref{conbfG},
			\begin{equation}\label{V1blinf}
				\begin{aligned}
					|\mathbf{V}_\b^1(x,z)|&\lesssim\left|\int_{\mathbb{R}^d}\bar {\mathbf{G}}(x,z,y')\nabla \mathcal{G}(y')dy'\right|\\
					&\lesssim\left|\int_{|x-y'|\leq 1}\bar {\mathbf{G}}(x,z,y')(\nabla \mathcal{G}(y')-\nabla \mathcal{G}(x))dy'\right|+c\|\nabla \mathcal{G}\|_{L^\infty}\\
					&\quad\quad+\left|\int_{|x-y'|\geq 1}\bar {\mathbf{G}}(x,z,y')\nabla \mathcal{G}(y')dy'\right|\\
					&\lesssim \|\nabla \mathcal{G}\|_{\dot C^{1,\log^\varkappa}}+c\|\nabla \mathcal{G}\|_{L^\infty}+\|\nabla \mathcal{G}\|_{L^2}\lesssim \|\nabla \mathcal{G}\|_{\dot C^{1,\log^\varkappa}}+\|\nabla \mathcal{G}\|_{L^2}.
				\end{aligned}
			\end{equation}
			For $\mathbf{V}_\b^2$, we further decompose 
			\begin{align*}
				&\int_{\mathbb{R}^{d+1}_+}\mathbf{G}(x,z,y',w)\tilde{\mathfrak{L}} \mathcal{G}(y')dy'dw=\int_{\mathbb{R}^{d}}\int_{\mathbb{R}^{d}}\frac{\bar {\mathbf{G}}(x,z,y')\delta_\alpha\mathcal{G}(y')}{\langle\hat \alpha\cdot \b\rangle^{d+1}}\frac{d\alpha dy'}{|\alpha|^{d+1}}\\
				&=\int_{\mathbb{R}^d}\int_{ |\alpha|\leq|x-y'|/2}\frac{\bar {\mathbf{G}}(x,z,y')\delta_\alpha\mathcal{G}(y')}{\langle\hat \alpha\cdot \b\rangle^{d+1}}\frac{d\alpha dy'}{|\alpha|^{d+1}}+\int_{\mathbb{R}^d}\int_{ |\alpha|\geq|x-y'|/2}\frac{\bar {\mathbf{G}}(x,z,y')\delta_\alpha\mathcal{G}(y')}{\langle\hat \alpha\cdot \b\rangle^{d+1}}\frac{d\alpha dy'}{|\alpha|^{d+1}}\\
				&=:I_1(x,z)+I_2(x,z).
			\end{align*}
			We first deal with $I_1$. We first write 
			\begin{align*}
				\delta_\alpha \mathcal{G}(y')=\int_0^1\alpha\cdot \nabla \mathcal{G}(y'-\alpha \tau)d\tau.
			\end{align*}
			By a change of variable, one has 
			\begin{align*}
				I_1(x,z)&=\int_{\mathbb{R}^{d}}\int_{ |\alpha|\leq|x-y'|/2}\int_0^1\frac{\bar {\mathbf{G}}(x,z,y')\alpha\cdot (\nabla \mathcal{G}(y'-\alpha \tau)-\nabla \mathcal{G}(x))}{\langle\hat \alpha\cdot \b\rangle^{d+1}}\frac{d\tau d\alpha dy'}{|\alpha|^{d+1}}\\
				&=\int_{\mathbb{R}^{d}}\int_{ |\alpha|\leq|x-y'|/2}\int_0^1\frac{\bar {\mathbf{G}}(x,z,y'+\alpha\tau)\alpha\cdot (\nabla \mathcal{G}(y')-\nabla \mathcal{G}(x))}{\langle\hat \alpha\cdot \b\rangle^{d+1}}\frac{d\tau d\alpha dy'}{|\alpha|^{d+1}}\\
				&=\int_{\mathbb{R}^{d}}\int_{ |\alpha|\leq|x-y'|/2}\left(\int_0^1(\bar {\mathbf{G}}(x,z,y'+\alpha\tau)-\bar {\mathbf{G}}(x,z,y'-\alpha\tau))d\tau\right)\frac{\alpha\cdot (\nabla \mathcal{G}(y')-\nabla \mathcal{G}(x))}{\langle\hat \alpha\cdot \b\rangle^{d+1}}\frac{ d\alpha dy'}{|\alpha|^{d+1}}
			\end{align*}
			By the condition \eqref{condiffy}, we have 
			\begin{align*}
				\left|\bar{\mathbf{G}}(x,z,y'+\alpha\tau)-\bar{\mathbf{G}}(x,z,y'-\alpha\tau)\right|\lesssim \frac{|\alpha|}{|x-y'|^{d+1}},\ \ \forall \tau\in(0,1),\ |\alpha|\leq|x-y'|/2.
			\end{align*}
			Hence,
			\begin{equation}\label{i1xz}
				\begin{aligned}
					|I_1(x,z)|&\lesssim\int_{\mathbb{R}^d} \frac{|\nabla \mathcal{G}(y')-\nabla \mathcal{G}(x))|}{|x-y'|^{d}}dy'\\
					&\lesssim 
					\int_{|x-y'|\leq 1} \frac{1}{|x-y'|^{d}\log^\varkappa(2+|x-y'|^{-1})} dy'\|\nabla\mathcal{G}\|_{\dot C^{\log^\varkappa}}+\left(\int_{|x-y'|\geq 1} \frac{1}{|x-y'|^{2d}}dy'\right)^\frac{1}{2}\|\nabla\mathcal{G}\|_{L^2}\\
					&\lesssim \|\nabla\mathcal{G}\|_{\dot C^{\log^\varkappa}}+\|\nabla\mathcal{G}\|_{L^2}.
				\end{aligned}
			\end{equation}
			Then we consider $I_2(x,z)$. We have 
			\begin{align*}
				I_2(x,z)&=\int_{\mathbb{R}^d}\int_{ |\alpha|\geq|x-y'|/2}\frac{\bar {\mathbf{G}}(x,z,y')(\delta_\alpha\mathcal{G}(y')-\delta_\alpha\mathcal{G}(x))}{\langle\hat \alpha\cdot \b\rangle^{d+1}}\frac{d\alpha dy'}{|\alpha|^{d+1}}+\int_{\mathbb{R}^d}\int_{ |\alpha|\geq|x-y'|/2}\frac{\bar {\mathbf{G}}(x,z,y')\delta_\alpha\mathcal{G}(x)}{\langle\hat \alpha\cdot \b\rangle^{d+1}}\frac{d\alpha dy'}{|\alpha|^{d+1}}\\
				&:=I_{2,1}(x,z)+I_{2,2}(x,z).
			\end{align*}
			Note that 
			\begin{align*}
				|\delta_\alpha\mathcal{G}(y')-\delta_\alpha\mathcal{G}(x)|\lesssim \min\{ |x-y'|\log^{-\varkappa}(2+|\alpha|^{-1})\|\mathcal{G}\|_{\dot C^{1,\log^\varkappa}},\|\mathcal{G}\|_{L^\infty}\}.
			\end{align*}
			This implies that 
			\begin{align*}
				|I_{2,1}(x,z)|\lesssim &\int_{|x-y'|\leq 1}\int_{|\alpha|\geq |x-y'|/2}\frac{1}{|x-y'|^{d-1}|\alpha|^{d+1}\log^{\varkappa}(2+|\alpha|^{-1})}d\alpha dy'\|\mathcal{G}\|_{\dot C^{1,\log^\varkappa}}\\
				&+\int_{|x-y'|\geq 1}\int_{|\alpha|\geq |x-y'|/2}\frac{1}{|x-y'|^d|\alpha|^{d+1}}d\alpha dy' \|\mathcal{G}\|_{L^\infty}\\
				\lesssim& \|\mathcal{G}\|_{ C^{1,\log^\varkappa}}.
			\end{align*}
			Moreover, by \eqref{conbfG}, following the proof of \eqref{gggg1} we obtain
			\begin{align*}
				|I_{2,2}|\lesssim c\sup_{r>0}\left|\int_{|\alpha|\geq r}\frac{\delta_\alpha \mathcal{G}(x)}{\langle\hat \alpha\cdot \b\rangle^{d+1}}\frac{d\alpha dy'}{|\alpha|^{d+1}}\right|\lesssim \|\mathcal{G}\|_{ C^{1,\log^\varkappa}}. 
			\end{align*}
			hence we get 
			\begin{align*}
				\|I_2\|_{L^\infty}\lesssim \|\mathcal{G}\|_{ C^{1,\log^\varkappa}}. 
			\end{align*}
			Combining this with \eqref{i1xz} and Sobolev interpolation inequality, we get 
			\begin{align*}
				\|\mathbf{V}_\b^2\|_{L^\infty}\lesssim  \|I_1\|_{L^\infty} +  \|I_2\|_{L^\infty}\lesssim \|\mathcal{G}\|_{\dot C^{1,\log^\varkappa}}+\|\mathcal{G}\|_{ L^2}. 
			\end{align*}
			Combining  this with \eqref{V1blinf}, we complete the proof.
		\end{proof}\vspace{0.1cm}\\
		The following lemma indicates that when considering any domain far from the interface $\{z=0\}$, the contribution of $\tilde v_\b^\pm$ are lower order terms.
		\begin{lemma}\label{estilv1low}
			Let $\tilde v_\b^\pm[f]$ be as defined in \eqref{fortvb}. Then for any $\sigma>0$, $n\in\mathbb{N}$,
			\begin{align*}
				\sup_\b\|\nabla_x^n\nabla_{x,z}\tilde v_\b^\pm[f]\|_{L^\infty( \mathbb{R}^{d+1}_\pm\cap\{|z|\geq \sigma\})}&+		\sup_\b\|\nabla_x^n\nabla_{x,z}\tilde v_\b^\pm[f]\|_{\cdot C^a( \mathbb{R}^{d+1}_\pm\cap\{|z|\geq \sigma\})}\\
				&\leq C(\sigma)(1+\|\nabla\phi\|_{L^\infty})^{d+n+10}\|\nabla f\|_{\dot C^1}(1+\|\nabla f\|_{L^\infty})^5.
			\end{align*}
		\end{lemma}
		\begin{proof}
			It follows from \eqref{deriK} that 
			\begin{align*}
				\sup_{|z|\geq \sigma}|\nabla_x^{n_1}\partial_z^{n_2}K_\b(x,z)|\lesssim \sigma^{-(d+n_1+n_2)}\langle\b\rangle^{d+1+2(n_1+n_2)},\ \ \ \forall n_1,n_2\in\mathbb{N}.
			\end{align*}
			Then by the formula of $\tilde v_\b^\pm[f]$ in \eqref{fortvb}, one has 
			\begin{align*}
				\sup_{|z|\geq \sigma}|\nabla_x^n\nabla_{x,z}\tilde v_\b^\pm|\lesssim \|\nabla_x^n\nabla_{x,z}K_\b(x,z)\|_{L^1_x} \|\mathcal{G}[f]\|_{L^\infty}\lesssim C(\sigma)\langle\b\rangle^{d+n+10}\|\mathcal{G}[f]\|_{L^\infty}.
			\end{align*}
			Combining this with the fact that 
			\begin{align*}
				\|\mathcal{G}[f]\|_{L^\infty}\lesssim \|\nabla f\|_{\dot C^1}(1+\|\nabla f\|_{L^\infty})^5,
			\end{align*}
			we complete the proof.
		\end{proof}
		
		~\vspace{0.3cm}\\
		With Lemma \ref{lemtvbhol}, \ref{lemsiG} and \ref{estilv1low} in hand, we are now ready to deal with the elliptic system \eqref{elliw}.
		For readers' convenience, we rewrite the system \eqref{elliw} of $\omega_\b^\pm=\omega_\b^\pm[f]$  here.
		\begin{equation*}
			\begin{aligned}
				&\operatorname{div}_{x,z}(\M(\nabla f)\nabla_{x,z} \omega^\pm_{\b})=\operatorname{div}_{x,z} F_{1,\b}^\pm,\ \ \text{in}\ \tilde\Omega_\eta^\pm,\\
				&\omega^-_{\b}-\omega^+_{\b}=F_2, \quad \text { on } \{z=0\},\\
				&e_{d+1}\cdot (\M(\nabla f)\nabla_{x,z} (\mu_+^{-1}\omega_\mathbf{b}^+-\mu_-^{-1}\omega^-_{\b}))=F_{3,\b}, \quad \text { on } \{z=0\},\\
				&\tilde \nu^\pm\cdot (\M(\nabla f)\nabla _{x,z} \omega^\pm_{\b}) =F_{4,\b},\ \ \quad\quad\quad\quad\ \text{on}\  \tilde \Gamma_f^\pm,
			\end{aligned}
		\end{equation*}	
		where 
		\begin{equation}\label{fortems}
			\begin{aligned}
				&F_{1,\b}^\pm=\left(\M(\mathbf{b})-\M(\nabla f)\right)\nabla_{x,z} \tilde v^\pm_\b[f],\\
				&F_2=\varrho_0f,\\
				&F_{3,\b}=\left.e_{d+1}\cdot \left((\mu_+^{-1}F_{1,\b}^+-\mu_-^{-1}F_{1,\b}^-)\right)\right|_{z=0},\\
				&F_{4,\b}^\pm=\left.\tilde \nu^\pm\cdot (\M(\nabla f)(\nabla_{x,z}  v^\pm_f-\nabla_{x,z}\tilde v_{\mathbf{b}}^\pm[f])) \right|_{ \tilde \Gamma_f^\pm}.
			\end{aligned}
		\end{equation}
		Recalling the definition of norm $|||\cdot|||_{a,\mathbb{K}}$ in \eqref{defnormbb},	we have the following result concerning the Holder estimates of the system \eqref{elliw}.
		\begin{lemma}\label{lemH2w}Consider $\|\cdot\|_T$, $\|\cdot\|_{X_T}$ defined in \eqref{defnorgm}, and 
			$\mathcal{X}^\sigma_{T,\phi}$ defined in \eqref{defsp}.	Let $\omega^\pm_\b[f]$ be the weak solution to system \eqref{elliw}. Then there exist $\sigma, T \ll 1$ depending on $\|\nabla\phi\|_{C^{d+m+4}}$ and $\mathbf{r}$ such that \\
			i) For any $f\in \mathcal{X}^\sigma_{T,\phi}$,
			\begin{align}\label{omeb0}
				& \sum_{+,-} \sup_{t\in[0,T]}\sup_\b t^\frac{2}{3}\|\nabla_{x,z} \omega_\b^\pm[f]\|_{L^\infty(\mathbb{R}^{d+1}_{\pm2 r})}\lesssim (\|f-\phi\|_{X_T}+T^\frac{1}{10})(1+\|f\|_{X_T})^5,\\
				&\sum_{+,-}\sup_{t\in[0,T]}|\log t|^\varkappa t^\frac{m+\kappa}{3} \|\nabla_{x}^m\nabla_{x,z} \omega_\b^\pm[f]\|_{\kappa-2,\mathbb{R}^{d+1}_{\pm r}}	\lesssim \| f-\phi\|_{T}( \| f-\phi\|_{T}+T^\frac{1}{10})+T^\frac{1}{100}\| f\|_{X_T}.\label{omeba}
			\end{align}
			ii) For $f_1,f_2\in \mathcal{X}^\sigma_{T,\phi}$, 
			\begin{align*}
				&\sum_{+,-} \sup_{t\in[0,T]}\sup_\b t^\frac{2}{3}\|\nabla_{x,z} (\omega_\b^\pm[f_1]-\omega_\b^\pm[f_2])\|_{L^\infty(\mathbb{R}^{d+1}_{\pm2 r})}\lesssim \|f_1-f_2\|_{X_T}(1+\|(f_1,f_2)\|_{X_T})^5,\\	&\sum_{+,-}\sup_{t\in[0,T]}|\log t|^\varkappa t^\frac{m+\kappa}{3} |||\nabla_x^m\nabla_{x,z} (\omega_\b^\pm[f_1]-\omega_\b^\pm[f_2])
				|||_{\kappa-2,\mathbb{R}^{d+1}_{\pm r}}	\lesssim ( T^\frac{1}{10}+\|f_1-\phi\|_{T}+\|f_2-\phi\|_{T})\|f_1-f_2\|_{X_T},
			\end{align*}
			where $r=\frac{\mathbf{r}}{50}$.
		\end{lemma}
		\begin{proof} We denote $a=\kappa-2$ and $\omega^\pm_\b=\omega^\pm_\b[f]$ for short. 
			Moreover, we use bold letters $\x=(x,z)$, $\y=(y,w),$ to denote $d+1$ dimensional variables.
			
			We first estimate the Lipschitz norm $\|\nabla_{x,z} \omega_\b^\pm[f]\|_{L^\infty(\mathbb{R}^{d+1}_{\pm2 r})}$. To avoid dealing with remainder term and the issue of losing derivatives, we go back to the equation of $v_f^\pm$ rather than estimating $\omega_\b$ directly. From \eqref{elleqeta}, we have 
			\begin{align*}
				&\operatorname{div}_{x,z}(\M(\nabla f)\nabla_{x,z} v^\pm)=0,\ \ \text{in}\ \mathbb{R}^{d+1}_{-4r}\cup \mathbb{R}^{d+1}_{4r},\\
				&v^--v^+=\mathcal{G}[f]+\varrho_0f, \quad \text { on } \{z=0\},\\
				&e_{d+1}\cdot (\M(\nabla f)\nabla_{x,z} (\mu_+^{-1}v^+-\mu_-^{-1}v^-))=0, \quad \text { on } \{z=0\}.
			\end{align*}
			Applying Lemma \ref{leminterf} with $A(x)=\M(\nabla f(x))$, $(g,h)=(0,\mathcal{G}[f]+\varrho_0f)$, we obtain that for any $x_0\in\mathbb{R}^d$,
			\begin{equation}\label{forvvv}
				\begin{aligned}
					v(x,z)=&\frac{\mu_+\mathbf{1}_{z\in(0,2r)}+\mu_-\mathbf{1}_{z\in(-2r,0)}}{\mu_++\mu_-}\left(\int_{\mathbb{R}^d}e_{d+1}\cdot(\M(\nabla f(x_0))\nabla_{y,w}\tilde\G_{x_0}(x,z,y,0)(\mathcal{G}[f]+\varrho_0f)(y))dy\right.\\
					&\quad\quad\quad\quad\quad\ \left.-\frac{1}{2r}\int_{2r}^{4r}\int_{\mathbb{R}^{d+1}_{\tau}\cup \mathbb{R}^{d+1}_{-\tau}} ((\M(\nabla f(x_0))-\M(\nabla f(y)))\nabla_{y,w}v(y,w))\nabla_{y,w}\tilde\G_{x_0}(x,z,y,w)dydwd\tau\right)\\
					&-\frac{\mu_+\mu_-}{\mu_++\mu_-}\sum_{+,-}\left(\frac{1}{2r\mu_\pm}\int_{2r}^{4r}\int_{\mathbb{R}^{d+1}_{\pm\tau}} ((\M(\nabla f(x_0))-\M(\nabla f(y)))\nabla_{y,w}v(y,w))\cdot\nabla_{y,w}\G_{x_0}(x,z,y,w)dydwd\tau\right)\\
					&+\mathrm{LO}(v, 0),
				\end{aligned}
			\end{equation}
			with the remainder term $\mathrm{LO}(v, 0)$ satisfying 
			\begin{align*}
				\|\nabla_{x,z}^n\mathrm{LO}(v, 0)\|_{L^\infty(\mathbb{R}^{d+1}_{2r})\cup\mathbb{R}^{d+1}_{-2r})}\lesssim \|v\|_{H^1(\mathbb{R}^{d+1}_{8r}\cup\mathbb{R}^{d+1}_{-8r})}.
			\end{align*}
			By Lemma \ref{lemLip}, we get
			\begin{align*}
				\|\nabla_{x,z} v\|_{L^\infty(\mathbb{R}^{d+1}_{-2r}\cup\mathbb{R}^{d+1}_{2r})}\lesssim \|\mathcal{G}[f]\|_{ C^{1,\log^\varkappa}}+\|f\|_{ C^{1,\log^\varkappa}}+c(r)\|\nabla_{x,z}v\|_{ L^\infty(\mathbb{R}^{d+1}_{-4r}\cup\mathbb{R}^{d+1}_{4r})}+C(r)\|\nabla_{x,z}v\|_{ L^2(\mathbb{R}^{d+1}_{-4r}\cup\mathbb{R}^{d+1}_{4r})}.
			\end{align*}
			Then we can take $r$ small enough, which implies $c(r)$ small enough to get 
			\begin{align*}
				\|\nabla_{x,z} v\|_{L^\infty(\mathbb{R}^{d+1}_{-2r}\cup\mathbb{R}^{d+1}_{2r})}\lesssim \|\mathcal{G}[f]\|_{ C^{1,\log^\varkappa}}+\|f\|_{ C^{1,\log^\varkappa}}+\|\nabla_{x,z}v\|_{ L^\infty(\mathbb{R}^{d}\times([-4r,-2r]\cup[2r,4r]))}+\|\nabla_{x,z}v\|_{ L^2(\mathbb{R}^{d+1}_{-4r}\cup\mathbb{R}^{d+1}_{4r})}.
			\end{align*}
			The right hand side involves the Lipschitz estimate of the solution in interior domain. To deal with this term, we go back to the function $q^\pm(x,z)=v^\pm(x,z-f(x))$, which is the solution to system \eqref{sysq}. Observe that $q^\pm$ in harmonic in interior domain. We have the following classical estimates of harmonic function  
			\begin{align*}
				\|\nabla_{x,z}q^\pm\|_{L^\infty(\Omega^\pm_r(t))}\lesssim r^{-\frac{d+1}{2}} \|\nabla_{x,z}q^\pm\|_{L^2(\Omega^\pm(t))},
			\end{align*}
			where $\Omega^\pm_r(t)=\{(x,z):B_r(x,z)\subset\Omega^\pm(t)\}$ denotes the interior of $\Omega^\pm(t)$.
			This yields that
			\begin{align*}
				\|\nabla_{x,z}v\|_{ L^\infty(\mathbb{R}^{d}\times([-4r,-2r]\cup[2r,4r]))}&\lesssim (1+\|\nabla f\|_{L^\infty})\|\nabla_{x,z}q^\pm\|_{{L^\infty(\Omega^\pm_r(t))}}\\
				&\lesssim r^{-\frac{d+1}{2}}(1+\|\nabla f\|_{L^\infty}) \|\nabla_{x,z}q^\pm\|_{L^2(\Omega^\pm(t))}\\
				&\lesssim r^{-\frac{d+1}{2}}(1+\|\nabla f\|_{L^\infty})^2 \|\nabla_{x,z}v^\pm\|_{L^2(\tilde \Omega^\pm)}.
			\end{align*}
			Hence, we get 
			\begin{align*}
				\|\nabla_{x,z} v\|_{L^\infty(\mathbb{R}^{d+1}_{-2r}\cup\mathbb{R}^{d+1}_{2r})}\lesssim \|\mathcal{G}[f]\|_{ C^{1,\log^\varkappa}}+\|f\|_{ C^{1,\log^\varkappa}}+(1+\|\nabla f\|_{L^\infty})^2\|\nabla_{x,z}v^\pm\|_{ L^2(\tilde \Omega^\pm)}.
			\end{align*}
			Combining this with Lemma \ref{lemtvbhol} and the definition $\omega_\b^\pm=v^\pm-\tilde v_\b^\pm$, we obtain 
			\begin{align}\label{vlinfty}
				\sup_\b\|\nabla_{x,z} \omega_\b^\pm[f]\|_{L^\infty(\mathbb{R}^{d+1}_{\pm2 r})}\lesssim  \|\mathcal{G}[f]\|_{ C^{1,\log^\varkappa}}+\|f\|_{ C^{1,\log^\varkappa}}+\|\nabla_{x,z}v\|_{L^2(\mathbb{R}^{d+1}_{-4r}\cup\mathbb{R}^{d+1}_{4r})}.
			\end{align}
			From this and \eqref{vl2}, together with the definition of $\|\cdot\|_{X_T}$ in \eqref{defnorgm}, we obtain \eqref{omeb0}.
			
			Then we prove \eqref{omeba}. For simplicity, we only consider 
			$|||\nabla_{x,z}\omega_\b^\pm[f]|||_{a,\mathbb{R}^{d+1}_{r}}$, the estimate of $|||\nabla_x^m\nabla_{x,z}\omega_\b^\pm[f]|||_{a,\mathbb{R}^{d+1}_{r}}$ follows by taking $m$ horizontal derivative on the equation.
			We consider the equation of $\omega_\b$ near the interface:
			\begin{align*}
				&\operatorname{div}_{x,z}(\M(\nabla f)\nabla_{x,z} \omega^\pm_{\b})=\operatorname{div}_{x,z} F_{1,\b}^\pm,\ \ \text{in}\ \mathbb{R}^{d+1}_{-2r}\cup \mathbb{R}^{d+1}_{2r},\\
				&\omega^-_{\b}-\omega^+_{\b}=F_2, \quad \text { on } \{z=0\},\\
				&e_{d+1}\cdot (\M(\nabla f)\nabla_{x,z} (\mu_+^{-1}\omega_\mathbf{b}^+-\mu_-^{-1}\omega^-_{\b}))=F_{3,\b}, \quad \text { on } \{z=0\},
			\end{align*}
			Applying Lemma \ref{leminterf} with $A(x)=\M(\nabla f(x))$, $(g,h)=(F_{1,\b},F_2)$, we obtain that for any $x_0\in\mathbb{R}^d$,
			\begin{align*}
				\omega_\b(x,z)=\T_{1,x_0}F_{1,\b}(x,z)+ \T_{2,x_0} F_2(x,z)+\R_{x_0}\omega_\b(x,z)+\mathrm{LO}(\omega_\b, F_{1,\b}),
			\end{align*}
			with
			\begin{align*}
				&\T_{1,x_0}F_{1,\b}(x,z)=-\frac{\mu_+\mathbf{1}_{z\in(0,r)}+\mu_-\mathbf{1}_{z\in(-r,0)}}{\mu_++\mu_-}\frac{1}{r}\int_{r}^{2r}\int_{R^{d+1}_{\tau}\cup R^{d+1}_{-\tau}} F_{1,\b}(y,w)\cdot\nabla_{y,w}\tilde\G_{x_0}(x,z,y,w)dydwd\tau\\
				&\quad\quad\quad-\frac{\mu_+\mu_-}{\mu_++\mu_-}\sum_{+,-}\left(\frac{1}{r\mu_\pm}\int_{r}^{2r}\int_{R^{d+1}_{\pm\tau}} F_{1,\b}(y,w)\cdot\nabla_{y,w}\G_{x_0}(x,z,y,w)dydwd\tau\right),\\
				&\T_{2,x_0}F_{2}(x,z)=\frac{\mu_+\mathbf{1}_{z\in(0,r)}+\mu_-\mathbf{1}_{z\in(-r,0)}}{\mu_++\mu_-}\int_{\mathbb{R}^d}e_{d+1}\cdot(\M(\nabla f(x_0))\nabla_{y,w}\tilde\G_{x_0}(x,z,y,0)F_2(y))dy,\\
				&\R_{x_0}\omega_\b(x,z)=\T_{1,x_0}((\M(\nabla f(x_0))-\M(\nabla f))\nabla_{x,z}\omega_\b)(x,z).
			\end{align*}
			and the lower order term $\mathrm{LO}(\omega_\b,F_{1,\b})$ satisfies
			\begin{align}\label{rereLO}
				\|\nabla_{x,z}^n\mathrm{LO}(\omega_\b,F_{1,\b})\|_{L^\infty(\mathbb{R}^{d+1}_{r/2}\cup\mathbb{R}^{d+1}_{-r/2})}\lesssim C(r,n) (\|\omega_\b\|_{H^1(\mathbb{R}^{d+1}_{2r}\cup \mathbb{R}^{d+1}_{-2r})}+\| F_{1,\b}\|_{L^2(\mathbb{R}^{d+1}_{2r}\cup\mathbb{R}^{d+1}_{-2r})}),\ \ \forall n\in\mathbb{N}.
			\end{align}
			By Remark \ref{remF2} and \eqref{hhhkf} in Lemma \ref{lemKf}, we have 
			\begin{align}\label{rereF2}
				\sup_{x_0}\|\nabla_{x,z}\T_{2,x_0}F_{2}\|_{\dot C^a(\mathbb{R}^{d+1}_{r}\cup\mathbb{R}^{d+1}_{-r})}\lesssim (1+\|\nabla f\|_{L^\infty})^2\|F_2\|_{\dot C^{1+a}}\lesssim (1+\|\nabla f\|_{L^\infty})^2\|f\|_{\dot C^{1+a}}.
			\end{align}
			To control $\nabla_{x,z}\T_{1,x_0}F_{1,\b}$, it suffices to consider  
			\begin{align*}
				&\T_{1,x_0}^i F_{1,\b}(x,z)=\int_{\mathbb{R}^{d+1}}K_i(x,z,y,w)F_{1,\b}(y,w)dydw,\ \ i=1,2,
			\end{align*}
			with 
			\begin{align*}
				&K_1(x,z,y,w)=\nabla_{x,z}\nabla_{y,w}\tilde \G_{x_0}(x,z,y,w)\mathbf{1}_{\mathbb{R}^{d+1}_{[
						-\tau,\tau]}}(y,w),\ \ K_2(x,z,y,w)=\nabla_{x,z}\nabla_{y,w} \G_{x_0}(x,z,y,w)\mathbf{1}_{\mathbb{R}^{d+1}_{[
						-\tau,\tau]}}(y,w).
			\end{align*}
			By Lemma \ref{lemsgit}, we obtain for any $x_0\in\mathbb{R}^d$,
			\begin{equation}\label{Tf1}
				\begin{aligned}
					\sum_{i=1,2}\sup_{\substack{\x,\y\in\mathbb{R}^{d+1}_r\\
							|\x-\y|\leq r}}\frac{|(\T_{1,x_0}^i F_{1,\b}(\x)-\T_{1,x_0}^i F_{1,\b}(\y))|_{\b=\nabla\phi(x)}|}{|\x-\y|^a}\lesssim \sup_{\substack{\x,\y\in\mathbb{R}^{d+1}_r,x'\in\mathbb{R}^d\\
							|\x-\y|\leq {10r}\\
							|x'-x|\leq |\x-\y|}}\frac{|(F_{1,\b}(\x)-F_{1,\b}(\y))|_{\b=\nabla\phi(x')}|}{|\x-\y|^a}.
				\end{aligned}
			\end{equation}
			Recalling the definition of $F_{1,\b}$ in \eqref{fortems}, we have for $\x=(x,z), \y=(y,w)\in\mathbb{R}^{d+1}_r$,
			\begin{align*}
				F_{1,\b}(\x)-F_{1,\b}(\y)=& (\M(\nabla f(y))-\M(\nabla f(x)))\nabla_{x,z}\tilde v_\b[f](\y)\\
				&+(\M(\b)-\M(\nabla f(x)))(\nabla_{x,z}\tilde v_\b[f](\x)-\nabla_{x,z}\tilde v_\b[f](\y)).
			\end{align*}
			Note that 
			\begin{align*}
				|\M(\nabla f(x))-\M(\nabla \phi(x'))|&\leq |\M(\nabla f(x))-\M(\nabla \phi(x))|+|\M(\nabla \phi(x))-\M(\nabla \phi(x'))|\\
				&\lesssim (\|\nabla (f-\phi)\|_{L^\infty}+|x-x'|^a\|\nabla \phi\|_{\dot C^a})(1+\|(\nabla f,\nabla\phi)\|_{L^\infty}).
			\end{align*}
			Then applying Lemma \ref{lemtvbhol} to get 
			\begin{align*}
				&\sup_{\substack{\x,\y\in\mathbb{R}^{d+1}_r,x'\in\mathbb{R}^d\\
						|\x-\y|\leq {10r}\\
						|x'-x|\leq |\x-\y|}}\frac{|(F_{1,\b}(\x)-F_{1,\b}(\y))|_{\b=\nabla\phi(x')}|}{|\x-\y|^a}\\
				&\lesssim \left((\|\nabla f\|_{\dot C^a}+\|\nabla\phi\|_{\dot C^a})\sup_\b\|\nabla_{x,z}\tilde v_\b[f]\|_{L^\infty}+\|\nabla (f-\phi)\|_{L^\infty}\sup_\b\|\nabla_{x,z}\tilde v_\b[f]\|_{\dot C^a}\right)(1+\|(\nabla f,\nabla\phi)\|_{L^\infty})\\
				&\lesssim \left((\|\nabla f\|_{\dot C^a}+\|\nabla\phi\|_{\dot C^a})\|\mathcal{G}[f]\|_{\dot C^{1,\log^\varkappa}}+\|\nabla (f-\phi)\|_{L^\infty}\|\mathcal{G}[f]\|_{\dot C^{1+a}}\right)(1+\|(\nabla f,\nabla\phi)\|_{L^\infty})^3.
			\end{align*}
			Combining this with \eqref{Tf1}, we obtain 
			\begin{align}\label{rereF1}
				&\sup_{x_0}|||\nabla_{x,z}\T_{1,x_0}F_{1,\b}|||_{a,\mathbb{R}^{d+1}_{\pm r}}\lesssim \sum_{i=1,2}\sup_{x_0}|||\T_{1,x_0}^iF_{1,\b}|||_{a,\mathbb{R}^{d+1}_{\pm r}}\nonumber\\
				&\lesssim\left((\|\nabla f\|_{\dot C^a}+\|\nabla\phi\|_{\dot C^a})\|\mathcal{G}[f]\|_{\dot C^{1,\log^\varkappa}}+\|\nabla (f-\phi)\|_{L^\infty}\|\mathcal{G}[f]\|_{\dot C^{1+a}}\right)(1+\|(\nabla f,\nabla\phi)\|_{L^\infty})^3.
			\end{align}
			Finally, we consider the remainder term $\nabla_{x,z}\R_{x_0}\omega_\b$. We have 
			\begin{equation}\label{rx0}
				\begin{aligned}
					&{\left|\left(\nabla_{x,z}\R_{x_0}\omega_\b(\x)-\nabla_{x,z}\R_{x_0}\omega_\b(\x')\right)|_{x_0=x}\right|}\\
					&\quad\quad\quad\lesssim \sum_{i=1,2}\int_{\mathbb{R}^{d+1}}|(\tilde K_{i,x_0}(\x,\y)-\tilde K_{i,x_0}(\x',\y))|_{x_0=x}|d\y\|\nabla_{x,z}\tilde \omega_\b\|_{L^\infty(\mathbb{R}^{d+1}_{2r}\mathbb{R}^{d+1}_{-2r})},
				\end{aligned}
			\end{equation}
			with 
			\begin{align*}
				&\tilde K_{1,x_0}(\x,\y)=\nabla_{x,z}\nabla_{y,w}\tilde \G_{x_0}(\x,\y)(\M(\nabla f(x_0))-\M(\nabla f(y)))\mathbf{1}_{\mathbb{R}^{d+1}_{\pm\tau}}(\y),\\ &\tilde K_{2,x_0}(\x,\y)=\nabla_{x,z}\nabla_{y,w} \G_{x_0}(\x,\y)(\M(\nabla f(x_0))-\M(\nabla f(y)))\mathbf{1}_{\mathbb{R}^{d+1}_{\pm\tau}}(\y).
			\end{align*}
			By the definition of $\tilde \G$ and $\G$ in , and 
			\begin{align*}
				|\M(\nabla f(x))-\M(\nabla f(y))|\lesssim |x-y|^a\|\nabla (f-\phi)\|_{\dot C^a},
			\end{align*}
			one has 
			\begin{align*}
				\sum_{i=1,2}\int_{\mathbb{R}^{d+1}}|(\tilde K_{i,x_0}(\x,\y)-\tilde K_{i,x_0}(\x',\y))|_{x_0=x}|d\y\lesssim |\x-\x'|^a\|\nabla f\|_{\cdot C^a}(1+\|\nabla f\|_{L^\infty}).
			\end{align*}
			Combining this with \eqref{rx0} to obtain
			\begin{equation}\label{rereR}
				\begin{aligned}
					&\sup_{\x,\x'}\frac{|\nabla_{x,z}\R_{x_0}\omega_\b(\x)-\nabla_{x,z}\R_{x_0}\omega_\b(\x')|_{x_0=x}|}{|\x-\x'|^a}\lesssim \|\nabla f\|_{\dot C^a}(1+\|\nabla f\|_{L^\infty})\|\nabla_{x,z}\tilde \omega_\b\|_{L^\infty(\mathbb{R}^{d+1}_{2r}\cup\mathbb{R}^{d+1}_{-2r})}\\
					&\quad\quad\quad\quad\quad\quad\lesssim\|\nabla f\|_{\dot C^a}(1+\|\nabla f\|_{L^\infty})(\|\mathcal{G}[f]\|_{\dot C^{1,\log^\varkappa}}+\|f\|_{\dot C^{1,\log^\varkappa}}+\|\nabla_{x,z}v\|_{L^2(\tilde \Omega_f)}),
				\end{aligned}
			\end{equation}
			where we applied \eqref{vlinfty} in the last inequality. We conclude from \eqref{rereLO}, \eqref{rereF2}, \eqref{rereF1}, and \eqref{rereR} that 
			\begin{align*}
				&\sum_{+,-}|||\nabla_{x,z}\omega_\b^\pm|||_{a,\mathbb{R}^{d+1}_{\pm r/2}}\\
				&\ \ \lesssim\left((\|\nabla f\|_{\dot C^a}+\|\nabla\phi\|_{\dot C^a})\|\mathcal{G}[f]\|_{\dot C^{1,\log^\varkappa}}+\|\nabla (f-\phi)\|_{L^\infty}\|\mathcal{G}[f]\|_{\dot C^{1+a}}\right)(1+\|(\nabla f,\nabla\phi)\|_{L^\infty})^3.\\
				&\quad\quad\quad+\|\nabla f\|_{\dot C^a}(1+\|\nabla f\|_{L^\infty})(\|\mathcal{G}[f]\|_{\dot C^{1,\log^\varkappa}}+\|f\|_{\dot C^{1,\log^\varkappa}}+\|\nabla_{x,z}v\|_{L^2(\tilde \Omega_f)})
			\end{align*}
			From this and \eqref{vl2}, together with the definition of $\|\cdot\|_{X_T}$ in \eqref{defnorgm}, we obtain \eqref{omeba}. This completes the proof of \textit{i)}. 
			
			The proof for \textit{ii)} is basically the same as for \textit{i)}. We give a brief explanation instead of a detailed proof to avoid redundancy. The only difference we have to mention is the estimate of Lipschitz norm, namely $\|\nabla_{x,z} (\omega_\b^\pm[f_1]-\omega_\b^\pm[f_2])\|_{L^\infty(\mathbb{R}^{d+1}_{\pm2 r})}$. As we know (see \eqref{gaes}), the Lipschitz estimate of solution to a second order elliptic equation will lose $\log^\varkappa$ derivative with $\kappa>1$. 
			In estimate \eqref{omeb0}, we essentially proved that  $\|\nabla_{x,z}\omega_\b^\pm[f]\|_{L^\infty}$ has the same estimate with  $\|\nabla_{x,z}\tilde v_\b^+\|_{L^\infty}$ by taking advantage of the fact that the equation for $v$ has no external force term, hence we did not lose derivative when estimating $\|\nabla_{x,z} v^\pm_f\|_{L^\infty}$. And  $\|\nabla_{x,z}\omega_\b^\pm[f]\|_{L^\infty}$ inherits the loss of $\log^\varkappa$ from $\|\nabla_{x,z}\tilde v_\b^\pm[f]\|_{L^\infty}$. However, when considering the contraction estimates, we have to make difference of two elliptic systems with different variable coefficients, then we cannot avoid the force term that arising from the difference of coefficients. More precisely, we have 
			\begin{align*}
				&\operatorname{div}_{x,z}(\M(\nabla f_1)\nabla_{x,z}(v^\pm_{f_1}-v_{f_2}^\pm))=\operatorname{div}_{x,z}\tilde F_{1,\b}^\pm,\ \ \ \text{in}\ \mathbb{R}^{d+1}_{\pm 2r},\\
				&(v_{f_1}^--v_{f_2}^-)-(v_{f_1}^+-v_{f_2}^+)=\mathcal{G}[f_1]-\mathcal{G}[f_2]+\rho_0(f_1-f_2),\ \ \text{on}\ \{z=0\},\\
				&e_{d+1}\cdot (\M(\nabla f_1)\nabla_{x,z}(\mu_+^{-1}(v_{f_1}^+-v_{f_2}^+-\tilde F_{1,\b}^+)-\mu_-^{-1}(v_{f_1}^--v_{f_2}^--\tilde F_{1,\b}^-))=0,
			\end{align*}
			where 
			\begin{align*}
				\tilde F_{1,\b}^\pm=(\M(\nabla f_2)-\M(\nabla f_1))\nabla_{x,z}v_{f_2}^\pm.
			\end{align*}
			Then by Lemma \ref{leminterf}, we need to deal with terms in the following type
			\begin{equation}\label{escontr}
				\begin{aligned}
					&\int_{\mathbb{R}^{d+1}_{\pm\tau}}\nabla_{x,z}\nabla_{y,w}\tilde \G_{x_0}(x,z,y,w)\tilde F_{1,\b}^\pm(y,w)dydw\\
					&=\int_{\mathbb{R}^{d+1}_{\pm\tau}}\nabla_{x,z}\nabla_{y,w}\tilde \G_{x_0}(x,z,y,w)(\M(\nabla f_2(x))-\M(\nabla f_1(x)))\nabla_{x,z}v_{f_2}^\pm(y,w)dydw+R(x,z),
				\end{aligned}
			\end{equation}
			here we just write the trouble term explicitly, and the remainder term can be controlled by 
			\begin{align*}
				\|R\|_{L^\infty}\lesssim \|\nabla (f_1-f_2)\|_{\dot C^{\log ^\varkappa}}\|\nabla_{x,z}v_{f_2}^\pm\|_{L^\infty(\mathbb{R}^{d+1}_{\pm\tau})}.
			\end{align*}
			To estimate the trouble term, we need to substitute the formula of $v_{f_2}^\pm$ in \eqref{forvvv} to \eqref{escontr}. Then we have to estimate the convolution of data with  double singular integral. In Lemma   \ref{lemsiG}, we proved that the convolution with  double singular integral loses $\log^\varkappa$ derivatives instead of $\log^{2\varkappa}$ derivatives. Then we obtain the contraction estimate in \textit{ii)}. This completes the proof.
		\end{proof}~\vspace{0.3cm}\\
		\begin{lemma}
			\label{lemh2}
			Let $\H_\b^{im}[f]$ be as defined in \eqref{defHim}, the norm $\|\cdot\|_{T}$ be as defined in \eqref{defnorgm} with constants $\varkappa,m,\kappa$ fixed in \eqref{consgm}. Then there exist $\sigma,T\ll 1$ depending on $\|\nabla\phi\|_{C^{m+4}}, \|\nabla \underline{b}^\pm\|_{C^1}$ and $\mathbf{r}$ such that  for any $f,f_1,f_2\in \mathcal{X}^\sigma_{T,\phi}$, 
			\begin{align}
				&\sup_{t\in[0,T]}|\log t|^\varkappa  (t^\frac{\kappa}{3}|||\H_{\b}^{im}[f](t)|||_{{\kappa-2},\mathbb{R}^{d+1}_{\pm r}}+t^\frac{m+\kappa}{3}|||\nabla_x^m\H_\b^{im}[f](t)|||_{\kappa-2,\mathbb{R}^{d+1}_{\pm r}})\nonumber\\
				&\quad\quad\quad\lesssim\|f-\phi\|_T(\|f-\phi\|_T+T^\frac{1}{20}\|\nabla\phi\|_{C^{m+3}})(1+\|(f,\phi)\|_T)^{2m+5}+T^\frac{1}{5}(1+\|\nabla\phi\|_{C^{n+3}})^{m+3}, \label{m2}\\
				&\sup_{t\in[0,T]}|\log t|^\varkappa  (t^\frac{\kappa}{3}||| (\H_\b^{im}[f_1]-\H_{\b}^{im}[f_2])(t)|||_{{\kappa-2},\mathbb{R}^{d+1}_{\pm r}}+t^\frac{m+\kappa}{3}\|\nabla_x^m (\H_\b^{im}[f_1]-\H_\b^{im}[f_2])(t)\|_{\kappa-2,\mathbb{R}^{d+1}_{\pm r}})\nonumber\\
				&\quad\quad\quad\lesssim \|f_1-f_2\|_{T}(\|(f_1-\phi,f_2-\phi)\|_{T}+T^\frac{1}{20}\|\nabla \phi\|_{C^{m+3}})(1+\|( f_1,f_2,\phi)\|_{T})^{2m+5},\label{dif2}
			\end{align}
			where $r=\frac{\mathbf{r}}{50}$.
		\end{lemma}
		\begin{proof}
			By the definition \eqref{defHim}, and Lemma \ref{Hol}, we have for $j=0,m$,
			\begin{align*}
				&|||\nabla_x^j\H_\b^{im}[f]|||_{\kappa-2,\tilde\Omega_t}\\
				&\ \ \ \ \lesssim \sum_{j_1+j_2=j}\left(\|\M(\nabla f)\|_{\dot C^{j_1}}|||\nabla_x^{j_2}\nabla_{x,z}\omega_\b^-[f]|||_{\kappa-2,\mathbb{R}^{d+1}_{\pm r}}+\|\M(\nabla f)\|_{\dot C^{j_1+\kappa-2}}|||\nabla_x^{j_2}\nabla_{x,z}\omega_\b^-[f]|||_{0,\mathbb{R}^{d+1}_{\pm r}}\right),
			\end{align*}
			and
			\begin{align*}
				|||\nabla_x^j(\H_\b^{im}[f_1]-
				\H_\b^{im}[f_2])|||_{ \kappa-2,\mathbb{R}^{d+1}_{\pm r}}&\lesssim \sum_{j_1+j_2=j}\left(\|\M(\nabla f_1)-\M(\nabla f_2)\|_{\dot C^{j_1}}|||\nabla_x^{j_2}\nabla_{x,z}\omega_\b^-[f_1]|||_{\kappa-2,\mathbb{R}^{d+1}_{\pm r}}\right.\\
				&\left.\quad+\|\M(\nabla f_1)-\M(\nabla f_2)\|_{\dot C^{j_1+\kappa-2}}|||\nabla_x^{j_2}\nabla_{x,z}\omega_\b^-[f_1]|||_{0,\mathbb{R}^{d+1}_{\pm r}}\right)\\
				&\quad+\sum_{j_1+j_2=j}\left(\|\M(\nabla f_2)\|_{\dot C^{j_1}}|||\nabla_x^{j_2}(\nabla_{x,z}\omega_\b^-[f_1]-\nabla_{x,z}\omega_\b^-[{f_2}])|||_{\kappa-2,\mathbb{R}^{d+1}_{\pm r}}\right.\\
				&   \left.\quad+\|\M(\nabla f_2)\|_{\dot C^{j_1+\kappa-2}}|||\nabla_x^{j_2}(\nabla_{x,z}\omega_\b^-[f_1]-\nabla_{x,z}\omega_\b^-[f_2])|||_{0,\mathbb{R}^{d+1}_{\pm r}}\right).
			\end{align*}
			Applying Lemma \ref{lemcom} to obtain 
			\begin{align*}
				&\|\M(\nabla f)\|_{\dot C^{j_1}}\lesssim \|\nabla f\|_{\dot C^1}^{j_1}+\|\nabla f\|_{\dot C^{j_1}},\\
				&\|\M(\nabla f)\|_{\dot C^{j_1+\kappa-2}}\lesssim\|\nabla f\|_{\dot C^{\kappa-2}}^{\frac{j_1+\kappa-2}{\kappa-2}}+\|\nabla f\|_{\dot C^{j_1+\kappa-2}},\\
				&\|\M(\nabla f_1)-\M(\nabla f_2)\|_{\dot C^{j_1}}\lesssim \sum_{l=0}^{j_1}\|\nabla(f_1-f_2)\|_{\dot C^l}(\|(\nabla f_1,\nabla f_2)\|_{\dot C^1}^{j_1-l}+\|(\nabla f_1,\nabla f_2)\|_{\dot C^{j_1-l}}),\\
				&\|\M(\nabla f_1)-\M(\nabla f_2)\|_{\dot C^{j_1+a}}\lesssim \sum_{l=0}^{j_1}\left\{\|\nabla(f_1-f_2)\|_{\dot C^{l+a}}(\|(\nabla f_1,\nabla f_2)\|_{\dot C^1}^{j_1-l}+\|(\nabla f_1,\nabla f_2)\|_{\dot C^{j_1-l}})\right.\\
				&\quad\quad\quad\quad\quad\quad\quad\quad\quad\quad\quad\quad\quad\quad +\left.\|\nabla(f_1-f_2)\|_{\dot C^{l}}(\|(\nabla f_1,\nabla f_2)\|_{\dot C^a}^{\frac{j_1-l+a}{a}}+\|(\nabla f_1,\nabla f_2)\|_{\dot C^{j_1-l+a}})\right\}.
			\end{align*}
			Combining this with the estimates of $\nabla_{x,z}\omega_\b^-[f]$ in Lemma \ref{lemH2w}, we conclude that 
			\begin{align*}
				\sup_{t\in[0,T]}|\log t|^\varkappa t^\frac{j+\kappa}{3}|||\nabla_x^j\H_\b^{im}[f](t)|||_{\kappa-2,\mathbb{R}^{d+1}_{\pm r}}\lesssim (1+\|f\|_T)^{j+2}\left(\|f-\phi\|_T(\|f-\phi\|_T+T^\frac{1}{10})+T^\frac{1}{100}\|f\|_{X_T}\right),
			\end{align*}
			and 
			\begin{align*}
				& \sup_{t\in[0,T]}|\log t|^\varkappa t^\frac{j+\kappa}{3}\|\nabla_x^j(\H_\b^{im}[f_1]-\H_\b^{im}[f_2])(t)|||_{\kappa-2,\mathbb{R}^{d+1}_{\pm r}}\\
				&\quad\quad\lesssim (1+\|(f_1,f_2)\|_T)^{j+2}\|f_1-f_2\|_{X_T}(\|(f_1-\phi,f_2-\phi)\|_T+T^\frac{1}{100})(1+\|(f_1,f_2)\|_T)^{j+2}.
			\end{align*}
			This completes the proof.
		\end{proof}\vspace{0.1cm}\\
		Finally, we control the lower order $L^2$ norm.
		The following lemmas helps to control the lower frequency of solution to system \eqref{elliw}.
		\begin{lemma}\label{tilv1L2}
			Let $ \mathcal{V}_\b$ be  defined in \eqref{defvb}. Then 
			\begin{align*}
				\|\nabla_{x,z}\mathcal{V}_\b[f]\|_{L^2(\mathbb{R}^{d+1}_\pm)}\lesssim  \langle \b\rangle
				\|f\|_{\dot H^\frac{1}{2}}.
			\end{align*}
		\end{lemma}
		\begin{proof}
			For simplicity, we only estimate $	\|\nabla_{x,z}\mathcal{V}_\b[f]\|_{L^2(\mathbb{R}^{d+1}_+)}$.
			By definition, taking Fourier transform in $x$ to obtain 
			\begin{align*}
				&\mathcal{F}\left(\nabla_ x\mathcal{V}_\b[f]\right)(\xi,z)=i\xi\hat f(\xi) e^{\lambda^\mp(\xi,\b) z},\\
				&\mathcal{F}\left(\partial_z\mathcal{V}_\b[f]\right)(\xi,z)=\lambda^\mp(\xi,\b)\hat f(\xi) e^{\lambda^\mp(\xi,\b) z},
			\end{align*}
			where $\lambda^\mp(\xi,\b)$ is defined in \eqref{deflanb}. 
			By Parseval's identity,
			\begin{align*}
				\|\nabla_{x}\mathcal{V}_\b[f]\|_{L^2(\mathbb{R}^{d+1}_\pm)}^2+\|\partial_z\mathcal{V}_\b[f]\|_{L^2(\mathbb{R}^{d+1}_\pm)}^2&=\int_{\mathbb{R}^\pm}\int_{\mathbb{R}^d}\left(\left|\mathcal{F}\left(\nabla_ x\mathcal{V}_\b[f]\right)(\xi,z)\right|^2+\left|\mathcal{F}\left(\partial_z\mathcal{V}_\b[f]\right)(\xi,z)\right|^2\right)d\xi dz\\
				&\quad\lesssim \int_{\mathbb{R}^\pm}\int_{\mathbb{R}^d}\left(|\xi|^2+|\lambda^\mp(\xi,\b)|^2\right)|\hat f(\xi)|^2|e^{2\lambda^\mp(\xi,\b) z}| d\xi dz.
			\end{align*}
			By the definition \eqref{deflanb}, one has
			\begin{align*}
				\frac{|\xi|}{\langle\b \rangle^2}\leq |\lambda^\mp(\xi,\b)|\leq |\xi|,\ \ \ \forall \b\in\mathbb{R}^d.
			\end{align*}
			Hence 
			\begin{align*}
				\|\nabla_{x,z}\mathcal{V}_\b[f]\|_{L^2(\mathbb{R}^{d+1}_\pm)}^2&\lesssim \int_0^\infty\int_{\mathbb{R}^d} |\xi|^2 |\hat f(\xi)|^2e^{-\frac{2|\xi|z}{\langle\b\rangle^2}}d\xi dz\\
				&\lesssim \langle\b\rangle^2\int_{\mathbb{R}^d} |\xi| |\hat f(\xi)|^2d\xi\\
				&\lesssim \langle\b\rangle^2\|f\|_{\dot H^\frac{1}{2}}^2.
			\end{align*}
			This completes the proof.
		\end{proof}
		\begin{lemma}\label{wL2}
			Let $v^\pm=v_f^\pm$, $\omega^\pm_\b=\omega^\pm_\b[f]$ be the weak solution to system \eqref{elleqeta} and  \eqref{elliw}, respectively, then 
			\begin{align}
				\label{vl2}	&	\sum_{+,-}		\|\nabla _{x,z}v^\pm\|_{L^2(\tilde \Omega_t^\pm)}
				\lesssim (1+\|\nabla f\|_{L^\infty})^5(\|\nabla f\|_{\dot C^\frac{3}{2}}+\|f\|_{L^2}),\\
				&	\sum_{+,-}	\sup_\b\|\nabla _{x,z}\omega^\pm_\b\|_{L^2(\tilde \Omega_t^\pm)}
				\lesssim (1+\|\nabla f\|_{L^\infty}+\|\nabla\phi\|_{C^1})^5(\|\nabla f\|_{\dot C^\frac{3}{2}}+\|f\|_{L^2}).\label{wbl2}
			\end{align}	
		\end{lemma}
		\begin{proof}
			Recall that $v^\pm=v^\pm_f$ is the weak solution to the elliptic system 
			\begin{equation}\label{eqmap}
				\begin{aligned}
					&\operatorname{div}_{x,z}(\M(\nabla f)\nabla_{x,z} v^\pm )=0,\ \ \text{in}\ \tilde \Omega^\pm(t),\\
					&v^--v^+=\mathcal{G}[f]+\varrho_0f, \quad \text { on } \{z=0\},\\
					&e_{d+1}\cdot (\M(\nabla f)(\mu_+^{-1}\nabla_{x,z}v^+-\mu_-^{-1}\nabla_{x,z}v^-))=0, \quad \text { on } \{z=0\},\\
					&\tilde \nu^\pm\cdot (\M(\nabla f)\nabla_{x,z} v^\pm) =0,\ \ \ \ \ \ \ \ \ \text{on} \ \tilde \Gamma^\pm_t,
				\end{aligned}
			\end{equation}
			Using $v^\pm$ as the test function for \eqref{eqmap}, we obtain  
			\begin{align*}
				0=&\mu_+^{-1}\int_{\tilde \Omega_t^+}\operatorname{div}_{x,z}(\M(\nabla f)\nabla_{x,z} v^+ ) v^+dxdz+\mu_-^{-1}\int_{\tilde \Omega_t^-}\operatorname{div}_{x,z}(\M(\nabla f)\nabla_{x,z} v^-) v^-dxdz\\
				=&-\mu_+^{-1}\int_{\tilde \Omega_t^+}(\nabla_{x,z} v^+)^\top   \M(\nabla f)\nabla_{x,z} v^+dxdz-\mu_-^{-1}\int_{\tilde \Omega_t^-}(\nabla_{x,z} v^-)^\top   \M(\nabla f)\nabla_{x,z} v^-dxdz\\
				&-\mu_+^{-1}\int_{\mathbb{R}^d}e_{d+1}\cdot (\M(\nabla f)\nabla_{x,z}v^+(x,0))v^+(x,0)dx+\mu_-^{-1}\int_{\mathbb{R}^d}e_{d+1}\cdot (\M(\nabla f)\nabla_{x,z}v^-(x,0))v^-(x,0)dx.
			\end{align*}
			By the ellipticity of $\M(\nabla f)$, we have 
			\begin{align*}
				&	\|\nabla _{x,z}v^+\|_{L^2(\tilde \Omega_t^+)}^2+\|\nabla _{x,z}v^-\|_{L^2(\tilde \Omega_t^-)}^2\\
				&\quad\leq \int_{\tilde \Omega_t^+}(\nabla_{x,z} v^+)^\top   \M(\nabla f)\nabla_{x,z} v^+dxdz+\int_{\tilde \Omega_t^-}(\nabla_{x,z} v^-)^\top   \M(\nabla f)\nabla_{x,z} v^-dxdz.
			\end{align*}
			On the other hand, with the boundary condition in \eqref{eqmap}, we obtain that 
			\begin{align*}
				&-\mu_+^{-1}\int_{\mathbb{R}^d}e_{d+1}\cdot (\M(\nabla f)\nabla_{x,z}v^+(x,0))v^+(x,0)dx+\mu_-^{-1}\int_{\mathbb{R}^d}e_{d+1}\cdot (\M(\nabla f)\nabla_{x,z}v^-(x,0))v^-(x,0)dx\\
				&\quad=\mu_+^{-1}\int_{\mathbb{R}^d}e_{d+1}\cdot (\M(\nabla f)\nabla_{x,z}v^+(x,0))\left(\mathcal{G}[f]+\varrho_0f\right)dx.
			\end{align*}
			By the trace theorem, there exists $\mathbf{f}:\mathbb{R}^{d+1}_+\to \mathbb{R}$ such that $\mathbf{f}(x,0)=\mathcal{G}[f](x)+\varrho_0f(x)$, and 
			\begin{align*}
				\|\mathbf{f}\|_{ H^1(\mathbb{R}^{d+1}_+)}\lesssim \left\|\mathcal{G}[f]+\varrho_0f\right\|_{  H^\frac{1}{2}},
			\end{align*}
			then apply Green's formula and \eqref{eqmap} we obtain
			\begin{align*}
				\int_{\mathbb{R}^d}e_{d+1}\cdot (\M(\nabla f)\nabla_{x,z}v^+(x,0))\left(\mathcal{G}[f]+\varrho_0f\right)(x)dx&= \int_{\tilde \Omega_t^+}\operatorname{div}_{x,z}\left(\M(\nabla f)\nabla_{x,z}v^+(x,z)\mathbf{f}(x,z)\right)dx\\
				&=\int_{\tilde \Omega_t^+} (\nabla _{x,z}\mathbf{f}(x,z))^\top  \M(\nabla f)\nabla_{x,z}v^+(x,z)dxdz.
			\end{align*}
			Hence, by H\"{o}lder's inequality, we get 
			\begin{align*}
				&	\|\nabla _{x,z}v^+\|_{L^2(\tilde \Omega_t^+)}^2+\|\nabla _{x,z}v^-\|_{L^2(\tilde \Omega_t^-)}^2\\
				&\quad\leq C\|\nabla_{x,z}v^+\|_{L^2(\tilde \Omega_t^+)}\|\nabla _{x,z}\mathbf{f}\|_{L^2(\mathbb{R}^{d+1})}\leq \frac{1}{10}\|\nabla_{x,z}v^+\|_{L^2(\tilde \Omega_t^+)}^2+C\left\|\mathcal{G}[f]+\varrho_0f\right\|_{  H^\frac{1}{2}}^2.
			\end{align*}
			Then we obtain 
			\begin{align*}
				\|\nabla _{x,z}v^+\|_{L^2(\tilde \Omega_t^+)}+\|\nabla _{x,z}v^-\|_{L^2(\tilde \Omega_t^-)}&\lesssim (1+\|\nabla f\|_{L^\infty})(\|\mathcal{G}[f]\|_{H^\frac{1}{2}}+\|f\|_{H^\frac{1}{2}})\\
				&\lesssim (1+\|\nabla f\|_{L^\infty})^5(\|\nabla f\|_{\dot C^\frac{3}{2}}+\|f\|_{L^2}).
			\end{align*}
			This completes the proof of \eqref{vl2}.
			On the other hand, Lemma \ref{tilv1L2} implies that 
			\begin{align}\label{tvL2}
				\|\nabla_{x,z}\tilde {\mathcal{V}}_\b^+\|_{L^2(\mathbb{R}^{d+1}_+)}+\|\nabla _{x,z}\tilde {\mathcal{V}}_\b^-\|_{L^2(\mathbb{R}^{d+1}_-)}\lesssim\langle \b\rangle \|\mathcal{G}[f]\|_{\dot H^\frac{1}{2}}. 
			\end{align}
			By  \eqref{vl2}, \eqref{tvL2} and the Sobolev interpolation inequality, we get 
			\begin{equation*}
				\begin{aligned}
					\|\nabla _{x,z}\omega^+_\b\|_{L^2(\tilde \Omega_t^+)}+\|\nabla _{x,z}\omega^-_\b\|_{L^2(\tilde \Omega_t^-)}	&\lesssim (1+\|\nabla f\|_{L^\infty}+\langle \b\rangle)^3(\|\mathcal{G}[f]\|_{H^\frac{1}{2}}+\|f\|_{H^\frac{1}{2}})\\
					&\lesssim (1+\|\nabla f\|_{L^\infty}+\langle \b\rangle)^5(\|\nabla f\|_{\dot C^\frac{3}{2}}+\|f\|_{L^2}).
				\end{aligned}
			\end{equation*}
			Then we obtain \eqref{wL2} and  complete the proof.
		\end{proof}\vspace{0.3cm}\\	
		
		\begin{lemma}\label{mapl2}
			Let  $T\in(0,1)$, $f\in \mathcal{X}^\sigma_{T,\phi}$ and $g=\mathcal{S}f$ be as defined in \eqref{infaeq1}. Then 
			\begin{align}\label{L2}
				\sup_{t\in[0,T]}\|g(t)-\phi\|_{L^2}\leq  \|\eta_0-\phi\|_{L^2}+C_1T^\frac{1}{50}.
			\end{align}
			Moreover, for $f_1,f_2\in \mathcal{X}^\sigma_{T,\phi}$, there holds
			\begin{align}\label{l2contr}
				\sup_{t\in[0,T]}\|\mathcal{S}f_1(t)-\mathcal{S}f_2(t)\|_{L^2}\leq C_2T^\frac{1}{50}\|f_1-f_2\|_{X_T}.
			\end{align}
			Here $C_1,C_2$ depends on $\|\eta_0\|_{L^2}$, $\|\nabla \eta_0\|_{\dot C^{\log ^\varkappa}}$ and $\|\nabla \phi\|_{C^1}$.
		\end{lemma}
		\begin{proof}
			For simplicity, we only show the proof of \eqref{L2}. And \eqref{l2contr} can be proved with minor changes of the proof.\\
			By \eqref{infaeq1}, and the fact that 
			\begin{align*}
				\H[f](x)=-\frac{1}{\mu_-}e_{d+1}\cdot  (\M(\nabla f)\nabla  _{x,z}v^-_f(x,z)|_{z=0})-\int_{\mathbb{R}^d} \A(\nabla \phi(x),\alpha):(\delta_\alpha \nabla^2 f(x))\frac{d\alpha}{|\alpha|^{d+1}},
			\end{align*}
			we have 
			\begin{align*}
				&\partial_t (g-\phi)(x)= \int_{\mathbb{R}^d} \A(\nabla \phi(x),\alpha):(\delta_{\alpha} \nabla^2 (g-f)(x))\frac{d\alpha}{|\alpha|^{d+1}}-\frac{1}{\mu_-}e_{d+1}\cdot  (\M(\nabla f)\nabla  _{x,z}v^-_f(x,z)|_{z=0}).
			\end{align*}
			Using $g-\phi$ as a test function to the above equation, we get
			\begin{align*}
				\frac{1}{2}\partial_t \|g(t)-\phi\|_{L^2}^2&=\int_{\mathbb{R}^d}\int_{\mathbb{R}^d} \A(\nabla \phi(x),\alpha):\mathcal{O}_{\alpha} \nabla^2 (g-f)(x))(g-\phi)(x)\frac{d\alpha dx}{|\alpha|^{d}}\\
				&\ \ -\frac{1}{\mu_-}\int_{\mathbb{R}^d}e_{d+1}\cdot  (\M(\nabla f)\nabla  _{x,z}v^-_f(x,z)|_{z=0})(g-\phi)(x)dx\\
				&:=I_1+I_2.
			\end{align*}
			Note that $|\A(\nabla \phi(x),\alpha)|\lesssim (1+\|\nabla\phi\|_{L^\infty})^2 $. By H\"{o}lder's inequality, we have 
			\begin{align*}
				I_1	&\lesssim (1+\|\nabla\phi\|_{L^\infty})^2\|g(t)-\phi\|_{L^2}\int_{\mathbb{R}^d}\|\mathcal{O}_{\alpha} \nabla^2 (g-f)\|_{L^2}\frac{d\alpha}{|\alpha|^{d}}.
			\end{align*}
			Let $\lambda>0$. It is easy to check that 
			\begin{align}\label{hf}
				\int_{|\alpha|\geq \lambda}\|\mathcal{O}_{\alpha} \nabla^2 (g-f)\|_{L^2}\frac{d\alpha}{|\alpha|^{d}}\lesssim \lambda^{-1}\|g-f\|_{\dot H^2}.
			\end{align}
			On the other hand, for high frequency part, by Holder's inequality we have 
			\begin{align*}
				\int_{|\alpha|\leq \lambda}\|\mathcal{O}_{\alpha} \nabla^2 (g-f)\|_{L^2}\frac{d\alpha}{|\alpha|^{d}}\lesssim	\lambda^{\frac{1}{2}}\left(\int_{|\alpha|\leq  \lambda}\|\mathcal{O}_{\alpha} \nabla^2 (g-f)\|_{L^2}^2\frac{d\alpha}{|\alpha|^{d+1}}\right)^{\frac{1}{2}}.
			\end{align*}
			By Parseval's identity, for any function $h$,
			\begin{align*}
				\|\mathcal{O}_\alpha h \|_{L^2}^2=\frac{1}{|\alpha|^2}\int_{\mathbb{R}^d}|\hat h(\xi)|^2|2-e^{i\alpha\cdot \xi}-e^{-i\alpha\cdot \xi}|^2d\xi\lesssim \frac{1}{|\alpha|^2} \int_{\mathbb{R}^d}|\hat h(\xi)|^2\min\{1,|\alpha|^4|\xi|^4\}d\xi.
			\end{align*}
			Hence we obtain 
			\begin{align*}
				\int_{|\alpha|\leq \lambda}\|\mathcal{O}_{\alpha} \nabla^2 (g-f)\|_{L^2}\frac{d\alpha}{|\alpha|^{d}}&\lesssim 	\lambda^\frac{1}{2}\left(\int_{|\alpha|\leq \lambda}\int_{\mathbb{R}^d}|\xi|^4|(\hat g-\hat f)(\xi)|^2\min\{1,|\alpha|^4|\xi|^4\}\frac{d\xi d\alpha}{|\alpha|^{d+3}}\right)^\frac{1}{2}\\
				&\lesssim \lambda^\frac{3}{4}\left(\int_{\mathbb{R}^d}|\hat g(\xi)-\hat f(\xi)|^2|\xi|^\frac{15}{2}d\xi\right)^\frac{1}{2}=\lambda^\frac{3}{4}\|g-f\|_{\dot H^\frac{15}{4}}.
			\end{align*}
			Combining this with \eqref{hf}, 	by taking $\lambda= \|g-f\|_{\dot H^2}^\frac{4}{7}\|g-f\|_{\dot H^\frac{15}{4}}^{-\frac{4}{7}}$ we obtain 
			\begin{align*}
				\int_{\mathbb{R}^d}\|\mathcal{O}_{\alpha} \nabla^2 (g-f)\|_{L^2}\frac{d\alpha}{|\alpha|^{d}}\lesssim \|g-f\|_{\dot H^2}^\frac{3}{7}\|g-f\|_{\dot H^\frac{15}{4}}^\frac{4}{7}\lesssim \|g-f\|_{L^2}+\|g-f\|_{\dot C^\frac{19}{5}},
			\end{align*}
			where we applied the Gagliardo--Nirenberg interpolation inequality in the last inequality. 
			Hence we obtain 
			\begin{equation}\label{iii11}	\begin{aligned}
					I_1	&\lesssim (1+\|\nabla\phi\|_{L^\infty})^2\|g-\phi\|_{L^2}(\|g-f\|_{L^2}+\|g-f\|_{\dot C^\frac{19}{5}})\\
					&\lesssim  (1+\|\nabla\phi\|_{L^\infty})^2\|g-\phi\|_{L^2}(1+t^{-\frac{14}{15}})(\|f\|_{X_T}+\|g\|_{X_T}).
				\end{aligned}
			\end{equation}
			Then we estimate $I_2$. 
			By the trace extension theorem, there exists $\bar g:\mathbb{R}^{d+1}_-\to \mathbb{R}$ such that 
			\begin{align*}
				\|\bar g\|_{H^1(\mathbb{R}^{d+1}_-)}\lesssim \|g-\phi\|_{H^\frac{1}{2}(\mathbb{R}^d)}, \ \ \ \ \ \  \bar g(x,0)=g(x)-\phi(x).
			\end{align*}
			Note that $\operatorname{div}_{x,z}  \left(\M(\nabla f)\nabla_{x,z} v^-_f(x,z)\right)=0$. Hence by Holder's inequality and \eqref{vl2},
			\begin{equation}\label{v21}
				\begin{aligned}
					&\left|I_2\right|=\frac{1}{\mu_-}\left|\int_{\tilde \Omega_t^-}\operatorname{div}_{x,z}  \left(\M(\nabla f)\nabla  _{x,z}v^-_f(x,z)\bar{g}(x,z)\right)dxdz\right| \\
					&\lesssim \int_{\tilde \Omega_t^-}\left|(\nabla _{x,z}\bar{g})^\top  \M(\nabla f)\nabla  _{x,z}v^-_f(x,z)\right|dxdz\lesssim (1+\|\nabla f\|_{L^\infty})\|\nabla_{x,z}\bar g\|_{L^2(\mathbb{R}^{d+1}_-)}\|\nabla_{x,z}v_f^-\|_{L^2(\mathbb{R}^{d+1}_-)}\\
					&\lesssim (1+\|\nabla f\|_{L^\infty})\|g-\phi\|_{H^\frac{1}{2}(\mathbb{R}^d)}\|\nabla_{x,z}v_f^-\|_{L^2(\mathbb{R}^{d+1}_-)}\lesssim  t^{-\frac{49}{50}}(1+\|\nabla f\|_{L^\infty})\|g-\phi\|_{X_T}\|f\|_{X_T}.
				\end{aligned}
			\end{equation}
			Combining this with \eqref{iii11}, and using the bounds $\|f\|_{X_T}+\|\nabla g\|_T+\|\nabla (g-\phi)\|_T+\|\phi\|_{L^2}\leq C$, where $C$ is a constant depending only on $\|\eta_0\|_{L^2}$ and $\|\nabla \eta_0\|_{\dot C^{\log ^\varkappa}}$, we obtain that 
			\begin{align*}
				\partial_t \|g(t)-\phi\|_{L^2}^2\lesssim |I_1|+|I_2|\lesssim \|g(t)-\phi\|_{L^2}^2+t^{-\frac{49}{50}}\|g(t)-\phi\|_{L^2}+1,
			\end{align*}
			where the implicit constant depends on $\|\eta_0\|_{L^2}$, $\|\nabla \eta_0\|_{\dot C^{\log ^\varkappa}}$ and $\|\nabla \phi\|_{C^1}$.
			Integrating  in time, and  using Gronwall's inequality,
			\begin{align*}
				\|g-\phi\|_{L^\infty_TL^2}^2\leq \|\eta_0-\phi\|_{L^2}^2+CT^\frac{1}{50}.
			\end{align*}
			This completes the proof.
		\end{proof}~\vspace{0.5cm}\\
		In the following we study the Schauder-type estimates for second order elliptic equations in divergence form with coefficients in $\dot C^{\log^\varkappa}, \varkappa>1$. The result is known in literature(see \cite{Don} for example). Here we give a short self-contained proof using the new freezing coefficient method introduced in Section \ref{secproof}. This  method can be extended to fractional  elliptic equations, but we will not discuss in more detail.	\vspace{0.1cm}\\
		We start from the  Laplace's equation, for which the fundamental solution satisfies 
		\begin{align}
			\label{fundlap}
			-\Delta \Phi(x)=\delta(x), \ \ x\in\mathbb{R}^d.
		\end{align}
		The fundamental solution has explicit formula
		\begin{equation}\label{defPhi}
			\Phi(x):=\begin{cases}
				&-\frac{1}{2\pi}\log|x|,\ \ \ d=2\\
				&\frac{1}{d(d-2)V(d)}\frac{1}{|x|^{d-2}},\ \ \ d\geq 3.
			\end{cases}\ \ \ \forall \ x\in\mathbb{R}^d\backslash\{0\}.
		\end{equation}
		Here $V(d)$ is the volume of the unit ball in $\mathbb{R}^d$. 
		We will sometimes slightly abuse notation and write $\Phi(x)=\Phi(|x|)$ to emphasize that the fundamental solution is radial.  We have the estimates
		\begin{align}\label{esfdphi}
			&	|\nabla^k \Phi(x)|\lesssim \frac{1}{|x|^{d-2+k}},\ \ \ \forall k\in\mathbb{N}^+,\\
			&	|\delta_\alpha \nabla^k \Phi(x)|\lesssim \left(\frac{1}{|x|^{d-2+k}}+\frac{1}{|x-\alpha|^{d-2+k}}\right)\min\left\{1,\frac{|\alpha|}{|x|}\right\},\quad\quad\quad \forall \alpha\in\mathbb{R}^d.\label{delphi}
		\end{align}
		Moreover, the fundamental solution has the following  cancellation property 
		\begin{align*}
			\int_{|x|=r}\partial_i\partial_j\Phi(x)dS(x)=0,\ \ \ \forall \ r>0,\ \ \forall \ i,j=1,\cdots,d.
		\end{align*}
		For the Poisson equation 
		\begin{align}\label{poisson}
			\Delta u=\operatorname{div}f.
		\end{align}
		The solution reads $u(x)=\int_{\mathbb{R}^d}\Phi(x-y)\operatorname{div}f(y)dy$. 
		By the classical Schauder theory, we have for $0<\alpha<1$,
		\begin{align}\label{ellhol}
			\|\nabla u\|_{\dot C^a}\lesssim\frac{1}{a(1-a)}\|f\|_{\dot C^a}.
		\end{align} 
		However, the $L^\infty\to L^\infty$ estimate does not hold. To control the $L^\infty$ norm of $\nabla u$, we need the $\dot C^{\log^\varkappa  }$ norm of $f$. In fact, integrating by parts to obtain 
		\begin{align*}
			\nabla u(x)&=\int_{\mathbb{R}^d}\nabla ^2\Phi(x-y)f(y)dy\\
			&=\int_{|x-y|\leq 1}\nabla ^2\Phi(x-y)(f(y)-f(x))dy+\int_{|x-y|\geq 1}\nabla ^2\Phi(x-y)f(y)dy\\
			:&=I_1+I_2.
		\end{align*} 
		It is easy to check that 
		\begin{align*}
			|I_1|\lesssim \int_{|x|\leq 1} \frac{1}{|x|^d\log^\varkappa  (|x|^{-1}+10)}dx\|f\|_{\dot C^{\log ^\varkappa}}\lesssim \|f\|_{\dot C^{\log ^\varkappa}},
		\end{align*}
		provided $\varkappa>1$. For lower frequency, we have 
		\begin{align*}
			|I_2|\lesssim \int_{|x-y|\geq 1}\frac{1}{|x-y|^d}|f(y)|dy\lesssim \|f\|_{L^2}.
		\end{align*}
		Hence we obtain the Lipschitz estimate to the Poisson equation \eqref{poisson}
		\begin{align}\label{gaes}
			\|\nabla u\|_{L^\infty}\lesssim  \|f\|_{\dot C^{\log ^\varkappa}}+\|f\|_{L^2}.
		\end{align}
		We remark that,  when considering elliptic equation with variable coefficients, the classical method to estimate equations with non-constant coefficient is to combine the standard freezing coefficient method together with  the estimate of constant coefficient equations. But this fails when we estimate the endpoint $L^\infty$ norm. Because the norm $\|\cdot\|_{\dot C^{\log^\varkappa  }}$ on the right hand side and the $\|\cdot\|_{L^\infty}$ norm on the left hand side of \eqref{gaes} do not match. Hence the remainder term for freezing coefficient cannot be absorbed by the main term on the left hand side. To obtain the $L^\infty$ estimate of solution to \eqref{poisson}, we need to find a suitable way to freeze the coefficient and estimate carefully. 
		\vspace{0.1cm}\\
		In the following lemma, we give a representation formula of solution to elliptic system with variable coefficient, on a domain with interface. To state the results clearly, we denote 
		\begin{align*}
			&B_r(x)=\{y\in\mathbb{R}^d:|y-x|\leq r\},\quad\quad B_r=B_r(0),\\ &\mathbb{R}^d_\tau=\mathbb{R}^{d-1}\times (0,\tau)\ \text{if}\ \tau>0,\ \ \mathbb{R}^d_\tau=\mathbb{R}^{d-1}\times (\tau,0)\ \text{if}\ \tau<0.
		\end{align*}
		Consider the system
		\begin{equation}\label{ellinte}
			\begin{aligned}
				&\operatorname{div}_{x,z}(A(x)\nabla_{x,z} u(x,z))=\operatorname{div}_{x,z}g(x,z),\ \ \text{in}\ \mathbb{R}^{d+1}_{2r}\cup \mathbb{R}^{d+1}_{-2r},\\
				&u(x,0^-)-u(x,0^+)=h(x),\ \ \text{on} \ \{z=0\},\\
				&e_{d+1}\cdot(A(x)(\mu_+^{-1}(\nabla_{x,z} u(x,0^+)-g(x,0^+))-\mu_-^{-1}\nabla_{x,z} (u(x,0^-)-g(x,0^-))))=0,\ \ \text{on} \ \{z=0\}.
			\end{aligned}
		\end{equation}
		Suppose the variable coefficient $A(x)\in \dot C^{\log^s}(\mathbb{R}^d,\mathbb{R}^{(d+1)\times(d+1)}_{sym})$ is uniformly elliptic:
		$$
		c_0	\mathrm{Id}\leq A(x)\leq c_0^{-1} \mathrm{Id},\ \ \forall x\in\mathbb{R}^d.
		$$
		For any $x_0\in\mathbb{R}^d$, there exists an invertible matrix satisfying $Q_{x_0}\in\mathbb{R}^{(d+1)\times(d+1)}$ such that  $Q_{x_0}Q_{x_0}^\top=A(x_0)$, and $(Q_{x_0})_{d+1,i}=0$ for $i=1,\cdots,d$. Let $\Phi_{x_0}(x,z)=\Phi(Q_{x_0}^{-1}(x,z)^\top)$, where $\Phi$ is defined in  \eqref{defPhi}.
		For $(x,z),(y,w)\in\mathbb{R}^d\times \mathbb{R}$. Define
		\begin{equation}\label{defGtiG}
			\begin{aligned}
				&\G_{x_0}(x,z,y,w)=|\operatorname{det}Q_{x_0}|^{-1}(\Phi_{x_0}(x-y,z-w)+\Phi_{x_0}(x-y,z+w)),\\
				&\tilde \G_{x_0}(x,z,y,w)=|\operatorname{det}Q_{x_0}|^{-1}(\Phi_{x_0}(x-y,z-w)-\Phi_{x_0}(x-y,z+w)).
			\end{aligned}
		\end{equation}
		Then one has 
		\begin{align}\label{Gbdy}
			\tilde \G_{x_0}(x,z,y,w)|_{w=0}=0,\quad\quad \ e_{d+1}\cdot(A(x_0)\nabla_{y,w}\G_{x_0}(x,z,y,w)|_{w=0}=0,
		\end{align}
		and
		\begin{align}
			&\operatorname{div}_{y,w}(A(x_0)\nabla_{y,w}\G_{x_0}(x,z,y,w))=\delta(x-y,z-w)+\delta(x-y,z+w),\label{eqG}\\
			&\operatorname{div}_{y,w}(A(x_0)\nabla_{y,w}\tilde\G_{x_0}(x,z,y,w))=\delta(x-y,z-w)-\delta(x-y,z+w).\label{eqtilG}
		\end{align}
		\begin{lemma}\label{leminterf}
			Let $u\in W^{1,2}(\mathbb{R}^{d+1}_{2r}\cup\mathbb{R}^{d+1}_{-2r})$ be a weak solution to \eqref{ellinte}. 
			Then for any $x_0\in\mathbb{R}^d$, there holds
			\begin{equation}\label{foruintf}
				\begin{aligned}
					u(x,z)=&\frac{\mu_+\mathbf{1}_{z\in(0,r)}+\mu_-\mathbf{1}_{z\in(-r,0)}}{\mu_++\mu_-}\left(\int_{\mathbb{R}^d}e_{d+1}\cdot(A(x_0)\nabla_{y,w}\tilde\G_{x_0}(x,z,y,0)h(y))dy\right.\\
					&\quad\quad\quad\quad\quad-\frac{1}{r}\int_{r}^{2r}\int_{R^{d+1}_{\tau}\cup R^{d+1}_{-\tau}} g(y,w)\cdot\nabla_{y,w}\tilde\G_{x_0}(x,z,y,w)dydwd\tau\\
					&\quad\quad\quad\quad\quad\ \left.-\frac{1}{r}\int_{r}^{2r}\int_{R^{d+1}_{\tau}\cup R^{d+1}_{-\tau}} ((A(x_0)-A(y))\nabla_{y,w}u(y,w))\cdot\nabla_{y,w}\tilde\G_{x_0}(x,z,y,w)dydwd\tau\right)\\
					&-\frac{\mu_+\mu_-}{\mu_++\mu_-}\sum_{+,-}\left(\frac{1}{r\mu_\pm}\int_{r}^{2r}\int_{R^{d+1}_{\pm\tau}} g(y,w)\cdot\nabla_{y,w}\G_{x_0}(x,z,y,w)dydwd\tau\right.\\
					&\quad\quad\quad\quad\quad\quad\left.+\frac{1}{r\mu_\pm}\int_{r}^{2r}\int_{R^{d+1}_{\pm\tau}} ((A(x_0)-A(y))\nabla_{y,w}u(y,w))\cdot\nabla_{y,w}\G_{x_0}(x,z,y,w)dydwd\tau\right)\\
					&+\mathrm{LO}(u, g),
				\end{aligned}
			\end{equation}
			with remainder term $\mathrm{LO}(u, g)$ satisfying \begin{align*}
				\|\nabla_{x,z}^n\mathrm{LO}(u,g)\|_{L^\infty(\mathbb{R}^{d+1}_{r/2})}\lesssim C(r,n) (\|u\|_{H^1(\mathbb{R}^{d+1}_{2r}\cup \mathbb{R}^{d+1}_{-2r})}+\| g\|_{L^2(\mathbb{R}^{d+1}_{2r}\cup\mathbb{R}^{d+1}_{-2r})}),\ \ \forall n\in\mathbb{N}.
			\end{align*}
		\end{lemma}
		\begin{proof}
			For $\tau\in(r,2r)$,  using integration by parts to obtain 
			\begin{equation}\label{ibp1}
				\begin{aligned}
					0&=\int_{R^{d+1}_\tau}\operatorname{div}_{y,w}(A(y)\nabla_{y,w}u- g)(y,w)\G_{x_0}(x,z,y,w)dydw\\
					&=\left. \int_{\mathbb{R}^d}e_{d+1}\cdot(A(y)\nabla_{y,w}u- g)(y,w)\G_{x_0}(x,z,y,w)dy\right|_{w=0}^{w=\tau}+\int_{R^{d+1}_\tau} g(y,w)\cdot\nabla_{y,w}\G_{x_0}(x,z,y,w)dydw\\
					&\quad\quad\quad-\int_{R^{d+1}_\tau}A(y)\nabla_{y,w}u(y,w)\cdot\nabla_{y,w}\G_{x_0}(x,z,y,w)dydw.
				\end{aligned}
			\end{equation}
			For the last term, we further have 
			\begin{equation}\label{ibp2}
				\begin{aligned}
					&-\int_{R^{d+1}_\tau}A(y)\nabla_{y,w}u(y,w)\cdot\nabla_{y,w}\G_{x_0}(x,z,y,w)dydw\\
					&=-\int_{R^{d+1}_\tau}A(x_0)\nabla_{y,w}u(y,w)\cdot\nabla_{y,w}\G_{x_0}(x,z,y,w)dydw\\
					&\quad\quad\quad+\int_{R^{d+1}_\tau}(A(x_0)-A(y))\nabla_{y,w}u(y,w)\cdot\nabla_{y,w}\G_{x_0}(x,z,y,w)dydw.
				\end{aligned}
			\end{equation}
			Here the last term is remainder term because the regularity of variable coefficient $A$ helps to compensate the singularity of kernel. For the first term, we integrate by parts again, and apply \eqref{Gbdy}, \eqref{eqG} to get
			\begin{align*}
				&-\int_{R^{d+1}_\tau}A(x_0)\nabla_{y,w}u(y,w)\cdot\nabla_{y,w}\G_{x_0}(x,z,y,w)dydw\\
				&\quad=\left.-\int_{\mathbb{R}^{d}}e_{d+1}\cdot(A(x_0)\nabla_{y,w}\G_{x_0}(x,z,y,w)u(y,w))dy\right|_{w=0}^{w=\tau}+(\mathbf{1}_{\mathbb{R}^{d+1}_\tau}u)(x,z)+(\mathbf{1}_{\mathbb{R}^{d+1}_\tau}u)(x,-z)\\
				&\quad=-\int_{\mathbb{R}^{d}}e_{d+1}\cdot(A(x_0)\nabla_{y,w}\G_{x_0}(x,z,y,\tau)u(y,\tau))dy+(\mathbf{1}_{\mathbb{R}^{d+1}_\tau}u)(x,z)+(\mathbf{1}_{\mathbb{R}^{d+1}_\tau}u)(x,-z).
			\end{align*}
			Combining this with \eqref{ibp1} and \eqref{ibp2} to obtain 
			\begin{equation}\label{u1+}
				\begin{aligned}
					&(\mathbf{1}_{\mathbb{R}^{d+1}_\tau}u)(x,z)+(\mathbf{1}_{\mathbb{R}^{d+1}_\tau}u)(x,-z)\\
					&=-\int_{R^{d+1}_\tau} g(y,w)\cdot\nabla_{y,w}\G_{x_0}(x,z,y,w)dydw-\int_{R^{d+1}_\tau}(A(x_0)-A(y))\nabla_{y,w}u(y,w)\cdot\nabla_{y,w}\G_{x_0}(x,z,y,w)dydw\\
					&\quad\quad +\int_{\mathbb{R}^{d}}e_{d+1}\cdot(A(x_0)\nabla_{y,w}\G_{x_0}(x,z,y,\tau)u(y,\tau))dy- \int_{\mathbb{R}^d}e_{d+1}\cdot(A(y)\nabla_{y,w}u- g)(y,\tau)\G_{x_0}(x,z,y,\tau)dy\\
					&\quad\quad+\int_{\mathbb{R}^d}e_{d+1}\cdot(A(y)\nabla_{y,w}u- g)(y,0^+)\G_{x_0}(x,z,y,0)dy.
				\end{aligned}
			\end{equation}
			Similarly, we get 
			\begin{equation}\label{u1-}
				\begin{aligned}
					&(\mathbf{1}_{\mathbb{R}^{d+1}_{-\tau}}u)(x,z)+(\mathbf{1}_{\mathbb{R}^{d+1}_{-\tau}}u)(x,-z)\\
					&=-\int_{R^{d+1}_{-\tau}} g(y,w)\cdot\nabla_{y,w}\G_{x_0}(x,z,y,w)dydw-\int_{R^{d+1}_{-\tau}}(A(x_0)-A(y))\nabla_{y,w}u(y,w)\cdot\nabla_{y,w}\G_{x_0}(x,z,y,w)dydw\\
					&\quad\quad -\int_{\mathbb{R}^{d}}e_{d+1}\cdot(A(x_0)\nabla_{y,w}\G_{x_0}(x,z,y,-\tau)u(y,-\tau))dy+ \int_{\mathbb{R}^d}e_{d+1}\cdot(A(y)\nabla_{y,w}u- g)(y,-\tau)\G_{x_0}(x,z,y,-\tau)dy\\
					&\quad\quad-\int_{\mathbb{R}^d}e_{d+1}\cdot(A(y)\nabla_{y,w}u- g)(y,0^-)\G_{x_0}(x,z,y,0)dy.
				\end{aligned}
			\end{equation}
			Then with the boundary condition $$
			e_{d+1}\cdot(A(x)(\mu_+^{-1}(\nabla_{x,z} u(x,0^+)-g(x,0^+))-\mu_-^{-1}\nabla_{x,z} (u(x,0^-)-g(x,0^-))))=0,
			$$
			we compute \eqref{u1+}$\times\mu_+^{-1}+$\eqref{u1-}$\times \mu_-^{-1}$ to obtain
			\begin{equation}\label{ueq1}
				\begin{aligned}
					&(\mu_+^{-1}u(x,z)+\mu_-^{-1}u(x,-z))\mathbf{1}_{\mathbb{R}^{d+1}_\tau}(x,z)+(\mu_-^{-1}u(x,z)+\mu_+^{-1}u(x,-z))\mathbf{1}_{\mathbb{R}^{d+1}_{-\tau}}(x,z)\\
					&=-\sum_{+,-}\left(\pm \mu_\pm^{-1}\int_{\mathbb{R}^d}e_{d+1}\cdot(A(y)\nabla_{y,w}u-g)(y,\pm\tau)\G_{x_0}(x,z,y,\pm\tau)dy\right)\\
					&\quad\quad\quad+\sum_{+,-}\left(\pm\mu_\pm^{-1} \int_{\mathbb{R}^d}e_{d+1}\cdot(A(x_0)\nabla_{y,w}\G(x,z,y,\pm\tau))u(y,\pm\tau)dy\right)\\
					&\quad\quad\quad-\sum_{+,-}\left(\mu_\pm^{-1}\int_{R^{d+1}_{\pm\tau}} g(y,w)\cdot\nabla_{y,w}\G_{x_0}(x,z,y,w)dydw\right)\\
					&\quad\quad\quad-\sum_{+,-}\left(\mu_\pm^{-1}\int_{R^{d+1}_{\pm\tau}} (A(x_0)-A(y))\nabla_{y,w}u(y,w)\cdot\nabla_{y,w}\G_{x_0}(x,z,y,w)dydw\right)\\
					&=-\sum_{+,-}\left(\mu_\pm^{-1}\int_{R^{d+1}_{\pm\tau}} g(y,w)\cdot\nabla_{y,w}\G_{x_0}(x,z,y,w)dydw\right)\\
					&\ \ -\sum_{+,-}\left(\mu_\pm^{-1}\int_{R^{d+1}_{\pm\tau}} (A(x_0)-A(y))\nabla_{y,w}u(y,w)\cdot\nabla_{y,w}\G_{x_0}(x,z,y,w)dydw\right)+\mathrm{LO}_\tau^1(u, g).
				\end{aligned}
			\end{equation}
			Here we denote $\mathrm{LO}_\tau^1(u, g)$ the boundary terms that only depend on $u,\nabla_{x,z}u, g$ on the boundary $\mathbb{R}^d\times \{-\tau,\tau\}$.\\
			On the other hand, using $\tilde {\G}_{x_0}$ as a test function, one has for 
			\begin{equation}\label{ibp11}
				\begin{aligned}
					0&=\int_{R^{d+1}_\tau}\operatorname{div}_{y,w}(A(y)\nabla_{y,w}u- g)(y,w)\tilde\G_{x_0}(x,z,y,w)dydw\\
					&=\left. \int_{\mathbb{R}^d}e_{d+1}\cdot(A(y)\nabla_{y,w}u- g)(y,w)\tilde\G_{x_0}(x,z,y,w)dy\right|_{w=0}^{w=\tau}+\int_{R^{d+1}_\tau} g(y,w)\cdot\nabla_{y,w}\tilde\G_{x_0}(x,z,y,w)dydw\\
					&\quad\quad\quad-\int_{R^{d+1}_\tau}A(y)\nabla_{y,w}u(y,w)\cdot\nabla_{y,w}\tilde\G_{x_0}(x,z,y,w)dydw.
				\end{aligned}
			\end{equation}
			Similar to \eqref{ibp2}, we further split the last term to the following two parts, 
			\begin{equation}\label{ibp22}
				\begin{aligned}
					&-\int_{R^{d+1}_\tau}A(y)\nabla_{y,w}u(y,w)\cdot\nabla_{y,w}\tilde\G_{x_0}(x,z,y,w)dydw\\
					&=-\int_{R^{d+1}_\tau}A(x_0)\nabla_{y,w}u(y,w)\cdot\nabla_{y,w}\tilde\G_{x_0}(x,z,y,w)dydw\\
					&\quad\quad\quad+\int_{R^{d+1}_\tau}(A(x_0)-A(y))\nabla_{y,w}u(y,w)\cdot\nabla_{y,w}\tilde\G_{x_0}(x,z,y,w)dydw.
				\end{aligned}
			\end{equation}
			For the first term, we integrate by parts again, and apply \eqref{Gbdy}, \eqref{eqtilG} to get
			\begin{align*}
				&-\int_{R^{d+1}_\tau}A(x_0)\nabla_{y,w}u(y,w)\cdot\nabla_{y,w}\tilde\G_{x_0}(x,z,y,w)dydw\\
				&=\left.-\int_{\mathbb{R}^{d}}e_{d+1}\cdot(A(x_0)\nabla_{y,w}\tilde\G_{x_0}(x,z,y,w)u(y,w))dy\right|_{w=0}^{w=\tau}+(\mathbf{1}_{\mathbb{R}^{d+1}_\tau}u)(x,z)+(\mathbf{1}_{\mathbb{R}^{d+1}_\tau}u)(x,-z)\\
				&=-\int_{\mathbb{R}^{d}}e_{d+1}\cdot(A(x_0)\nabla_{y,w}\tilde\G_{x_0}(x,z,y,\tau)u(y,\tau))dy+\int_{\mathbb{R}^{d}}e_{d+1}\cdot(A(x_0)\nabla_{y,w}\tilde\G_{x_0}(x,z,y,0)u(y,0^+))dy\\
				&\ \ \ +(\mathbf{1}_{\mathbb{R}^{d+1}_\tau}u)(x,z)-(\mathbf{1}_{\mathbb{R}^{d+1}_\tau}u)(x,-z).
			\end{align*}
			Combining this with \eqref{ibp11} and \eqref{ibp22} to obtain\begin{equation}
				\label{u2+}
				\begin{aligned}
					&(\mathbf{1}_{\mathbb{R}^{d+1}_\tau}u)(x,z)-(\mathbf{1}_{\mathbb{R}^{d+1}_\tau}u)(x,-z)\\
					&=-\int_{R^{d+1}_\tau} g(y,w)\cdot\nabla_{y,w}\tilde\G_{x_0}(x,z,y,w)dydw-\int_{R^{d+1}_\tau}(A(x_0)-A(y))\nabla_{y,w}u(y,w)\cdot\nabla_{y,w}\tilde\G_{x_0}(x,z,y,w)dydw\\
					&\ \ -\int_{\mathbb{R}^d}e_{d+1}\cdot(A(y)\nabla_{y,w}u- g)(y,\tau)\tilde\G_{x_0}(x,z,y,\tau)dy+\int_{\mathbb{R}^{d}}e_{d+1}\cdot(A(x_0)\nabla_{y,w}\tilde\G_{x_0}(x,z,y,\tau)u(y,\tau))dy\\
					&\ \ -\int_{\mathbb{R}^{d}}e_{d+1}\cdot(A(x_0)\nabla_{y,w}\tilde\G_{x_0}(x,z,y,0)u(y,0^+))dy.
				\end{aligned}
			\end{equation}
			Similarly, 
			\begin{equation}\label{u2-}
				\begin{aligned}
					&(\mathbf{1}_{\mathbb{R}^{d+1}_{-\tau}}u)(x,z)-(\mathbf{1}_{\mathbb{R}^{d+1}_{-\tau}}u)(x,-z)\\
					&=-\int_{R^{d+1}_{-\tau}} g(y,w)\cdot\nabla_{y,w}\tilde\G_{x_0}(x,z,y,w)dydw-\int_{R^{d+1}_{-\tau}}(A(x_0)-A(y))\nabla_{y,w}u(y,w)\cdot\nabla_{y,w}\tilde\G_{x_0}(x,z,y,w)dydw\\
					&\ \ +\int_{\mathbb{R}^d}e_{d+1}\cdot(A(y)\nabla_{y,w}u- g)(y,-\tau)\tilde\G_{x_0}(x,z,y,-\tau)dy-\int_{\mathbb{R}^{d}}e_{d+1}\cdot(A(x_0)\nabla_{y,w}\tilde\G_{x_0}(x,z,y,-\tau)u(y,-\tau))dy\\
					&\ \ +\int_{\mathbb{R}^{d}}e_{d+1}\cdot(A(x_0)\nabla_{y,w}\tilde\G_{x_0}(x,z,y,0)u(y,0^-))dy.
				\end{aligned}
			\end{equation}
			Taking a sum of \eqref{u2+} and \eqref{u2-}, and applying the boundary condition
			$u(y,0^+)=u(y,0^-)$, we obtain
			\begin{equation}\label{ueq2}
				\begin{aligned}
					&(u(x,z)-u(x,-z))\mathbf{1}_{\mathbb{R}^{d+1}_\tau\cup\mathbb{R}^{d+1}_{-\tau}}(x,z)\\
					&=-\int_{R^{d+1}_{\tau}\cup R^{d+1}_{-\tau}} g(y,w)\cdot\nabla_{y,w}\tilde\G_{x_0}(x,z,y,w)dydw+\int_{\mathbb{R}^{d}}e_{d+1}\cdot(A(x_0)\nabla_{y,w}\tilde\G_{x_0}(x,z,y,0)h(y))dy\\
					&\ \ -\int_{R^{d+1}_{\tau}\cup R^{d+1}_{-\tau}}(A(x_0)-A(y))\nabla_{y,w}u(y,w)\cdot\nabla_{y,w}\tilde\G_{x_0}(x,z,y,w)dydw+\mathrm{LO}_\tau^2(u,g),
				\end{aligned}
			\end{equation}
			where we put terms related to boundaries $\mathbb{R}^{d}\times \{-\tau,\tau\}$
			in $\mathrm{LO}_\tau^2(u,g)$.
			With
			\eqref{ueq1} and \eqref{ueq2}, we can obtain a representation formula of $u(x,z)$ for any $ (x,z)\in\mathbb{R}^{d+1}_{r}$ by solving the elementary equation
			\begin{align*}
				\begin{cases}
					&  \mu_+^{-1}u(x,z)+ \mu_-^{-1}u(x,-z)=B_1,\\
					&u(x,z)-u(x,-z)=B_2.
				\end{cases}   \Rightarrow \begin{cases}
					& u(x,z)=\frac{\mu_+\mu_-}{\mu_++\mu_-}B_1+\frac{\mu_+}{\mu_++\mu_-}B_2,\\
					&u(x,-z)=\frac{\mu_+\mu_-}{\mu_++\mu_-}B_1-\frac{\mu_-}{\mu_++\mu_-}B_2.
				\end{cases}
			\end{align*}
			Here $B_1,B_2$ are the right hand side terms in \eqref{ueq1} and \eqref{ueq2}, respectively,  which only relate to data $ g, h$ and boundary values $u|_{\mathbb{R}^d\times\{-\tau,\tau\}}$. Specially, after taking an average in $\tau\in[r,2r]$, one has for any $(x,z)\in\mathbb{R}^{d+1}_{r}$,
			\begin{align*}
				u(x,z)=&\frac{\mu_+}{\mu_++\mu_-}\left(\int_{\mathbb{R}^d}e_{d+1}\cdot(A(x_0)\nabla_{y,w}\tilde\G_{x_0}(x,z,y,0)h(y))dy\right.\\
				&\quad\quad\quad\quad\quad-\frac{1}{r}\int_{r}^{2r}\int_{R^{d+1}_{\tau}\cup R^{d+1}_{-\tau}} g(y,w)\cdot\nabla_{y,w}\tilde\G_{x_0}(x,z,y,w)dydwd\tau\\
				&\quad\quad\quad\quad\quad\ \left.-\frac{1}{r}\int_{r}^{2r}\int_{R^{d+1}_{\tau}\cup R^{d+1}_{-\tau}} ((A(x_0)-A(y))\nabla_{y,w}u(y,w))\cdot\nabla_{y,w}\tilde\G_{x_0}(x,z,y,w)dydwd\tau\right)\\
				&-\frac{\mu_+\mu_-}{\mu_++\mu_-}\sum_{+,-}\left(\frac{1}{r\mu_\pm}\int_{r}^{2r}\int_{R^{d+1}_{\pm\tau}} g(y,w)\cdot\nabla_{y,w}\G_{x_0}(x,z,y,w)dydwd\tau\right.\\
				&\quad\quad\quad\quad\quad\quad\left.+\frac{1}{r\mu_\pm}\int_{r}^{2r}\int_{R^{d+1}_{\pm\tau}} ((A(x_0)-A(y))\nabla_{y,w}u(y,w))\cdot\nabla_{y,w}\G_{x_0}(x,z,y,w)dydwd\tau\right)\\
				&+\mathrm{LO}(u, g),
			\end{align*}
			\begin{align*}
				u(x,-z)=&-\frac{\mu_-}{\mu_++\mu_-}\left(\int_{\mathbb{R}^d}e_{d+1}\cdot(A(x_0)\nabla_{y,w}\tilde\G_{x_0}(x,z,y,0)h(y))dy\right.\\
				&\quad\quad\quad\quad\quad-\frac{1}{r}\int_{r}^{2r}\int_{R^{d+1}_{\tau}\cup R^{d+1}_{-\tau}} g(y,w)\cdot\nabla_{y,w}\tilde\G_{x_0}(x,z,y,w)dydwd\tau\\
				&\quad\quad\quad\quad\quad\ \left.-\frac{1}{r}\int_{r}^{2r}\int_{R^{d+1}_{\tau}\cup R^{d+1}_{-\tau}} ((A(x_0)-A(y))\nabla_{y,w}u(y,w))\cdot\nabla_{y,w}\tilde\G_{x_0}(x,z,y,w)dydwd\tau\right)\\
				&-\frac{\mu_+\mu_-}{\mu_++\mu_-}\sum_{+,-}\left(\frac{1}{r\mu_\pm}\int_{r}^{2r}\int_{R^{d+1}_{\pm\tau}} g(y,w)\cdot\nabla_{y,w}\G_{x_0}(x,z,y,w)dydwd\tau\right.\\
				&\quad\quad\quad\quad\quad\quad\left.+\frac{1}{r\mu_\pm}\int_{r}^{2r}\int_{R^{d+1}_{\pm\tau}} ((A(x_0)-A(y))\nabla_{y,w}u(y,w))\cdot\nabla_{y,w}\G_{x_0}(x,z,y,w)dydwd\tau\right)\\
				&+\mathrm{LO}(u, g),
			\end{align*}
			here we denote 
			$\mathrm{LO}(u,g)$, which may vary from line to line, some lower order terms satisfying
			\begin{align*}
				\|\nabla_{x,z}^n\mathrm{LO}(u,g)\|_{L^\infty(\mathbb{R}^{d+1}_{r/2})}\lesssim C(r,n) (\|u\|_{H^1(\mathbb{R}^{d+1}_{2r}\cup \mathbb{R}^{d+1}_{-2r})}+\| g\|_{L^2(\mathbb{R}^{d+1}_{2r}\cup\mathbb{R}^{d+1}_{-2r})}),\ \ \forall n\in\mathbb{N}.
			\end{align*}
			Finally, observe that 
			\begin{align*}
				&\nabla_{y,w}\tilde \G_{x_0}(x,-z,y,w)=-\nabla_{y,w}\tilde \G_{x_0}(x,z,y,w),\ \ \ \nabla_{y,w} \G_{x_0}(x,-z,y,w)=\nabla_{y,w} \G_{x_0}(x,z,y,w).
			\end{align*}
			Hence, we obtain 
			\begin{align*}
				u(x,z)=&\frac{\mu_+\mathbf{1}_{z\in(0,r)}+\mu_-\mathbf{1}_{z\in(-r,0)}}{\mu_++\mu_-}\left(\int_{\mathbb{R}^d}e_{d+1}\cdot(A(x_0)\nabla_{y,w}\tilde\G_{x_0}(x,z,y,0)h(y))dy\right.\\
				&\quad\quad\quad\quad\quad-\frac{1}{r}\int_{r}^{2r}\int_{R^{d+1}_{\tau}\cup R^{d+1}_{-\tau}} g(y,w)\cdot\nabla_{y,w}\tilde\G_{x_0}(x,z,y,w)dydwd\tau\\
				&\quad\quad\quad\quad\quad\ \left.-\frac{1}{r}\int_{r}^{2r}\int_{R^{d+1}_{\tau}\cup R^{d+1}_{-\tau}} ((A(x_0)-A(y))\nabla_{y,w}u(y,w))\cdot\nabla_{y,w}\tilde\G_{x_0}(x,z,y,w)dydwd\tau\right)\\
				&-\frac{\mu_+\mu_-}{\mu_++\mu_-}\sum_{+,-}\left(\frac{1}{r\mu_\pm}\int_{r}^{2r}\int_{R^{d+1}_{\pm\tau}} g(y,w)\cdot\nabla_{y,w}\G_{x_0}(x,z,y,w)dydwd\tau\right.\\
				&\quad\quad\quad\quad\quad\quad\left.+\frac{1}{r\mu_\pm}\int_{r}^{2r}\int_{R^{d+1}_{\pm\tau}} ((A(x_0)-A(y))\nabla_{y,w}u(y,w))\cdot\nabla_{y,w}\G_{x_0}(x,z,y,w)dydwd\tau\right)\\
				&+\mathrm{LO}(u, g).
			\end{align*}
			This completes the proof.
		\end{proof}
		\begin{remark}\label{remF2}
			If we take $A(x)=\M(\nabla f(x))$, where $\M$ is defined in \eqref{ndefM}, then 
			\begin{align*}
				e_{d+1}\cdot A(x_0)\nabla_{y,w}\tilde\G_{x_0}(x,z,y,0)= K_\b(x-y,z)|_{\b=\nabla f(x_0)},
			\end{align*}
			where $K_\b$ is the kernel defined in \eqref{defK}.
		\end{remark}
		We have the following lemma considering the Lipschitz estimate of the solution in Lemma \ref{leminterf}.
		\begin{lemma}\label{lemLip}
			Let $u\in W^{1,2}(\mathbb{R}^{d+1}_{2r}\cup\mathbb{R}^{d+1}_{-2r})$ be a weak solution to \eqref{ellinte}. Then 
			\begin{align*}
				\|\nabla_{x,z}u\|_{L^\infty(\mathbb{R}^{d+1}_{r/2}\cup\mathbb{R}^{d+1}_{-r/2})}\lesssim& \|h\|_{ C^{1,\log^\varkappa}(\mathbb{R}^d)}+\|g\|_{\dot C^{\log^\varkappa}\cap L^2(\mathbb{R}^{d+1}_{2r}\cup\mathbb{R}^{d+1}_{-2r})}+c(r)\|\nabla_{x,z}u\|_{L^2(\mathbb{R}^{d+1}_{2r}\cup\mathbb{R}^{d+1}_{-2r})}\\
				&\ \ +C(r)\|\nabla_{x,z}u\|_{L^\infty(\mathbb{R}^{d+1}_{2r}\cup\mathbb{R}^{d+1}_{-2r})}.
			\end{align*}
			where $\varkappa>1$,  $c(r)=\log^{1-\varkappa}(2+r^{-1})$, $C(r)=r^{-\frac{2d+1}{2}}$.
		\end{lemma}
		\begin{proof}
			By \eqref{foruintf}, it suffices to consider 
			\begin{align*}
				&I_{1,x_0}=\int_{\mathbb{R}^d}e_{d+1}\cdot(A(x_0)\nabla_{y,w}\tilde\G_{x_0}(x,z,y,0)h(y))dy,\\
				&I_{2,x_0}^1=\int_{\mathbb{R}^{d+1}_{\pm\tau}} g(y,w)\cdot\nabla_{y,w}\G_{x_0}(x,z,y,w)dydw,\\ &I_{2,x_0}^2=\int_{\mathbb{R}^{d+1}_{\pm\tau}} g(y,w)\cdot\nabla_{y,w}\tilde \G_{x_0}(x,z,y,w)dydw, \\
				&I_{3,x_0}^1=\int_{\mathbb{R}^{d+1}_{\pm\tau}} ((A(x_0)-A(y))\nabla_{y,w}u(y,w))\cdot\nabla_{y,w}\G_{x_0}(x,z,y,w)dydw,\\ &I_{3,x_0}^2=\int_{\mathbb{R}^{d+1}_{\pm\tau}} ((A(x_0)-A(y))\nabla_{y,w}u(y,w))\cdot\nabla_{y,w}\bar\G_{x_0}(x,z,y,w)dydw,
			\end{align*}
			with $\tau\in(r,2r)$ and $(x,z)\in\mathbb{R}^{d+1}_\pm r/2$.
			By Remark \ref{remF2} and Lemma \ref{lemKf}, we obtain 
			\begin{align}\label{i1x0}
				\|\nabla_{x,z}I_{1,x_0}\|_{L^\infty}\lesssim \|h\|_{C^{1,\log^\varkappa}}.
			\end{align}
			Moreover, by the definition \eqref{defGtiG}, we have 
			\begin{align*}
				&|\nabla_{x,z}\nabla_{y,w}\G_{x_0}(x,z,y,w)|+|\nabla_{x,z}\nabla_{y,w}\tilde \G_{x_0}(x,z,y,w)|\lesssim \frac{1}{(|x-y|+|z-w|)^{d+1}},\\
				&\left|\int_{\mathbb{R}^{d+1}_{\pm \tau}}\nabla_{x,z}\nabla_{y,w}\G_{x_0}(x,z,y,w)dydw\right|+\left|\int_{\mathbb{R}^{d+1}_{\pm \tau}}\nabla_{x,z}\nabla_{y,w}\tilde \G_{x_0}(x,z,y,w)dydw\right|\lesssim C_r,
			\end{align*}
			for any $(x,z)\in\mathbb{R}^{d+1}_{\pm r/2}, \tau\in(r,2r)$.
			Hence, we obtain
			\begin{equation}\label{i2x0}
				\begin{aligned}
					& |\nabla_{x,z} I_{2,x_0}^1(x,z)|+ |\nabla_{x,z} I_{2,x_0}^2(x,z)|\\
					&\ \ \ \lesssim \int_{\mathbb{R}^{d+1}_{\pm \tau}} |g(y,w)-g(x,z)|(|\nabla_{x,z}\nabla_{y,w}\G_{x_0}(x,z,y,w)|+|\nabla_{x,z}\nabla_{y,w}\tilde \G_{x_0}(x,z,y,w)|)dydw+C(r) |g(x,z)|\\
					&\ \ \ \lesssim \|g\|_{\dot C^{\log ^\varkappa}\cap L^2(\mathbb{R}^{d+1}_{2r}\cup\mathbb{R}^{d+1}_{-2r})}.
				\end{aligned}
			\end{equation}
			Finally, for the remainder term, by taking $x_0=x$, we get 
			\begin{equation}\label{i3x0}
				\begin{aligned}
					&|(\nabla_{x,z}I_{3,x_0}^1)|_{x_0=x}|+|(\nabla_{x,z}I_{3,x_0}^2)|_{x_0=x}|\\
					&\ \ \ \lesssim \int_{\mathbb{R}^{d+1}_{\pm\tau}}\frac{|A(x)-A(y)||\nabla_{y,w}u(y,w)|}{(|x-y|+|z-w|)^{d+1}}dydw\\
					&\ \ \ \lesssim c(r)\|A\|_{ C^{\dot \log^\varkappa}}\|\nabla_{x,z}u\|_{  L^\infty(\mathbb{R}^{d+1}_{2r}\cup\mathbb{R}^{d+1}_{-2r})}+C(r)\|\nabla_{x,z}u\|_{  L^2(\mathbb{R}^{d+1}_{2r}\cup\mathbb{R}^{d+1}_{-2r})},
				\end{aligned}
			\end{equation}
			where $c(r)=\log^{1-\varkappa}(2+r^{-1})$, $C(r)=r^{-\frac{2d+1}{2}}$.
			Combining \eqref{i1x0}, \eqref{i2x0} and \eqref{i3x0}, we get 
			\begin{align*}
				&\|\nabla_{x,z}u\|_{L^\infty(\mathbb{R}^{d+1}_{r/2}\cup\mathbb{R}^{d+1}_{-r/2})}\\
				&\lesssim \sup_{x_0}\|\nabla_{x,z}I_{1,x_0}\|_{L^\infty(\mathbb{R}^{d+1}_{r/2}\cup\mathbb{R}^{d+1}_{-r/2})}+\sup_{x_0}\sum_{i=1,2}\|\nabla_{x,z}I_{2,x_0}^i\|_{L^\infty(\mathbb{R}^{d+1}_{r/2}\cup\mathbb{R}^{d+1}_{-r/2})}\\
				&\quad\quad\quad+\sum_{i=1,2}\|(\nabla_{x,z}I_{3,x_0}^i)|_{x_0=x}\|_{L^\infty(\mathbb{R}^{d+1}_{r/2}\cup\mathbb{R}^{d+1}_{-r/2})}\\
				&\lesssim  \|h\|_{ C^{1,\log^\varkappa}(\mathbb{R}^d)}+\|g\|_{\dot C^{\log^\varkappa}\cap L^2(\mathbb{R}^{d+1}_{2r}\cup\mathbb{R}^{d+1}_{-2r})}+c(r)\|\nabla_{x,z}u\|_{ L^\infty(\mathbb{R}^{d+1}_{2r}\cup\mathbb{R}^{d+1}_{-2r})}\\
				&\quad\quad\quad +C(r)\|\nabla_{x,z}u\|_{ L^\infty(\mathbb{R}^{d+1}_{2r}\cup\mathbb{R}^{d+1}_{-2r})}.
			\end{align*}
			This completes the proof.
		\end{proof}\vspace{0.3cm}\\
		Finally,		we present the following lemmas concerning estimates of singular integrals. 
		Let $K:\mathbb{R}^d\to \mathbb{R}$ such that 
		\begin{equation}\label{defczop}
			\begin{cases}
				&|K(x,y)|\leq C|x-y|^{-d},\ \ \ \forall x,y\in\mathbb{R}^d, x\neq y,\\
				&\int_{|x-y|\geq 10 |\alpha|}|\delta_\alpha^x K(x,y)||x-y|^a dy\lesssim |\alpha|^a,\ \ \ \forall \alpha\in \mathbb{R}^d, a\in(0,1).
			\end{cases}
		\end{equation}
		For any $h:\mathbb{R}^d\times \mathbb{R}^d\to \mathbb{R}$, define the integral operator $\T $ as following
		\begin{align}
			\T h(x,x'):=\mathrm{P.V.}\int_{\mathbb{R}^d}K(x,y) h(y,x') dy.
		\end{align}
		We have the following results.
		\begin{lemma}\label{lemsgit}
			Let $K$ satisfy \eqref{defczop} and 
			\begin{align}\label{kcanc}
				\sup_{0<r_1<r_2<\infty} \left|\int_{r_1<|x-y|<r_2}K(x,y)dy\right|\leq c_1,
			\end{align}
			for a constant $c_1\geq 0$.
			Then for $h:\mathbb{R}^d\times \mathbb{R}^d\to \mathbb{R}$, $a\in(0,1)$, $\tau>0$,
			\begin{align}
				&       \sup_x\sup_{|\alpha|\leq \tau} \frac{|\T h(x,x)-\T h(x-\alpha,x)|}{|\alpha|^a}\lesssim \frac{1}{a(1-a)}\lceil h\rfloor_{a,\tau},\label{holdTf},
			\end{align}
			where 
			$$
			\lceil h\rfloor_{a,\tau}:=    \sup_{x}\sup_{\substack{|\alpha|\leq 10\tau,\\|z-x|\leq |\alpha|}}\frac{| h(x,z)-h(x-\alpha,z)|}{|\alpha|^a}.
			$$
		\end{lemma}
		\begin{proof}
			For any $x_1,x_2\in\mathbb{R}^d$, denote $\delta=|x_1-x_2|$. By \eqref{defczop} and \eqref{kcanc} we have 
			\begin{align*}
				&	|\T h(x_1,x_1)-\T h(x_2,x_1)|\\
				&\leq \sum_{j=1,2}\left|\int_{|x_j-y|\leq 10\delta}K(x_j,y)(h(y,x_1)-h(x_j,x_1))dy\right|+ 2c_1|h(x_1,x_1)-h(x_2,x_1)|\\
				&+\int_{\mathbb{R}^d} \left|\mathbf{1}_{|x_1-y|>10\delta}K(x_1,y)-\mathbf{1}_{|x_2-y|>10\delta}K(x_2,y)\right||h(y,x_1)-h(x_1,x_1)|dy\\
				&\lesssim \left(\int_{y\leq 10\delta}\frac{|y|^a}{|y|^d}dy+\delta^a+\int_{|x_1-y|>10\delta}|K(x_1,y)-K(x_2,y)||y-x_1|^ady\right.\\
				&\quad\left.\quad\quad+\int_{9\delta\leq |x_1-y|,|x_2-y|\leq 10\delta} (|K(x_1,y)|+|K(x_2,y)|)|y-x_1|^ady\right)\lceil h\rfloor_a\\
				&\lesssim \frac{\delta^a\lceil h\rfloor_{a,\tau}}{a(1-a)},
			\end{align*}
			which leads to \eqref{holdTf}. 
			This completes the proof.
		\end{proof}
		
	\end{document}